\documentclass[reqno,twoside,a4paper]{amsbook}
\usepackage{amssymb}
\usepackage{latexsym}
\usepackage{amsmath}
\input{dumacros}	
\usepackage{curves}

\newdimen\ybox\ybox=5pt 
\setbox249=\vbox{\hbox{\vrule{\vbox to\ybox{}}\kern\ybox}\hrule}
\setbox248=\hbox{\copy249}

\newcount\n
\def\R#1,{%
\n=#1%
\nointerlineskip
\hbox{%
\loop\ifnum\n>0\copy248\advance\n by-1\repeat%
\vrule}%
}

\def\YD#1.{\,\vcenter{\hrule#1}\,}

\renewcommand {\Box}{\YD\R1,.}
\newcommand {\sym}{\YD\R2,.}
\newcommand {\symm}{\YD\R3,.}
\newcommand {\alt}{\YD\R1,\R1,.}
\newcommand {\altt}{\YD\R1,\R1,\R1,.}
\newcommand {\twone}{\YD\R2,\R1,.}

\newcommand {\Hbox}{\H_{\Box}}
\newcommand {\Hsym}{\H_{\sym}}
\newcommand {\Hsymm}{\H_{\symm}}
\newcommand {\Halt}{\H_{\alt}}
\newcommand {\Haltt}{\H_{\altt}}
\newcommand {\Htwone}{\H_{\twone}}

\newcommand {\IN}{\mathbb{N}}                          
\newcommand {\IZ}{\mathbb{Z}}

\newcommand {\IR}{\mathbb{R}}                           
\newcommand {\IC}{\mathbb{C}}
\newcommand {\T}{\mathbb T}

\newcommand {\A}{\mathcal A}
\newcommand {\B}{\mathcal B}
\newcommand {\D}{\mathcal D}
\newcommand {\E}{\mathcal E}
\newcommand {\F}{\mathcal F}
\newcommand {\G}{\mathcal G}
\renewcommand {\H}{\mathcal H}
\newcommand {\II}{\mathcal I}
\newcommand {\K}{\mathcal K}
\renewcommand {\L}{\mathcal L}
\newcommand {\Nu}{\mathcal V}
\newcommand {\U}{\mathfrak U}

\newcommand {\g}{\mathfrak{g}}
\newcommand {\gc}{\mathfrak{g}_{_{\IC}}}
\renewcommand {\t}{\mathfrak{t}}
\newcommand {\tc}{\t_{\IC}}                  
\newcommand {\Lie}{\operatorname{Lie}}

\newcommand {\half}[1]{\frac{#1}{2}}

\newcommand {\rot}{\operatorname{Rot}S^{1}}

\newcommand {\fuse}{\boxtimes}

\newtheorem {theorem}{Theorem}[subsection]
\newtheorem {proposition}[theorem]{Proposition}
\newtheorem {corollary}[theorem]{Corollary}
\newtheorem {lemma}[theorem]{Lemma}

\renewcommand {\thechapter}{\Roman{chapter}}
\newcommand {\ssection}[1]{\section{#1}\setcounter{theorem}{0}\setcounter{equation}{0}}
\newcommand {\ssubsection}[1]{\subsection{#1}\hfill\vspace{0.5\baselineskip}}
\renewcommand {\thetheorem}
{\ifcase\arabic{section}\else\arabic{section}.\ifcase\arabic{subsection}\else\arabic{subsection}.\fi\fi\arabic{theorem}}

\numberwithin {equation}{subsection}
\renewcommand {\theequation}
{\ifcase\arabic{section}\else\arabic{section}.\ifcase\arabic{subsection}\else\arabic{subsection}.\fi\fi\arabic{equation}}

\renewcommand {\proof}{{\sc Proof.}\ }

\newcommand {\remark}{{\sc Remark.\ }}
\newcommand {\definition}{{\sc Definition.\ }}
\newcommand {\halmos}{$\Diamond$}

\newcommand {\ltwo}[1]{L^{2}(S^{1},{#1})}
\newcommand {\cchi}{\raisebox{2pt}{$\chi$}}

\newcommand {\fermno}[1]{{_{\circ}^{\circ}}#1{_{\circ}^{\circ}}}

\newcommand {\ffin}{{\mathcal F}_{\operatorname{fin}}}
\newcommand {\lpol}{L^{\scriptscriptstyle{\operatorname{pol}}}}

\newcommand {\Hom}{\operatorname{Hom}}
\newcommand {\End}{\operatorname{End}}
\newcommand {\Ker}{\operatorname{Ker}}
\renewcommand {\Im}{\operatorname{Im}}
\newcommand {\Det}{\operatorname{Det}}
\newcommand {\rank}{\operatorname{rank}}
\newcommand {\<}{\langle}
\renewcommand {\>}{\rangle}
\newcommand {\sqnm}[1]{\<#1,#1\>}

\DeclareMathOperator*{\medsqcup}{\textstyle\bigsqcup}
\DeclareMathOperator*{\medsum}{\textstyle{\sum}}

\DeclareMathOperator*{\medplus}{\textstyle{\bigoplus}}

\newcommand {\lig}{\mathcal L_{I}G}
\newcommand {\licg}{\mathcal L_{I^{c}}G}
\newcommand {\Lig}{L_{I}G}
\newcommand {\Licg}{L_{I^{c}}G}
\newcommand {\Lis}{L_{I}\Spin_{2n}}
\newcommand {\Lics}{L_{I^{c}}\Spin_{2n}}

\newcommand {\XX}{\mathfrak X}
\newcommand {\YY}{\mathfrak Y}
\newcommand {\RR}{\mathcal R_{0}}

\newcommand {\pl}{\mathcal P_{\ell}}

\newcommand {\Vbox}{V_{\Box}}
\newcommand {\Vsym}{V_{\sym}}
\newcommand {\Vsymm}{V_{\symm}}
\newcommand {\Valt}{V_{\alt}}
\newcommand {\Valtt}{V_{\altt}}
\newcommand {\Vtwone}{V_{\twone}}

\renewcommand {\AA}{\mathcal A}

\newcommand {\fin}{^{\scriptscriptstyle{\operatorname{fin}}}}
\newcommand {\Hfin}{\H\fin}
\newcommand {\hfin}{\H\fin}
\newcommand {\hsmooth}{\H^{\infty}}
\newcommand {\Homl}{\Hom^{\ell}}

\newcommand {\vertex}[3]{\left(\begin{array}{c}#1\\#2\thickspace #3\end{array}\right)}
\newcommand {\field}[3]{\phi_{#2\medspace #3}^{#1}}

\newcommand {\diff}{\operatorname{Diff}(S^{1})}
\newcommand {\Exp}{\operatorname{Exp}}

\newcommand {\tad}{\A d}
\newcommand {\phs}[1]{\alpha_{\xi}(#1)}
\newcommand {\der}{\left.\frac{d}{dt}\right|_{t=0}}
\newcommand {\ders}{\left.\frac{d}{ds}\right|_{s=0}}
\newcommand {\C}[2]{\begin{pmatrix}#1\\#2\end{pmatrix}}
\newcommand {\norm}[1]{|#1|_{\half{3}}}
		
\renewcommand {\root}{\Lambda_{R}}
\newcommand {\weight}{\Lambda_{W}}
\newcommand {\domweight}{\weight^{+}}
\newcommand {\vvee}{^{\raisebox{.6pt}{$\scriptscriptstyle{\vee}$}}}
\newcommand {\coroot}{\Lambda_{R}\vvee}
\newcommand {\coweight}{\Lambda_{W}\vvee}

\newcommand {\Ie}{I^{0}}

\newcommand {\primed}[1]{{#1}^{\prime}}
\newcommand {\cor}[1]{\alpha_{#1}^{\vee}}	
\newcommand {\cow}[1]{\lambda_{#1}^{\vee}}	
\newcommand {\cowp}[1]{\primed{\cow{#1}}}	

\newcommand {\lzg}{L_{Z}G}

\newcommand {\Ad}{\operatorname{Ad}}
\newcommand {\ad}{\operatorname{ad}}

\newcommand {\wt}{\widetilde}
\newcommand {\wh}{\widehat}

\newcommand {\al}{\mathcal A_{\ell}}

\newcommand {\deltabar}{\overline{\Delta}}

\newcommand {\SU}{\operatorname{SU}}
\newcommand {\Spin}{\operatorname{Spin}}
\newcommand {\Sp}{\operatorname{Sp}}
\newcommand {\SO}{\operatorname{SO}}
\newcommand {\PSO}{\operatorname{PSO}}

\newcommand {\so}{\mathfrak{so}}

\input{diagrams}

\newcommand{\sizex}{480}
\newcommand{\sizey}{80}
\newcommand{\picturebegin}{\begin{picture}(\sizex,\sizey)(0,-10)}
\newcommand{\pictureend}{\end{picture}}

\def\axis{\thinlines
            \put(40,0){\line(100,0){400}}
            \put(40,0){\circle*{10}}\put(40,15)
               {\makebox(0,0){${\scriptscriptstyle \infty}$}}
            \put(140,0){\circle*{10}}\put(140,15)
	       {\makebox(0,0){${\scriptscriptstyle z}$}}
            \put(240,0){\circle*{10}}\put(240,15)
               {\makebox(0,0){${\scriptscriptstyle 0}$}}
            \put(340,0){\circle*{10}}\put(340,15)
               {\makebox(0,0){${\scriptscriptstyle 1}$}}
            \put(440,0){\circle*{10}}\put(440,15)
               {\makebox(0,0){${\scriptscriptstyle \infty}$}}}

\def\arc #1,#2,#3,#4{
  \countdef\x=50
  \x=100
  \multiply\x by #3
  \advance\x by 40
  \countdef\y=51
  \countdef\mid=52
  \mid=\x
  \countdef\z=53
  \countdef\height=54
  \ifcase#2 
    \y = 30
    \ifcase#1 
      \advance\mid by 40
      \z = 25
      \height = 47
      \scaleput(\mid,\height){${\textstyle t_{#4}}$}
      \scaleput(\mid,\z){${\scriptscriptstyle \blacktriangleleft}$}
    \or       
      \advance\mid by 42
      \z = 24
      \height = 40
      \scaleput(\mid,\height){${\textstyle t_{#4}}$}
      \scaleput(\mid,\z){${\scriptscriptstyle \blacktriangleright}$}
    \fi
  \or       
    \y = -30
    \ifcase#1 
      \advance\mid by 40
      \z = -34
      \height = -67
      \scaleput(\mid,\height){${\textstyle t_{#4}}$}
      \scaleput(\mid,\z){${\scriptscriptstyle \blacktriangleleft}$}
    \or       
      \advance\mid by 42
      \z = -33
      \height = -67
      \scaleput(\mid,\height){${\textstyle t_{#4}}$}
      \scaleput(\mid,\z){${\scriptscriptstyle \blacktriangleright}$}
    \fi
  \fi
  \put(\x,0){\curve(0,0,50,\y,100,0)}}

\def\leftuparc #1,#2{\arc 0,0,#1,#2}
\def\leftdownarc #1,#2{\arc 0,1,#1,#2}
\def\rightuparc #1,#2{\arc 1,0,#1,#2}
\def\rightdownarc #1,#2{\arc 1,1,#1,#2}

\def\I#1#2
  {\picturebegin\axis
  \ifcase#1 
    \ifcase#2
    \or \rightuparc 3,1\rightdownarc 3,2 
    \or \rightuparc 3,1\leftdownarc 1,2  
    \or \leftuparc  1,1\leftdownarc 1,2  
    \fi
  \or 
    \ifcase#2
    \or \leftuparc 0,1\leftdownarc 0,2 
    \or \leftuparc 2,1\leftdownarc 0,2 
    \or \leftuparc 2,1\leftdownarc 2,2 
    \fi
  \or 
    \ifcase#2
    \or \rightuparc 2,1\rightdownarc 2,2 
    \or \rightuparc 2,1\leftdownarc 0,2  
    \or \leftuparc 0,1\leftdownarc 0,2   
    \fi
  \fi
  \pictureend}

\begin{document}



\thispagestyle{empty}
\vspace*{2cm}

\begin{center}

{\huge\bf
Fusion of Positive Energy Representations of\\
LSpin$_{\mathbf{2n}}$}

\vskip 6cm

{\LARGE
Valerio Toledano Laredo\\[.5 ex]
St. John's College\\[.5 ex]
Cambridge}

\vskip 8 cm

{\Large
A thesis submitted for the degree of\\
Doctor of Philosophy\\
at the University of Cambridge\\
June 1997}

\end{center}
\newpage
\thispagestyle{empty}
\vspace*{2cm}
\hspace{.65\textwidth}
\begin{minipage}{.35\textwidth}
{\it
To Corinne,\\[1.2 ex]
and my parents\\[1.2 ex]
for their biblical patience}
\end{minipage}
\newpage
\thispagestyle{empty}

{\global\topskip 7.5pc\relax
  \begingroup
  \LARGE\bfseries\centering
  Preface
  \par \endgroup
  \vskip 1.2 cm}

The contents of this dissertation are original, except where explicit
reference is made to the work of others. The original material is the
result of my own work and includes nothing which is the outcome of
a collaborative effort. No part of this dissertation has been, or
is currently being submitted for a degree, diploma or qualification
at this or any other University.\\

I am deeply grateful to my research supervisor, Dr.~Antony Wassermann
for suggesting the problems discussed in this thesis, his continual
guidance and for sharing his ideas on conformal field theory and
operator algebras. I would also like to thank Professor Graeme Segal
for his interest in my work and many enlightening and stimulating
conversations, Hans Wenzl for sending me some unpublished notes on
braid group representations the use of which was central to the main
computation of this thesis and Alex Selby, for his expertise on local
systems at a time when I was struggling with the Knizhnik--Zamolodchikov
equations.\\

I am indebted in countless ways to Giovanni Forni and Alessio
Corti. It is a pleasure to thank them both for their friendship,
advice and much needed encouragement. Grazie.\\

The pleasant working conditions in the Department of Pure Mathematics
were made even more so by the company of Jochen Theis, to whom
I am grateful for putting up with me over so many lunches, and
the assistance of Sally Lowe, whom I thank for her friendly
help in all matters administrative.\\

This research could not have been carried out without the generous
financial support of the European Commission through the Human Capital and
Mobility and the Travel and Mobility of Researchers programs. I also
gratefully acknowledge a Science and Engineering Research Council
fees award and financial assistance from St.~John's College during
the tenth term of my research.

\vskip10ex
\hspace{0.5\textwidth}
Valerio Toledano Laredo
\thispagestyle{empty}
\newpage


\setcounter{tocdepth}{1}
\tableofcontents







\chapter*{Introduction}

Let $G$ be a compact, connected and simply--connected simple Lie group.
The {\it loop group} $LG=C^{\infty}(S^{1},G)$ admits a class of projective
unitary representations, those of {\it positive energy}, whose theory
parallels that of the representations of $G$. In particular, they are
completely reducible and the irreducible ones are described by the same
data as those of $G$ together with an additional parameter, the {\it level}
$\ell\in\IN$ which classifies the underlying projective cocycle. Morevoer,
only finitely many irreducibles exist at a given level $\ell$.\\

In their pioneering work, Belavin--Polyakov--Zamolodchikov \cite{BPZ} and
Knizhnik--Zamolodchikov \cite{KZ} associate to each positive energy
representation a quantum {\it primary field}. They introduce the notion of
{\it fusion} of primary fields, based on their short range
{\it operator product expansion}. Their computations implicitly suggest
the existence of a tensor product operation, differing from the usual one,
between positive energy representations at a fixed level which should give
them a structure similar to the representation ring of a finite group.
The fundamental problem of giving a mathematically sound definition of this
product which is associative and reflects the computations of the physicists
has proved to be difficult so far. Several independent, and at present
unrelated definitions have been suggested by Segal, Kazhdan--Lusztig,
Jones--Wassermann, Borcherds, Beilinson--Drinfeld and Huang--Lepowsky.\\

The aim of the present dissertation is to give such a definition for the loop
group of $G=\Spin_{2n}$, $n\geq 3$ and to characterise the resulting algebraic
structure on the category $\pl$ of positive energy representations at level
$\ell$. Conjecturally, the latter is described by the {\it Verlinde rules}
\cite{Ve}, that is as a quotient of the representation ring of $\Spin_{2n}$ and
we establish a number of results in agreement with this. In particular, we compute
explicitly the fusion of the vector representation of $L\Spin_{2n}$, {\it i.e.}~the
one corresponding to the defining representation of $\SO_{2n}$ with the positive
energy representations corresponding to single--valued representations of $\SO_{2n}$
and, at odd level with all representations. Moreover, we prove that the level 1
fusion ring of $L\Spin_{2n}$ is isomorphic to the group algebra of the centre of
$\Spin_{2n}$. As A.~Wassermann has recently indicated, our techniques should
extend to yield the complete structure of the fusion ring of $L\Spin_{2n}$.
Together with a study of the category theoretic properties of $\pl$, this would
lead to the first rigorous definition of the three--manifold invariants defined
by link surgery corresponding to the group $\Spin_{2n}$.\\

Our study is modelled on the von Neumann algebra approach to fusion, first
introduced by Wassermann in relation to the fusion ring of $L\SU_{n}$
\cite{Wa1,Wa2,Wa3,Wa4,J4}. The latter originates in a joint program with V. Jones
aimed at understanding unitary conformal field theories from the point of view of
operator algebras which arose from Jones' suggestion that there might be a subfactor
explanation for the coincidence of certain braid group representations that
appeared in statistical mechanics and conformal field theory \cite{J2,J3}.
It relies on the use of {\it Connes fusion}, a tensor product operation on bimodules
over von Neumann algebras originally developed by Connes \cite{Co,Sa}. The link with
loop groups arises by regarding positive energy representations as bimodules over the
groups $L_{I}\Spin_{2n},L_{I^{c}}\Spin_{2n}$ of loops supported in a proper interval
$I\subset S^{1}$ and its complement and leads to a definition of fusion which is
manifestly unitary and associative.
The use of Connes fusion requires one to check that positive energy representations
at a given level satisfy the axioms of locality, Haag duality and local equivalence.
This amounts in essence to showing that $\pl$ constitutes a quantum field theory in
the sense of Doplicher--Haag--Roberts and provides interesting examples of such
theories in 1+1 dimensions since, as pointed out by Goddard, Nahm and Olive \cite{GNO}
they are distinct from free field theories for all but finitely many values of $\ell$.\\

Wassermann showed that the intertwiners necessary to compute Connes fusion could
be obtained by smearing the primary fields of Knizhnik and Zamolodchikov,
provided their continuity as operator--valued distributions was established. With
this analysis done, the actual computation of fusion depends on the commutation
or {\it braiding} properties of the primary fields through {\it braiding--fusion
duality}. More precisely, as established in \cite{KZ,TK1}, the structure constants
governing the braiding of primary fields arise as the entries of the analytic continuation
matrix of a vector--valued fuchsian ODE, the {\it Knizhnik--Zamolodchikov equation},
from the singular point $0$ to $\infty$. The computability of fusion therefore rests
on the solvability of this ODE. Wassermann's method was successfully applied by Loke
\cite{Lo} to the positive energy representations of the diffeomorphism group of the
circle.\\

Another important input in this study is the theory of superselection sectors
developed within the context of algebraic quantum field theory by Doplicher,
Haag and Roberts \cite{DHR1,DHR2}.
After the laborious analysis involved in explicitly checking a few fusion rules,
it follows fairly easily that the properties required to apply the DHR theory
hold. One is therefore guaranteed the existence of a quantum or statistical
dimension, of a Markov trace and an action of the braid group whose use plays
an important r\^ole in deducing the remaining fusion rules.\\

As mentioned above, an important precedent for the present study is Wassermann's
computation of the fusion ring of $L\SU_{n}$. This constituted an invaluable
conceptual framework and the backbone on which this dissertation rests. We wish
however to stress that, much
as the invariant theory of $\Spin_{2n}$ differs from that of $\SU_{n}$, the study
of the loop groups of the spin groups poses a whole new layer of conceptual and
technical difficulties. We sketch them here, in technical jargon and refer to the
outline of our results given in the next section or to the main text for more
details.\\

With a degree of oversimplification, one may say that most analytic results needed for
the positive energy representations of $L\SU_{n}$ can be derived from their fermionic
construction. This realises the level 1
representations as summands of the basic representation of $LU_{n}$ obtained from that
of the Clifford algebra of $\H=L^{2}(S^{1},\IC^{n})$ on the Fock space $\F=\Lambda\H_{+}$
of the Hardy space $\H_{+}$. The level $\ell$ representations are then realised inside
the $\ell$--fold tensor product $\F^{\otimes\ell}$, itself a Fock space representation
of the Clifford algebra of $L^{2}(S^{1},\IC^{n\ell})$. This uniform construction simplifies
the proof of the following points

\begin{enumerate}
\item All positive energy representations of $L\SU_{n}$ at level $\ell$ give rise to
the same projective cocycle. This is obvious since they are summands of a single
representation.
\item Positive energy representations at equal level are unitarily equivalent for the
local loop groups $L_{I}\SU_{n}$. This follows from the (non--trivial) fact that these
generate type III factors in $\F^{\otimes\ell}$ and that all projections in such
factors are equivalent.
\item All primary fields extend to unbounded operator--valued distributions. This simply
follows because they may be realised as products of Fermi fields compressed by projections,
although checking that this procedure yields all primary fields is a subtler matter.
\end{enumerate}

By contrast,
\begin{enumerate}
\item The level 1 representations of $L\Spin_{2n}$ appear, grouped by two as
summands of two distinct Fock spaces, the Neveu--Schwarz and the Ramond sector. 
This separation obscures (i) above, the proof of which requires the use of the spin primary
field which links the two sectors, or of analytic vector techniques.
\item The local equivalence at level 1 must be established by the use of the outer
automorphic action of the centre of $\Spin_{2n}$ on $L\Spin_{2n}$ via conjugation by
discontinuous loops.
\item The fermionic construction for $L\Spin_{2n}$ only gives information on
the vector primary field but little or none on the spin fields.
We remedy this problem by making use of the vertex operator model of Segal and
Kac--Frenkel to give an explicit construction of all the level 1 primary fields
in the bosonic picture. A technique due to Wassermann
then allows to establish their continuity properties as well as those of a number of
higher level primary fields.\\
\end{enumerate}

Another substantial difference lies in the study of the relevant Knizhnik--Zamolodchikov
equations. As well--known, solutions of the KZ equations for all simple Lie algebras have
been given by Schechtman and Varchenko \cite{SV} but their combinatorial complexity makes
them intractable for computational purposes. Somewhat miraculously, the computation of
the fusion of a general representation of $L\SU_{n}$ with the vector representation may
be obtained from a detailed study of the monodromy properties of the generalised
hypergeometric function studied by Thomae at the end of the 19th century \cite{Wa3}. No
such classical treatment was available to us, but a similar
miracle occurs for $L\Spin_{2n}$. The 3rd order ODE required to compute the fusion of the
vector representation with its symmetric powers reduces to the Dotsenko--Fateev equation
\cite{DF} discovered by the physicists in the context of minimal models.
The analysis of the first part of the dissertation then enables us to establish a small
number of explicit fusion rules together with upper bounds for the multiplicities of the
general fusion rules with the vector representation. The computation is concluded by
resorting to the algebraic techniques of quantum invariant theory, and in particular to
the action of the braid group coming from the Doplicher--Haag--Roberts theory. The
corresponding representations can be handled using some algebraic techniques of Wenzl
connected with the Birman--Wenzl algebra \cite{We2,We3}. The fusion rules are finally
pinned down by interpreting the quantum dimension of DHR in the context of
Perron--Frobenius theory.

\vskip 2 cm

\ssection{Outline of contents}

We outline below the contents of this dissertation, developing only what is minimally
necessary to sketch a proof of our main results. A detailed account of each chapter
may be found in its opening paragraphs. \S \ref{ss:representations} describes the
classification of positive energy representations of the loop group of $\Spin_{2n}$.
In \S \ref{ss:connes fusion}, we define the Connes fusion of two such representations
and indicate a scheme to compute it based upon the commutation properties of certain
intertwiners. The latter are derived in \S \ref{ss:primary fields} from the braiding
properties of primary fields. These are the building blocks of conformal field theory
and we outline their definition and main properties. Our main results are given in \S
\ref{ss:main results}.

\ssubsection{Positive energy representations of $LG$}\label{ss:representations}

Let $G=\Spin_{2n}$, $n\geq 3$ be the universal covering group of $\SO_{2n}$
and $\g=\so_{2n}$ its Lie algebra. The irreducible representations of $G$
are classified by their highest weight, a sequence
$\zeta_{1}\geq\cdots\geq\zeta_{n-1}\geq|\zeta_{n}|$ with $\zeta_{i}\in\IZ$
for any $i$ or $\zeta_{i}\in\half{1}+\IZ$ for any $i$ corresponding
respectively to the single and two--valued representations of $\SO_{2n}$.\\

Consider now the {\it loop group} $LG=C^{\infty}(S^{1},G)$ of $G$ with Lie
algebra $L\g=C^{\infty}(S^{1},\g)$. Both are acted upon by $\rot$, the group
of rotations of the circle, via reparametrisation. The representations of
$LG$ have a similar classification, if one restricts one's attention to
those of {\it positive energy}.
These are projective unitary representations $\pi:LG\rightarrow PU(\H)$
extending to the semi--direct product $LG\rtimes\rot$ in such a way that
the infinitesimal generator of rotations is bounded below and has
finite--dimensional eigenspaces. In other words, $\H=\bigoplus_{n\geq 0}\H(n)$
where $\H(n)=\{\xi\in\H|\pi(R_{\theta})\xi=e^{in\theta}\xi\}$, the {\it subspace
of energy $n$}, supports a finite--dimensional representation of the subgroup
of constant loops $G\subset LG$.\\

Positive energy representations are completely reducible and their classification
is obtained in the following way. As explained in chapter \ref{ch:classification},
the {\it finite energy subspace} $\Hfin$ of a positive energy representation $\H$,
that is the algebraic direct sum $\bigoplus_{n\geq 0}\H(n)$ supports a projective
representation of $\lpol\g$, the dense subalgebra of $L\g$ consisting of
$\g$--valued trigonometric polynomials, and
therefore one of the affine Kac--Moody algebra $\wh{\g}$ corresponding to $\g$.
The latter is the semi--direct product $\wt\g\rtimes\IC d$ of the universal central
extension $\wt\g$ of $\lpol\g$ by the infinitesimal action of rotations.
Moreover, the classification of $\H$ as an $LG$--module is equivalent to that
of the $\wt\g$--module $\Hfin$. In particular, an irreducible $\H$ is uniquely
determined by an integer $\ell$ called the {\it level}, which classifies the Lie
algebra cocyle on $\lpol\g$ and its {\it lowest energy subspace} $\H(0)$, an
irreducible $G$--module whose highest weight $\zeta$ is bound by the requirement
that $\zeta_{1}+\zeta_{2}\leq\ell$. An irreducible $G$--module satisfying this
constraint is called {\it admissible at level $\ell$} and always arises as the
lowest energy subspace of a (necessarily unique) irreducible positive energy
representation at level $\ell$.\\

The centre $Z(G)$ of $G$ acts by outer automorphisms on $LG$ via conjugation
by {\it discontinuous loops}, {\it i.e.}~lifts to $G$ of loops in $G/Z(G)$. As
shown in chapter \ref{ch:classification}, this induces a level preserving
action of $Z(G)$ on the positive energy representations of $LG$. In fact,
if $z\in Z(G)$ and $\H$ is an irreducible positive energy representation,
the lowest energy subspace of the conjugated representation $z\H$ may be
characterised explicitly in terms of the level $\ell$, the lowest energy
subspace $\H(0)$ and $z$ by realising $Z(G)$ as a distinguished subgroup
of the automorphisms of the extended Dynkin diagram of $G$. This shows in
particular that at level 1, the action of $Z(G)$ is transitive and free.\\


Chapter \ref{ch:analytic} is devoted to the analytic properties of positive energy
representations and in particular to the construction of a dense subspace on which
both $LG$ and its Lie algebra $L\g$ act, which is required as a natural domain for
the smeared primary fields. These are densely defined $LG$--intertwiners mapping 
between different positive energy representations and do not in general extend to
bounded maps. If $\H$ is a positive energy representation, such a domain is provided
by the subspace of {\it smooth vectors} $\hsmooth$ for $\rot$, {\it i.e.}~those vectors
whose projection onto the subspace of energy $n$ decreases in norm faster than any
polynomial in $n$.

\ssubsection{Connes fusion of positive energy representations}\label{ss:connes fusion}

The notion of {\it Connes fusion} of positive energy representations of $LG$ arises
by regarding them as bimodules over the subgroups $L_{I}G,L_{I^{c}}G$ of loops
supported in a given interval $I\subset S^{1}$ and its complement and using a
tensor product operation on bimodules over von Neumann algebras developed by
Connes \cite{Co,Sa}.\\

Recall that a bimodule $\H$ over a pair $(M,N)$ of von Neumann algebras is a Hilbert
space supporting commuting representations of $M$ and $N$. To any two bimodules $X,Y$
over the pairs $(M,N)$, $(\wt N,P)$, Connes fusion associates an $(M,P)$--bimodule
denoted by $X\boxtimes Y$. The definition of $X\boxtimes Y$ relies on, but is ultimately
independent of the choice of a reference or {\it vacuum} $(N,\wt N)$--bimodule $\Nu$
with a cyclic vector $\Omega$ for both actions and for which {\it Haag duality} holds,
{\it i.e.}~the actions of $N$ and $\wt N$ are each other's commutant. Given $\Nu$, we form
the intertwiner spaces
\begin{xalignat}{3}
\XX&=\Hom_{N}(\Nu,X)&
&\text{and}&
\YY&=\Hom_{\wt N}(\Nu,Y)
\end{xalignat}
and consider the sequilinear form on the algebraic tensor product $\XX\otimes\YY$
given by
\begin{equation}\label{inner product}
\<x_{1}\otimes y_{1},x_{2}\otimes y_{2}\>=
(x_{2}^{*}x_{1} y_{2}^{*}y_{1}\Omega,\Omega)
\end{equation}
where the inner product on the right hand--side is taken in $\Nu$. If $x_{1}=x_{2}$
and $y_{1}=y_{2}$, Haag duality implies that $x_{2}^{*}x_{1}$ and $y_{2}^{*}y_{1}$
are commuting positive operators and therefore that $\<\cdot,\cdot\>$ is positive
semi--definite. By definition, the bimodule
$X\boxtimes Y$ is the Hilbert space completion of $\XX\otimes\YY$ with respect
to $\<\cdot,\cdot\>$, with $(M,P)$ acting as $(m,p)x\otimes y=mx\otimes py$.\\

Applying the above to the positive energy representations of $LG$ requires a number
of preliminary results which are established in chapters \ref{ch:fermionic} and
\ref{ch:loc loops}. Let $\pl$ be the set of positive energy representations at a
fixed level $\ell$. We wish to regard any $(\H,\pi)\in\pl$ as a bimodule over the
pair $\pi_{0}(L_{I}G)^{''},\pi_{0}(L_{I^{c}}G)^{''}$ where $\pi_{0}$ is the
{\it vacuum representation} at level $\ell$ whose lowest energy subspace is, by
definition the trivial $G$--module. The well--foundedness of this change of
perspective is justified by the following properties

\begin{enumerate}
\item {\it Locality} :
$\pi(L_{I}G)''\subset\pi(L_{I^{c}}G)'$ for any $(\pi,\H)\in\pl$. In other words,
$\H$ is a $(\pi(L_{I}G)^{''},\pi(L_{I^{c}}G)^{''})$--bimodule.

\item {\it Local equivalence} :
All $(\pi,\H)\in\pl$ are unitarily equivalent as $L_{I}G$--modules. Thus we may
unambiguously identify $\pi(L_{I}G)''$ with $\pi_{0}(L_{I}G)''$ and consider $\H$
as a $(\pi_{0}(L_{I}G)'',\pi_{0}(L_{I^{c}}G)'')$--bimodule.

\item {\it von Neumann density} :
$\pi(L_{I}G)\times\pi(L_{I^{c}}G)$ is strongly dense in $\pi(LG)$. Thus, inequivalent
irreducible positive energy representations of $LG$ remain so when regarded as bimodules.
\end{enumerate}

The r\^ole of the reference bimodule is played by the vacuum representation
$(\pi_{0},\H_{0})$. Two crucial facts need to be established in this respect

\begin{enumerate}
\item[(iv)] {\it Reeh-Schlieder theorem} : any finite energy vector of a positive
energy representation $\pi$ is cyclic under $\pi(L_{I}G)$. In particular, the
lowest energy vector $\Omega\in\H_{0}(0)$ is cyclic for $\pi_{0}(L_{I}G)''$
and $\pi_{0}(L_{I^{c}}G)''$.
\item[(v)] {\it Haag duality} : $\pi_{0}(L_{I}G)^{''}=\pi_{0}(L_{I^{c}}G)^{'}$.
\end{enumerate}

Finally, another technically crucial property of the algebras $\pi(L_{I}G)''$ is the
following
\begin{enumerate}
\item[(vi)] {\it Factoriality} : The algebras $\pi(L_{I}G)''$ with $I\subsetneq S^{1}$
are type III$_{1}$ factors.
\end{enumerate}

Once properties (i)-(vi) are established, the definition of Connes fusion may be
adapted to our setting. Let $\H_{i},\H_{j}\in\pl$ and form the intertwiner spaces
\begin{xalignat}{2}\label{intertwiner spaces}
\XX_{i}&=\Hom_{\Licg}(\H_{0},\H_{i})&
\YY_{j}&=\Hom_{\Lig}(\H_{0},\H_{j})
\end{xalignat}
Then, $\H_{i}\fuse\H_{j}$ is the completion of $\XX_{i}\otimes\YY_{j}$ with respect
to the inner product \eqref{inner product}. $\H_{i}\fuse\H_{j}$ is manifestly unitary
and the functorial properties of Connes fusion imply that $\fuse$ is associative.
It is not {\it a priori} clear however that $\H_{i}\fuse\H_{j}$ is of positive energy.
In fact there is no naturally defined action of $LG$ on it, let alone an intertwining
one of $\rot$. These facts will be checked a posteriori by explicitly computing the
fusion.\\

Consider now a first attempt at computing $\H_{i}\fuse\H_{j}$. Notice that, if
$y_{j0}\in\Hom_{\Lig}(\H_{0},\H_{j})$, then
\begin{equation}
y^{*}_{j0}y_{j0}\in\Hom_{\Lig}(\H_{0},\H_{0})=\pi_{0}(\Lig)'
\end{equation}
By Haag duality, $y_{j0}^{*}y_{j0}$ lies in $\pi_{0}(\Licg)''$ and may therefore
be represented, via local equivalence, on $\H_{i}$. Thus, if $x_{i0}\in\XX_{i}$
so that $x_{i0}\pi_{0}(\gamma)=\pi_{i}(\gamma)x_{i0}$ for any $\gamma\in\Licg$,
then
\begin{equation}\label{eq:transport -}
x_{i0}y_{j0}^{*}y_{j0}=\pi_{i}(y_{j0}^{*}y_{j0})x_{i0}
\end{equation}
Much of our efforts will be aimed at establishing a {\it transport formula} giving
an explicit expression for $\pi_{i}(y_{j0}^{*}y_{j0})$, that is the identity
\begin{equation}\label{transport}
\pi_{i}(y_{j0}^{*}y_{j0})=\sum_{k}\beta_{k}y_{ki}^{*}y_{ki}
\end{equation}
where the $\beta_{k}$ are some positive constants labelled by inequivalent positive
energy representations $\pi_{k}\in\pl$ and the $y_{ki}\in\Hom_{\Lig}(\H_{i},\H_{k})$
depend linearly on $y_{j0}$. Deferring to \S \ref{ss:primary fields} an explanation
of why \eqref{transport} should hold at all, notice that it allows to compute
$\H_{i}\fuse\H_{j}$. Indeed, by \eqref{eq:transport -} and \eqref{transport}
\begin{equation}
\begin{split}
\|x_{i0}\otimes y_{j0}\|^{2}
&=(x_{i0}^{*}x_{i0}y_{j0}^{*}y_{j0}\Omega,\Omega)\\[1.2 ex]
&=(x_{i0}^{*}\pi_{i}(y_{j0}^{*}y_{j0})x_{i0}\Omega,\Omega)\\[1.2 ex]
&=\sum_{k}\beta_{k}(x_{i0}^{*}y_{ki}^{*}y_{ki}x_{i0}\Omega,\Omega)\\
&=\|\bigoplus_{k}\beta_{k}^{\half{1}}y_{ki}x_{i0}\Omega\|^{2}
\end{split}
\end{equation}
and therefore the map
\begin{equation}
U:\XX_{i}\otimes\YY_{j}\rightarrow\bigoplus_{k}\H_{k},
\thickspace
x_{i0}\otimes y_{j0}\rightarrow
\bigoplus_{k}\beta_{k}^{\half{1}}y_{ki}x_{i0}\Omega
\end{equation}
extends to an isometry $\H_{i}\fuse\H_{j}\rightarrow\bigoplus_{k}\H_{k}$ which is
easily seen to be $\Lig\times\Licg$--equivariant. Since the $\H_{k}$ are irreducible
and inequivalent $\Lig\times\Licg$--modules, and the image of $U$ has non--zero
intersection with each of the $\H_{k}$, it is dense by Shur's lemma and we conclude
that
\begin{equation}
\H_{i}\fuse\H_{j}=\bigoplus_{k}\H_{k}
\end{equation}
In particular, $\H_{i}\fuse\H_{j}$ is of positive energy or, more precisely, arises
as the restriction to $\Lig\times\Licg$ of a positive energy $LG$--module
\footnote{In the foregoing, we have implicitly restricted our attention to the case
where the transport formula does not involve any multiplicities and therefore the
summands of $\H_{i}\fuse\H_{j}$ have multiplicity one. The general case follows in
a similar way.}.\\

To conclude, notice that the transport formula only needs to be established for one
unitary $y_{j0}\in\YY_{j}$. Such an element certainly exists by local equivalence.
Moreover, by Haag duality,
\begin{equation}
\YY_{j}=\Hom_{\Lig}(\H_{0},\H_{j})=y_{j0}\pi_{0}(\Lig)'=y_{j0}\pi_{0}(\Licg)''
\end{equation}
Thus, if \eqref{transport} holds for $y_{j0}$ and $a$ lies in the *--algebra generated
by $\pi_{0}(\Licg)$, we have for $y=y_{j0}a$
\begin{equation}
\pi_{i}(y^{*}y)=
\pi_{i}(a^{*}y_{j0}^{*}y_{j0}a)=
\pi_{i}(a)^{*}\pi_{i}(y_{j0}^{*}y_{j0})\pi_{i}(a)=
\sum_{k}\beta_{k}
\pi_{i}(a)^{*}y_{ki}^{*}y_{ki}\pi_{i}(a)
\end{equation}
and the transport formula holds for $y$. 

\ssubsection{Primary fields and their braiding properties}\label{ss:primary fields}

Let $\H_{i},\H_{j}$ be irreducible positive energy representations at level $\ell$
with lowest energy subspaces $V_{i},V_{j}$. Let $V_{k}$ be an irreducible
$G$--module admissible at level $\ell$ and $V_{k}[z,z^{-1}]$ the space of
$V_{k}$--valued trigonometric polynomials. $\lpol\g$ acts on $V_{k}[z,z^{-1}]$ by
multiplication and $\rot$ by $R_{\theta} v\otimes z^{n}=e^{-in\theta}v\otimes z^{n}$.
A {\it primary field} of charge $V_{k}$ is a linear map
\begin{equation}
\phi_{ji}^{k}:\Hfin_{i}\otimes V_{k}[z,z^{-1}]\rightarrow\Hfin_{j}
\end{equation}
intertwining the actions of $\lpol\g\rtimes\rot$. $\phi_{ji}^{k}$ may be regarded as
an endomorphism--valued algebraic distribution associating to any $f\in V_{k}[z,z^{-1}]$
the {\it smeared field}
$\phi_{ji}^{k}(f)\in\Hom(\Hfin_{i},\Hfin_{j})$ and is best represented by its generating
function
\begin{equation}
\phi_{ji}^{k}(z)=
\phi_{ji}^{k}(v,z)=\sum_{n\in\IN}\phi_{ji}^{k}(v,n)z^{-n}
\end{equation}
where $\phi_{ji}^{k}(v,n)=\phi_{ji}^{k}(v\otimes z^{n})$, $v\in V_{k}$
\footnote{For the {\it cognoscenti}, we shall abusively denote by $\phi_{ji}^{k}(z)$
both the integrally moded primary field and $\sum_{n\in\IZ}\phi_{ji}^{k}(v,n)z^{-n-\Delta}$
where $\Delta$ is the conformal dimension of $\phi_{ji}^{k}$.}.
By restriction, $\phi_{ji}^{k}$ defines a finite--dimensional $G$--intertwiner or {\it
initial term} $\varphi:V_{i}\otimes V_{k}\rightarrow V_{j}$ which determines $\phi_{ji}^{k}$
uniquely. Thus, the space of $\phi_{ji}^{k}$ is a subspace of
$\Hom_{G}(V_{i}\otimes V_{k},V_{j})$ and is therefore finite--dimensional.\\

The relevance of primary fields to the study of fusion lies in the fact that they
may be used to construct intertwiners for the local loop groups and ultimately
explicit elements in \eqref{intertwiner spaces} by smearing them against functions
supported in a disjoint interval. This involves extending them to operator--valued
distributions on $C^{\infty}(S^{1},V_{k})$ and requires a careful study of their
continuity properties which is carried out in chapters \ref{ch:vertex operator} and
\ref{ch:sobolev fields}. This shows that they extend to jointly continuous maps
$\hsmooth_{i}\otimes C^{\infty}(S^{1},V_{k})\rightarrow\hsmooth_{j}$, where
$\hsmooth_{i}\subset\H_{i},\hsmooth_{j}\subset\H_{j}$ are the subspaces of smooth
vectors satisfying as expected
\begin{equation}
\pi_{j}(\gamma)\phi_{ji}^{k}(f)\pi_{i}(\gamma)^{*}=
\phi_{ji}^{k}(\gamma f)
\end{equation}
for any $\gamma\in LG$.
Choosing $f$ supported in $I^{c}$ yields an operator $\hsmooth_{i}\rightarrow\hsmooth_{j}$
commuting with $\Lig$. Although in general $\phi(f)$ is unbounded, a bounded operator
may be obtained by taking the phase of its polar decomposition. This yields
an element of $\Hom_{\Lig}(\H_{i},\H_{j})$.\\

The transport formula \eqref{transport} for these elements is derived from the remarkable
commutation or {\it braiding} properties of primary fields to which we now turn. These
read
\begin{equation}\label{braiding identity}
\phi_{i_{4}j}^{i_{3}}(w)\phi^{i_{2}}_{ji_{1}}(z)=
\sum_{h}\lambda_{h}
\phi_{i_{4}h}^{i_{2}}(z)\phi^{i_{3}}_{hi_{1}}(w)
\end{equation}
and are a generalisation of the canonical commutation or anti--commutation relations
$\phi(z)\phi(w)=\pm\phi(w)\phi(z)$ to a setting where the symmetric group
$\mathfrak S_{2}$ is replaced
by the braid group on two strings, {\it i.e.}~$\IZ$. The sum on the right hand--side spans
all level $\ell$ irreducible positive energy representations $\H_{h}$ and all primary
fields
\begin{xalignat}{3}
\phi_{hi_{1}}^{i_{3}}&:\Hfin_{i_{1}}\otimes V_{i_{3}}[z,z^{-1}]\rightarrow\Hfin_{h}&
&\text{and}&
\phi_{i_{4}h}^{i_{2}}&:\Hfin_{h}\otimes V_{i_{2}}[z,z^{-1}]\rightarrow\Hfin_{i_{4}}
\end{xalignat}
Note that \eqref{braiding identity} is not a formal power series identity.
Rather, it should be interpreted as saying that the matrix coefficients
of the right hand--side define a holomorphic function on $|z|>|w|$. Similarly,
those of the left hand--side define a holomorphic function on $|w|>|z|$ which
may be analytically continued to a (multi--valued) function on $z\neq w$ and
coincides with the right--hand side on $|w|>|z|$.\\

Braiding identities such as \eqref{braiding identity} are established by means
of the {\it four--point function} of the primary fields
\begin{equation}\label{four point}
F=
(\phi_{i_{4}j}^{i_{3}}(v_{3},w)\phi_{ji_{1}}^{i_{2}}(v_{2},z)v_{1},v_{4})
\end{equation}
where $v_{1}\in V_{i_{1}}=\H_{i_{1}}(0)$ and $v_{4}\in V_{i_{4}}=\H_{i_{4}}(0)$.
$F$ is a formal power series in $z,w$ with coefficients in
$V_{i_{4}}\otimes V_{i_{3}}^{*}\otimes V_{i_{2}}^{*}\otimes V_{i_{1}}^{*}$,
and uniquely determines the product
$\phi_{i_{4}j}^{i_{3}}(w)\phi_{ji_{1}}^{i_{2}}(z)$. An important feature of $F$
is that it satisfies the {\it Knizhnik--Zamolodchikov equations} with respect to
the variable $\zeta=zw^{-1}$
\begin{equation}\label{kz equations}
\frac{dF}{d\zeta}=
\frac{1}{\kappa}\Bigl(\frac{\Omega_{12}}{\zeta}+\frac{\Omega_{23}}{\zeta-1}\Bigr)F
\end{equation}
where $\kappa=\ell+2(n-1)$ and the $\Omega_{ij}$ are matrices canonically
associated to the tensor product 
$V_{i_{4}}\otimes V_{i_{3}}^{*}\otimes V_{i_{2}}^{*}\otimes V_{i_{1}}^{*}$.
The above is a fuchsian ODE with regular singular points at $0,1,\infty$ and,
as explained in chapter \ref{ch:algebraic fields}, the four--point functions of
products $\phi_{i_{4}j}^{i_{3}}(w)\phi_{ji_{1}}^{i_{2}}(z)$ with variable $j$ or
products $\phi_{i_{4}h}^{i_{2}}(w)\phi_{hi_{1}}^{i_{3}}(z)$ with variable $h$
form a basis of solutions of these equations
\footnote{In fact, of a canonical subspace of these.} so that
\eqref{braiding identity} holds. In fact, the solutions corresponding to the
products  $\phi_{i_{4}j}^{i_{3}}(w)\phi_{ji_{1}}^{i_{2}}(z)$ diagonalise the
monodromy of \eqref{kz equations} about the singular point $0$ while those
corresponding to $\phi_{i_{4}h}^{i_{2}}(w)\phi_{hi_{1}}^{i_{3}}(z)$ diagonalise
the monodromy about $\infty$. It follows that the braiding constants $\lambda_{h}$
may be computed by analytically continuing the four--point function \eqref{four point}
from $0$ to $\infty$ and re--expressing it as a sum of functions diagonalising
the monodromy at $\infty$.\\

The smeared primary fields obey the same braiding relations as their algebraic
counterparts. Fix for definiteness $I=(0,\pi)$ so that $I^{c}=(\pi,2\pi)$ and
let $f\in C^{\infty}(S^{1},V_{i_{3}})$ and $g\in C^{\infty}(S^{1},V_{i_{2}})$
be supported in $I$ and $I^{c}$ respectively. Then, as shown in chapter
\ref{ch:algebraic fields}, the following holds on $\hsmooth_{i_{1}}$
\begin{equation}
\phi_{i_{4}j}^{i_{3}}(f)\phi^{i_{2}}_{ji_{1}}(g)=
\sum_{h}\beta_{h}
\phi_{i_{4}h}^{i_{2}}(ge_{\alpha_{h}})\phi^{i_{3}}_{hi_{1}}(fe_{-\alpha_{h}})
\end{equation}
where $e_{\mu}(\theta)=e^{i\mu\theta}$ and the $\alpha_{jh}$ are inessential
phase corrections.\\

We now outline the proof of the transport formula \eqref{transport}. We start from
the braiding relation
\begin{equation}\label{eq:first}
\phi_{i0}^{i}(w)\phi_{0j}^{\overline\jmath}(z)=
\sum_{k}\lambda_{k}\phi_{ik}^{\overline\jmath}(z)\phi_{kj}^{i}(w)
\end{equation}
Here, $\phi_{i0}^{i}(w)$ and $\phi_{0j}^{\overline\jmath}(z)$ are the primary fields
whose initial terms are the canonical intertwiners $V_{i}\otimes\IC\rightarrow V_{i}$
and $V_{j}\otimes V_{j}^{*}\rightarrow\IC$. $k$ labels the irreducible summands
of $V_{j}\otimes V_{i}$ and we assume for simplicity that these have multiplicity one
so that for any such $k$ there are at most a primary field of the form $\phi_{kj}^{i}$
and one of the form $\phi_{ik}^{\overline\jmath}$. We shall need another braiding
relation, namely
\begin{equation}\label{second}
\phi_{kj}^{i}(w)\phi_{j0}^{j}(z)=
\epsilon_{k}\phi_{ki}^{j}(z)\phi_{i0}^{i}(w)
\end{equation}
where the sum on the right hand--side of \eqref{second} contains only one term
since there are only one primary field of the form $\phi_{i0}^{i}$ and one of
the form $\phi_{ki}^{j}$. Let now $f\in C^{\infty}(S^{1},V_{i})$ and
$g\in C^{\infty}(S^{1},V_{j})$
so that $\overline{g}\in C^{\infty}(S^{1},V_{j}^{*})$ be supported in $I$
and $I^{c}$ respectively. Then, smearing \eqref{eq:first}--\eqref{second}, we
find
\begin{align}
\phi_{i0}^{i}(f)\phi_{0j}^{\overline\jmath}(\overline{g})
&=
\sum_{k}\lambda_{k}
\phi_{ik}^{\overline\jmath}(\overline{g}e_{\alpha_{k}})
\phi_{kj}^{i}(fe_{-\alpha_{k}})\\
\phi_{kj}^{i}(fe_{-\alpha_{k}})\phi_{j0}^{j}(g)
&=
\epsilon_{k}\phi_{ki}^{j}(ge_{-\alpha_{k}})\phi_{i0}^{i}(f)
\end{align}
It is easy to see that
$\phi_{0j}^{\overline\jmath}(\overline{g})\subseteq\phi_{j0}^{j}(g)^{*}$
and similarly that
$\phi_{ik}^{\overline\jmath}(\overline{g}e_{\alpha_{k}})
 \subseteq\phi_{ki}^{j}(ge_{-\alpha_{k}})^{*}$.
Granted this, we obtain by alleviating notations
\begin{align}
x_{i0}y_{j0}^{*}&=\sum_{k}\lambda_{k}y_{ki}^{*}x_{kj} \label{one}\\
x_{kj}y_{j0}&=\epsilon_{k}y_{ki}x_{i0}		      \label{two}
\end{align}
where the $x_{qp},y_{qp}$ commute with $\Licg$ and $\Lig$ respectively. The above
is an unbounded version of the transport formula \eqref{transport} since
$x_{i0}y_{j0}^{*}y_{j0}=\lambda_{k}\epsilon_{k}y_{ki}^{*}y_{ki}x_{i0}$.\\

A key lemma due to Wassermann \cite{Wa2} asserts that the $x_{qp}$ and $y_{qp}$
above may be replaced by {\it bounded} operators
$x_{qp}\in\Hom_{\Licg}(\H_{p},\H_{q})$, $y_{qp}\in\Hom_{\Lig}(\H_{p},\H_{q})$
without altering the relations \eqref{one}--\eqref{two} in such a way that $x_{i0}$
and $y_{j0}$ are replaced by their phases
\footnote{At present, the lemma only applies when all primary fields with charge $V_{i}$
(or $V_{j}$) define {\it bounded} operator--valued distributions. We shall however only
be concerned with this case in the present dissertation.}.
These do not quite give the transport formula yet but a further alteration of
the operators based on the fact that the $\pi_{p}(L_{I}G)''$ are type III factors
yields unitary $x_{i0}$ and $y_{j0}$. Moreover, each $\lambda_{k}\epsilon_{k}$ 
is necessarily non--negative and vanishes iff $\lambda_{k}$ does since
$\epsilon_{k}$ turns out to be a root of unity. Thus,

\begin{equation}
x_{i0}y_{j0}^{*}y_{j0}=
\sum_{k}\lambda_{k}\epsilon_{k}y_{ki}^{*}y_{ki}x_{i0}
\end{equation}
As pointed out in \S \ref{ss:connes fusion} however,
$x_{i0}y_{j0}^{*}y_{j0}=x_{i0}\pi_{i}(y_{j0}^{*}y_{j0})$
and therefore simplifying by $x_{i0}$ we find \eqref{transport}.
Combining the above with the discussion of \S \ref{ss:connes fusion}, we obtain

\begin{proposition}\label{upper bound}
Let $\H_{i},\H_{j}$ be irreducible positive energy representations at level $\ell$
with lowest energy subspaces $V_{i},V_{j}$. Assume in addition that the irreducible
summands of $V_{i}\otimes V_{j}$ have multiplicity 1 and that all primary fields
with charge $V_{i}$ extend to bounded operator--valued distributions. Then
\begin{equation}
\H_{i}\fuse\H_{j}=\bigoplus_{k}\H_{k}
\end{equation}
where the sum spans the positive energy representations at level $\ell$ whose
lowest energy subspace $V_{k}$ is contained in $V_{i}\otimes V_{j}$ and such
that the corresponding braiding coefficient $\lambda_{k}$ in \eqref{eq:first}
does not vanish.
\end{proposition}

\ssubsection{Statement of main results}\label{ss:main results}

We adopt a more graphical notation and label the simplest representations of
$\SO_{2n}$ by the corresponding Young diagrams. Thus, we denote the defining
representation of $\SO_{2n}$ by $\Vbox$ and its second exterior and symmetric
(traceless) powers by $\Valt$ and $\Vsym$ respectively. Moreover, if $U$ is an
irreducible representation of $G=\Spin_{2n}$ admissible at level $\ell$, we
denote by $\H_{U}$ the corresponding positive energy representation of
$L\Spin_{2n}$ whose lowest energy subspace is $U$.

\begin{theorem}\label{th:main}
Let $U$ be an irreducible representation of $\SO_{2n}$ admissible at
level $\ell$. Then,
\begin{equation}\label{fuse with box}
\Hbox\fuse\H_{U}=
\bigoplus_{W\subset \Vbox\otimes U}N_{\Box U}^{W}\H_{W}
\end{equation}
where $N_{\Box U}^{W}=1$ is $W$ is admissible at level $\ell$ and $0$
otherwise.
\end{theorem}

The proof of theorem \ref{th:main} rests on the following important
special case

\begin{proposition}\label{pr:rule}
The following fusion rules hold at level $\ell$
\begin{align}
\Hbox\fuse\Hbox&=\Hsym\oplus\Halt\oplus\H_{0}\\
\intertext{if $\ell\geq 2$ and}
\Hbox\fuse\Hbox&=\H_{0}
\end{align}
if $\ell=1$.
\end{proposition}

{\sc Proof of proposition \ref{pr:rule}.}
We use proposition \ref{upper bound} and the tensor product rule
\begin{equation}
\Vbox\otimes\Vbox=\Vsym\oplus\Valt\oplus V_{0}
\end{equation}
The primary fields with charge $\Vbox$ are shown to extend to bounded
operator--valued distributions in chapter
\ref{ch:sobolev fields}. The proof relies on the simple observation that they
are essentially Fermi fields. The relevant braiding coefficients are computed
in chapter \ref{ch:box/kbox braiding} by explicitly solving the corresponding
Knizhnik--Zamolodchikov equations and shown to be non--zero. The $\ell=1$ fusion
rule for $\Hbox\fuse\Hbox$ differs from those for $\ell\geq 2$ because $\Vsym$
and $\Valt$ are not admissible at level 1 \halmos\\

{\sc Proof of theorem \ref{th:main}.}
We begin by using proposition \ref{upper bound} to obtain an upper bound on
$\Hbox\fuse\H_{U}$. Since $\Vbox$ is a minimal representation of $\Spin_{2n}$, all
irreducible summands
of $\Vbox\otimes U$ have multiplicity one. Moreover, the primary fields with charge
$\Vbox$ extend to bounded operator--valued distributions and therefore
\eqref{fuse with box} holds with the $N_{\Box\medspace U}^{W}$ replaced by some
$0\leq\wt N_{\Box\medspace U}^{W}\leq N_{\Box\medspace U}^{W}$ since part of the
braiding coefficients involved in the computation of fusion might vanish.
The matrix $N_{\Box}$ whose entries are the $N_{\Box\medspace U}^{W}$, where $U$ and $W$
are single--valued representations of $\SO_{2n}$ is non--negative and irreducible
\footnote{This is the reason for demanding that $U$ be a single--valued representation of
$\SO_{2n}$ in theorem \ref{th:main}. Without this restriction, $N_{\Box}$ would have two
irreducible diagonal blocks, corresponding respectively to the single and two--valued
representations of $\SO_{2n}$, and zero off--diagonal blocks and would therefore fail
to be irreducible.}.
By the Perron--Frobenius theorem \cite{GM}, $N_{\Box}$ has, up to a multiplicative
constant, a unique eigenvector with non--negative entries. On the other hand, if
$2\rho$ is the sum of the positive roots of $\Spin_{2n}$ and $\cchi_{U}$ the
character of the irreducible $\Spin_{2n}$--module $U$, the
'dimensions' $\delta_{U}=\cchi_{U}(\exp(\frac{2\pi i\rho}{\ell+2(n-1)}))\in\IR_{+}$
given in \cite[\S 13.8]{Ka1} obey, by a simple computation
\begin{equation}
\sum_{W\subset\Vbox\otimes U}N_{\Box\medspace U}^{W}\delta_{W}=\delta_{\Box}\delta_{U}
\end{equation}
and it follows that the Perron--Frobenius eigenvalue of $N_{\Box}$ is equal to
$\delta_{\Box}$.\\

A similar statement about the matrix $\wt N_{\Box}$ may be obtained by using some
results of Doplicher--Haag--Roberts on the algebraic theory of superselection
sectors \cite{DHR1,DHR2}.
In essence, the fact that $\Hbox\boxtimes\Hbox$ contains the vacuum representation
$\H_{0}$,
by virtue of proposition \ref{pr:rule} implies that, on the ring $\RR$ generated
by the irreducible summands of the iterated powers $\Hbox^{\boxtimes m}$, $m\in\IN$
there exists a unique positive character or {\it quantum dimension} function $d$
which is additive under direct
sums, multiplicative under fusion and normalised by $d(\H_{0})=1$.
Applying $d$ to \eqref{fuse with box} with $N_{\Box\medspace U}^{W}$ replaced by
$\wt N_{\Box\medspace U}^{W}$ yields
\begin{equation} \label{almost pf}
\sum_{W\subset\Vbox\otimes U}
\wt N_{\Box\medspace U}^{W}d(\H_{W})=
d(\Hbox)d(\H_{U})
\end{equation}
whenever $U$ is a single--valued representation of $\SO_{2n}$ which appears in
$\RR$ \footnote{In principle these are all the single--valued representations
but that is part of what we are trying to establish.}. We do not know $\wt N_{\Box}$
to be irreducible, nor is \eqref{almost pf} a statement about $\wt N_{\Box}$ since
some $\H_{U}$ may not lie in $\RR$ but if $M_{\Box}$ is the matrix whose entries
are $\wt N_{\Box\medspace U}^{W}$ if $\H_{U},\H_{W}\in\RR$ and zero otherwise,
we clearly have $M_{\Box}\leq N_{\Box}$ entry--wise and therefore, by Perron--Frobenius
theory, $d(\Hbox)\leq \delta_{\Box}$ with equality only if $M_{\Box}=N_{\Box}$,
{\it i.e.}~only if \eqref{fuse with box} holds. Finally, a computation based on
\eqref{pr:rule} and paralleling one of Wenzl \cite{We3} shows that
$d(\Hbox)=\delta_{\Box}$
and we have therefore established \eqref{fuse with box} \halmos

\begin{theorem}\label{th:closed}
The positive energy representations of $L\Spin_{2n}$ whose lowest energy subspace is a
single--valued representation of $\SO_{2n}$ are closed under fusion.
\end{theorem}
\proof This follows at once because they are, by theorem \ref{th:main} exactly the
representations that appear in the iterated fusion powers $\Hbox^{\fuse m}$ \halmos\\

The restriction to single--valued representations of $\SO_{2n}$ in theorems
\ref{th:main} and \ref{th:closed} is technical rather than conceptual and
these results conjecturally hold for all positive energy representations of
$L\Spin_{2n}$. At present, we can show this at {\it odd} level by resorting
to the action of the centre $Z(\Spin_{2n})$ via conjugation by discontinuous
loops. Denote by $z\H$ the representation obtained by conjugating $\H$ by
$z\in Z(\Spin_{2n})$ as in \S \ref{ss:representations}. We need the following

\begin{lemma}\label{le:trick}
Let $\H$ be a positive energy representation of $L\Spin_{2n}$ and $\H_{0}$ the
vacuum representation at the same level. Then, for any $z\in Z(\Spin_{2n})$,
\begin{equation}
z\H_{0}\fuse\H\cong z\H\cong \H\fuse z\H_{0}
\end{equation}
\end{lemma}
\proof
We prove the first identity only, the second follows in a similar way.
Let $\zeta\in C^{\infty}([0,2\pi],\Spin_{2n})$ be the lift of a loop in
$\Spin_{2n}/Z(\Spin_{2n})$ with $\zeta(2\pi)\zeta(0)^{-1}=z$. By definition,
the conjugate $z\K$ of any positive energy representation $(\pi,\K)$ is
given by the isomorphism class of the representation
$\gamma\rightarrow\pi(\zeta^{-1}\gamma\zeta)$ of $L\Spin_{2n}$ on $\K$.
Form now the intertwiner spaces
\begin{xalignat}{2}
\XX&=\Hom_{\Lics}(\H_{0},z\H_{0})&
\YY&=\Hom_{\Lis} (\H_{0}, \H)
\end{xalignat}
Choosing a $\zeta$ equal to one on $I$ so that
$\pi_{0}(\zeta^{-1}\gamma\zeta)=\pi_{0}(\gamma)$ for any $\gamma\in\Lics$,
we have
\begin{equation}
\XX=
\Hom_{\Lics}(\H_{0},\H_{0})=
\pi_{0}(\Lis)''
\end{equation}
where the last identity follows by Haag duality. We claim now that the
map $U:\XX\otimes\YY\rightarrow\H$, $x\otimes y\rightarrow yx\Omega$ is 
norm--preserving.
Indeed,
\begin{equation}
\|x\otimes y\|^{2}
=(x^{*}xy^{*}y\Omega,\Omega)
=(x^{*}y^{*}yx\Omega,\Omega)
=\|yx\Omega\|^{2}
\end{equation}
where we have used locality so that $y^{*}y\in\Hom_{\Lis}(\H_{0},\H_{0})$
commutes with $x$. $U$ is equivariant for the actions of
$\Lis\times\Lics$ on $\XX\otimes\YY$ and $z\H$ and extends to an isometry
$z\H_{0}\fuse\H\rightarrow z\H$ whose range is dense by the Reeh--Schlieder
theorem. Thus, $U$ is the required unitary equivalence \halmos

\begin{theorem}\hfill
\begin{enumerate}
\item Let $\H_{U}$ be the irreducible positive energy representations of
$L\Spin_{2n}$ at odd level $\ell$ whose lowest energy subspace is $U$.
Then
\begin{equation}\label{eq:again}
\Hbox\fuse\H_{U}=\bigoplus_{W\subset \Vbox\otimes U}N_{\Box U}^{W}\H_{W}
\end{equation}
where $N_{\Box U}^{W}$ is 1 if $W$ is admissible at level $\ell$ and zero
otherwise.
\item The positive energy representations of $L\Spin_{2n}$ at odd level
$\ell$ are closed under fusion.
\end{enumerate}
\end{theorem}
\proof
(i)
By theorem \ref{th:main}, \eqref{eq:again} holds if $U$ is a single--valued
representation of $\SO_{2n}$. If that is the case, we find by conjugating both
sides by an element $z\in Z(\Spin_{2n})$ and using lemma \ref{le:trick} that
\begin{equation}
\Hbox\fuse z\H_{U}=
\bigoplus_{W\subset\Vbox\otimes U}N_{\Box U}^{W}\thickspace z\H_{W}
\end{equation}
By definition, $z\H_{U}=\H_{zU}$ and $z\H_{W}=\H_{zW}$ where the notation
refers to the induced action of $Z(\Spin_{2n})$ on the irreducible
$\Spin_{2n}$--modules admissible at level $\ell$. The explicit form of this
action is computed in chapter \ref{ch:classification} and shows in particular
that $W\subset U\otimes V_{\Box}$ if, and only if $zW\subset\Vbox\otimes zU$.
Thus,
\begin{equation}
\Hbox\fuse\H_{zU}=
\bigoplus_{W'\subset\Vbox\otimes zU}N_{\Box zU}^{W'}\thickspace\H_{W'}
\end{equation}
Moreover, when $\ell$ is odd, any two--valued representation $U'$ of $\SO_{2n}$
may be written as $z U$ where $z\in Z(\Spin_{2n})$ and $U$ is a single--valued
representation of $\SO_{2n}$. This yields (i).\\

(ii)
By theorem \ref{th:closed}, $\H_{U_{1}}\fuse\H_{U_{2}}$ is of positive energy
if $U_{1}$ and $U_{2}$ are single--valued $SO_{2n}$--modules. Conjugating
successively by $z_{2}$ and $z_{1}\in Z(\Spin_{2n})$ and using lemma \ref{le:trick},
we find that
\begin{equation}
\H_{z_{1}U_{1}}\fuse\H_{z_{2}U_{2}}\cong
z_{1}\H_{0}\fuse(\H_{U_{1}}\fuse\H_{U_{2}})\fuse z_{2}\H_{0}\cong
z_{1}z_{2}(\H_{U_{1}}\fuse\H_{U_{2}})
\end{equation}
is of positive energy.
Since any irreducible $\Spin_{2n}$--module admissible at odd level
$\ell$ may be written as $z U$ for some $z\in Z(\Spin_{2n})$ and
$U$ a single--valued $\SO_{2n}$--module, the result follows
\halmos

\begin{theorem}
The level 1 representations of $L\Spin_{2n}$ are closed under fusion.
Moreover, if $\H_{0}$ is the level 1 vacuum representation, then
\begin{equation}
z\longrightarrow z\H_{0}
\end{equation}
yields an isomorphism of the group algebra $\IC[Z(\Spin_{2n})]$ and
the level 1 fusion ring of $L\Spin_{2n}$.
\end{theorem}
\proof
The above map is bijective since, as noted in \S \ref{ss:representations},
the action of $Z(\Spin_{2n})$ on the level 1 representations is transitive
and free. Moreover, by lemma \ref{le:trick},
\begin{equation}
z_{1}\H_{0}\fuse z_{2}\H_{0}\cong z_{1}(z_{2}\H_{0})=z_{1}z_{2}\H_{0}
\end{equation}
\halmos


\part{Positive energy representations and primary fields}




\chapter{Positive energy representations and primary fields}
\label{ch:classification}

This chapter introduces the basic objects of study of this dissertation, the
{\it positive energy representations} of the loop group $LG=C^{\infty}(S^{1},G)$
of the Lie group $G=\Spin_{2n}$, $n\geq 3$. These are projective unitary
representations supporting an intertwining action of $\rot$ with
finite--dimensional eigenspaces and positive spectrum. As explained
in section \ref{se:classification of per}, they are completely reducible
and the irreducible ones are classified by their {\it level} $\ell\in\IN$,
which determines the underlying projective cocycle, and their {\it lowest
energy subspace}, an irreducible $G$--module. Moreover, only finitely many
irreducibles exist at any given level $\ell$.
In section \ref{se:level 1 per}, we restrict our attention to the level 1
representations of $LG$ and their lowest energy subspaces, the {\it minimal}
$G$--modules. These correspond in a natural way to elements in the dual
of the centre of $G$ and play an important r\^ole in view of the fact that
any level $\ell$ representation is contained in the $\ell$--fold tensor
product $\H_{i_{1}}\otimes\cdots\otimes\H_{i_{\ell}}$ of level 1 representations.\\

In section \ref{se:disc loops}, we show that the outer action of the centre
$Z(G)$ on $LG$ via conjugation by discontinuous loops, {\it i.e.}~lifts to $G$ of
loops in $G/Z(G)$, induces a level preserving one on the positive energy
representations of $LG$. At level 1, the action is transitive and free, a fact
which will be used repeatedly, most notably in chapter \ref{ch:connes fusion}
to prove the isomorphism of the level 1 fusion ring of $LG$ and the group
algebra of $Z(G)$.
Finally, in section \ref{se:classification of pf}, we outline the classification
of the {\it primary fields} of $LG$. These are algebraic, operator--valued
distributions mapping between positive energy representations of equal level.
The study of their analytic and algebraic properties, both of which are 
essential ingredients in the computation of fusion, will be carried out
in chapters \ref{ch:vertex operator}--\ref{ch:sobolev fields} and
\ref{ch:algebraic fields}--\ref{ch:box/kbox braiding} respectively.

\ssection{Definition and classification of positive energy representations}
\label{se:classification of per}

In this section and the rest of this chapter, $G$ denotes a compact, connected
and simply--connected simple Lie group with Lie algebra $\g$.

\ssubsection{The loop group $LG$}\label{ss:LG}

The {\it loop group} of $G$ is, by definition $LG=C^{\infty}(S^{1},G)$.
When endowed with the $C^{\infty}$--topology and pointwise multiplication, $LG$ is
a real analytic Fr\'echet Lie group. Its Lie algebra is $L\g=C^{\infty}(S^{1},\g)$
with pointwise bracket. It has complexification $L\gc=C^{\infty}(S^{1},\gc)$ and
is often more conveniently regarded as the $-1$ eigenspace of the anti--linear
anti--involution $X\rightarrow -\overline{X}$ determined by the canonical
conjugation on $\gc$. As Lie groups, $LG\cong\Omega G\rtimes G$ where $\Omega G$
is the space of based loops on which $G$ acts by conjugation and the map is given
by $\gamma\rightarrow(\gamma\gamma(0)^{-1},\gamma(0))$.
Since $\pi_{0}(G)=\pi_{1}(G)=\pi_{2}(G)=0$, $LG$ is connected and
simply--connected. In particular, it is generated by the image of the exponential
map since the latter is locally one--to--one.\\

$LG$ admits a smooth automorphic action of $\rot$ given by $R_{\theta}\gamma=
\gamma_{\theta}$ where $\gamma_{\theta}(\phi)=\gamma(\phi-\theta)$.
The corresponding semi-direct product $LG\rtimes\rot$ is therefore a Fr\'echet
Lie group, which however fails to be analytic, since for a fixed $\gamma\in LG$,
$\theta\rightarrow\gamma_{\theta}$ is analytic iff $\gamma$ is an analytic loop
in $G$. Identifying the Lie algebra of $\rot$ and $i\IR$ with generator $id$ so
that $R_{\theta}=\exp(i\theta d)$, the Lie algebra of $LG\rtimes\rot$ is
$L\g\rtimes i\IR$ and the following relations hold for $\gamma\in LG$ and
$X\in L\g$
\begin{align}
\Ad(R_{\theta})X&=X_{\theta}\\
[id,X]&=
\left.\frac{d}{d\theta}\right|_{\theta=0}\Ad(R_{\theta})X=
-\dot X\\
\intertext{and}
\Ad(\gamma) id&=
\left.\frac{d}{d\theta}\right|_{\theta=0}\gamma R_{\theta}\gamma^{-1}=
id-\gamma\dot{(\gamma^{-1})}=
id+\dot\gamma\gamma^{-1} \label{Ad rot}
\end{align}
where we have used
$0=\dot{(\gamma\gamma^{-1})}=\dot\gamma\gamma^{-1}+\gamma\dot{(\gamma^{-1})}$.\\

We shall be concerned with projective representations of $LG\rtimes\rot$ and will
therefore need to consider continuous Lie algebra cocycles on $L\g\rtimes i\IR$,
{\it i.e.}~skew--symmetric, $i\IR$--valued maps $\beta$ satisfying the Jacobi
identity
\begin{equation}
\beta([X,Y],Z)+\beta([Z,X],Y)+\beta([Y,Z],X)=0
\end{equation}
Up to the addition of coboundaries $\beta(X,Y)=d\alpha(X,Y)=\alpha([X,Y])$, these
vanish on the Lie algebra of $\rot$ and their restriction to $L\g$ is a multiple
of the fundamental cocycle \cite[4.2.4]{PS}
\begin{equation}\label{fundamental cocycle}
iB(X,Y)=i\int_{0}^{2\pi}\<X,\dot Y\>\frac{d\theta}{2\pi}
\end{equation}
where $\<\cdot,\cdot\>$ is the {\it basic inner product} on $\gc$, {\it i.e.}~the
unique multiple of the Killing form for which the highest root $\theta$ has squared
length 2.\\

$L\g$ contains a distinguished dense sub--algebra $\lpol\g$ consisting of all
$\g$--valued trigonometric polynomials. Its complexification $\lpol\gc$ is spanned
by the elements $x(n)=x\otimes e^{in\theta}$, $x\in\gc$, $n\in\IZ$ with bracket
$[x(n),y(m)]=[x,y](n+m)$. The restriction of the cocycle \eqref{fundamental cocycle}
to $\lpol\gc$ reads
\begin{equation}\label{combinatorial cocycle}
iB(x(n),y(m))=n\delta_{n+m,0}\<x,y\>
\end{equation}
The action of $\rot$ leaves $\lpol\gc$ invariant and satisfies $[d,x(n)]=-nx(n)$.
The central extension of $\lpol\gc$ with central term $c$ and cocycle
given by \eqref{combinatorial cocycle} is usually denoted by $\wt\gc$.
The semi--direct product $\wt\gc\rtimes\IC d$ is the affine, untwisted
Kac--Moody algebra $\wh\gc$ corresponding to $\gc$ \cite{Ka1}\footnote{The
$d$ above is the opposite of that used in the theory of affine Kac--Moody
algebras.}.
$\lpol\gc$ possesses two decompositions we shall make use of. The first is
the triangular decomposition determined by a maximal torus $T\subset G$
\begin{equation}\label{eq:triangle 1}
\lpol\gc=\g_{\leq}\oplus\tc\oplus\g_{\geq}
\end{equation}
where $\tc$ is the complexified Lie algebra of $T$ and $\g_{\leq}$ (resp. $\g_{\geq}$)
is the nilpotent Lie algebra spanned by the $x(n)$ with $n<0$ (resp. $n>0$) and
$x\in\gc$ or $n=0$ and $x$ lying in a negative (resp. positive) root space of $\gc$.
The second decomposition is given by
\begin{equation}\label{eq:triangle 2}
\lpol\gc=\g_{-}\oplus\gc\oplus\g_{+}
\end{equation}
where $\g_{\pm}$ are spanned by the $x(n)$, $x\in\gc$, $n\gtrless 0$.

\ssubsection{Positive energy representations of $LG$}\label{ss:per}

We outline below the classification of positive energy representations of $LG$
following \cite{Wa2} to which we refer for further details. The basic properties
of projective representations and central extensions relevant to the present
discussion may be found in \S \ref{ss:projective}.
Let $\pi$ be a projective unitary representation of $LG\rtimes\rot$ on
a complex Hilbert space $\H$, {\it i.e.}~a strongly continuous homomorphism
\begin{equation}
\pi:LG\rtimes\rot\longrightarrow PU(\H)=U(\H)/\T
\end{equation}
Over $\rot$, $\pi$ possesses a continuous lift to a unitary representation
which we denote by the same symbol. By definition,
$\H$ is of {\it positive energy} if $\H=\bigoplus_{n\geq n_{0}}\H(n)$ with
$\dim\H(n)<\infty$ and  $\left.\pi(R_{\theta})\right|_{\H(n)}=e^{in\theta}$.
The lift is clearly only determined up to multiplication by a character of
$\rot$ and we normalise it by choosing $n_{0}=0$ and $\H(0)\neq 0$.
Since $G$ is simple, connected and simply--connected, the restriction of $\pi$
to $G$ lifts uniquely to a strongly continuous unitary representation which we
also denote by $\pi$. It commutes with the unitary action of $\rot$
since, projectively $\pi(g)\pi(R_{\theta})\pi(g)^{*}\pi(R_{\theta})^{*}=1$. Thus,
lifting each representation, the following holds in $U(\H)$
\begin{equation}
\pi(g)\pi(R_{\theta})\pi(g)^{*}\pi(R_{\theta})^{*}=\cchi(g,\theta)
\end{equation}
where $\cchi(g,\theta)\in\T$ is continuous and multiplicative in either variable
and therefore defines a continuous homomorphism $G\rightarrow\Hom(\rot,\T)\cong\IZ$.
Since $G$ is connected, $\cchi\equiv 1$ as claimed and it follows in particular
that each $\H(n)$ is a finite--dimensional $G$--module.\\

The spectral assumption on the action of $\rot$ implies that a positive energy
representation is completely reducible. Since it is projective however, some
care is needed in defining the direct sum of its irreducible components. Indeed,
as explained in \S \ref{ss:projective}, the direct sum of two projective
representations $(\pi_{i},\H_{i})$, $i=1,2$ of a topological group $\Gamma$ may
be defined only if the corresponding central extensions
\begin{equation}
\pi_{i}^{*}U(\H_{i})=\{(g,u)\in\Gamma\times U(\H_{i})|\pi_{i}(g)=p_{i}(u)\}
\end{equation}
obtained by pulling back to $\Gamma$ the canonical central extensions
\begin{equation}
1\rightarrow\T\rightarrow U(\H_{i})\xrightarrow{p_{i}}PU(\H_{i})
\end{equation}
are isomorphic. When that is the case, the definition of $\pi_{1}\oplus\pi_{2}$
depends upon the choice of an identification
$\pi_{1}^{*}U(\H_{1})\cong\pi_{2}^{*}U(\H_{2})$ which is unique only
up to multiplication by an element of $\Hom(\Gamma,\T)$.\\

Since $\Hom(LG,\T)=\{1\}$ \cite[prop. 3.4.3]{PS} and $\Hom(LG\rtimes\rot,\T)=\IZ$,
the direct sum of projective representations of $LG\rtimes\rot$ is ill--defined
and we shall therefore forgetfully regard positive energy representations as
$LG$--modules. This does not affect the irreducible ones since a positive energy
representation is irreducible under $LG\rtimes\rot$ if, and only if it remains so
when restricted to $LG$ \cite[prop. 9.2.3]{PS} and leads to a canonically, if
partially defined notion of direct sum. Thus, we regard two positive energy
representations $(\pi_{i},\H_{i})$ as unitarily equivalent if they are so as
projective $LG$--modules, {\it i.e.}~if there exists a unitary
$U:\H_{1}\rightarrow\H_{2}$ such that
$U\pi_{1}(\gamma)U^{*}=\pi_{2}(\gamma)$ in $PU(\H_{2})$ for all $\gamma\in LG$.
If both are irreducible, this is easily seen to imply their unitary equivalence
as projective $LG\rtimes\rot$--modules.\\

The classification of a positive energy representation $(\pi,\H)$ is obtained via
the associated infinitesimal action of $\lpol\g$ in the following way. Consider
the subspace $\hfin\subset\H$ of {\it finite energy vectors} for $\rot$, that
is the algebraic direct sum $\sum_{n\geq 0}\H(n)$. The latter is a core for the
normalised self--adjoint generator of rotations which we abusively denote by $d$.
Thus
\begin{xalignat}{3}
\left.d\right|_{\H(n)}&=n&&\text{and}&\pi(\exp_{\rot}(i\theta d))&=e^{i\theta d}
\end{xalignat}
For any $X\in\lpol\g$, the one--parameter projective group $\pi(\exp_{LG}(tX))$
possesses a continuous lift to $U(\H)$, unique up to multiplication by a character
of $\IR$. It is therefore given, via Stone's theorem by $e^{t\pi(X)}$ where
$\pi(X)$ is a skew--adjoint operator determined up to an additive constant.

\begin{theorem}\label{core}
The subspace $\Hfin$ of finite energy vectors is an invariant core for the
operators $\pi(X)$, $X\in\lpol\g$. The operators $\pi(X)$ may be chosen
uniquely so as to satisfy $[d,\pi(X)]=i\pi(\dot X)$ on $\Hfin$ and then
$X\rightarrow\pi(X)$ gives a projective representation of $\lpol\g$ on
$\Hfin$ such that
\begin{equation}\label{eq:KM reln}
[\pi(X),\pi(Y)]=\pi([X,Y])+i\ell B(X,Y)
\end{equation}
where $iB(X,Y)$ is given by \eqref{fundamental cocycle} and $\ell$ is
a non--negative integer called the level of the representation.
\end{theorem}

When no confusion arises, we denote the restriction of the operators $\pi(X)$,
$X\in\lpol\g$ to $\Hfin$ by the same symbol and extend the resulting projective
representation $\pi:\lpol\g\rightarrow\End(\Hfin)$ to one of $\lpol\gc$ satisfying
\eqref{eq:KM reln} as well as the formal adjunction property
$\pi(X)^{*}=-\pi(\overline{X})$.\\

By theorem \ref{core}, the operators $\pi(X)$ and $d$ give rise to a unitarisable
representation of the Kac--Moody algebra $\wh\gc$ at level $\ell$ such that
$d$ is diagonal with finite--dimensional eigenspaces and spectrum in $\IN$.
Such representations split into a direct sum of irreducibles, each of which
is necessarily an integrable highest weight representation, that is a module
generated over the enveloping algebra $\U(\g_{\leq})$ by a vector $v$ annihilated
by $\g_{\geq}$ and diagonalising the action of $T\rtimes\rot$. Thus, for any
$h\in\tc$
\begin{align}
d v&=n v\\
\pi(h) v&=\lambda(h) v
\end{align}
for some $n\in\IN$ and dominant integral weight $\lambda$ of $G$ satisfying
$\<\lambda,\theta\>\leq\ell$ where $\theta$ is the highest root. The
pair $(\ell,\lambda)$ classifies the representation as an $\lpol\g$--module
uniquely \cite{Ka1}.
The highest weight representation is also generated over the enveloping algebra
$\U(\g_{-})$ by its {\it lowest energy subspace}, {\it i.e.}~the $d$--eigenspace with
lowest eigenvalue which in fact coincides with the irreducible $G$--module with
highest weight $\lambda$ generated by $v$. We shall interchangeably adopt either
point of view. When the
level is understood, we refer to $\lambda$ as the highest weight of the
representation. The collection of dominant integral weights of $G$ satisfying
$\<\lambda,\theta\>\leq\ell$ is a finite set called the level $\ell$
{\it alcove} and is denoted by $\al$.\\

Propositions \ref{simple} and \ref{equivalent} below imply that the classification
of positive energy representations of $LG$ is equivalent to that of their finite
energy subspaces as $\lpol\g$--modules.

\begin{proposition}\label{simple}
Let $(\pi,\H)$ be a positive energy representation of $LG$ with finite energy subspace
$\Hfin$. If $\K\subset\H\fin$ is invariant under $\lpol\g\rtimes\IC d$, then $\overline{\K}$
is invariant under $LG\rtimes\rot$. In particular, $\H$ is topologically irreducible under
$LG\rtimes\rot$ if, and only if $\H\fin$ is algebraically irreducible for 
$\lpol\g\rtimes\IC d$.
\end{proposition}
\proof
By continuity of $\pi$, it is sufficient to prove the invariance of $\overline{\K}$ under
the dense subgroup generated by $\exp_{LG}(\lpol\g)$. Let $X\in\lpol\g$ and $\pi(X)$ the
corresponding skew--adjoint operator on $\H$ with domain $\D(\pi(X))$ and invariant core
$\Hfin\subset\D(\pi(X))$. It is sufficient to prove that $P\pi(X)\subset\pi(X)P$ where $P$
is the orthogonal projection on $\overline{\K}$ for then, by the spectral theorem, $P$
commutes with bounded functions of $\pi(X)$ and in particular with $e^{t\pi(X)}=
\pi(\exp_{LG}(tX))$.
Notice that $\K$ is invariant under $d$ and is therefore a graded subspace of $\hfin$
so that it coincides with the finite energy subspace of $\overline{\K}$. Since $P$
commutes with $\rot$, it follows that its restriction to $\Hfin$ is the orthogonal
projection onto $\K$ and therefore $\pi(X)P\eta=P\pi(X)\eta$ for any $\eta\in\Hfin$.
Let now $\xi\in\D(\pi(X))$ and $\xi_{n}\in\Hfin$ a sequence such that
$\xi_{n}\rightarrow\xi$ and $\pi(X)\xi_{n}\rightarrow\pi(X)\xi$. Then
$\eta_{n}=P\xi_{n}\in\Hfin\subset\D(\pi(X))$ converges to $P\xi$ and
$\pi(X)\eta_{n}=P\pi(X)\xi_{n}\rightarrow P\pi(X)\xi$ whence $P\pi(X)\subset\pi(X)P$ as
claimed. To conclude, notice that $\overline{\K}=\H$ iff $\K\fin=\H\fin$ and therefore
topological irreducibility of $\H$ implies algebraic irreducibility of $\H\fin$.
The converse holds by complete reducibility of positive energy representations
and functoriality \halmos

\begin{proposition}\label{equivalent}
Two positive energy representations $\H_{1}$ and $\H_{2}$ are unitarily equivalent
as $LG$--modules if, and only if their finite energy subspaces are isomorphic as
$\lpol\g$--modules.
\end{proposition}
Evidently, the unitarily equivalence of $\H_{1},\H_{2}$ as $LG$--modules implies that
of $\Hfin_{1},\Hfin_{2}$ as $\lpol\g$--modules. By complete reducibility, we need only
prove the converse when the $\H_{i}$ and therefore the $\Hfin_{i}$ are irreducible.
Let $U:\Hfin_{1}\rightarrow\Hfin_{2}$ be an isomorphism of $\lpol\g$--modules. Up to
multiplication by a scalar, $U$ is an isometry since highest weight representations
admit a unique $\lpol\g$--invariant inner product.
Thus U extends to a unitary $\H_{1}\rightarrow\H_{2}$ such that
$U\pi_{1}(X)U^{*}=\pi_{2}(X)$ on $\H_{2}\fin$ for any $X\in\lpol\g$. Since $\Hfin_{i}$
is a core for $\pi_{i}(X)$, it follows that $U\pi_{1}(X)U^{*}=\pi_{2}(X)$ holds as an
operator identity and therefore that $U$ intertwines with the one parameter groups
$\exp_{LG}(tX)$, $X\in\lpol\g$. Since these generate a dense subgroup in $LG$, $U$
is an $LG$ intertwiner \halmos\\

\begin{corollary}\label{co:data}
An irreducible, positive energy representation of $LG$ is uniquely determined
by its level $\ell\in\IN$ and its lowest energy subspace, an irreducible
$G$--module whose highest weight $\lambda$ satisfies $\<\lambda,\theta\>\leq\ell$.
\end{corollary}

\definition An irreducible $G$--module $V$ whose highest weight $\lambda$ satisfies
$\<\lambda,\theta\>\leq\ell$ will be called {\it admissible at level $\ell$}.
The corresponding positive energy representation will be denoted by $\H_{V}$
or $\H_{\lambda}$.\\

\remark We will show in chapter \ref{ch:analytic} (proposition \ref{classify extension})
that the central extensions of $LG$ corresponding to positive energy representations
of levels $\ell_{1},\ell_{2}$ are isomorphic if, and only if $\ell_{1}=\ell_{2}$. In
particular, the direct sum of two positive energy representations of equal level is
unambiguously defined.\\

Corollary \ref{co:data} settles the uniqueness of irreducible positive energy
representations of $LG$ corresponding to a given level $\ell$ and irreducible
$G$--module $V$. The corresponding existence problem may be solved by one of
the following methods
\begin{enumerate}
\item
It is easy to construct the corresponding highest weight representation $\Nu$
of $\wh\gc$ using Verma modules and to prove that $\Nu$ is unitarisable \cite{Ka1}.
Some simple estimates of Goodman and Wallach \cite{GoWa} then show that the
action of $\wt\gc$ extends to one of the Lie algebra $L\gc$ on a suitable
Fr\'echet completion of $\Nu$. This action may then be exponentiated to
a projective unitary representation of $LG$ on the Hilbert space completion
of $\Nu$. This is carried out for the dense subgroup of analytic loops
$S^{1}\rightarrow G$ in \cite{GoWa} and, by different methods for the
full loop group $LG$ in \cite{TL1}.
\item
One may alternatively consider the space of holomorphic section of a suitable
vector bundle with fibre $V$ over the flag manifold $LG/G$. This space possesses
an $LG$--invariant inner product and its Hilbert space completion is $\H_{V}$
\cite[chap. 11]{PS}.
\end{enumerate}
The positive energy representations of the loop groups relevant to this thesis,
namely those corresponding to $G=\Spin_{2n}$, $n\geq 3$ possess an alternative
Fermionic construction which will be given in chapter \ref{ch:fermionic}.

\ssubsection{Appendix : projective representations and central extensions}
\label{ss:projective}

We give the definitions and elementary properties of central extensions and
projective representations and discuss in particular the problems associated
with defining the direct sum of the latter. All groups below are assumed to
be topological and the associated homomorphisms continuous.\\

\definition
Let $\Gamma,A$ be groups with $A$ abelian. A {\it central extension} of $\Gamma$
by $A$ is a short exact sequence
\begin{equation}\label{central}
1\rightarrow A\xrightarrow{i}\wt\Gamma\xrightarrow{p}\Gamma\rightarrow 1
\end{equation}
with $i(A)<Z(\wt\Gamma)$. A homomorphism of central extensions is a
map $\phi:\wt\Gamma_{1}\rightarrow\wt\Gamma_{2}$ making the following
a commutative diagram
\vskip -1.5em
\begin{equation}
\begin{diagram}[height=1.2em,width=1.7em]
 &    & &     &\wt\Gamma_{1}	&     &		&    & \\
 &    & &\ruTo&         	&\rdTo&		&    & \\
1&\rTo&A&     &\dTo>\phi	&     &\Gamma	&\rTo&1\\
 &    & &\rdTo&         	&\ruTo&		&    & \\
 &    & &     &\wt\Gamma_{2}	&     &		&    & \\
\end{diagram}
\end{equation}

Given two central extensions $\wt\Gamma_{1}$, $\wt\Gamma_{2}$ we may form their
{\it product} $\wt\Gamma=\wt\Gamma_{1}\star\wt\Gamma_{2}$ as the quotient of
\begin{equation}
\{(\wt g_{1},\wt g_{2})\in\wt\Gamma_{1}\times\wt\Gamma_{2}|
    \thinspace p_{1}(\wt g_{1})=p_{2}(\wt g_{2})\}
\end{equation}
by the image of the diagonal embedding $A\rightarrow\wt\Gamma_{1}\times\wt\Gamma_{2}$,
$a\rightarrow (i_{1}(a),i_{2}(a^{-1}))$. The associated maps are given by
$p(\wt g_{1},\wt g_{2})=p_{1}(\wt g_{1})=p_{2}(\wt g_{2})$ and
$i(a)=(i_{1}(a),1)=(1,i_{2}(a))$. The set of isomorphism classes of central
extensions of $\Gamma$ by $A$ endowed with this product operation is easily seen
to be an abelian group, usually denoted by $H^{2}(\Gamma,A)$ with
the trivial central extension $\Gamma\times A$ as identity element. The inverse
$\overline{\wt\Gamma}$ of the central extension \eqref{central} is given by
$\overline{\wt\Gamma}=\wt \Gamma$, $\overline{p}=p$, $\overline{\imath}(a)=i(a^{-1})$.\\

Let $\H$ be a complex Hilbert space. We endow the unitary group $U(\H)$
with the strong operator topology and the projective unitary group $PU(\H)=U(\H)/\T$
with the corresponding quotient topology.\\

\definition
A projective unitary representation of $\Gamma$ on $\H$ is a continuous homomorphism
$\pi:\Gamma\rightarrow PU(\H)$. Two such representations $(\H_{i},\pi_{i})$, $i=1,2$
are {\it unitarily equivalent} if there exists a unitary $V:\H_{1}\rightarrow\H_{2}$
such that for any $g\in\Gamma$,
\begin{equation}
V\pi_{1}(g)V^{*}=\pi_{2}(g)\thickspace\text{in $PU(\H_{2})$}
\end{equation}
A projective unitary representation is {\it irreducible} if it leaves no proper,
non--trivial subspace of $\H$ invariant or equivalently if $T\in B(\H)$ commutes
with $\pi$ iff $T$ is a scalar.\\

\remark
\begin{enumerate}
\item Two inequivalent unitary representations of a group $\Gamma$ may become
equivalent when regarded as projective representations. For example, all characters of
$\T=U(1)$ are projectively equivalent. However, this is the case iff the original
representations differ by a character of $\Gamma$.
\item The irreducibility of projective representations is a well--defined concept
since the map $U(\H)\times\B(\H)\rightarrow\B(\H)$,
$(V,T)\rightarrow VTV^{*}$ descends to $PU(\H)\times\B(\H)$ and it therefore
makes sense for an operator, and in particular a projection to commute with
$\pi(g)\in PU(\H)$ for any $g\in \Gamma$.\\
\end{enumerate}

By keeping track of the phases of the operators $\pi(g)$, a projective unitary
representation $(\H,\pi)$ may be viewed as a unitary representation of the group
\begin{equation}
\pi^{*}U(\H)=\{(g,V)\in \Gamma\times U(\H)|\thinspace\pi(g)=[V]\}
\end{equation}
where $[V]$ is the equivalence class of $V$ in $PU(\H)$, by defining $\wt\pi(g,V)=V$.
$\pi^{*}U(\H)$ is in fact the central extension of $\Gamma$ by $\T$ obtained by pulling
back the canonical central extension
\begin{equation}\label{central pu}
1\rightarrow\T\rightarrow U(\H)\rightarrow PU(\H)\rightarrow 1
\end{equation}
The isomorphism class of $\pi^{*}U(\H)$ depends on the particular representation
and it is easy to see that $(\H_{1},\pi_{1})$ and $(\H_{2},\pi_{2})$ are unitarily
equivalent iff the corresponding central extensions are isomorphic and the Hilbert
spaces are unitarily equivalent as representations of either of the extensions.\\

The homomorphism $U(\H_{1})\times U(\H_{2})\rightarrow U(\H_{1}\otimes\H_{2})$,
$(V_{1},V_{2})\rightarrow V_{1}\otimes V_{2}$ descends to the corresponding projective
quotients, and we may therefore define the tensor product $\pi_{1}\otimes\pi_{2}$ of
two projective unitary representations $\pi_{1}$, $\pi_{2}$. The map
$((g,V_{1}),(g,V_{2}))\rightarrow(g,V_{1}\otimes V_{2})$ then yields an isomorphism
\begin{equation}
\pi_{1}^{*}U(\H_{1})\star\pi_{2}^{*}U(\H_{2})\cong
(\pi_{1}\otimes\pi_{2})^{*}U(\H_{1}\otimes\H_{2})
\end{equation}

Similarly, the map $U(\H)\rightarrow U(\overline{\H})$,
$V\rightarrow IVI^{-1}$ where $\overline{\H}$ is $\H$ endowed with the
opposite complex structure and $I:\H\rightarrow\overline{\H}$ is the
canonical anti--linear identification, descends to the projective quotients. We may
therefore define the conjugate $\overline{\pi}$ of a projective representation $\pi$
and the map $(g,V)\rightarrow (g,IVI^{-1})$ gives an isomorphism
\begin{equation}
\overline{\pi^{*}U(\H)}\cong
\overline{\pi}^{*}U(\overline{\H})
\end{equation}

On the other hand, the map
\begin{equation}
U(\H_{1})\times U(\H_{2})\rightarrow U(\H_{1}\oplus\H_{2}),
\quad
(V_{1},V_{2})\rightarrow\begin{pmatrix}V_{1}&0\\0&V_{2}\end{pmatrix}
\end{equation}
does {\bf not} descend to the projective quotients and it is therefore impossible
in general to define the direct sum of projective unitary representations. In fact,
if $(\H_{i},\pi_{i})$ are two projective representations of $\Gamma$ possessing
a direct sum, {\it i.e.}~a projective representation $\pi$ on $\H_{1}\oplus\H_{2}$
leaving the summands invariant and restricting on each $\H_{i}$ to $\pi_{i}$,
then $\pi_{1}^{*}U(\H_{1})\cong\pi_{2}^{*}U(\H_{2})$. The isomorphism is explicitly
given by $(g,V_{1})\rightarrow (g,V_{2})$ where $V_{2}\in U(\H_{2})$ is uniquely
determined by
\begin{xalignat}{3}
[V_{2}]&=\pi_{2}(g)&
&\text{and}&
\left[\begin{pmatrix}V_{1}&0\\0&V_{2}\end{pmatrix}\right]&=\pi(g)
\end{xalignat}
When $\pi_{1}^{*}U(\H_{1})\cong\pi_{2}^{*}U(\H_{2})$, {\it a} direct sum may be
defined as the projectivisation of the direct sum representation
$\wt\pi_{1}\oplus\wt\pi_{2}$ of $\pi_{1}^{*}U(\H_{1})\cong \pi_{2}^{*}U(\H_{2})$.
However, the definition depends upon the isomorphism 
$\pi_{1}^{*}U(\H_{1})\cong \pi_{2}^{*}U(\H_{2})$ the choice of which is
unique only up to a character $\chi\in\Hom(\Gamma,\T)$ and is therefore canonical
only when the latter group is trivial. The following examples should illustrate 
our discussion\\

\begin{enumerate}
\item The trivial and spin $\half{1}$ representations $V_{0},V_{1}$ of $\SU_{2}$ may be
regarded as projective representations of $\SO_{3}=\SU_{2}/\{\pm 1\}$. However, their
direct sum does not factor through $\SO_{3}$ since $\pi_{0}(-1)=1$, $\pi_{1}(-1)=-1$.
The difficulty is that
$\pi_{0}^{*}U(V_{1})\cong\SO_{3}\times\T\ncong
 \SU_{2}\times\T/(-1,-1)\cong\pi_{1}^{*}U(V_{2})$.
\item The projectivisations of
$z\rightarrow\begin{pmatrix}z^{n}&0\\0&1\end{pmatrix}$ give, for varying $n\in\IZ$
inequivalent projective representations of $\T$ which are a direct sum of two copies
of the unique irreducible projective representation of $\T$.
\end{enumerate}

\ssection{Level 1 representations of $LG$}
\label{se:level 1 per}

We consider in this section the lowest energy subspaces of the level 1
positive energy representations of $LG$. When $G=\Spin_{2n}$ or is more
generally {\it simply--laced}, {\it i.e.}~all roots have equal length, these
are exactly the {\it minimal} $G$--modules. Their simple weight structure
and tensor product rules are described in \S \ref{ss:minimal}. They are
used in \S \ref{ss:level 1} to prove that any level $\ell$ representation
of $L\Spin_{2n}$ occurs as a summand in an $\ell$--fold tensor product
$\H_{i_{1}}\otimes\cdots\otimes\H_{i_{\ell}}$ of level 1 representations.

\ssubsection{Lattices and Lie groups}\label{ss:lattices}

We begin by gathering some elementary properties of the lattices canonically
associated to $G$. The present discussion follows \cite{GO1}. Let $T\subset G$
be a maximal torus with Lie algebra $\t\subset\g$. By the {\it roots} of $G$
we shall always mean its infinitesimal roots, namely the set $R$ of linear
forms $\alpha\in i\t^{*}=\Hom(\t,i\IR)$ such that the subspace
\begin{equation}
\g_{\alpha}=\{x\in\gc|\thinspace[h,x]=\alpha(h)x\thickspace\forall h\in\tc\}
\end{equation}
is non--zero. Let $\Delta=\{\alpha_{1},\ldots,\alpha_{n}\}$ be a basis of
$R$ and $\theta$ the corresponding highest root.
The basic inner product $\<\cdot,\cdot\>$ is positive definite on $i\t$ and gives
an identification $i\t^{*}\cong i\t$ of which we shall make implicit use. The
{\it coroots}  of $G$ are the elements of $i\t$ given by
$\alpha^{\vee}=\frac{2\alpha}{\<\alpha,\alpha\>}$. They form the dual root system
$R^{\vee}$.\\

The {\it root} and {\it coroot} lattices $\root\subset i\t^{*}$, $\coroot\subset i\t$ are
the lattices spanned by $R$ and $R^{\vee}$ respectively. They have $\IZ$-basis given by
$\Delta$ and $\Delta^{\vee}=\{\cor{1},\ldots,\cor{n}\}$. Since $\theta$ is a long root and
there are at most two root lengths in $R$ with the ratio of the squared length of a long
root by that of a short root equal to 2 or 3, rewriting
$\alpha^{\vee}=\frac{\<\theta,\theta\>}{\<\alpha,\alpha\>}\alpha$ we see that
$\coroot\subset\root$.
Notice that 
$\<\alpha^{\vee},\alpha^{\vee}\>= \frac{4}{\<\alpha,\alpha\>}=
2\frac{\<\theta,\theta\>}{\<\alpha,\alpha\>}$
so that $\coroot$ is an even, and therefore integral lattice. The {\it weight} and
{\it coweight} lattices $\weight\subset i\t^{*}$,
$\coweight\subset i\t$ are the lattices dual to $\coroot$ and $\root$ respectively.
They have $\IZ$-basis given by the {\it fundamental (co)weights} $\lambda_{i}$,
$\cow{i}$ which are defined by
\begin{equation}
\<\lambda_{i},\cor{j}\>=\<\cow{i},\alpha_{j}\>=\delta_{ij}
\end{equation}
Clearly, $\coweight\subset\weight$. Moreover, by the integrality properties of root systems,
$\<\alpha,\beta^{\vee}\>\in\IZ$ for any root $\alpha$ and coroot $\beta^{\vee}$ so that
$\root\subset\weight$ and, dually, $\coroot\subset\coweight$. Graphically,

\begin{equation}\label{eq:summa}
\begin{diagram}[height=2em,width=3em]
\root	&\subset&\weight	&\subset	&i\t^{*}\\
\cup	&	&\cup		&		&	\\
\coroot	&\subset&\coweight	&\subset	&i\t	\\
\end{diagram}
\end{equation}

Let $Z(G)$ be the centre of $G$ and $\wh{Z(G)}=\Hom(Z(G),\T)$ its Pontriagin dual. Then,

\begin{lemma}\label{centre}\hfill
\begin{enumerate}
\item The map $e(h)=\exp_{T}(-2\pi ih)$ induces an isomorphism $\coweight/\coroot\cong Z(G)$.
\item The pairing $\mu(\exp_{T}(h))=e^{\<\mu,h\>}$ induces an isomorphism
$\weight/\root\cong\wh{Z(G)}$.
\end{enumerate}
\end{lemma}
\proof
(i) Since $G$ is connected, $T$ is maximal abelian and therefore $Z(G)\subset T$.
It follows that $Z(G)\cong e^{-1}(Z(G))/\Ker e$ where the integral lattice
$\Ker e\cong\Hom(\T,T)$ is equal to $\coroot$ since $G$ is simply--connected
\cite[thm. 5.47] {Ad}.
To show that $e(h)\in Z(G)$ iff $h\in\coweight$, we use $Z(G)=\Ker(\Ad_{G})$ and the
fact that if $0\neq x_{\alpha}\in\g_{\alpha}$,
then $\Ad(e(h))x_{\alpha}=\exp(-2\pi i\ad(h))x_{\alpha}=e^{-2\pi i\alpha(h)}x_{\alpha}$
is equal to $x_{\alpha}$ iff $\alpha(h)\in\IZ$.

(ii) The map $\weight/\root\rightarrow\wh{\coweight/\coroot}$,
$\mu\rightarrow e^{-2\pi i\<\mu,\cdot\>}$ is readily seen to be an isomorphism and
coincides with the given pairing under the identification $Z(G)\cong\coweight/\coroot$
\halmos\\

\remark When $G$ is simply--laced, {\it i.e.}~with roots of equal length, the basic inner
product identifies roots and coroots and the vertical inclusions in \eqref{eq:summa}
are equalities. Moreover, lemma \ref{centre} yields a a canonical isomorphism
$\wh{Z(G)}\cong Z(G)$.\\

The Weyl group $W$ of $G$ is the finite group generated in $\End(i\t^{*})$ by the
orthogonal reflections $\sigma_{\alpha}$ corresponding to the roots $\alpha\in R$.
Since
\begin{equation}
 \sigma_{\alpha}(\mu)=
 \mu-2\frac{\<\mu,\alpha\>}{\<\alpha,\alpha\>}\alpha=
 \mu-\<\mu,\alpha^{\vee}\>\alpha=
 \mu-\<\mu,\alpha\>\alpha^{\vee}
\end{equation}
the action of $W$ preserves $\coroot$--cosets in $\coweight$ and $\root$--cosets
in $\weight$. Call $\mu\in\coweight$ (resp. $\mu\in\weight$) {\it minimal} if it
is of minimal length in its $\coroot$ (resp. $\root$)--coset. The following gives
a characterisation of minimal (co)weights.

\begin{proposition}\label{pr:minimal orbit}
There is, in each $\coweight/\coroot$--coset (resp. $\weight/\root$--coset) a unique
$W$--orbit of elements of minimal length. These may equivalently be characterised as
those $\lambda$ such that
\begin{equation}\label{minimal}
\<\lambda,\alpha\>\in\{0,\pm 1\}\qquad
\text{(resp. $\<\lambda,\alpha^{\vee}\>\in\{0,\pm 1\}$)}
\end{equation} for any root $\alpha$ (resp. coroot $\alpha^{\vee}$).
\end{proposition}
\proof
It is sufficient to consider the case of $\coweight/\coroot$ since $\root,\weight$
are the coroot and coweight lattices of the dual root system $R^{\vee}$.
Let $\mu\in\coweight$ be of minimal length in its $\coroot$--coset. Then, for any root
$\beta$ and corresponding coroot $\beta^{\vee}=\frac{2\beta}{\<\beta,\beta\>}$, we have
$\|\mu\pm\beta^{\vee}\|^{2}\geq\|\mu\|^{2}$ and, expanding $|\<\mu,\beta\>|\leq 1$.
Assume that $\lambda\in\coweight$ satisfies \eqref{minimal} and $\nu=\lambda$ mod $\coroot$
is of minimal length in its coset. We claim that $w\lambda=\nu$ for an appropriate
$w\in W$. To see this, write $\nu=\lambda+\beta^{\vee}_{1}+\cdots+\beta^{\vee}_{r}$
where the $\beta^{\vee}_{i}$ are (possibly repeated) coroots. Clearly, one cannot
have $\<\lambda,\beta_{i}\>\geq 0$ for all $i$ otherwise
\begin{equation}
\<\nu,\nu\>=
\<\lambda,\lambda\>+\<\sum\beta^{\vee}_{i},\sum\beta^{\vee}_{i}\> +
2\<\lambda,\sum\beta^{\vee}_{i}\> >
\<\lambda,\lambda\>
\end{equation}
in contradiction with the minimality of $\nu$. Thus, by \eqref{minimal} there exists an
$i\in\{1,\ldots,r\}$ such that $\<\lambda,\beta_{i}\>=-1$ and therefore
$\lambda_{1}:=\sigma_{\beta_{i}}\lambda=\lambda+\beta^{\vee}_{i}$. Moreover, $\lambda_{1}$
satisfies \eqref{minimal} since $W$ permutes the roots and preserves $\<\cdot,\cdot\>$.
We may therefore iterate the above step to find a permutation $\tau$ of $\{1,\ldots,r\}$
such that
\begin{equation}
\lambda_{i}:=
\lambda+\beta^{\vee}_{\tau(1)}+\cdots+\beta^{\vee}_{\tau(i)}=
\sigma_{\beta_{\tau(i)}}\cdots\sigma_{\beta_{\tau(1)}}\lambda
\end{equation}
In particular, $\lambda_{r}=\nu$ and therefore $\nu\in W\lambda$ whence
$\|\lambda\|=\|\nu\|$ \halmos\\

Recall that a weight $\mu\in\weight$ is {\it dominant} if it lies in the cone
\begin{equation}
 \domweight=
 \{\nu\in\weight|\thinspace\<\nu,\cor{i}\>\geq 0
 \thickspace\forall\cor{i}\in\Delta^{\vee}\}=
 \bigoplus\lambda_{i}\IN
\end{equation}
Since $\domweight$ is a fundamental domain for the action of $W$ on $\weight$,
lemma \ref{centre} and proposition \ref{pr:minimal orbit} establish a bijective
correspondence between elements in $\wh{Z(G)}$ and minimal dominant weights.
Dually, the elements of $Z(G)$ correspond to the minimal dominant coweighs,
{\it i.e.}~those $\mu\in\coweight$ of minimal length in their $\coroot$--coset lying
in $\bigoplus\cow{i}\cdot\IN$.

\ssubsection{Minimal $G$--modules}\label{ss:minimal}

An irreducible $G$--module is {\it minimal} if its highest weight $\mu$ is a
minimal (dominant) weight. The following is proved in \cite{GO1}

\begin{proposition}\label{pr:weights of minimal}
Let $V$ be an irreducible $G$--module. Then $V$ is minimal if, and only if its
weights lie in a single orbit of the Weyl group and therefore have multiplicity
one.
\end{proposition}
\proof
We use the fact that the set of weights $\Pi(V)$ of $V$ is the union of the
$W$--orbits of its highest weight $\lambda$ and any other dominant weight
$\mu$ differing from $\lambda$ by a sum of negative roots \cite[\S 21.3]{Hu}.
If $\lambda$ is of minimal length in its $\root$--coset and $\pi\neq 0$ is a
sum of positive roots such that $\lambda-\pi$ is dominant, then
$\<\lambda-\pi,\pi\>\geq 0$ and in particular $\<\lambda,\pi\>>0$. Thus
\begin{equation}
\|\lambda-\pi\|^{2}=\|\lambda\|^{2}-\<\lambda-\pi,\pi\>-\<\pi,\lambda\><
\|\lambda\|^{2}
\end{equation}
a contradiction.
Conversely, assume that $\Pi(V)=W\lambda$ and let $\mu\in\lambda+\root$
be dominant and of minimal length. We claim that $\mu=\lambda$. To see 
this, write $\mu-\lambda=\sum k_{i}\alpha_{i}=\pi-\nu$ where the
$\alpha_{i}$ are simple roots, $\pi=\sum_{i:k_{i}>0} k_{i}\alpha_{i}$
and $\nu=\sum_{i:k_{i}<0}-k_{i}\alpha_{i}$. Notice that $\<\pi,\nu\>\leq 0$
since distinct simple roots form an obtuse angle. It follows that
\begin{equation}
\|\mu\|^{2}=
\|\lambda-\nu\|^{2}+2\<\lambda-\nu,\pi\>+\|\pi\|^{2}\geq
\|\lambda-\nu\|^{2}+\|\pi\|^{2}
\end{equation}
whence $\pi=0$ by minimality of $\mu$. Thus, $\mu=\lambda-\nu$ and therefore
$\nu=0$ since $\mu$ is dominant and $\Pi(V)=W\lambda$ \halmos\\

The following is an instance of the Brauer--Weyl rules for computing
tensor products \mbox{\cite[ex. 9, \S 24.4]{Hu}}

\begin{proposition}\label{pr:tensor with minimal}
Let $V_{\mu},V_{\lambda}$ be the irreducible $G$--modules with highest weights
$\mu,\lambda$. If $\lambda$ is minimal, then
\begin{equation}
V_{\mu}\otimes V_{\lambda}=\bigoplus_{\nu}V_{\mu+\nu}
\end{equation}
where $\nu$ ranges over those weights of $V_{\lambda}$ such that $\mu+\nu$ is
dominant.
\end{proposition}
\proof
Let $\rho=\half{1}\sum_{\alpha>0}\alpha$ be the half--sum of the positive
roots so that $\<\rho,\cor{i}\>=1$ for any simple coroot $\cor{i}$. By the
Weyl character formula, the character of $V_{\mu}$ is
\begin{equation}
\cchi_{\mu}=\frac{A(\mu+\rho)}{\delta}
\end{equation}
where $A(\beta)=\sum_{w\in W}(-1)^{w}e^{w\beta}$ with
$e^{\beta}(\exp_{T}(h))=e^{\beta(h)}$ and
$\delta=\prod_{\alpha>0}(e^{\half{\alpha}}-e^{-\half{\alpha}})$.
On the other hand, by proposition \ref{pr:weights of minimal}
\begin{xalignat}{3}
\cchi_{\lambda}&=\frac{1}{N_{\lambda}}\sum_{w\in W}e(w\lambda)&
&\text{where}&
N_{\lambda}&=|\{w\in W|w\lambda=\lambda\}|
\end{xalignat}
It follows that
\begin{equation}
\cchi_{\mu}\cchi_{\lambda}=
\frac{1}{\delta N_{\lambda}}
\sum_{w'\in W}(-1)^{w'}e(w'(\mu+\rho))\sum_{w\in W}e(w'w\lambda) =
\frac{1}{N_{\lambda}}\sum_{w\in W}\frac{A(\mu+w\lambda+\rho)}{\delta}
\end{equation}
If $\mu+w\lambda$ is not dominant then, for some simple coroot $\cor{i}$,
\begin{equation}
0 > \<\mu+w\lambda,\cor{i}\>=\<\mu,\cor{i}\>+\<w\lambda,\cor{i}\>
\end{equation}
so that $\<\mu,\cor{i}\>=0$ and $\<w\lambda,\cor{i}\>=-1$ since $\mu$ is 
dominant and $w\lambda$ minimal and therefore satisfies \eqref{minimal}.
Thus, $\<\mu+w\lambda+\rho,\cor{i}\>=0$ and if
$\sigma_{i}\in W$ is the simple reflection corresponding to $\cor{i}$, then
\begin{equation}
A(\mu+w\lambda+\rho)=-A(\sigma_{i}(\mu+w\lambda+\rho))=-A(\mu+w\lambda+\rho)=0
\end{equation}
whence, denoting the set of weight of $V_{\lambda}$ by $\Pi(V_{\lambda})$
\begin{equation}
\cchi_{\mu}\cchi_{\lambda} =
\frac{1}{N_{\lambda}}\sum_{w\in W:\mu+w\lambda\in\domweight}\cchi_{\mu+w\lambda} =
\sum_{\nu\in\Pi(V_{\lambda}):\mu+\nu\in\domweight}\cchi_{\mu+\nu}
\end{equation} \halmos

\begin{corollary}\label{co:one dim}
Let $V_{i},V_{k},V_{j}$ be irreducible $G$--modules one of which is minimal.
Then, $\Hom_{G}(V_{i}\otimes V_{k},V_{j})$ is at most one--dimensional.
\end{corollary}
\proof
This follows from proposition \ref{pr:tensor with minimal} if $V_{i}$ or $V_{k}$ are
minimal. If $V_{j}$ is minimal, so is $V_{j}^{*}$ by proposition \ref{pr:weights of 
minimal} and therefore
$\Hom_{G}(V_{i}\otimes V_{k},V_{j})\cong\Hom_{G}(V_{i}\otimes V_{j}^{*},V_{k}^{*})$
is at most one--dimensional \halmos

\ssubsection{Level 1 representations of $L\Spin_{2n}$}\label{ss:level 1}

Consider now $G=\Spin_{2n}$, $n\geq 3$. Let $\H$ be an irreducible, level 1 positive
energy representation of $LG$. Its lowest energy subspace $\H(0)$ is an irreducible
$\Spin_{2n}$--module whose highest weight $\lambda$ satisfies $\<\lambda,\theta\>\leq 1$.
Denote by $\theta_{i}$, $i=1\ldots n$ an orthonormal basis of $\IR^{n}$ and identify
the simple roots of $\Spin_{2n}$ with the vectors 
$\alpha_{i}=\theta_{i}-\theta_{i+1}$, $i=1\ldots n-1$ and
$\alpha_{n}=\theta_{n-1}+\theta_{n}$. The corresponding highest root is
$\theta=\theta_{1}+\theta_{2}$ and, by inspection
\footnote{see also the tables in \S \ref{ss:Z on alcove 2} or \cite{Bou}.}, $\lambda$
is one of $\{0,\lambda_{1},\lambda_{n-1},\lambda_{n}\}$ where $\lambda_{i}$ is the
fundamental weight corresponding to $\alpha_{i}$ so that
\begin{xalignat}{3}
\lambda_{1}  &=v=         	\theta_{1}&
\lambda_{n-1}&=s_{+}=\half{1}(	\theta_{1}+\cdots+\theta_{n-1}-\theta_{n})&
\lambda_{n}  &=s_{-}=\half{1}(	\theta_{1}+\cdots+\theta_{n})
\end{xalignat}
The corresponding irreducible representations are the trivial, vector and spin
representations. Notice that the above weights are minimal. Indeed, since $\Spin_{2n}$
is simply--laced, $\theta^{\vee}=\theta$ is the highest coroot and therefore, for
any coroot $\alpha^{\vee}$
\begin{equation}
\<\lambda,\alpha^{\vee}\>=\<\lambda,\theta\>-\<\lambda,\theta-\alpha^{\vee}\>\leq 1
\end{equation}
since $\theta-\alpha^{\vee}$ is a sum of positive coroots. By proposition
\ref{pr:minimal orbit}, $\lambda$ is a minimal dominant weight. Since the Weyl group
acts by permutation and even numbers of sign changes of the $\theta_{i}$, it follows
by proposition \ref{pr:weights of minimal} that the weights of the corresponding
irreducible representations are
\begin{align}
\Pi(V_{0})&=\{0\}\\
\Pi(V_{v})&=\{\pm\theta_{i}\}\\
\Pi(V_{s_{+}})&=
\{\half{1}\sum\epsilon_{i}\theta_{i}|
  \thinspace\epsilon_{i}\in\{\pm 1\}
  \thickspace\text{and $\epsilon_{i}=-1$ for an odd number of $i$}\}\\
\Pi(V_{s_{-}})&=
\{\half{1}\sum\epsilon_{i}\theta_{i}|
  \thinspace\epsilon_{i}\in\{\pm 1\}
  \thickspace\text{and $\epsilon_{i}=-1$ for an even number of $i$}\}
\end{align}
and that their tensor product rules are governed by proposition \ref{pr:tensor with minimal}.

\begin{lemma}\label{le:l lemma}
An irreducible $\Spin_{2n}$--module $V$ is admissible at level $\ell$
if, and only if it arises as a summand of an $\ell$--fold tensor product
$V_{i_{1}}\otimes\cdots\otimes V_{i_{\ell}}$ whose factors are admissible
at level 1.
\end{lemma}
\proof
Let $\lambda$ be the highest weight of $V$ and denote by $V_{\mu}$ the
irreducible $G$--module with highest weight $\mu$. If $V$ is contained
in $V_{\lambda_{1}}\otimes\cdots\otimes V_{\lambda_{\ell}}$
where $\<\lambda_{i},\theta\>\leq 1$, then $\lambda=\sum_{i}\lambda_{i}-\pi$
where $\pi$ is a sum of positive roots and therefore
$\<\lambda,\theta\>\leq\ell$.
Conversely, assume that $\<\lambda,\theta\>\leq\ell$. If the inequality
is strict, then by induction,
$V\subset
 V_{\lambda_{1}}\otimes\cdots\otimes V_{\lambda_{\ell-1}}\otimes\IC$
with $\<\lambda_{i},\theta\>\leq 1$.
Suppose therefore that $\<\lambda,\theta\>=\ell\geq 2$ and write
$\lambda=\sum_{i}\lambda_{i}\theta_{i}$. If $\lambda_{2}=0$, then
$\mu=\lambda-\theta_{1}$ is dominant and satisfies
$\<\mu,\theta\>\leq\ell-1$. Moreover, $V \subset V_{\mu}\otimes V_{\theta_{1}}$.
If, on the other hand, $\lambda_{2}>0$ and the coordinates of $\lambda$ are of the
form $\lambda_{1}=\lambda_{2}=\cdots=\lambda_{m}>\lambda_{m+1}\geq\cdots\geq\lambda_{n}$
then $\mu=\lambda-\half{1}(\theta_{1}+\cdots+\theta_{m}-\theta_{m+1}-\cdots-\theta_{n})$
is dominant and satisfies $\<\mu,\theta\>\leq\ell-1$. Since
$\half{1}(\theta_{1}+\cdots+\theta_{m}-\theta_{m+1}-\cdots-\theta_{n})$ is a weight
of one of the two spin representations, call it $V_{\sigma}$, proposition
\ref{pr:tensor with minimal} yields $V\subset V_{\mu}\otimes V_{\sigma}$.
The result now follows by induction \halmos

\begin{lemma}\label{le:generation}
Let $(\pi,\H)$ be a positive energy representation of $LG$ and $V\subset\H(0)$
an irreducible $G$-module. Then, the closure of the span of $\pi(LG)V$ is
unitarily equivalent to $\H_{V}$.
\end{lemma}
\proof
Let $\K\subset\Hfin$ be the subspace generated by $V$ under the action of $\lpol\gc$.
Its lowest energy subspace is $V$ since any element in the enveloping algebra of
$\lpol\gc$ may be written as a sum of monomials of the form $X^{-}X^{0}X^{+}$ where
the $X^{\mp}$ are a product of $x(n)$, with $n\lessgtr 0$ and $X^{0}\in\mathfrak{U}\g$,
and $V$ is annihilated by $X^{+}$ and left invariant by $X^{0}$. It follows that
$\K$ is an irreducible $\lpol\g\rtimes\IC d$--module. Indeed, any submodule $\Nu\subset\K$
is necessarily graded and the corresponding orthogonal projection $P$ commutes with
$\rot$. It therefore maps $V=\K(0)$ into $V$ and since $P$ commutes with $\g$ we have
$PV=0$ or $V$. Thus, $\Nu=P\U\lpol\gc V=\U\lpol\gc PV$ is $0$ or $\K$.
By proposition \ref{simple}, $\overline{\K}$ is invariant and irreducible under $LG$
since its finite energy subspace is $\K$ and it follows that it is unitarily equivalent
to $\H_{V}$. Since $\overline{\K}$ contains the closure of the span of $\pi(LG)V$ it
coincides with it \halmos\\

The following useful result is due to Pressley and Segal \cite[prop. 9.3.9]{PS}

\begin{proposition}\label{pr:per l lemma}
Any irreducible, level $\ell$ positive energy representation $\H$ of $L\Spin_{2n}$
is a summand in an $\ell$--fold tensor product
$\H_{i_{1}}\otimes\cdots\otimes\H_{i_{\ell}}$ of level 1 representations.
\end{proposition}
\proof
By lemma \ref{le:l lemma}, the lowest energy subspace $\H(0)$ of $\H$ is
contained in some tensor product $V_{i_{1}}\otimes\cdots\otimes V_{i_{\ell}}$
of minimal representations of $\Spin_{2n}$. Let $\H_{i_{k}}$ be the irreducible
level 1 representations whose lowest energy subspace are the $V_{i_{k}}$. 
The lowest energy subspace of the level $\ell$ representation
$\H_{i_{1}}\otimes\cdots\otimes\H_{i_{\ell}}$ contains $\H(0)$ and therefore,
by lemma \ref{le:generation}, the closure of the linear span of $L\Spin_{2n}\H(0)$
inside $\H_{i_{1}}\otimes\cdots\otimes\H_{i_{\ell}}$ is isomorphic to $\H$ \halmos\\

\remark
Lemma \ref{le:l lemma} holds for any simply--connected classical Lie group $G$.
This may be established for $G=\Spin_{2n+1}$ by an almost identical proof to
that of lemma \ref{le:l lemma} and for $\SU_{n}$ and $\Sp_{n}$ by simply noticing
that all fundamental weights satisfy $\<\lambda,\theta\>\leq 1$. Thus, proposition
\ref{pr:per l lemma} holds for the loop groups of all classical Lie groups.
Notice that lemma \ref{le:l lemma} does not hold for $G=E_{8}$ since the only
irreducible module admissible at level 1 is the trivial representation.

\ssection{Discontinuous loops and outer automorphisms of $LG$}
\label{se:disc loops}

We consider in this section the automorphic action on $LG$ of the group of
{\it discontinuous loops}
\begin{equation}\label{eq:disc}
\lzg=\{f\in C^{\infty}(\IR,G)|\thinspace f(x+2\pi)f(x)^{-1}\in Z(G)\}
\end{equation}
We show in \S \ref{ss:Z on per} that the category $\pl$ of positive energy
representations of $LG$ at level $\ell$ is closed under conjugation by
$\lzg$. The corresponding abstract action of $Z(G)=\lzg/LG$ on the level
$\ell$ alcove of $G$ which parametrises the irreducibles in $\pl$ coincides 
with the geometric one obtained by realising $Z(G)$ as a distinguished
subgroup of the automorphisms of the extended Dynkin diagram of $G$. We
begin by studying the latter. An explicit description of this action
according to the Lie type of $G$ may be found in \S \ref{ss:Z on alcove 2}.

\ssubsection{Geometric action of $Z(G)$ on the level $\ell$ alcove}
\label{ss:Z on alcove}

This subsection is essentially an expanded version of \cite[ch. VI, \S 2.3]{Bou}.
The notation follows that of the section \ref{se:level 1 per}.

\begin{lemma}\label{special roots}
There is a bijective correspondence between elements of $Z(G)\backslash\{1\}$ and
fundamental coweights corresponding to special roots, {\it i.e.}~the
$\alpha_{i}\in\Delta$ bearing the coefficient 1 in the expansion
\begin{equation}\label{highest root}
\theta=\sum m_{i}\alpha_{i}
\end{equation}
\end{lemma}
\proof
By proposition \ref{pr:minimal orbit}, $\mu\in(\coweight)^{+}$ is minimal iff
$\<\mu,\theta\>\leq 1$. Indeed, for any positive root $\alpha$, we get
$0\leq\<\mu,\alpha\>\leq\<\mu,\theta\>-\<\mu,\theta-\alpha\>\leq\<\mu,\theta\>$.
Since $\<\mu,\theta\>=0$ implies $\mu=0$, the non--zero minimal dominant coweights
are those $\mu\in(\coweight)^{+}$ such that $\<\mu,\theta\>=1$. Writing
$\mu=\sum_{i}k_{i}\cow{i}$, $k_{i}\geq 0$ and using \eqref{highest root}, we find
$\<\mu,\theta\>=\sum k_{i}m_{i}$. Since $\theta-\alpha_{i}$ is a sum of positive
roots, $m_{i}\geq 1$ for any $i$ and result follows \halmos\\

Denote $-\theta$ by $\alpha_{0}$, then

\begin{lemma}
For any special root $\alpha_{i}$, the set
$\Delta_{i}=\Delta\backslash\{\alpha_{i}\}\cup\{\alpha_{0}\}$ is a basis of $R$
with highest root $-\alpha_{i}$ and dual basis
\begin{align}
\cowp{0}&=-\cow{i} \label{dual one}\\
\cowp{j}&= \cow{j}-\<\theta,\cow{j}\>\cow{i} \label{dual two}
\end{align}
\end{lemma}
\proof
Let $x\in\tc$, then
\begin{equation}
 x=\sum\<x,\cow{j}\>\alpha_{j}
  =\<x,\cow{i}\>\theta+\sum_{j\neq i}
   (\<x,\cow{j}\>-\<x,\cow{i}\>\<\theta,\cow{j}\>)\alpha_{j}
\end{equation}
so that $\Delta_{i}$ is a vector space basis of $\tc$ with dual basis given by
\eqref{dual one}--\eqref{dual two}.
If $0<\beta\in R$ then either $\<\beta,\cow{i}\>=0$ in which case 
$\<\beta,\cowp{0}\>$ and $\<\beta,\cowp{j}\>$ are all non-negative or
$\<\beta,\cow{i}\>=1$ since $\cow{i}$ is a minimal dominant coweight. In the latter
case $\<\beta,\cowp{0}\>=-1$ and $\<\beta,\cowp{j}\>=\<\beta-\theta,\cow{j}\>\leq 0$.
Thus, $\Delta_{i}$ is a basis of $R$.
Next, for any $\beta\in R$, \eqref{minimal} yields
$\<-\alpha_{i}-\beta,\cowp{0}\>=1+\<\beta,\cow{i}\>\geq 0$ since $\cow{i}$ is minimal.
Moreover, for $j\neq i$
\begin{equation}\label{positive}
 \<-\alpha_{i}-\beta,\cowp{j}\>=
 \<\theta,\cow{j}\>-\<\beta,\cow{j}\>+\<\theta,\cow{j}\>\<\beta,\cow{i}\>=
 \<\theta-\beta,\cow{j}\>+\<\theta,\cow{j}\>\<\beta,\cow{i}\>
\end{equation}
The above is clearly non-negative if $\<\beta,\cow{i}\>\geq 0$. If, on the other hand
$\<\beta,\lambda_{i}\vvee\>=-1$ so that $\beta<0$ then \eqref{positive} is equal to
$-\<\beta,\lambda_{j}\vvee\>\geq 0$ and $-\alpha_{i}$ is the highest root relative to
$\Delta_{i}$ \halmos

\begin{proposition}\label{centre weyl}
Let $\deltabar=\Delta\cup\{\alpha_{0}\}$. Then, for any special root $\alpha_{i}$,
there exists a unique
\begin{equation}
w_{i}\in W_{0}=\{w\in W|\thinspace w\deltabar=\deltabar\}
\end{equation}
such that $w\alpha_{0}=\alpha_{i}$. The resulting map $\i:Z(G)\rightarrow W_{0}$
obtained by identifying $Z(G)\setminus\{1\}$ with the set of special roots is an
isomorphism.
\end{proposition}
\proof
The existence of $w_{i}$ follows from the previous lemma since $W$ acts transitively
on the set of basis of $R$ and maps highest roots to highest roots.
$w_{i}$ is unique because an element $w\in W_{0}$ is determined by $w\alpha_{0}$.
Indeed, if $w_{1}\alpha_{0}=\alpha_{j}=w_{2}\alpha_{0}$, then $w_{2}^{-1}w_{1}$ is a
permutation of $\Delta$ and is therefore the identity since $W$ acts simply on basis.
$\i$ is injective because $w_{i}\alpha_{0}=\alpha_{i}$. Let now $w\in W_{0}$. We claim
that $\alpha_{i}=w\alpha_{0}$ is a special root. It then follows by uniqueness that
$w=w_{i}$ and therefore that $\i$ is surjective. To see this, we apply $w$ to
\eqref{highest root} and get $-\alpha_{i}=\sum_{j}m_{j}w\alpha_{j}$ while at the same
time $-\alpha_{i}=m_{i}^{-1}(\alpha_{0}+\sum_{j\neq i}m_{j}\alpha_{j})$. Comparing
the coefficients of $\alpha_{0}$ we get $m_{i}=1$.
To prove that $\i$ is a homomorphism, let $\alpha_{i}$ and $\alpha_{j}$ be special
roots. Then either $w_{i}w_{j}=1$ or $w_{i}w_{j}=w_{k}$ where $\alpha_{k}$ is another
special root. In the former case, $w_{i}\alpha_{j}=\alpha_{0}$ and therefore, by
\eqref{dual one}
\begin{equation}\label{ccompo}
w_{i}\cow{j}=\cowp{0}=-\cow{i}
\end{equation}
so that $\cow{j}=-\cow{i}$ mod $\coroot$ since $W$ leaves $\coweight/\coroot$ cosets
invariant. In the latter, $w_{i}\alpha_{j}=\alpha_{k}$ and therefore, using
\eqref{dual two}
\begin{equation}\label{compo}
w_{i}\cow{j}=\cowp{k}=\cow{k}-\<\theta,\cow{k}\>\cow{i}=\cow{k}-\cow{i}
\end{equation}
whence $\cow{i}+\cow{j}=\cow{k}$ mod $\coroot$ \halmos\\

The following is well--known and often rediscovered \cite{OT}

\begin{corollary}
$Z(G)$ is canonically isomorphic to the group of automorphisms of the extended
Dynkin diagram of $G$ induced by Weyl group elements.
\end{corollary}

\begin{proposition}\label{geometric action}
For any $\ell\in\IN$, there is a canonical action of $Z(G)$ on the level $\ell$ alcove
$\al$ given by
\begin{equation}\label{eq:action of Z}
z\longrightarrow A_{i}=\tau(\ell\cow{i})w_{i}
\end{equation}
where $i$ is the index of the special root corresponding to $z$ via lemma
\ref{special roots}, $\tau$ denotes translation and $w_{i}=\i(z)$ corresponds
to $z$ via proposition \ref{centre weyl}.
\end{proposition}
\proof
The level $\ell$ alcove is given by
\begin{equation}
\al=\{\lambda\in\weight|\<\lambda,\alpha_{i}\>\geq 0,\medspace\<\lambda,\theta\>\leq\ell\}
\end{equation}
If $\lambda\in\al$ then for $j\neq i$,
$\<A_{i}\lambda,\alpha_{j}\>=\<\lambda,w_{i}^{-1}\alpha_{j}\>\geq 0$ since
$w_{i}^{-1}\alpha_{j}\neq \alpha_{0}$. On the other hand,
$\<A_{i}\lambda,\alpha_{i}\>=\ell+\<\lambda,\alpha_{0}\>\geq 0$. Finally,
$\<A_{i}\lambda,\theta\>=\ell-\<\lambda,w_{i}^{-1}\theta\>\leq\ell$ so that the $A_{i}$
leave $\al$ invariant.
Next, $A_{i}A_{j}=\tau(\ell(\cow{i}+w_{i}\cow{j}))w_{i}w_{j}$. If $w_{i}w_{j}=1$, we get
by \eqref{ccompo} and the previous proposition $A_{i}A_{j}=1$. If on the other hand
$w_{i}w_{j}=w_{k}$, \eqref{compo} yields $A_{i}A_{j}=A_{k}$ \halmos

\begin{corollary}\label{co:level 1 alcove}
If $G$ is simply--laced, the canonical action of $Z(G)$ on the level 1 alcove
\begin{equation}
{\mathcal A}_{1}\cong\wh{Z(G)}\cong Z(G)
\end{equation}
coincides with left multiplication. In particular, it is transitive and free.
\end{corollary}
\proof
Let $\lambda\in\AA_{1}$. Since $G$ is simply--laced, $\theta^{\vee}=\theta$
is the highest coroot and therefore, for any coroot $\alpha^{\vee}$,
\begin{equation}
\<\lambda,\alpha^{\vee}\>=
\<\lambda,\theta\>-\<\lambda,\theta-\alpha^{\vee}\>\leq
\<\lambda,\theta\>\leq
1
\end{equation}
since $\theta-\alpha^{\vee}$ is a sum of positive coroots. Thus, by proposition
\ref{pr:minimal orbit} the points in $\AA_{1}$ are the minimal dominant weights
and therefore are in one--to--one correspondence with elements in $\wh{Z(G)}$.
As previously remarked, the basic inner product identifies $\coroot,\coweight$
with $\root,\weight$ respectively and therefore determines a natural isomorphism
$Z(G)\cong\wh{Z(G)}$. Let now $z\in Z(G)\setminus\{1\}$
correspond to the fundamental coweight $\cow{i}$. Then, if $\lambda\in\AA_{1}$,
we have by \eqref{eq:action of Z},
\begin{equation}
A_{i}\lambda=
\cow{i}+w_{i}\lambda=
\cow{i}+\lambda+(w_{i}\lambda-\lambda)
\end{equation}
Since $W$ preserves $\root$--cosets in $\weight$, $A_{i}\lambda=\cow{i}+\lambda
\mod\root$ as claimed \halmos

\ssubsection{Action of $Z(G)$ on the positive energy representations}
\label{ss:Z on per}

The group of discontinuous loops $\lzg$ defined by \eqref{eq:disc} acts
automorphically on $LG$ by conjugation. This induces an action on the set
of equivalence classes of projective unitary representations of $LG$ given by
$\zeta_{*}\pi(\gamma)=\pi(\zeta^{-1}\gamma\zeta)$ for any $\zeta\in\lzg$. Clearly,
$(\zeta_{1}\zeta_{2})_{*}\pi={\zeta_{1}}_{*}{\zeta_{2}}_{*}\pi$ and this action
factors through $Z=\lzg/LG$ since
$\gamma_{*}\pi(\cdot)=\pi(\gamma)^{*}\pi(\cdot)\pi(\gamma)$ whenever $\gamma\in LG$.
Let $(\pi,\H)$ be irreducible and of positive energy and denote by $U_{\theta}$
the corresponding action of $\rot$. For a fixed $\zeta\in\lzg$, any intertwining
action of $\rot$ for $\zeta_{*}\pi$, whether of positive energy or not is necessarily
given by 
\begin{equation}\label{twisted rotations}
V_{\theta}=\pi(\zeta^{-1}\zeta_{\theta})U_{\theta}
\end{equation}
Indeed, $\zeta^{-1}\zeta_{\theta}\in LG$ and $V_{\theta}$ yields a projective action of
$\rot$ satisfying
$V_{\theta}\zeta_{*}\pi(\gamma)V_{\theta}^{*}=\zeta_{*}\pi(\gamma_{\theta})$.
Moreover, if $V_{\theta}^{i}$, $i=1,2$ are two intertwining actions of $\rot$, then
$W_{\theta}=(V_{\theta}^{1})^{*}V_{\theta}^{2}$ commutes projectively with the
action of $LG$ so that the following holds in $U(\H)$ for any $\theta$ :
$W_{\theta}\pi(\gamma)W_{\theta}^{*}\pi(\gamma)^{*}=\chi(\gamma,\theta)$,
where $\chi(\gamma,\theta)\in\T$ and depends multiplicatively on $\gamma\in LG$.
Since $LG$ is equal to its commutator subgroup \cite[propn. 3.4.1.]{PS},
$\chi\equiv 1$ and it follows by Shur's lemma that  $W_{\theta}=1$ in $PU(\H)$.
Thus, $\zeta_{*}\pi$ is of positive energy iff $V_{\theta}$ is a positive energy
representation of $\rot$.

\begin{proposition}
If $(\pi,\H)$ is a positive energy representation of $LG$ and $\zeta\in\lzg$,
the conjugated representation $\zeta_{*}\pi$ is of positive energy.
\end{proposition}
\proof
It is sufficient to prove the above for a given set of representatives of $LG$--cosets
in $\lzg$. A particularly convenient choice is obtained via lemma \ref{centre} by
considering the discontinuous loops
$\zeta_{\mu}(\theta)=\exp_{T}(-i\theta\mu)$ where $\mu\in\coweight\subset i\t$ is
a coweight. If $\mu\in\coroot$, then $\zeta_{\mu}\in LG$ and the
corresponding action of $\rot$ may be rewritten, by \eqref{twisted rotations} as
$\pi(\zeta^{-1}_{\mu}{\zeta_{\mu}}_{\theta})U_{\theta}=
 \pi(\zeta_{\mu})^{*}U_{\theta}\pi(\zeta_{\mu})$ which is clearly of positive
energy. The general case $\mu\in\coweight$ is settled by the following simple
observation. Notice first that
$\zeta_{\mu}^{-1}{\zeta_{\mu}}_{\theta}=\zeta_{\mu}(-\theta)=\exp_{T}(i\theta\mu)$
since $\zeta_{\mu}$ is a homomorphism and write $\mu$ as a convex combination of
elements in the coroot lattice,
$\mu=\sum_{i=1}^{m} t_{i}\alpha_{i}$,
$t_{i}\in(0,1]$, $\sum_{i} t_{i}=1$, $\alpha_{i}\in\coroot$.
Recall that on $\H$ we have a unitary lift $\wt\pi$ of $\pi$ over $T\times\rot$.
Therefore,
\begin{equation}\label{convex}
\wt\pi(\exp_{T}(i\theta\mu))\wt\pi(\exp_{\rot}(i\theta d))=
\prod_{j}
\wt\pi(\exp_{T}(i\theta t_{j}\alpha_{j}))
\wt\pi(\exp_{\rot}(i\theta t_{j} d))
\end{equation}
is a lift of $\pi(\zeta_{\mu}^{-1}{\zeta_{\mu}}_{\theta})U_{\theta}$ and the
product of $m$ commuting representations of $\IR$ which by our previous argument
are of positive energy. It follows that
$\pi(\zeta_{\mu}^{-1}{\zeta_{\mu}}_{\theta})U_{\theta}$ is of positive energy \halmos

\begin{proposition}\label{conjugated action}
Let $(\pi,\H)$ be a positive energy representation of $LG$ of level $\ell$ and
$\zeta\in\lzg$. Then,
\begin{enumerate}
\item $\zeta_{*}\pi$ is of level $\ell$.
\item If $\zeta_{\mu}(\phi)=\exp_{T}(-i\phi\mu)$ is the discontinuous
loop corresponding to $\mu\in\coweight$, the subspaces of finite energy vectors
of $\pi$ and $\zeta_{*}\pi$ coincide.
\item If $\Ad(\zeta)\lpol\g=\lpol\g$ and the finite energy subspaces of $\pi$ and
$\zeta_{*}\pi$ coincide, the conjugated action of $\lpol\g$ on $\hfin$ is given by
\begin{equation}\label{twisted lg}
\zeta_{*}\pi(X)=
\pi(\zeta^{-1} X\zeta)+i\ell\int_{0}^{2\pi}\<\dot\zeta\zeta^{-1},X\>\frac{d\theta}{2\pi}
\end{equation}
\end{enumerate}
\end{proposition}
\proof
It is sufficient to check (i) for the discontinuous loops $\zeta_{\mu}$,
$\mu\in\coweight$. This will be done in the course of the proof of (iii).

(ii)
It was remarked in the proof of the previous proposition that the conjugated
action of rotations \eqref{twisted rotations} corresponding to $\zeta_{\mu}$
is given by
\begin{equation}
 \pi(\zeta^{-1}_{\mu}{\zeta_{\mu}}_{\theta})U_{\theta}=
 \pi(\exp_{T}(i\theta\mu))\pi(\exp_{\rot}(i\theta d))
\end{equation}
which commutes with the original action of $\rot$ given by
$U_{\theta}=\pi(\exp_{\rot}(i\theta d))$. Since both are of
positive energy, their finite energy subspaces coincide.

(iii)
Let $h$ be the level of $\zeta_{*}\pi$ and denote by $\pi$ and $\zeta_{*}\pi$
the projective representations of $\lpol\g$ on $\hfin$ given by theorem \ref{core}
so that
\begin{align}
[\pi(X),\pi(Y)]&=\pi([X,Y])+i\ell B(X,Y)\\
[\zeta_{*}\pi(X),\zeta_{*}\pi(Y)]&=\zeta_{*}\pi([X,Y])+ihB(X,Y) \label{new level}
\end{align}
where $B(X,Y)=\int_{0}^{2\pi}\<X,\dot Y\>\frac{d\theta}{2\pi}$ is the fundamental
2-cocycle on $L\g$. Evidently, $\zeta_{*}\pi(X)=\pi(\zeta^{-1}X\zeta)+iF(X)$
for some $F(X)\in\IR$ since
$\zeta_{*}\pi(\exp_{LG}(X))= \pi(\exp_{LG}(\zeta^{-1}X\zeta))=
 e^{\pi(\zeta^{-1}X\zeta)}$ in $PU(\H)$. It follows that
\begin{equation}
\begin{split}
[\zeta_{*}\pi(X),\zeta_{*}\pi(Y)]
&=
[\pi(\zeta^{-1}X\zeta),\pi(\zeta^{-1}Y\zeta)]\\
&=
\pi(\zeta^{-1}[X,Y]\zeta)+
i\ell\int_{0}^{2\pi}
\<\zeta^{-1}X\zeta,
  \zeta^{-1}\dot Y\zeta\>
\frac{d\theta}{2\pi}+
i\ell\int_{0}^{2\pi}
\<\zeta^{-1}X\zeta,
  \dot{(\zeta^{-1})}Y\zeta+\zeta^{-1}Y\dot\zeta\>
\frac{d\theta}{2\pi}\\
&=
\zeta_{*}\pi([X,Y])-iF([X,Y])+
i\ell B(X,Y)+
i\ell\int_{0}^{2\pi}
\<\zeta^{-1}X\zeta,
  \zeta^{-1}[Y,\dot\zeta\zeta^{-1}]\zeta\>
\frac{d\theta}{2\pi}\\
&=
\zeta_{*}\pi([X,Y])-idF(X,Y)+i\ell B(X,Y)+
i\ell\int_{0}^{2\pi}
\<\dot\zeta\zeta^{-1},[X,Y]\>
\frac{d\theta}{2\pi}
\end{split}
\end{equation}
where we used the $\Ad_{G}$ invariance of $\<\cdot,\cdot\>$ and the fact
that $\dot{(\zeta^{-1})}\zeta+\zeta^{-1}\dot\zeta=
      \dot{(\zeta^{-1}\zeta)}=0$.
Since $hB$ and $\ell B$ define the same cohomology class iff $h=\ell$, we find by
equating the above with \eqref{new level} that the level of $\zeta_{*}\pi$ is $\ell$.
Moreover, \eqref{twisted lg} holds since $[\lpol\g,\lpol\g]=\lpol\g$ \halmos

\begin{theorem}\label{th:Z on per}
Let $(\pi,\H)$ be an irreducible positive energy representation of $LG$ of
level $\ell$ and highest weight $\lambda$ and $\zeta\in\lzg$. Then, the conjugated
representation $\zeta_{*}\pi(\gamma)=\pi(\zeta^{-1}\gamma\zeta)$ on $\H$ is of positive
energy, level $\ell$ and highest weight $\zeta\lambda$ where the notation refers to the
geometric action of $Z(G)=\lzg/LG$ on the level $\ell$ alcove defined by proposition
\ref{geometric action}.
\end{theorem}
\proof
As customary, it is sufficient to prove the result for a given choice of representatives
of $LG$--cosets in $\lzg$. Let $z\in Z(G)\backslash\{1\}$ correspond to the special root
$\alpha_{j}$ by lemma \ref{special roots} and consider the discontinuous loop
$\zeta=\zeta_{\cow{j}}w_{j}$
where $\cow{j}$ is the associated fundamental coweight and $w_{j}\in G$ a representative
of the Weyl group element corresponding to $z$ by proposition \ref{centre weyl}.
Since the action of $G$ commutes with $\rot$ on $\H$, the subspace of finite energy
vectors of $\zeta_{*}\pi$ coincides with that of ${\zeta_{\cow{j}}}_{*}\pi$ and, by the
previous proposition, with that of $\pi$. We may therefore compare the infinitesimal
actions of $\lpol\gc$ corresponding to $\pi$ and $\zeta_{*}\pi$ on $\Hfin$.\\

If $\alpha$ is a root and $x_{\alpha}\in\g_{\alpha}$, we have
$[\cow{j},x_{\alpha}]=\<\cow{j},\alpha\>x_{\alpha}$. Since
$\zeta_{\cow{j}}(\theta)=\exp_{T}(-i\cow{j}\theta)$, this gives
\begin{equation}
 \zeta_{\cow{j}}^{-1}x_{\alpha}(n)\zeta_{\cow{j}}(\theta)=
 x_{\alpha}\otimes e^{i\theta(n+\<\cow{j},\alpha\>)}=
 x_{\alpha}(n+\<\cow{j},\alpha\>)(\theta)
\end{equation}
Thus, $\zeta x_{\alpha}(n)\zeta^{-1}$ lies in the root space
$\g_{w_{j}^{-1}\alpha,n+\<\cow{j},\alpha\>}$ and therefore, up to a non--zero
multiplicative constant
\begin{equation}\label{new e alpha}
\zeta_{*}\pi(e_{\alpha}(n))=\pi(e_{w_{j}^{-1}\alpha}(n+\<\cow{j},\alpha\>))
\end{equation}
since no additional term arise from \eqref{twisted lg} because $\dot\zeta\zeta^{-1}=-i\cow{j}$
lies in $\tc$ which is orthogonal to $\g_{\alpha}$.
If, on the other hand $h\in\tc$, then $\zeta^{-1}h(n)\zeta=w_{j}^{-1}h(n)$ and \eqref{twisted lg}
reads
\begin{equation}\label{new hn}
\zeta_{*}\pi(h(n))=\pi(w_{j}^{-1}h(n))+\ell\delta_{n,0}\<h,\cow{j}\>
\end{equation}
Let $\Omega\in\hfin$ be the highest weight vector for $\zeta_{*}\pi$. We claim that, up to a
scalar factor, $\Omega=\Upsilon$, the highest weight vector for $\pi$. To see this, recall
that $\Omega$ is the unique element of $\hfin$ annihilated by the subalgebra spanned by the
$x(n)$, $x\in\gc$ and $n>0$ and the $x_{\alpha}(0)$ with
$\alpha>0$. This in turn is generated by the elements corresponding the simple affine roots,
namely $e_{\alpha_{i}}(0)$ and $e_{\alpha_{0}}(1)$ where $\alpha_{0}=-\theta$. Recalling from
proposition \ref{centre weyl} that $w_{j}^{-1}$ acts as a permutation of
$\deltabar=\{\alpha_{0},\ldots,\alpha_{n}\}$ and maps $\alpha_{j}$ to $\alpha_{0}$, we get,
using \eqref{new e alpha}
\begin{align}
\zeta_{*}\pi(e_{\alpha_{k}}(0))
&=
\left\{\begin{array}{ll}
\pi(e_{w_{j}^{-1}\alpha_{k}}(0))&\text{if $k\neq j$}\\
\pi(e_{\alpha_{0}}(1))		&\text{if $k=j$}
\end{array}\right.\\
\zeta_{*}\pi(e_{\alpha_{0}}(1))
&=\pi(e_{w_{j}^{-1}\alpha_{0}}(0))
\end{align}
whence $\Omega=\Upsilon$.
To find the weight of $\Omega$ and therefore the highest weight of $\zeta_{*}\pi$, we use
\eqref{new hn} and the fact that $\pi(h(0))\Upsilon=\<\lambda,h\>\Upsilon$ whenever $h\in\tc$
\begin{equation}
\zeta_{*}\pi(h(0))\Omega=
\<h,w_{j}\lambda+\ell\cow{j}\>\Omega=
\<h,\zeta\lambda\>\Omega
\end{equation}
\halmos\\


We now derive a number of simple corollaries of the above results. First, the level
of a positive energy representation may be detected globally in view of the following 

\begin{corollary}\label{co:global level}
Let $(\pi,\H)$ be a level $\ell$ positive energy representation of $LG$. Then, for
any $\tau\in T$ and coroot $\alpha\in\coroot$ 
\begin{equation}\label{eq:global level}
\pi(\zeta_{\alpha})\pi(\tau)\pi(\zeta_{\alpha})^{*}\pi(\tau)^{*}=\alpha(\tau)^{-\ell}
\end{equation}
where $\zeta_{\alpha}(\phi)=\exp_{T}(-i\alpha\phi)$ and the right hand--side of
\eqref{eq:global level} refers to the canonical pairing $\coroot\times
T\rightarrow \T$, $(\alpha,\exp_{T}(h))\rightarrow e^{\<\alpha,h\>}$.
\end{corollary}
\proof
Since $\zeta_{\alpha}$ and $T$ commute in $LG$, the following holds in $PU(\H)$ for
any $h\in\t$ and $s\in\IR$
\begin{equation}
\pi(\zeta_{\alpha})e^{s\pi(h)}\pi(\zeta_{\alpha})^{*}=
\pi(\zeta_{\alpha})\pi(\exp_{T}(sh))\pi(\zeta_{\alpha})^{*}=
\pi(\exp_{T}(sh))=
e^{s\pi(h)}
\end{equation}
Thus, in $U(\H)$,
\begin{equation}\label{eq:global coc}
\pi(\zeta_{\alpha})e^{s\pi(h)}\pi(\zeta_{\alpha})^{*}=
\lambda(s)e^{s\pi(h)}
\end{equation}
where $\lambda:\IR\rightarrow\T$ is a continuous homomorphism and is therefore of the 
form $e^{s c}$ for some $c\in i\IR$. Let now $\xi\in\Hfin$. By (ii) of proposition
\ref{conjugated action}, $\pi(\zeta_{\alpha})^{*}\xi\in\Hfin$ so that applying both
sides of \eqref{eq:global coc} to $\xi$ and differentiating at $s=0$, we find
\begin{equation}\label{eq:infi coc}
\pi(\zeta_{\alpha})\pi(h)\pi(\zeta_{\alpha})^{*}\xi=\pi(h)\xi+c\xi
\end{equation}
The result now follows by comparison with (iii) of proposition \ref{conjugated action}.
Indeed, $\zeta_{\alpha}\in LG$ and therefore
${\zeta_{\alpha}^{-1}}_{*}\pi(X)=\pi(\zeta_{\alpha})\pi(X)\pi(\zeta_{\alpha})^{*}$
for any $X\in\lpol\g$. Since $\dot{(\zeta_{\alpha}^{-1})}\zeta_{\alpha}=i\alpha$, we
find by equating \eqref{eq:infi coc} and \eqref{twisted lg}, that $c=-\ell\<\alpha,h\>$
\halmos

\begin{corollary}\label{co:level 1 Z}
If $G$ is simply--laced, the action of $Z(G)$ via conjugation by discontinuous
loops on the irreducible level 1 positive energy representations of $LG$ is
transitive and free.
\end{corollary}
\proof
This follows at once from theorem \ref{th:Z on per} and corollary
\ref{co:level 1 alcove} \halmos

\begin{corollary}\label{co:exchange}
If $\ell$ is odd, the action of each of the two elements of $Z(\Spin_{2n})$
corresponding to $-1\in Z(\SO_{2n})$ on the positive energy representations
of level $\ell$ maps those whose lowest energy subspace is a single--valued
$\SO_{2n}$--module to those whose lowest energy subspace is a two--valued
$\SO_{2n}$--module and viceversa.
\end{corollary}
\proof
This follows from proposition \ref{th:Z on per} and inspection of the tables
in \S \ref{ss:Z on alcove 2} as follows. The highest weights of irreducible
$\Spin_{2n}$--modules are given
by sequences $\mu_{1}\geq\cdots\geq\mu_{n-1}\geq|\mu_{n}|$ where the $\mu_{j}$
are either all integral or half--integral, the latter being the case iff
the corresponding representation is a two--valued $\SO_{2n}$--module.
The elements of $Z(\Spin_{2n})$ mapping to $-1\in\SO_{2n}$ correspond to the
fundamental (co)weights 
$\cow{n-1}$ and $\cow{n}$ whose associated representations are the spin
modules. The result now follows from equations \eqref{eq:n-1}--\eqref{eq:n}
\halmos

\ssubsection{Appendix : explicit action of $Z(G)$ on the level $\ell$ alcove}
\label{ss:Z on alcove 2}

We give the explicit action of $Z(G)$ on $\al$ for each $G$. For simply--laced $G$, the
coroot and coweight lattices are identified with the root and weight lattices respectively.
We denote by $\theta_{i}$, $i=1\ldots n$ and $\<\cdot,\cdot\>$ the standard basis
and inner product in $\IR^{n}$, by $I$ the self-dual lattice $\bigoplus_{i}\theta_{i}\IZ$
and by $\Ie=\{\lambda\in I|\thinspace|\lambda|=\sum_{i}\lambda_{i}\in 2\IZ\}$. Unless
otherwise indicated, the basic inner product is the standard one.\\

{\bf SU$_{\bf n}$, ${\bf n\geq 2}$ (simply-laced)}\\
roots : $\theta_{i}-\theta_{j}$, $i\neq j$.\\
root lattice : $\root=\{\xi \in I|\thinspace \sum_{i}\xi_{i}=0\}$.\\
simple roots : $\alpha_{i}=\theta_{i}-\theta_{i+1}$, $i=1\ldots n-1$.\\
highest root : $\theta=\theta_{1}-\theta_{n}=\alpha_{1}+\cdots+\alpha_{n-1}$.\\
fundamental weights :
$\lambda_{i}=\theta_{1}+\cdots+\theta_{i}-\frac{i}{n}\sum_{j}\theta_{j}$.\\
weight lattice : generated by the $\lambda_{i}$ but more conveniently identified
with $I/\xi\IZ$ where $\xi=\sum\theta_{j}$\\
with inner product
$\<\lambda,\mu\>=\<\lambda-\frac{\<\lambda,\xi\>}{\<\xi,\xi\>}\xi,
                   \mu-\frac{\<\mu,\xi\>}{\<\xi,\xi\>}\xi\>$.\\
centre : $\lambda_{k}=k\lambda_{1}$ mod $\root$ and therefore
$\weight/\root\cong\IZ_{n}$ is generated by $\lambda_{1}$.\\
Weyl group : $\mathfrak S_{n}$ acting by permutation of the $\theta_{i}$.\\
$W_{0}$ : $w_{k}$ is the cyclic permutation
$(\theta_{1}\cdots\theta_{n})^{k}=(\alpha_{0}\cdots\alpha_{n-1})^{k}$.\\
level $\ell$ alcove :
$\al=\{\lambda\in I|\thinspace\lambda_{1}\geq\cdots\geq\lambda_{n},\thickspace
                              \lambda_{1}-\lambda_{n}\leq\ell\}/(\sum_{j}\theta_{j})$.\\
action of the centre :
$A_{k}(\mu_{1},\ldots,\mu_{n})=
 (\ell+\mu_{n+1-k},\ldots,\ell+\mu_{n},\mu_{1},\ldots,\mu_{n-k})$.\\

{\bf Spin$_{\bf 2n+1}$}\\
Since $\Spin_{3}\cong\SU_{2}$, we assume $n\geq 2$.\\
roots : $\pm\theta_{i}\pm\theta_{j}$, $i\neq j$ and $\pm\theta_{i}$.\\
root lattice : $\root=I$.\\
simple roots : $\alpha_{i}=\theta_{i}-\theta_{i+1}$, $i=1\ldots n-1$ and
$\alpha_{n}=\theta_{n}$.\\
highest root : $\theta=\theta_{1}+\theta_{2}=\alpha_{1}+2(\alpha_{2}+\cdots+\alpha_{n})$.\\
coroots : $\pm\theta_{i}\pm\theta_{j}$, $i\neq j$ and $\pm 2\theta_{i}$.\\
coroot lattice : $\coroot=\Ie$.\\
simple coroots : $\cor{i}=\theta_{i}-\theta_{i+1}$, $i=1\ldots n-1$ and
$\cor{n}=2\theta_{n}$.\\
highest coroot : $\wh\theta=2\theta_{1}=2(\cor{1}+\cdots+\cor{n-1})+\cor{n}$.\\
coweight lattice : $\coweight=I^{*}=I$.\\
fundamental coweights : $\cow{i}=\theta_{1}+\cdots+\theta_{i}$, $i=1\ldots n$.\\
centre : $\coweight/\coroot\cong\IZ_{2}$ generated by $\cow{1}$.\\
weight lattice : $\weight={\coroot}^{*}=I+\half{1}(\theta_{1}+\cdots+\theta_{n})\IZ$.\\
fundamental weights : $\lambda_{i}=\theta_{1}+\cdots+\theta_{i}$, $i=1\ldots n-1$ and
$\lambda_{n}=\half{1}(\theta_{1}+\cdots+\theta_{n})$.\\
dual of centre : $\weight/\root\cong\IZ_{2}$ generated by $\lambda_{n}$.\\
Weyl group : $\mathfrak S_{n}\ltimes\IZ_{2}^{n}$ acting by permutations and sign changes of
the $\theta_{i}$.\\
$W_{0}$ : $w_{1}$ is the sign change $\theta_{1}\rightarrow-\theta_{1}$ permuting
$\alpha_{0}$ and $\alpha_{1}$.\\
level $\ell$ alcove :
$\al=\{\mu\in I+\half{1}(\theta_{1}+\ldots+\theta_{n})\IZ|
     \thinspace\mu_{1}\geq\cdots\geq\mu_{n}\geq 0,\medspace \mu_{1}+\mu_{2}\leq\ell\}$.\\
action of the centre :
$A_{1}(\mu_{1},\ldots,\mu_{n})=(\ell-\mu_{1},\mu_{2},\ldots,\mu_{n})$.\\
level 1 weights : $\lambda_{1},\lambda_{n}$.\\

{\bf Sp$_{\bf n}$}\\
Since $\Sp_{1}\cong\SU_{2}$, we assume $n\geq 2$.\\
roots : $\pm\theta_{i}\pm\theta_{j}$, $i\neq j$ and $\pm 2\theta_{i}$.\\
root lattice : $\root=\Ie$.\\
simple roots : $\alpha_{i}=\theta_{i}-\theta_{i+1}$, $i=1\ldots n-1$ and
$\alpha_{n}=2\theta_{n}$.\\
highest root : $\theta=2\theta_{1}=2(\alpha_{1}+\cdots+\alpha_{n-1})+\alpha_{n}$.\\
basic inner product : half the standard one on $\IR^{n}$.\\
coroots : $\pm 2(\theta_{i}\pm\theta_{j})$, $i\neq j$ and $\pm 2\theta_{i}$.\\
coroot lattice : $\coroot=2I$.\\
simple coroots : $\cor{i}=2(\theta_{i}-\theta_{i+1})$, $i=1\ldots n-1$ and
$\cor{n}=2\theta_{n}$.\\
highest coroot : $\wh\theta=2(\theta_{1}+\theta_{2})=\cor{1}+2(\cor{2}+\cdots+\cor{n})$.\\
coweight lattice : $\coweight=2I+(\theta_{1}+\cdots+\theta_{n})\IZ$.\\
fundamental coweights : $\cow{i}=2(\theta_{1}+\cdots+\theta_{i})$, $i=1\ldots n-1$ and
$\cow{n}=\theta_{1}+\ldots+\theta_{n}$.\\
centre : $\coweight/\coroot\cong\IZ_{2}$ generated by $\cow{n}$.\\
weight lattice : $\weight=I$.\\
fundamental weights : $\lambda_{i}=\theta_{1}+\cdots+\theta_{i}$, $i=1\ldots n$.\\
dual of centre : $\weight/\root\cong\IZ_{2}$ generated by $\lambda_{1}$.\\
Weyl group : $\mathfrak S_{n}\ltimes\IZ_{2}^{n}$ acting by permutations and sign changes of
the $\theta_{i}$.\\
$W_{0}$ : $w_{n}$ is the transformation $\theta_{i}\rightarrow -\theta_{n+1-i}$.\\
level $\ell$ alcove :
$\al=\{\mu\in I|\thinspace \ell\geq\mu_{1}\geq\cdots\geq\mu_{n}\geq 0\}$.\\
action of the centre :
$A_{n}(\mu_{1},\ldots,\mu_{n})=(\ell-\mu_{n},\ldots,\ell-\mu_{1})$.\\
level 1 weights : $\lambda_{i}$, $i=1\ldots n$.\\

{\bf Spin$_{\bf 2n}$, ${\bf n\geq 3}$ (simply laced)}\\
roots : $\pm\theta_{i}\pm\theta_{j}$, $i\neq j$.\\
root lattice : $\root=\coroot=\Ie$.\\
simple roots : $\alpha_{i}=\theta_{i}-\theta_{i+1}$, $i=1\ldots n-1$ and
$\alpha_{n}=\theta_{n-1}+\theta_{n}$.\\
highest root : $\theta=
\theta_{1}+\theta_{2}=\alpha_{1}+2(\alpha_{2}+\cdots+\alpha_{n-2})+\alpha_{n-1}+\alpha_{n}$.\\
weight lattice : $\weight=\coweight=I+\half{1}(\theta_{1}+\cdots+\theta_{n})\IZ$.\\
fundamental weights :
\begin{align}
\lambda_{i}  &= \theta_{1}+\cdots+\theta_{i}\qquad\qquad i=1\ldots n-2\\
\lambda_{n-1}&= \half{1}(\theta_{1}+\cdots+\theta_{n-1}-\theta_{n})\\
\lambda_{n}  &= \half{1}(\theta_{1}+\cdots+\theta_{n-1}+\theta_{n})
\end{align}
centre : We have $2\lambda_{1}=2\theta_{1}\in\root$. Moreover, for $n$ even,
$2\lambda_{n-1}=(\theta_{1}+\theta_{2})+\cdots+(\theta_{n-1}\pm\theta_{n})=0$
mod $\root$ and $2\lambda_{n}=0$ mod $\root$. On the other hand, for $n$ odd,
$2\lambda_{n}=\theta_{1}+(\theta_{2}+\theta_{3})+\cdots+(\theta_{n-1}\pm\theta_{n})=
  \lambda_{1}$ mod $\root$ and similarly $2\lambda_{n-1}=\lambda_{1}$ mod $\root$.
Thus, $Z(\Spin_{2n})$ is isomorphic to $\IZ_{2}\times\IZ_{2}$ for $n$ even and to
$\IZ_{4}$ for $n$ odd with $\lambda_{n-1}$ and $\lambda_{n}$ of order 4.\\
Weyl group : $\mathfrak S_{n}\ltimes\IZ_{2}^{n-1}$ acting by permutations and even numbers of
sign changes of the $\theta_{i}$.\\
$W_{0}$ : $w_{1}$ is the sign change $\theta_{1}\rightarrow -\theta_{1}$,
$\theta_{n}\rightarrow -\theta_{n}$ and permutes $\{\alpha_{0},\alpha_{1}\}$
and $\{\alpha_{n-1},\alpha_{n}\}$.
For $n$ even, $w_{n-1}$ is given by $\theta_{i}\rightarrow -\theta_{n+1-i}$,
$2\leq i\leq n-1$ and $\theta_{1}\leftrightarrow\theta_{n}$ and permutes
$\{\alpha_{0},\alpha_{n-1}\}$ and $\{\alpha_{1},\alpha_{n}\}$ while
$w_{n}$ is given by $\theta_{i}\rightarrow -\theta_{n+1-i}$ and permutes
$\{\alpha_{k},\alpha_{n-k}\}$.
For $n$ odd, $w_{n-1}$ is given by $\theta_{1}\rightarrow\theta_{n}$ and
$\theta_{i}\rightarrow -\theta_{n+1-i}$, $i=2\ldots n$ and acts as the cyclic permutation
$\begin{pmatrix}\alpha_{1}&\alpha_{n}&\alpha_{0}&\alpha_{n-1}\end{pmatrix}$ while
$w_{n}$ is given by $\theta_{i}\rightarrow -\theta_{n+1-i}$, $i=1\ldots n-1$ and
$\theta_{n}\rightarrow\theta_{1}$ and acts as
$\begin{pmatrix}\alpha_{1}&\alpha_{n}&\alpha_{0}&\alpha_{n-1}\end{pmatrix}^{-1}$.\\
level $\ell$ alcove :
$\al=\{\mu\in I+\half{1}(\theta_{1}+\ldots+\theta_{n})\IZ|
     \thinspace\mu_{1}\geq\cdots\geq\mu_{n-1}\geq|\mu_{n}|,
     \medspace\mu_{1}+\mu_{2}\leq\ell\}$.\\
action of the centre :
\begin{align}
A_{1}(\mu_{1},\ldots,\mu_{n})&=(\ell-\mu_{1},\mu_{2},\ldots,\mu_{n-1},-\mu_{n}).\\
A_{n-1}(\mu_{1},\ldots,\mu_{n})&=
\begin{cases}
(\half{\ell}+\mu_{n},\half{\ell}-\mu_{n-1},\ldots,\half{\ell}-\mu_{2},-\half{\ell}+\mu_{1})
&\text{$n$ even}\\[1.2 ex]
(\half{\ell}-\mu_{n},\ldots,\half{\ell}-\mu_{2},-\half{\ell}+\mu_{1})
&\text{$n$ odd}
\end{cases}\label{eq:n-1}\\
A_{n}(\mu_{1},\ldots,\mu_{n})&=
\begin{cases}
(\half{\ell}-\mu_{n},\ldots,\half{\ell}-\mu_{1})
&\text{$n$ even}\\[1.2 ex]
(\half{\ell}+\mu_{n},\half{\ell}-\mu_{n-1},\ldots,\half{\ell}-\mu_{1})
&\text{$n$ odd}\label{eq:n}
\end{cases}
\end{align}

\ssection{Definition and classification of primary fields}
\label{se:classification of pf}

For any $G$--module $V$, $\lpol\gc$ acts on the space of $V$--valued Laurent
polyomials $V[z,z^{-1}]$ by multiplication. We extend this to an action of
$\lpol\gc\rtimes\rot$ by setting $R_{\theta}f(z)=f(e^{-i\theta}z)$. In terms
of the basis $v(n)=v\otimes z^{n}$ of $V[z,z^{-1}]$, the action is given by
\begin{xalignat}{2}
X(m)v(n)&=Xv(n+m)&dv(n)&=-nv(n)
\end{xalignat}

Let now $\ell\in\IN$ and consider the set $\pl$ of positive energy
representations of $LG$ at level $\ell$.\\

\definition
Let $\H_{i},\H_{j}\in\pl$ be irreducibles and $V_{i},V_{j}$ the corresponding
lowest energy subspaces. If $V_{k}$ is an irreducible $G$--module, a
{\it primary field} of charge $V_{k}$ is a linear map $\phi:\hfin_{i}\otimes
V_{k}[z,z^{-1}]\rightarrow\hfin_{j}$ intertwining the action of
$\lpol\g\rtimes\rot$, so that
\begin{xalignat}{2}
[X,\phi(f)]&=\phi(Xf)&
[d,\phi(f)]&=i\phi(\dot f)
\end{xalignat}
We represent $\phi$ by the (formal) operator--valued distribution
$\phi(v,z)=\sum_{n\in\IZ}\phi(v,n)z^{-n-\Delta_{\phi}}$ where
$\phi(v,n)=\phi(v\otimes z^{n}):\hfin_{i}\rightarrow\hfin_{j}$. Here,
$\Delta_{\phi}=\Delta_{i}+\Delta_{k}-\Delta_{j}$ is the {\it conformal weight}
of $\phi$, the summands being the Casimirs of the corresponding $G$-modules
divided by $2\kappa$ where $\kappa=\ell+\half{C_{\g}}$ and $C_{\g}$ is the
Casimir of the adjoint representation. The intertwining properties of $\phi$
may then be rephrased as
\begin{xalignat}{2}
[X(m),\phi(v,n)]&=\phi(Xv,m+n)&
[d,\phi(v,n)]&=-n\phi(v,n)\\
\intertext{or as}
[X(m),\phi(v,z)]&=\phi(Xv,z)z^{m}&
[d,\phi(v,z)]&=(z\frac{d}{dz}+\Delta_{\phi})\phi(v,z) \label{intertwine}
\end{xalignat}

By restriction, $\phi(\cdot,0)$ gives rise to a $G$-intertwiner or {\it initial
term} $\varphi:V_{i}\otimes V_{k}\rightarrow V_{j}$ which determines $\phi$
uniquely. In fact, the description of $\hfin_{i}$ as the quotient of the module
freely generated by the action of $\wh\g_{-}=\oplus_{n<0}\gc\otimes z^{n}$ on $V_{i}$
modulo the single relation $e_{\theta}(-1)^{\ell-\<\lambda_{i},\theta\>+1}v_{i}=0$
where $\lambda_{i}$ is the highest weight of $V_{i}$, $v_{i}$ the corresponding
eigenvector and $\theta$ the highest root of $G$ \cite{Ka1}, entails the following

\begin{proposition}[Tsuchiya--Kanie]\label{existence}
There exists a (necessarily unique) primary field
\begin{equation}
\phi:\hfin_{i}\otimes V_{k}[z,z^{-1}]\rightarrow\hfin_{j}
\end{equation}
with given initial term $\varphi\in\Hom_{G}(V_{i}\otimes V_{k},V_{j})$ if, and only if
the restriction of $\varphi$ to any triple $U_{i}\subset V_{i}$, $U_{k}\subset V_{k}$,
$U_{j}\subset V_{j}$ of irreducible $\{e_{\theta},f_{\theta},h_{\theta}\}$--submodules
with highest weights $s_{i},s_{k},s_{j}$ is zero whenever $s_{i}+s_{k}+s_{j}\geq 2\ell+1$.
If $\varphi$ is non--zero and satisfies this criterion, then $V_{k}$ is admissible at
level $\ell$. In particular, the charge of a non--zero primary field is necessarily
$\ell$-admissible.
\end{proposition}
\proof The above is proved in \cite[thm. 2.3]{TK1} for $G=\SU_{2}$ and in \cite{Wa2}
for $G=\SU_{n}$. The extension to any $G$ is immediate \halmos\\

\remark Following Tsuchiya and Kanie, we shall represent the $G$--modules attached
to a primary field $\phi:\hfin_{i}\otimes V_{k}[z,z^{-1}]\rightarrow\hfin_{j}$ by
the {\it vertex} $\vertex{V_{k}}{V_{j}}{V_{i}}$ or, more succintly $\vertex{k}{j}{i}$.\\

\definition For any triple $V_{i},V_{k},V_{j}$ of irreducible representations of $G$,
let $\Homl_{G}(V_{i}\otimes V_{k},V_{j})\subset\Hom_{G}(V_{i}\otimes V_{k},V_{j})$ be
the subspace of intertwiners satisfying the condition of proposition \ref{existence}.
By symmetry of the latter, the isomorphisms
$\Hom_{G}(V_{i}\otimes V_{k},V_{j})\cong\Hom_{G}(V_{k}\otimes V_{i},V_{j})$
and
$\Hom_{G}(V_{i}\otimes V_{k},V_{j})\cong\Hom_{G}(V_{j}^{*}\otimes V_{k},V_{i}^{*})$
restrict to isomorphisms of the corresponding $\Homl_{G}$ spaces and the same is true
for any permutation of $V_{i},V_{j},V_{k}$.

\begin{lemma}\label{lemmino}
Let $V_{i},V_{k},V_{j}$ be irreducible $G$-modules which are admissible at level $\ell$.
If one of them is minimal, then
$\Homl_{G}(V_{i}\otimes V_{k},V_{j})=\Hom_{G}(V_{i}\otimes V_{k},V_{j})$.
\end{lemma}
\proof
Up to a permutation, we may assume that $V_{k}$ is minimal. By corollary \ref{co:one dim},
$\Hom_{G}(V_{i}\otimes V_{k},V_{j})$ is at most one-dimensional. Assume the generator
$\varphi$ restricts to a non--zero intertwiner $U_{i}\otimes U_{k}\rightarrow U_{j}$ for
some triple of $\{e_{\theta},f_{\theta},h_{\theta}\}$--submodules with highest weights
$s_{i},s_{k},s_{j}$.
Since $V_{k}$ is minimal, $s_{k}\in\{0,1\}$ and therefore, by the Clebsch-Gordan rules
$\min(s_{i},s_{j})=\max(s_{i},s_{j})-s_{k}$ so that
$s_{i}+s_{j}+s_{k}=2\max(s_{i},s_{j})\leq 2\ell$ \halmos\\

\definition Let $\phi:\hfin_{i}\otimes V_{k}[z,z^{-1}]\rightarrow\hfin_{j}$ be a
primary field with charge $V_{k}$. The {\it adjoint field}
$\phi^{*}:\hfin_{j}\otimes V_{k}^{*}[z,z^{-1}]\rightarrow\hfin_{i}$ is the unique
primary field with charge $V_{k}^{*}$ satisfying
\begin{equation}\label{eq:adj field}
(\phi(f)\xi,\eta)=(\xi,\phi^{*}(\overline{f})\eta)
\end{equation}
for any $\xi\in\hfin_{i},\eta\in\hfin_{j}$ and $f\in V_{k}[z,z^{-1}]$. It is defined
in the following way. Using the anti--linear identification
$\overline{\cdot}:V_{k}\rightarrow V_{k}^{*}$, set
\begin{equation}
\phi^{*}(\overline{v},n)=\phi(v,-n)^{*}
\end{equation}
then $\phi^{*}$ satisfies \eqref{eq:adj field} and defines the required primary field.









\chapter{Analytic properties of positive energy representations}\label{ch:analytic}

Using the Segal--Sugawara formula, we show in section \ref{se:smooth vectors} that
the infinitesimal action of $\lpol\g$ on the finite energy subspace of a positive
energy representation $\H$ extends to one of $L\g$ on the space $\hsmooth$ of {\it
smooth vectors} for $\rot$. When exponentiated, this action coincides with the
original representation of $LG$. We prove that $\hsmooth$ is invariant under $LG$
and find that the crossed homomorphism for the joint projective action of $LG$ and
$L\g\rtimes i\IR d$ on $\hsmooth$ agrees with that given by Pressley and Segal
\cite[4.9.4]{PS}.\\

In section \ref{se:central ext}, we show that the topological central extensions of
$LG$ arising from positive energy representations are smooth, in fact real--analytic
and compute their Lie algebra cocycle. Their classification follows at once from
that of Pressley and Segal \cite[4.4.1]{PS} and shows in particular that the central
extensions corresponding to positive energy representations of the same level are
canonically isomorphic. Moreover, every positive energy representation possesses a
dense subspace of analytic, and {\it a fortiori} smooth vectors for $LG$, a result
conjectured by Pressley and Segal \cite[\S 9.3]{PS}.

\ssection{The subspace of smooth vectors}
\label{se:smooth vectors}

Let $(\pi,\H)$ be a positive energy representation of $LG$ and $d$ the self--adjoint
generator of rotations normalised by $\left.d\right|_{\H(n)}=n$, $n\in\IN$. As
pointed out by Goodman and Wallach \cite{GoWa}, the Segal--Sugawara formula
(\cite[\S 9.4]{PS} and \S\ref{quadratic} below) implies that $d$ plays a r\^ole
analogous to that of the laplacian of $LG$ in the representation. In particular,
we shall prove in \S\ref{LG on smooth} that the action of $LG$ preserves the
abstract Sobolev scale corresponding to $\rot$ which we presently define.
Let $\|\cdot\|_{s}$, $s\in\IR$ be the Sobolev norm on $\Hfin$ given by
\begin{equation}\label{sobolev norm}
\|\xi\|_{s}^{2}=\|(1+d)^{s}\xi\|^{2}=((1+d)^{2s}\xi,\xi)
\end{equation}
where the powers $(1+d)^{s}$ are defined by the spectral theorem or simply by
expanding $\xi$ as a sum of eigenvectors of $d$. The {\it scale} $\H^{s}$ is
the Hilbert space completion of $\Hfin$ with respect to $\|\cdot\|_{s}$. Thus,
$\H^{s}=\D((1+d)^{s})$ for $s\geq 0$ and $(1+d)$ gives a unitary
$\H^{t}\rightarrow\H^{t-1}$ for any $t$. The space $\hsmooth$ of {\it smooth
vectors} for $\rot$ is, by definition, $\bigcap_{s}\H^{s}$. When endowed with
the corresponding direct limit topology, $\hsmooth$ is a Fr\'echet space which
may equivalently be described as the space of $\xi\in\H$ such that the function
$\theta\rightarrow\pi(R_{\theta})\xi$ is smooth, topologised as a closed subspace
of $C^{\infty}(S^{1},\H)$.

\ssubsection{The Banach Lie algebras $L\g_{t}$}

We shall need an alternative description of the Fr\'echet topology on $L\g$.
The restriction of the basic inner product $\<\cdot,\cdot\>$ to $i\g\subset\gc$
is positive definite and we extend it to a $G$-invariant hermitian form on $\gc$
with associated norm $\|\cdot\|$. For $X=\sum_{k}a_{k}e^{ik\theta}\in\lpol\gc$
and $t\geq 0$, define
\begin{equation}\label{Lg sobolev}
|X|_{t}=\sum_{k}(1+|k|)^{t}\|a_{k}\|
\end{equation}
Let $L\g_{t}$ be the completion of $\lpol\g$ with respect
to $|\cdot|_{t}$. If $a_{k}$, $b_{l}$ are the Fourier coefficents of
$X,Y\in\lpol\g$ and $C_{\g}>0$ is such that
$\|[x,y]\|\leq C_{\g}\|x\|\|y\|$ for any $x,y\in\g$, then
\begin{equation}
|[X,Y]|_{t}
\leq C_{\g}\sum_{k,l}\|a_{k}\|\|b_{l}\|(1+|k+l|)^{t}
\leq C_{\g}|X|_{t}|Y|_{t}
\end{equation}
since $(1+|k+l|)^{t}\leq(1+|k|)^{t}(1+|l|)^{t}$ for $t\geq 0$. $L\g_{t}$
is therefore a Banach Lie algebra. If $\|.\|_{\infty}$ is the supremum
norm on $C(S^{1},\g)$, and $X=\sum_{k}a_{k}e^{ik\theta}\in\lpol\g$, then
for any $t\geq n$
\begin{equation}
\|X^{(n)}\|_{\infty}\leq\sum_{k}\|a_{k}\||k|^{n}\leq|X|_{t}
\end{equation}

Moreover, if $2\geq s>1$,
\begin{equation}
\begin{split}
|X|_{t}
&=\sum_{k}(1+|k|)^{-\half{s}}(1+|k|)^{t+\half{s}}\|a_{k}\|\\
&\leq\biggl\{\sum_{k}(1+|k|)^{-s}\biggr\}^{\half{1}}
     \biggl\{\sum_{k}(1+|k|)^{2t+s}\|a_{k}\|^{2}\biggr\}^{\half{1}}\\
&=C_{s}\|(1+|\frac{d}{d\theta}|)^{t+\half{s}}X\|_{L^{2}(S^{1},\g)}\\[1.2ex]
&\leq C'_{s}\|X\|_{C^{\lceil t\rceil+1}(S^{1},\g)}
\end{split}
\end{equation}

Consequently, we have the following norm-continuous embeddings with dense image
\begin{equation}
C^{\lceil t\rceil+1}(S^{1},\g)\hookrightarrow 
L\g_{t}\hookrightarrow C^{\lfloor t\rfloor}(S^{1},\g)
\end{equation}
showing that the Fr\'echet topology on $L\g$ given by the $C^{k}(S^{1},\g)$ norms
may equivalently be given by the norms $|\cdot|_{t}$.

\ssubsection{The Segal--Sugawara formula and Sobolev estimates on the action of $\lpol\g$}
\label{quadratic}

Let $\H$ be a positive energy representation of $LG$ at level $\ell$. The
restriction of the infinitesimal generator of rotations to $\Hfin$ is given
via the Segal--Sugawara formula \cite[\S 9.4]{PS}
\begin{equation}\label{sugawara}
L_{0}=
\frac{1}{\kappa}\Bigl(\half{1}X_{i}(0)X^{i}(0)+\sum_{m>0} X_{i}(-m)X^{i}(m)\Bigr)
\end{equation}
where $\{X_{i}\},\{X^{i}\}$ are dual basis of $\gc$ with respect to the basic
inner product and $\kappa=\ell+g$ with $g$ the dual Coxeter number of $\g$ (half
the Casimir of the adjoint representation). More precisely, $L_{0}$ satisfies
$[L_{0},x(n)]=-nx(n)$ and therefore differs from $d$ by an additive constant on
each irreducible summand of $\Hfin$. If $\H$, and therefore $\H(0)$, is irreducible,
$L_{0}$ acts on $\H(0)$ as the scalar $\Delta=\frac{C}{2\kappa}$, where $C$ is
the Casimir of $\H(0)$ so that
\begin{equation}\label{energy shift}
L_{0}=d+\Delta
\end{equation}
More generally, only finitely many irreducible $G$-modules appear as lowest energy
subspaces of positive energy representations at level $\ell$. It follows that for
any level $\ell$ representation $\H$, $d\leq L_{0}\leq d+C_{\ell}$ for some constant
$C_{\ell}$ depending only on $\ell$ and therefore that $d$ and $L_{0}$ define the same
scales. For convenience, we shall often supersede definition \eqref{sobolev norm} and
refer to $\|(1+L_{0})^{s}\xi\|$ as $\|\xi\|_{s}$. The following estimates are due to
Goodman and Wallach \cite[\S 3.2]{GoWa}

\begin{proposition}\label{sobolev estimates}
Let $X\in\lpol\gc$ and $\xi\in\Hfin$. Then, for any $s\in\IR$
\begin{equation}
\|\pi(X)\xi\|_{s}\leq\sqrt{2\kappa}|X|_{|s|+\half{1}}\|\xi\|_{s+\half{1}}
\end{equation}
\end{proposition}
\proof
It is clearly sufficient to show that for any $x\in\g$ of norm 1 and $p\in\IZ$
\begin{equation}\label{easier}
\|x(p)\xi\|_{s}^{2}\leq 2\kappa(1+|p|)^{2|s|+1}\|\xi\|^{2}_{s+\half{1}}
\end{equation}
Let $X_{i}$ be an orthonormal basis of $\g$ so that $X_{i}$ and $-X_{i}$ are dual basis
with respect to $\<\cdot,\cdot\>$. By formal adjunction, $X_{i}(n)^{*}=-X_{i}(-n)$ and
therefore, by \eqref{sugawara}
\begin{equation}
\begin{split}
(L_{0}\xi,\xi)
&=\frac{1}{2\kappa}
\sum_{i}\Bigl(-(X_{i}(0)X_{i}(0)\xi,\xi)-2\sum_{p>0}(X_{i}(-p)X_{i}(p)\xi,\xi)\Bigr)\\
&=\frac{1}{2\kappa}
\sum_{i}\Bigl(\|X_{i}(0)\xi\|^{2}+2\sum_{p>0}\|X_{i}(p)\xi\|^{2}\Bigr)
\end{split}
\end{equation}
whence, for any $p\geq 0$,
$\|X_{i}(p)\xi\|^{2}\leq 2\kappa\|L_{0}^{\half{1}}\xi\|^{2}$.
If $p<0$, $[X_{i}(p),X_{i}(-p)]=-p$ and therefore
$\|X_{i}(p)\xi\|^{2}=\|X_{i}(-p)\xi\|^{2}-p\|\xi\|^{2}$. Thus, for any $p$
\begin{equation}
\|X_{i}(p)\xi\|^{2}\leq 2\kappa\|L_{0}^{\half{1}}\xi\|^{2}+|p|\|\xi\|^{2}
\end{equation}
Since the eigenspaces of $L_{0}$ are orthonormal for any of the Sobolev norms, is it 
sufficient, to show \eqref{easier}, to restrict to the case $L_{0}\xi=m\xi$.
Therefore, assuming $p\leq m$ (else the left hand side vanishes),
\begin{equation}
\begin{split}
\|X_{i}(p)\xi\|_{s}^{2}
&=\|(1+L_{0})^{s}X_{i}(p)\xi\|^{2}\\[1.2 ex]
&=(1+m-p)^{2s}\|X_{i}(p)\xi\|^{2} \\[1.2 ex]
&\leq(1+m-p)^{2s}(2\kappa m+|p|)\|\xi\|^{2} \\[1.2 ex]
&\leq 2\kappa\frac{(1+m-p)^{2s}}{(1+m)^{2s}}(1+|p|)\|\xi\|^{2}_{s+\half{1}}             
\end{split}
\end{equation}
If $s\geq 0$,
$\frac{(1+m-p)^{2s}}{(1+m)^{2s}}\leq g(m)=\frac{(1+m+|p|)^{2s}}{(1+m)^{2s}}$
and the bound is decreasing in $m$ so that $g(m)\leq g(0)=(1+|p|)^{2s}$.
If, on the contrary, $s<0$ and $p\leq 0$ we have
$\frac{(1+m-p)^{2s}}{(1+m)^{2s}}=(1-\frac{p}{m+1})^{2s}\leq 1$.
Lastly, if $0<p\geq m$, $h(m)=\frac{(1+m-p)^{2s}}{(1+m)^{2s}}$ is decreasing
in $m$ and therefore $h(m)\leq h(p)=(1+p)^{2|s|}$. Either way, we obtain
\begin{equation}
\|X_{i}(p)\xi\|_{s}^{2}
\leq 2\kappa(1+|p|)^{2|s|+1}\|\xi\|_{s+\half{1}}^{2}
\end{equation}
\halmos

\ssubsection{The action of $L\g$ on $\hsmooth$}

Let $(\pi,\H)$ be a level $\ell$ positive energy representation of $LG$ with
finite energy subspace $\Hfin$. Then,

\begin{corollary}
The action of $\lpol\gc$ on $\Hfin$ given by theorem \ref{ch:classification}.\ref{core}
extends to a jointly continuous, projective action
$L\gc\otimes\hsmooth\rightarrow\hsmooth$ satisfying
\begin{align}
[d,\pi(X)]&=i\pi(\dot X) \label{commutation 1}\\
[\pi(X),\pi(Y)]&=\pi([X,Y])+i\ell B(X,Y)\label{commutation 2}
\end{align}
where, as customary $B(X,Y)=\int_{0}^{2\pi}\<X,\dot Y\>\frac{d\theta}{2\pi}$. Moreover,
the operators $\pi(X)$, $X\in L\g$ are essentially skew--adjoint on $\hsmooth$.
\end{corollary}
\proof The first claim is an immediate consequence of the estimates of proposition
\ref{sobolev estimates}. The second follows from the estimates
\begin{align}
\|\pi(X)\xi\|_{s}&\leq\sqrt{2\kappa}|X|_{|s|+\half{1}}\|\xi\|_{s+\half{1}}\\
\|[1+d,\pi(X)]\xi\|_{s}
&\leq\sqrt{2\kappa}|\dot X|_{|s|+\half{1}}\|\xi\|_{s+\half{1}}
 \leq\sqrt{2\kappa}|X|_{|s|+\half{3}}\|\xi\|_{s+\half{1}}
\end{align}
and Nelson's commutator theorem \cite{Ne2} \halmos\\

Clearly, $\pi$ extends to a jointly continuous map $L\gc\otimes\H^{s}\rightarrow\H^{s-\half{1}}$.
It follows that the restriction of the operators $\pi(X)$, $X\in L\g$ is essentially
skew--adjoint on any dense subspace of $\hsmooth$ for the Fr\'echet topology by Nelson's
commutator theorem \cite{Ne2}. We shall usually consider the operators $\pi(X)$ as being
defined on the invariant core $\hsmooth$. When no confusion arises, we denote by
the same symbol their restriction to any of the $\H^{s}$ and their skew--adjoint closure.\\

The following shows that the exponentiation of the representation of $L\g$ on $\hsmooth$
coincides with the original representation of $LG$.

\begin{proposition}\label{lift one}
For any $X\in L\g$, the following holds in $PU(\H)$
\begin{equation}\label{lift}
\pi(\exp_{LG}X)=e^{\pi(X)}
\end{equation}
\end{proposition}
\proof
When $X\in\lpol\g$, \eqref{lift} holds by the definition of the operators $\pi(X)$.
Assume now $X\in L\g$ and let $X_{n}\in\lpol\g$, $X_{n}\rightarrow X$ in $L\g$, and
{\it a fortiori} in the $|\cdot|_{\half{1}}$ norm. By proposition \ref{sobolev estimates},
\begin{equation}
\|(\pi(X)-\pi(X_{n}))\xi\|\leq
\sqrt{2\kappa}|X-X_{n}|_{\half{1}}\|(1+L_{0})^{\half{1}}\xi\|
\end{equation}
for any $\xi\in\hsmooth$. Since $\hsmooth$ is a core for $\pi(L\g)$, we deduce that
$\pi(X_{n})\rightarrow\pi(X)$ in the strong-resolvent sense \cite[thm. VIII.25 (a)]{RS}
and therefore, by Trotter's theorem $e^{\pi(X_{n})}\rightarrow e^{\pi(X)}$ strongly in
$U(\H)$ \cite[thm. VIII.21]{RS}. Thus,
\begin{equation}
e^{\pi(X)}=
\lim_{n\rightarrow\infty}e^{\pi(X_{n})}=
\lim_{n\rightarrow\infty}\pi(\exp_{LG}(X_{n}))=
\pi(\exp_{LG}(X))
\end{equation}
\halmos

\ssubsection{The exponential map of $LG\rtimes\rot$}

We study below the well--posedness of the exponential map of $LG\rtimes\rot$ and
obtain an explicit formula for $\exp_{LG\rtimes\rot}(i\theta d+X)$, $X\in L\g$.
This will be needed to extend proposition \ref{lift one} to $X\in L\g\rtimes id\IR$.\\

The integration of one--parameter groups in $LG\rtimes\rot$ reduces, after some
simple manipulations to that of a linear, time--dependent ordinary differential
equation for a function with values in $LG$. The latter is easily solved using
Volterra product integrals whose salient features we now recall. Details may be 
found in \cite[chap. 2]{Ne1}. Let $E$ be a Banach space and $\B(E)$ the Banach
algebra of all bounded linear maps $E\rightarrow E$ endowed with the operator
norm. For $A\in\B(E)$, define the operator $e^{A}$ by the convergent power
series $\sum_{n}\frac{A^{n}}{n!}$. Clearly,
$\|e^{A}\|\leq\sum_{n}\frac{\|A\|^{n}}{n!}\leq e^{\|A\|}$. Moreover, from
$A^{n}-B^{n}=\sum_{k=0}^{n-1}A^{k}(A-B)B^{n-1-k}$ we deduce that
\begin{equation}\label{exp inequality}
\|e^{A}-e^{B}\|\leq
\sum_{n\geq 0}\frac{\|A^{n}-B^{n}\|}{n!}\leq
\|A-B\|\sum_{n\geq 1}\frac{\max(\|A\|,\|B\|)^{n-1}}{(n-1)!}=
\|A-B\|e^{\max(\|A\|,\|B\|)}
\end{equation}
For any $A\in C(\IR,\B(E))$ and $a<b\in\IR$,
we define the product integral $\prod_{b\geq\tau\geq a}\Exp(A(\tau)d\tau)$ as
follows. Consider first step functions $A:[a,b]\rightarrow\B(E)$,
$A(\tau)=A_{j}$ for $\tau_{j}>\tau>\tau_{j-1}$ corresponding to subdivisions
$b=\tau_{n}>\tau_{n-1}>\cdots>\tau_{1}>\tau_{0}=a$ and set
\begin{equation}
\prod_{b\geq\tau\geq a}\Exp(A(\tau)d\tau):=
e^{\Delta_{n}A_{n}}\cdots e^{\Delta_{1}A_{1}}
\end{equation}
where $\Delta_{j}=\tau_{j}-\tau_{j-1}$. If $A$, $B$ are two step functions, which
we may take as defined on a common subdivision, the identity
\begin{equation}
 e^{\Delta_{n}A_{n}}\cdots e^{\Delta_{1}A_{1}}-
 e^{\Delta_{n}B_{n}}\cdots e^{\Delta_{1}B_{1}}=\sum_{k=1}^{n}
 e^{\Delta_{n}A_{n}}\cdots e^{\Delta_{k+1}A_{k+1}}
 (e^{\Delta_{k}A_{k}}-e^{\Delta_{k}B_{k}})
 e^{\Delta_{k-1}B_{k-1}}\cdots e^{\Delta_{1}B_{1}}
\end{equation}
and \eqref{exp inequality} imply that
\begin{equation}
\|\prod_{b\geq\tau\geq a}\Exp(A(\tau)d\tau)-
  \prod_{b\geq\tau\geq a}\Exp(B(\tau)d\tau)\|\leq
e^{(b-a)\max(\|A\|_{[a,b]},\|B\|_{[a,b]})}(b-a)\|A-B\|_{[a,b]}
\end{equation}
where the norms on the right hand-side refer to the supremum norm on functions
$[a,b]\rightarrow\B(E)$. We may therefore define the product integral for any
$A\in C([a,b],\B(E))$ by using a sequence of approximating step functions.
Notice that product integrals are invertible operators, in fact
$\prod_{b\geq\tau\geq a}\Exp(A(\tau)d\tau)^{-1}=
 \prod_{b\geq\tau\geq a}\Exp(-\check A(\tau)d\tau)$
where $\check A(\tau)=A(a+b-\tau)$ so that we may define, for $a>b$,
$\prod_{b\geq\tau\geq a}\Exp(A(\tau)d\tau):=
 \prod_{a\geq\tau\geq b}\Exp(A(\tau)d\tau)^{-1}$.

\begin{theorem}\label{product exp}
Let $E$ be a Banach space and $A:\IR\rightarrow\B(E)$ a norm-continuous map. For any
$\xi_{0}\in E$, the time-dependent, linear ordinary differential equation in
$\xi:\IR\rightarrow E$
\begin{align}
\dot\xi(t)&=A(t)\xi(t)\\
\xi(0)&=\xi_{0}
\end{align}
possesses a unique solution $\xi\in C^{1}(\IR,E)$ given by
\begin{equation}
\xi(t)=\prod_{t\geq\tau\geq0}\Exp(A(\tau)d\tau)\xi_{0}
\end{equation}
Moreover, if $A,\widetilde A\in C(\IR,\B(E))$ and $\xi,\widetilde\xi$ are
the corresponding solutions with initial condition $\xi_{0}$, then, for any $t\in\IR$,
\begin{equation}\label{exp growth}
\|\xi(t)-\widetilde\xi(t)\|\leq 
e^{|t|\max(\|A\|_{[0,t]},\|\widetilde A\|_{[0,t]})}|t|
\|A-\widetilde A\|_{[0,t]}\|\xi_{0}\|
\end{equation}
\end{theorem}
\proof See \cite[Thm. 1, page 17]{Ne1} \halmos

\begin{corollary}[\cite{Wa2}]\label{semidir exp}
The exponential map $L\g\rtimes i\IR d\rightarrow LG\rtimes\rot$ is well--defined and
continuous. It is given by 
\begin{equation}\label{explicit exp}
\begin{split}
\exp_{LG\rtimes i\IR}(X+i\alpha d)
&=\prod_{1\geq t\geq 0}\Exp(X_{\alpha(1-t)}dt)R_{\alpha}\\
&=
\lim_{n\rightarrow\infty}
\exp_{LG}\Bigl(\frac{X}{n}\Bigr)
\exp_{LG}\Bigl(\frac{X_{\alpha\cdot 1/n}}{n}\Bigr)
\cdots
\exp_{LG}\Bigl(\frac{X_{\alpha\cdot(n-1)/n}}{n}\Bigr)
R_{\alpha}
\end{split}
\end{equation}
\end{corollary}
\proof
To compute the exponential map, fix $X+i\alpha d\in L\g\rtimes i\IR d$ and
consider $f=\gamma R_{\phi}:\IR\rightarrow LG\rtimes\rot$ satisfying
$\dot f=(X+i\alpha d)f$ and $f(0)=1$. As a manifold, $LG\rtimes\rot$ is the
product of the two factors and therefore $s\rightarrow\exp_{LG}(sX)R_{s\alpha}$ is
an integral curve for $X+i\alpha d$ at 1. Thus
\begin{equation}
(X+i\alpha d)f=
\left.\frac{d}{ds}\right|_{s=0}
\exp_{LG}(sX)R_{s\alpha}\gamma R_{\phi}=
\left.\frac{d}{ds}\right|_{s=0}
\exp_{LG}(sX)\gamma_{s\alpha}R_{s\alpha+\phi}=
X\gamma-\alpha\partial_{\theta}\gamma+i\alpha dR_{\phi}
\end{equation}
whence $\phi(t)=\alpha t$
and we must solve $\partial_{t}\gamma=X\gamma-\alpha\partial_{\theta}\gamma$ with
boundary condition $\gamma(\cdot,0)\equiv 1$. The corresponding homogeneous equation
$\partial_{t}\gamma=-\alpha\partial_{\theta}\gamma$ is easily solved by setting
$\gamma(\theta,t)=\gamma_{0}(\theta-\alpha t)$ for some $\gamma_{0}\in LG$. Varying the
constants, we set $\gamma(\theta,t)=\gamma_{0}(\theta-\alpha t,t)$ and the original
equation yields, in terms of $\gamma_{0}$, 
$\partial_{t}\gamma_{0}(\theta,t)=X(\theta+\alpha t)\gamma_{0}(\theta,t)$ {\it i.e.}~in the
notation of theorem \ref{product exp}, $\dot\gamma_{0}(t)=A(t)\gamma_{0}(t)$, 
$\gamma_{0}(0)=1$, where $A:\IR\rightarrow L\g$ is given by $A(t)=X_{-\alpha t}$.
Thus, if we embed $G$ in a space of matrices $M_{m}(\IC)$ and $LG$ as a closed
subspace of the Fr\'echet space
$C^{\infty}(S^{1},M_{m}(\IC))=\bigcap_{k}C^{k}(S^{1},M_{m}(\IC))$, then
\begin{equation}
\gamma_{0}(1)=
\prod_{1\geq\tau\geq 0}\Exp(X_{-\alpha\tau}d\tau)=
\lim_{n\rightarrow\infty}
\exp\Bigl(\frac{X_{-\alpha}}{n}\Bigr)\cdots
\exp\Bigl(\frac{X_{-2\alpha/n}}{n}\Bigr)
\exp\Bigl(\frac{X_{-\alpha/n}}{n}\Bigr)
\end{equation}
where the right-side converges in each Banach space $C^{k}(S^{1},M_{m}(\IC))$ and hence in
$LG$. \eqref{explicit exp} now follows from $\gamma(\theta,1)=\gamma_{0}(\theta-\alpha,1)$.
The continuity of the exponential map is a direct corollary of \eqref{exp growth} \halmos

\begin{corollary}\label{lie two}
For any $X\in L\g$ and $\theta\in\IR$, we have
\begin{equation}
\exp_{LG\rtimes\rot}(X+i\theta d)=\lim_{n\rightarrow\infty}
\biggl(\exp_{LG}\Bigl(\frac{X}{n}\Bigr)
       \exp_{\rot}\Bigl(\frac{i\theta d}{n}\Bigr)\biggr)^{n}
\end{equation}
\end{corollary}
\proof
We have
\begin{equation}
\biggl(\exp_{LG}\Bigl(\frac{X}{n}\Bigr)
       \exp_{\rot}\Bigl(\frac{i\theta d}{n}\Bigr)\biggr)^{n}=
\exp_{LG}\Bigl(\frac{X}{n}\Bigr)
\exp_{LG}\Bigl(\frac{X_{\theta\cdot 1/n}}{n}\Bigr)
\cdots
\exp_{LG}\Bigl(\frac{X_{\theta\cdot(n-1)/n}}{n}\Bigr)
R_{\theta}
\end{equation}
which, by corollary \ref{semidir exp}, tends to $\exp_{LG\rtimes\rot}(X+i\theta d)$
as $n\rightarrow\infty$ \halmos\\

\remark It is fairly easy to prove that the integration of one--parameter groups in
$LG\rtimes\rot$ is smoothly well--posed and therefore that the exponential map of this
group is smooth. The following result will not be used elsewhere.

\begin{proposition}
The exponential map of $L\IC^{*}\rtimes\rot$ is not locally surjective.
\end{proposition}
\proof
Let $u+i\alpha d\in L\IC\rtimes i\IR=\Lie(L\IC^{*}\rtimes\rot)$. Since $L\IC^{*}$ is
abelian, \eqref{explicit exp} yields
\begin{equation}\label{explicit}
\exp_{L\IC^{*}\rtimes\rot}(u+i\alpha d)=
e^{\int_{0}^{1}u_{\alpha\tau}d\tau}R_{\alpha}
\end{equation}
Taking logarithms on each of the factors of $L\IC^{*}\rtimes\rot$, we see that 
\eqref{explicit} is locally surjective iff the following map is
\begin{equation}
\E:L\IC\rtimes i\IR\rightarrow L\IC\rtimes i\IR,\thickspace
u+i\alpha d\rightarrow \int_{0}^{1}u(\cdot-\alpha\tau)d\tau+i\alpha d
\end{equation}
Write $u=\sum_{k}a_{k}e^{ik\theta}$, then, for $\alpha\neq 0$,
\begin{equation}
\E(u+i\alpha d)-i\alpha d=
a_{0}+\sum_{k\neq 0}a_{k}e^{ik\theta}\frac{1-e^{-ik\alpha}}{ik\alpha}
\end{equation}
and it follows that, for any $b\in\IC^{*}$ and $k\in\IZ^{*}$, the elements
$be^{ik\theta}+\frac{2\pi}{k}i d$, which may be taken arbitrarily close to
the origin, are not in the image of $\E$ \halmos

\ssubsection{The action of $LG$ on $\hsmooth$}\label{LG on smooth}

The invariance of $\hsmooth$ under the action of $LG$ depends upon the following
considerations. Let $\gamma\in LG$ and consider the smooth one--parameter group
in $LG\rtimes\rot$ given by
$\theta\rightarrow\gamma R_{\theta}\gamma^{-1}=\gamma\gamma_{\theta}^{-1}R_{\theta}$.
Since $0=\dot{(\gamma\gamma^{-1})}=\dot\gamma\gamma^{-1}+\gamma\dot\gamma^{-1}$, its
derivative at 0 is $id+\dot\gamma\gamma^{-1}$ and therefore, by uniqueness
\begin{equation}\label{useful}
\gamma R_{\theta}\gamma^{-1}=\exp_{LG\rtimes\rot}(\theta(id+\dot\gamma\gamma^{-1}))
\end{equation}
As shown below, this implies that
$\pi(\gamma)e^{id\theta}\pi(\gamma)^{*}=e^{\theta(id+\pi(\dot\gamma\gamma^{-1}))}$
in $PU(\H)$ and therefore, up to an additive term,
$\pi(\gamma)d\pi(\gamma)^{*}=d-i\pi(\dot\gamma\gamma^{-1})$.
Thus, the perturbed operator $\pi(\gamma)d\pi(\gamma)^{*}$ is equal to $d$ up to
terms which, by proposition \ref{sobolev estimates} are of lower order with respect
to the scale $\H^{s}$ and it follows from this that $\hsmooth$ is invariant under
$LG$. A detailed argument follows.\\

Extending our former notation, we denote by $\pi(X)\in\End(\hsmooth)$,
$X=Y+i\theta d\in L\g\rtimes i\IR$, the operator $\pi(Y)+i\theta d$ so that we get a
projective representation of $L\g\rtimes i\IR$ on $\hsmooth$ satisfying
\begin{equation}
[\pi(X_{1}),\pi(X_{2})]=\pi([X_{1},X_{2}])+i\ell B(X_{1},X_{2})
\end{equation}
where
\begin{equation}
B(Y_{1}+i\phi_{1} d,Y_{2}+i\phi_{2} d)=B(Y_{1},Y_{2})=
\int_{0}^{2\pi}\frac{d\theta}{2\pi}\<Y_{1},\dot Y_{2}\>
\end{equation}
Extending the Sobolev norms $|\cdot|_{t}$ to $L\g\rtimes i\IR$ by
$|Y+i\theta d|_{t}=|Y|_{t}+|\theta|$ and using proposition \ref{sobolev estimates},
we get
\begin{equation}
\|\pi(X)\xi\|_{s}=
\|(1+L_{0})^{s}\pi(X)\xi\|\leq
C|Y|_{|s|+\half{1}}\|\xi\|_{s+\half{1}}+|\theta|\|\xi\|_{s+1}\leq
\max(1,C)|X|_{|s|+\half{1}}\|\xi\|_{s+1}
\end{equation}
so that the $\pi(X)$ extend to bounded operators $\H^{s}\rightarrow\H^{s-1}$
and are essentially skew--adjoint on $\hsmooth$ by Nelson's commutator theorem \cite{Ne2}.
The following extends proposition \ref{lift one} to the Lie algebra of $LG\rtimes\rot$

\begin{lemma}\label{lifting}
For any $X\in L\g\rtimes i\IR d$, the following holds in $PU(\H)$
\begin{equation}
\pi(\exp_{LG\rtimes\rot}X)=e^{\pi(X)}
\end{equation}
\end{lemma}
\proof Let $X=Y+i\theta d\in L\g\rtimes i\IR$. Then, using Trotter's formula \cite{Ch},
proposition \ref{lift one} and corollary \ref{lie two}
\begin{equation}
e^{\pi(Y)+i\theta  d}=
\lim_{n\rightarrow\infty}
(e^{\frac{\pi(Y)}{n}}e^{i\frac{\theta}{n} d})^{n}=
\lim_{n\rightarrow\infty}
\pi(\exp_{LG}\Bigl(\frac{Y}{n}\Bigr)
    \exp_{\rot}\Bigl(\frac{i\theta d}{n}\Bigr))^{n}=
\pi(\exp_{LG\rtimes\rot}(Y+i\theta d))
\end{equation}
\halmos

\begin{lemma}\label{abstract commutation}
The following operator identities hold for any $\gamma\in LG$ and $X\in L\g$
\begin{align}
\pi(\gamma)\pi(X)\pi(\gamma)^{*}&=
\pi(\Ad(\gamma)X)+ic_{1}(\gamma,X)\\
\pi(\gamma) d\pi(\gamma)^{*}&=
d-i\pi(\dot\gamma\gamma^{-1})+c_{2}(\gamma)
\end{align}
for some real-valued functions $c_{1},c_{2}$.
\end{lemma}
\proof
Let $\gamma\in LG$ and $X\in L\g\rtimes i\IR$. Then, by lemma \ref{lifting}, the
following holds in $PU(\H)$
\begin{equation}
\begin{split}
\pi(\gamma)e^{t\pi(X)}\pi(\gamma)^{*}
&=\pi(\gamma)\pi(\exp_{LG\rtimes\rot}(tX))\pi(\gamma^{-1})\\
&=\pi(\gamma\exp_{LG\rtimes\rot}(tX)\gamma^{-1})\\
&=\pi(\exp_{LG\rtimes\rot}(t\Ad(\gamma)X))\\
&=e^{t\pi(\Ad(\gamma)X)}
\end{split}
\end{equation}
and consequently
$\pi(\gamma)e^{t\pi(X)}\pi(\gamma)^{*}=
 \lambda(t)e^{t\pi(\Ad(\gamma)X)}$
where $\lambda$ is a continuous homomorphism $\IR\rightarrow\T$ so that
$\lambda(t)=e^{itc(\gamma,X)}$. The claimed formulae now follow from Stone's
theorem \halmos\\

The following result is due to A. Wassermann \cite{Wa2} and parallels the corresponding
one for positive energy representations of $\diff$ \cite[thm. 7.4]{GoWa}

\begin{proposition}\label{invariance}
The subspaces $\H^{n}$ are invariant under the action of $LG$ and the corresponding
map $LG\times\H^{n}\rightarrow\H^{n}/\T$ is jointly continuous.
\end{proposition}
\proof We proceed in several steps.

{\it Invariance of $\H^{n}$}. From lemma \ref{abstract commutation} and the fact
that $\D(UAU^{*})=U\D(A)$ for any operator $A$ and unitary $U$, we have
\begin{equation}
\begin{split}
\H^{n}
&=\D((1+d)^{n})\\
&\subset\D(\left.(
1+d-i\pi(\dot\gamma\gamma^{-1})+c_{2}(\gamma))\right|_{\D(1+d)}^{n})\\
&=\pi(\gamma)\D((1+d)^{n})\\
&=\pi(\gamma)\H^{n}
\end{split}
\end{equation}
and therefore $\pi(\gamma^{-1})\H^{n}\subset\H^{n}$ for any
$\gamma\in LG$.

{\it Continuity of the cocycle $c_{2}$}. Let $\gamma_{m}\rightarrow\gamma\in LG$ and
choose lifts $\pi(\gamma_{m})$ 
and $\pi(\gamma)$ such that $\pi(\gamma_{m})\rightarrow\pi(\gamma)$ in
$U(\H)$. Let $\xi\in\H^{\infty}$, then from lemma
\ref{abstract commutation},
\begin{equation}
\pi(\gamma_{m})(1+d)\xi=
(1+d-i\pi(\dot\gamma_{m}\gamma_{m}^{-1}))\pi(\gamma_{m})\xi+
c_{2}(\gamma_{m})\pi(\gamma_{m})\xi
\end{equation}
To deduce from this that $c_{2}(\gamma_{m})\rightarrow c_{2}(\gamma)$, we
regularise by multiplying both sides by the compact operator $(1+d)^{-1}$
\begin{equation}
(1+d)^{-1}\pi(\gamma_{m})(1+d)\xi=
\pi(\gamma_{m})\xi-i(1+d)^{-1}\pi(\dot\gamma_{m}\gamma_{m}^{-1})\pi(\gamma_{m})\xi+
c_{2}(\gamma_{m})(1+d)^{-1}\pi(\gamma_{m})\xi
\end{equation}
The first terms on the right and left hand-sides manifestly tend to the corresponding
terms for $\gamma$. To see that this is the case for the second term on the right
hand-side, we use proposition \ref{sobolev estimates}
\begin{equation}
\begin{split}
    &\|(1+d)^{-1}\pi(\dot\gamma_{m}\gamma_{m}^{-1})\pi(\gamma_{m})\xi-
       (1+d)^{-1}\pi(\dot\gamma\gamma^{-1})\pi(\gamma)\xi\|\\
\leq&
\|(1+d)^{-1}\pi(\dot\gamma_{m}\gamma_{m}^{-1})(\pi(\gamma_{m})-\pi(\gamma))\xi\|\\
   +&
\|(1+d)^{-1}(\pi(\dot\gamma_{m}\gamma_{m}^{-1})-\pi(\dot\gamma\gamma^{-1}))
                 \pi(\gamma))\xi\|\\
\leq&
C|\dot\gamma_{m}\gamma_{m}^{-1}|_{\half{3}}
\|(1+d)^{-\half{1}}(\pi(\gamma_{m})-\pi(\gamma))\xi\|\\
   +&
C|\dot\gamma_{m}\gamma_{m}^{-1}-\dot\gamma\gamma^{-1}|_{\half{3}}
\|(1+d)^{-\half{1}}\pi(\gamma))\xi\|
\end{split}
\end{equation}
which tends to zero. Thus,
$c_{2}(\gamma_{m})(1+d)^{-1}\pi(\gamma_{m})\xi\rightarrow
 c_{2}(\gamma)(1+d)^{-1}\pi(\gamma)\xi$ and therefore
$c_{2}(\gamma_{m})\rightarrow c_{2}(\gamma)$ since
$(1+d)^{-1}\pi(\gamma_{m})\xi\rightarrow(1+d)^{-1}\pi(\gamma)\xi$.

{\it Norm estimates on the $\pi(\gamma)$.}
Set $A(\gamma)=-i\pi(\gamma^{-1}\dot\gamma)+c_{2}(\gamma^{-1})$ so
that, on $\hsmooth$, $\pi(\gamma)^{*}(1+d)\pi(\gamma)=1+d+A(\gamma)$. Notice that,
by proposition \ref{sobolev estimates}, $\|A(\gamma)\xi\|_{p}\leq M_{p}\|\xi\|_{p+1}$
where $M_{p}=(C|\gamma^{-1}\dot\gamma|_{p+\half{1}}+|c_{2}(\gamma^{-1})|)$.
If $\xi\in\H^{\infty}$, we have
\begin{equation}
\|\pi(\gamma)\xi\|_{n}=
\|(1+d)^{n}\pi(\gamma)\xi\|=
\|\pi(\gamma)^{*}(1+d)^{n}\pi(\gamma)\xi\|=
\|(1+d+A(\gamma))^{n}\xi\|
\end{equation}
and therefore
\begin{equation}\label{norm estimate}
\begin{split}
\|\pi(\gamma)\xi\|_{n}
&\leq\sum_{\substack{0\leq p_{i}+q_{i}\leq 1\\p_{i}+q_{i}=1}}
 \|(1+d)^{p_{1}}A(\gamma)^{q_{1}}\cdots
   (1+d)^{p_{n}}A(\gamma)^{q_{n}}\xi\| \\
&\leq\sum_{\substack{0\leq p_{i}+q_{i}\leq 1\\p_{i}+q_{i}=1}}
M_{p_{1}}^{q_{1}}\|(1+d)^{1+p_{2}}A(\gamma)^{q_{2}}\cdots
   (1+d)^{p_{n}}A(\gamma)^{q_{n}}\xi\| \\
&\leq\sum_{\substack{0\leq p_{i}+q_{i}\leq 1\\p_{i}+q_{i}=1}}
M_{p_{1}}^{q_{1}}M_{1+p_{2}}^{q_{2}}\cdots M_{(n-1)+p_{n}}^{q_{n}}\|(1+d)^{n}\xi\|\\
&\leq\sum_{\substack{0\leq p_{i}+q_{i}\leq 1\\p_{i}+q_{i}=1}}
M_{n-1}^{\sum q_{i}}\|(1+d)^{n}\xi\|\\
&=(1+M_{n-1})^{n}\|\xi\|_{n}\\
\end{split}
\end{equation}

{\it Joint continuity.} Fix $\gamma\in LG$. Let $U\subset LG$
be a neighborhood of the identity on which $\pi$ has a continuous lift to $U(\H)$ which
we denote by the same letter. Define a continuous lift of $\pi$ over $\gamma U$ by
$\pi(\zeta)=\pi(\gamma)\pi(\gamma^{-1}\zeta)$ where $\pi(\gamma)$ is arbitrary, so that
$\pi(\gamma)^{*}\pi(\zeta)=\pi(\gamma^{-1}\zeta)$. Now for $\xi,\eta\in\H^{n}$ and
$\zeta\in\gamma U$, we have
\begin{equation}
\begin{split}
\|\pi(\gamma)\xi-\pi(\zeta)\eta\|_{n}
&\leq
\|\pi(\gamma)(\xi-\eta)\|_{n}+\|(\pi(\gamma)-\pi(\zeta))\eta\|_{n}\\
&\leq
(1+M_{n-1}(\gamma))^{n}\|\xi-\eta\|_{n}+
(1+M_{n-1}(\gamma))^{n}\|(1-\pi(\gamma)^{*}\pi(\zeta))\eta\|_{n}\\
\end{split}
\end{equation}
and
\begin{equation}
\begin{split}
\|(1-\pi(\gamma)^{*}\pi(\zeta))\eta\|_{n}
&=
\|(1+d)^{n}(\pi(\gamma^{-1}\zeta)-1)\eta\|\\
&=
\|\pi(\gamma^{-1}\zeta)^{*}(1+d)^{n}(\pi(\gamma^{-1}\zeta)-1)\eta\|\\
&\leq
\|(1+d+A(\gamma^{-1}\zeta))^{n}\eta-(1+d)^{n}\eta\|+
\|(1-\pi(\gamma^{-1}\zeta)^{*})(1+d)^{n}\eta\|\\
&\leq
((1+M_{n-1}(\gamma^{-1}\zeta))^{n}-1)\|\eta\|_{n}+
\|(1-\pi(\gamma^{-1}\zeta)^{*})(1+d)^{n}\eta\|\\
\end{split}
\end{equation}
where the last inequality is proved like \eqref{norm estimate}. Since $c_{2}$ is
continuous and $c_{2}(1)=0$ the above tends to zero as $\zeta\rightarrow\gamma$
and $\eta\rightarrow\xi$ \halmos

\ssubsection{The crossed homomorphism corresponding to $\hsmooth$}

We compute below the crossed homomorphism corresponding to the joint
action of $LG$ and $L\g\rtimes i\IR d$ on $\hsmooth$, namely the map
$c:LG\times(L\g\rtimes i\IR d)\rightarrow \IR$ satisfying
\begin{equation}
\pi(\gamma)\pi(X)\pi(\gamma)^{*}=\pi(\gamma X\gamma^{-1})+ic(\gamma,X)
\end{equation}
as show that it agrees with the one given by Pressley and Segal
\cite[4.9.4]{PS}. We shall need two preliminary results.

\begin{proposition}\label{diff of exp}
Let $\G$ be a Fr\'echet Lie group with smooth exponential map. For
any $X,Y\in\Lie\G$, the left and right logarithmic derivatives of
the exponential map at $X$ are given by
\begin{align}
D^{L}_{X}\exp_{\G}\thinspace Y &:=
\left.\frac{d}{ds}\right|_{s=0}\exp_{\G}(-X)\exp_{\G}(X+sY) =
\int_{0}^{1}\exp_{\G}(-tX)Y\exp_{\G}(tX)dt \label{left log}\\
D^{R}_{X}\exp_{\G}\thinspace Y &:=
\left.\frac{d}{ds}\right|_{s=0}\exp_{\G}(X+sY)\exp_{\G}(-X) =
\int_{0}^{1}\exp_{\G}(tX)Y\exp_{\G}(-tX)dt \label{right log}
\end{align}
\end{proposition}
\proof We use Duhamel's formula \footnote{This idea is due to A. Selby}. Let
$A,B\in\Lie\G$ and set $f(t)=\exp_{\G}(-tA)\exp_{\G}(tB)$ so that
$\dot f=\exp_{\G}(-tA)(B-A)\exp_{\G}(tB)$. If $\phi$ is a smooth, real-valued
function defined on a neighborhood of $f([0,1])$, then
\begin{equation}
\phi(\exp_{\G}(-A)\exp_{\G}(B))-\phi(1)=
 \int_{0}^{1}\frac{d}{dt}\phi\circ f\thinspace dt=
 \int_{0}^{1}\exp_{\G}(-tA)(B-A)\exp_{\G}(tB)\phi\thinspace dt
\end{equation}
Set now $A=X$, $B=X+sY$ and let $\phi$ be defined and smooth near the identity.
Then,
\begin{equation}
\begin{split}
D^{L}_{X}\exp_{\G}\thinspace Y\phi
&=\left.\frac{d}{ds}\right|_{s=0}\phi(\exp_{\G}(-X)\exp_{\G}(X+sY))\\
&=\left.\frac{d}{ds}\right|_{s=0}s\int_{0}^{1}\exp_{\G}(-tX)Y\exp_{\G}(t(X+sY))
  \phi\thinspace dt\\
&=\int_{0}^{1}\exp_{\G}(-tX)Y\exp_{\G}(tX)\phi\thinspace dt
\end{split}
\end{equation}
and  \eqref{left log} follows since the above holds for any $\phi$ smooth near $1$.
The formula for the right logarithmic derivative follows similarly \halmos

\begin{corollary}\label{log diff exp}
Let $X\in L\g$ and $\gamma=\exp_{LG}(X)$. Then, the left and right logarithmic
derivatives of $\gamma$, seen as a smooth path in $G$ are given by
\begin{align}
\gamma^{-1}\dot\gamma&=
\int_{0}^{1}\exp_{LG}(-tX)\dot X\exp_{LG}(tX)dt \label{left exp log}\\
\dot\gamma\gamma^{-1}&=
\int_{0}^{1}\exp_{LG}(tX)\dot X\exp_{LG}(-tX)dt \label{right exp log}
\end{align}
\end{corollary}
\proof
Recall that if $p:I\rightarrow G$ is a smooth path, the left and right
logarithmic derivatives are the paths $I\rightarrow\g$ defined by
\begin{xalignat}{3}
p^{-1}\dot p(t)&:=\left.\frac{d}{dh}\right|_{h=0}p^{-1}(t)p(t+h)&
&\text{and}&
\dot pp^{-1}(t)&:=\left.\frac{d}{dh}\right|_{h=0}p(t+h)p^{-1}(t)
\end{xalignat}
so that $\dot pp^{-1}=\Ad(p)p^{-1}\dot p$. Therefore, if
$\gamma(\theta)=\exp_{G}(X(\theta))$, the chain rule and proposition 
\ref{diff of exp} for $\G=G$ show that the left and right hand-sides of
\eqref{left exp log} and \eqref{right exp log} are equal pointwise \halmos

\begin{theorem}\label{th:segal formulae}
For any $X\in L\g$ and $\gamma\in LG$, the following identities hold on $\H^{\infty}$
\begin{align}
\pi(\gamma)\pi(X)\pi(\gamma)^{*}&=
\pi(\Ad(\gamma)X)-
i\ell\int_{0}^{2\pi}\<\gamma^{-1}\dot\gamma,X\>\frac{d\theta}{2\pi}
\label{segal one}\\
\pi(\gamma)d\pi(\gamma)^{*}&=
d-i\pi(\dot\gamma\gamma^{-1})-
\half{\ell}\int_{0}^{2\pi}
\<\gamma^{-1}\dot\gamma,\gamma^{-1}\dot\gamma\>\frac{d\theta}{2\pi}
\label{segal two}\end{align}
\end{theorem}
\proof
Specialising lemma \ref{abstract commutation} to $\H^{\infty}$, we find for any
$\gamma\in LG$ and $X\in L\g\rtimes i\IR d$
\begin{equation}\label{basic}
\pi(\gamma)\pi(X)\pi(\gamma)^{*}=\pi(\Ad(\gamma)X)+ic(\gamma,X)
\end{equation}
from which \eqref{segal one} follows by the same proof as that of (iii) of proposition
\ref{ch:classification}.\ref{conjugated action}. \eqref{segal two} is more difficult
to establish because $d$ does not lie in the derived subalgebra of $L\g\rtimes i\IR d$.
To compute the values $c(\gamma,id)$, we shall show that the restriction of $c$ to
one--parameter subgroups in $LG$ satisfies a simple linear, first order ODE which may
be solved explicitly. It then easily follows that the
right hand-sides of \eqref{segal one}--\eqref{segal two} and $c$ coincide on the image of
$\exp_{LG}$ and hence on $LG$ because both are crossed homomorphisms
$LG\rightarrow(L\g\rtimes i\IR d)^{*}$ and $\exp_{LG}(L\g)$ generates $LG$.
We proceed in several steps.\\
{\it Continuity of $c$}. We claim that $c$ is jointly continuous in its arguments
or, equivalently that it can be viewed as a continuous map
$c:LG\rightarrow (L\g\rtimes i\IR d)^{*}$ where the latter is the locally convex
space of all continuous, linear forms on $L\g\rtimes i\IR d$ endowed with
the weak-* topology. {\it i.e.}~the topology of pointwise convergence.
To see this, let $(\gamma,X)\in LG\times(L\g\rtimes i\IR d)$ and choose a neighborhood
$\gamma\in U$ such that $\pi:U\rightarrow PU(\H)$ possesses a continuous lift to
$U(\H)$. Fix $\xi\in\H^{\infty}$ and for $\zeta\in U$ and
$Y\in L\g\rtimes i\IR d$, write
\begin{equation}
ic(\zeta,Y)\pi(\zeta)\xi=\pi(\zeta)\pi(Y)\xi-\pi(\Ad(\zeta)Y)\pi(\zeta)\xi
\end{equation}
so that
\begin{equation}
\begin{split}
|c(\zeta,Y)-c(\gamma,X)|
&\leq
\|c(\zeta,Y)\pi(\zeta)\xi-c(\gamma,X)\pi(\gamma)\xi\|+
 |c(\gamma,X)|\|\pi(\gamma)\xi-\pi(\zeta)\xi\|\\
&\leq
\|\pi(\zeta)\pi(Y)\xi-\pi(\gamma)\pi(X)\xi\|+
\|\pi(\Ad(\zeta)Y)\pi(\zeta)\xi-\pi(\Ad(\gamma)X)\pi(\gamma)\xi\|\\
&+
|c(\gamma,X)|\|\pi(\gamma)\xi-\pi(\zeta)\xi\|\\
&\leq
\|\pi(\zeta)(\pi(Y)-\pi(X))\xi\|+
\|(\pi(\zeta)-\pi(\gamma))\pi(X)\xi\|\\
&+
\|\pi(\Ad(\zeta)Y)(\pi(\zeta)-\pi(\gamma))\xi\|+
\|(\pi(\Ad(\zeta)Y)-\pi(\Ad(\gamma)X))\pi(\gamma)\xi\|\\
&+
|c(\gamma,X)|\|\pi(\gamma)\xi-\pi(\zeta)\xi\|\\
&\leq
C|Y-X|_{\half{1}}\|(1+d)\xi\|+
\|(\pi(\zeta)-\pi(\gamma))\pi(X)\xi\|\\
&+
C|\Ad(\zeta)Y|_{\half{1}}\|(1+d)(\pi(\zeta)-\pi(\gamma))\xi\|+
C|\Ad(\zeta)Y-\Ad(\gamma)X|_{\half{1}}\|(1+d)\pi(\gamma)\xi\|\\
&+
|c(\gamma,X)|\|\pi(\gamma)\xi-\pi(\zeta)\xi\|
\end{split}
\end{equation}
which, in view of the continuity of the action of $LG$ on $\H^{1}$, can be made
arbitrarily small.\\
{\it Functional Equation.} The following is an immediate corollary of \eqref{basic}
\begin{equation}\label{functional}
c(\gamma_{2}\gamma_{1},X)=c(\gamma_{1},X)+c(\gamma_{2},\Ad(\gamma_{1})X)
\end{equation}
in other words, $c$ is a continuous crossed homomorphism
$LG\rightarrow(L\g\rtimes i\IR d)^{*}$.\\
{\it Differentiability at $1$.} We claim that the restriction of $c$ to one-parameter
subgroups $\gamma_{t}=\exp_{LG}(tY)$ is differentiable at $t=0$. To see this, notice
that by proposition \ref{lifting}, the unitary $e^{t\pi(Y)}$ is a lift of $\pi(\gamma_{t})$.
Applying both sides of \eqref{basic} to $\xi\in\hsmooth$ and developing in Taylor series,
we get
\begin{equation}
\pi(\Ad(\gamma_{t})X)\xi=\pi(X+t[Y,X]+o(t))\xi=\pi(X)+t\pi([Y,X])\xi+o(t)
\end{equation}
where $\pi(o(t))\xi=o(t)$ follows from
$\|\pi(Z)\xi\|\leq C|Z|_{\half{1}}\|(1+d)\xi\|$. Next
\begin{equation}
\begin{split}
e^{t\pi(Y)}\pi(X)e^{-t\pi(Y)}\xi
&=e^{t\pi(Y)}\pi(X)(\xi-t\pi(Y)\xi+R(t))\\
&=e^{t\pi(Y)}(\pi(X)\xi-t\pi(X)\pi(Y)\xi+o(t))\\
&=\pi(X)\xi+t[\pi(Y),\pi(X)]\xi+o(t)
\end{split}
\end{equation}
since $\pi(X)R(t)=o(t)$. Indeed, the remainder term $R(t)$ is equal to
$t^{2}\int_{0}^{1}du u\int_{0}^{1}dv\pi(Y)^{2}e^{-tuv\pi(Y)}\xi$
and therefore
\begin{equation}
\|\pi(X)R(t)\|\leq Ct^{2}|X|_{\half{1}}
\|(1+d)\int_{0}^{1}du u\int_{0}^{1}dv\pi(Y)^{2}e^{-tuv\pi(Y)}\xi\|
\end{equation}
which is an $o(t)$ since the unitaries $e^{-tuv\pi(Y)}$ are uniformly
bounded in the $\|\cdot\|_{2}$ norm by proposition \ref{invariance}. Summarising,
\begin{equation}
ic(\gamma_{t},X)\xi=
t([\pi(Y),\pi(X)]-\pi([X,Y]))\xi+o(t)=it\ell B(X,Y)\xi+o(t)
\end{equation}
Thus, $c(\gamma_{t})$ is differentiable at $0$ and
\begin{equation}\label{diff at 0}
\left.\partial_{t}\right|_{t=0}c(\gamma(t),X)=\ell B(Y,X)
\end{equation}
{\it Ordinary Differential Equation.}
We shall now show that $c(\gamma_{t})$ is differentiable at any $t\in\IR$ and that
\begin{equation}\label{crossed ODE}
\partial_{t}c(\gamma_{t},X)=\ell B(Y,X)+c(\gamma_{t},[Y,X])
\end{equation}
Using \eqref{functional}, we get
$c(\gamma_{t+h},X)=c(\gamma_{h},X)+c(\gamma_{t},\Ad(\gamma_{h})X)$ and therefore
\begin{equation}
\frac{1}{h}\Bigl(c(\gamma_{t+h},X)-c(\gamma_{t},X))=
\frac{1}{h}c(\gamma_{h},X)+c(\gamma_{t},\frac{1}{h}(\Ad(\gamma_{h})X-X))
\end{equation}
Letting $h\rightarrow 0$ and using \eqref{diff at 0} and continuity in the second
variable yields \eqref{crossed ODE}.\\
{\it Formal solution of the ODE.} Somewhat changing notation, we write \eqref{crossed ODE}
as the inhomogeneous ODE
\begin{equation}\label{crossed ODE2}
\dot c_{t}=i_{Y}B-Yc_{t}
\end{equation}
where $i_{Y}B(X)=B(Y,X)$ and $Yc_{t}$ denotes the coadjoint action of $LG$ on
$(L\g\rtimes i\IR d)^{*}$. The underlying homogeneous equation is formally solved by
$c_{t}=e^{-tY}c_{0}$. Varying the constants, we set $c_{t}=e^{-tY}c_{0}(t)$, so that
\eqref{crossed ODE2} reads, in terms of $c_{0}$, $e^{-tY}\dot c_{0}=i_{Y}B$, $c_{0}(0)=0$
whence
$c_{0}(t)=\int_{0}^{t}e^{\tau Y}i_{Y}B\thinspace d\tau$.
In other words, reverting to our
former notation
\begin{equation}\label{solution}
c(\gamma_{t},X)=
\ell B(Y,\int_{0}^{t}\Ad(\gamma_{t-\tau})X\thinspace d\tau)=
\ell B(Y,\int_{0}^{t}\Ad(\gamma_{\tau})X\thinspace d\tau)
\end{equation}
{\it Existence and uniqueness of solution.} Using
$\frac{d}{dt}\Ad(\gamma_{t})X=[Y,\Ad(\gamma_{t})X]$, it is easy to verify that
\eqref{solution} defines a $C^{1}(\IR,(Lg\rtimes i\IR d)^{*})$ solution of
\eqref{crossed ODE} with initial condition $c_{0}=0$. If $c^{\prime}$ is another solution,
$\kappa_{t}=c_{t}-c^{\prime}_{t}$  satisfies $\dot\kappa_{t}(X)=\kappa_{t}([Y,X])$,
$\kappa_{0}=0$. Setting
$K_{t}(X)=\kappa_{t}(\Ad(\gamma_{-t})X)$, we see that $\dot K_{t}=0$ whence $K$ is
locally constant since it takes values in a locally convex vector space.
Thus, $\kappa\equiv 0$ and $c=c^{\prime}$.\\
{\it Reinterpretation of the solution.} We derive now a formula for $c$ which naturally
extends to the whole of $LG$. Let $Y,X\in L\g$ and set $\gamma_{t}=\exp_{LG}(tY)$.
By \eqref{solution} we have, in view of the $Ad_{G}$ invariance of $\<\cdot,\cdot\>$ and
corollary \ref{log diff exp}
\begin{equation}
\begin{split}
c(\exp_{LG}Y,X)
&=
-\ell\int_{0}^{2\pi}
\int_{0}^{1}dt\<\dot Y,\gamma_{t}X\gamma_{-t}\>
\frac{d\theta}{2\pi}\\
&=
-\ell\int_{0}^{2\pi}\frac{d\theta}{2\pi}
\int_{0}^{1}dt\<\gamma_{-t}\dot Y\gamma_{t},X\>
\frac{d\theta}{2\pi}\\
&=
-\ell\int_{0}^{2\pi}\frac{d\theta}{2\pi}\<\gamma_{1}^{-1}\dot\gamma_{1},X\>
\frac{d\theta}{2\pi}
\end{split}
\end{equation}
Similarly, using $\Ad(\gamma) id=id+\dot\gamma\gamma^{-1}$, corollary \ref{log diff exp}
and the symmetry of $\<\cdot,\cdot\>$,
\begin{equation}
\begin{split}
c(\exp_{LG}Y,i d)
&=
-\ell\int_{0}^{2\pi}
\int_{0}^{1}dt\<\dot Y,\dot\gamma_{t}\gamma_{t}^{-1}\>\frac{d\theta}{2\pi}\\
&=
-\ell\int_{0}^{2\pi}
\int_{0}^{1}\int_{0}^{t}\<\dot Y,\gamma_{u}\dot Y\gamma_{-u}\>dudt\frac{d\theta}{2\pi}\\
&=
-\half{\ell}\int_{0}^{2\pi}\biggl(
\int_{0}^{1}\int_{0}^{t} \<\dot Y,\gamma_{u}\dot Y\gamma_{-u}\>dudt+
\int_{0}^{1}\int_{-t}^{0}\<\dot Y,\gamma_{u}\dot Y\gamma_{-u}\>dudt
\biggr)\frac{d\theta}{2\pi}\\
&=
-\half{\ell}\int_{0}^{2\pi}\biggl(
\iint\limits_{0\leq v\leq t\leq 1}\<\dot Y,\gamma_{t-v}\dot Y\gamma_{v-t}\>dvdt+
\iint\limits_{0\leq t\leq v\leq 1}\<\dot Y,\gamma_{t-v}\dot Y\gamma_{v-t}\>dvdt
\biggr)\frac{d\theta}{2\pi}\\
&=
-\half{\ell}\int_{0}^{2\pi}
\<\gamma^{-1}_{1}\dot\gamma_{1},\gamma^{-1}_{1}\dot\gamma_{1}\>\frac{d\theta}{2\pi}
\end{split}
\end{equation}
Therefore, for any $\gamma$ in the image of $\exp_{LG}$, we have
\begin{equation}
c(\gamma,X+i\phi d)=
-\ell\int_{0}^{2\pi}
\<\gamma^{-1}\dot\gamma,X\>\frac{d\theta}{2\pi}
-\half{\ell\phi}\int_{0}^{2\pi}
\<\gamma^{-1}\dot\gamma,\gamma^{-1}\dot\gamma\>\frac{d\theta}{2\pi}
\end{equation}
The right hand-side extends in an obvious way to a function
$\widetilde c:LG\rightarrow(L\g\rtimes i\IR d)^{*}$ which is easily seen to satisfy
\eqref{functional}. Since any solution of that functional equation is uniquely 
determined by its restriction to a generating set in $LG$ and the image of $\exp_{LG}$
is such a set, we conclude that $c\equiv\widetilde c$ \halmos

\ssection{Central extensions of $LG$ arising from positive energy representations}
\label{se:central ext}

We show below that positive energy representations of equal level induce canonically
isomorphic central extensions of $LG$. As explained in chapter \ref{ch:classification}, 
this basic fact is needed to define the direct sum of these representations and
ultimately to show that the category of all positive energy representations at
a given level is abelian.
Our method shows in fact that these topological central extensions possess the structure
of real analytic Fr\'echet Lie groups, a fact which, surprisingly is easier to establish
than their smoothness.
It relies on proving that their local multiplication is analytic, which in turn is
derived from Nelson's analytic domination theorem \cite{Ne3} and the Sobolev estimates
of proposition \ref{sobolev estimates}.
Once their analyticity is established, the classification of the central extensions
reduces, by an elegant argument of Pressley and Segal \cite[4.4.1]{PS} to that of their
Lie algebra cocycle which we compute explicitly. As expected, the latter differs from
the fundamental cocycle by a factor equal to the level of the representation.
An alternative method which proves the smoothness of the central extensions only but
also applies to Diff$(S^{1})$ may be found in \cite{TL1}.

\ssubsection{The infinite dimensional projective unitary group}
\label{ss:PU}

Let $\H$ be a separable, infinite dimensional complex Hilbert space. We endow the unitary
group $U(\H)$ with the strong operator topology so that it becomes a complete, metrisable
topological group which moreover is contractible \cite[Lemme 3, page 251]{DiDo}. 
When the projective unitary group $PU(\H)=U(\H)/\T$ is endowed with the quotient
topology, the short exact sequence
\begin{equation}\label{extension}
1\rightarrow\T\rightarrow U(\H)\xrightarrow{p} PU(\H)\rightarrow 1
\end{equation}
possesses local continuous sections (see below) and may therefore be regarded as a
continuous, principal $\T$-bundle over $PU(\H)$. The associated long exact sequence of
homotopy groups shows that $PU(\H)$ is connected and simply--connected and therefore
that \eqref{extension} does not possess local homomorphic sections, contrary to its
finite--dimensional analogues. Indeed, any such would extend, by the local monodromy
principle, to a global section yielding an isomorphism $U(\H)\cong PU(\H)\times\T$ in
contradiction with the contractibility of $U(\H)$.\\

Local continuous sections of \eqref{extension} may be exhibited as follows \cite{Ba}.
Fix $0\neq\xi\in\H$ and consider the open set
\begin{equation}
U_{\xi}=\{[u]\in PU(\H)|\thinspace |(u\xi,\xi)|>0\}
\end{equation}
where $[u]$ is the equivalence class of $u\in U(\H)$ in $PU(\H)$, so that
$p^{-1}(U_{\xi})=\{u\in U(\H)|\thinspace (u\xi,\xi)\neq 0\}$. Define a function
\begin{equation}
\alpha_{\xi}:p^{-1}(U_{\xi})\longrightarrow\T,
\quad
u\longrightarrow\alpha_{\xi}(u)=\frac{(u\xi,\xi)}{|(u\xi,\xi)|}
\end{equation}
and notice that $\alpha_{\xi}(e^{i\theta}u)=e^{i\theta}\alpha_{\xi}(u)$ so that
the map $\phi:p^{-1}(U_{\xi})\rightarrow U_{\xi}\times\T$,
$\phi(u)=([u],\alpha_{\xi}(u))$
is a $\T$-equivariant local trivialisation with inverse
$\phi^{-1}([u],z)=u\alpha_{\xi}(u)^{-1}z$.
If $1\in V_{\xi}\subset U_{\xi}$ is a neighborhood such that
$V_{\xi}\cdot V_{\xi}\subset U_{\xi}$ and $V_{\xi}^{-1}=V_{\xi}$, we may use $\phi$
to define a local multiplication and inversion on $V_{\xi}\times\T$ by
\begin{align}
x\star y&=\phi(\phi^{-1}x\cdot\phi^{-1}y)\\
\II(x) &=\phi\bigl((\phi^{-1}x)^{-1}\bigr)
\end{align}
where $\cdot$ is the multiplication in $U(\H)$. Explicitly, using
$\phs{u^{*}}=\overline{\phs{u}}$,
\begin{align}
([u],z)\star([v],w)&=\Bigl([u\cdot v],zw
\frac{\alpha_{\xi}(uv)}{\alpha_{\xi}(u)\alpha_{\xi}(v)}\Bigr)\\
\II([u],z)&=([u^{-1}],z^{-1})
\end{align}
where the quotient $\alpha_{\xi}(uv)\alpha_{\xi}(u)^{-1}\alpha_{\xi}(v)^{-1}$ is
independent of the choice of the lifts $u,v$ of $[u],[v]\in PU(\H)$.

\newcommand {\LG}{{\mathcal L}G}

\ssubsection{The central extensions $\LG$ and their local multiplication}
\label{localmult}

Let $\pi:LG\rightarrow PU(\H)$ be a positive energy representation of $LG$ and
$\LG=\pi^{*}U(\H)$ the topological central extension
\begin{equation}\label{eq:LG extension}
1\rightarrow\T\rightarrow\LG\xrightarrow{p} LG\rightarrow{1}
\end{equation}
obtained by pulling back \eqref{extension} by $\pi$. Explicitly,
\begin{equation}
\LG=\{(g,u)\in LG\times U(\H)|\thinspace\pi(g)=[u]\}
\end{equation}
As in \S \ref{ss:PU}, we may trivialise $\LG$ over $\pi^{-1}(U_{\xi})$ to obtain
a homeomorphism $p^{-1}\pi^{-1}(U_{\xi})\rightarrow\pi^{-1}(U_{\xi})\times\T$ with
inverse $(\gamma,z)\rightarrow(\gamma,\pi(\gamma)\alpha_{\xi}(\pi(\gamma))^{-1}z)$
and corresponding local multiplication and inversion
\begin{align}
(\zeta,z)\star(\gamma,w)&=\Bigl(\zeta\gamma,zw
\frac{\alpha_{\xi}(\pi(\zeta)\pi(\gamma))}{\alpha_{\xi}(\pi(\zeta))\alpha_{\xi}(\pi(\gamma))}
\Bigr)\label{localone}\\
\II(\zeta,z)&=(\zeta^{-1},z^{-1})\label{localtwo}
\end{align}
We will prove in \S \ref{ss:real analytic} that $\xi$ may be chosen so that
$\star$ is (real) analytic near the identity. This in turn shows that the
above local trivialisation of $\LG$ extends to an analytic atlas and therefore
that \eqref{eq:LG extension} is a real analytic central extension of $LG$
in view of the following elementary

\begin{lemma}\label{local}
Let $\G$ be a connected topological group homeomorphic near 1 to an open
subset of a Fr\'echet space $E$. If the transported local multiplication
$\star:E\times E\rightarrow E$ and inversion $\II:E\rightarrow E$ are real
analytic maps, $\G$ possesses the structure of a real analytic Fr\'echet
Lie group modelled on $E$.
\end{lemma}
\proof
Let $\phi:U\subset\G\rightarrow\phi(U)\subset E$ be the local homeomorphism and
$1\in V\subset U$ a neighborhood such that $V\cdot V\subset U$ and $V^{-1}=V$.
Assume for simplicity that $\phi(1)=0$. The transported multiplication
$\star:\phi(V)\times\phi(V)\rightarrow\phi(U)$ and inversion
$\II:\phi(V)\rightarrow\phi(V)$ are defined by
\begin{align}
x\star y&=\phi(\phi^{-1}x\cdot\phi^{-1}y)\\
\II(x) &=\phi\bigl((\phi^{-1}x)^{-1}\bigr)
\end{align}
and are, by assumption, analytic.
Consider a neighborhood $1\in W\subset V$ such that $W\cdot W\subset V$
and $W^{-1}=W$. We define an atlas $\A=\Bigl\{(W_{g},\phi_{g})\Bigr\}$ indexed
by $g\in\G$ by setting $W_{g}=gW$, $\phi_{g}:W_{g}\rightarrow\phi(W)$,
$\phi_{g}(h)=\phi(g^{-1}h)$ so that $\phi_{g}^{-1}(x)=g\phi^{-1}(x)$. We claim
that $(\G,\A)$ has the structure of an analytic Fr\'echet Lie group. 

{\it Analyticity of the atlas}. Assume $W_{g_{1}}\cap W_{g_{2}}\neq 0$.
Then, $g_{2}^{-1}g_{1}\in W\cdot W^{-1}\subset V$ may be written as $\phi^{-1}(y)$
with $y\in\phi(V)$ and therefore for $x\in E$ close to the origin, the transition
map
\begin{equation}
\tau_{g_{2}g_{1}}(x)=
\phi_{g_{2}}\circ\phi_{g_{1}}^{-1}(x)=
\phi(g_{2}^{-1}g_{1}\phi^{-1}(x))=
\phi(\phi^{-1}(y)\phi^{-1}(x))=
y\star x
\end{equation}
is analytic.

{\it Analyticity of the adjunction}. Let $g\in\G$ and define the transported adjunction
$\tad(g)$ about $0\in E$ by $\tad(g)x=\phi(g\phi^{-1}(x)g^{-1})$. If $g$ lies in a
suitably small neighborhood $1\in S\subset W$, $\tad(g)$ is smooth
(analytic) near 0 since $\tad(g)x=y\star x\star\II(y)$ with $y=\phi(g)$. Similarly, if
$g=g_{1}\cdots g_{k}\in S^{k}$, then $\tad(g)=\tad(g_{1})\circ\cdots\circ\tad(g_{k})$
is analytic near 1. Since $\G$ is connected, $\bigcup_{k}S^{k}=\G$ and 
$\tad(g)$ is analytic near 0 for any $g\in\G$.

{\it Analyticity of multiplication}. Let $g,h\in\G$ and $x,y\in E$ small enough,
then
\begin{equation}
\phi_{gh}\Bigl(\phi_{g}^{-1}(x)\phi_{h}^{-1}(y)\Bigr)=
\phi(h^{-1}\phi^{-1}(x)h\phi^{-1}(y))=
(\tad(h^{-1}) x)\star y
\end{equation}
is analytic.

{\it Analyticity of inversion}. Let $g\in\G$ and $x\in E$ small enough, then
\begin{equation}
\phi_{g^{-1}}(\phi_{g}^{-1}(x)^{-1})=
\phi(g\phi^{-1}(x)^{-1}g^{-1})=
\tad(g)\circ\II\thinspace(x)
\end{equation}
is analytic \halmos

\ssubsection{Analytic Domination of $L\g$ in Positive Energy Representations}

We will need Goodman's refinement of Nelson's fundamental analytic domination
theorem \cite[Thm. 1]{Ne3}

\begin{theorem}\label{domination}
Let $\{X_{i}\}_{i\in I}$ be a family of endomorphisms of a normed vector space
$V$ and $A\in\End(V)$ such that, for any $\xi\in V$
\begin{align}
\|X_{i}\xi\|&\leq\|A\xi\| \label{domi}\\
\|\ad X_{i_{n}}\cdots\ad X_{i_{1}}(A)\xi\|&\leq n!\|A\xi\| \label{ad domi}
\end{align}
for any $i\in I$, $n\in\IN$ and finite subset $\{i_{1},\ldots,i_{n}\}\subset I$.
Then there exists a constant $\widetilde M>0$ such that the inequalities
\begin{equation}
\|X_{i_{n}}\cdots X_{i_{1}}\xi\|\leq n!\widetilde M^{n}
\end{equation}
hold whenever $\xi\in V$ satisfies $\|A^{m}\xi\|\leq m!$ for any $m\in\IN$.
\end{theorem}
\proof The above is simply a reformulation of lemma $2^{\prime}$ of \cite{Go}
(where the author simply omitted to state explicitly the independence of
$\widetilde M$ on $\xi$). The proof is obtained by combining the sketch proof
of lemma $2^{\prime}$ and the proof of lemma 2 in \cite{Go} \halmos

\begin{corollary}\label{estimate}
Let $\H$ be a positive energy representation of $LG$ with finite energy subspace $\Hfin$.
Then, there exists a constant $\infty>M>0$ such that for any $\xi\in\Hfin$
and $X_{1},\ldots,X_{n}\in L\g$, the following inequalities hold
\begin{equation}
\|\pi(X_{1})\cdots\pi(X_{n})\xi\|\leq
\lambda_{\xi}\norm{X_{1}}\cdots\norm{X_{n}}M^{n}n!
\end{equation}
for some constant $0<\lambda_{\xi}<\infty$ depending on $\xi$.
\end{corollary}
\proof
Set $V=\H^{\infty}$  and $A=1+d$ in theorem \ref{domination} where
$d$ is the infinitesimal generator of rotations on $\H$. We will
show that the inequalities \eqref{domi}--\eqref{ad domi} hold for a suitable neighborhood
of $0$ in $L\g$.
From the estimates of proposition \ref{sobolev estimates} and the commutation
relations \eqref{commutation 1}--\eqref{commutation 2}, we find that
$\|\pi(X)\xi\|\leq C|X|_{\half{1}}\|A\xi\|$ for any $\xi\in V$ and, for $n\geq 2$
\begin{equation}
\begin{split}
\|\ad\pi(X_{n})\cdots\ad\pi(X_{1})A\xi\|
&\leq C|\ad X_{n}\cdots\ad X_{2}\dot X_{1}|_{\half{1}}\|A\xi\|+
      \ell|B(X_{n},\ad X_{n-1}\cdots\ad X_{2}\dot X_{1})|\thinspace\|\xi\| \\
&\leq CC_{\g}^{n-1}
      |X_{n}|_{\half{1}}\cdots|X_{2}|_{\half{1}}|\dot X_{1}|_{\half{1}}\|A\xi\|+
      \ell C_{\g}^{n-2}
      |X_{n}|_{1}|X_{n-1}|_{0}\cdots|X_{2}|_{0}|\dot X_{1}|_{0}\|\xi\| \\
&\leq C_{g}^{n-2}(CC_{\g}+\ell)|X_{n}|_{\half{3}}\cdots|X_{1}|_{\half{3}}\|A\xi\|
\end{split}
\end{equation}
where $\ell$ is the level of $\H$ and we have used $|\dot X|_{t}\leq|X|_{t+1}$,
$|[X,Y]|_{t}\leq C_{\g}|X|_{t}|Y|_{t}$ and $|B(X,Y)|\leq|X|_{1}|Y|_{0}$.
Thus, for $r>0$ small enough, the family
$\{\pi(X)|\thinspace X\in L\g, |X|_{\half{3}}\leq r\}$ satisfies the
inequalities \eqref{domi}--\eqref{ad domi} of theorem \ref{domination} for any
$\xi\in V$.
Let now $\Hfin\ni\xi=\sum_{k}\xi_{k}$ be a finite sum of eigenvectors of
$A$ with eigenvalues $\mu_{k}$. Set $\sigma=\max_{k}\mu_{k}$ and
$\lambda_{\xi}=\sup_{m}\frac{\sigma^{m}}{m!}$, then
$\|A^{m}\lambda_{\xi}^{-1}\xi\|\leq m!$ for all $m\in\IN$ and therefore,
for any $X_{1},\ldots X_{n}\in L_{\g}$
\begin{equation}
\|\pi(X_{n})\cdots\pi(X_{1})\xi\|
=\lambda_{\xi}r^{-n}|X_{n}|_{\half{3}}\cdots|X_{1}|_{\half{3}}
 \|r\frac{\pi(X_{n})}{|X_{n}|_{\half{3}}}\cdots
   r\frac{\pi(X_{1})}{|X_{1}|_{\half{3}}}\frac{\xi}{\lambda_{\xi}}\|
\leq\lambda_{\xi}|X_{n}|_{\half{3}}\cdots|X_{1}|_{\half{3}}M^{n}n!
\end{equation}
where $M=\widetilde Mr^{-1}$ \halmos

\ssubsection{$\LG$ as a real analytic Fr\'echet Lie group}\label{ss:real analytic}

The following is a simple consequence of the spectral theorem

\begin{lemma}\label{exponential}
Let $X,Y$ be two essentially skew--adjoint operators on a Hilbert space $\H$
defined on an invariant domain $\D$. Denote by $\overline{X}$, $\overline{Y}$ the
closures of $X$ and $Y$ and let the unitaries $e^{\overline{X}}$, $e^{\overline{Y}}$
be defined by the spectral theorem. If $\xi\in\D$, the identities
\begin{enumerate}
\item[(i)]$\displaystyle{\sum_{n\geq 0}\frac{Y^{n}}{n!}\xi=e^{\overline{Y}}\xi}$
\item[(ii)]$\displaystyle{\sum_{m,n\geq 0}\frac{X^{m}Y^{n}}{m!n!}\xi=
			  e^{\overline{X}}e^{\overline{Y}}\xi}$
\end{enumerate}
hold whenever their left-hand sides are absolutely convergent.
\end{lemma}
\proof
We may suppose that $X$ and $Y$ are closed with invariant subdomain $\D$.

(i) By the spectral theorem, we may assume $\H=L^{2}(\Omega,\mu)$ with $Y$
acting as multiplication by a measurable, $i\IR$-valued, $\mu$-almost everywhere
finite function $f$. Then,
$\D(Y^{n})=\{g\in\H|\int_{\Omega}|f|^{2n}|g|^{2}d\mu<\infty\}$.
If $\xi\in\bigcap_{n}\D(Y^{n})$ is such that the left-hand side of (i) is absolutely
convergent, the sequence of functions $\eta_{N}=\sum_{n=0}^{N}\frac{1}{n!}f^{n}\xi$
converges to some $\eta\in\H$. Therefore \cite[Theorem 3.12]{Ru}, for almost
every $\omega\in\Omega$, $|f(\omega)|,|\xi(\omega)|<\infty$ and
$\eta_{N}(\omega)\rightarrow\eta(\omega)$.
However,
\begin{equation}
\eta_{N}(\omega)=
\sum_{n=0}^{N}\frac{1}{n!}f^{n}(w)\xi(w)\rightarrow
e^{f(\omega)}\xi(\omega)
\end{equation}
whence $\eta=e^{f}\xi=e^{Y}\xi$.\\
(ii) Assume that $\xi\in\D$ is such that the left hand-side of (ii) is absolutely
convergent and set $\eta=e^{Y}\xi$. We will show inductively that $\eta\in\D(X^{m})$
and that $X^{m}\eta=\sum_{n\geq 0}X^{m}\frac{Y^{n}}{n!}\xi$ for any $m\in\IN$. The
case $m=0$ is settled by (i). If $m\geq 1$, set
$\theta_{N}=\sum_{n=0}^{N}X^{m-1}\frac{Y^{n}}{n!}\xi$, then by assumption
$\theta_{N}\in\D\subset\D(X)$ and by induction $\theta_{N}\rightarrow X^{m-1}\eta$.
Moreover, by assumption $X\theta_{N}$ is Cauchy and hence by closedness of $X$,
$X^{m-1}\eta\in\D(X)$, {\it i.e.}~$\eta\in\D(X^{m})$ and
$X^{m}\eta=\sum_{n\geq 0}X^{m}\frac{Y^{n}}{n!}\xi$. Therefore,
$\eta\in\bigcap_{m}D(X^{m})$ and the series $\sum_{m\geq 0}\frac{X^{m}}{m!}\eta$
is absolutely convergent and therefore equal to $e^{X}\eta$. Finally,
\begin{equation}
\sum_{m,n\geq 0}\frac{X^{m}Y^{n}}{m!n!}\xi=
\sum_{m\geq 0}\sum_{n\geq 0}\frac{X^{m}Y^{n}}{m!n!}\xi=
\sum_{m\geq 0}\frac{X^{m}}{m!}e^{Y}\xi=
e^{X}e^{Y}\xi
\end{equation}
\halmos

\begin{theorem}\label{analytic}
Let $(\pi,\H)$ be a positive energy representation of $LG$ and
\begin{equation}\label{pull back}
1\rightarrow\T\rightarrow\LG\rightarrow LG\rightarrow 1
\end{equation}
the corresponding topological central extension induced by $\pi$. Then
\begin{enumerate}
\item $\LG$ possesses the structure of a real analytic Fr\'echet Lie group modelled
locally on $L\g\oplus i\IR$ such that \eqref{pull back} is a real analytic central
extension of $LG$.
\item $\H$ possesses a dense set of analytic vectors for the action of $LG$.
\end{enumerate}
\end{theorem}
\proof
(i) By lemma \ref{local}, we need only prove that the vector $\xi$ used in the local
trivialisation of $\LG$ may be chosen so that the local multiplication \eqref{localone}
is analytic. Fix $\xi\in\Hfin$ and $\delta\leq(2M)^{-1}$ and let $X,Y\in L\g$,
$|X|_{\half{3}},|Y|_{\half{3}}<\delta$. Then, by corollary \ref{estimate},
\begin{equation}\label{cvgce}
\sum_{m,n\geq 0}\frac{\|\pi(X)^{m}\pi(Y)^{n}\xi\|}{m!n!}
\leq\lambda_{\xi}\sum_{m,n\geq 0}\C{m}{n}|X|_{\half{3}}^{m}|Y|_{\half{3}}^{n}M^{m+n}
=\lambda_{\xi}\sum_{k\geq 0}(M|X|_{\half{3}}+M|Y|_{\half{3}})^{k}
<\infty
\end{equation}
Therefore, if $B_{\delta}=\{X\in L\g|\thinspace|X|_{\half{3}}<\delta\}$ then,
in view of lemma \ref{exponential}, the function
$f_{\xi}:L\g\times L\g\rightarrow\H$,
$(X,Y)\rightarrow e^{\pi(X)}e^{\pi(Y)}\xi$ is analytic on
$B_{\delta}\times B_{\delta}$. Moreover, since
\begin{equation}
\|e^{\pi(X)}e^{\pi(Y)}\xi-\xi\|\leq
\lambda_{\xi}\sum_{k\geq 1}(M|X|_{\half{3}}+M|Y|_{\half{3}})^{k}
\end{equation}
we may take $\delta$ small enough so that
$(e^{\pi(X)}e^{\pi(Y)}\xi,\xi)\neq 0$ for any
$(X,Y)\in B_{\delta}\times B_{\delta}$ and it follows that the function
\begin{equation}
(X,Y)\rightarrow
\frac{\phs{e^{\pi(X)}e^{\pi(Y)}}}
     {\phs{e^{\pi(X)}}\phs{e^{\pi(Y)}}}
\end{equation}
is analytic on $B_{\delta}\times B_{\delta}$.
Since $LG$ has a local analytic logarithm and for any $X\in L\g$, $e^{\pi(X)}$ is
a lift of $\pi(\exp_{LG}X)$ by proposition \ref{lift one}, the local multiplication
\eqref{localone} is analytic.

(ii) Let $\xi$ as in (i) and $\eta\in\Hfin$. We claim that $\eta$ is analytic for the
canonical action of $LG$ given by $(\gamma,u)\xi=u\xi$. It is sufficient to check this
near 1 where we may trivialise $\LG$ using $\xi$ as in \S\ref{localmult}. Then
\begin{equation}
(\gamma,z)\eta=
(\gamma,\pi(\gamma)\alpha_{\xi}(\pi(\gamma))^{-1}z)\eta=
\pi(\gamma)\alpha_{\xi}(\pi(\gamma))^{-1}z\eta
\end{equation}
For $\gamma\in LG$ near 1, this is equal to
\begin{equation}
z\alpha_{\xi}(e^{-\pi(X)})e^{\pi(X)}\eta
\end{equation}
where $X=\log(\gamma)$ and is therefore analytic \halmos

\newpage
\begin{proposition}\label{classify extension}
Let $(\pi_{i},\H_{i})$ be positive energy representations of levels $\ell_{i}$.
Then,
\begin{enumerate}
\item The Lie algebra cocycle corresponding to $\L^{i}G=\pi_{i}^{*}U(\H_{i})$
is $\ell_{i}$ times the fundamental cocycle
\begin{equation}
iB(X,Y)=i\int_{0}^{2\pi}\<X,\dot Y\>\frac{d\theta}{2\pi}
\end{equation}
\item $\L^{1}G$ and $\L^{2}G$ are (canonically) isomorphic if, and only if
$\ell_{1}=\ell_{2}$.
\end{enumerate}
\end{proposition}
\proof
(i) Fix $\xi\in\Hfin$ and trivialise $\L^{i}G$ near 1 as in \S\ref{localmult}.
We begin by computing the corresponding local adjoint action.
By \eqref{localone}--\eqref{localtwo}
\begin{equation}\label{adjunction}
\begin{split}
(g,z)\star(h,w)\star\II(g,z)
&=
\biggl(gh,zw\frac{\phs{\pi(g)\pi(h)}}{\phs{\pi(g)}\phs{\pi(h)}}\biggr)
\star(g^{-1},z^{-1})\\
&=
\biggl(ghg^{-1},w\frac{\phs{\pi(g)\pi(h)\pi(g)^{*}}}{\phs{\pi(h)}}\biggr)
\end{split}
\end{equation}
where the lifts $\pi(g),\pi(h)\in U(\H)$ are arbitrary and we have
chosen $\pi(g)\pi(h)$ as a lift of $[\pi(gh)]$ and $\pi(g)^{*}$ as
a lift of $[\pi(g)^{-1}]$.
Set $h:=h_{t}=\exp_{LG}(tX)$ for some $X\in L\g$ and choose $e^{t\pi(X)}$ as a lift of
$\pi(h_{t})$. We will compute the derivative of the right hand-side of \eqref{adjunction}
at $t=0$.
Since $\xi\in\Hfin\subset\hsmooth$,
\begin{equation}
\der e^{t\pi(X)}\xi=\pi(X)\xi
\end{equation}
and therefore, assuming that $\|\xi\|=1$
\begin{equation}
\der\phs{e^{t\pi(X)}}=
\der\frac{(e^{t\pi(X)}\xi,\xi)}
{\sqrt{(e^{t\pi(X)}\xi,\xi)
 (\xi,e^{t\pi(X)}\xi)}}=
(\pi(X)\xi,\xi)
\end{equation}
since, by skew--adjointness, the function at the denominator is even.
Similarly, $\pi(g)^{*}$ leaves $\H^{\infty}$ invariant by proposition \ref{invariance}
and therefore
\begin{align}
\der\pi(g)e^{t\pi(X)}\pi(g)^{*}\xi&=\pi(g)\pi(X)\pi(g)^{*}\xi\\
\intertext{whence}
\der\phs{\pi(g)e^{t\pi(X)}\pi(g)^{*}}&=(\pi(g)\pi(X)\pi(g)^{*}\xi,\xi)
\end{align}
Thus
\begin{equation}
\der\frac{\phs{\pi(g)\pi(h_{t})\pi(g)^{*}}}{\phs{\pi(h_{t})}}=
(\pi(g)\pi(X)\pi(g)^{*}\xi,\xi)-(\pi(X)\xi,\xi)
\end{equation}
and the local adjoint action of $\L^{i}G$ is
\begin{equation}
(g,z)(X,ix)(g,z)^{-1}=
(gXg^{-1},ix+(\pi(g)\pi(X)\pi(g)^{*}\xi,\xi)-(\pi(X)\xi,\xi))
\end{equation}
where we have identified the Lie algebra of $\T$ with $i\IR$.
To compute the Lie bracket, set $g:=g_{s}=\exp_{LG}(sY)$, $Y\in L\g$ and choose the
lift $\pi(g_{s})=e^{s\pi(Y)}$. Then $\ders\pi(g_{s})^{*}\xi=-\pi(Y)\xi$ and therefore
\begin{equation}
\ders(\pi(g_{s})\pi(X)\pi(g_{s})^{*}\xi,\xi)=
\ders(\pi(X)\pi(g_{s})^{*}\xi,\pi(g_{s})^{*}\xi)=
([\pi(Y),\pi(X)]\xi,\xi)
\end{equation}
and it follows by the commutation relations \eqref{commutation 2} that
the bracket on the Lie algebra of $\L^{i}G$ is given by
\begin{equation}
[(X,ix),(Y,iy)]=([X,Y],i\ell_{i}B(X,Y)+d\beta_{\xi}(X,Y))
\end{equation}
where $\beta_{\xi}$ is the linear form on $L\g$ given by $\beta_{\xi}(Z)=(\pi(Z)\xi,\xi)$
and $d\beta_{\xi}(X,Y)=\beta_{\xi}([X,Y])$. Notice that $\beta_{\xi}$ is continuous
since $\|\pi(Z)\xi\|\leq C|Z|_{\half{1}}\|\xi\|_{1}$.
Thus, the cocycle classifying the Lie algebra of $\L^{i}G$ as a central extension of
$L\g$ is $i\ell_{i}B+d\beta_{\xi}$ which is cohomologous to $i\ell_{i}B$.

(ii)
follows from (i) and the classification of central extensions of $LG$
\cite[4.4.1 (ii)]{PS}. The isomorphism $\L^{1}G\cong\L^{2}G$ is unique
since any two differ by an element of $\Hom(LG,\T)=\{1\}$ \cite[3.4.1]{PS}
\halmos\\

\remark
Denote by $\wt{LG}$ the central extension of $LG$ whose Lie algebra cocycle is
$iB(\cdot,\cdot)$. By \cite[4.4.6]{PS}, $\wt{LG}$ is the universal central extension
of $LG$ and it is easy to deduce from (i) of proposition \ref{classify extension}
that $\LG\cong\wt{LG}/Z_{\ell}$ where $\ell$ is the level of $\H$ and
$Z_{\ell}\subset\T\subset\wt{LG}$ is the group of $\ell$ roots of unity. It
follows that any positive energy representation of $LG$ at level $\ell$ gives rise
to a continuous, unitary representation of $\wt{LG}$ such that the centre of $\wt{LG}$
acts by the character $z\rightarrow z^{\ell}$. Moreover, as conjectured by Pressley and
Segal \cite[\S 9.3]{PS}, $\H$ possesses by theorem \ref{analytic} a dense subspace of
analytic, and {\it a fortiori} smooth vector for the action of $\wt{LG}$, 




\renewcommand {\AA}{\mathfrak A}
\newcommand {\Vj}{V_{j}}
\newcommand {\VJ}{V_{j}}
\newcommand {\VC}{V_{\IC}}
\newcommand {\VCH}{\VC^{1,0}}
\newcommand {\VCA}{\VC^{0,1}}
\newcommand {\HJ}{\H_{J}}
\newcommand {\HC}{\H_{\IC}}
\newcommand {\HIR}{\H_{\IR}}
\newcommand {\HCH}{\HC^{1,0}}
\newcommand {\HCA}{\HC^{0,1}}
\newcommand {\HS}{\operatorname{\scriptscriptstyle{HS}}}
\newcommand {\HR}{\H_{\operatorname{\scriptscriptstyle{R}}}}
\newcommand {\HNS}{\H_{\operatorname{\scriptscriptstyle{NS}}}}
\newcommand {\FJ}{\mathcal F_{J}}
\newcommand {\FP}{\mathcal F_{P}}
\newcommand {\FR}{\F_{R}}
\newcommand {\FRE}{\F_{R,0}}
\newcommand {\FRO}{\F_{R,1}}
\newcommand {\FNS}{\F_{NS}}
\newcommand {\FNSE}{\F_{NS,0}}
\newcommand {\FNSO}{\F_{NS,1}}
\newcommand {\ures}{U_{\operatorname{res}}}
\newcommand {\ores}{O_{\operatorname{res}}}

\newcommand {\ext}{\Lambda}
\newcommand {\evenext}{\ext_{0}}
\newcommand {\oddext}{\ext_{1}}
\newcommand {\tr}{\operatorname{tr}}
\newcommand {\ffine}{\F^{0}}
\newcommand {\ffino}{\F^{1}}
\newcommand {\cv}{C(V)}
\newcommand {\cvo}{C(V)^0}
\newcommand {\calg}{C^{\operatorname{alg}}(\H)}

\newcommand {\NS}{Neveu--Schwarz }


\chapter{Fermionic construction of level 1 representations}
\label{ch:fermionic}

This chapter is devoted to the well--known {\it quark model} or Fermionic construction
of the level 1 positive energy representations of $L\Spin_{2n}$. This parallels the
Clifford algebra construction of the spin modules of $\Spin_{2n}$ and realises the
level 1 representations as summands of two distinct Fermionic Fock spaces, the Ramond
and Neveu--Schwarz sectors $\FR,\FNS$.\\

Our main interest in this construction comes from Algebraic Quantum Field Theory. We
shall prove in chapter \ref{ch:loc loops} that the von Neumann algebras generated by
the groups $L_{I}\Spin_{2n}$ of loops supported in a proper interval $I\subset S^{1}$
in positive energy representations are type III$_{1}$ factors and satisfy Haag duality
in the vacuum sector. These results are necessary for the very definition of
fusion and will be derived from their well--known analogues for the free Fermi field.
Another feature of the quark model is that the vector primary field is realised as a
Fermionic field and therefore obeys $L^{2}$ bounds. The continuity properties of other
primary fields, most notably the spin ones, are not as easily derived in the fermionic
picture and will be proved in chapter \ref{ch:sobolev fields} using the equivalent
bosonic construction of level 1 representations of $L\Spin_{2n}$.\\

In section \ref{se:fd spinors}, we review the Clifford algebra construction of the
$\SO_{2n}$ spin modules.
We do this in some detail as an explicit description of these modules will be needed
elsewhere. In section \ref{se:quark model}, we construct the level 1 representations
of $L\Spin_{2n}$ inside the Ramond and \NS Fock spaces $\FR,\FNS$ by using a global
version of the quark model. 
In section \ref{se:currents}, we identify the abstract action of $\lpol\so_{2n}$
on $\FR,\FNS$ corresponding to that of $L\Spin_{2n}$ with well--known bilinear
expressions in the Fermi field and use this to prove the finite reducibility of
$\FR$ and $\FNS$ under $L\Spin_{2n}$.
Finally, in section \ref{se:vector fields}, we prove that the vector primary
field defines a bounded operator--valued distribution by identifying it with
a Fermi field.

\ssection{The finite--dimensional spin representations}\label{se:fd spinors}

We describe below the construction of the $\SO_{2n}$ spin modules following Brauer
and Weyl's original lines \cite{BrWe}. This stems from the simple observation that
the Clifford algebra $C(V)$ of $V=\IR^{2n}$ is a matrix algebra and therefore
possesses a unique irreducible representation $\F$. Conjugating it by the
automorphic action of $\SO_{2n}$ on $C(V)$ thus leads to
unitarily equivalent ones and therefore to a projective representation of $\SO_{2n}$
on $\F$. The latter lifts to an ordinary representation of the universal cover
$\Spin_{2n}$ of $\SO_{2n}$ which decomposes as the sum of the two spin modules.

\ssubsection{Clifford algebras and CAR algebras}\label{ss:clifford and CAR}

Let $V=\IR^{2n}$ with inner product $B(\cdot,\cdot)$ and $C(V)$ the corresponding Clifford
algebra, that is the complex *--algebra generated by the self--adjoint symbols $\psi(v)$,
$v\in V$ satisfying
\begin{equation}\label{eq:clifford}
\{\psi(u),\psi(v)\}=\psi(u)\psi(v)+\psi(v)\psi(u)=2B(u,v)
\end{equation}
$C(V)$ is naturally $\IZ_{2}$--graded by decreeing that the $\psi(v)$ are of degree
one. We claim that $C(V)$ is isomorphic to the matrix algebra $M_{2^{n}}(\IC)$. This
may be seen by using the factorisation property
$C(V_{1}\oplus V_{2})=C(V_{1})\wh\otimes C(V_{2})$, where $\wh\otimes$ is the
graded tensor product and the fact that $C(\IR^{2})=M_{2}(\IC)$
but is more conveniently obtained by giving a different presentation of $C(V)$
which naturally suggests a module for it.\\

Introduce for this purpose an orthogonal complex structure $j$ on $V$ and denote by
$\Vj$ the corresponding complex vector space with hermitian inner product
$(u,v)=B(u,v)+iB(u,jv)$. The {\it canonical anticommutation relations} or CAR algebra
$\AA(\Vj)$ is the $\IZ_{2}$--graded, complex *--algebra generated by the degree 1,
$\IC$--linear symbols $c(v)$, $v\in\Vj$ subject to the relations
\begin{xalignat}{2}\label{eq:CAR}
\{c(u),c(v)\}&=0&
\{c(u),c(v)^{*}\}&=(u,v)
\end{xalignat}
These suggest that a natural $\AA(\Vj)$--module is the exterior algebra
$\ext\Vj=\bigoplus_{k=0}^{n}\ext^{k}\Vj$ with the $c(v)$ acting as creation
operators, namely
\begin{align}
c(v)\medspace v_{1}\wedge\cdots\wedge v_{m}&=v\wedge v_{1}\wedge\cdots\wedge v_{m}
\label{eq:repr 1}\\
\intertext{so that}
c(v)^{*}\medspace v_{1}\wedge\cdots\wedge v_{m}&=\sum_{j=1}^{m}
(-1)^{j-1}(v_{j},v)v_{1}\wedge\cdots\wedge\wh{v_{i}}\wedge\cdots\wedge v_{m}
\label{eq:repr 2}
\end{align}

\begin{lemma}\label{le:fd irred}
The representation of $\AA(\Vj)$ defined by \eqref{eq:repr 1}--\eqref{eq:repr 2}
is irreducible and yields an isomorphism $\AA(\Vj)\cong\End(\ext\Vj)$.
\end{lemma}
\proof Let $\pi$ be the above representation. 
Notice that the number operator $N$ acting as multiplication by $k$ on $\ext^{k}\Vj$
may be written as
\begin{equation}
N=\sum_{k}c(v_{k})c(v_{k})^{*}
\end{equation}
where $v_{k}$ is any complex basis of $\Vj$. Thus, if $T\in\End(\ext\Vj)$ commutes
with $\AA(\Vj)$, then $[T,N]=0$ and $T$ leaves $\ext^{0}\Vj$ invariant. Thus,
$T\Omega=\lambda\Omega$ where
$\Omega\in\ext^{0}\Vj$ is a generator and therefore $T\equiv\lambda$ by cyclicity of
$\Omega$. Thus, $\pi$ is irreducible and, by the double commutant theorem,
$\pi(\AA(\Vj))=\End(\ext\Vj)$ is of dimension $2^{2n}$. Since $\AA(\Vj)$
is spanned by the monomials
\begin{equation}
c(v_{i_{1}})\cdots c(v_{i_{k}})c(v_{\ell_{1}})^{*}\cdots c(v_{\ell_{m}})^{*}
\end{equation}
where $v_{1}\ldots v_{n}$ is a complex basis of $\VC$ and
$1\leq i_{1}<\cdots<i_{k}\leq n$, $1\leq\ell_{1}<\cdots<\ell_{m}\leq n$, its
dimension is bounded above by $2^{2n}$ and $\pi$ is an isomorphism \halmos\\

Returning to the Clifford algebra $C(V)$, we notice that it is canonically
isomorphic to $\AA(\Vj)$ by
\begin{xalignat}{2}\label{eq:fd isom}
\psi(v)&\rightarrow c(v)+c(v)^{*}&
   c(v)&\rightarrow \half{1}\Bigl(\psi(v)-i\psi(jv)\Bigr)
\end{xalignat}
and is therefore a matrix algebra whose unique irreducible module may be
identified with $\ext\Vj$. 

\begin{proposition}
There exists a projective unitary representation $\Gamma$ of $\SO_{2n}$ on
$\ext\Vj$ satisfying, and uniquely determined by
\begin{equation}\label{eq:fd covariance}
\Gamma(R)\psi(v)\Gamma(R)^{*}=\psi(Rv)
\end{equation}
and extending the canonical action of the unitary group $U(\Vj)\subset\SO_{2n}$ given by
\begin{equation}
\Gamma(U)\medspace v_{1}\wedge\cdots\wedge v_{k}=Uv_{1}\wedge\cdots\wedge Uv_{k}
\end{equation}
Moreover, $\Gamma$ leaves the even and odd subspaces of $\ext\Vj$, namely
$\evenext\Vj=\bigoplus_{k}\ext^{2k}\Vj$ and $\oddext\Vj=\bigoplus_{k}\ext^{2k+1}\Vj$
invariant.
\end{proposition}
\proof
The natural action of $\SO_{2n}$ on $V$ induces an automorphic one on $C(V)$ given
by $\psi(v)\rightarrow\psi(Rv)$, $R\in\SO_{2n}$. Conjugating the representation of
$C(V)$ on $\ext\Vj$ by $R\in\SO_{2n}$ gives of necessity a unitary equivalent one.
It follows that there exists a unitary $\Gamma(R)$ on $\ext\Vj$ satisfying
\eqref{eq:fd covariance}. By irreducibility, $\Gamma(R)$ is unique up to a phase
so that $\Gamma(R_{1})\Gamma(R_{2})=\Gamma(R_{1}R_{2})$ and we get a
projective unitary representation of $\SO_{2n}$ on $\ext\Vj$. By uniqueness, $\Gamma$
extends the action of $U(\Vj)$.
Let $P\in U(\Vj)$ be multiplication by $-1$ on $\Vj$ so that $\evenext\Vj,\oddext\Vj$
are the $\pm 1$ eigenspaces of $\Gamma(P)$. Since $P$ commutes with $\SO_{2n}$ on
$V_{j}$, we get $\Gamma(R)\Gamma(P)\Gamma(R)^{*}=\epsilon(R)\Gamma(P)$ where
$\epsilon(R)\in\T$ has square 1. Since $\SO_{2n}$ is connected,
$\epsilon(R)\equiv\epsilon(1)$ and therefore $\evenext\Vj$ and $\oddext\Vj$ are
invariant under $\SO_{2n}$ \halmos\\

The projective representation of $\SO_{2n}$ on $\ext\Vj$ lifts to an ordinary,
unitary representation of the universal covering group $\Spin_{2n}$ of $\SO_{2n}$
which we denote by the same symbol $\Gamma$. The lift is unique since any two
differ by a character and $\Spin_{2n}$ is simple. We will show in \S
\ref{ss:infinitesimal} that $\evenext\Vj$ and $\oddext\Vj$ are  irreducible under
$\Spin_{2n}$ by using the infinitesimal action of $\so_{2n}$. For this purpose,
we give below a more convenient description of $\ext\Vj$ borrowed from
\cite[\S 12.1]{PS}.

\ssubsection{Holomorphic spinors}\label{ss:holomorphic}

The complex structure $j$ on $V$ induces a splitting $\VC=\VCH\oplus\VCA$ of
the complexification $\VC$ where the summands are the $\pm i$ eigenspaces of
$j\otimes 1$. Both are isotropic for the $\IC$--bilinear extension
of the inner product $B$ on $V$, which we denote by the same symbol, since
$B(u,v)=B(ju,jv)=-B(u,v)$ if $v\in\VCH$ or $\VCA$. They are
unitarily and anti--unitarily isomorphic to $\Vj$ by the maps
\begin{equation}\label{eq:identifications}
U_{\pm}:v\rightarrow\frac{1}{\sqrt{2}}(v\otimes 1\mp jv\otimes i)
\end{equation}
and we shall identify $\Vj$ with $\VCH$ via $U_{+}$. The latter induces an
isomorphism of CAR algebras $\AA(\Vj)\cong\AA(\VCH)$ by $c(f)\rightarrow c(U_{+}f)$
and a unitary $\ext U_{+}:\ext\Vj\rightarrow\ext\VCH$ which is equivariant for
the corresponding canonical actions of $\AA(\Vj)$ and $\AA(\VCH)$. Transporting
the action of $C(V)$ via $\ext U_{+}$, we see by \eqref{eq:fd isom} that $\psi(v)$
acts as
\begin{equation}\label{eq:fd transport}
c(U_{+}v)+c(U_{+}v)^{*}=
c(U_{+}v)+c(\overline{U_{-}v})^{*}
\end{equation}
where $\overline{u}$ is the canonical conjugation on $\VC$.
We may extend by linearity the definition of the symbols $\psi(v)$, $v\in V$ to
$\IC$--linear symbols $\psi(v)$, $v\in\VC$ satisfying
\begin{xalignat}{3}\label{eq:complex clifford}
\{\psi(u),\psi(v)\}&=2B(u,v)&
&\text{and}&
\psi(v)^{*}&=\psi(\overline{v})
\end{xalignat}
Similarly, we extend the maps $U_{\pm}$ to $\IC$--linear
maps $\VC\rightarrow\VCH,\VCA$. As is readily verified, $U_{\pm}=\sqrt{2}P_{\pm}$
where the $P_{\pm}$ are the orthogonal projections onto $\VCH$ and $\VCA$
respectively. It follows from \eqref{eq:fd transport} that the $\psi(v)$, $v\in\VC$
act on $\ext\VCH$ by $\sqrt{2}c(P_{+}v)+\sqrt{2}c(\overline{P_{-}v})^{*}$.
In particular, since the inner product $(\cdot,\cdot)$ on $\VC$ is given by
$(u,v)=B(u,\overline{v})$, we get
\begin{xalignat}{2}
\psi(v)\medspace v_{1}\wedge\cdots\wedge v_{m}&=
\sqrt{2}\medspace v\wedge v_{1}\wedge\cdots\wedge v_{m}&
&\text{if $v\in\VCH$}\label{eq:create}\\
\intertext{and}
\psi(v)\medspace v_{1}\wedge\cdots\wedge v_{m}&=\sum_{j=1}^{m}
(-1)^{j-1}\sqrt{2}\medspace B(v,v_{j})
v_{1}\wedge\cdots\wedge\wh{v_{i}}\wedge\cdots\wedge v_{m}&
&\text{if $v\in\VCA$}\label{eq:annihilate}
\end{xalignat}

Let now $R\in\SO_{2n}$. We give below an explicit formula for the action of
$\Gamma(R)$ on $\ext\VCH$. By the cyclicity of $\Omega\in\ext^{0}\VCH$ under
$C(V)$ and \eqref{eq:fd covariance}, it is sufficient to give a formula for
$\Gamma(R)\Omega$. 
Notice that the action of $R=R\otimes 1$ on $\VC$ leaves $B$ invariant and
therefore $R\VCA$ is an isotropic subspace for $B$. Thus, if $U\subseteq R\VCA$
is a subspace, the expression
\begin{equation}
\Det(U)=\psi(u_{1})\cdots\psi(u_{k})
\end{equation}
where $u_{1}\cdots u_{k}$ is an orthonormal basis of $U$ is, up to multiplication
by a complex number of modulus one, a well--defined element of $C(V)$ by
\eqref{eq:complex clifford}. Let now $U$ be a complementary subspace of
$R\VCA\cap\VCA$ in $R\VCA$ and set
\begin{equation}
\Det(R\VCA/R\VCA\cap\VCA)=\Det(U)
\end{equation}
This is well--defined only up to elements lying in the two--sided ideal
generated by the exterior algebra of $R\VCA\cap\VCA$. Since the latter
annihilates $\Omega$, $\Det(R\VCA/R\VCA\cap\VCA)\Omega$ is a well--defined
ray in $\ext\VCH$. We now have

\begin{proposition}\label{pr:explicit repr}
Let $R\in\SO_{2n}$, then
\begin{equation}\label{eq:formula}
\IC\cdot\Gamma(R)\Omega=\IC\cdot\Det(R\VCA/R\VCA\cap\VCA)\Omega
\end{equation}
\end{proposition}
\proof
The vector $\xi=\Gamma(R)\Omega$ differs from zero and, by \eqref{eq:fd covariance}
is annihilated by $\psi(R\VCA)$. Conversely, we claim that any $\eta$ satisfying
these two requirements is a multiple of $\xi$. Indeed, assuming $\|\xi\|=\|\eta\|$,
it is easy to see that the map $x\xi\rightarrow x\eta$, $x\in C(V)$ is
norm--preserving and therefore well--defined. Since it commutes with $C(V)$,
it is a scalar and therefore $\eta=\alpha\xi$, $\alpha\in\IC$.
We claim now that the right--hand side of \eqref{eq:formula} is non--zero. Let
$v_{1}\ldots v_{k}$ be a basis of a complementary subspace of $R\VCA\cap\VCA$
in $R\VCA$ so that the projections $v_{i}^{+}$ of the $v_{i}$ on $\VCH$ are linearly
independent. Then, by \eqref{eq:create}--\eqref{eq:annihilate}, the leading term of
$\Det(R\VCA/R\VCA\cap\VCA)\Omega$ is
\begin{equation}
v_{1}^{+}\wedge\cdots\wedge v_{k}^{+}\Omega
\end{equation}
and therefore does not vanish.
We check now that $\Det(R\VCA/R\VCA\cap\VCA)\Omega$ is annihilated by the elements
$\psi(R\VCA)$. This is clear for $\psi(v)$, $v\in R\VCA\cap\VCA$ since by isotropy
$\psi(v)$ commutes with the $\Det$ term and therefore annihilates $\Omega$. If, on
the other hand $v$ does not lie in $R\VCA\cap\VCA$ then we may pick a complementary
subspace of $R\VCA\cap\VCA$ containing $v$. The $\Det$ term may therefore be taken
of the form $\psi(v)\wedge\cdots$ and therefore is annihilated by $\psi(v)$. The
identity \eqref{eq:formula} now follows \halmos

\ssubsection{The infinitesimal action of $\so_{2n}$}\label{ss:infinitesimal}

We now characterise $\Gamma$ infinitesimally. The bilinear form $B$ on $\VC$ yields
an isomorphism $\so_{2n,\IC}\cong\VC\wedge\VC$ where the latter acts on $\VC$ by
$u\wedge v\thickspace\thickspace w=B(v,w)u-B(u,w)v$ with commutation relations
\begin{equation}
[u_{1}\wedge v_{1},u_{2}\wedge v_{2}]=
 (u_{1},v_{2})v_{1}\wedge u_{2}+(v_{1},u_{2})u_{1}\wedge v_{2}
-(u_{1},u_{2})v_{1}\wedge v_{2}-(v_{1},v_{2})u_{1}\wedge u_{2}
\end{equation}
Using this identification, it is easy to see that the map
$\rho:\so_{2n,\IC}\rightarrow C(V)$
\begin{equation}\label{eq:fd currents}
u\wedge v\rightarrow
\frac{1}{4}(\psi(u)\psi(v)-\psi(v)\psi(u))=
\frac{1}{2}(\psi(u)\psi(v)-B(u,v))
\end{equation}
satisfies $[\rho(X),\rho(Y)]=\rho([X,Y])$ for any $X,Y\in\so_{2n,\IC}$ and
\begin{equation}\label{eq:exponentiation relation}
[\rho(X),\psi(v)]=\psi(Xv) 
\end{equation}
for any $X\in\so_{2n,\IC}, v\in\VC$. If $X\in\so_{2n}$, the exponentiation of
\eqref{eq:exponentiation relation} in $U(\ext\VCH)$ yields
\begin{equation}
e^{\rho(X)}\psi(v)e^{-\rho(X)}=\psi(e^{X}v)
\end{equation}
and therefore, by irreducibility,
$e^{\rho(X)}=\Gamma(e^{X})$ in $PU(\ext\VCH)$ so that $\rho$ is the representation
of $\so_{2n}$ corresponding to the projective action of $\SO_{2n}$ on $\ext\VCH$.
We now use this to classify $\evenext\VCH$ and $\oddext\VCH$

\begin{proposition}\label{pr:fd classification}
Let $V_{s_{\pm}}$ be the irreducible $\Spin_{2n}$--modules with highest weights
$s_{\pm}=\half{1}(\theta_{1}+\cdots+\theta_{n-1}\pm\theta_{n})$. Then,
\begin{xalignat}{3}\label{eq:fd classification}
V_{s_{+}}&\cong\ext_{\epsilon}\VCH&
&\text{and}&
V_{s_{-}}&\cong\ext_{1-\epsilon}\VCH
\end{xalignat}
where $\epsilon=0$ or $1$ according to whether $n$ is even or odd.
\end{proposition}
\proof
Let $e_{1},\ldots,e_{2n}$ be an orthonormal basis of $V$ satisfying $j e_{2k-1}=e_{2k}$.
By \eqref{eq:identifications}, the vectors $f_{\pm k}=\sqrt{2}^{-1}(e_{2k-1}\mp ie_{2k})$,
$k=1\ldots n$ are orthonormal basis of $\VCH$ and $\VCA$ respectively and
$B(f_{k},f_{l})=\delta_{k+l,0}$.
If the maximal torus $T$ of $\SO_{2n}$ is chosen as consisting of those elements
whose matrices in the basis $e_{k}$ are of the form
\begin{xalignat}{3}\label{eq:standard torus}
&\begin{pmatrix} R(t_{1})&&0\\&\ddots&\\0&&R(t_{n})\end{pmatrix}&
&\text{with}&
R(t_{k})&=
\left(\begin{array}{rr}\cos t_{k}&-\sin t_{k}\\\sin t_{k}&\cos t_{k}\end{array}\right)
\end{xalignat}
the Lie algebra $\t$ of $T$ has a basis given by
\begin{equation}\label{eq:basis of torus}
\Theta_{k}=e_{2k}\wedge e_{2k-1}=if_{k}\wedge f_{-k}
\end{equation}
so that $\Theta_{k}$ is a matrix all of whose entries are zero except for the $k$th
diagonal block which is of the form
$\left(\begin{array}{rr}0&-1\\1&0\end{array}\right)$.
By \eqref{eq:fd currents}, $\Theta_{k}$ is represented on $\ext\VCH$ by
$\half{i}(\psi(f_{k})\psi(f_{-k})-1)$ and therefore, by
\eqref{eq:create}--\eqref{eq:annihilate},
\begin{equation}
\Theta_{k} f_{k_{1}}\wedge\cdots\wedge f_{k_{\ell}}=
\half{i}\left\{\begin{array}{rl}
 f_{k_{1}}\wedge\cdots\wedge f_{k_{\ell}}
&\text{if $k\in\{k_{1},\ldots,k_{\ell}\}$}\\[1.2em]
-f_{k_{1}}\wedge\cdots\wedge f_{k_{\ell}}
&\text{if $k\notin\{k_{1},\ldots,k_{\ell}\}$}
\end{array}\right.
\end{equation}
so that the vector $f_{k_{1}}\wedge\cdots\wedge f_{k_{\ell}}$ corresponds to the weight
$\sum_{p}\theta_{k_{p}}-\half{1}\sum_{k}\theta_{k}$ where the $\theta_{k}$ are the dual
basis elements to the $-i\Theta_{k}$. It follows that the part of $\ext\VCH$ containing
the top exterior power $\ext^{n}\VCH=\IC f_{1}\wedge\cdots\wedge f_{n}$ has highest
weight $s_{+}$ and the other has highest weight $s_{-}$ with corresponding weight vector
$f_{1}\wedge\cdots\wedge f_{n-1}$. Thus,
$V_{s_{+}}\subset\ext_{\epsilon}\VCH$ and $V_{s_{-}}\subset\ext_{1-\epsilon}\VCH$.
To conclude, notice
that the parity subspaces of $\ext\VCH$ are irreducible under the even part of the
Clifford algebra and therefore $\so_{2n}$ since the former is generated by the
$\psi(u)\psi(v)$ \halmos\\

\remark By the tensor product rules of proposition
\ref{ch:classification}.\ref{pr:tensor with minimal},
$\Hom_{\Spin_{2n}}(\VC\otimes V_{s_{\pm}},V_{s_{\mp}})\cong\IC$. The
corresponding intertwiners may be constructed using the Clifford multiplication map
\begin{xalignat}{2}\label{eq:vector intertwiner}
\VC\otimes\ext\VCH&\rightarrow\ext\VCH&
v\otimes w&\rightarrow\psi(v)w
\end{xalignat}
which, by \eqref{eq:exponentiation relation} commutes with $\Spin_{2n}$. Its restrictions
to $\VC\otimes\evenext\VCH$ and $\VC\otimes\oddext\VCH$ are the required intertwiners.

\ssection{The infinite--dimensional spin representations}\label{se:quark model}

In this section, we construct the level 1 representations of $L\Spin_{2n}$ by using
an analogue of the Clifford algebra construction of section \ref{se:fd spinors} and
realise them, grouped in pairs as summands of two distinct exterior algebras or
Fermionic Fock spaces, the \NS and Ramond sectors $\FR,\FNS$.\\

We begin by discussing the representations of the Clifford algebra $C(\H)$ of a real,
infinite--dimensional Hilbert space. As in section \ref{se:fd spinors}, these may be
obtained by introducing an orthogonal complex structure $J$ on $\H$ and regarding
the latter as a complex Hilbert space $\H_{J}$. An irreducible representation of
$C(\H)$ is then got via the canonical action of the isomorphic CAR algebra $\AA(\HJ)$
on the Fock space $\ext\HJ$.
Unlike the finite--dimensional case however, different complex structures $J_{1},J_{2}$
lead in general to inequivalent representations and we give below a necessary and
sufficient criterion due to I. Segal \cite{BSZ} for that to be the case. The discussion
is technically simpler when $J_{1}$ and $J_{2}$ commute with a given complex structure
$i$ on $\H$ and we first formulate the criterion as the equivalence of representations
of the reference CAR algebra $\AA(\H_{i})$. This is the context of {\it complex fermions}
and is treated in \S \ref{ss:CAR}. The more general case of {\it real fermions} is studied
in \S \ref{ss:clifford}.\\

The complex Fermionic criterion leads at once to a projective representation of a
distinguished subgroup of the unitary group of $\H_{i}$ on the Fock space $\ext\HJ$.
In \S \ref{ss:LU}, we derive from it the {\it basic representation} of the loop group
of $U_{n}$ starting from its unitary action on $\ltwo{\IC^{n}}$, as done in \cite{Wa3}.
Similarly, the real Fermionic criterion yields a projective representation of a
subgroup of the orthogonal group of $\H$ extending the previous one. When applied
to the orthogonal action of $L\SO_{2n}$ on the Hilbert spaces of periodic and
anti--periodic functions on $S^{1}$ with values in $\IR^{2n}$, this yields positive
energy representations of $L\Spin_{2n}$ and is carried out in \S \ref{ss:R and NS}.
These are classified in \S \ref{ss:level 1} and shown to contain all level 1
representations.

\ssubsection{Complex fermions and CAR algebras}\label{ss:CAR}

Consider a {\it complex}, infinite--dimensional Hilbert space $\H$ with hermitian form
$(\cdot,\cdot)$ and the corresponding CAR algebra $\AA(\H)$ defined as in section
\ref{se:fd spinors} by the $\IC$--linear symbols $c(f)$, $f\in\H$ subject to the
anticomutation relations
\begin{xalignat}{2}\label{eq:CAR 2}
\{c(f),c(g)\}&=0&
\{c(f),c(g)^{*}\}&=(f,g)
\end{xalignat}
$\AA(\H)$ possesses a canonical representation on the exterior algebra or Fermionic Fock
space $\ext\H$ given by
\begin{align}
c(f)\thickspace g_{1}\wedge\cdots\wedge g_{m}&=f\wedge g_{1}\wedge\cdots\wedge g_{m}
\label{eq:repr 2-1}\\
\intertext{so that}
c(f)^{*}\thickspace g_{1}\wedge\cdots\wedge g_{m}&=\sum_{j=1}^{m}
(-1)^{j-1}(g_{j},f)g_{1}\wedge\cdots\wedge\wh{g_{i}}\wedge\cdots\wedge g_{m}
\label{eq:repr 2-2}
\end{align}
which is irreducible \cite[cor. 3.3.1]{BSZ}.
We also consider the underlying real Hilbert space $\HIR$ with inner product
$B(\cdot,\cdot)=\Re(\cdot,\cdot)$ and the associated Clifford algebra $C(\HIR)$
generated by the real--linear, self--adjoint symbols $\psi(f)$, $f\in\HIR$ satisfying
\begin{equation}\label{eq:clifford 2}
\{\psi(f),\psi(g)\}=2B(f,g)
\end{equation}

The maps
\begin{xalignat}{2}\label{eq:isom 2}
\psi(f)&\rightarrow c(f)+c(f)^{*}&
   c(f)&\rightarrow \half{1}\Bigl(\psi(f)-i\psi(if)\Bigr)
\end{xalignat}
extend to an isomorphism $C(\HIR)\cong\AA(\H)$ which we shall however use in a somewhat
opposite way to that of section \ref{se:fd spinors}. Indeed, choosing an orthogonal
complex structure $J$ on $\H_{\IR}$ differing from the original one and denoting by $\HJ$
the corresponding complex Hilbert space, we get $\AA(\H)\cong C(\H_{\IR})\cong\AA(\HJ)$
and therefore a representation $\pi_{J}$ of $\AA(\H)$ via the canonical action of
$\AA(\HJ)$ on the Fock space $\FJ=\ext\HJ$.\\

The equivalence class of $\pi_{J}$ depends upon $J$ in a way explained below.
The discussion is simpler
if $J$ commutes with the original complex structure on $\H$ so that it is unitary.
We may then diagonalise it and write $\H=\H_{+}\oplus\H_{-}$ where the $\H_{\pm}$
are the $\pm i$ eigenspaces of $J$.
Then, $\H_{J}=\H_{+}\oplus\overline{\H_{-}}$ where
$\overline{\H_{-}}\cong\H_{-}^{*}$ is $\H_{-}$ with a reversed complex structure and
$\FJ=\ext\H_{+}\wh\otimes\ext\H_{-}^{*}$.
The corresponding action of $\AA(\H)$ on $\FJ$ is explicitly given by
\begin{equation}\label{eq:pi J}
\pi_{J}(c(f))=c(f_{+})+c(f_{-})^{*}
\end{equation}
where $f_{\pm}$ are the projections of $f$ on $\H_{\pm}$.\\

In marked contrast to the finite--dimensional case, the automorphic action of the full
unitary group of $\H$ is not implemented on $\FJ$ in the sense that, for a given
$u\in U(\H)$ there does not exist in general a unitary operator $\Gamma(u)$ on $\FJ$
such that $\Gamma(u)\pi_{J}(c(f))\Gamma(u)^{*}=\pi_{J}(c(uf))$. Notice though that
the operator $\Gamma(u)$, when it exists is unique up to a phase since $\pi_{J}$ is
irreducible. The existence of $\Gamma(u)$ is in fact easily seen to be equivalent
to the unitary equivalence of $\pi_{J}$ and $\pi_{\wt J}$ where $\wt J=uJu^{*}$.
The equivalence criterion given below determines a projective representation of a
distinguished subgroup of $U(\H)$ which we presently define. Consider the
{\it restricted unitary group} of $\H$, relative to the unitary complex structure
$J$ given by
\begin{equation}\label{eq:ures}
\ures(\H,J)=\{u\in U(\H)|\thinspace \|[u,J]\|_{\HS}<\infty\}
\end{equation}
where $\|\cdot\|_{\HS}$ is the Hilbert--Schmidt norm. We topologise $\ures(\H,J)$
by endowing it with the strong operator topology combined with the metric topology
determined by $d(u,v)=\|[u-v,J]\|_{\HS}$. It contains the intersection of the unitary
groups of $\H$ and $\HJ$. The automorphic action of the latter is implemented on
$\FJ$ by
\begin{equation}
\Gamma(u)g_{1}\wedge\cdots\wedge g_{m}=ug_{1}\wedge\cdots\wedge ug_{m}
\end{equation}
we refer to this as the {\it canonical quantisation} of $U(\H)\cap U(\HJ)$.
The following result is due to I. Segal \cite{BSZ}. A short proof as well as the details
of (iii), may be found in \cite{Wa3}.

\begin{theorem}\label{th:complex segal}
Let $\H$ be a complex Hilbert space and $\AA(\H)$ the corresponding CAR algebra. Then
\begin{enumerate}
\item Any unitary complex structure $J$ on $\H$ determines by \eqref{eq:pi J} an irreducible
representation $\pi_{J}$ of $\AA(\H)$ on $\FJ=\ext\H_{+}\wh\otimes\ext\overline\H_{-}$ where
the $\H_{\pm}$ are the $\pm i$ eigenspace of $J$.
\item The representations $\pi_{J_{1}}$ and $\pi_{J_{2}}$ are unitarily equivalent if,
and only if $J_{1}-J_{2}$ is a Hilbert--Schmidt operator.
\item For any unitary complex structure $J$, there exists a strongly continuous, projective
unitary representation $\Gamma$ of the restricted unitary group \eqref{eq:ures} on $\FJ$
satisfying and uniquely determined by
\begin{equation}\label{eq:complex implementation}
\Gamma(u)\pi_{J}(c(f))\Gamma(u)^{*}=\pi_{J}(c(uf))
\end{equation}
\end{enumerate}
\end{theorem}

\ssubsection{The basic representation of $LU_{n}$}\label{ss:LU}

Let $LU_{n}=C^{\infty}(S^{1},U_{n})$ and consider the action of $LU_{n}\rtimes\rot$
on the Hilbert space $\H=\ltwo{\IC^{n}}$ where $LU_{n}$ acts by multiplication and
$\rot$ by $R_{\theta}f=f_{\theta}$ where $f_{\theta}(\phi)=f(\phi-\theta)$. We 
denote the standard complex structure on $\H$ by $i$. Thus, $R_{\theta}=e^{i\theta d}$
where the self--adjoint generator of rotations $d=i\frac{d}{d\theta}$ is non--negative
on the subspace $\H_{+}$ spanned by the functions $e^{im\theta}\otimes v$, $v\in\IC^{n}$,
$m\leq 0$ and negative on $\H_{-}=\H_{+}^{\perp}$.
Let $J$ be the complex structure acting as $\pm i$ on $\H_{\pm}$. Up to a
finite--dimensional ambiguity on the space of constant functions $\IC^{n}\subset\H$,
$J$ is the unique complex structure on $\H$ such that $R_{\theta}$ is a unitary action
on $\H_{J}$ decomposing into a sum of non--negative characters. In other words, $J$
commutes with $R_{\theta}$ and the corresponding self--adjoint generator on $\H_{J}$,
namely $J\frac{d}{d\theta}$ is non--negative.\\

Consider the representation $\pi_{J}$ of the CAR algebra $\AA(\H)$ on
$\F=\ext\HJ=\ext\H_{+}\wh\otimes\ext\overline{\H_{-}}$ given as in \eqref{eq:pi J}.
$R_{\theta}$ commutes with $J$ and is therefore canonically quantised on $\F$, call
the corresponding action $U_{\theta}=\Gamma(R_{\theta})$.
Since $R_{\theta}$ has non--negative spectrum and finite--dimensional
eigenspaces on $\HJ$, $U_{\theta}$ is a positive energy representation of $\rot$.
Let $\gamma\in LU_{n}$ with Fourier expansion
$\sum_{m\in\IZ}\wh\gamma(m)e^{im\theta}$. It is readily verified, using the basis 
$e^{im\theta}\otimes v$ of $\H$ that
\begin{equation}
\|[\gamma,J]\|_{\HS}^{2}=4\sum_{m}|m|\|\wh\gamma(m)\|^{2}
\end{equation}
where $\|\wh\gamma(m)\|$ is the Hilbert--Schmidt norm on $\End(\IC^{n})$. On the other
hand, the operator norm of $\gamma$ is equal to the supremum norm $\|\gamma\|_{\infty}$
and by theorem \ref{th:complex segal} we get a continuous projective unitary
action of $LU_{n}$ of positive energy on $\F$ called the {\it basic representation}
of $LU_{n}$.

\ssubsection{Real fermions and Clifford algebras}\label{ss:clifford}

Let now $\H$ be a {\it real} Hilbert space with inner product $B(\cdot,\cdot)$ and
Clifford algebra $C(\H)$. Every orthogonal complex structure $J$ on $\H$ gives rise
to a representation $\pi_{J}$ of $C(\H)$ on the exterior algebra $\ext\HJ$ via the
isomorphism $C(\H)\cong\AA(\HJ)$. Explicitly,
\begin{equation}
\pi_{J}(\psi(f))=c(f)+c(f)^{*}
\end{equation}
where $c(f)$ and $c(f)^{*}$ act as in \eqref{eq:repr 2-1}--\eqref{eq:repr 2-2}.
We give below a criterion for the representations of $C(\H)$ determined by two
orthogonal complex structures $J_{1},J_{2}$ to be unitary equivalent. 
Its proof proceeds by noting that, up to the replacement of $\H$ by
$\H\oplus\H=\H\otimes\IC$ and $J_{k}$ by $J_{k}\otimes 1$, one may always assume
that $J_{1},J_{2}$ commute with a reference complex structure on $\H$. The criterion
is then obtained from that of theorem \ref{th:complex segal} in view of the following
simple observation 

\begin{lemma}\label{le:doubling}
Let $\H$ be a real Hilbert space with Clifford algebra $C(\H)$ and $J_{1},J_{2}$ two
orthogonal complex structures.
Then, the representations of $C(\H)$ on $\ext\H_{J_{1}}$, $\ext\H_{J_{2}}$
are unitarily equivalent if, and only if those of $C(\H\oplus\H)$ on
$\ext(\H_{J_{1}}\oplus\H_{J_{1}})$,
$\ext(\H_{J_{2}}\oplus\H_{J_{2}})$ are.
\end{lemma}
\proof We proceed as in \cite[lemma 1]{SS}. Let $V:\ext\H_{J_{1}}\rightarrow\ext\H_{J_{2}}$
be a unitary $C(\H)$--intertwiner. Then,
\begin{equation}
 V\otimes V:
 \ext(\H_{J_{1}}\oplus\H_{J_{1}})\cong\ext\H_{J_{1}}\wh\otimes\ext\H_{J_{1}}
 \longrightarrow
 \ext\H_{J_{2}}\wh\otimes\ext\H_{J_{2}}\cong\ext(\H_{J_{2}}\oplus\H_{J_{2}})
\end{equation}
intertwines the action of $C(\H\oplus\H)=C(\H)\wh\otimes C(\H)$. Conversely, let
$U:\ext\H_{J_{1}}\wh\otimes\ext\H_{J_{1}}
   \rightarrow
   \ext\H_{J_{2}}\wh\otimes\ext\H_{J_{2}}$
be a $C(\H)\wh\otimes C(\H)$--intertwiner and
$\Omega_{1}\otimes\Omega_{1}\in\ext\H_{J_{1}}\wh\otimes\ext\H_{J_{1}}$ be the vacuum
vector. We claim that $U(\Omega_{1}\otimes\Omega_{1})=\xi\otimes\eta$ for some
$\xi,\eta\in\ext\H_{J_{2}}$ of norm one. Assuming this
for the moment, one sees that the map $V:\ext\H_{J_{1}}\rightarrow\ext\H_{J_{2}}$
given by $\pi_{J_{1}}(x)\Omega_{1}\rightarrow\pi_{J_{2}}(x)\xi$ is well--defined and
norm preserving since, for any $x\in C(\H)$
\begin{equation}
\|\pi_{J_{2}}(x)\xi\|=
\|\pi_{J_{2}}(x)\otimes 1(\xi\otimes\eta)\|=
\|U(\pi_{J_{1}}(x)\otimes 1)(\Omega_{1}\otimes\Omega_{1})\|=
\|\pi_{J_{1}}(x)\Omega_{1}\|
\end{equation}
and therefore yields a unitary map intertwining $C(\H)$. Returning to our claim, let
$U(\Omega_{1}\otimes\Omega_{1})=\sum_{k}\xi_{k}\otimes\eta_{k}$ where
$(\eta_{k},\eta_{k'})=\delta_{k,k'}$. Then,
\begin{equation}
\begin{split}
(\pi_{J_{1}}(x)\Omega_{1},\Omega_{1})
&=(\pi_{J_{1}}(x)\otimes 1(\Omega_{1}\otimes\Omega_{1}),
   \Omega_{1}\otimes\Omega_{1})\\
&=(U\pi_{J_{1}}(x)\otimes 1U^{*}U(\Omega_{1}\otimes\Omega_{1}),
   U(\Omega_{1}\otimes\Omega_{1}))\\
&=\sum_{k}(\pi_{J_{2}}(x)\xi_{k},\xi_{k})
\end{split}
\end{equation}
However, $\pi_{J_{1}}$ is the GNS representation of $C(\H)\cong\AA(\H_{J_{1}})$ for the
vector state $\phi(x)=(\pi_{J_{1}}(x)\Omega_{1},\Omega_{1})$. By theorem
\ref{th:complex segal}, $\pi_{J_{1}}$ is irreducible so that $\phi$ is pure.
It follows that all non--zero states $(\pi_{J_{2}}(x)\xi_{k},\xi_{k})$
must be equal. By irreducibility of $\pi_{J_{2}}$, this is the case if, and only if all
$\xi_{k}$ are proportional and therefore $U\Omega_{1}\otimes\Omega_{1}=\xi\otimes\eta$ 
as claimed \halmos\\

It follows almost at once that $\pi_{J_{1}}$ is unitarily equivalent to $\pi_{J_{2}}$
if, and only if $J_{1}-J_{2}$ is a (real) Hilbert--Schmidt operator on $\H$. As in
theorem \ref{th:complex segal}, this criterion leads to a projective representation
of a distinguished subgroup of the orthogoanl group $O(\H)$. More precisely, consider
the {\it restricted orthogonal group} of $\H$ relative to $J$ defined by
\begin{equation}\label{eq:ores}
\ores(\H,J)=\{R\in O(\H)|\thinspace\|[R,J]\|_{\HS}<\infty\}
\end{equation}
where the Hilbert--Schmidt norm refers to $[R,J]$ as a real operator. We topologise
$\ores(\H,J)$ by the strong operator topology combined with the metric $d(R,S)=\|[R-S,J]\|_{HS}$.
This gives, for any other reference complex structure $i$ on $\H$ commuting with $J$ a
continuous inclusion $\ures(\H_{i},J)\subset\ores(\H,J)$. The following theorem is due
to Shale and Stinespring \cite{BSZ,SS}

\begin{theorem}\label{th:real segal}
Let $\H$ be a real Hilbert space with Clifford algebra $C(\H)$. Then
\begin{enumerate}
\item Any orthogonal complex structure $J$ on $\H$ determines an irreducible
representation $\pi_{J}$ of $C(\H)$ on $\FJ=\ext\H_{J}$ given by
\begin{equation}
\pi_{J}(\psi(f))=c(f)+c(f)^{*}
\end{equation}
\item The representations $\pi_{J_{1}}$ and $\pi_{J_{2}}$ are unitarily equivalent
if, and   only if $J_{1}-J_{2}$ is a Hilbert--Schmidt operator.
\item For any unitary $J$, there exists a strongly continuous, projective unitary
representation $\Gamma$ of the restricted orthogonal group \eqref{eq:ores} on $\FJ$
satisfying and uniquely determined by
\begin{equation}\label{eq:real implementation}
\Gamma(R)\pi_{J}(\psi(f))\Gamma(R)^{*}=\pi_{J}(\psi(Rf))
\end{equation}
\item If $J$ commutes with a reference complex structure $i$ on $\H$, $\Gamma$ extends the
representation of $\ures(\H_{i},J)\subset\ores(\H,J)$ given by theorem \ref{th:complex segal}.
\end{enumerate}
\end{theorem}
\proof
(i) follows at once from the isomorphism $C(\H)\cong\AA(\H_{J})$ given by
\eqref{eq:isom 2} and (i) of theorem \ref{th:complex segal}.

(ii) By lemma \ref{le:doubling}, $\ext\H_{J_{1}}$ and $\ext\H_{J_{2}}$ are unitarily equivalent
as $C(\H)$--modules iff $\ext(\H_{J_{k}}\oplus\H_{J_{k}})$, $k=1,2$ are unitarily equivalent
as $C(\H\oplus\H)$--modules.
On the other hand, $\H\oplus\H=\H\otimes_{\IR}\IC$ possesses a natural complex
structure commuting with $J_{k}\otimes 1$ and since $C(\H\oplus\H)\cong\AA(\H\otimes_{\IR}\IC)$,
theorem \ref{th:complex segal} implies that the original representations are equivalent
if, and only $J_{1}\otimes 1-J_{2}\otimes 1$ is a complex Hilbert--Schmidt operator on
$\H\otimes_{\IR}\IC$ and therefore if, and only if $J_{1}-J_{2}$ is a real Hilbert--Schmidt
operator on $\H$.

(iii) Let $R\in O(\H)$ and regard it as a unitary operator $\HJ\rightarrow\H_{\wt J}$
where $\wt J=RJR^{-1}$. As such it extends to a unitary $\ext R:\ext\HJ\rightarrow\ext\H_{\wt J}$
satisfying $\ext R\pi_{J}(c(f)){\ext R}^{*}=\pi_{\wt J}(c(Rf))$. Notice that $[R,J]$ is
Hilbert--Schmidt
if, and only if $J-\wt J$ is. Thus, if $R\in\ores(\H,J)$, there exists by (i) a unitary
$V:\ext\HJ\rightarrow\ext\H_{\wt J}$ intertwining $C(\H)$. The operator
$\Gamma(R)=V^{*}RV$ is easily seen to have all the required properties. Moreover, by
irreducibility, $\Gamma(R)$ is uniquely determined, up to a phase by \eqref{eq:real implementation}
and in particular, $\Gamma(R_{1}R_{2})=\Gamma(R_{1})\Gamma(R_{2})$. As is readily verified,
$\Gamma$ is such that the following diagram commutes
\newarrow{Ddot}.....
\begin{equation}
\begin{diagram}[height=2.2em,width=3.0em]
\ores(\H,J)&\rTo&\ures(\H\otimes_{\IR}\IC,J\otimes 1)&\rTo{\Gamma_{\otimes}}
								&PU(\ext\HJ\wh\otimes\ext\HJ)\\
	   &\rdTo{\Gamma}&			     &\ruTo{\Delta}&\\				   
	   &	&		PU(\ext\HJ)	     &		   &\\				   
\end{diagram}
\end{equation}
where $\Gamma_{\otimes}$ is the representation of $\ures(\H\otimes_{\IR}\IC,J\otimes 1)$
given by theorem \ref{th:complex segal} and $\Delta(u)=u\otimes u$. It follows that
$\Gamma$ is continuous since $\Delta$ is a homeomorphism onto its image.

(iv) is clear in view of the uniqueness of the operators $\Gamma(R)$,
\eqref{eq:complex implementation} and the fact that the isomorphism $C(\H)\cong\AA(\H_{i})$
is equivariant for the automorphic actions of $U(\H_{i})\subset O(\H)$ \halmos

\ssubsection{The Ramond and \NS Fock spaces}\label{ss:R and NS}

We construct below positive energy representations of $L\Spin_{2n}$ by mapping it to
the restricted orthogonal groups of suitable real Hilbert spaces and using theorem
\ref{th:real segal}. These representations factor through the loop group of $\SO_{2n}$
and extend the basic representation of $LU_{n}\subset L\SO_{2n}$ obtained in \S
\ref{ss:LU}.\\

Consider the Ramond and \NS Hilbert spaces $\HR,\HNS$ of square integrable,
$\IC^{n}$--valued periodic and anti--periodic functions on $S^{1}$ spanned by the
functions $e^{im\theta}\otimes v$ with $m\in\IZ$ or $m\in\half{1}+\IZ$ respectively.
They are naturally complex Hilbert spaces with standard complex structure $i$ but
we shall regard them as real ones. As such, they support orthogonal actions of
$L\SO_{2n}\rtimes\rot$ where the
first factor acts by multiplication on $\IC^{n}\cong\IR^{2n}$--valued functions and
the second by $R_{\theta}f=f_{\theta}$
\footnote{The action of $\rot$ on $\HNS$ is two--valued and is therefore really one
of the non--trivial double cover of $\rot$.}.
Let $\H$ be $\HR$ or $\HNS$ and split is as a direct sum $\H_{+}\oplus\H_{-}$ where
the $\H_{\pm}$ are the (complex) subspaces spanned by the functions
$e^{im\theta}\otimes v$,
$m\leq 0$ and $m>0$ respectively with $m\in\IZ$ for $\HR$ and $m\in\half{1}+\IZ$ for
$\HNS$. Let $J$ be the complex structure acting as multiplication by $\pm i$ on
$\H_{\pm}$. If $\H$ is endowed with the complex structure $J$, the action of $\rot$
is unitary and has non--negative spectrum and finite--dimensional eigenspaces. Its
canonical quantisation on $\F=\ext\H_{J}$ is therefore of positive energy. Let
$\ores(\H,J)$ be defined by \eqref{eq:ores}, then

\begin{lemma}\label{le:HS cond}
Let $\H=\HR$ or $\HNS$ and $J$ the complex structure defined above. Then,
$L\SO_{2n}\subset\ores(\H,J)$ and the inclusion is continuous.
\end{lemma}
\proof
Let $\gamma\in L\SO_{2n}$. To compute $\|[\gamma,J]\|_{\HS}$, it is more convenient
to complexify $\H$ and consider the splitting
\begin{equation}
\HC=\HC^{1,0}\oplus\HC^{0,1}
\end{equation}
where $\HCH$, $\HCA$ are the $\pm 1\otimes i$ eigenspaces of $J\otimes 1$ on $\HC$.
If $\H=\HNS$, these are spanned by the functions $v(m)=e^{(1\otimes i)m\theta}v$,
$v\in\VC$ and $m\lessgtr 0$ respectively. For $\H=\HR$, $\HCH$ is spanned by the
$v(m)$ with $m<0$ or $m=0$ and $v\in\VCH$, the $1\otimes i$ eigenspace of $i\otimes 1$
on $\VC$ while $\HCA$ is spanned by $v(m)$ with $m>0$ or $m=0$ and $v\in\VCA$, the
$-1\otimes i$ eigenspace of $i\otimes 1$ on $\VC$. Using these basis, a simple
computation shows that for $\gamma\in L\SO_{2n}$,
\begin{equation}
\|[\gamma\otimes 1,J\otimes 1]\|^{2}_{\HS(\HC)}\leq
 4\sum_{m\in\IZ}(|m|+1)\|\wh\gamma(m)\|^{2}
\end{equation}
where the $\wh\gamma(m)\in\End(\VC)$ are the Fourier coefficients of $\gamma$ and
$\|\wh\gamma(m)\|$ is the Hilbert--Schmidt norm on $\End(\VC)$. The lemma follows
since $\|[\gamma,J]\|_{\HS(\H)}=\|[\gamma\otimes 1,J\otimes 1]\|_{\HS(\HC)}$ \halmos\\

By lemma \ref{le:HS cond} and theorem \ref{th:real segal}, we get a continuous projective
representation of $L\SO_{2n}$ on the Fock space $\F=\ext\HJ$ with $\H=\HR$ or $\HNS$ which,
as remarked above is of positive energy.
Notice that by identifying $U_{n}$ with the subgroup of $\SO_{2n}$ commuting
with the complex structure on $\IR^{2n}\cong\IC^{n}$, we get a natural inclusion
$LU_{n}\subset L\SO_{2n}$. The restriction of $\HR$ to $LU_{n}$ is the standard
representation of $LU_{n}$ on $\ltwo{\IC^{n}}$ considered in \S \ref{ss:LU},
and the restriction of $\HNS$ to $LU_{n}$ is unitary equivalent to $\ltwo{\IC^{n}}$
by the map $f\rightarrow e^{-i\half{\theta}}f$.
Since these identifications carry the complex structure $J$ on $\HR$ or $\HNS$
considered above to the one on $\ltwo{\IC^{n}}$ given in \S \ref{ss:LU}, the
restrictions of $\FR,\FNS$ to $LU_{n}$ are equivalent to its basic representation.

\begin{lemma}\label{le:setup}
Let $\H=\HR$ or $\HNS$ and $\F_{0}=\bigoplus_{k\geq 0}\ext^{2k}\HJ$ and
$\F_{1}=\bigoplus_{k\geq 0}\ext^{2k+1}\HJ$ the even and odd subspaces of $\F$.
Then, each is invariant under $L\Spin_{2n}\rtimes\rot$
and therefore defines a positive energy representation of $L\Spin_{2n}$.
\end{lemma}
\proof
Let $P$ be multiplication by $-1$ on $\H$. $P$ is canonically quantised on $\F=\ext\HJ$
and $\F_{0},\F_{1}$ are the $\pm 1$ eigenspaces of $\Gamma(P)$.
The action of $\rot$ on $\H$ commutes with $P$ and $J$. It is canonically quantised on
$\F$ and commutes with $\Gamma(P)$ and therefore leaves each $\F_{k}$ invariant.
$L\Spin_{2n}$ commutes with $P$ on $\H$ hence
\begin{equation}
\Gamma(\gamma)\Gamma(P)\Gamma(\gamma)^{*}=
\epsilon(\gamma)\Gamma(P)
\end{equation}
for some $\epsilon(\gamma)\in\T$ depending continuously on $\gamma$. Since
$L\Spin_{2n}$ is connected and $\epsilon(\gamma)^{2}=1$, we have
$\epsilon(\gamma)=\epsilon(1)=1$ so that the $\F_{k}$ are invariant under $L\Spin_{2n}$.
We now distinguish two cases. If $\H=\HR$, the action of $\rot$ on $\H$ is a genuine
one and we therefore get a positive energy representation of $L\Spin_{2n}\rtimes\rot$ on
$\F_{0},\F_{1}$.
On the other hand, the action of $\rot$ on $\H=\HNS$ is really one of its non--trivial
double cover $\wt{\rot}$. This leads to a positive energy representation
of $\wt{\rot}$ on $\F$ with characters $z\rightarrow z^{n}$, $n\in\IN$ on $\F^{0}$
and $z^{n}$, $n\in\half{1}+\IN$ on $\F^{1}$. Thus, multiplying the action on $\F^{1}$
by the character $z^{-\half{1}}$ we get a positive energy representation of $\rot$
on $\F_{0},\F_{1}$ and therefore one of $L\Spin_{2n}\rtimes\rot$ \halmos

\ssubsection{The level 1 representations of $L\Spin_{2n}$}\label{ss:level 1 of spin}

We prove below that the level of $\FR$ and $\FNS$ as positive energy representations
of $L\Spin_{2n}$ is equal to one by computing the commutator cocycle for the action
of the torus and the coroot lattice of $\Spin_{2n}$ and using the criterion of
corollary \ref{ch:classification}.\ref{co:global level}. An alternative infinitesimal
proof will be given in section \ref{se:currents}. We begin by giving an
equivalent description of the Fock spaces $\FR$ and $\FNS$ parallel to that of the
finite--dimensional spin modules obtained in \S \ref{ss:holomorphic}.\\

Let $\H=\HR$ or $\HNS$ and $J$ the complex structure of \S \ref{ss:R and NS}. Split
the complexification $\HC$ as
\begin{equation}\label{eq:decomposition}
\HC=\HCH\oplus\HCA
\end{equation}
where the summands are the $\pm 1\otimes i$ eigenspaces of $J\otimes 1$. 
As in the finite--dimensional case, they are isotropic for the $\IC$--bilinear
extension of the inner product $B$ on $\H$ and are unitarily and anti--unitarily
isomorphic to $\HJ$ via the maps
\begin{equation}
U_{\pm}:f\rightarrow\frac{f\otimes 1-Jf\otimes i}{\sqrt{2}}
\end{equation}
If $\H=\HNS$, $\HCH$ and $\HCA$ are spanned by the functions
$v(m)=e^{(1\otimes i)m\theta}\otimes v$ where $v\in\VC$ and $m\in\half{1}+\IZ$ is
$\lessgtr0$ respectively. If, on the other hand $\H=\HR$, $\HCH$ is spanned by the
$v(m)$, with  $m\in\IN_{-}$ and $v\in\VC$ or $m=0$ and $v\in\VCH$. The latter is a
summand of the splitting $\VC=\VCH\oplus\VCA$ considered in \S \ref{ss:holomorphic}
and corresponds to the $1\otimes i$--eigenspace of $i\otimes 1$ on $\VC$. Similarly,
$\HCA$ is spanned by the $v(m)$ with $m\in\IN_{+}$ and $v\in\VC$ or $m=0$ and $v\in\VCA$.\\

Using $U_{+}$, we identify $\HJ$ with $\HCH$ and $\F=\ext\HJ$ with $\ext\HCH$. Transporting
the representation of the Clifford algebra $C(\H)$ to the latter Fock space, we find as in
\S \ref{ss:holomorphic} that it acts as 
\begin{xalignat}{2}
\psi(f)\medspace g_{1}\wedge\cdots\wedge g_{m}&=
\sqrt{2}\medspace f\wedge g_{1}\wedge\cdots\wedge g_{m}&
&\text{if $f\in\HCH$}\label{eq:complex create}\\
\intertext{and}
\psi(f)\medspace g_{1}\wedge\cdots\wedge g_{m}&=\sum_{j=1}^{m}
(-1)^{j-1}\sqrt{2}\medspace B(f,g_{j})
g_{1}\wedge\cdots\wedge\wh{g_{i}}\wedge\cdots\wedge g_{m}&
&\text{if $f\in\HCA$}\label{eq:complex annihilate}
\end{xalignat}
where we have extended the definition of the symbols $\psi(\cdot)$ to $\IC$--linear
symbols $\psi(f)$, $f\in\HC$ satisfying
\begin{xalignat}{3}\label{eq:new anticommutation}
\{\psi(f),\psi(g)\}&=2B(f,g)&
&\text{and}&
\psi(f)^{*}&=\psi(\overline{f})
\end{xalignat}

We now characterise explicitly the action of the coroot lattice of $\Spin_{2n}$ on $\F$.
Let $e_{1}\cdots e_{2n}$ be an orthonormal basis of $V$ such that $ie_{2k-1}=e_{2k}$.
Let $T$ be the torus of $\SO_{2n}$ consisting
of block diagonal matrices of the form \eqref{eq:standard torus} and $\wt T$ its pre--image
in $\Spin_{2n}$.
The coroot lattice $\coroot\cong\Hom(\T,\wt T)\subset L\wt T$ of $\Spin_{2n}$ embeds in the
integral lattice $I=\Hom(\T,T)$ of $\SO_{2n}$. The latter has a $\IZ$--basis given, in the
notation of \eqref{eq:basis of torus} by $i\Theta_{k}=i e_{2k}\wedge e_{2k-1}$ where to
each $\lambda\in\bigoplus_{k}i\Theta_{k}\cdot\IZ$, we associate the loop
$\zeta_{\lambda}(\phi)=\exp_{T}(-i\phi\lambda)$. 
Let
$f_{\pm k}=\sqrt{2}^{-1}(e_{2k-1}\otimes 1\mp e_{2k}\otimes i)$,
$k=1\ldots n$ be the orthonormal basis of $\VCH,\VCA$ corresponding to
$e_{1}\cdots e_{2n}$ via \eqref{eq:identifications} so that
$B(f_{k},f_{l})=\delta_{k+l,0}$. Then,

\begin{lemma}\label{le:explicit coroot}
In $\FR$, we have, projectively,
\begin{equation}
U(\zeta_{i\Theta_{k}})\Omega=
\frac{1}{\sqrt{2}}\psi(f_{-k}(-1))\Omega=
f_{-k}(-1)
\end{equation}
\end{lemma}
\proof
The proof of proposition \ref{pr:explicit repr} shows that if $R\in\ores(\H,J)$
is such that $R\HCA\cap\HCA$ is of finite codimension in $R\HCA$, then
$\Gamma(R)\Omega$ is given, as in \eqref{eq:formula} by
$\Det(R\HCA/R\HCA\cap\HCA)\Omega$. We choose now $R=\zeta_{i\Theta_{k}}$ and
compute $R\HCH$. From $\Theta_{k}=e_{2k}\wedge e_{2k-1}=if_{k}\wedge f_{-k}$,
we find $\Theta_{k}f_{\pm l}=\pm i\delta_{kl} f_{\pm l}$ so that
\begin{equation}\label{eq:T on fk}
\exp(\sum t_{k}\Theta_{k})f_{\pm l}=e^{\pm it_{l}}f_{\pm l}
\end{equation}
and therefore
\begin{equation}
\zeta_{i\Theta_{k}}f_{\pm l}(n)\medspace(\phi)=
e^{(1\otimes i)(n\pm\delta_{kl})\phi}\otimes f_{\pm l}=
f_{\pm l}(n\pm\delta_{kl})\medspace(\phi)
\end{equation}
Since
\begin{equation}
\HCA=\bigoplus_{m>0,l}\IC f_{l}(m)\oplus\IC f_{-l}(m)
     \bigoplus_{l}\IC f_{-l}(0)
\end{equation}
we find
\begin{equation}
\zeta_{i\Theta_{k}}\HCA=
\bigoplus_{m> 0,l}\IC f_{l}(m+\delta_{kl})\oplus\IC f_{-l}(m-\delta_{kl})
\bigoplus_{l}\IC f_{-l}(-\delta_{kl})
\end{equation}
so that
\begin{equation}
\zeta_{i\Theta_{k}}\HCA/(\zeta_{i\Theta_{k}}\HCA\cap\HCA)\cong\IC f_{-k}(-1)
\end{equation}
whence the conclusion \halmos

\begin{lemma}\label{le:level 1}
$\FR$ and $\FNS$ are positive energy representations of $L\Spin_{2n}$ of level 1.
\end{lemma}
\proof
By corollary \ref{ch:classification}.\ref{co:global level}, it is
sufficient to compute the following cocycle on $\coroot\times\wt T$ where $\wt T$
is a maximal torus in $\Spin_{2n}$ and $\coroot\cong\Hom(\T,\wt T)$ the corresponding
coroot lattice. Define
\begin{equation}
\kappa(\alpha,\tau)=
\Gamma(\zeta_{\alpha})\Gamma(\tau)\Gamma(\zeta_{\alpha})^{*}\Gamma(\tau)^{*}
\end{equation}
Since $\Gamma$ factors through $L\SO_{2n}$, $\kappa$
extends to a bilinear map $I\times T\rightarrow\T$ where $T=\wt T/\IZ_{2}\subset\SO_{2n}$
and $I=\Hom(\T,T)$ is the integral lattice spanned by the $i\Theta_{k}$.
Moreover, the inclusion $U_{n}\subset\SO_{2n}$ identifies the diagonal torus of $U_{n}$
with $T$ and the corresponding integral lattices. Thus, $I,T\subset LU_{n}$ 
and since the Ramond and \NS Fock spaces are unitarily
equivalent as $LU_{n}$--modules, we need only compute $\kappa$ in one of them, say
$\FR$. Notice that $T\subset U_{n}$ commutes with the complex structure $J$ on $\HR$
so that it is canonically quantised and fixes $\Omega$. By lemma \ref{le:explicit coroot}
and \eqref{eq:T on fk}
\begin{equation}
\begin{split}
\Gamma(\zeta_{i\Theta_{k}})^{*}\Gamma(\exp_{T}(h))
\Gamma(\zeta_{i\Theta_{k}})\Gamma(\exp_{T}(h))^{*}\Omega
&=
\Gamma(\zeta_{\Theta_{k}})^{*}\Gamma(\exp_{T}(h))\frac{1}{\sqrt{2}}\psi(f_{-k}(-1))\Omega\\
&=
\Gamma(\zeta_{\Theta_{k}})^{*}\frac{1}{\sqrt{2}}\psi(\exp_{T}(h)f_{-k}(-1))\Omega\\
&=
e^{\<i\Theta_{k},h\>}\Omega
\end{split}
\end{equation}
where $\<\cdot,\cdot\>$ is the basic inner product so that
$\<\Theta_{k},\Theta_{l}\>=-\delta_{kl}$.
It follows by linearity that
\begin{equation}
\kappa(\alpha,\exp_{T}(h))=e^{-\<\alpha,h\>}
\end{equation}
and therefore that $\FR$ and $\FNS$ are of level 1 \halmos

\begin{proposition}\label{pr:level 1 reps}
Let $\H_{0},\H_{v},\H_{s_{\pm}}$ be the level 1 representations of $L\Spin_{2n}$ and
$\F=\FNS$ or $\FR$ the \NS and Ramond Fock spaces. Then, denoting by $\F_{0},\F_{1}$
the even and odd subspaces of $\F$, we have
\begin{xalignat}{2}
\H_{0}&\subseteq\FNSE&\H_{v}&\subseteq\FNSO
\label{eq:first set}\\
\intertext{and}
\H_{s_{+}}&\subseteq\F_{R,\epsilon}&\H_{s_{-}}&\subseteq\F_{R,1-\epsilon}
\label{eq:second set}
\end{xalignat}
where $\epsilon=0$ or $1$ according to whether $n$ is even or odd.
\end{proposition}
\proof
Let $\F=\FNS$.
The action of $\SO_{2n}\subset L\SO_{2n}$ on the complexification of $\H=\HNS$
preserves the decomposition \eqref{eq:decomposition} and therefore commutes with
the complex structure $J$. Thus, $\SO_{2n}$ is canonically quantised on $\F$ and
acts as
\begin{xalignat}{3}
\Gamma(R)\Omega&=\Omega&&\text{and}&
\Gamma(R)v_{1}(-n_{1})\wedge\cdots\wedge v_{k}(-n_{k})&=
Rv_{1}(-n_{1})\wedge\cdots\wedge Rv_{k}(-n_{k})
\label{canonical rot action}
\end{xalignat}
The lowest energy subspace of $\F_{0}$ is spanned by the vacuum vector $\Omega$
and therefore, by \eqref{canonical rot action}, $\F_{0}\cong\IC\cong V_{0}$ as
$\SO_{2n}$--modules. It follows by lemma
\ref{ch:classification}.\ref{le:generation} and lemma \ref{le:setup} that the
closure of the linear span of $\Gamma(L\Spin_{2n})\Omega\subset\F_{0}$ is isomorphic
to $\H_{0}$. Similarly, the lowest energy subspace of $\F_{1}$ is spanned by the
vectors $v(-\half{1})$, $v\in\VC$ and by \eqref{canonical rot action} is isomorphic,
as $\SO_{2n}$--module to $V_{\IC}$ so that $\H_{v}\subseteq\F_{1}$.
Let now $\F=\FR$. The operators $\psi(v)$, $v\in\VC$ satisfy
$U_{\theta}\psi(v)U_{\theta}^{*}=\psi(R_{\theta}v)=\psi(v)$ and therefore map
$\F(0)=\ext\VCH$ to itself. This action coincides in fact with the Clifford
algebra multiplication \eqref{eq:vector intertwiner}.
The projective action of $\Spin_{2n}\subset L\Spin_{2n}$ on $\F$ leaves $\F(0)$
invariant and, by \eqref{eq:real implementation} satisfies
$\Gamma(S)\psi(v)\Gamma(S)^{*}=\psi(Sv)$. It follows by the irreducibility of
$\ext\VCH$ under $C(V)$ that, as $\Spin_{2n}$--modules $\F(0)\cong\ext\VCH$.
The inclusions \eqref{eq:second set} now follow by proposition
\ref{pr:fd classification} and lemma \ref{ch:classification}.\ref{le:generation}
\halmos\\

\remark We shall prove in the next section that the inclusions
\eqref{eq:first set}--\eqref{eq:second set} are in fact equalities
by using the infinitesimal action of $\lpol\so_{2n}$ on $\FR,\FNS$.\\

\ssection{The infinitesimal action of $\lpol\so_{2n}$ on $\FR,\FNS$}
\label{se:currents}

In this section, we prove that the even and odd parts of the Ramond and \NS Fock
spaces are irreducible under the action of $L\Spin_{2n}$ and therefore that
the inclusions of proposition \ref{pr:level 1 reps} are equalities.
We proceed by studying the abstract action of $\lpol\so_{2n}$ on $\FR$ and $\FNS$
induced by that of $L\Spin_{2n}$ and identify it with
bilinear expressions in the Fermi field similar to those of \S \ref{ss:infinitesimal}.
The irreducibility result then follows from the well--known fact that these give a
level 1 representation of the Kac--Moody algebra $\wh{\so_{2n}}$ for which the even
and odd parts of $\FR,\FNS$ are irreducible \cite{IgF,GO2}.\\

In \S \ref{ss:alg fermions}, we consider the algebraic Clifford algebra of $\H=\HR$
or $\HNS$ generated by the $\psi(f)$ where $f\in\H$ is a trigonometric polynomial
and show that it acts irreducibly on the finite energy subspace of the corresponding
Fock space. The bilinear expressions giving the action of $\wh{\so_{2n}}$ are derived
in \S \ref{ss:bilinear}. Their identification with the abstract action of $\lpol\so_{2n}$
then follows by Shur's lemma because both representations have the same commutation
relations with the algebraic $\psi(f)$. In \S \ref{ss:irreducibility} we derive
from this the irreducibility of the even and odd subspaces of $\FR,\FNS$ under
$L\Spin_{2n}$.

\ssubsection{Algebraic Fermions}\label{ss:alg fermions}

Let $\H=\HR$ or $\HNS$ and $\F$ the corresponding Fock space. The subspace $\F\fin\subset\F$
of finite energy vectors for the canonically quantised action of $\rot$ given by
$U_{\theta}=\Gamma(R_{\theta})$ is spanned by the monomials
$v_{1}(-n_{1})\wedge\cdots\wedge v_{k}(-n_{k})$
where $n_{i}\in\IN$ for $\FR$ and $\half{1}+\IN$ for $\FNS$.
Consider the {\it algebraic Clifford algebra} $\calg\subset C(\H)$ generated by the
$\psi(v,n):=\psi(v(n))$, $v\in\VC$. By \eqref{eq:real implementation}
\begin{equation}\label{rot cov}
U_{\theta}\psi(v,n)U_{\theta}^{*}=e^{-in\theta}\psi(v,n)
\end{equation}
so that $\ffin$ is invariant under the action of $\calg$ and, differentiating
\begin{equation}
[d,\psi(v,n)]=-n\psi(v,n)
\end{equation}

We prove below that $\calg$ acts irreducibly on $\ffin$. The proof is similar to that of
lemma \ref{le:fd irred} with the replacement of the grading corresponding to the number
operator by that given by the infinitesimal generator of rotations $d$. It relies on the
fact that $d$ can be written as a bilinear expression in the $\psi(v,n)$
which we presently derive. To fix notations, if 
$e_{i}$ is an orthonormal basis of $V$, we denote $\psi(e_{i},n)$ by $\psi_{i}(n)$ so that,
by \eqref{eq:new anticommutation}
\begin{xalignat}{3}\label{psivn}
\{\psi_{i}(n),\psi_{j}(m)\}&=2\delta_{i,j}\delta_{m+n,0}&
&\text{and}&
\psi_{i}(n)^{*}&=\psi_{i}(-n)
\end{xalignat}
Thus, if we consider the $\psi_{i}(n)$ as anticommuting variables, then
$\displaystyle{\ad_{+}(\psi_{i}(-n))=2\frac{\partial}{\partial\psi_{i}(n)}}$ where
$\ad_{+}(a)b=\{a,b\}$. Recall the definition of {\it fermionic normal ordering}
\begin{equation}\label{eq:fermno}
\fermno{a(n)b(m)}=
\left\{\begin{array}{cl}
 a(n)b(m)&n<m \\
\half{1}(a(n)b(m)-b(m)a(n))&m=n \\
-b(m)a(n)&m>n
\end{array}\right.
=a(n)b(m)-H(n-m)\{a(n),b(m)\}
\end{equation}
where
\begin{equation}
H(p)=\begin{cases}0&p<0\\\half{1}&p=0\\1&p>0\end{cases}
\end{equation}
and the quantum field theoretic prescription that products
$\psi_{i_{1}}(n_{1})\cdots\psi_{i_{k}}(n_{k})$ should be written in normal order
to prevent the occurrence of infinities. With this in mind, we have using
$[ab,c]=a\{b,c\}-\{a,b\}c$
\begin{equation}\label{eq:sugawara derivation}
\begin{split}
-4n\psi_{i}(n)
&=\sum_{j,p}\Bigl(
  p\psi_{j}(-p)\{\psi_{j}(p),\cdot\}-p\{\psi_{j}(-p),\cdot\}\psi_{j}(p)
  \Bigr)\psi_{i}(n)\\
&=\sum_{j,p}[p\psi_{j}(-p)\psi_{j}(p),\psi_{i}(n)]\\
&=\sum_{j,p}[p\fermno{\psi_{j}(-p)\psi_{j}(p)},\psi_{i}(n)]
\end{split}
\end{equation}
where $p\in\IN$ for $\FR$ and $\half{1}+\IZ$ for $\FNS$. The operator
$\sum_{j,p}p\fermno{\psi_{j}(-p)\psi_{j}(p)}$ is well--defined on $\F\fin$ since the sum is
locally finite. In fact,

\begin{lemma}\label{le:sugawara}
The infinitesimal generator of rotations on $\F\fin$ is given by
\begin{equation}\label{eq:sugawara}
d=\frac{1}{4}\sum_{j,p}p\fermno{\psi_{j}(-p)\psi_{j}(p)}
     =\frac{1}{2}\sum_{j,p>0}p\psi_{j}(-p)\psi_{j}(p)
\end{equation}
\end{lemma}
\proof
Let $D$ be the operator defined by the right hand--side of \eqref{eq:sugawara}. By
\eqref{eq:sugawara derivation}, $[D,\psi(v,n)]=-n\psi(v,n)$. Moreover, $D\Omega=0$
where $\Omega\in\F$ is the vacuum vector and therefore by cyclicity of $\Omega$
under the action of the $\psi(v,n)$, $D=d$ \halmos

\begin{proposition}\label{pr:algebraic Schur}
$\F\fin$ is an irreducible $\calg$--module. Moreover, if $T\in\End(\F\fin)$ commutes with
$\calg$, then $T=\lambda\cdot 1$ for some $\lambda\in\IC$.
\end{proposition}
\proof
Notice first that any vector $\xi$ in $\F\fin(0)$ is cyclic for $\calg$, {\it i.e.}~
$\calg\xi=\F\fin$. This is obvious if $\F=\FNS$ since $\F\fin(0)=\IC\medspace\Omega$
and follows for $\FR$ from the irreducibility of $\F\fin(0)=\ext\VCH$ as
$C(V)\subset\calg$-module which implies $\cv\xi=\ext\VJ$.
If $T\in\End(\F\fin)$ commutes with $\calg$, it commutes with $d$ by lemma
\ref{le:sugawara} and therefore leaves $\F\fin(0)$ invariant. Thus, $T$ must have
an eigenvector $\xi\in\F\fin(0)$ with eigenvalue $\lambda$, and by cyclicity of $\xi$,
$T=\lambda\cdot 1$.
If $0\neq\Nu\subset\F\fin$ is invariant under $\calg$, then $d\Nu\subset\Nu$ by
lemma \ref{le:sugawara} so that $\Nu$ is graded, {\it i.e.}~$\Nu=\bigoplus_{n}\Nu(n)$.
The corresponding orthogonal projection $P_{\Nu}:\F\fin\rightarrow\Nu$ commutes
with $\calg$ and therefore $P_{\Nu}=1$ whence $\Nu=\F\fin$ \halmos

\ssubsection{The action of $\lpol\so_{2n}$ via bilinear expressions in Fermions}
\label{ss:bilinear}

Let $\H=\HR$ or $\HNS$ and $\F$ the corresponding Fock space. Denote by
$\Gamma:\lpol\so_{2n}\rightarrow\End(\F\fin)$ the projective representation
determined by the positive energy action of $L\Spin_{2n}$ on $\F$ via theorem
\ref{ch:classification}.\ref{core}.

\begin{lemma}\label{le:abstract covariance}
Let $f\in\HC$ be a trigonometric polynomial so that $\psi(f)\in\End(\F\fin)$
and $X\in\lpol\so_{2n}$. Then
\begin{equation}\label{eq:abstract covariance}
[\Gamma(X),\psi(f)]=\psi(Xf)
\end{equation}
\end{lemma}
\proof
By theorem \ref{ch:classification}.\ref{core}, $\Gamma(X)$ is essentially skew--adjoint
on $\F\fin$ and $e^{t\Gamma(X)}=\Gamma(e^{tX})$ in $PU(\F)$. By the covariance relation
\eqref{eq:real implementation}, we have
$e^{t\Gamma(X)}\psi(f)e^{-t\Gamma(X)}=\psi(e^{tX}f)$ and applying both sides to
$\xi\in\F\fin$, we find $e^{t\Gamma(X)}\psi(f)\xi=\psi(e^{tX}f)e^{t\Gamma(X)}\xi$.
The relations \eqref{eq:abstract covariance} now follow by taking derivatives at $t=0$.
Specifically,
$\psi(f)\xi\in\F\fin\subset\D(\Gamma(X))$
and therefore
\begin{equation}
e^{t\Gamma(X)}\psi(f)\xi=\psi(f)\xi+t\Gamma(X)\psi(f)\xi+o(t)
\end{equation}
Similarly,
\begin{equation}
\begin{split}
\psi(e^{tX}f)e^{t\Gamma(X)}\xi
&=\psi(f+tXf+o(t))(\xi+t\Gamma(X)\xi+o(t))\\
&=\psi(f)\xi+t(\psi(Xf)\xi+\psi(f)\Gamma(X)\xi)+o(t)
\end{split}
\end{equation}
since, by \eqref{eq:new anticommutation}, $\|\psi(g)\|\leq\sqrt{2}\|g\|$. Letting
$t\rightarrow 0$ yields \eqref{eq:abstract covariance} \halmos\\

Specialising the above, we find  
\begin{equation}\label{eq:to solve}
[E_{ij}(n),\psi_{k}(m)]=\delta_{j,k}\psi_{i}(m+n)-\delta_{i,k}\psi_{j}(m+n)
\end{equation}
where we denote $\Gamma(X)$ and $X$ by the same symbol and 
$E_{ij}\in\so_{2n}\cong\ext^{2}V$ is the basis given by $e_{i}\wedge e_{j}$.
Because of the irreducibility of the action of $\calg$, the above equations have essentially
unique solutions in the $E_{ij}(n)$. To solve them, notice that the right hand--side
of \eqref{eq:to solve} may be rewritten using the identity $a\{b,c\}-\{a,c\}b=[ab,c]$ as
\begin{equation}\label{eq:ansatz}
\begin{split}
&\half{1}\sum_{p}\Bigl(
\psi_{i}(p+n)\{\psi_{j}(-p),\cdot\}-\{\psi_{i}(-p),\cdot\}\psi_{j}(p+n)
\Bigr)\psi_{k}(m)=\\
&\half{1}\sum_{p}\Bigl(
\psi_{i}(p+n)\{\psi_{j}(-p),\cdot\}-\{\psi_{i}(p+n),\cdot\}\psi_{j}(-p)
\Bigr)\psi_{k}(m)=\\
&\half{1}\sum_{p}
[\psi_{i}(p+n)\psi_{j}(-p),\psi_{k}(m)]=\\
&\half{1}\sum_{p}
[\fermno{\psi_{i}(p+n)\psi_{j}(-p)},\psi_{k}(m)]
\end{split}
\end{equation}
{\it Define} now 
\begin{equation}\label{eq:currents}
E_{ij}(n)=\half{1}\sum_{p}\fermno{\psi_{i}(n+p)\psi_{j}(-p)}
\end{equation}
and use the above to define operators $X(n)$ for any $X\in\so_{2n}$. The following result is
well--known, see for example \cite{IgF}

\begin{proposition}\label{pr:level 1 commutation}
For $X,Y\in\so_{2n}$, the operators $X(n)$ satisfy
\begin{enumerate}
\item[(i)] $[X(n),\psi(v,m)]=\psi(Xv,m+n)$
\item[(ii)] $[X(n),Y(m)]=[X,Y](n+m)+n\delta_{n+m,0}\<X,Y\>$
\end{enumerate}
where $\<\cdot,\cdot\>$ is the basic inner product on $\so_{2n}$.
\end{proposition}
\proof
(i) follows from \eqref{eq:ansatz}. From (i) we deduce that $[X(n),Y(m)]-[X,Y](n+m)$
commutes with the $\psi(v,n)$ and hence by proposition \ref{pr:algebraic Schur} is equal to
a scalar $\lambda(X(n),Y(m))$. Since $[d,\psi(v,n)]=-n\psi(v,n)$, the operators $Z(p)$
are homogeneous of degree $-p$ and therefore $\lambda$ vanishes unless $m+n=0$. We may
compute it in that case by evaluating
$(\Omega,[X(n),Y(-n)]\Omega)-(\Omega,[X,Y](0)\Omega)$ and choosing for $X,Y$ basis
elements since $\lambda$ is bilinear in its arguments.
From
\begin{equation}
[ab,cd]=\{b,c\}ad-\{a,d\}cb+\{b,d\}ca-\{a,c\}bd
\end{equation}
which holds if all anti--commutators are scalars, we find for $n\geq 0$
\begin{equation}\label{eq:comm one}
\begin{split}
(\Omega,[E_{ij}(n),E_{kl}(-n)]\Omega)
&=\sum_{0\leq p,q\leq n}
(\Omega,[\half{1}\psi_{i}(n-p)\psi_{j}(p),\half{1}\psi_{k}(-n+q)\psi_{l}(-q)]\Omega)\\
&=\sum_{0\leq p,q\leq n}
\delta_{j,k}\delta_{p,n-q}(\Omega,\half{1}\psi_{i}(q)\psi_{l}(-q)\Omega)-
\delta_{i,l}\delta_{q,n-p}(\Omega,\half{1}\psi_{k}(-p)\psi_{j}(p)\Omega) \\
&+\sum_{0\leq p,q\leq n}
\delta_{j,l}\delta_{p,q}(\Omega,\half{1}\psi_{k}(-p)\psi_{i}(p)\Omega)-
\delta_{i,k}\delta_{p,q}(\Omega,\half{1}\psi_{j}(p)\psi_{l}(-p)\Omega)
\end{split}
\end{equation}
where $p,q\in\IZ$ for $\FR$ and $\half{1}+\IZ$ for $\FNS$. Using \eqref{eq:fermno},
we may rewrite \eqref{eq:comm one} as
\begin{equation}
\begin{split}
 &\sum_{0\leq p\leq n}
\delta_{j,k}(\Omega,\half{1}\fermno{\psi_{i}(p)\psi_{l}(-p)}\Omega)+
\delta_{i,l}(\Omega,\half{1}\fermno{\psi_{j}(p)\psi_{k}(-p)}\Omega) \\
-&\sum_{0\leq p\leq n}
\delta_{j,l}(\Omega,\half{1}\fermno{\psi_{i}(p)\psi_{k}(-p)}\Omega)+
\delta_{i,k}(\Omega,\half{1}\fermno{\psi_{j}(p)\psi_{l}(-p)}\Omega)\\
+&\sum_{0\leq p\leq n}
\delta_{i,l}\delta_{j,k}\biggl(H(2p)-H(-2p)\biggr)-
\delta_{i,k}\delta_{j,l}\biggl(H(2p)-H(-2p)\biggr)\\[1.2em]
=&
\delta_{j,k}(\Omega,E_{il}(0)\Omega)+\delta_{i,l}(\Omega,E_{jk}(0)\Omega)-
\delta_{j,l}(\Omega,E_{ik}(0)\Omega)-\delta_{i,k}(\Omega,E_{jl}(0)\Omega)
+n(\delta_{i,l}\delta_{j,k}-\delta_{i,k}\delta_{j,l})\\[1.0em]
=&
(\Omega,[E_{ij},E_{kl}](0)\Omega)+n(\delta_{i,l}\delta_{j,k}-\delta_{i,k}\delta_{j,l})
\end{split}
\end{equation}
The following lemma then shows that (ii) holds \halmos

\begin{lemma}\label{dynkin}
Let $\<\cdot,\cdot\>$ be the basic inner product on $\so_{2n}$. Then,
$\<E_{ij},E_{kl}\>=\delta_{j,k}\delta_{i,l}-\delta_{j,l}\delta_{i,k}$.
\end{lemma}
\proof
Let $V=\IR^{2n}$. The form $\tau(X,Y)=\tr_{V}(XY)$ is symmetric, bilinear and
ad--invariant and is therefore a multiple $\beta$ of $\<\cdot,\cdot\>$.
To compute $\beta$, let $h_{i}$ and $h^{i}$ be dual basis of the Cartan subalgebra
$\t\subset\so_{2n}$ with respect to $\<\cdot,\cdot\>$.
Denoting by $\Pi(\VC)$ the weights of $\VC$, and computing with a basis of $\VC$
diagonal for the action of $t$, we get
\begin{equation}
\beta\dim\t=
\tau(h_{i},h^{i})=
\sum_{\mu\in\Pi(V)}\mu(h_{i})\mu(h^{i})=
\sum_{\mu\in\Pi(V)}\|\mu\|^{2}
\end{equation}
Since the weights of $\VC$ are $\pm\theta_{i}$, we get $\beta=2$. To conclude,
notice that in $V$,
$E_{ij}E_{kl}=
\delta_{jk}\varepsilon_{il}+\delta_{il}\varepsilon_{jk}-
\delta_{ik}\varepsilon_{jl}-\delta_{jl}\varepsilon_{ik}$
where $\epsilon_{pq}e_{r}=\delta_{qr}e_{p}$ and therefore
$\tr_{V}(E_{ij}E_{kl})=2(\delta_{jk}\delta_{il}-\delta_{ik}\delta_{jl})$ \halmos

\ssubsection{Finite reducibility of $\FR,\FNS$ under $L\Spin_{2n}$}\label{ss:irreducibility}

\begin{proposition}\label{pr:irreducibility}
Let $\F$ be the \NS or Ramond sector. Then, the even and odd subspaces
of $\F$ are irreducible under $L\Spin_{2n}$ and therefore
\begin{xalignat}{3}
\FNS&=\H_{0}\oplus\H_{v}&
&\text{and}&
\FR&=\H_{s_{+}}\oplus\H_{s_{-}}
\end{xalignat}
\end{proposition}
\proof
Let $\Gamma$ and $\pi$ be the projective representations of $\lpol\so_{2n}$ on $\F\fin$
determined by the action of $L\Spin_{2n}$ and the bilinears \eqref{eq:currents} respectively.
By lemma \ref{le:abstract covariance} and proposition \ref{pr:level 1 commutation}, $\Gamma$
and $\pi$ have the same commutation relations with the $\psi(f)$, $f$ a trigonometric
polynomial and it follows by proposition \ref{pr:algebraic Schur} that
$\Gamma(X)=\pi(X)+\alpha(X)$ for some $\alpha(X)\in\IC$. Since the $\pi(X)$ act irreducibly
on the even and odd subspaces of $\F\fin$ \cite[thm. I.3.21]{IgF}, \cite[\S 5.6]{GO2} so do
the $\Gamma(X)$ and therefore, by proposition \ref{ch:classification}.\ref{simple}
$\F_{0},\F_{1}$ are irreducible under $L\Spin_{2n}$. The conclusion now follows from
proposition \ref{pr:level 1 reps} \halmos\\

\remark The fact that the abstract action of $\lpol\so_{2n}$ and that given by the
bilinears \eqref{eq:currents} coincide gives, by proposition \ref{pr:level 1 commutation}
another way of showing that $\F$ is a level 1 representation of $L\Spin_{2n}$.

\ssection{The level 1 vector primary fields}\label{se:vector fields}

We show below that the level 1 vector primary fields define bounded operator--valued
distributions by identifying them with the Fermi fields within the construction
of the level 1 representations of $L\Spin_{2n}$ given in sections
\ref{se:quark model}--\ref{se:currents}.
The same result holds in fact at any level and will be established in this generality
in chapter \ref{ch:sobolev fields}.\\

Let $\VC$ be the complexified defining representation of $\SO_{2n}$ and $V_{0}$,
$V_{s_{\pm}}$ the trivial and spin modules. By the tensor product rules of proposition
\ref{ch:classification}.\ref{pr:tensor with minimal}, the spaces
$\Hom_{\Spin_{2n}}(V_{i}\otimes\VC,V_{j})$ where $V_{i},V_{j}$ are admissible at level
1 are non--zero if, and only if $\{V_{i},V_{j}\}=\{V_{0},\VC\}$ or
$\{V_{i},V_{j}\}=\{V_{s_{+}},V_{s_{-}}\}$. Thus, by proposition \ref{pr:irreducibility},
the non--zero vector primary fields at level 1 are endomorphisms of the spaces of finite
energy vectors of the Ramond and \NS Fock spaces. We shall construct them as such.

\begin{proposition}\label{pr:L2 vector}
Let $\phi:\hfin_{i}\otimes\VC[z,z^{-1}]\rightarrow\hfin_{j}$ be a vector primary field
at level 1. Then, $\phi$ extends to a bounded map $\ltwo{\VC}\rightarrow\B(\H_{i},\H_{j})$
intertwining the action of $L\Spin_{2n}\rtimes\rot$ and satisfying
\begin{equation}\label{eq:ltwo}
\|\phi(f)\|\leq \sqrt{2}\|f\|
\end{equation}
\end{proposition}
\proof
We begin with the Ramond sector $\H=\HR$. Since $\HC=\ltwo{\VC}$ we may define a map
$\phi:\ltwo{\VC}\rightarrow\B(\FR)$ by $f\rightarrow\psi(f)$ where the latter acts
by \eqref{eq:complex create}--\eqref{eq:complex annihilate}.
If $P_{\epsilon}$, $\epsilon=0,1$ are the orthogonal
projections onto the even and odd subspaces of $\FR$, we claim that the restrictions
of $P_{\epsilon}\phi(\cdot)P_{1-\epsilon}$ to $\VC[z,z^{-1}]$ are primary fields of
charge $\VC$. We have already noted in \S \ref{ss:alg fermions} that if $f$ is a
polynomial, $\psi(f)$ defines an endomorphism of the finite energy subspace $\F\fin$.
Moreover, by lemma \ref{le:abstract covariance} and the fact that $\F_{0},\F_{1}$ are
invariant under $L\Spin_{2n}$, $P_{\epsilon}\phi(\cdot)P_{1-\epsilon}$ intertwines the
action of $\lpol\so_{2n}$ as claimed. Notice in passing that the initial terms of
$P_{\epsilon}\phi(\cdot)P_{1-\epsilon}$, that is the restriction to $\Spin_{2n}$
intertwiners
\begin{equation}
\ext_{1-\epsilon}\VCH=
\F_{1-\epsilon}(0)\rightarrow\F_{\epsilon}(0)=\ext_{\epsilon}\VCH
\end{equation}
coincided with the Clifford multiplication map given by
\eqref{eq:vector intertwiner}. By construction, the $P_{\epsilon}\phi(\cdot)P_{1-\epsilon}$
extend to $\ltwo{\VC}$ and, by \eqref{eq:real implementation} intertwine
$L\Spin_{2n}\rtimes\rot$. Moreover, by \eqref{eq:new anticommutation}
\begin{equation}
\|\psi(f)\omega\|^{2}+\|\psi(f)^{*}\omega\|=2\|f\|^{2}\|\omega\|^{2}
\end{equation} and \eqref{eq:ltwo} follows.
The construction of the vector primary fields in the \NS sector $\H=\HNS$ proceeds
in exactly the same way. Consider unitaries
$S_{\pm}:\ltwo{\VC}\rightarrow\HC$ given by
$f\rightarrow e^{\pm\half{1}(1\otimes i)\theta}f$. These intertwine the actions of
$L\Spin_{2n}$ and satisfy
$S^{\pm}R_{\theta}{S^{\pm}}^{*}=e^{\pm\half{1}(1\otimes i)\theta}R_{\theta}$.
If $P_{0}$, $P_{1}$ are the projections onto the even and odd subspaces of $\FNS$,
then $f\rightarrow P_{1}\psi(S_{-}f)P_{0},P_{0}\psi(S_{+}f)P_{1}$ are easily seen
to be the required primary fields \halmos







\newcommand {\Out}{\operatorname{Out}}
\newcommand {\Aut}{\operatorname{Aut}}
\newcommand {\Inn}{\operatorname{Inn}}
\newcommand {\mob}
{\begin{pmatrix}\alpha&\beta\\\overline{\beta}&\overline{\alpha}\end{pmatrix}}
\newcommand {\mobinv}
{\begin{pmatrix}\overline{\alpha}&-\beta\\-\overline{\beta}&\alpha\end{pmatrix}}
\newcommand {\J}{\mathcal J}

\chapter{Local loop groups and their associated von Neumann algebras}
\label{ch:loc loops}

This chapter assembles various results required for the definition of Connes fusion
to be given in chapter \ref{ch:connes fusion}.\\

The notion of {\it Connes fusion} of positive energy representations of $LG$, $G=\Spin_{2n}$
arises by regarding them as bimodules over the subgroups $L_{I}G,L_{I^{c}}G$ of loops
supported in a given interval $I\subset S^{1}$ and its complement and using a tensor
product operation on bimodules over von Neumann algebras due to Connes \cite{Co,Sa}.
Recall that a bimodule $\H$ over a pair $(M,N)$ of von Neumann algebras is a Hilbert
space supporting commuting representations of $M$ and $N$. To any two bimodules $X,Y$
over the pairs $(M,N)$, $(\wt N,P)$, Connes fusion associates an $(M,P)$--bimodule
denoted by $X\boxtimes Y$.
The definition of $X\boxtimes Y$ relies on, but is ultimately independent of the choice
of a reference or {\it vacuum} $(N,\wt N)$--bimodule $\Nu$ with a cyclic vector $\Omega$
for both actions and for which {\it Haag duality} holds, {\it i.e.}~the actions of $N$ and
$\wt N$ are each other's commutant.
Given $\Nu$, we form the intertwiner spaces $\XX=\Hom_{N}(\Nu,X)$ and
$\YY=\Hom_{\wt N}(\Nu,Y)$ and consider the sequilinear form on the algebraic
tensor product $\XX\otimes\YY$ given by
\begin{equation*}
\<x_{1}\otimes y_{1},x_{2}\otimes y_{2}\>=
(x_{2}^{*}x_{1} y_{2}^{*}y_{1}\Omega,\Omega)
\end{equation*}
where the inner product on the right hand--side is taken in $\Nu$. If $x_{1}=x_{2}$
and $y_{1}=y_{2}$, Haag duality implies that $x_{2}^{*}x_{1}$ and $y_{2}^{*}y_{1}$
are commuting positive operators and therefore that $\<\cdot,\cdot\>$ is positive
semi--definite. By definition, the bimodule $X\boxtimes Y$ is the Hilbert space
completion of $\XX\otimes\YY$ with respect to $\<\cdot,\cdot\>$, with $(M,P)$ acting
as $(m,p)x\otimes y=mx\otimes py$.\\

Applying the above to the positive energy representations of $LG$ requires a number
of preliminary results which are established in this chapter. Let $\pl$ be the set
of positive energy representations at a fixed level $\ell$. We wish to regard any
$(\H,\pi)\in\pl$ as a bimodule over the pair
$\pi_{0}(L_{I}G)^{''},\pi_{0}(L_{I^{c}}G)^{''}$
where $\pi_{0}$ is the vacuum representation at level $\ell$ whose lowest energy
subspace is, by definition the trivial $G$--module.
The well--foundedness of this change of perspective is justified by the following
properties

\begin{enumerate}
\item {\it Locality} :
$\pi(L_{I}G)''\subset\pi(L_{I^{c}}G)'$ for any $(\pi,\H)\in\pl$. In other words,
$\H$ is a $(\pi(L_{I}G)^{''},\pi(L_{I^{c}}G)^{''})$ bimodule.

\item {\it Local equivalence} :
All $(\pi,\H)\in\pl$ are unitarily equivalent as $L_{I}G$--modules. Thus we may
unambiguously identify $\pi(L_{I}G)''$ with $\pi_{0}(L_{I}G)''$ where $\pi_{0}\in\pl$
is the vacuum representation and consider $\H$ as a
$(\pi_{0}(L_{I}G)'',\pi_{0}(L_{I^{c}}G)'')$--bimodule.

\item {\it von Neumann density} :
$\pi(L_{I}G)\times\pi(L_{I^{c}}G)$ is strongly dense in $\pi(LG)$. Thus, inequivalent
irreducible positive energy representations of $LG$ remain so when regarded as bimodules.
\end{enumerate}

The r\^ole of the reference bimodule is played by the vacuum representation
$(\pi_{0},\H_{0})$. Two crucial facts need to be established in this respect

\begin{enumerate}
\item[(iv)] {\it Reeh-Schlieder theorem} : any finite energy vector of a positive
energy representation $\pi$ is cyclic under $\pi(L_{I}G)$. In particular, the
lowest energy vector $\Omega\in\H_{0}(0)$ is cyclic for $\pi_{0}(L_{I}G)''$
and $\pi_{0}(L_{I^{c}}G)''$.
\item[(v)] {\it Haag duality} : $\pi_{0}(L_{I}G)^{''}=\pi_{0}(L_{I^{c}}G)^{'}$.
\end{enumerate}

Finally, another technically crucial property of the algebras $\pi(L_{I}G)''$ is the
following
\begin{enumerate}
\item[(vi)] {\it Factoriality} : The algebras $\pi(L_{I}G)''$ with $I\subsetneq S^{1}$
are type III$_{1}$ factors.
\end{enumerate}

Properties (i)-(vi) were established by Jones and Wassermann for $G=SU_{n}$ in \cite{Wa3}
and, in unpublished notes for all simple, simply--connected compact Lie groups using the
quark model and conformal inclusions \cite{Wa1}. In section \ref{se:local groups}, we give
a model independent proof of locality and von Neumman density in the latter generality.
The remaining properties are proved in section \ref{se:local algebras} 
by following Jones and Wassermann's original lines.

\ssection{Local loop groups}\label{se:local groups}

For any open, possibly improper, interval $I\subseteq S^{1}$, define the local loop
group $L_{I}G=\{f\in LG|\left.f\right|_{S^{1}\backslash I}\equiv 1\}$. The
 complementary interval $I^{c}$ is, by definition $S^{1}\backslash\overline{I}$.

\ssubsection{Locality}\label{ss:locality}

The following is an immediate generalisation of proposition 3.4.1. of \cite{PS}.

\begin{lemma}\label{local is perfect}
Assume $G$ is simple. Then for any open interval $I\subseteq S^{1}$, $L_{I}G$
is perfect {\it i.e.}~it is equal to its commutator subgroup. In particular,
$\Hom(L_{I}G,\T)=\{1\}$.
\end{lemma}
\proof Consider the commutator map $C:G\times G\rightarrow G$,
$(g,h)\rightarrow ghg^{-1}h^{-1}$. Its differential at $(1,1)$ is given
by the bracket $[\cdot,\cdot]:\g\times\g\rightarrow\g$ and is therefore
surjective since $G$ is simple. By the implit function theorem, $C$
has a smooth right inverse defined on a neighborhood $U$ of the identity
and mapping $1$ to $(1,1)$. It follows that $L_{I}U\subset L_{I}G$
is contained in $[L_{I}G,L_{I}G]$ and therefore
$L_{I}G=\bigcup_{n}L_{I}U^{n}\subset[L_{I}G,L_{I}G]$ since $L_{I}G$ is
connected \halmos

\begin{proposition}[Locality]
If $I,J\subset S^{1}$ are disjoint open intervals and $(\pi,\H)$ is a
positive energy representation of $LG$, then
\begin{equation}
\pi(L_{I}G)^{\prime\prime}\subset\pi(L_{J}G)^{\prime}
\end{equation}
\end{proposition}
\proof
For any $\gamma_{I},\gamma_{J}\in L_{I}G\times L_{J}G$ the following holds in
$U(\H)$ :
$\pi(\gamma_{I})\pi(\gamma_{J})\pi(\gamma_{I})^{*}\pi(\gamma_{J})^{*}=
\chi(\gamma_{I},\gamma_{J})$, 
where $\chi(\gamma_{I},\gamma_{J})\in\T$ does not depend upon the choices of the
particular lifts and is multiplicative in either variable. By
lemma \ref{local is perfect}, $\chi\equiv 1$ and $\pi(L_{I}G)$ commutes with
$\pi(L_{J}G)$ \halmos

\ssubsection{The Sobolev $\half{1}$ space}\label{ss:sobolev half}

We establish below the density of the space of smooth, complex--valued functions
vanishing to all orders on a finite subset of $S^{1}$ in $C^{\infty}(S^{1})$ for
the Sobolev $\half{1}$--norm. We shall use this result in \S \ref{ss:von Neumann} to
prove von Neumann density by showing that the topology on $LG$ obtained by pulling
back the strong operator topology on the projective unitary group of a positive
energy representation is weaker than the Sobolev $\half{1}$--topology. The discussion
below is taken from \cite[Chap. VIII]{Wa1}.\\

Let $A\subset S^{1}$ be a finite subset and denote by $C_{A}^{k}(S^{1})$,
$k\in\IN\cup\{\infty\}$, the space of smooth functions vanishing up to order $k$ on $A$.
If $f:S^{1}\rightarrow\IC$ is a trigonometric polynomial with Fourier series
$\sum_{k}a_{k}e^{ik\theta}$ and $s\in\IR$, define
$\|f\|_{s}^{2}=\sum_{k}(1+|k|)^{2s}|a_{k}|^{2}$ and
$|f|_{s}=\sum_{k}(1+|k|)^{s}|a_{k}|$. Denote by $H^{s}(S^{1})$ the completion of
the space of trigonometric polynomials with respect to the norm $\|\cdot\|_{s}$. It
contains $C^{\infty}(S^{1})$ since the Fourier coefficents of a smooth function decrease
faster than any polynomial.

\begin{lemma}
If $f,g\in C^{\infty}(S^{1})$ and $s\geq 0$, then $\|fg\|_{s}\leq |f|_{s}\|g\|_{s}$.
\end{lemma}
\proof Consider first the case $f=a_{k}e^{ik\theta}$ and let
$g=\sum_{n}b_{n}e^{in\theta}$. Then
\begin{equation}
\|fg\|_{s}^{2}=
|a_{k}|^{2}\sum_{n}(1+|n|)^{2s}|b_{n-k}|^{2}\leq
|a_{k}|^{2}(1+|k|)^{2s}\sum_{n}(1+|n-k|)^{2s}|b_{n-k}|^{2}=
|a_{k}|^{2}(1+|k|)^{2s}\|g\|_{s}^{2}
\end{equation}
since $(1+|p+q|)\leq(1+|p|)(1+|q|)$. Thus, if $f=\sum_{k}a_{k}e^{ik\theta}$,
we have $\|fg\|_{s}\leq\sum_{k}|a_{k}|(1+|k|)^{s}\|g\|_{s}=|f|_{s}\|g\|_{s}$
\halmos

\begin{lemma}
$C^{1}_{A}(S^{1})$ is dense in $C^{\infty}(S^{1})$ for the $\|\cdot\|_{\half{1}}$
norm.
\end{lemma}
\proof Consider first the case $A=\{0\}$.
Let $f_{n}=\sum_{k=2}^{n}\frac{\cos(k\theta)}{k\log k}$ and set
$\chi_{n}=f_{n}f_{n}(0)^{-1}$. Then $f=\lim_{n}f_{n}$ exists in
$H^{\half{1}}$ but $\lim_{n}f_{n}(0)=\infty$ and therefore $\chi_{n}(0)=1$,
$\chi_{n}^{\prime}(0)=0$ and $\chi_{n}\rightarrow 0$ in $H^{\half{1}}$. Thus,
if $h\in C^{\infty}(S^{1})$, then $h_{n}=h(1-\chi_{n})\in C^{1}_{A}(S^{1})$ and by
the above lemma, $\|h-h_{n}\|_{\half{1}}\leq|h|_{\half{1}}\|\chi_{n}\|_{\half{1}}$
which tends to zero. The general case $A=\{\theta_{1},\ldots,\theta_{k}\}$ follows
easily \halmos

\begin{proposition}\label{pr:half density}
$C^{\infty}_{A}(S^{1})$ is dense in $C^{\infty}(S^{1})$ for the $\|\cdot\|_{\half{1}}$
norm.
\end{proposition}
\proof
It is sufficient to establish the density of $C^{\infty}_{A}(S^{1})$ in
$C^{1}_{A}(S^{1})$. Moreover, since
\begin{equation}
\|f\|_{\half{1}}^{2}\leq
 \sum_{k}|a_{k}|^{2}(1+k^{2})=
 \|f\|_{L^{2}}^{2}+\|\dot f\|_{L^{2}}^{2}\leq
 \|f\|_{\infty}^{2}+\|\dot f\|_{\infty}^{2}\leq
 (4\pi^{2}+1)\|\dot f\|_{\infty}^{2}
\end{equation}
for $f\in C^{1}_{A}(S^{1})$, density for the $C^{1}$ norm will do.
Let $f\in C^{1}_{A}(S^{1})$. By the Stone-Weierstrass theorem, there exists a sequence
$g_{n}\in C^{\infty}_{c}(S^{1}\backslash A)$ such that
$\|\dot f-g_{n}\|_{\infty}\rightarrow 0$. If $A=\{\theta_{1},\ldots,\theta_{n}\}$ with
the points arranged in incresing order, choose $\rho_{i}$ smooth and supported in
$(\theta_{i},\theta_{i+1})$ with $\int\rho_{i}=1$. Set
$g_{n}^{\prime}=g_{n}-\sum_{i}\rho_{i}\int_{\theta_{i}}^{\theta_{i+1}}g_{n}(t)dt$ and
$G_{n}(\theta)=\int_{\theta_{1}}^{\theta}g_{n}^{\prime}(t)dt$ so that
$G_{n}\in C^{\infty}_{c}(S^{1}\backslash A)$. Since
$\int_{\theta_{i}}^{\theta_{i+1}}\dot f=0$, we have
\begin{equation}
\|\dot f-\dot G_{n}\|_{\infty}\leq
\|\dot f-g_{n}\|_{\infty}+
\sum_{i}\|\rho_{i}\|_{\infty}\int_{\theta_{i}}^{\theta_{i+1}}|g_{n}-\dot f|dt
\rightarrow 0
\end{equation}
\halmos

\ssubsection{Von Neumann density}\label{ss:von Neumann}

Our proof of von Neumann density follows Jones and Wassermann's original lines \cite{Wa1}.
It bypasses however the use of conformal inclusions to obtain the general case from
$SU_{n}$ or $\Spin_{2n}$ by making use of the Sobolev estimates of Goodman and Wallach
\cite[\S 3.2]{GoWa}.\\

Let $I\subseteq S^{1}$ be an open interval and $A\subset I$ a finite set.
Denote by $L_{I}^{A}G\subset L_{I}G$ the normal subgroup of loops equal to 1 
to all orders on $A$. Clearly, if $I\backslash A=I_{1}\sqcup\ldots\sqcup I_{k}$,
then $L_{I}^{A}G=L_{I_{1}}G\times\cdots\times L_{I_{k}}G$. We will prove that
$\pi(L_{I}^{A}G)$ is dense in $\pi(L_{I}G)$ for any positive energy representation
$(\pi,\H)$. We shall need for this purpose a different version of the Sobolev
estimates governing the action of $L\gc$ on the space of smooth vectors $\hsmooth$
obtained in proposition \ref{ch:analytic}.\ref{sobolev estimates}.
For $X=\sum_{k}X(k)\in L\gc$ and $s\in\IR$ define 
$\|X\|_{s}^{2}=\sum_{k}(1+|k|)^{2s}|X(k)|^{2}$ and recall that
$|X|_{s}=\sum_{k}(1+|k|)^{s}|X(k)|$ where $|\cdot|$ is a norm on $\gc$. Denote as
customary by $L_{0}$ the infinitesimal generator of rotations given by the
Segal--Sugawara formula (\ref{ch:analytic}.\ref{sugawara}).

\begin{lemma}
For any $X\in L\gc$, $\xi\in\H^{\infty}$ and $\epsilon>0$ the following holds
\begin{equation}
\|\pi(X)\xi\|\leq C_{\epsilon}\|X\|_{\half{1}}\|(1+L_{0})^{1+\epsilon}\xi\|
\end{equation}
\end{lemma}
\proof
Let $X=X(k)\in\gc\otimes e^{ik\theta}$ and $\xi\in\H$ be a finite energy vector.
By the estimates of proposition \ref{ch:analytic}.\ref{sobolev estimates},
\begin{equation}
\|\pi(X)\xi\|^{2}\leq C^{2}(1+|k|)|X(k)|^{2}\|(1+L_{0})^{\half{1}}\xi\|^{2}
\end{equation}
Thus, if $X=\sum_{k}X(k)$ and $\xi$ is an eigenvector of $L_{0}$ so that the
$\pi(X(k))\xi$ are orthogonal, we get
$\|\pi(X)\xi\|\leq C\|X\|_{\half{1}}\|(1+L_{0})^{\half{1}}\xi\|$.
More generally, if $\xi=\sum_{m}\xi_{m}$ with $L_{0}\xi_{m}=m\xi_{m}$ and
$\epsilon>0$, then
\begin{equation}
\begin{split}
\|\pi(X)\xi\|
&\leq
C\|X\|_{\half{1}}\sum_{m}\|(1+L_{0})^{\half{1}}\xi_{m}\|\\
&\leq
C\|X\|_{\half{1}}
\Bigl(\sum_{m}(1+m)^{-1-2\epsilon}\Bigr)^{\half{1}}
\Bigl(\sum_{m}\|(1+L_{0})^{1+\epsilon}\xi_{m}\|^{2}\Bigr)^{\half{1}}\\
&\leq
C_{\epsilon}\|X\|_{\half{1}}\|(1+L_{0})^{1+\epsilon}\xi\|\\
\end{split}
\end{equation}
\halmos

\begin{proposition}\label{pr:sobolev density}
Let $\pi$ be a positive energy representation of $LG$. Then $\pi(L_{I}^{A}G)$ is
strongly dense in $\pi(L_{I}G)$.
\end{proposition}
\proof
Since $L_{I}G$ is connected, is suffices to show that
$\pi(V)\subset\overline{\pi(L_{I}^{A}G)}$ for a suitable neighborhhod
$1\in V\subset L_{I}G$ for then
$\pi(L_{I}G)=\bigcup_{n}\pi(V)^{n}\subset\overline{\pi(L_{I}^{A}G)}$.
Let $\mathcal U$ be a neighborhood of the identity in $G$ on which the logarithm is
well-defined and set $V=L_{I}\mathcal U$. If $\gamma=\exp_{LG}(X)\in L_{I}G$, we claim
that for any $\xi\in\H$, $\pi(\gamma)\xi=\exp(\pi(X))\xi$ may be
arbitrarily well approximated by $\exp(\pi(Y))\xi$ with $Y\in L_{I}^{A}\g$.
By the norm-boundedness of the unitaries $\pi(\gamma)$ it is sufficent to prove this
on the dense subspace of vectors $\xi\in\hsmooth\subset\H$. Since these are
invariant under $LG$ by proposition \ref{ch:analytic}.\ref{invariance}, the function
$F(t)=e^{-t\pi(Y)}e^{t\pi(X)}\xi$ is smooth and
$\dot F(t)=e^{-t\pi(Y)}(\pi(X)-\pi(Y))e^{t\pi(X)}\xi$. Therefore,
\begin{equation}
\begin{split}
\|e^{\pi(Y)}\xi-e^{\pi(X)}\xi\|
&=\|e^{\pi(Y)}(1-e^{-\pi(Y)}e^{\pi(X)})\xi\|\\
&\leq\int_{0}^{1}\|e^{\pi(Y)}\dot F(t)\|dt  \\
&\leq C_{\epsilon}M_{X}\|X-Y\|_{\half{1}}\|(1+L_{0})^{1+\epsilon}\xi\|
\end{split}
\end{equation}
where we have used that $LG$ acts continuously on each Sobolev space $\H^{s}$.
The result now follows from the density of $L_{I}^{A}\g$ in $L_{I}\g$ with
respect to the $\|\cdot\|_{\half{1}}$ norm given by proposition
\ref{pr:half density} \halmos

\begin{corollary}[von Neumann density]
On any positive energy representation $\pi$, we have
\begin{equation}
\pi(L_{I}G)''\vee\pi(L_{I^{c}}G)''=\pi(LG)''
\end{equation}
\end{corollary}

\remark The Sobolev estimates of proposition \ref{ch:analytic}.\ref{sobolev estimates},
namely $\|\pi(X)\xi\|\leq C|X|_{\half{1}}\|(1+L_{0})^{\half{1}}\xi\|$ where
$|X|_{s}=\sum_{k}|X(k)|(1+|k|)^{s}$ cannot be used directly to establish the above density 
result. Indeed, $|\cdot|_{s}$ dominates the $C(S^{1},\g)$ norm for any $s\geq 0$ and
therefore $L_{I}^{A}\g$ is {\bf not} dense in $L_{I}\g$ for the corresponding topology.

\ssubsection{The Reeh--Schlieder theorem}\label{ss:Reeh}

\begin{proposition}[Reeh--Schlieder theorem]
Let $(\pi,\H)$ be an irreducible positive energy representation of $LG$ and $I\subset S^{1}$
an open interval. Then, any finite energy vector $\xi\in\H$ is cyclic for $L_{I}G$ {\it i.e.}~
the linear span of $\pi(L_{I}G)\xi$ is dense in $\H$.
\end{proposition}
\proof This is proved in \cite[\S 17]{Wa3} \halmos

\begin{corollary}
If $I\subset S^{1}$ is an open non-dense interval and $(\pi,\H)$ an irreducible positive
energy representation of $LG$ then any finite energy-vector $\xi\in\H$ is cyclic and
separating for $M_{I}=\pi(L_{I}G)''$, {\it i.e.}~
$\overline{M_{I}\xi}=\overline{M_{I}'\xi}=\H$. 
\end{corollary}
\proof The result follows at once from the Reeh--Schlieder theorem and locality since
$M_{I}'\supset\pi(L_{I^{c}}G)$ \halmos

\ssection{Von Neumann algebras generated by local loops groups}\label{se:local algebras}

In this section, we establish the Haag duality, local equivalence and factoriality
properties for the von Neumann algebras generated by local loop groups $L_{I}G$ in
positive energy representations.\\

The proof of Haag duality relies on the modular theory of Tomita and Takesaki which
we outline in \S \ref{ss:modular theory}. This allows in principle to compute the
commutant of a von Neumann algebra $M$ acting on a Hilbert space $\H$ with a cyclic
and separating vector $\Omega$, in that it yields the existence of a canonical
conjugation $J$ on $\H$ such that $M'=JMJ$. When $M$ is the von Neumann algebra
generated by a local loop group $L_{I}G$ in the vacuum representation $\H_{0}$
and $I$ is the upper semi--circle, the action of $J$ corresponds to the map
$z\rightarrow z^{-1}$ on $S^{1}$ exchanging $I$ and $I^{c}$. It follows from this
that $\pi_{0}(L_{I}G)'=\pi_{0}(L_{I^{c}}G)''$ and therefore that Haag duality holds.
This explicit characterisation of $J$ is obtained in \S \ref{ss:Haag} by using the
fermionic realisation of $\H_{0}$ obtained in chapter \ref{ch:fermionic} and a
computation of the modular conjugation for local Clifford algebras due to Jones
and Wassermann \cite{Wa3} which we review in \S \ref{ss:modular fermions}. The local
equivalence and factoriality properties are proved in \S \ref{ss:local equivalence}.

\ssubsection{Modular theory}\label{ss:modular theory}

We outline Tomita and Takesaki's modular theory, details may be found in \cite{BR}.
Let $M\subset\B(\H)$ be a von Neumann algebra with commutant $M'$. Notice that
$\xi\in\H$ is {\it cyclic} for $M$, {\it i.e.}~$\overline{M\xi}=\H$ iff it is {\it separating}
for $M'$ {\it i.e.}~$x\xi=0$, $x\in M'$ implies $x=0$. Indeed, if $\xi$ is cyclic and
annihilated by $x\in M'$, then $xM\xi=Mx\xi=0$ and therefore $x=0$. Conversely,
the orthogonal projection $p$ on $\overline{M\xi}$ lies in $M'$ and $(1-p)\xi=0$
whence $p=1$.\\

Assume now that $\Omega\in\H$ is cyclic and separating for $M$ so that
$\overline{M\Omega}=\H=\overline{M'\Omega}$ and $\H$ is the GNS representation
corresponding to the (faithful and normal) state $\phi_{\Omega}(x)=(x\Omega,\Omega)$,
$x\in M$. Let $\kappa$ be the densely defined conjugation given by
$\kappa x\Omega=x^{*}\Omega$
on $\D(\kappa)=M\Omega$. If $\phi_{\Omega}$ is a tracial state {\it i.e.}~
$\phi_{\Omega}(ab)=\phi_{\Omega}(ba)$ then
$\kappa$ is an isometry and therefore extends to the whole of $\H$. As is easily
verified, $\kappa M\kappa'\subset M'$ and von Neumann's commutation theorem asserts that in
that case $\kappa M\kappa=M'$.
More generally, $\kappa$ is closeable with closure $\overline{\kappa}=J\Delta^{\half{1}}$
where $J$ is a conjugate linear isometry,
the phase of the polar decomposition of $\overline{\kappa}$ and the modulus
$\Delta^{\half{1}}$ is a positive self-adjoint operator which measures the failure of
$\phi_{\Omega}$ to be tracial. Tomita's fundamental theorem asserts that $JMJ=M'$ and that,
in addition $\Delta^{it}M\Delta^{-it}=M$.
$J$ and $\Delta^{it}$ are called the modular conjugation and group.\\

The modular group is a fundamental tool in the study of $M$ in view of Connes' result that 
the class of $\Delta^{it}$ in $\Out(M)=\Aut(M)/\Inn(M)$ is independent of the particular
pair $(\Omega,\H)$. This implies that if $\Delta^{it}$ acts ergodically, {\it i.e.}~without
fixed points on $M$ then $M$ is a factor (for $Z(M)$ is fixed by $\Ad(\Delta^{it})$)
of type III$_{1}$ for in all other cases the modular group is inner, as a
suitable choice of $\H$ shows. An important special case arises when $\Delta^{it}$
leaves no vectors in $\H$ fixed aside from $\Omega$. Then $\Delta^{it}x\Delta^{-it}=x$
for some $x\in M$ implies that $x\Omega$ is fixed by the modular group and therefore
$x\Omega=\alpha\Omega$ for some $\alpha\in\IC$. Since $\Omega$ is separating,
$x\equiv\alpha$ and the modular group acts ergodically on $M$.\\

The modular operators are {\it hereditary} in the following sense \cite{Ta}. If
$N\subset M$ is a sub--von Neumann algebra such that $\Delta^{it}N\Delta^{-it}=N$,
then $J$ and $\Delta^{it}$ restrict to the modular operators for $N$ acting on
$\K=\overline{N\Omega}$.
More precisely, if $p\in N'$ is the orthogonal projection on $\K$, then $y\rightarrow yp$
is an isomorphism of $N$ onto $pNp=Np\subset\B(\K)$ for $yp=0$ implies
$y\Omega=y(1-p)\Omega=0$ and therefore $y=0$ since $\Omega$ is separating. $\Omega\in\K$
is cyclic and separating for $pNp$ and $J$ and $\Delta^{it}$ leave $\K$ invariant and
restrict to the modular operators for $pNp$ relative to $\Omega$. Following Jones and
Wassermann \cite{Wa1,Wa2,Wa3}, we shall refer to the hereditarity of the modular operators
as {\it Takesaki devissage}.

\ssubsection{Modular operators for complex fermions on $S^{1}$}\label{ss:modular fermions}

We describe, following Jones and Wassermann \cite{Wa1,Wa3}, the modular operators for
the local CAR algebras $\AA(L^{2}(I,\IC^{n}))$, $I\subset S^{1}$ in the Fock space
corresponding to the basic representation of $LU_{n}$. Their explicit knowledge will
be used in \S \ref{ss:Haag} to deduce Haag duality in the vacuum representation of
$L\Spin_{2n}$ from its fermionic realisation.\\

Recall that the group
\begin{equation}
\SU(1,1)_{\pm}=\{\mob|\thinspace |\alpha|^{2}-|\beta|^{2}=\pm 1\}
\end{equation}
acts on $S^{1}=\{z|\thinspace |z|=1\}$ by M\"obius transformations. $\SU(1,1)_{\pm}$ is
the semi--direct product of the connected component of the identity $\SU(1,1)$ by its 
group of components $\IZ_{2}$ generated by the matrix
\begin{equation}\label{eq:j}
j=\left(\begin{array}{rr}0&-1\\-1&0\end{array}\right)
\end{equation}
acting on $S^{1}$ as the flip $z\rightarrow z^{-1}$.
$\SU(1,1)$ acts by orientation preserving diffeomorphisms of the circle and maps
any three points $z_{1},z_{2},z_{3}\in S^{1}$ to any other three points $w_{1},w_{2},w_{3}$
provided $z_{2}$ and $w_{2}$ lie on the anticlockwise circular arcs going from $z_{1}$
to $z_{3}$ and $w_{1}$ to $w_{3}$ respectively.
In particular, there is a one--parameter group in $\SU(1,1)_{+}$ of elements fixing the
endpoints of the upper semi--circle $\{z\in S^{1}|\Im z\geq 0\}$ given by
\begin{equation}\label{mob flow}
t\rightarrow d^{t}=
\begin{pmatrix}\cosh\pi t&\sinh\pi t\\\sinh\pi t&\cosh\pi t\end{pmatrix}
\end{equation}

Let now $\H=L^{2}(S^{1},\IC^{n})$ with CAR algebra $\AA(\H)$ generated by the
$\IC$--linear symbols $c(f)$ subject to
\begin{xalignat}{2}
\{c(f),c(g)\}&=0&
\{c(f),c(g)^{*}\}&=(f,g)
\end{xalignat}
Denote by $\H_{+}\subset\H$ the space of functions with vanishing positively moded
Fourier coefficents, {\it i.e.}~the boundary values of holomorphic functions on $|z|>1$
and by $\H_{-}$ its orthogonal complement. Let $\J$ be the complex structure acting
as multiplication by $\pm i$ on $\H_{\pm}$ and denote by $\H_{\J}$ the Hilbert space
$\H$ with complex multiplication given by $\J$. We consider, as in \S \ref{ss:LU} of
chapter \ref{ch:fermionic}, the corresponding representation of $\AA(\H)$ on the Fock
space $\F=\ext\H_{\J}$. The subgroup of the unitary group $U(\H)$ commuting with $\J$
is canonically quantised on $\F$ by
\begin{equation}
\Gamma(u)f_{1}\wedge\cdots\wedge f_{m}=uf_{1}\wedge\cdots\wedge uf_{m}
\end{equation}
Moreover, unitaries of $\H$ anticommuting with $\J$ are canonically represented on
$\F$ by anti--unitaries.\\

$SU(1,1)_{\pm}$ acts unitarily on $\H$ by
\begin{equation}
\mob f(z)=\left(\frac{z}{\overline{\alpha}z-\beta}\right)f(\mob^{-1}z)
\end{equation}
When restricted to $SU(1,1)$, this action leaves $\H_{+}$ invariant since
$|\alpha|>|\beta|$ on this subgroup and therefore 
$z\rightarrow z(\overline{\alpha}z-\beta)^{-1}$ is holomorphic on $|z|>1$. Thus,
$SU(1,1)$ commutes with $\J$ and is therefore canonically quantised on $\F$.
On the other hand, the flip \eqref{eq:j} acts by $jf(z)=zf(z^{-1})$ and therefore
anticommutes with $\J$ so that $SU(1,1)_{-}=jSU(1,1)$ acts on $\F$ by antiunitaries.\\

For any open, non dense interval $I\subset S^{1}$, denote by $\AA(I)\subset B(\F)$
the von Neumann algebra generated by the operators $c(f)$ with $f$ supported in $I$.
Let now $I=(0,\pi)$ be the upper semi--circle.

\begin{theorem}[Jones--Wassermann \cite{Wa1,Wa3}]\label{modular}\hfill
\begin{enumerate}
\item The vacuum vector $\Omega\in\ext^{0}\H_{\J}\subset\F$ is cyclic and separating for
$\AA(I)$.
\item The corresponding modular group $\Delta_{I}^{it}$ is the canonical quantisation of
the action of the M\"obius flow \eqref{mob flow} on $\H$.
\item The corresponding modular conjugation $J_{I}$ is given by $\kappa^{-1}\Gamma(j)$
where $\kappa=i^{\epsilon}$ is the Klein transform corresponding to the natural
$\IZ_{2}$--grading $\epsilon$ on $\F$ and $j$ is the action of the flip \eqref{eq:j} on
$\H$.
\item $\Delta_{I}^{it}$ leaves no vectors in $\F$ fixed aside from $\Omega$.
\end{enumerate}
\end{theorem}

\remark Similar results are obtained for a general $I\subsetneq S^{1}$ by using the
action of $\SU(1,1)$ on $\F$.\\

\remark It follows at once from Tomita's theorem and the above that
$\AA(I)'=J\AA(I)J=\kappa^{-1}\AA(I^{c})\kappa$. Moreover, the $\AA(I)$
are type III$_{1}$ factors since $\Delta_{I}^{it}$ acts ergodically.

\ssubsection{Haag duality in the vacuum sector}\label{ss:Haag}

Recall from chapter \ref{ch:fermionic} that the Neveu--Schwarz Hilbert space $\HNS$ of
anti--periodic, $\IC^{n}$--valued functions on $S^{1}$ carries an orthogonal action of
$L\Spin_{2n}\rtimes\rot$. If $\wt\J$ is the complex structure acting as multiplication
by $i$ on the subspace of functions with vanishing positively--moded Fourier coefficients
and by $-i$ on its orthogonal complement, then by lemma
\ref{ch:fermionic}.\ref{le:HS cond}, $L\Spin_{2n}$ commutes with $\wt\J$ up to
Hilbert--Schmidt operators and therefore acts on the Fock space
$\FNS=\ext\H_{\operatorname{NS},\wt\J}$.
By proposition \ref{ch:fermionic}.\ref{pr:irreducibility}, the corresponding positive
energy representation is the sum of the level 1 vacuum and vector representations.\\

$\HNS$ is unitarily isomorphic to $\H=\ltwo{\IC^{n}}$ considered in
\S \ref{ss:modular fermions} via the map $f\rightarrow z^{\half{1}}f$ which
identifies the complex structures $\wt\J$ and $\J$. Transporting the action of
$L\Spin_{2n}\rtimes\rot$ via this identification,
we see that $L\Spin_{2n}$ acts on $\H$ by
$\gamma f(z)=z^{\half{1}}\gamma(z)z^{-\half{1}}f(z)$
and $\rot$ by $R_{\theta}f(z)=e^{i\half{\theta}}f(e^{-i\theta}z)$
\footnote{This is not unitarily equivalent to the action on $\H$ given by
$\gamma f(z)=\gamma(z)f(z)$ since $\Spin_{2n}$ does not commute with the complex structure
on $\IC^{n}$. The latter action is easily recognised to be that of the Ramond sector
whose quantisaton leads to the two level 1 spin representations of $L\Spin_{2n}$}.
Denote by $\Gamma_{1}$ the corresponding projective representation of $L\Spin_{2n}$ on
$\F=\ext\H_{J}$. More generally, $L\Spin_{2n}$ acts on
$\H\otimes\IC^{\ell}$ by $u\rightarrow u\otimes 1$ and commutes, up to Hilbert--Schmidt
operators with $\J\otimes 1$. Since
$\ext(\H\otimes\IC^{\ell}_{\J\otimes 1})\cong\ext\H_{\J}^{\otimes\ell}=
 \F^{\otimes\ell}$,
the resulting projective representation $\Gamma_{\ell}$ is equivalent to
$\Gamma_{1}^{\otimes\ell}$.\\

Let, as in \ref{ss:modular fermions}, $\AA(I)$ be the von Neumann algebra generated
in $\B(\F^{\otimes\ell})$ by the $c(f)$ with $f$ supported in $I$.
Notice that $\Gamma_{\ell}(L_{I}\Spin_{2n})$ commutes with $\AA(I^{c})$ since, by
(\ref{ch:fermionic}.\ref{eq:isom 2}) and (\ref{ch:fermionic}.\ref{eq:real implementation}),
\begin{equation}\label{eq:commutation}
\Gamma_{\ell}(\gamma)c(f)\Gamma_{\ell}(\gamma)^{*}=
\half{1}\Gamma_{\ell}(\gamma)(\psi(f)-i\psi(if))\Gamma_{\ell}(\gamma)^{*}=
\half{1}(\psi(\gamma f)-i\psi(\gamma if))=
c(f)
\end{equation}
whenever $\gamma\in L_{I}\Spin_{2n}$ and $f$ is supported in $I^{c}$. Let
$\Delta_{I}^{it}$ be the modular group for $\AA(I)$ relative to the vacuum
$\Omega^{\otimes\ell}$, then

\begin{lemma}\label{normalised}
$\Gamma_{\ell}(L_{I}\Spin_{2n})''\subset \AA(I)$ and is normalised by $\Delta_{I}^{it}$.
\end{lemma}
\proof
We claim that $L\Spin_{2n}$ acts by even operators, {\it i.e.}~that $\Gamma_{\ell}$
commutes with the $\IZ_{2}$-grading $\epsilon$ on $\F^{\otimes\ell}$. Indeed,
$\epsilon$ is the canonical quantisation of multiplication by $-1$ on $\H\otimes\IC^{\ell}$
and therefore
\begin{equation}
\Gamma_{\ell}(\gamma)\epsilon\Gamma_{\ell}(\gamma)^{*}=
\chi(\gamma)\epsilon
\end{equation}
where $\chi(\gamma)=\pm 1$ and depends continuously on $\gamma$. By connectedness 
of $L\Spin_{2n}$, $\chi\equiv 1$. Now, by \eqref{eq:commutation},
$\Gamma_{\ell}(L_{I}\Spin_{2n})$ commutes with $\AA(I^{c})$
and therefore, by theorem \ref{modular} lies in
$\AA(I^{c})'=\kappa^{-1}\AA(I)\kappa$. Thus, since
$[\epsilon,\Gamma_{\ell}(L_{I}\Spin_{2n})]=0$,
we get
\begin{equation}
\Gamma_{\ell}(L_{I}\Spin_{2n})''=
\kappa\Gamma_{\ell}(L_{I}\Spin_{2n})''\kappa^{-1}\subset
\AA(I)
\end{equation} 
To see that $\Gamma_{\ell}(L_{I}\Spin_{2n})''$ is normalised by $\Delta_{I}^{it}$,
notice that if $\gamma\in L_{I}\Spin_{2n}$ and $A^{-1}\in SU(1,1)$ leaves
$I$ invariant, then in $PU(\F^{\otimes\ell})$,
$\Gamma(A)\Gamma_{\ell}(\gamma)\Gamma(A)^{*}=\Gamma(A\gamma A^{-1})$. Now $A\gamma A^{-1}$
is multiplication by
\begin{equation}
\left(\frac{z}{\overline{\alpha}z-\beta}\right)
w^{\half{1}}\gamma(A^{-1}z)w^{-\half{1}}
\left(\frac{z}{\overline{\alpha}z-\beta}\right)^{-1}
=
z^{\half{1}}
\left(\frac{z}{\overline{\alpha}z-\beta}\right)
\left(\frac{w}{z}\right)^{\half{1}}
\gamma(A^{-1}z)
\left(\frac{w}{z}\right)^{-\half{1}}
\left(\frac{z}{\overline{\alpha}z-\beta}\right)^{-1}
z^{-\half{1}}
\end{equation}
where $w=A^{-1}z$. This is the action of a loop in $L_{I}\Spin_{2n}$ since $SO_{2n}$
is a normal subgroup of $M_{2n}(\IR)$ and therefore, by the explicit geometric form
of the modular group given by theorem \ref{modular}, $\Gamma_{\ell}(L_{I}\Spin_{2n})''$
is invariant under $\Delta_{I}^{it}$ \halmos\\

\remark Notice that the inclusion $\pi(L_{I}\Spin_{2n})''\subset \AA(I)$ and the evenness
of the operators $\Gamma_{\ell}(\gamma)$, $\gamma\in L\Spin_{2n}$ imply that
$\Gamma_{\ell}(L_{I^{c}}\Spin_{2n})''\subset\Gamma_{\ell}(L_{I}\Spin_{2n})'$ and therefore
give another proof of locality.

\begin{lemma}\label{vacuum}
$\overline{\Gamma_{\ell}(L_{I}\Spin_{2n})''\Omega^{\otimes\ell}}\subset\F^{\otimes\ell}$
is the vacuum representation of $L\Spin_{2n}$ at level $\ell$.
\end{lemma}
\proof
Let $\K=\overline{\Gamma_{\ell}(L\Spin_{2n})''\Omega^{\otimes\ell}}$. $\K$ is
an irreducible $L\Spin_{2n}\rtimes\rot$ module since if $\mathcal W\subset K$ is a submodule
and $P$ is the corresponding orthogonal projection then
$\mathcal W=\overline{\Gamma_{\ell}(L\Spin_{2n})''P\Omega^{\otimes\ell}}$. Since $P$ commutes
with $\rot$, $P\Omega^{\otimes\ell}=\lambda\Omega^{\otimes\ell}$ with $\lambda\in\{0,1\}$
since the latter is the only vector in $\F^{\otimes\ell}$ fixed by $\rot$
\footnote{We are referring to the action of $\rot$ transported from $\HNS$ {\bf not} to the
natural one on $\H$}. It follows that $\mathcal W=0$ or $\K$. 
Since $\K$ is of positive energy, it is uniquely characterised by its lowest energy
subspace and therefore isomorphic to the vacuum representation of level $\ell$. The claim
now follows from the Reeh--Schlieder theorem \halmos

\begin{proposition}[Haag duality in the vacuum sector]\label{haag duality}
Let $(\pi,\H_{0})$ be the vacuum representation at level $\ell$. Then,
\begin{equation}
\pi_{0}(L_{I}\Spin_{2n})'=
\pi_{0}(L_{I^{c}}\Spin_{2n})''
\end{equation}
\end{proposition}
\proof
By lemma \ref{normalised} and Takesaki devissage, the modular conjugation $J_{I}$
of $\AA(I)$ restricts to the one for $\Gamma_{\ell}(L_{I}\Spin_{2n})''$ on
$\overline{\Gamma_{\ell}(L_{I}\Spin_{2n})''\Omega^{\otimes\ell}}\subset\F^{\otimes\ell}$
which, by lemma \ref{vacuum} is isomorphic to $\H_{0}$. When $I$ is the upper
semi--circle, Haag duality follows from (iii) of theorem \ref{modular} because,
in $PU(\F^{\otimes\ell})$,
\begin{equation}
J\Gamma_{\ell}(\gamma)J=
 \kappa^{-1}\Gamma(j)\Gamma_{\ell}(\gamma)\Gamma(j)^{*}\kappa=
 \kappa^{-1}\Gamma(j\gamma j)\kappa
\end{equation}
and $j\gamma j$ is multiplication by $z^{\half{1}}\gamma(z^{-1})z^{-\half{1}}$
which lies in $L_{I^{c}}\Spin_{2n}$. For a general $I$, the result follows by
conjugating by an appropriate element of $SU(1,1)$ \halmos\\

\remark Since conjugation by discontinuous loops normalises each $L_{I}\Spin_{2n}$,
Haag duality holds for any representation obtained by conjugating $\pi_{0}$ by an
element of $L_{Z(\Spin_{2n})}\Spin_{2n}$ and in particular for all level 1
representations of $L\Spin_{2n}$.

\ssubsection{Local equivalence and factoriality}\label{ss:local equivalence}

\begin{proposition}[Local equivalence]\label{local equivalence}
All irreducible positive energy representations of level $\ell$ are unitarily
equivalent for the local loop groups $L_{I}\Spin_{2n}$.
\end{proposition}
\proof
For $\ell=1$, local equivalence follows by conjugating by localised discontinuous
loops. Indeed, by corollary \ref{ch:classification}.\ref{co:level 1 Z}, the action 
of $L_{Z(\Spin_{2n})}\Spin_{2n}$ on the irreducible representations of level 
1 is transitive since $\Spin_{2n}$ is simply--laced. Thus, if $\H_{1}$ is
obtained from $\H_{2}$ by conjugating by discontinuous loops in a given
$L\Spin_{2n}$--coset, choosing a representative equal to 1 on $I$ yields a
unitary equivalence of $\H_{1}$ and $\H_{2}$ as $L_{I}\Spin_{2n}$--modules.

For a general $\ell$, we shall need the fact that $M=\Gamma_{\ell}(L_{I}\Spin_{2n})''$
is a factor of type III$_{1}$. Indeed, by lemma \ref{normalised} and Takesaki devissage,
$\Delta^{it}_{I}$ is the modular group of $M$ relative to $\Omega^{\otimes\ell}$
which, by (iv) of theorem \ref{modular} acts ergodically on $M$ thus proving our
claim. Now, if $\pi_{0}$ is the level 1 vacuum representation of $L\Spin_{2n}$,
then $N=\pi_{0}^{\otimes\ell}(L_{I}\Spin_{2n})''$ is also a type III$_{1}$ factor.
Indeed, $\pi_{0}^{\otimes\ell}$ is a subrepresentation of $\Gamma_{\ell}$ so that
if $p\in M'$ is the corresponding orthogonal projection, $N=pMp=Mp$. Moreover, the
map $x\rightarrow xp$ is an isomorphism since, by the Reeh--Schlieder theorem, any
finite energy vector  $\xi\in\pi_{0}^{\otimes\ell}$ is separating for $M$ and
therefore $xp=0$ implies $x\xi=x(1-p)\xi=0$. Thus,
$\pi_{0}^{\otimes\ell}(L_{I}\Spin_{2n})''$ is a factor of type III and therefore
so is $N'$. It follows that all subrepresentations of $\pi_{0}^{\otimes\ell}$ are
isomorphic since they correspond to the projections in $N'$ and these are all
equivalent.

To conclude, note that by local equivalence at level 1, $\pi_{0}^{\otimes\ell}$
is locally equivalent to $\pi_{i_{1}}\otimes\cdots\otimes\pi_{i_{\ell}}$ where the
$\pi_{i_{j}}$ are any level 1 representations. By proposition
\ref{ch:classification}.\ref{pr:per l lemma}, any irreducible positive energy
representation of $L\Spin_{2n}$ at level $\ell$ is a summand of some such tensor
product and our claim is established \halmos\\

The proof of proposition \ref{local equivalence} has the following important

\begin{corollary}[Factoriality]\label{Factoriality}
The von Neumann algebras generated in a positive energy representation by
the local loop groups $L_{I}\Spin_{2n}$  are factors of type III$_{1}$.
\end{corollary}

\remark\begin{enumerate}
\item
Proposition \ref{local equivalence} implies that the local von Neumann algebras
$\pi_{i}(L_{I}\Spin_{2n})''$ generated in two positive energy representations
$(\pi_{i},\H_{i})$ of equal level are spatially and canonically isomorphic.
Indeed, it yields the existence of a unitary
$U:\H_{1}\rightarrow\H_{2}$ satisfying $U\pi_{1}(\gamma)U^{*}=\pi_{2}(\gamma)$ in
$PU(\H_{2})$ for any $\gamma\in L_{I}\Spin_{2n}$. If $U'$ is another such unitary
and $V={U'}^{*}U$, then, in $U(\H_{1})$,
\begin{equation}
V\pi_{1}(\gamma)V^{*}=\cchi(\gamma)\pi_{1}(\gamma)
\end{equation}
where $\chi(\gamma)\in\T$ is a character of $L_{I}\Spin_{2n}$ and is therefore 
trivial by lemma \ref{local is perfect}. Thus, in $U(\H_{2})$
$U\pi_{1}(\gamma)U^{*}=U'V\pi_{1}(\gamma)V^{*}{U'}^{*}=U'\pi_{1}(\gamma){U'}^{*}$
and
\begin{equation}\label{eq:spatial}
x\longrightarrow UxU^{*}
\end{equation}
is a canonical spatial isomorphism
$\pi_{1}(L_{I}\Spin_{2n})''\cong\pi_{2}(L_{I}\Spin_{2n})''$.
\item
Notice that \eqref{eq:spatial} yields an isomorphism of the central extensions of
$L_{I}\Spin_{2n}$ induced by $\pi_{1}$ and $\pi_{2}$. Since $\pi_{i}(L_{I}\Spin_{2n})$
and $\pi_{i}(L_{I^{c}}\Spin_{2n})$ commute by locality, this isomorphism extends to
one of the central extensions of $L_{I}\Spin_{2n}\times L_{I^{c}}\Spin_{2n}$
corresponding to $\pi_{1}$ and $\pi_{2}$. Using proposition \ref{pr:sobolev density},
it is easy to show that it extends futher to one of central extensions of $L\Spin_{2n}$
thus giving an alternative proof of (ii) of proposition
\ref{ch:analytic}.\ref{classify extension}.
\end{enumerate}



\renewcommand {\SS}{\mathcal S}
\newcommand {\reg}{\mathcal R}
\newcommand {\ns}{\operatorname{NS}}
\newcommand {\ra}{\operatorname{R }}
\renewcommand {\l}{\mathfrak l}

\chapter{Bosonic construction of level 1 representations and \\primary fields}
\label{ch:vertex operator}

This chapter is devoted to the bosonic or Frenkel--Kac--Segal vertex operator
construction of the level 1 representations of $L\Spin_{2n}$. This stems from
the {\it a posteriori} observation that they remain irreducible when restricted
to the abelian group $LT$ of loops in a maximal torus $T\subset\Spin_{2n}$. A
suitable Stone--von Neumann theorem may then be used to reconstruct them from
the Heisenberg representations of $LT$.
Our interest in this construction lies in the fact that it gives a particularly
convenient description of all level 1 primary fields which will be used in
chapter \ref{ch:sobolev fields} to establish their continuity properties as
operator--valued distributions.\\

The vertex operator construction applies in fact to the loop groups of
all compact, simply--connected, simple and simply--laced Lie groups $G$ and we
shall treat it in this generality.
In section \ref{se:restriction}, we study the restriction of positive energy
representations of $LG$ to the loop group $LT=C^{\infty}(S^{1},T)$ of a maximal
torus $T\subset G$ and compute the corresponding projective cocycle on $LT$. As
is well--known, the associated central extension factors as a product of two
Heisenberg groups, one corresponding to the identity component of $LT/T$, the
other to the product of the constant loops $T\subset LT$ by the integral lattice
$\check T=\Hom(\T,T)$. Suitable versions of the Stone--von Neumann theorem classifying
the irreducible representations of these groups are obtained in section
\ref{se:level 1 LT}.\\

We also compute in section \ref{se:restriction} the commutation relations of $LT$
with the infinitesimal action of $\lpol\g$ in a positive energy representation.
In section \ref{se:level 1 LT}, we construct operators in the irreducible
representations of $LT$ which mimick $\lpol\g$ in that they have the same commutation
relations with $LT$. We then show in section \ref{se:vertex operators} that they give
level 1 representations of $\lpol\g$. The reader familiar with the vertex operator
construction may skip to section \ref{se:level 1 primary} where the construction
of the level 1 primary fields is carried out by following a similar scheme. The
equivariance properties of the primary fields with respect to $LT$ are obtained
and operators satisfying these constructed explicitly which, by inspection have
the required commutation relations with $\lpol\g$.

\ssection{Formal variables}
\label{se:formal}

We begin by briefly reviewing the formal variable approach of Frenkel, Lepowski and
Meurman \cite[chap. 2]{FLM}. This gives a convenient way to encode the infinitesimal
action of $\lpol\g$ on the finite energy subspace $\hfin$ of a positive energy
representation of $LG$.\\

Let $\ell$ be the level of $\H$.
The commutation relations satisfied by $\lpol\g$ on $\hfin$, namely
\begin{equation}
[X(m),Y(n)]=[X,Y](m+n)+\ell m\delta_{m+n,0}\<X,Y\>
\end{equation}
may equivalently be written in terms of the generating
function $X(z)\in\End(\hfin)[[z,z^{-1}]]$ defined by
\begin{equation}
X(z)=\sum_{n\in\IZ}X(n)z^{-n}\label{eq:X(z)}
\end{equation}
as
\begin{equation}
\begin{split}
[X(z),Y(\zeta)]
&=\medsum_{m,n}[X,Y](m+n)\zeta^{-(m+n)}\Bigl(\frac{\zeta}{z}\Bigr)^{m}+
\ell\<X,Y\>\medsum_{m}m\Bigl(\frac{\zeta}{z}\Bigr)^{m}\\
&=[X,Y](\zeta)\delta\Bigl(\frac{\zeta}{z}\Bigr)+
\ell\<X,Y\>\delta'\Bigl(\frac{\zeta}{z}\Bigr)
\end{split}
\end{equation}
where
\begin{xalignat}{3}
\delta(u) &=\sum_{n\in\IZ} u^{n}=\delta(u^{-1})&
&\text{and}&
\delta'(u)&=\sum_{n\in\IZ} nu^{n}=-\delta'(u^{-1})
\end{xalignat}
are the formal Dirac delta function and its first derivative.
Similarly, the relation $[d,X(n)]=-nX(n)$ and the formal adjunction property
$X(n)^{*}=-\overline{X}(-n)$ are equivalent to
\begin{align}
[d,X(z)]&=z\frac{d}{dz}\label{eq:homogeneous}\\
X(z)^{*}&=-\overline{X}(z)
\end{align}
where the adjoint of $A(z)=\sum_{n\in\IZ}A(n)z^{-n}$ is defined whenever the modes
$A(n)$ possess formal adjoints by $\sum_{n\in\IZ}A(n)^{*}z^{n}$ so that it satisfies
\eqref{eq:homogeneous} if, and only if $A(z)$ does. We will need, for later reference,
the following elementary 

\begin{lemma}\label{le:delta}
Let $u$ be a formal variable and define, for any $a\in\IC$
\begin{equation}\label{eq:binomial}
(1-u)^{a}=\sum_{n\geq 0}\begin{pmatrix}a\\n\end{pmatrix}(-u)^{n}
\end{equation}
Then
\begin{align}
\delta(u)&=(1-u)^{-1}+(1-u^{-1})^{-1}u^{-1}
\label{eq:delta}\\
\delta'(u)&=u(1-u)^{-2}-(1-u^{-1})^{-2}u^{-1}
\label{eq:delta'}
\end{align}
\end{lemma}

\begin{lemma}\label{le:delta identity}
Set $u=\frac{\zeta}{z}$ for two formal variables $\zeta$ and $z$ and regard $\delta(u)$
and $\delta'(u)$ as formal Laurent series in $\zeta,z$.
Then, for any $f(z,\zeta)\in\IC[[z,\zeta,z^{-1},\zeta^{-1}]]$ such that
$f(\zeta,\zeta)$ is a well--defined formal Laurent series,
\begin{align}
f(z,\zeta)\delta\Bigl(\frac{\zeta}{z}\Bigr)&=
f(\zeta,\zeta)\delta\Bigl(\frac{\zeta}{z}\Bigr)=
f(z,z)\delta\Bigl(\frac{\zeta}{z}\Bigr)
\label{eq:delta identity}\\
f(z,\zeta)\delta'\Bigl(\frac{\zeta}{z}\Bigr)&=
f(\zeta,\zeta)\delta'\Bigl(\frac{\zeta}{z}\Bigr)+
z\frac{\partial f}{\partial z}(\zeta,\zeta)\delta\Bigl(\frac{\zeta}{z}\Bigr)=
f(z,z)\delta'\Bigl(\frac{\zeta}{z}\Bigr)-
\zeta\frac{\partial f}{\partial\zeta}(z,z)\delta\Bigl(\frac{\zeta}{z}\Bigr)
\label{eq:delta' identity}
\end{align}
\end{lemma}

\remark It is important to assume that $f(z,\zeta)$ only contains integral powers of
$z,\zeta$ in the above lemma since $f(z,\zeta)=\bigl(\frac{\zeta}{z}\bigr)^{\half{1}}$
clearly does not satisfy \eqref{eq:delta identity}.

\ssection{Restriction of positive energy representations to $LT$}
\label{se:restriction}

Unlike the quark model of chapter \ref{ch:fermionic}, the vertex operator construction
is purely infinitesimal and gives an explicit description of the action of $\lpol\g$
on the finite energy subspaces of level 1 representations only. These subspaces are
studied as $LT$--modules in this section and will be reconstructed as such in section
\ref{se:vertex operators}. Strictly speaking, they are not invariant under the identity
component of $LT/T$ and we shall trade the latter for its algebraic Lie algebra
$\lpol\t/\t$ consisting of polynomial loops in $\t=\Lie(T)$ with zero average.

\ssubsection{The loop group $LT$}\label{ss:LT}

Let $T\subset G$ be a maximal torus with Lie algebra $\t$ and integral lattice
$I=\{h\in\t|\thinspace\exp_{T}(2\pi h)=1\}$. We denote by $LT=C^{\infty}(S^{1},T)$
the loop group of $T$. Since 
\begin{equation}
0\rightarrow 2\pi I\rightarrow\t\xrightarrow{\exp} T\rightarrow 1
\end{equation}
is the universal covering group of $T$, we have
\begin{equation}
LT\cong\{f\in C^{\infty}([0,2\pi],\t) | f(2\pi)-f(0)\in 2\pi I\}/2\pi I=
{\medsqcup_{\lambda\in I}}LT_{\lambda}
\end{equation}
where $LT_{\lambda}$ is the connected component of loops $f$ whose
winding number $\nu_{f}=\frac{1}{2\pi}\Bigl(f(2\pi)-f(0)\Bigr)$
equals $\lambda\in I$. We therefore have an isomorphism
$LT\cong L\t/\t\times T\times\check T$, where $L\t/\t$
is the additive group of smooth loops in $\t$ of zero average and
$\check T=\Hom(\T,T)$, given by
\begin{equation}
f\rightarrow\Bigl(f-f_{0}-\nu_{f}(\theta-\pi),\exp_{T}(f_{0}),
 \exp_{T}(\theta\nu_{f})\Bigr)
\end{equation}
where $f_{0}=\frac{1}{2\pi}\int_{0}^{2\pi}f(\theta)d\theta$.
Since $G$ is simply--connected, multiplication by $i$ yields an identification of
the integral lattice $I$ with the coroot lattice $\coroot\subset i\t$
\cite[thm. 5.47]{Ad}. The simple--lacedness of $G$ implies in turn that the latter
is isomorphic to the root lattice $\root\subset i\t^{*}$ via the identification
$i\t\cong i\t^{*}$ defined by the basic inner product $\<\cdot,\cdot\>$. Using these
identifications, we shall parametrise elements of $\check T$ by associating to each
$\alpha\in\root$ the homomorphism
\begin{equation}
\zeta_{\alpha}(\theta)=\exp_{T}(-i\alpha\theta)
\end{equation}

\ssubsection{Commutation relations of $LT$ and $\lpol\g$}
\label{ss:LT on finite energy}

Let $(\pi,\H)$ be a level $\ell$ positive energy representation of $LG$ and denote
by $U_{\theta}=e^{id\theta}$ the corresponding integrally--moded action of $\rot$. 
$\pi$ lifts uniquely over $G\subset LG$ to a unitary representation commuting with
$U_{\theta}$ which we denote by the same symbol.
Recall from chapter \ref{ch:analytic} that the action of $LG$ leaves the subspace
of smooth vectors $\hsmooth$ invariant and, by theorem
\ref{ch:analytic}.\ref{th:segal formulae} satisfies
\begin{align}
\pi(\gamma)\pi(X)\pi(\gamma)^{*}&=
\pi(\gamma X\gamma^{-1})-
i\ell\int_{0}^{2\pi}\<\gamma^{-1}\dot\gamma,X\>\frac{d\theta}{2\pi}
\label{eq:crossed one}\\
\pi(\gamma)d\pi(\gamma)^{*}&=
d-i\pi(\dot\gamma\gamma^{-1})-
\half{\ell}\int_{0}^{2\pi}
\<\gamma^{-1}\dot\gamma,\gamma^{-1}\dot\gamma\>\frac{d\theta}{2\pi}
\label{eq:crossed two}
\end{align}
for any $\gamma\in LG$ and $X\in L\g$.

\begin{proposition}\label{pr:zeromodes energy}
The action of $T\times\check T$ leaves the subspace of finite energy vectors $\hfin$
invariant and satisfies, for any $\tau\in T$ and $\alpha\in\root$
\begin{align}
[d,\pi(\tau)]&=0
\label{eq:torus energy}\\
[d,\pi(\zeta_{\alpha})]&=\pi(\zeta_{\alpha})(\pi(\alpha)+\half{\<\alpha,\alpha\>})
\label{eq:root energy}
\end{align}
\end{proposition}
\proof
The action of $T$ leaves $\hfin$ invariant since $\pi(g)U_{\theta}\pi(g)^{*}=U_{\theta}$
for any $g\in G$. If $\alpha\in\root$ then,
projectively
\begin{equation}\label{eq:projective commutation}
\pi(\zeta_{\alpha})U_{\theta}\pi(\zeta_{\alpha})^{*}=
\pi(\zeta_{\alpha}(\zeta_{\alpha})_{\theta}^{-1})U_{\theta}=
\pi(\exp_{T}(-i\alpha\theta))U_{\theta}
\end{equation}
Thus, $\theta\rightarrow\pi(\exp_{T}(-i\alpha\theta))U_{\theta}$ is a positive energy
action of $\rot$. Since it commutes with $U_{\theta}$, their finite energy subspaces
coincide and therefore $\pi(\zeta_{\alpha})\hfin=\hfin$. The relations
\eqref{eq:torus energy}--\eqref{eq:root energy} follow at once from
\eqref{eq:crossed two} \halmos

\begin{proposition}\label{pr:zeromodes affine}
The action of $T\times\check T$ on $\hfin$ satisfies, for any $X\in\lpol\g$
\begin{align}
\pi(\tau)\pi(X)\pi(\tau)^{*}&=\pi(\tau X\tau^{-1})
\label{eq:torus affine}\\
\pi(\zeta_{\alpha})\pi(X)\pi(\zeta_{\alpha})^{*}&=
\pi(\zeta_{\alpha}X\zeta_{\alpha}^{-1})-
\ell\int_{0}^{2\pi}\<\alpha,X\>\frac{d\theta}{2\pi}
\label{eq:root affine}
\end{align}
In particular,
\begin{equation}
\pi(\tau)\pi(\zeta_{\alpha})\pi(\tau)^{*}=
\alpha(\tau)^{\ell}\pi(\zeta_{\alpha})
\label{eq:root torus}
\end{equation}
where $\alpha(\exp_{T}(h))=e^{\<\alpha,h\>}$.
Moreover, if $X_{\beta}\in\g_{\beta}$ and $X_{\beta}(z)$ is given by
\eqref{eq:X(z)}, then
\begin{align}
\pi(\tau)X_{\beta}(z)\pi(\tau)^{*}&=\beta(\tau)X_{\beta}(z)
\label{eq:torus vertex}\\
\pi(\zeta_{\alpha})X_{\beta}(z)\pi(\zeta_{\alpha})^{*}&=
z^{-\<\alpha,\beta\>}X_{\beta}(z)
\label{eq:root vertex}
\end{align}
\end{proposition}
\proof
The relations \eqref{eq:torus affine}, \eqref{eq:root affine} and
\eqref{eq:torus vertex} follow at once from \eqref{eq:crossed one}. Since $T$ and
$\check T$ commute projectively, the following holds in $U(\H)$ for any $h\in\t$,
\begin{equation}
\pi(\zeta_{\alpha})e^{t\pi(h)}\pi(\zeta_{\alpha})^{*}=e^{tc}e^{t\pi(h)}
\end{equation}
for some $c\in i\IR$. Applying both sides to $\xi\in\hfin$ and taking derivatives
at $t=0$ yields, by comparison with \eqref{eq:root affine}, $c=-\ell\<\alpha,h\>$
and therefore \eqref{eq:root torus}. Let now
$X_{\beta}\in\g_{\beta}$ so that $[h,X_{\beta}]=\<h,\beta\>X_{\beta}$ for any
$h\in\t$. Then, $\Ad(\exp_{T}(h))X_{\beta}=e^{\<h,\beta\>}X_{\beta}$ and therefore,
in $L\g$,
\begin{equation}
 \zeta_{\alpha}X_{\beta}(n)\zeta_{\alpha}^{-1}(\theta)=
 e^{n-i\<\alpha,\beta\>\theta}X_{\beta}=
 X_{\beta}(n-\<\alpha,\beta\>)
\end{equation}
Since $\zeta_{\alpha}^{-1}\dot\zeta_{\alpha}=-i\alpha\in\t$ and this subspace is
orthogonal to $\g_{\beta}$ with respect to the Killing form, it follows from
\eqref{eq:crossed one} that 
\begin{equation}
\pi(\zeta_{\alpha})\pi(X_{\beta}(n))\pi(\zeta_{\alpha})^{*}=
\pi(X_{\beta}(n-\<\alpha,\beta\>))
\end{equation}
and therefore that \eqref{eq:root vertex} holds \halmos

\ssubsection{The projective cocycle on $\check T$}\label{ss:cocycle on LT}

The following result is due to Pressley and Segal \cite[prop. 4.8.1]{PS}

\begin{proposition}\label{pr:coroot cocycle}
Let $(\pi,\H)$ be a positive energy representation of $LG$ at level $\ell$. Then, for
any $\alpha,\beta\in\root$,
\begin{equation}\label{eq:coroot cocycle}
\pi(\zeta_{\alpha})\pi(\zeta_{\beta})\pi(\zeta_{\alpha})^{*}\pi(\zeta_{\beta})^{*}=
(-1)^{\ell\<\alpha,\beta\>}
\end{equation}
\end{proposition}

We shall derive proposition \ref{pr:coroot cocycle} from a somewhat more general result
which will be needed in section \ref{se:level 1 primary}.
Recall that the isomorphism $i\t^{*}\cong i\t$ corresponding to the basic inner product
$\<\cdot,\cdot\>$ identifies the weight and coweight lattices $\weight,\coweight$ since
$G$ is simply--laced. In particular, any $\mu\in\weight$ defines, in the terminology of
section \ref{ch:classification}.\ref{se:disc loops}, a discontinuous
loop $\zeta_{\mu}(\theta)=\exp_{T}(-i\mu\theta)$ which acts on $L\g$ by conjugation.

\begin{proposition}\label{pr:extended cocycle}
Let $(\pi,\H)$ be a level $\ell$ positive energy representation of $LG$. Assume that
the finite energy subspace $\hfin\subset\H$ supports a projective unitary representation
$\rho$ of an intermediate lattice $\root\subset\Lambda\subset\weight$ satisfying,
for any $X\in\lpol\g$
\begin{equation}\label{eq:infi cocycle}
\rho(\mu)\pi(X)\rho(\mu)^{*}=
\pi(\zeta_{\mu}X\zeta_{\mu}^{-1})-\ell\int_{0}^{2\pi}\<\mu,X\>\frac{d\theta}{2\pi}
\end{equation}
Then, for any $\alpha\in\root$ and $\mu\in\Lambda$
\begin{equation}\label{eq:extended cocycle}
\rho(\mu)\pi(\zeta_{\alpha})\rho(\mu)^{*}\pi(\zeta_{\alpha})^{*}=
(-1)^{\ell\<\mu,\alpha\>}
\end{equation}
\end{proposition}

{\sc Proof of proposition \ref{pr:coroot cocycle}}.
Let $\Lambda=\root$ in proposition \ref{pr:extended cocycle} and $\rho$ the restriction
of $\pi$ to $\root\cong\check T$. By proposition \ref{pr:zeromodes energy}, $\rho$ leaves
$\hfin$ invariant and, by \eqref{eq:root affine} satisfies \eqref{eq:infi cocycle}. The
conclusion therefore follows from proposition \ref{pr:extended cocycle} \halmos\\

{\sc Proof of proposition \ref{pr:extended cocycle}}.
Notice first that \eqref{eq:extended cocycle} holds if $\mu=-\alpha$. Indeed, by
proposition \ref{pr:zeromodes energy} and \eqref{eq:root affine}, the unitary
$T_{\alpha}=\pi(\zeta_{\alpha})\rho(-\alpha)$
leaves $\hfin$ invariant and commutes with the action of $\lpol\g$. By proposition
\ref{ch:classification}.\ref{equivalent} it projectively commutes with $LG$ and therefore,
for any $\gamma\in LG$, the following holds in $U(\H)$
\begin{equation}
T_{\alpha}\pi(\gamma)T_{\alpha}^{*}=\cchi(\gamma)\pi(\gamma)
\end{equation}
for some $\cchi(\gamma)\in\T$ which defines a character of $LG$. By lemma
\ref{ch:loc loops}.\ref{local is perfect}, $\Hom(LG,\T)=1$ and therefore $\cchi\equiv 1$.
In particular,
\begin{equation}
\rho(-\alpha)\pi(\zeta_{\alpha})\rho(-\alpha)^{*}\pi(\zeta_{\alpha})^{*}=
\rho(-\alpha)\pi(\zeta_{\alpha})T_{\alpha}^{*}=
\rho(-\alpha)T_{\alpha}^{*}\pi(\zeta_{\alpha})=
1
\end{equation}
which coincides with the right--hand side of \eqref{eq:extended cocycle} since $\root$
is an even lattice when $G$ simply--laced.\\

We now establish \eqref{eq:extended cocycle} when $\alpha$ is a positive root and
$\mu\in\Lambda$ is such that $\<\mu,\alpha\>\in\{0,1\}$.
The homomorphism $\zeta_{\alpha}(\theta)=\exp(-i\theta\alpha)$ may
be written as a product of two exponentials in $LG$ \cite[4.8.1]{PS}, namely
\begin{equation}\label{product of exp}
\zeta_{\alpha}=
\exp_{LG}\Bigl(-\half{\pi}(e_{\alpha}(0)-f_{\alpha}(0) )\Bigr)
\exp_{LG}\Bigl( \half{\pi}(e_{\alpha}(1)-f_{\alpha}(-1))\Bigr)
\end{equation}

To see this, consider the loop $\sigma_{\alpha}:SU_{2}\rightarrow G$ corresponding to
the $\mathfrak{sl}_{2}(\IC)$--subalgebra of $\gc$ spanned by
$e_{\alpha},f_{\alpha},h_{\alpha}=\alpha$. This
induces a homomorphism $LSU_{2}\rightarrow LG$ in the obvious way mapping
\begin{equation}
\theta\rightarrow\begin{pmatrix}e^{-i\theta}&0\\0&e^{i\theta}\end{pmatrix}
\end{equation}
to $\zeta_{\alpha}$. Using $\sigma_{\alpha}$ and the standard basis of
$\mathfrak{sl}_{2}(\IC)$ given by
\begin{xalignat}{3}
e_{\alpha}&=\begin{pmatrix}0&1\\0&0 \end{pmatrix}&
f_{\alpha}&=\begin{pmatrix}0&0\\1&0 \end{pmatrix}&
h_{\alpha}&=\begin{pmatrix}1&0\\0&-1\end{pmatrix}
\end{xalignat}
\eqref{product of exp} reduces to a matrix check.\\

If $h\in\t$, then $[h,e_{\alpha}]=\<h,\alpha\>e_{\alpha}$ whence
$\Ad(\exp_{T}(h))e_{\alpha}=e^{\<h,\alpha\>}e_{\alpha}$.
Therefore, since $\zeta_{\mu}(\theta)=\exp_{T}(-i\theta\mu)$, we have
$\zeta_{\mu}e_{\alpha}(n)\zeta_{\mu}^{-1}(\theta)=
 e_{\alpha}\otimes e^{i\theta(n-\<\alpha,\mu\>)}=
 e_{\alpha}(n-\<\alpha,\mu\>)(\theta)$.
In other words,
\begin{align}
\zeta_{\mu}e_{\alpha}(n)\zeta_{\mu}^{-1}&=e_{\alpha}(n-\<\alpha,\mu\>)\\
\zeta_{\mu}f_{\alpha}(n)\zeta_{\mu}^{-1}&=f_{\alpha}(n+\<\alpha,\mu\>)
\end{align}
Since $-i\mu\in\t$ and this subspace is orthogonal to $\g_{\alpha}\oplus\g_{-\alpha}$
with respect to the Killing form, \eqref{eq:infi cocycle} yields
\begin{align}\label{infi formula}
\rho(\mu)\pi(e_{\alpha}(n))\rho(\mu)^{*}&=\pi(e_{\alpha}(n-\<\alpha,\mu\>))\\
\rho(\mu)\pi(f_{\alpha}(n))\rho(\mu)^{*}&=\pi(f_{\alpha}(n+\<\alpha,\mu\>))
\end{align}

Thus, \eqref{eq:extended cocycle} holds if $\<\alpha,\mu\>=0$. If, on the other hand
$\<\alpha,\mu\>=1$, then by \eqref{product of exp} 
\begin{equation}\label{temp}
\begin{split}
\rho(\mu)    \pi(\zeta_{\alpha})
\rho(\mu)^{*}\pi(\zeta_{\alpha})^{*}
=&
\exp\Bigl(-\half{\pi}\pi(e_{\alpha}(-1)-f_{\alpha}(1) )\Bigr)
\exp\Bigl( \half{\pi}\pi(e_{\alpha}(0) -f_{\alpha}(0) )\Bigr)\\
\cdot&
\exp\Bigl(-\half{\pi}\pi(e_{\alpha}(1) -f_{\alpha}(-1))\Bigr)
\exp\Bigl( \half{\pi}\pi(e_{\alpha}(0) -f_{\alpha}(0) )\Bigr)\\
=&
\exp\Bigl(-\half{\pi}\pi(e_{\alpha}(-1)-f_{\alpha}(1) )\Bigr)\\
\cdot&
\exp\Bigl( \half{\pi}\pi(e_{\alpha}(0) -f_{\alpha}(0) )\Bigr)
\exp\Bigl(-\half{\pi}\pi(e_{\alpha}(1) -f_{\alpha}(-1))\Bigr)
\exp\Bigl( \half{\pi}\pi(e_{\alpha}(0) -f_{\alpha}(0) )\Bigr)^{*}\\
\cdot&
\exp\Bigl(       \pi \pi(e_{\alpha}(0) -f_{\alpha}(0) )\Bigr)
\end{split}
\end{equation}
where all exponentials are defined using the spectral theorem. As is easily checked
using $\sigma_{\alpha}$, we have
\begin{equation}
\exp_{LG}\Bigl( \half{\pi}(e_{\alpha}(0) -f_{\alpha}(0) )\Bigr)
(e_{\alpha}(1)-f_{\alpha}(-1))
\exp_{LG}\Bigl( \half{\pi}(e_{\alpha}(0) -f_{\alpha}(0) )\Bigr)^{-1}=
e_{\alpha}(-1)-f_{\alpha}(1)
\end{equation}
Moreover, since we are conjugating by a constant loop, no correction term arises from
\eqref{eq:crossed one} and it follows that \eqref{temp} is equal to
\begin{equation}\label{intermediate}
\exp\Bigl(-\pi\pi(e_{\alpha}(-1)-f_{\alpha}(1))\Bigr)
\exp\Bigl( \pi\pi(e_{\alpha}(0) -f_{\alpha}(0))\Bigr)
\end{equation}

To proceed, we seek to diagonalise the above unitaries. This is best done in $LSU_{2}$ using
the identity
\begin{equation}
e_{\alpha}(-m)-f_{\alpha}(m)=V(m)ih_{\alpha}(0)V(m)^{*}
\end{equation}
where $V(m)\in LSU_{2}$ is given by
$\theta\rightarrow\frac{1}{\sqrt{2}}
 \begin{pmatrix}1&ie^{-im\theta}\\ie^{im\theta}&1\end{pmatrix}$.
Since
\begin{equation}
V^{-1}(m)\dot V(m)=\half{m}\Bigl(ih_{\alpha}(0)+e_{\alpha}(-m)-f_{\alpha}(m)\Bigr)
\end{equation}
we have
\begin{equation}
\int_{0}^{2\pi}\<V^{-1}(m)\dot V(m),ih_{\alpha}(0)\>\frac{d\theta}{2\pi}=
-\half{m}\|h_{\alpha}\|^{2}=-m
\end{equation}
since $G$ is simply--laced. Therefore, using \eqref{eq:crossed one}, \eqref{intermediate} is
equal to
\begin{equation}
\begin{split}
 &
e^{i\pi\ell }
\pi(V(1))\exp\Bigl(-\pi\pi(ih_{\alpha}(0))\Bigr)\pi(V(1))^{*}
\pi(V(0))\exp\Bigl( \pi\pi(ih_{\alpha}(0))\Bigr)\pi(V(0))^{*}\\
=&
e^{i\pi\ell }
\pi(V(1))\pi(\exp_{LG}(-\pi ih_{\alpha}(0)))\pi(V(1))^{*}
\pi(V(0))\pi(\exp_{LG}( \pi ih_{\alpha}(0)))\pi(V(0))^{*}
\end{split}
\end{equation}

Since $\exp_{SU_{2}}(-i\pi h_{\alpha}(0))=-1$ lies in the centre of any central extension
of $L\SU_{2}$, the above is equal to $(-1)^{\ell}=(-1)^{\ell\<\alpha,\mu\>}$.\\

To summarise, \eqref{eq:extended cocycle} holds if $\alpha$ is a simple root and
$\mu\in\Lambda$ equals $-\alpha$ or satisfies $\<\mu,\alpha\>=0,1$ and therefore
if $\mu\in\Lambda$ is the opposite of a simple root or a minimal dominant weight.
Since each $\root$ coset in $\Lambda$ contains a minimal dominant weight, such
$\alpha$, $\mu$ span $\root$ and $\Lambda$ respectively and it follows that $\rho$
and the restriction of $\pi$ to $\check T$ projectively commute. Thus,
\begin{equation}
\rho(\mu)\pi(\zeta_{\alpha})\rho(\mu)^{*}\pi(\zeta_{\alpha})^{*}=
\omega(\mu,\alpha)
\end{equation}
where $\omega(\mu,\alpha)\in\T$ is easily seen to be $\IZ$--bilinear in each of its
arguments. Since $\omega(\mu,\alpha)=(-1)^{\ell\<\mu,\alpha\>}$ for a spanning set
of $\mu$ and $\alpha$, \eqref{eq:extended cocycle} follows \halmos

\ssection{Level 1 representations of $LT$}
\label{se:level 1 LT}

We classify below the irreducible projective unitary representations of $LT$,
or more precisely $\lpol\t/\t\times(T\times\check T)$ satisfying the level 1 relations
computed in the previous section and show that, like the level 1 representations of $LG$
they are parametrised by the dual of $Z(G)$. Notice that, by proposition
\ref{pr:zeromodes affine},
\begin{align}
\pi(\tau)\pi(h(n))\pi(\tau)^{*}&=\pi(h(n))
\label{eq:desiderata 4}\\
\pi(\zeta_{\alpha})\pi(h(n))\pi(\zeta_{\alpha})^{*}&=
\pi(h(n))-\delta_{n,0}\<\alpha,h\>
\label{eq:desiderata 5}
\end{align}
so that $\lpol\t/\t$ and $T\times\check T$ commute and we may therefore study the
factors separately. In each case, we introduce the vertex operators of Segal and
Kac--Frenkel. These are formal Laurent series having the same commutation relations
with $\lpol\t/\t\times T\times\check T$ as the generating functions for a Cartan--Weyl
basis of $\lpol\gc$. In section \ref{se:vertex operators}, we shall show that their
modes give, together with the action of $\lpol\t$ all level 1 representations of
$\lpol\g$.

\ssubsection{Stone--von Neumann theorem and vertex operators $T\times\check T$}
\label{ss:Stone for zero}

When restricted to $T\times\check T\subset LT$, a positive energy representation has a
somewhat hybrid nature. On the one hand, it possesses a lift to a unitary representation
of $T$ which is uniquely determined by the requirement that it should extend to one of
$G$. On the other, proposition \ref{pr:coroot cocycle} implies that it is genuinely
projective on $\check T$. We shall accordingly consider projective unitary representations
$\pi$ of $T\times\check T$ with a preferred lift over $T$ which we denote by the same
symbol. These will be required to satisfy \eqref{eq:root torus} and \eqref{eq:coroot
cocycle} for $\ell=1$, namely
\begin{align}
\pi(\tau)\pi(\zeta_{\alpha})\pi(\tau)^{*}&=\alpha(\tau)\pi(\zeta_{\alpha})
\label{eq:desiderata 2}\\
\pi(\zeta_{\alpha})\pi(\zeta_{\beta})\pi(\zeta_{\alpha})^{*}\pi(\zeta_{\beta})^{*}&=
(-1)^{\<\alpha,\beta\>}
\label{eq:desiderata 3}
\end{align}
Two such $(\pi_{i},\H_{i})$, $i=1,2$ will be regarded as unitarily equivalent if there
exists a unitary $U:\H_{1}\rightarrow\H_{2}$ such that
\begin{xalignat}{3}\label{eq:hybrid equivalence}
U\pi_{1}(\tau)U^{*}&=\pi_{2}(\tau)\thickspace\text{in $U(\H_{2})$}&
&\text{and}&
U\pi_{1}(\zeta_{\alpha})U^{*}&=\pi_{2}(\zeta_{\alpha})\thickspace\text{in $PU(\H_{2})$}
\end{xalignat}
Notice that if each $\pi_{i}$ is the restriction of a positive energy representation of $LG$,
any unitary equivalence $U:\H_{1}\rightarrow\H_{2}$ as $LG$--modules necessarily satisfies
\eqref{eq:hybrid equivalence}. Indeed, for any $g\in G$ with canonical lifts
$\pi_{i}(g)\in U(\H_{i})$,
\begin{equation}
U\pi_{1}(g)U^{*}=\cchi(g)\pi_{2}(g)
\end{equation}
for some $\cchi(g)\in\T$ which defines a character of $G$. Since $G$ is simple, $\cchi\equiv 1$
and \eqref{eq:hybrid equivalence} holds.\\

Since our goal is to reconstruct the finite energy subspaces of the level 1 positive energy
representations of $LG$, we shall in fact be interested in projective representations of
$T\times\check T$ on {\it pre}--Hilbert spaces $\H$ and demand accordingly that they be
algebraically irreducible. We shall further assume that they are the algebraic sum of
their $T$--eigenspaces. To construct representations of $T\times\check T$ satisfying
\eqref{eq:desiderata 2} and \eqref{eq:desiderata 3},
we need some elementary facts about cocycles on $\root$. Let $\Gamma,A$ be groups
with $A$ abelian and recall that an $A$--valued {\it cocyle} on $\Gamma$ is a map
$\epsilon:\Gamma\times\Gamma\rightarrow A$ satisfying
\begin{equation}\label{eq:cocycle id}
\epsilon(\lambda,\mu)\epsilon(\lambda\mu,\rho)=
\epsilon(\mu,\rho)\epsilon(\lambda,\mu\rho)
\end{equation}
A {\it coboundary} $b:\Gamma\times\Gamma\rightarrow A$ is a map of the form
\begin{equation}\label{eq:coboundary id}
b(\lambda,\mu)=
\frac{a(\lambda)a(\mu)}{a(\lambda\mu)}
\end{equation}
for some function $a:\Gamma\rightarrow A$. As is readily verified, coboundaries always
satisfy \eqref{eq:cocycle id}. Two cocyles $\epsilon,\epsilon'$ are {\it cohomologous} if they
differ by a coboundary $b$, {\it i.e.}~ if $\epsilon'=\epsilon b$.
If $\Gamma$ is abelian, we may attach to any cocycle $\epsilon$ a {\it commutator map}
$\omega:\Gamma\times \Gamma\rightarrow A$ defined by
$\omega(\lambda,\mu)=\epsilon(\lambda,\mu)\epsilon(\mu,\lambda)^{-1}$. $\omega$ depends
only upon the cohomolgy class of $\epsilon$ and, by \eqref{eq:cocycle id} is bilinear and
skew--symmetric. The following is well--known, see for example \cite[propn. 5.2.3]{FLM}
and \cite[\S 2.3]{FK}

\begin{proposition}\label{pr:commutator map}
If $\Gamma\cong\IZ^{r}$, then
\begin{enumerate}
\item The cocyles $\epsilon$ and $\epsilon'$ are cohomologous if, and only if they have
the same commutator map.
\item For any skew--symmetric, bilinear form $\omega:\Gamma\times\Gamma\rightarrow A$
there exists a cocycle $\epsilon$ whose commutator map is $\omega$. Moreover, $\epsilon$
may be chosen with the following normalisations
\begin{xalignat}{3}\label{eq:normalised cocycle}
\epsilon(\lambda,0)&=1&\epsilon(0,\lambda)&=1&\epsilon(\lambda,-\lambda)&=1
\end{xalignat}
\end{enumerate}
\end{proposition}

Let now $\epsilon$ be a $\T$--valued cocycle on $\root$ with commutator map
$\omega(\alpha,\beta)=(-1)^{\<\alpha,\beta\>}$. We consider the action of $T\times\check T$
on the group algebra $\IC[\root]$ and more generally on $\IC[\nu+\root]$, $\nu\in\weight$,
given by
\begin{align}
\pi(\tau)f\medspace(\mu)&=\mu(\tau)f(\mu)
\label{eq:model action 1}\\
\pi(\zeta_{\alpha})f&=L_{\alpha}\epsilon_{\alpha}f
\label{eq:model action 2}
\end{align}
where $L_{\alpha}$ is the operator of left translation by $\alpha$ and $\epsilon_{\alpha}$
acts as multiplication by the function $\mu\rightarrow\epsilon(\alpha,\mu-\nu)$,
$\mu\in\nu+\root$. Using the natural basis $\delta_{\mu}$ of $\IC[\nu+\root]$, one readily
verifies that $\pi(\tau)\pi(\zeta_{\alpha})\pi(\tau)^{*}=\alpha(\tau)\pi(\zeta_{\alpha})$.
Moreover, \eqref{eq:cocycle id} implies that
\begin{equation}
\pi(\zeta_{\alpha})\pi(\zeta_{\beta})=\epsilon(\alpha,\beta)\pi(\zeta_{\alpha+\beta})
\end{equation}
so that $\pi$ is a projective unitary representation of $T\times\check T$ satisfying
\eqref{eq:desiderata 2} and \eqref{eq:desiderata 3}. Notice that the choice of the cocycle
is immaterial since (i) of proposition \ref{pr:commutator map} implies that any two choices
$\epsilon,\eta$ satisfy $\epsilon=\eta\cdot b$ for some coboundary $b$ of the form
\eqref{eq:coboundary id}.
Thus, if $M_{a}$ acts on $\IC[\nu+\root]$ as multiplication by the function
$\mu\rightarrow a(\mu-\nu)$, $M_{a}$ commutes with $T$ and satisfies
\begin{equation}
M_{a}\pi_{\epsilon}(\zeta_{\alpha})M_{a}^{*}=a(\alpha)\pi_{\eta}(\zeta_{\alpha})
\end{equation}
and therefore gives the unitary equivalence of $\pi_{\epsilon}$ and $\pi_{\eta}$.
We now have

\begin{proposition}\label{pr:SvN for zeromodes}\hfill
\begin{enumerate}
\item For any $\nu\in\weight$, the representation of $T\times\check T$ on $\IC[\nu+\root]$
given by \eqref{eq:model action 1}--\eqref{eq:model action 2} is irreducible. Any operator
on $\IC[\nu+\root]$ commuting with $\pi$ is a scalar.
\item $\IC[\nu+\root]$ and $\IC[\nu'+\root]$ are unitarily equivalent iff $\nu-\nu'\in\root$.
\item Let $\H$ be a pre--Hilbert space supporting a projective unitary representation $\pi$
of $T\times\check T$ with a preferred unitary lift over $T$ satisfying
\eqref{eq:desiderata 2} and \eqref{eq:desiderata 3}. If $\H$ is the algebraic direct sum of its
$T$--weight spaces and is irreducible, then for some $\nu\in\weight$, $\H\cong\IC[\nu+\root]$.
\end{enumerate}
\end{proposition}
\proof
(i)
The irreducibility of $\IC[\nu+\root]$ is clear since any non--zero submodule $\K$ is the direct
sum of its weight spaces and therefore contains a basis vector $\delta_{\mu}$, $\mu\in\nu+\root$.
Thus, $\H=\check T\IC\delta_{\mu}\subseteq\K$. Any operator $S$ commuting with $T$ necessarily
acts as multiplication by a function $f_{S}$ since $\pi(\tau)S\delta_{\mu}=\mu(\tau)S\delta_{\mu}$
and therefore $S\delta_{\mu}=f_{S}(\mu)\delta_{\mu}$ for some $f_{S}(\mu)\in\T$. If $S$ commutes
with $\check T$, $f_{S}$ is invariant under translations and therefore $S\equiv f_{S}(\nu)\cdot 1$.

(ii)
If $\nu-\nu'\in\root$, then by \eqref{eq:cocycle id}
\begin{equation}
\epsilon(\alpha,\mu-\nu')=
\epsilon(\alpha,\mu-\nu)
\frac{\epsilon(\alpha+\mu-\nu,\nu-\nu')}{\epsilon(\mu-\nu,\nu-\nu')}
\end{equation}
for any $\alpha\in\root$ and $\mu\in\nu+\root$. It follows that the operator $M$ acting on
$\IC[\nu+\root]$ as multiplication by $\mu\rightarrow\epsilon(\mu-\nu,\nu-\nu')$ satisfies
$M\pi_{\nu}(\zeta_{\alpha})M^{*}=\pi_{\nu'}(\zeta_{\alpha})$ and clearly commutes with
$T$ thus giving the unitary equivalence of $\pi_{\nu}$ and $\pi_{\nu'}$. The converse is clear.

(iii)
Let $\H=\bigoplus_{\mu\in\weight}\H_{\mu}$ be the weight space decomposition of $\H$ for
the action of $T$ and $0\neq v_{\nu}\in\H_{\nu}$ a vector of unit length.
$\K_{\nu}=\bigoplus_{\alpha\in\root}\IC\cdot\pi(\zeta_{\alpha})v_{\nu}$ is invariant
under $T$ by \eqref{eq:desiderata 2} and \eqref{eq:desiderata 3} and under $\check T$. By
irreducibility,
$\H=\K_{\nu}\cong\IC[\nu+\root]$ where the $T$--equivariant equivalence is given by
$\pi(\zeta_{\alpha})v_{\nu}\rightarrow\delta_{\alpha+\nu}$ and we may therefore assume
that $T\times\check T$ is acting on $\IC[\nu+\root]$.
Denoting the left--translation action of $\root$ on $\H$ by $\alpha\rightarrow L_{\alpha}$,
we notice that, by \eqref{eq:desiderata 2} and \eqref{eq:desiderata 3} the operator
$L_{\alpha}^{*}\pi(\zeta_{\alpha})$ commutes with $T$ and therefore acts as multiplication
by a $\T$--valued function $\mu\rightarrow\eta(\alpha,\mu)$. Normalising the lifts
$\pi(\zeta_{\alpha})\in U(\H)$, we may assume that $\eta(\alpha,\nu)=1$ for any $\alpha$.
The $\pi(\zeta_{\alpha})$ give a projective representation of $\check T$ and therefore
satisfy
\begin{equation}\label{eq:epsilon}
\pi(\zeta_{\alpha})\pi(\zeta_{\beta})=\epsilon(\alpha,\beta)\pi(\zeta_{\alpha+\beta})
\end{equation}
for some function $\epsilon:\root\times\root\rightarrow\T$ which, by the associativity
of the multiplication in $U(\H)$ and \eqref{eq:desiderata 2} and \eqref{eq:desiderata 3}
is a cocycle with commutator map $\omega$. Since
$\pi(\zeta_{\alpha})=L_{\alpha}\eta(\alpha,\cdot)$,
\eqref{eq:epsilon} yields, for any $\mu\in\nu+\root$
\begin{equation}
\eta(\alpha,\beta+\mu)\eta(\beta,\mu)=\epsilon(\alpha,\beta)\eta(\alpha+\beta,\mu)
\end{equation}
Evaluating at $\mu=\nu$ yields $\eta(\alpha,\mu)=\epsilon(\alpha,\mu-\nu)$ \halmos\\

Recall now from section \ref{se:restriction} that, on the finite energy subspace of a
positive energy representation, the infinitesimal generator of rotations $d$ and the
generating function $X_{\alpha}(z)$ satisfy
\begin{xalignat}{2}
[d,\pi(\tau)]&=0&
[d,\pi(\zeta_{\alpha})]&=\pi(\zeta_{\alpha})(\pi(\alpha)+\half{\<\alpha,\alpha\>})
\label{eq:zeromodes requirements 2}\\
\pi(\tau)X_{\alpha}(z)\pi(\tau)^{*}&=\alpha(\tau)X_{\alpha}&
\pi(\zeta_{\beta})X_{\alpha}(z)\pi(\zeta_{\beta})^{*}&=z^{-\<\beta,\alpha\>}X_{\alpha}(z)
\label{eq:zeromodes requirements 3}
\end{xalignat}
We now give an explicit construction of operators possessing the same commutation
relations in any of the irreducible $T\times\check T$--modules constructed above.

\begin{proposition}
Let $\H=\IC[\nu+\root]$ be the irreducible representation of $T\times\check T$ given by
\eqref{eq:model action 1}--\eqref{eq:model action 2}. Then
\begin{enumerate}
\item Any operator $d$ on $\H$ satisfying \eqref{eq:zeromodes requirements 2} is given,
up to addition by a constant, by
\begin{equation}\label{eq:mini sugawara}
d=\half{1}\sum_{k}\pi(h_{k})\pi(h^{k})
\end{equation}
where $h_{k},h^{k}$ are dual basis of $\tc$ with respect to the basic inner product.
\item For $\alpha\in\root$, the unique solution in $X_{\alpha}(z)\in\End(\H)[[z^{-1},z]]$
to the commutation relations \eqref{eq:zeromodes requirements 3} and
\begin{equation}
[d,X_{\alpha}(z)]=z\frac{d}{dz}X_{\alpha}(z)\label{with energy}
\end{equation}
is given, up to a multiplicative constant by
\begin{equation}\label{eq:vertex for zeromodes}
X_{\alpha}(z)=V_{\alpha}z^{\alpha+\half{\<\alpha,\alpha\>}}
\end{equation}
where $V_{\alpha}=L_{\alpha}\epsilon^{\dagger}_{\alpha}$ and $\epsilon^{\dagger}_{\alpha}$
acts as multiplication by the function $\mu\rightarrow\epsilon(\mu-\nu,\alpha)$ so that
\begin{align}
V_{\alpha}V_{\beta}&=\epsilon(\beta,\alpha)V_{\alpha+\beta}
\label{eq:transposed cocycle}\\
V_{\alpha}z^{\beta}&=z^{\beta-\<\alpha,\beta\>}V_{\alpha}
\label{eq:V and z}
\end{align}
\item
If the cocycle $\epsilon$ giving the action \eqref{eq:model action 2} is chosen with the
normalisations \eqref{eq:normalised cocycle}, then $X_{0}(z)=1$ and $X_{\alpha}(z)$ has
the formal adjunction property $X_{\alpha}(z)^{*}=X_{-\alpha}(z)$.
\end{enumerate}
\end{proposition}
\proof
(i)
The infinitesimal action of $\tc$ on $\H$ corresponding to $\pi$ is given
by $\pi(h)\delta_{\mu}=\<h,\mu\>\delta_{\mu}$ and therefore satisfies
\begin{equation}\label{eq:torus roots}
[h,\pi(\zeta_{\alpha})]=\<h,\alpha\>\pi(\zeta_{\alpha})
\end{equation}
It follows that the operator $d$ given by \eqref{eq:mini sugawara} satisfies
\eqref{eq:zeromodes requirements 2}. The uniqueness of $d$ follows from (i)
of proposition \ref{pr:SvN for zeromodes}.

(ii) By \eqref{eq:torus roots}, the operator $z^{\alpha}$ satisfies
\begin{equation}
\pi(\zeta_{\beta})z^{\alpha}\pi(\zeta_{\beta})^{*}=z^{\alpha-\<\alpha,\beta\>}
\end{equation}
As is readily verified,
\begin{xalignat}{3}
\pi(\tau)V_{\alpha}\pi(\tau)^{*}&=\alpha(\tau)V_{\alpha}&
&\text{and}&
\pi(\zeta_{\beta})V_{\alpha}\pi(\zeta_{\beta})^{*}&=V_{\alpha}
\end{xalignat}
and it follows that $X_{\alpha}(z)=V_{\alpha}z^{\alpha+\half{\<\alpha,\alpha\>}}$ satisfies
\eqref{eq:zeromodes requirements 3}. Moreover, using \eqref{eq:zeromodes requirements 2}
\begin{equation}
[d,X_{\alpha}(z)]=z\frac{d}{dz}X_{\alpha}(z)
\end{equation}
so that $X_{\alpha}(z)$ does indeed satisfy the required commutation relations. The relations
\eqref{eq:transposed cocycle} and \eqref{eq:V and z} follow at once from \eqref{eq:cocycle id}
and the fact that $L_{\alpha}\beta L_{\alpha}^{*}=\beta-\<\alpha,\beta\>$ respectively.
The uniqueness of the solution follows from irreducibility. Indeed, if $Y_{\alpha}(z)$
satisfies the required commutation relations, then
$Z_{\alpha}(z)=Y_{\alpha}(z)X_{-\alpha}(z)$ commutes with $T\times\check T$ and must therefore
be a scalar function of $z$. \eqref{with energy} then implies that $Z_{\alpha}(z)$ is a scalar.

(iii)
If \eqref{eq:normalised cocycle} hold, then $V_{0}=1$ and
$V_{\alpha}V_{-\alpha}=\epsilon(-\alpha,\alpha)V_{0}=1$ hence $V_{\alpha}^{*}=V_{-\alpha}$
from which the claimed adjunction property follows.
\halmos\\

\remark Notice that the operator defined by \eqref{eq:vertex for zeromodes} only involves
integral powers of $z$ since $\root$ is even and $\<\alpha,\mu\>\in\IZ$ for any $\alpha\in
\root$ and $\mu\in\weight$.\\

We conclude this subsection by addressing a minor technical point, namely the uniqueness
of the direct sum of projective representations of $T\times\check T$ satisfying
\eqref{eq:desiderata 2}--\eqref{eq:desiderata 3}.
Recall from \S \ref{ss:projective} of chapter \ref{ch:classification}
that a direct sum of projective representations $(\pi_{i},\H_{i})$ of a group $\Gamma$
is a projective representation of $\Gamma$ on $\bigoplus\H_{i}$ restricting on each
$\H_{i}$ to $\pi_{i}$. One such exists iff the pull--backs of the central extensions
\begin{equation}
1\rightarrow\T\rightarrow U(\H_{i})\rightarrow PU(\H_{i})\rightarrow 1
\end{equation}
to $\Gamma$ are isomorphic and depends upon the choice of identifications
$\pi_{i}^{*}U(\H_{i})\cong\pi_{j}^{*}U(\H_{j})$. Since these are unique only up to
multiplication by a character of $\Gamma$, a canonically defined direct sum does not
exist in general. When $\Gamma=T\times\check T$, the following holds however

\begin{proposition}\label{pr:existence of direct sums}
Let $(\pi_{i},\H_{i})$ be projective unitary representations of $T\times\check T$
satisfying \eqref{eq:desiderata 2} and \eqref{eq:desiderata 3}. Then, there exists a
direct sum representation of $T\times\check T$ on $\bigoplus_{i}\H_{i}$, unique up to
unitary equivalence.
\end{proposition}
\proof
Let $\epsilon$ be a cocycle on $\root$ with associated commutator map $\omega$ and
choose lifts $U_{\alpha}^{i}\in U(\H_{i})$ of $\pi_{i}(\zeta_{\alpha})$ satisfying
$U_{\alpha}^{i}U_{\beta}^{i}=\epsilon(\alpha,\beta)U_{\alpha+\beta}^{i}$. The
representation of $T\times\check T$ on $\bigoplus\H_{i}$ obtained by letting
$\tau\in\T$ act as $\bigoplus_{i}\pi_{i}(\tau)$ and $\zeta_{\alpha}$ as
$\bigoplus_{i}U_{\alpha}^{i}$ clearly is a direct sum representation of the $\H_{i}$.
Let now $\pi,\rho$ be two representations on $\H=\bigoplus\H_{i}$ leaving each $\H_{i}$
restricting on each $\H_{i}$ to unitarily equivalent representations.
Let $\epsilon$ be a cocycle with commutator map $\omega$ and choose lifts
$U_{\alpha},V_{\alpha}\in U(\H)$ of $\pi(\zeta_{\alpha}),\rho(\zeta_{\alpha})$
respectively satisfying
\begin{xalignat}{3}\label{eq:choice of lifts}
U_{\alpha}U_{\beta}&=\epsilon(\alpha,\beta)U_{\alpha+\beta}&
&\text{and}&
V_{\alpha}V_{\beta}&=\epsilon(\alpha,\beta)V_{\alpha+\beta}
\end{xalignat}
By assumption, there exist unitaries $I_{i}:\H_{i}\rightarrow\H_{i}$ such that
\begin{xalignat}{3}\label{eq:by assumption}
I_{i}P_{i}\pi(\tau)P_{i}I_{i}^{*}&=P_{i}\rho(\tau)P_{i}&
&\text{and}&
I_{i}P_{i}U_{\alpha}P_{i}I_{i}^{*}&=\cchi_{i}(\alpha)P_{i}V_{\alpha}P_{i}
\end{xalignat}
where $P_{i}$ is the orthogonal projection onto $\H_{i}$ and $\cchi_{i}(\alpha)\in\T$.
Comparing \eqref{eq:choice of lifts} and \eqref{eq:by assumption} shows that each
$\cchi_{i}$ is a character of $\root$ and is therefore given by
$\alpha\rightarrow\alpha(\tau_{i})$ for some $\tau_{i}\in T$. It follows from
\eqref{eq:desiderata 2} and \eqref{eq:desiderata 3} that
\begin{equation}
\pi(\tau_{i})I_{i}P_{i}U_{\alpha}P_{i}I_{i}^{*}\pi(\tau_{i})^{*}=P_{i}V_{\alpha}P_{i}
\end{equation}
Since $\pi(\tau_{i})I_{i}$ intertwines the action of $T$, the unitary
$\bigoplus_{i}\pi(\tau_{i})I_{i}$ gives the equivalence of $\pi$ and $\rho$ \halmos

\ssubsection{Stone--von Neumann theorem and vertex operators for $\lpol\t/\t$}
\label{ss:Stone for nonzero}

This subsection follows \cite[chap. 3]{FLM}. Consider the action of $\lpol\tc/\tc$
on the finite energy subspace of a level 1 positive energy representation of $LG$.
It satisfies
\begin{align}
[h_{1}(n),h_{2}(m)]&=n\delta_{n+m,0}\<h_{1},h_{2}\>
\label{eq:nonzero nonzero}\\
[d,h(n)]&=-nh(n)
\label{eq:nonzero energy}
\end{align}
as well as the formal adjunction property
\begin{equation}\label{eq:formal adjunction}
h(n)^{*}=-\overline{h}(-n)
\end{equation}
where $\overline{h}$ is the canonical conjugation on $\tc$.
An explicit representation of $\lpol\tc/\tc$ satisfying the above may be obtained
as follows. Write $\lpol\tc/\tc=V_{+}\bigoplus V_{-}$ where the $V_{\pm}$ are
spanned by the $h(n)=h\otimes e^{in\theta}$ with $n\lessgtr 0$. Define an action
of $\lpol\tc/\tc\rtimes\IC d$ on the symmetric algebra
\begin{equation}
\SS=SV_{+}=\bigoplus_{k} S^{k}V_{+}
\end{equation}
by
\begin{align}
h(-n)\medspace h_{1}(-n_{1})\otimes\cdots\otimes h_{k}(-n_{k})&=
h(-n)\otimes h_{1}(-n_{1})\otimes\cdots\otimes h_{k}(-n_{k})
\label{eq:nonzero action 1}\\
h(n) \medspace h_{1}(-n_{1})\otimes\cdots\otimes h_{k}(-n_{k})&=
\sum_{j}n\delta_{n,n_{j}}\<h,h_{j}\>
h_{1}(-n_{1})\otimes\cdots\otimes\widehat{h_{j}(-n_{j})}\otimes\cdots\otimes h_{k}(-n_{k})
\label{eq:nonzero action 2}\\
d    \medspace h_{1}(-n_{1})\otimes\cdots\otimes h_{k}(-n_{k})&=
(n_{1}+\cdots+n_{k}) h_{1}(-n_{1})\otimes\cdots\otimes h_{k}(-n_{k})
\end{align}
where $n,n_{1},\ldots n_{k}>0$ throughout. It is readily verified that the relations
\eqref{eq:nonzero nonzero}--\eqref{eq:nonzero energy} hold.
Moreover, the hermitian form $(\cdot,\cdot)$ defined on $V_{+}$ by
$(h_{1}(n),h_{2}(m))=n\delta_{n,m}\<h_{1},\overline{h_{2}}\>$ is positive definite since
$h_{1}\otimes h_{2}\rightarrow -\<h_{1},\overline{h_{2}}\>$ is a positive definite inner
product on $\tc$. $(\cdot,\cdot)$ yields an inner product on $\SS$ for which
\eqref{eq:formal adjunction} holds.\\

Let $\H$ be a pre--Hilbert space supporting a representation of $\lpol\tc/\tc$ satisfying
\eqref{eq:nonzero nonzero} and \eqref{eq:formal adjunction}. By definition, $\H$ is
of positive energy if $\H$ has an $\IN$--grading $\H=\bigoplus_{n\geq 0}\H(n)$, with
$\dim H(n)<\infty$ and the operator $d$ acting as multiplication by $n$ on $\H(n)$
satisfies \eqref{eq:nonzero energy}.

\begin{proposition}\label{pr:nonzero SvN}\hfill
\begin{enumerate}
\item The following holds on $\SS$,
\begin{equation}\label{eq:boson sugawara}
d=\sum_{m>0}h_{k}(-m)h^{k}(m)
\end{equation}
where $h_{k},h^{k}$ are dual basis of $\tc$ for the basic inner product $\<\cdot,\cdot\>$.
\item $\SS=SV_{+}$ is algebraically irreducible under the action of $\lpol\tc/\tc$.
Moreover, $T\in\End(\SS)$ commutes with the $h(n)$ iff $T$ is a scalar.
\item Let $\H$ be a unitary, positive energy representation of $\lpol\tc/\tc$ and define
\begin{equation}\label{eq:vacuum}
\Nu_{\H}=\{\nu\in\H|\thinspace h(n)\nu=0\thickspace\text{for any $h\in\tc$ and $n>0$}\}
\end{equation}
Then $\H\cong\SS\otimes\Nu_{\H}$ where the unitary equivalence is given
by $p\cdot 1\otimes\nu\rightarrow p\nu$ where $p$ is any polynomial in
the $h(n)$ and $\nu\in\Nu_{\H}$.
\end{enumerate}
\end{proposition}
\proof
(i) Let $D$ be the operator given by the right hand--side of \eqref{eq:boson sugawara}
and notice that $D$ is well--defined since the $h(m)$, $m>0$ are locally nilpotent on
$\SS$. It follows from \eqref{eq:nonzero nonzero} that
\begin{equation}
[D,h(n)]=-n h(n)=[d,h(n)]
\end{equation}
Denoting the lowest energy vector in $\SS$ by $1\in S^{0}V_{+}$, we have $D1=0=d1$ and
it follows by cyclicity of 1 under the action of $\lpol\tc/\tc$ that $D=d$.

(ii) If $\Nu\subset\SS$ is invariant under the $h(n)$, it is invariant under
$d$ by \eqref{eq:boson sugawara}. Consequently, it contains
a non--zero lowest energy vector $\nu$ necessarily annihilated by the energy decreasing
$h(n)$, $n>0$. This can only be if $\nu\in S^{0}V_{+}$ and then
$\SS=\U(\lpol\tc/\tc)S^{0}V_{+}\subset\Nu$, where $\U$ denotes the enveloping algebra.
Similarly, any $T$ commuting with the $h(n)$ commutes with $d$ and consequently leaves
$S^{0}V_{+}$ invariant. It therefore acts on it as multiplication by a scalar $\lambda$
and by cyclicity, $T\equiv\lambda$.

(iii) Let $\nu\in\Nu_{\H}$. We claim that the map $\SS\rightarrow\H$ given by $p\cdot 1
\rightarrow p\nu$ is norm--preserving and therefore injective. Let $p,q$ be polynomials
in $\U(\lpol\tc/\tc)$ so that $(p\nu,q\nu)=(q^{*}p\nu,\nu)$. Using
the commutation relations of the $h(n)$'s, the product $q^{*}p$ may be written as the
sum of a constant term $\tau_{0}$ and terms of the form
$h_{1}(-n_{1})\cdot h_{k}(-n_{k})h'_{1}(m_{1})\cdot h'_{l}(m_{l})$ where the $n_{i}$
and $m_{j}$ are positive. However,
\begin{equation}
(h_{1}(-n_{1})\cdots h_{k}(-n_{k})h'_{1}(m_{1})\cdots h'_{l}(m_{l})\nu,\nu)=
(-1)^{k}
(h'_{1}(m_{1})\cdots h'_{l}(m_{l})\nu,
\overline{h}_{k}(n_{k})\cdots\overline{h}_{1}(n_{1})\nu)=0
\end{equation}
and therefore $(p\nu,q\nu)=\tau_{0}\|\nu\|^{2}$. Similarly, $(p\cdot 1,q\cdot 1)=\tau_{0}$
and it follows that  $p\cdot 1\otimes\nu\rightarrow p\nu$ is an $\lpol\tc/\tc$--equivariant
isometry $\SS\otimes\Nu_{\H}\rightarrow\H$. To see that it is surjective, notice
that $\Nu_{\H}$ is invariant under $d$ so that $\SS\otimes\Nu_{\H}$ is a
graded submodule of $\H$. Since the eigenspaces of $d$ are finite--dimensional, it possesses
a graded orthogonal complement $\K$ invariant under $\lpol\tc/\tc$. If $\K\neq\{0\}$, it
possesses a lowest energy
vector $\xi$ since $d$ is bounded below. However, $\xi$ is necessarily annihilated by the
$h(n)$, $n>0$ and therefore lies in $\Nu_{\H}$, a contradiction \halmos\\

Using $\<\cdot,\cdot\>$ to identify $\tc^{*}$ with $\tc$, we consider, for any
$\alpha\in\tc^{*}$ the elements $\alpha(m)\in\lpol\tc/\tc$ and introduce the
exponential operators 
\begin{equation}
E^{\pm}(\alpha,z)=
\exp\Bigl(-\sum_{m\gtrless 0} \frac{\alpha(m)}{m}z^{-m}\Bigr)
\end{equation}
These are to be viewed as formal Laurent series with coefficents in the
endomorphisms of $\SS$. We will need the following elementary

\begin{lemma}\label{expformula}
If $[a,b]$ commutes with both $a$ and $b$, then
$[a,e^{b}]=[a,b]e^{b}$ and $e^{a}e^{b}=e^{[a,b]}e^{b}e^{a}$.
\end{lemma}
\proof
Induction shows that $[a,b^{n}]=n[a,b]b^{n-1}$ and  therefore $[a,e^{b}]=[a,b]e^{b}$.
It follows that $ae^{b}=e^{b}(a+[a,b])$ thus $a^{n}e^{b}=e^{b}(a+[a,b])^{n}$
whence $e^{a}e^{b}=e^{[a,b]}e^{b}e^{a}$ \halmos

\begin{proposition}\label{pr:exponential commutators}
 The following commutation relations hold
\begin{align}
E^{\pm}(\alpha,z)E^{\pm}(\beta,z)&=E^{\pm}(\alpha+\beta,z)\\
E^{\pm}(\alpha,z)E^{\pm}(-\alpha,z)&=1\\
E^{\pm}(\alpha,z)^{*}&=E^{\mp}(\overline{\alpha},z)\\
E^{\pm}(\alpha,z)E^{\pm}(\beta,\zeta)&=E^{\pm}(\beta,\zeta)E^{\pm}(\alpha,z)\\
E^{+}(\alpha,z)E^{-}(\beta,\zeta)&=
\Bigl(1-\frac{\zeta}{z}\Bigr)^{\<\alpha,\beta\>}E^{-}(\beta,\zeta)E^{+}(\alpha,z)
\label{eq:ope}
\end{align}
where $\Bigl(1-\frac{\zeta}{z}\Bigr)^{\<\alpha,\beta\>}$ is defined by \eqref{eq:binomial}.
\end{proposition}
\proof All but the last relation are trivial. The last follows from lemma \ref{expformula}
and the fact that
\begin{equation}
[\sum_{m>0}\frac{\alpha(m)}{m}z^{-m},\sum_{n<0}\frac{\beta(n)}{n}\zeta^{-n}]
=-\<\alpha,\beta\>\sum_{m>0}\frac{1}{m}\Bigl(\frac{\zeta}{z}\Bigr)^{m}
=\<\alpha,\beta\>\log\Bigl(1-\frac{\zeta}{z}\Bigr)
\end{equation}
\halmos
\newpage

\begin{proposition}\label{pr:nonzero vertex}\hfill
\begin{enumerate}
\item For any $\alpha\in\tc$, the formal Laurent series
$X_{\alpha}(z)=E^{-}(\alpha,z)E^{+}(\alpha,z)\in\End(\SS)[[z,z^{-1}]]$ is well--defined,
satisfies
\begin{align}
[h(n),X_{\alpha}(z)]&=\<h,\alpha\>z^{n}X_{\alpha}(z)
\label{eq:covariance}\\
[d,X_{\alpha}(z)]&=z\frac{d}{dz}X_{\alpha}(z)
\label{eq:homogeneity}
\end{align}
and the formal adjunction property $X_{\alpha}(z)^{*}=X_{\overline{\alpha}}(z)$.
\item The above commutation relations characterise $X_{\alpha}(z)$ uniquely, up to
a multiplicative constant.
\end{enumerate}
\end{proposition}
\proof 
(i)
$X_{\alpha}(z)$ is well--defined on $\SS$ since, for any $\psi\in\SS$, only finitely many
terms in $E^{+}(\alpha,z)\psi$ do not vanish. Moreover,
$[h(n),-\sum_{m\gtrless 0}\frac{\alpha(m)}{m}z^{-m}]=z^{n}\<h,\alpha\>\delta_{n\lessgtr 0}$
and therefore, by lemma \ref{expformula},
$[h(n),X_{\alpha}(z)]=z^{n}\<h,\alpha\>X_{\alpha}(z)$. The second relation follows from
the general fact that if $\phi(z)=\sum_{n}\phi_{n}z^{-n}$ satisfies
$[d,\phi(z)]=z\frac{d}{dz}\phi(z)$ or equivalently
$[d,\phi_{n}]=-n\phi_{n}$, then so do $p(\phi(z))$ and $\exp(\phi(z))$ for any polynomial
$p$, whenever these are well--defined. The adjunction property of $X_{\alpha}(z)$ follows
from those of the $E^{\pm}(\alpha,z)$.

(ii)
Let $Y_{\alpha}\in\End(\SS)[[z,z^{-1}]]$ satisfy \eqref{eq:covariance}--\eqref{eq:homogeneity}.
The operator $\phi(z)=E^{-}(-\alpha,z)Y_{\alpha}(z)E^{+}(-\alpha,z)$ is easily seen to be a
well--defined element of $\End(\SS)[[z,z^{-1}]]$ commuting with the $h(n)$. By (i), it follows
that $\phi(z)=\sum_{n}a_{n}z^{-n}$ for some $a_{n}\in\IC$. The relation \eqref{eq:homogeneity}
however implies that $a_{n}=0$ if $n\neq 0$ \halmos

\ssubsection{Irreducible representations of $\lpol\t/\t\times T\times\check T$}
\label{ss:LT irreducibles}

\begin{proposition}\label{pr:LT irreducibles}\hfill
\begin{enumerate}
\item The tensor product representations of $\lpol\t/\t\times T\times\check T$ on
the modules $\SS\otimes\IC[\nu+\root]$ given by
\eqref{eq:model action 1}--\eqref{eq:model action 2} and
\eqref{eq:nonzero action 1}--\eqref{eq:nonzero action 2} are irreducible. Moreover,
$\SS\otimes\IC[\nu+\root]$ and $\SS\otimes\IC[\nu'+\root]$ are equivalent if, and
only if $\nu-\nu'\in\root$.
\item Let $\H$ be a unitary positive energy representation of $\lpol\t/\t\times T\times
\check T$ satisfying \eqref{eq:desiderata 2}--\eqref{eq:desiderata 3}, \eqref{eq:nonzero
nonzero} and \eqref{eq:desiderata 4}--\eqref{eq:desiderata 5}. If $\H$ is irreducible,
it is unitarily equivalent to $\SS\otimes\IC[\nu+\root]$ for some $\nu\in\weight$.
\end{enumerate}
\end{proposition}
\proof
(i)
Let $0\neq V\subset\SS\otimes\IC[\nu+\root]$ be a submodule. $V$ is the direct sum
of its $T$--weight spaces which are necessarily of the form $V_{\mu}\otimes\delta_{\mu}$
where $V_{\mu}\subset\SS$ is invariant under $L\t/\t$ and therefore equal to $0$ or
$\SS$ by proposition \ref{pr:nonzero SvN}. Thus, $V_{\mu}=\SS$ for some $\mu$ and
therefore, using the action of $\check T$, $V_{\mu}=\SS$ for all $\mu\in\nu+\root$
whence $V=\SS$. The last statement follows from proposition \ref{pr:SvN for zeromodes}.

(ii)
By proposition \ref{pr:nonzero SvN}, $\H\cong\SS\otimes\Nu_{\H}$.
Clearly, $T\times\check T$ leaves $\Nu_{\H}$ invariant and the isomorphism is easily
seen to be equivariant for this action. By irreducibility of $\H$, $\Nu_{\H}$ is
irreducible under $T\times\check T$ and therefore, by proposition
\ref{pr:SvN for zeromodes} of the form $\IC[\nu+\root]$ for some $\nu\in\weight$
\halmos

\ssection{The vertex operator construction}
\label{se:vertex operators}

This section follows \cite[chap. 7]{FLM} and \cite[\S 6]{GO2}.
We show below that the vertex operators $X_{\alpha}(z)$ give rise, in any
of the irreducible representations of $\lpol\t/\t\times T\times\check T$ classified
in section \ref{se:level 1 LT}, to level 1 irreducible positive energy
representations of $\lpol\g$ and that the latter may all be obtained in this way.\\

It will in fact be more convenient to work with the direct sum
\begin{equation}
\SS\otimes\IC[\weight]=\bigoplus_{\nu\in\weight/\root}\SS\otimes\IC[\nu+\root]
\end{equation}
of all irreducible $\lpol\t/\t\times T\times\check T$--modules. By proposition
\ref{pr:existence of direct sums}, the action of $T\times\check T$ on the factor
$\IC[\weight]$ is unambiguously defined and we may, by uniqueness assume that it
extends to one of $T\times\weight$ given by
\begin{align}
\pi(\tau)f(\mu)&=\mu(\tau)f(\mu)\\
\pi(\lambda)f&=L_{\lambda}\epsilon_{\lambda}f
\end{align}
As customary, $L_{\lambda}$ is the left--translation operator and $\epsilon_{\lambda}$
acts as multiplication by the function $\mu\rightarrow\epsilon(\lambda,\mu)$ where
$\epsilon$ is a $\T$--valued cocycle on $\weight$ whose commutator map $\wt\omega$
satisfies $\wt\omega(\alpha,\beta)=(-1)^{\<\alpha,\beta\>}$ whenever $\alpha,\beta\in\root$.
The choice of a particular $\wt\omega$ satisfying this requirement is clearly immaterial.\\

Choose $\epsilon$ to be normalised in the sense of \eqref{eq:normalised cocycle} and
let $\epsilon^{\dagger}(\lambda,\mu)=\epsilon(\mu,\lambda)$ be the transposed cocycle.
Let
\begin{equation}\label{eq:def U}
U_{\lambda}:=\pi(\lambda)=L_{\lambda}\epsilon_{\lambda}
\end{equation}
and $V_{\lambda}=L_{\lambda}\epsilon_{\lambda}^{\dagger}$ so that
\begin{align}
U_{\lambda}U_{\mu}&=\epsilon(\lambda,\mu)U_{\lambda+\mu}
\label{eq:U}\\
V_{\lambda}V_{\mu}&=\epsilon(\mu,\lambda)V_{\lambda+\mu}
\label{eq:V}\\
U_{\lambda}V_{\mu}U_{\lambda}^{*}V_{\mu}^{*}&=1
\label{eq:UV}
\end{align}
and define for any $\alpha\in\root$, the vertex operators
\begin{equation}\label{eq:definition of vo}
X_{\alpha}(z):=
E^{-}(\alpha,z)V_{\alpha}z^{\alpha+\half{\<\alpha,\alpha\>}}E^{+}(\alpha,z)=
\sum_{n\in\IZ}X_{\alpha}(n)z^{-n}
\end{equation}

We shall need the following

\begin{lemma}\label{le:norm 2}
Let $\alpha\in\root$. Then $\<\alpha,\alpha\>=2$ iff $\alpha$ is a root.
\end{lemma}
\proof
Let $\alpha\in\root$ of length $\sqrt{2}$. It is clearly sufficient to show that
$w\alpha$ is a root for some $w\in W$. Write $\alpha=\beta_{1}+\cdots+\beta_{r}$
where the $\beta_{i}$ are roots. Clearly, one cannot have
$\<\alpha,\beta_{i}\>\leq 0$ for all $i$ else
\begin{equation}
\|\alpha-\sum\beta_{i}\|^{2}=
\|\alpha\|^{2}-2\sum\<\alpha,\beta_{i}\>+\|\sum\beta_{i}\|^{2}>0
\end{equation}
Thus there exists a $\beta_{i_{1}}$ with $\<\alpha,\beta_{i_{1}}\>\in\{1,2\}$. If
$\<\alpha,\beta_{i_{1}}\>=2$ then $\alpha=\beta_{i_{1}}$ and $\alpha$ is a root.
Otherwise, let $\sigma_{\beta_{i_{1}}}\in W$ be the simple reflection corresponding
to $\beta_{i_{1}}$,
then
\begin{equation}
\sigma_{\beta_{i_{1}}}\alpha=\alpha-\beta_{i_{1}}=\sum_{i\neq i_{1}}\beta_{i}
\end{equation}
Iterating this argument shows that
$\sigma_{\beta_{i_{k}}}\cdots\sigma_{\beta_{i_{1}}}\alpha$
is a root for some sequence of distinct $1\leq i_{j}\leq r$ \halmos

\begin{theorem}\label{th:existence of representation}
There exists a basis of $\gc/\tc$ given by root vectors $x_{\alpha}$ such that
\begin{align}
\pi(x_{\alpha}(n))&=X_{\alpha}(n)\\
\pi(h(n))&=h(n)\thickspace n\neq 0\\
\pi(h)&=h\\
\pi(d)&=\half{1}h_{k}h^{k}+\sum_{m>0}h_{k}(-m)h^{k}(m)
\label{eq:pi(d)}
\end{align}
defines a level 1 positive energy representation of $\lpol\gc\rtimes\IC d$ on
$\SS\otimes\IC[\weight]$.
\end{theorem}
\proof
We begin by computing the commutation relations of the vertex operators
$X_{\alpha}(z)$. Using \eqref{eq:ope}, \eqref{eq:V} and
$V_{\beta}^{*}z^{\alpha}V_{\beta}=z^{\alpha+\<\alpha,\beta\>}$, we get
\begin{equation}
X_{\alpha}(z)X_{\beta}(\zeta)=
\epsilon(\beta,\alpha)
\Bigl(1-\frac{\zeta}{z}\Bigr)^{\<\alpha,\beta\>}
\Bigl(\frac{z}{\zeta}\Bigr)^{\half{\<\alpha,\beta\>}}
V_{\alpha+\beta}
z^{\alpha+\half{\<\alpha,\alpha\>}+\half{\<\alpha,\beta\>}}
\zeta^{\beta+\half{\<\beta,\beta\>}+\half{\<\alpha,\beta\>}}
\reg(\alpha,\beta,z,\zeta)
\end{equation}
where
\begin{equation}
\reg(\alpha,\beta,z,\zeta)=
E^{-}(\alpha,z)E^{-}(\beta,\zeta)E^{+}(\alpha,z)E^{+}(\beta,\zeta)=
\reg(\beta,\alpha,\zeta,z)
\end{equation}
Similarly,
\begin{equation}
X_{\beta}(\zeta)X_{\alpha}(z)=
\epsilon(\alpha,\beta)
\Bigl(1-\frac{z}{\zeta}\Bigr)^{\<\alpha,\beta\>}
\Bigl(\frac{\zeta}{z}\Bigr)^{\half{\<\alpha,\beta\>}}
V_{\alpha+\beta}
z^{\alpha+\half{\<\alpha,\alpha\>}+\half{\<\alpha,\beta\>}}
\zeta^{\beta+\half{\<\beta,\beta\>}+\half{\<\alpha,\beta\>}}
\reg(\alpha,\beta,z,\zeta)
\end{equation}
Since $\epsilon(\alpha,\beta)=(-1)^{\<\alpha,\beta\>}\epsilon(\beta,\alpha)$,
we get
\begin{equation}
[X_{\alpha}(z),X_{\beta}(\zeta)]=
\epsilon(\beta,\alpha)
\D(z,\zeta)
V_{\alpha+\beta}
    z^{\alpha+\half{\<\alpha,\alpha\>}+\half{\<\alpha,\beta\>}}
\zeta^{\beta +\half{\<\beta,\beta\>}  +\half{\<\alpha,\beta\>}}
\reg(\alpha,\beta,z,\zeta)
\end{equation}
where, using \eqref{eq:delta}--\eqref{eq:delta'}
\begin{equation}
\begin{split}
\D(z,\zeta)
&=
\Bigl(\frac{\zeta}{z}\Bigr)^{-\half{\<\alpha,\beta\>}}
\biggl[
\Bigl(1-\frac{\zeta}{z}\Bigr)^{\<\alpha,\beta\>}
-(-1)^{\<\alpha,\beta\>}
\Bigl(1-\frac{z}{\zeta}\Bigr)^{\<\alpha,\beta\>}
\Bigl(\frac{z}{\zeta}\Bigr)^{-\<\alpha,\beta\>}
\biggr]\\
&=
\begin{cases}
\delta'\Bigl(\frac{\zeta}{z}\Bigr)&\text{if $\<\alpha,\beta\>=-2$}\\[1.5ex]
\Bigl(\frac{\zeta}{z}\Bigr)^{\half{1}}
\delta\Bigl(\frac{\zeta}{z}\Bigr)&\text{if $\<\alpha,\beta\>=-1$}\\[1.5ex]
0&\text{if $\<\alpha,\beta\>\geq 0$}
\end{cases}
\end{split}
\end{equation}

Let $\Delta\subset\root$ be the set of elements of squared length 2. Then, for
$\alpha,\beta\in\Delta$, we get by \eqref{eq:delta identity}--\eqref{eq:delta'
identity} and \eqref{eq:normalised cocycle},
\begin{equation}\label{KM commutation}
[X_{\alpha}(z),X_{\beta}(\zeta)]=
\left\{\begin{array}{lll}
\alpha(\zeta)\delta\Bigl(\frac{\zeta}{z}\Bigr)+
\delta'\Bigl(\frac{\zeta}{z}\Bigr)
&\text{if $\alpha+\beta=0$}
&\text{{\it i.e.}~ iff $\<\alpha,\beta\>=-2$}\\[1.5ex]
\epsilon(\beta,\alpha)X_{\alpha+\beta}(\zeta)\delta\Bigl(\frac{\zeta}{z}\Bigr)
&\text{if $\alpha+\beta\in\Delta$}
&\text{{\it i.e.}~ iff $\<\alpha,\beta\>=-1$}\\[1.5ex]
0
&\text{if $\alpha+\beta\notin\Delta\cup\{0\}$}
&\text{{\it i.e.}~ iff $\<\alpha,\beta\>\geq 0$}
\end{array}\right.
\end{equation}

The commutation relations \eqref{KM commutation} show that the complex vector space
$\l\subset\End(\SS\otimes\IC[\weight])$ spanned by the operators $X_{\alpha}(0)$
and $h$, with $\alpha\in\Delta$ and $h\in\tc$ is closed under the commutator bracket
and is a simple Lie algebra.
By lemma \ref{le:norm 2}, $\l$ and $\gc$ have the same root system and it follows
that there exists a Lie algebra isomorphism $\pi:\gc\rightarrow\l$ acting as the
identity on the Cartan subalgebra $\tc$. Set $x_{\alpha}=\pi^{-1}(X_{\alpha}(0))$.
The Lie algebra spanned by the modes of the $X_{\alpha}(z)$ and the $h(n)$ is clearly
a central extension of $\l[z,z^{-1}]\cong\gc[z,z^{-1}]$.
By \eqref{KM commutation}, the corresponding cocycle is given by $a(m)\otimes b(n)
\rightarrow m\delta_{m+n,0}(a,b)$ where $(\cdot,\cdot)$ is an ad--invariant bilinear
form on $\l$ and is therefore a multiple of the basic inner product $\<\cdot,\cdot\>$
on $\gc\cong\l$.
By construction however, the two coincide on $\tc$ and therefore on $\gc$ and it
follows that $\SS\otimes\IC[\weight]$ is a level 1 representation of $\lpol\gc$.
To conclude, notice that the operator $d$ given by \eqref{eq:pi(d)} satisfies by
construction
$[d,X_{\alpha}(n)]=-nX_{\alpha}(n)$ and $[d,h(m)]=-mh(m)$, is bounded below and
has finite--dimensional eigenspaces. Finally, the unitarity of the representation
follows from the formal adjunction property $h(z)^{*}=-\overline{h}(z)$ and
$X_{\alpha}(z)^{*}=X_{-\alpha}(z)$ \halmos

\begin{theorem}\label{th:vo reducibility}
Each $\SS\otimes\IC[\nu+\root]\subset\SS\otimes\IC[\weight]$ is invariant and
irreducible under $\lpol\g$ and is the level 1 positive energy representation
whose lowest energy subspace is the minimal $G$--module with weights lying in
$\nu+\root$. The corresponding highest weight vector is $1\otimes\delta_{\mu}$
where $\mu\in\nu+\root$ is the unique minimal dominant weight. It follows that
$\SS\otimes\IC[\weight]$ is the direct sum of all level 1 representations of
$\lpol\g$, each with multiplicity one.
\end{theorem}
\proof
We follow \cite[\S 6.4]{GO2}.
By construction, $\SS\otimes\IC[\nu+\root]$ is invariant under
$\lpol\t/\t\times T\times\check T$ and therefore under the $X_{\alpha}(z)$.
Since $\SS\otimes\IC[\nu+\root]$ is a unitary representation, it is the sum
of its irreducible summands. These correspond exactly to and are generated
by highest weight vectors, {\it i.e.}~ elements $\Omega_{\mu}$ satisfying, for some
$\mu\in\nu+\root$
\begin{alignat}{2}
h(0)\Omega_{\mu}&=\<h,\mu\>\Omega_{\mu}&&\quad\text{$h\in\tc$}\\
h(n)\Omega_{\mu}&=0&&\quad\text{$h\in\tc$, $n>0$}\\
X_{\alpha}(n)\Omega_{\mu}&=0&&\quad\text{$n>0$ or $n=0$, $\alpha>0$}
\end{alignat}

The second condition implies that $\Omega_{\mu}\in 1\otimes\IC[\weight]$ since
this is precisely the vacuum subspace for the $h(n)$, $n>0$ and together with the
first is equivalent to $\Omega_{\mu}=1\otimes\delta_{\mu}$. Next,
\begin{equation}
X_{\alpha}(n)=\frac{1}{2\pi i}\oint\frac{dz}{z}z^{n}X_{\alpha}(z)
\end{equation}
where the formal contour integration is shorthand for taking the coefficent of
$z^{-1}$ in the power series expansion of the integrand. Since
$h(n)1\otimes\delta_{\mu}=0$ for $n>0$ so that $\Omega_{\mu}$ is fixed by
$E^{+}(\alpha,z)$, we have
\begin{equation}\label{highest weight}
X_{\alpha}(n)1\otimes\delta_{\mu}=
\frac{\epsilon(\alpha,\mu)}{2\pi i}\oint\frac{dz}{z}
z^{n+\<\alpha,\mu\>+\half{\<\alpha,\alpha\>}}
E^{-}(\alpha,z)1\otimes\delta_{\alpha+\mu}
\end{equation}

The coefficient of $z^{k}$ in the power series expansion of
$E^{-}(\alpha,z)=\exp\Bigl(\sum_{m>0}\frac{\alpha(-m)}{m}z^{m}\Bigr)$ is 0 if $k<0$
and otherwise of the form $\frac{\alpha(-1)^{k}}{k!}$ plus terms involving powers
of $\alpha(-1)$ strictly smaller than $k$ so that it does not vanish if $k\geq 0$.
Since the map $p\rightarrow p\cdot 1$, where $p$ is any polynomial in the $h(n)$
with $n>0$, is an isomorphism, we conclude that \eqref{highest weight} vanishes iff
\begin{equation}
n+\<\alpha,\mu\>+\half{\<\alpha,\alpha\>}\geq 1
\end{equation}
If $n>0$, this is equivalent to $\<\alpha,\mu\>\geq -1$ for any $\alpha$ {\it i.e.}~
(replacing $\alpha$ by $-\alpha$) to $|\<\alpha,\mu\>|\leq 1$ so that $\mu$ is
a minimal weight. Choosing now $n=0$ and $\alpha>0$ yields $\<\alpha,\mu\>\geq 0$
so that $\mu$ is dominant and is therefore the unique minimal dominant weight in
$\nu+\root$. It follows that $\SS\otimes\IC[\weight]$ is the direct sum of all
level 1 representations of $\lpol\g$ since these are in bijective correspondence
with $\root$ cosets in $\weight$ \halmos\\

\remark As remarked in the introduction, the vertex operator construction does
not provide one with a natural action of $LG$ on the Hilbert space completion
of $\SS\otimes\IC[\weight]$. The infinitesimal action of $\lpol\g$ can however
be exponentiated to $LG$ by using analytic methods \cite{GoWa,TL1}. When
$G=\Spin_{2n}$ (or $\SU_{n}$), one may alternatively resort to the Fermionic
realisation of the level 1 representations described in chapter \ref{ch:fermionic}.
The latter exponentiates to $LG$ by construction and the isomorphism of
$\lpol\g$--modules may be used to transport the action of $LG$
on the fermionic Fock spaces to one on the completion of $\SS\otimes\IC[\weight]$.

\ssection{The level 1 primary fields}
\label{se:level 1 primary}

Having realised the finite energy subspaces of the level 1 representations of
$LG$ as summands of $\SS\otimes\IC[\weight]$, we construct in this section the
level 1 primary fields by following a scheme similar to the one adopted for the
vertex operators $X_{\alpha}(z)$. We obtain in \S \ref{ss:LT equivariance}
their equivariance properties with respect to $LT$. The calculation relies upon
the assumption that the primary fields extend to continuous operator--valued
distributions and therefore has heuristic value only. We use it nonetheless to
deduce their explicit form and then show in \S \ref{ss:construction} that the
guesses have the correct commutation relations with $\lpol\g$.\\

The form of the primary fields resembles that of the $X_{\alpha}(z)$ with an
additional term accounting for the non--integrality of the weight lattice.
In \S \ref{ss:extension}, we give the complete list of simply--laced $G$ for
which this correction factor cannot altogether be dispensed with. It comprises
$\SU_{2}$ and $\Spin_{2n}$ with $n$ not divisible by $4$.\\

Finally, in \S \ref{ss:fermi field}, we study the level 1 primary fields for
$L\Spin_{2n}$. We consider the vector primary field and recover the result of
chapter \ref{ch:fermionic} that it is a Fermionic field. We also prove that
all level 1 primary fields for $L\Spin_{8}$ are Fermi fields.

\ssubsection{Action of $\check T$ and $L_{0}$ on $\SS\otimes\IC[\weight]$}
\label{ss:identify rot and root}

We begin by identifying explicitly the operators giving the action of $\check T\subset LT$
and $L_{0}$ on $\SS\otimes\IC[\weight]$. Notice that the latter supports two projective
representations of $\check T$. The first is given by \eqref{eq:def U} and extends to one of
$\weight$ while the second is obtained by restricting the representation $\pi$ of $LG$ to
$\check T\subset LT$. Surprisingly perhaps, these do not coincide. In fact the following
is true

\begin{proposition}\label{pr:discrepancy}
For any $\alpha\in\root$, the following holds in $PU(\SS\otimes\IC[\weight])$
\begin{equation}\label{eq:discrepancy}
\pi(\zeta_{\alpha})=U_{\alpha}\cchi_{\alpha}
\end{equation}
where $\cchi_{\alpha}\in\End(\IC[\weight])$ acts as multiplication by the function
$\displaystyle{\cchi_{\alpha}(\mu)=\frac{(-1)^{\<\alpha,\mu\>}}{\wt\omega(\alpha,\mu)}}$.
\end{proposition}
\proof
By construction, the operators $U_{\mu}$, $\mu\in\weight$ satisfy the hypothesis of
proposition \ref{pr:extended cocycle}. Indeed,
$U_{\mu}h(n)U_{\mu}^{*}=h(n)-\delta_{n0}\<\mu,h\>$ and, by \eqref{eq:definition of vo},
$U_{\mu}X_{\alpha}(z)U_{\mu}^{*}=z^{-\<\alpha,\mu\>}X_{\alpha}(z)$. They therefore satisfy
\begin{equation}\label{eq:uno}
U_{\mu}\pi(\zeta_{\alpha})U_{\mu}^{*}\pi(\zeta_{\alpha})^{*}=(-1)^{\<\mu,\alpha\>}
\end{equation}
On the other hand, by \eqref{eq:U},
\begin{equation}
U_{\mu}U_{\alpha}U_{\mu}U_{\alpha}^{*}=\wt\omega(\mu,\alpha)
\end{equation}
and therefore \eqref{eq:uno} continues to hold if $\pi(\zeta_{\alpha})$ is replaced by
$W_{\alpha}=U_{\alpha}\cchi_{\alpha}$. It follows that the operators
$C_{\alpha}=W_{\alpha}\pi(\zeta_{\alpha})^{*}$ commute with the $U_{\mu}$. By proposition
\ref{pr:zeromodes affine}, they also commute with the action of $\lpol\t$ and therefore
act as scalars \halmos\\

The following useful result may be found in \cite[\S 4.3]{GO2}

\begin{lemma}\label{le:GO lemma}
Let $V$ be an irreducible $G$--module with Casimir $C_{V}$. If $\Pi(V)$ is the set of
weights of $V$ counted with multiplicities, then
\begin{equation}\label{eq:Casimir identity}
\frac{\dim V}{\sum_{\mu\in\Pi(V)}\|\mu\|^{2}}C_{V}=\frac{\dim G}{\rank G}
\end{equation}
In particular, if $G$ is simply--laced and $V$ is a minimal representation with highest
weight $\lambda$, the level 1 conformal dimension of $V$ is equal to
\begin{equation}\label{eq:dimension identity}
\Delta_{V}:=\frac{C_{V}}{2+C_{\g}}=\half{1}\|\lambda\|^{2}
\end{equation}
\end{lemma}
\proof
The map $X\otimes Y\rightarrow\tr_{V}(XY)$ is a symmetric, bilinear and ad--invariant form
on $\gc$ so that
\begin{equation}
\tr_{V}(XY)=\alpha_{V}\<X,Y\>
\end{equation}
where $\<\cdot,\cdot\>$ is the basic inner product and $\alpha_{V}\in\IC$. The proportionality
constant $\alpha_{V}$ may be evaluated in two different ways. Choosing basis $X_{k},X^{k}$ of
$\gc$ dual for $\<\cdot,\cdot\>$, we get
\begin{equation}
\dim(V)C_{V}=\tr(X_{k}X^{k})=\alpha_{V}\dim(G)
\end{equation}
On the other hand, if $h_{i},h^{i}$ are dual basis of $\tc$, we find by evaluating the trace
in a basis of weight vectors that
\begin{equation}
\sum_{\mu\in\Pi(V)}\|\mu\|^{2}=\tr(h_{i}h^{i})=\alpha_{V}\rank(G)
\end{equation}
Eliminating $\alpha_{V}$ yields \eqref{eq:Casimir identity}. 
If $V$ is minimal, all its weights lie on the orbit of $\lambda$ under the Weyl group and
have multiplicity one by proposition \ref{ch:classification}.\ref{pr:weights of minimal}.
Thus, $\sum_{\mu\in\Pi(V)}\|\mu\|^{2}=\dim(V)\|\lambda\|^{2}$.
Moreover, if $G$ is simply--laced, we find from \eqref{eq:Casimir identity}
\begin{equation}
C_{\g}=2\frac{\dim(G)-\rank(G)}{\rank(G)}
\end{equation}
and therefore \eqref{eq:dimension identity} \halmos\\

Let now $L_{0}$ be the infinitesimal generator of rotations given by the Segal--Sugawara
formula
\begin{equation}
L_{0}=
\frac{1}{\kappa}\Bigl(\half{1}X_{i}(0)X^{i}(0)+\sum_{m>0} X_{i}(-m)X^{i}(m)\Bigr)
\end{equation}
where $\kappa=1+\half{C_{\g}}$. Then
\begin{corollary}\label{co:zero nonzero sugawara}
The action of the infinitesimal generator of rotations $L_{0}$ on $\SS\otimes\IC[\weight]$
is given by
\begin{equation}\label{eq:zero nonzero sugawara}
L_{0}=\half{1}h_{k}h^{k}+\sum_{m>0}h_{k}(-m)h^{k}(m)
\end{equation}
where $h_{k},h^{k}$ are dual basis of $\tc$ for the basic inner product. In particular,
\begin{equation}
L_{0}\medspace 1\otimes\delta_{\mu}=\half{1}\<\mu,\mu\> 1\otimes\delta_{\mu}
\end{equation}
\end{corollary}
\proof
Let $D$ be operator defined by the right--hand side of \eqref{eq:zero nonzero sugawara}. It
satisfies by construction $[D,X(n)]=-nX(n)$ and therefore differs from $L_{0}$ by an additive
constant on each irreducible summand $\SS\otimes\IC[\mu+\root]$. Choosing $\mu$ dominant
minimal, $L_{0}$ acts on the corresponding highest weight vector $1\otimes\delta_{\mu}$ as
multiplication by the conformal dimension of the lowest energy subspace of $\IC[\mu+\root]$
and therefore, by the previous lemma, by $\half{\<\mu,\mu\>}$. Since
$D 1\otimes\delta_{\mu}=\half{\<\mu,\mu\>}$, the two operators coincide \halmos

\ssubsection{$LT$--equivariance of primary fields}\label{ss:LT equivariance}

Assuming that the level 1 primary fields extend to operator--valued distributions,
we derive below their commutation relations with $LT$. Although the continuity
properties of these fields will only be established in chapter \ref{ch:sobolev
fields}, the present discussion serves as motivation for the next subsection.\\

By proposition \ref{ch:classification}.\ref{existence}, the charges of the level 1
primary fields are necessarily admissible at level 1 and are therefore the minimal
$G$--modules since $G$ is simply--laced. Fix one such $V_{k}$ and let
$\phi_{ji}:\hfin_{i}\otimes V_{k}[z,z^{-1}]\rightarrow\hfin_{j}$ be a primary field
of charge $V_{k}$.
Since $\hfin_{i},\hfin_{j}$ are summands of $\hfin=\SS\otimes\IC[\weight]$, we regard
$\phi_{ji}$ as a map $\hfin\otimes V_{k}[z,z^{-1}]\rightarrow\hfin$ by extending it by
zero. By definition, it satisfies
\begin{equation}\label{eq:Lg equivariance}
[\pi(X),\phi_{ji}(f)]=\phi_{ji}(Xf)
\end{equation}
for any $X\in\lpol\g$ and $f\in V_{k}[z,z^{-1}]$. We assume in this subsection that
$\phi_{ji}$ extends to a jointly continuous map
$\hsmooth\otimes C^{\infty}(S^{1},V_{k})\rightarrow\hsmooth$. It follows by continuity
that \eqref{eq:Lg equivariance} holds on $\hsmooth$ for any $X\in L\g$ and
$f\in C^{\infty}(S^{1},V_{k})$ and therefore that, for any $\gamma\in LG$
\begin{equation}\label{eq:LG equivariance}
\pi(\gamma)\phi_{ji}(f)\pi(\gamma)^{*}=\phi_{ji}(\gamma f)
\end{equation}
To see this, it is sufficient to consider the case $\gamma=\exp_{LG}(X)$ since $LG$
is generated by the image of the exponential map. Let
$F(t)=e^{-t\pi(X)}\phi_{ji}(e^{tX}f)e^{t\pi(X)}\xi$ where $\xi\in\hsmooth$.
The $LG$--invariance of $\hsmooth$ (proposition \ref{ch:analytic}.\ref{invariance})
and \eqref{eq:Lg equivariance} imply that $F$ is differentiable and that $\dot F=0$
whence $F(1)\equiv F(0)=\phi_{ji}(f)\xi$.\\

We now restrict our attention to $\gamma\in LT$ and rephrase \eqref{eq:LG equivariance}
in terms of the generating function
$\phi_{ji}(z)=\sum_{n\in\IZ}\phi_{ji}(v,n)z^{-n-(\Delta_{i}+\Delta_{k}-\Delta_{j})}$
where, as customary $\phi_{ji}(v,n)=\phi_{ji}(v\otimes e^{in\theta})$ and the
$\Delta_{\cdot}$ are the conformal dimensions of $\H_{i}(0),V_{k}$ and $\H_{j}(0)$
respectively.
Let $v_{\mu}\in V_{k}$ be a vector of weight $\mu\in\weight$ and, for $n\in\IN$ set
$v_{\mu}(n)=e^{in\theta}\otimes v_{\mu}\in C^{\infty}(S^{1},V_{k})$. Clearly, for
$h\in\tc$ and $\tau\in T$,
\begin{xalignat}{3}
h(m)v_{\mu}(n)&=\<h,\mu\>v_{\mu}(n+m)&
&\text{and}&
\tau v_{\mu}(n)&=\mu(\tau)v_{\mu}(n)
\end{xalignat}
Moreover, since $\zeta_{\alpha}(\theta)=\exp(-i\alpha\theta)$, we find
$\zeta_{\alpha}v_{\mu}(n)=v_{\mu}(n-\<\alpha,\mu\>)$. Thus, in terms of the formal power
series $\phi_{ji}(v_{\mu},z)=\sum_{n\in\IZ}\phi(v_{\mu}(n))z^{-n}$, we have

\begin{align}
[h(n),\phi_{ji}(v_{\mu},z)]&=\<h,\mu\>z^{n}\phi_{ji}(v_{\mu},z)
\label{eq:cov 1}\\[1.2 ex]
\pi(\tau)\phi_{ji}(v_{\mu},z)\pi(\tau)^{*}&=\mu(\tau)\phi_{ji}(v_{\mu},z)
\label{eq:cov 2}\\[1.2 ex]
\pi(\zeta_{\alpha})\phi_{ji}(v_{\mu},z)\pi(\zeta_{\alpha})^{*}&=
z^{-\<\alpha,\mu\>}\phi_{ji}(v_{\mu},z)
\label{eq:cov 3}\\
\intertext{and, by \eqref{eq:dimension identity}}
[L_{0},\phi_{ji}(v_{\mu},z)]&=
(z\frac{d}{dz}+\Delta_{k})\phi_{ji}(v_{\mu},z)=
(z\frac{d}{dz}+\half{\<\mu,\mu\>})\phi_{ji}(v_{\mu},z)
\label{eq:cov 4}
\end{align}

\ssubsection{The construction of level 1 primary fields}\label{ss:construction}

We turn now to the construction of the primary fields $\phi_{ji}(z)$. It will be
convenient to consider all fields with a given charge $V_{k}$ at once by working
with $\Phi(z)=\bigoplus_{j,i}\phi_{ji}(z)$. Since $V_{k}$ is minimal, the spaces
$\Hom(\H_{i}(0)\otimes V_{k},\H_{j}(0))$ are at most one--dimensional by
corollary \ref{ch:classification}.\ref{co:one dim} and the
individual $\phi_{ji}(z)$ may be recovered by sandwiching $\Phi(z)$ between the
orthogonal projections $P_{j},P_{i}$ onto
$\hfin_{j},\hfin_{i}\subset\SS\otimes\IC[\weight]$ respectively.
Clearly, if $v_{\mu}\in V_{k}$ is a weight vector of weight $\mu$, 
$\Phi_{\mu}(z)=\Phi(v_{\mu},z)$ satisfies \eqref{eq:cov 1}--\eqref{eq:cov 4}.\\

To determine the form of $\Phi_{\mu}(z)$ from these commutation relations, fix
for any coset $\mu+\root\subset\weight$ a representative $[\mu]\in\weight$ and
set $[\alpha]=0$ for any $\alpha\in\root$. Let $\mu\in\weight$ and set
\begin{equation}\label{eq:ansatz field}
\Phi_{\mu}(z)=
E^{-}(\mu,z)
V_{\mu}z^{\mu+\half{\<\mu,\mu\>}}
\frac{\wt\omega(\mu,\cdot-[\cdot])}{e^{i\pi\<\mu,\cdot-[\cdot]\>}}
E^{+}(\mu,z)
\end{equation}
It is easy to check, using proposition \eqref{pr:discrepancy} that $\Phi_{\mu}$
satisfies \eqref{eq:cov 1}--\eqref{eq:cov 3}. Moreover, by corollary
\ref{co:zero nonzero sugawara},
\begin{equation}
[L_{0},\Phi_{\mu}(z)]=z\frac{d}{dz}\Phi_{\mu}(z)
\end{equation}
and we should therefore be considering $\Phi_{\mu}(z)z^{-\half{\<\mu,\mu\>}}$
instead. By analogy with the expressions for the vertex operators $X_{\alpha}(z)$,
we shall retain the erroneous $z^{\half{\<\mu,\mu\>}}$ instead.
We show below that, as $\mu$ spans the weights of $V_{k}$, the $\Phi_{\mu}(z)$
describe the components of a primary field of charge $V_{k}$. We shall need for
this purpose an explicit description of the infinitesimal action of $\g$ on any
minimal $G$--module.

\begin{lemma}\label{le:explicit action}
Let $\{\mu\}\subset\weight$ be the collection of weights of minimal length in a 
given $\root$--coset. Then, the subspace of $\SS\otimes\IC[\weight]$ spanned by
the vectors $1\otimes\delta_{\mu}$ is invariant and irreducible under the action
of $\g$ given by the operators $X_{\alpha}(0)$ and $h$. It is therefore isomorphic to
the minimal $G$--module $V_{\lambda}$ whose highest weight $\lambda$ is the unique
dominant element in $\{\mu\}$. The $\g$--action is explicitly given by
\begin{align}
h		\medspace 1\otimes\delta_{\mu}&=\<h,\mu\> 1\otimes\delta_{\mu}\\
X_{\alpha}(0)	\medspace 1\otimes\delta_{\mu}&=
\left\{\begin{array}{cll}
\epsilon(\mu,\alpha)1\otimes\delta_{\alpha+\mu}&
\text{if $\<\alpha,\mu\>=-1$}	&\text{{\it i.e.}~ iff   $\|\alpha+\mu\|=\|\mu\|$}\\
0&
\text{if $\<\alpha,\mu\>\geq 0$}&\text{{\it i.e.}~ iff   $\|\alpha+\mu\|>\|\mu\|$}
\end{array}\right.
\end{align}
\end{lemma}
\proof
Let $\mu$ be of minimal length so that $\|\mu\pm\alpha\|^{2}\geq\|\mu\|^{2}$ for any
root $\alpha$ and therefore $\<\mu,\alpha\>\in\{-1,0,1\}$. Proceeding as in the proof
of theorem \ref{th:vo reducibility}, we find
\begin{equation}\label{minimal weight vectors}
\begin{split}
X_{\alpha}(0)\thickspace 1\otimes\delta_{\mu}
&=\frac{1}{2\pi i}\oint\frac{dz}{z}X_{\alpha}(z)1\otimes\delta_{\mu}\\[1.5ex]
&=\frac{\epsilon(\mu,\alpha)}{2\pi i}\oint\frac{dz}{z}
z^{\<\alpha,\mu\>+\half{\<\alpha,\alpha\>}}
E^{-}(\alpha,z)1\otimes\delta_{\alpha+\mu}\\[1.5ex]
&=\left\{\begin{array}{cll}
\epsilon(\mu,\alpha)1\otimes\delta_{\alpha+\mu}&
\text{if $\<\alpha,\mu\>=-1$}&\text{{\it i.e.}~ iff   $\|\alpha+\mu\|=\|\mu\|$}\\
0&
\text{if $\<\alpha,\mu\>\geq 0$}&\text{{\it i.e.}~ iff   $\|\alpha+\mu\|>\|\mu\|$}
\end{array}\right.
\end{split}
\end{equation}
Thus, the span of $1\otimes\delta_{\mu}$ with $\mu$ varying in the vectors of minimal
length in a given $\root$--coset is invariant under the $X_{\alpha}(0)$. Since its weights
coincide by proposition \ref{ch:classification}.\ref{pr:weights of minimal} with those
of $V_{\lambda}$ the two representations have the same character and are therefore
isomorphic \halmos

\begin{theorem}\label{th:level 1 fields}
Let $V$ be a minimal $G$--module. Then $V$ has a basis of weight vectors $v_{\mu}$ such
that the assignement
\begin{equation}\label{eq:def of primary}
v_{\mu}\longrightarrow\Phi_{\mu}(z)=
E^{-}(\mu,z)
V_{\mu}z^{\mu+\half{\<\mu,\mu\>}}
\frac{\wt\omega(\mu,\cdot-[\cdot])}{e^{i\pi\<\mu,\cdot-[\cdot]\>}}
E^{+}(\mu,z)
\end{equation}
is the direct sum of all level 1 primary fields of charge $V$.
\end{theorem}
\proof
Set
$\displaystyle{\eta_{\mu}(\nu)=
 \frac{\wt\omega(\mu,\nu-[\nu])}{e^{i\pi\<\mu,\nu-[\nu]\>}}}$
and notice that $\eta_{\mu+\alpha}=\eta_{\mu}$ for any $\alpha\in\root$ since
$\wt\omega(\alpha,\beta)=(-1)^{\<\alpha,\beta\>}$ for any $\alpha,\beta\in\root$. Moreover,
$[\nu+\alpha]=[\nu]$ and $[\alpha]=0$ for $\alpha\in\root$ implies that
$\eta_{\mu}(\nu+\alpha)=\eta_{\mu}(\nu)\eta_{\mu}(\alpha)$. We now have

\begin{equation}
X_{\alpha}(z)\Phi_{\mu}(\zeta)=
\epsilon(\mu,\alpha)
\Bigl(1-\frac{\zeta}{z}\Bigr)^{\<\alpha,\mu\>}
\Bigl(\frac{z}{\zeta}\Bigr)^{\half{\<\alpha,\mu\>}}
V_{\alpha+\mu}
z^{\alpha+\half{\sqnm{\alpha}}+\half{\<\alpha,\mu\>}}
\zeta^{\mu+\half{\sqnm{\mu}}+\half{\<\alpha,\mu\>}}
\eta_{\mu}
\reg(\alpha,\mu,z,\zeta)
\end{equation}
where
\begin{equation}
\reg(\alpha,\mu,z,\zeta)=
E^{-}(\alpha,z)E^{-}(\mu,\zeta)E^{+}(\alpha,z)E^{+}(\mu,\zeta)=
\reg(\mu,\alpha,\zeta,z)
\end{equation}
On the other hand,
\begin{equation}
\Phi_{\mu}(\zeta)X_{\alpha}(z)=
\epsilon(\alpha,\mu)\eta_{\mu}(\alpha)
\Bigl(1-\frac{z}{\zeta}\Bigr)^{\<\alpha,\mu\>}
\Bigl(\frac{\zeta}{z}\Bigr)^{\half{\<\alpha,\mu\>}}
V_{\mu+\alpha}
z^{\alpha+\half{\sqnm{\alpha}}+\half{\<\alpha,\mu\>}}
\zeta^{\mu+\half{\sqnm{\mu}}+\half{\<\alpha,\mu\>}}
\eta_{\mu}
\reg(\alpha,\mu,z,\zeta)
\end{equation}
Thus, since
\begin{equation}
\epsilon(\alpha,\mu)\eta_{\mu}(\alpha)=
\epsilon(\mu,\alpha)\wt\omega(\alpha,\mu)\frac{\wt\omega(\mu,\alpha)}{(-1)^{\mu,\alpha}}=
\epsilon(\alpha,\mu)(-1)^{\<\alpha,\mu\>}
\end{equation}
we get

\begin{equation}
[X_{\alpha}(z),\Phi_{\mu}(\zeta)]=
\epsilon(\mu,\alpha)\mathcal D(z,\zeta)
V_{\mu+\alpha}
z^{\alpha+\half{\sqnm{\alpha}}+\half{\<\alpha,\mu\>}}
\zeta^{\mu+\half{\sqnm{\mu}}+\half{\<\alpha,\mu\>}}
\eta_{\mu}
\reg(\alpha,\mu,z,\zeta)
\end{equation}
where
\begin{equation}
\begin{split}
\mathcal D(z,\zeta)
&=
\Bigl(\frac{\zeta}{z}\Bigr)^{-\half{\<\alpha,\mu\>}}
\biggl[\Bigl(1-\frac{\zeta}{z}\Bigr)^{\<\alpha,\mu\>}-
(-1)^{\<\alpha,\mu\>}\Bigl(1-\frac{z}{\zeta}\Bigr)^{\<\alpha,\mu\>}
\Bigl(\frac{\zeta}{z}\Bigr)^{\<\alpha,\mu\>}\biggr]\\
&=
\begin{cases}
0&
\text{if $\<\alpha,\mu\>\geq 0$}\\
\Bigl(\frac{\zeta}{z}\Bigr)^{\half{1}}\delta\Bigl(\frac{\zeta}{z}\Bigr)&
\text{if $\<\alpha,\mu\>=-1$}
\end{cases}
\end{split}
\end{equation}

In the latter case, we get by \eqref{eq:delta identity}, the fact that $\reg$
contains only integral powers of $z,\zeta$ and that $\alpha$ has integral
eigenvalues

\begin{equation}
\begin{split}
[X_{\alpha}(z),\Phi_{\mu}(\zeta)]
&=\epsilon(\mu,\alpha)
\Bigl(\frac{\zeta}{z}\Bigr)^{\half{1}}\delta\Bigl(\frac{\zeta}{z}\Bigr)
V_{\mu+\alpha}
z^{\alpha+\half{\sqnm{\alpha}}+\half{\<\alpha,\mu\>}}
\zeta^{\mu+\half{\sqnm{\mu}}+\half{\<\alpha,\mu\>}}
\eta_{\mu}
\reg(\alpha,\mu,z,\zeta)\\
&=\epsilon(\mu,\alpha)\delta\Bigl(\frac{\zeta}{z}\Bigr)
V_{\mu+\alpha}
\zeta^{\alpha+\mu+\half{\sqnm{\alpha+\mu}}}
\eta_{\mu+\alpha}
\reg(\alpha,\mu,\zeta,\zeta)\\
&=\epsilon(\mu,\alpha)\delta\Bigl(\frac{\zeta}{z}\Bigr)
\Phi_{\alpha+\mu}(\zeta)
\end{split}
\end{equation}
The theorem now follows from lemma \ref{le:explicit action} \halmos\\

\remark The additional factor
$\displaystyle{{\wt\omega(\mu,\cdot-[\cdot])}{e^{-i\pi\<\mu,\cdot-[\cdot]\>}}}$ used in
the definition of the primary field is required to ensure that $\Phi_{\mu}(z)$ has the
correct commutation relations with the $X_{\alpha}(z)$. It can only be dispensed with
if $\weight$, or less ambitiously the intermediate lattice $\root\subset\Lambda\subset
\weight$ containing the collection of weights $\{\mu\}$ of a given primary field
possesses a skew--symmetric form $\wt\omega$ satisfying $\omega(\mu,\alpha)=(-1)
^{\<\mu,\alpha\>}$ whenever $\mu\in\Lambda$ and $\alpha\in\root$. Such a form need
not exist. For example, when $G=\SU_{2}$, the weight lattice is cyclic and therefore
does not possess non--trivial skew--symmetric forms. On the other hand, a suitable
$\wt\omega$ would satisfy
$\wt\omega(\alpha,\half{\alpha})=(-1)^{\<\alpha,\half{\alpha}\>}=-1$ where
$\alpha$ is the positive root of $\SU_{2}$ and $\half{\alpha}$ the corresponding 
generator of the weight lattice. In \S \ref{ss:extension}, we give the complete
list of intermediate lattices $\Lambda$ possessing such an $\omega$. We simply
note here that it comprises integral lattices, since in that case
$\wt\omega(\lambda,\mu)=(-1)^{\<\lambda,\mu\>+\sqnm{\lambda}\sqnm{\mu}}$ has the
required properties. This fact will be exploited below to study the vector primary
field for $L\Spin_{2n}$.

\ssubsection{The level 1 primary fields of $L\Spin_{2n}$}\label{ss:fermi field}

We consider now the case $G=\Spin_{2n}$ and rederive the result obtained in chapter
\ref{ch:fermionic} that the level 1 vector primary field for $L\Spin_{2n}$ is a
Fermi field. Its construction within the vertex operator model is an instance of
the celebrated Boson--Fermion correspondence.\\

\begin{lemma}\label{le:norm one}
Let $\root\subset\Lambda\subset\weight$ be an intermediate lattice which is integral
and choose $\wt\omega$ such that
\begin{equation}\label{eq:integral cocycle}
\wt\omega(\lambda,\mu)=(-1)^{\<\lambda,\mu\>+\sqnm{\lambda}\sqnm{\mu}}
\end{equation}
Then, if $\lambda,\mu\in\Lambda$ are of norm one
\begin{equation}
\{\Phi_{\lambda}(z),\Phi_{\mu}(\zeta)\}=
\left\{\begin{array}{cl}
\bigl(\frac{\zeta}{z}\bigr)^{-\lambda+\half{\sqnm{\lambda}}}
\delta\bigl(\frac{\zeta}{z}\bigr)
&\text{if $\lambda+\mu=0$}\\[1.2ex]
0&\text{otherwise}
\end{array}\right.
\end{equation}
\end{lemma}
\proof
Notice that the skew--symmetric form given by \eqref{eq:integral cocycle} is well--defined
and bilinear since $\Lambda$ is integral and therefore
$\|\sum\lambda_{i}\|^{2}=\sum\|\lambda_{i}\|^{2}\mod 2$. Moreover, since $\root$ is even,
$\omega(\mu,\alpha)=(-1)^{\<\mu,\alpha\>}$ whenever $\mu\in\Lambda$ and $\alpha\in\root$
and it follows that \eqref{eq:ansatz field} reduces to
\begin{equation}
\Phi_{\mu}(z)=E^{-}(\mu,z)V_{\mu}z^{\mu+\half{\sqnm{\mu}}}E^{+}(\mu,z)
\end{equation}
Thus, 
\begin{align}
\Phi_{\lambda}(z)\Phi_{\mu}(\zeta)
&=
\epsilon(\mu,\lambda)
\Bigl(1-\frac{\zeta}{z}\Bigr)^{\<\lambda,\mu\>}
\Bigl(\frac{z}{\zeta}\Bigr)^{\half{\<\lambda,\mu\>}}
V_{\lambda+\mu}
z^{\lambda+\half{\<\lambda,\lambda\>}+\half{\<\lambda,\mu\>}}
\zeta^{\mu+\half{\<\mu,\mu\>}+\half{\<\lambda,\mu\>}}
\reg(\lambda,\mu,z,\zeta)\\
\intertext{where}
\reg(\lambda,\mu,z,\zeta)&=
E^{-}(\lambda,z)E^{-}(\mu,\zeta)E^{+}(\lambda,z)E^{+}(\mu,\zeta)=
\reg(\mu,\lambda,\zeta,z)\\
\intertext{so that}
\Phi_{\mu}(\zeta)\Phi_{\lambda}(z)&=
\epsilon(\lambda,\mu)
\Bigl(1-\frac{z}{\zeta}\Bigr)^{\<\lambda,\mu\>}
\Bigl(\frac{\zeta}{z}\Bigr)^{\half{\<\lambda,\mu\>}}
V_{\lambda+\mu}
z^{\lambda+\half{\<\lambda,\lambda\>}+\half{\<\lambda,\mu\>}}
\zeta^{\mu+\half{\<\mu,\mu\>}+\half{\<\lambda,\mu\>}}
\reg(\lambda,\mu,z,\zeta)\\
\intertext{Since
$\epsilon(\lambda,\mu)=\epsilon(\mu,\lambda)
(-1)^{\<\lambda,\mu\>+\sqnm{\lambda}\sqnm{\mu}}=-(-1)^{\<\lambda,\mu\>}\epsilon(\mu,\lambda)$,
we get}
\{\Phi_{\lambda}(z),\Phi_{\mu}(\zeta)\}&=
\epsilon(\mu,\lambda)\D(z,\zeta)
V_{\lambda+\mu}
z^{\lambda+\half{\<\lambda,\lambda\>}+\half{\<\lambda,\mu\>}}
\zeta^{\mu+\half{\<\mu,\mu\>}+\half{\<\lambda,\mu\>}}
\reg(\lambda,\mu,z,\zeta)
\end{align}
where, using \eqref{eq:delta}
\begin{equation}
\begin{split}
\D(z,\zeta)
&=
\Bigl(\frac{\zeta}{z}\Bigr)^{-\half{\<\lambda,\mu\>}}
\biggl[
\Bigl(1-\frac{\zeta}{z}\Bigr)^{\<\lambda,\mu\>}
-(-1)^{\<\lambda,\mu\>}
\Bigl(1-\frac{z}{\zeta}\Bigr)^{\<\lambda,\mu\>}
\Bigl(\frac{z}{\zeta}\Bigr)^{-\<\lambda,\mu\>}
\biggr]\\
&=
\begin{cases}
0
&\text{if $\<\lambda,\mu\>\geq 0$}\\[1.5ex]
\Bigl(\frac{\zeta}{z}\Bigr)^{\half{1}}
\delta\Bigl(\frac{\zeta}{z}\Bigr)
&\text{if $\<\lambda,\mu\>=-1$}
\end{cases}
\end{split}
\end{equation}
Since $I$ is an integral lattice and $\lambda,\mu$ are of norm one,
$\<\lambda,\mu\>\in\{-1,0,1\}$ with $\<\lambda,\mu\>=-1$ iff $\mu=-\lambda$.
In the latter case, we get using \eqref{eq:delta identity} and the fact that
$\reg$ only contains integral powers of $z$ and $\zeta$,
\begin{equation}
 \delta\Bigl(\frac{\zeta}{z}\Bigr)\reg(\lambda,\mu,z,\zeta)=
 \delta\Bigl(\frac{\zeta}{z}\Bigr)\reg(\lambda,-\lambda,\zeta,\zeta)=1
\end{equation}
Moreover, since $\epsilon$ is normalised, $\epsilon(\lambda,-\lambda)=1$ and
$V_{0}=1$. Thus, if $\lambda=-\mu$,
\begin{equation}
\{\Phi_{\lambda}(z),\Phi_{-\lambda}(\zeta)\}=
\Bigl(\frac{\zeta}{z}\Bigr)^{\half{1}}
\delta\Bigl(\frac{\zeta}{z}\Bigr)
z^{\lambda}
\zeta^{-\lambda}
\end{equation}
as claimed \halmos\\

The following is well--known \cite[\S 7.1]{GO2}
\newpage

\begin{proposition}\label{pr:vector in vertex}\hfill
\begin{enumerate}
\item The level 1 primary fields of $L\Spin_{2n}$ whose charge is the vector
representation are Fermi fields.
\item All level 1 primary fields of $L\Spin_{8}$ are Fermi fields.
\end{enumerate}
\end{proposition}
\proof
(i)
The root lattice of $\Spin_{2n}$ is the $\IZ$--span in $\IR^{n}$ of the vectors
$\theta_{i}-\theta_{i+1}$, $i=1\ldots n-1$ and $\theta_{n-1}+\theta_{n}$ and is
therefore the lattice of integral points with even sum of coordinates. The dual
lattice $\weight$ is $I+\half{1}(\theta_{1}+\cdots+\theta_{n})$ where
$I=\IZ^{n}$ is the integral lattice. The minimal dominant weights are
$v=\theta_{1}$ and $s_{\pm}=\half{1}(\theta_{1}+\cdots+\theta_{n-1}\pm\theta_{n})$
and correspond to the vector and spin modules respectively. Let now
$\lambda\in v+\root=I$ of norm one so that $\lambda=\pm\theta_{i}$ and the
corresponding operator on $\SS\otimes\IC[\weight]$ has integral eigenvalues on
\begin{equation}
\F\fin_{NS}=\SS\otimes\IC[I]=
\SS\otimes\IC[\root]\medplus\SS\otimes\IC[v+\root]
\end{equation}
and half integral eigenvalues ({\it i.e.}~ elements of $\half{1}+\IZ$) on
\begin{equation}
\F\fin_{R}=\SS\otimes\IC[s_{+}+I]=
\SS\otimes\IC[s_{+}+\root]\medplus\SS\otimes\IC[s_{-}+\root]
\end{equation}
The notation $\F\fin_{NS}$ and $\F\fin_{R}$ refers to the fact that, by theorem
\ref{th:vo reducibility} and and proposition \ref{ch:fermionic}.\ref{pr:irreducibility}
these subspaces are isomorphic, as $\lpol\gc$--modules to the finite energy subspaces
of the Neveu--Schwarz and Ramond Fermionic Fock spaces constructed in chapter
\ref{ch:fermionic}. By \eqref{eq:delta identity}
\begin{equation}\label{eq:semi-integral}
\Bigl(\frac{\zeta}{z}\Bigr)^{-\lambda+\half{\sqnm{\lambda}}}
\delta\Bigl(\frac{\zeta}{z}\Bigr)
=
\left\{\begin{array}{rl}
\Bigl(\frac{\zeta}{z}\Bigr)^{\half{1}}\delta\Bigl(\frac{\zeta}{z}\Bigr)
&\text{on $\F\fin_{NS}$}\\[1.5ex]
\delta\Bigl(\frac{\zeta}{z}\Bigr)&\text{on $\F\fin_{R}$}
\end{array}\right.
\end{equation}
Moreover, the modes of $\Phi_{\lambda}$ clearly
preserve the splitting $\SS\otimes\IC[\weight]=\F\fin_{NS}\bigoplus\F\fin_{R}$.
Thus, defining, for any $i$,
$\Phi_{\theta_{i}}(z)=\psi^{\ns}_{i}(z)+\psi^{\ra}_{i}(z)$ where
\begin{align}
\psi^{\ns}_{i}(z)&=\sum_{n\in\half{1}+\IZ}\Phi_{\theta_{i}}(n)z^{-n}&
\psi^{\ra}_{i}(z)&=\sum_{n\in\IZ}\Phi_{\theta_{i}}(n)z^{-n}
\end{align}
and noticing that the formal adjunction property
$\Phi_{\theta_{i}}(z)^{*}=\Phi_{-\theta_{i}}(z)$ is inherited by $\psi^{\ns}_{i}(z)$
and $\psi^{\ra}_{i}(z)$ so that $\psi_{-i}^{\ns}(-n)=\psi_{i}^{\ns}(n)^{*}$, lemma
\ref{le:norm one} and \eqref{eq:semi-integral} imply that
\begin{equation}
\{\psi^{\ns}_{i}(n),\psi^{\ns}_{j}(m)\}=
0=
\{\psi^{\ra}_{i}(n),\psi^{\ra}_{j}(m)\}
\end{equation}
and
\begin{xalignat}{2}
\{\psi^{\ns}_{i}(n),\psi^{\ns}_{j}(m)^{*}\}
&=
\begin{cases}
\delta_{nm}\delta_{ij}&\text{on $\F\fin_{NS}$}\\
0&\text{on $\F\fin_{R}$}
\end{cases}&
\{\psi^{\ra}_{i}(n),\psi^{\ra}_{j}(m)^{*}\}
&=
\begin{cases}
0&\text{on $\F\fin_{NS}$}\\
\delta_{nm}\delta_{ij}&\text{on $\F\fin_{R}$}
\end{cases}
\end{xalignat}
and in particular $\psi^{\ns}_{i}(z)=0$ on $\F\fin_{R}$ and $\psi^{\ra}_{i}(z)=0$
on $\F\fin_{NS}$

(ii)
A similar analysis holds for all primary fields of $L\Spin_{8}$ since in that each
of the cosets $v+\root$, $s_{\pm}+\root$ are integral \halmos \\

\remark The same computation as that performed in the proof of theorem \ref{th:existence
of representation} shows that the components of the spin primary fields for $L\Spin_{16}$
give rise to two (non--commuting) action of $\lpol\mathfrak e_{8}$. Indeed, in this case
the highest weights of the spin representations, namely
$s_{\pm}=\half{1}(\theta_{1}+\cdots+\theta_{n-1}\pm\theta_{n})$ with $n=8$ have squared
length 2 and it follows that $s_{\pm}+\root$ are both even lattices which, by inspection
coincide with the root lattice of $E_{8}$.

\ssubsection{Appendix : intermediate lattices with the extension property}\label{ss:extension}

Let $G$ be simply--laced. We give below the list of intermediate lattices
$\root\subset\Lambda\subset\weight$ possessing the {\it extension property},
{\it i.e.}~ a skew--symmetric, bilinear form $\omega\in\Hom(\Lambda\wedge\Lambda,\T)$
satisfying $\omega(\alpha,\mu)=(-1)^{\<\alpha,\mu\>}$ whenever $\alpha\in\root$.
The existence of such a form simplifies the formulae \eqref{eq:ansatz field} giving the
primary fields whose charges have their weights in $\Lambda$. It also removes the
sole obstruction to extending the representation of $LG$ on the Hilbert space completion
of $\SS\otimes\IC[\weight]$ to one of the group of discontinuous loops whose
endpoints differ by an element of $Z(G)$ lying in the image of $\Lambda$ \cite{TL2}.

\begin{lemma}
If $\Lambda/\root$ is cyclic of order $k$ then $\Lambda$ has the extension
property if and only if $k\<\lambda,\lambda\>\in 2\IZ$ where $\lambda\in\Lambda$
is a generator. In this case, the extension is unique and takes values in
$\{\pm 1\}$.
\end{lemma}
\proof If $\omega:\Lambda\wedge\Lambda\rightarrow\T$ is an extension then
$k\lambda\in\root$ and therefore
$1=\omega(k\lambda,\lambda)=(-1)^{k\<\lambda,\lambda\>}$. Conversely,
if $k\<\lambda,\lambda\>\in 2\IZ$ then the skew-symmetric form
$\wt\omega(\alpha\oplus m\lambda,\beta\oplus n\lambda)=
 (-1)^{\<\alpha,\beta\>+n\<\alpha,\lambda\>+m\<\lambda,\beta\>}$ defined on
$\root\oplus\IZ\lambda$ descends to one on
$\root\oplus\IZ\lambda/-k\lambda\oplus k\lambda\cong\Lambda$ satisfying the
extension property \halmos

\begin{proposition}
The following is the complete list of subgroups $\{1\}\neq Z\subset Z(G)$, $G$
simply--laced such that $(2\pi i)^{-1}\exp^{-1}(Z)$ has the extension property.
\begin{enumerate}
\item Subgroups of $Z(\SU_{2n})$ of even index and all subgroups of $Z(\SU_{2n+1})$.
\item The subgroup $Z\subset Z(\Spin_{2n})$ such that $\Spin_{2n}/Z=\SO_{2n}$ and all
subgroups of $Z(\Spin_{8n})$.
\item $Z(E_{6})$
\end{enumerate}
\end{proposition}
\proof We proceed by enumeration according to the Lie theoretic of $G$. In what
follows, $\theta_{i}$, $i=1\ldots n$ and $\<\cdot,\cdot\>$ are the standard basis and
inner product in $\IR^{n}$, $I$ the self-dual lattice $\bigoplus_{i}\theta_{i}\IZ$
and $\Ie=\{\lambda\in I|\thinspace|\lambda|=\sum_{i}\lambda_{i}\in 2\IZ\}$.\\

${\bf SU_{n}}$\\
The roots are $\theta_{i}-\theta_{j}$, $i\neq j$ and span
$\root=\{\alpha\in I|\thinspace|\alpha|=0\}$.
The weight $\lambda=\theta_{1}-\frac{1}{n}(\theta_{1}+\cdots+\theta_{n})$
does not lie in $\root$ nor do $k\lambda$, $k=0\ldots n-1$ and therefore $\lambda$
generates $\weight/\root\cong Z(SU_{n})=\IZ_{n}$. Thus, for any $m|n$, the index
$m$ sublattice of $\weight$ is generated by $\root$ and $m\lambda$
and possesses the extension property if, and only if
$n/m\|m\lambda\|^{2}=m(n-1)\in 2\IZ$ and therefore iff $n$ is odd or $n$ and $m$
are even.\\

{\bf{Spin$_{\bf 2n}$}, ${\bf n\geq 3}$}\\
The simple roots are $\theta_{i}-\theta_{i+1}$, $i=1\ldots n-1$ and
$\theta_{n-1}+\theta_{n}$ and span $\root=\Ie$. The weight lattice is easily seen to
be $I+\half{1}(\theta_{1}+\cdots+\theta_{n})\IZ$ with minimal dominant weights given
by $v=\theta_{1}$ and $s_{\pm}=\half{1}(\theta_{1}+\cdots+\theta_{n-1}\pm\theta_{n})$.
We must distinguish two cases :

{\it $n$ even}. Then
$2s_{\pm}=
 (\theta_{1}+\theta_{2})+\cdots+(\theta_{n-1}\pm\theta_{n})=0$
mod $\root$ and $2v=2\theta_{1}=0$ mod $\root$ so that
$Z(\Spin_{2n})\cong\IZ_{2}\times\IZ_{2}$ and each of the minimal weights generates
a subgroup of order 2. In $SO_{2n}$, $\exp_{T}(-2\pi iv)=1$ and
therefore $\IZ_{2}^{v}$ possesses the extension property since
$2\|v\|^{2}=2$. The groups $\IZ_{2}^{s_{\pm}}$ have it if and only
if $n$ is a multiple of 4 since $2\|s_{\pm}\|^{2}=\half{n}$. 
Lastly, if $Z(\Spin_{2n})$ has the extension property then so
do $\IZ_{2}^{s_{\pm}}$ and therefore $n$ is a multiple of $4$. Conversely,
if that is the case, the skew-symmetric form
$\omega(\alpha\oplus pv\oplus qs_{+},\beta\oplus p^{\prime}v\oplus q^{\prime}s_{+})=
 i^{pq^{\prime}-qp^{\prime}}
 (-1)^{\<\alpha,\beta\>+\<\alpha,p^{\prime}v+q^{\prime}s_{+}\>+\<pv+qs_{+},\beta\>}$
defined on $\root\oplus v\IZ\oplus s_{+}\IZ$ descends to one on
$\weight=\root\oplus v\IZ\oplus s_{+}\IZ/
 (-2v\oplus 2v\oplus 0)\IZ+(-2s_{+}\oplus 0\oplus 2s_{+})\IZ$. In this case the extension
is unique only up to multiplication by the non-trivial $\IZ_{2}$-valued skew-form on
$\IZ_{2}\times\IZ_{2}$.

{\it $n$ odd}. Then
$2s_{\pm}=
 \theta_{1}+(\theta_{2}+\theta_{3})+\cdots+(\theta_{n-1}\pm\theta_{n})=
 v$ mod $\root$ and therefore $Z(\Spin_{2n})\cong\IZ_{4}$ with $v$ of order $2$ and
$s_{\pm}$ of order 4. Therefore $\root+v$ possesses the extension
property because $2\|v\|^{2}=2$ but $Z(\Spin_{2n})$ doesn't for
$4\|s_{+}\|^{2}=n\in 2\IZ+1$.\\

${\bf E_{6}}$\\
The simple roots of $E_{6}$ are
$\alpha_{1}=\half{1}(\theta_{1}-\theta_{2}-\cdots-\theta_{7}+\theta_{8}),
\alpha_{2}=\theta_{2}+\theta_{1},\alpha_{3}=\theta_{2}-\theta_{1},\ldots,
\alpha_{6}=\theta_{5}-\theta_{4}$.
The fundamental weight $\lambda_{1}=-\frac{2}{3}(\theta_{6}+\theta_{7}-\theta_{8})$
does not lie in the root lattice and therefore generates $\weight/\root\cong\IZ_{3}$
\cite[Planche V]{Bou}. Moreover, $3\|\lambda_{1}\|^{2}=4$.\\

${\bf E_{7}}$\\
The simple roots of $E_{7}$ are
$\alpha_{1}=\half{1}(\theta_{1}-\theta_{2}-\cdots-\theta_{7}+\theta_{8}),
\alpha_{2}=\theta_{2}+\theta_{1},\alpha_{3}=\theta_{2}-\theta_{1},\ldots,
\alpha_{7}=\theta_{6}-\theta_{5}$.
The fundamental weight
$\lambda_{2}=\half{1}(\theta_{1}+\cdots+\theta_{6}-2\theta_{7}+2\theta_{8})$
does not lie in the root lattice and therefore generates $\weight/\root\cong\IZ_{2}$
\cite[Planche VI]{Bou}. Moreover, $2\|\lambda_{2}\|^{2}=5$.\\

${\bf E_{8}}$\\
$E_{8}$ has trivial centre \cite[Planches VII]{Bou} and therefore no intermediate
lattices \halmos\\

\remark Notice that the above list agrees with the isomorphisms $\Spin_{6}\cong SU_{4}$.




\newcommand {\sbar}{\overline{s}}


\chapter{Analytic properties of primary fields}
\label{ch:sobolev fields}

This chapter is devoted to the study of the continuity properties of primary fields
of $L\Spin_{2n}$. These are required to construct (at first unbounded) explicit
intertwiners for the local loop groups by smearing the fields on functions supported
in complementary intervals. We show that any primary field
$\phi:\hfin_{i}\otimes V_{k}[z,z^{-1}]\rightarrow\hfin_{j}$ such that one of the
$\Spin_{2n}$--modules $\H_{i}(0),V_{k},\H_{j}(0)$ is minimal extends to a jointly
continuous operator--valued distribution
$\hsmooth_{i}\otimes C^{\infty}(S^{1},V_{k})\rightarrow\hsmooth_{j}$ satisfying as
expected
\begin{equation}\label{eq:smeared eq}
\pi_{j}(\gamma)\phi(f)\pi_{i}(\gamma)^{*}=\phi(\gamma f)
\end{equation}
for any $\gamma\in L\Spin_{2n}$.
When the charge $V_{k}$ is the vector representation, $\phi$ satisfies stronger
continuity properties and extends to a bounded map
$L^{2}(S^{1},V_{k})\rightarrow\B(\H_{i},\H_{j})$.\\

The level 1 result is obtained in section \ref{se:sobolev 1} from the bosonic
construction of the primary fields given in chapter \ref{ch:vertex operator}.
The level $\ell$ result is proved in section \ref{se:sobolev l} and follows
because the primary fields belonging to the above class may be obtained as
$\ell$--fold tensor products of level 1 primary fields.
Some care is required in checking this and the corresponding finite--dimensional
analysis is carried out in section \ref{se:level l for intertwiners}. The identity
\eqref{eq:smeared eq} is proved in section \ref{se:smeared inter}.

\ssection{Continuity of level 1 primary fields}
\label{se:sobolev 1}

As observed by Wassermann \cite{Wa5}, the continuity of the level 1 primary
fields depends upon the fact that they may be written as the product of
generating functions whose modes are bounded operators. We begin by studying
these. The notation follows chapter \ref{ch:vertex operator}.

\ssubsection{Norm boundedness of the vertex operators for $T\times\check T$}

Recall that the action of $\tc$ on $\IC[\weight]$ is given by
\begin{equation}\label{eq:reminder}
h \delta_{\mu}=\<h,\mu\>\delta_{\mu}
\end{equation}

\begin{lemma}\label{zero modes boundedness}
Let $\eta$ be a $\T$--valued function on $\weight$ and
\begin{equation}
X_{\lambda}(z)=
V_{\lambda}z^{\lambda+\half{\sqnm{\lambda}}}\eta(\cdot)=
\sum_{m}X_{\lambda}(m)z^{-m}
\end{equation}
Then, $X_{\lambda}(m)$ is bounded in norm by one.
\end{lemma}
\proof
By \eqref{eq:reminder}, $X_{\lambda}(m)=V_{\lambda}\eta(\cdot)P_{m}$ where $P_{m}$
is the orthogonal projection onto
\begin{equation}
\bigoplus_{\substack{\mu\in\weight\\ \<\mu,\lambda\>+\half{\sqnm{\lambda}}=-m}}
\IC\cdot\delta_{\mu}
\end{equation}
The claimed boundedness follows since $V_{\lambda}\eta(\cdot)$ is unitary \halmos

\ssubsection{Norm boundedness of pre-vertex operators of small conformal dimension}

This subsection follows \cite{Wa5}. Recall from \S \ref{ss:Stone for nonzero} of
chapter \ref{ch:vertex operator} the definition of the exponential operators
\begin{align}
E^{-}(\alpha,z)&=
\exp\Bigl(-\sum_{n<0}\frac{\alpha(n)}{n}z^{-n}\Bigr)&
E^{+}(\alpha,z)&=
\exp\Bigl(-\sum_{n>0}\frac{\alpha(n)}{n}z^{-n}\Bigr)
\end{align}
which are formal power series with coefficents in $\End(\SS)$ and the fact that
\begin{equation}
E^{+}(\alpha,z)E^{-}(\beta,\zeta)=
\Bigl(1-\frac{\zeta}{z}\Bigr)^{\<\alpha,\beta\>}
E^{-}(\beta,\zeta)E^{+}(\alpha,z)
\end{equation}
We have now

\begin{proposition}\label{non zero norm boundedness}
For any $\alpha\in i\t\cong\IR^{n}$ such that $\sqnm{\alpha}\leq 1$, the
modes $Y(n)$ of the pre--vertex operator
\begin{equation}
Y_{\alpha}(z)=
E^{-}(\alpha,z)E^{+}(\alpha,z)=
\exp\Bigl(-\sum_{n<0}\frac{\alpha(n)}{n}z^{-n}\Bigr)
\exp\Bigl(-\sum_{n>0}\frac{\alpha(n)}{n}z^{-n}\Bigr)
\end{equation}
satisfy
\begin{equation}\label{eq:bounded modes}
\|Y(n)\psi\|\leq\|\psi\|
\end{equation}
for any $\n\in\IZ$ and $\psi\in\SS$.
\end{proposition}
\proof
Consider first the case $\sqnm{\alpha}=1$. Then,
\begin{equation}
Y_{\alpha}(z)Y_{-\alpha}(\zeta)+
\Bigl(\frac{\zeta}{z}\Bigr)^{-1}Y_{-\alpha}(\zeta)Y_{\alpha}(z)=
\biggl[\Bigl(1-\frac{\zeta}{z}\Bigr)^{-1}+
\Bigl(1-\frac{z}{\zeta}\Bigr)^{-1}\Bigl(\frac{\zeta}{z}\Bigr)^{-1}\biggr]
\reg(\alpha,z,\zeta)
\end{equation}
where
\begin{equation}
\reg(\alpha,z,\zeta)=
E^{-}(\alpha,z)E^{-}(-\alpha,\zeta)E^{+}(\alpha,z)E^{+}(-\alpha,\zeta)=
\reg(-\alpha,\zeta,z)
\end{equation}
By (\ref{ch:vertex operator}.\ref{eq:delta}), the bracketed term is equal to
$\delta\Bigl(\frac{\zeta}{z}\Bigr)$ and since $\reg(\alpha,\zeta,\zeta)=1$,
we find by (\ref{ch:vertex operator}.\ref{eq:delta identity})
\begin{equation}
Y_{\alpha}(z)Y_{-\alpha}(\zeta)+
\Bigl(\frac{\zeta}{z}\Bigr)^{-1}Y_{-\alpha}(\zeta)Y_{\alpha}(z)=
\delta\Bigl(\frac{\zeta}{z}\Bigr)
\end{equation}
Writing $Y_{\alpha}(z)=\sum_{n}a_{n}z^{-n}$,
$Y_{-\alpha}(\zeta)=\sum_{n}b_{n}z^{-n}$ and recalling that the formal adjunction
property $Y_{\alpha}(z)^{*}=Y_{-\alpha}(z)$ implies that $b_{n}^{*}=a_{-n}$, we get
by taking modes on both sides
\begin{equation}
a_{n}a_{m}^{*}+a_{m-1}^{*}a_{n-1}=\delta_{n,m}
\end{equation}
and in particular
\begin{equation}
\|a_{n}\psi\|\leq\|\psi\|
\end{equation}

The general case may be settled by the following factorisation trick. The
space $V_{+}=z^{-1}\tc[z^{-1}]$ splits as $V_{+}^{1}\bigoplus V_{+}^{2}$ for
any orthogonal decomposition $\tc=\tc^{1}\oplus\tc^{2}$. Correspondingly,
\begin{equation}
\SS=\bigoplus_{k}S^{k}V_{+}=\SS^{1}\bigotimes\SS^{2}
\end{equation}
and
\begin{equation}
Y_{\lambda\oplus\lambda^\prime}(z)=Y_{\lambda}(z)\otimes Y_{\lambda^\prime}(z)
\end{equation}
for any $\lambda\in\tc^{1}$, $\lambda'\in\tc^{2}$. If $\sqnm{\alpha}\leq 1$,
we may find, by possibly enlarging $\tc$ if it is one--dimensional,
$i\t\ni\alpha'\perp\alpha$ such that $\sqnm{\alpha\oplus\alpha'}=1$.
The modes $Y_{\alpha\oplus\alpha'}(n)$ satisfy \eqref{eq:bounded modes} and are
equal to
\begin{equation}
\sum_{p+q=n} Y_{\alpha}(p)\otimes Y_{\alpha'}(q)
\end{equation}
Let $d=d_{1}+d_{2}$ be the infinitesimal generators of rotations on $\SS^{1},\SS^{2}$.
Let $\eta\in\SS^{2}$ be the lowest energy vector so that
$Y_{\alpha'}(z)\eta=E^{-}(\alpha',z)\eta$ and therefore
$Y_{\alpha'}(0)\eta=\eta$. If $\xi\in\SS^{1}$ is any eigenvector of $d_{1}$, then
\begin{equation}
\|\xi\|^{2}\|\eta\|^{2}=\|\xi\otimes\eta\|^{2}\geq
\|Y_{\alpha\oplus\alpha'}(n)\xi\otimes\eta\|^{2}=
\sum_{p'+q'=n}\|Y_{\alpha}(p')\xi\|^{2}\|Y_{\alpha'}(q')\eta\|^{2}\geq
\|Y_{\alpha}(n)\xi\|^{2}\|\eta\|^{2}
\end{equation}
since the vectors $Y_{\alpha}(p)\xi$, $Y_{\alpha}(p')\xi$ have different energies
and are therefore orthogonal for $p\neq p'$ and the same holds for
$Y_{\alpha}(q)\eta$ and $Y_{\alpha}(q')\xi$ whenever $q\neq q'$. Thus
\begin{equation}
\|Y_{\alpha}(n)\xi\|\leq \|\xi\|
\end{equation}
Lastly, if $\xi=\sum \xi_{n}$ is a sum of eigenvectors of $d_{1}$ with
distinct eigenvalues, then
\begin{equation}
\|Y_{\alpha}(m)\xi\|^{2}=
\sum\|Y_{\alpha}(m)\xi_{n}\|^{2}\leq
\sum\|\xi_{n}\|^{2}=\|\xi\|^{2}
\end{equation}
\halmos

\ssubsection{Sobolev estimates for products of norm bounded homogeneous fields}

This subsection follows \cite{Wa5}.
Let $\F_{i}$, $i=1\ldots k$ be inner product spaces supporting positive energy
representations $U_{\theta}^{i}=e^{i\theta d_{i}}$ of (a cover of) $\rot$. Call
a field $Y_{i}(z)=\sum_{n}Y_{i}(n)z^{-n}\in\End(\F_{i})[[z,z^{-1}]]$
{\it homogeneous} if $[d_{i},Y_{i}(n)]=-nY_{i}(n)$.

\begin{proposition}\label{factorisation}
Let $Y_{1}(z)\cdots Y_{k}(z)$ be homogeneous fields with uniformly bounded modes
acting on $\F_{1}\cdots\F_{k}$. Then, the modes of 
$Y(z)=Y_{1}(z)\otimes\cdots\otimes Y_{k}(z)$ satisfy
\begin{equation}
\|Y(m)\xi\|\leq C(1+|m|)^{k-1}\|(1+d)^{k-1}\xi\|
\end{equation}
for any $\xi\in\F=\F_{1}\otimes\cdots\otimes\F_{k}$ where
$d=d_{1}\otimes1\otimes\cdots\otimes1+
\cdots+1\otimes\cdots\otimes 1\otimes d_{k}$
\end{proposition}
\proof We have
\begin{equation}
Y(m)=\sum_{p_{1}+\cdots+p_{k}=m}Y_{1}(p_{1})\otimes\cdots\otimes Y_{k}(p_{k})
\end{equation}
Assume the $\F_{i}$ support positive energy representations of an $s$--sheeted cover
of $\rot$. If the right hand--side is applied to $\xi\in\F$ of energy $n\geq 0$, the
sum reduces to one involving only terms for which $p_{i}\leq n$ for any $i$. Their
number is bounded by that of solutions of $\sum p_{i}=m$, $p_{i}\leq n$,
$p_{i}\in s^{-1}\IZ$. However, $p_{i}=m-\sum_{j\neq i}p_{j}\geq m-n(k-1)$ and hence
for any $i$,
\begin{equation}
m-n(k-1)\leq p_{i}\leq n
\end{equation}
so that there are at most $s^{k-1}(1+nk-m)^{k-1}\leq s^{k-1}(1+nk+|m|)^{k-1}$
solutions since $p_{k}$ is determined once $p_{1},\ldots p_{k-1}$ are fixed.
It follows from $\|A\otimes B\|\leq\|A\|\|B\|$ for operators $A$, $B$
(Cauchy--Schwarz) that
\begin{equation}
\begin{split}
\|Y(m)\xi\|
&\leq Cs^{k-1}(1+nk+m)^{k-1}\|\xi\|\\
&\leq C'(1+|m|)^{k-1}(1+n)^{k-1}\|\xi\|\\
&=C'(1+|m|)^{k-1}\|(1+d)^{k-1}\xi\|
\end{split}
\end{equation}
where $C,C'$ are constants independent of $m$ and $\xi$ and we have used
$(1+d)\xi=(1+n)\xi$. The claimed inequality therefore holds for eigenvectors
of $d$ and hence for any $\xi\in\F$ since it is stable under taking orthogonal
sums of eigenvectors \halmos

\ssubsection{Continuity of the level 1 spin primary fields}

\begin{theorem}\label{th:spin sobolev}
Let $\phi_{s}:\hfin_{i}\otimes V_{s}[z,z^{-1}]\rightarrow\hfin_{j}$ be a
level 1 primary field of $L\Spin_{2n}$ whose charge is one of the spin
modules. Then $\phi_{s}$ extend to a jointly continuous map
$\hsmooth_{i}\otimes C^{\infty}(S^{1},V_{s})\rightarrow\hsmooth_{j}$.
\end{theorem}
\proof
Denote as customary by $\H_{i}^{t},\H_{j}^{t}$ the completion of $\hfin_{i},
\hfin_{j}$ with respect to the norm $\|\xi\|_{t}=\|(1+d)^{t}\xi\|$ where
$d$ is the infinitesimal generator of rotations.
The weights of the spin representations are of the form
$\mu=\half{1}(\epsilon_{1}\theta_{1}+\cdots+\epsilon_{n}\theta_{n})$
where $\epsilon_{i}\in\{\pm 1\}$ and may therefore be decomposed as orthogonal
sums of $\Delta=\left\lceil\frac{n}{4}\right\rceil$ vectors $\lambda_{i}$ with
$\sqnm{\lambda_{i}}\leq 1$. By theorem
\ref{ch:vertex operator}.\ref{th:level 1 fields}, the corresponding component
of $\phi_{s}$ factorises as
\begin{equation}
\Phi_{\mu}(z)=
Y_{\lambda_{1}}(z)\otimes\cdots\otimes Y_{\lambda_{\Delta}}(z)
\otimes V_{\mu}z^{\mu+\half{\sqnm{\mu}}}\eta(\cdot)
\end{equation}
for some $\T$--valued function $\eta$ on $\weight$. By lemma
\ref{zero modes boundedness} and proposition \ref{non zero norm boundedness},
all factors have modes bounded in norm by 1. Thus, if $\xi\in\hfin_{i}$ is of
energy $n\geq m$, proposition \ref{factorisation} yields
\begin{align}
\|\Phi_{\mu}(m)\xi\|_{t}
&=(1+n-m)^{t}\|\Phi_{\mu}(m)\xi\|\notag\\
 &\leq C(1+|m|)^{\Delta}(1+n-m)^{t}\|\xi\|_{\Delta}\notag\\
&=C(1+|m|)^{\Delta}\frac{(1+n-m)^{t}}{(1+n)^{t}}\|\xi\|_{\Delta+t}\\
&\leq C(1+|m|)^{\Delta+\half{|t|}}\|\xi\|_{\Delta+t}\notag\\
\intertext{whence, for any $\xi\in\hfin_{i}$}
\|\Phi_{\mu}(m)\xi\|_{t}
&\leq C(1+|m|)^{\Delta+|t|}\|\xi\|_{{\Delta+t}}
\end{align}
since $\xi$ may be written as an orthogonal sum of eigenvectors of $d$.
Next, if
\begin{equation}
f=\sum_{\substack{m\in\IZ\\ \mu\in\Pi(s)}}
  a_{m,\mu}v_{\mu}(m)\in V_{s}[z,z^{-1}]
\end{equation}
where $\Pi(s)$ is the set of weights of $V_{s}$, we have, by definition
$\phi_{s}(f)=\sum_{m,\mu}a_{m,\mu}\Phi_{\mu}(m)$. Then,
\begin{equation}\label{eq:spin estimate}
\|\phi_{s}(f)\xi\|_{t}\leq
C\sum_{m,\mu}|a_{m,\mu}|(1+|m|)^{\Delta+|t|}\|\xi\|_{\Delta+t}\leq
C' |f|_{\Delta+|t|}\|\xi\|_{\Delta+t}
\end{equation}
and $\phi_{s}$ extends to a continuous map
$C^{\infty}(S^{1},V_{s})\rightarrow\B(\H_{i}^{\Delta+t},\H_{j}^{t})$
\halmos\\

\remark Notice that the estimates \eqref{eq:spin estimate} are not quite
optimal. For example, for $L\Spin_{8}$, the spin primary fields are Fermi
fields by proposition \ref{ch:vertex operator}.\ref{pr:vector in vertex}
and therefore extend to bounded maps
$L^{2}(S^{1},V_{s_{\pm}})\rightarrow\B(\H_{i},\H_{j})$.

\ssection{Finite--dimensional $\Spin_{2n}$--intertwiners}
\label{se:level l for intertwiners}

By lemma \ref{ch:classification}.\ref{le:l lemma}, any $\Spin_{2n}$--module
admissible at level $\ell$ is contained in an $\ell$--fold tensor product
$V_{p_{1}}\otimes\cdots\otimes V_{p_{\ell}}$ where the $V_{p}$ are admissible
at level 1 and therefore minimal. In this section, we prove an analogous
factorisation for intertwiners $\phi:V_{i}\otimes V_{k}\rightarrow V_{j}$
when $V_{i},V_{k}, _{j}$ are admissible at level $\ell$ and one of them is
minimal. This will be used in the next section to show that level $\ell$
primary fields corresponding to the vertex $\vertex{V_{k}}{V_{j}}{V_{i}}$
can be written as the tensor product of level 1 primary fields.

\ssubsection{Minimal representations}

Recall from proposition \ref{ch:classification}.\ref{pr:weights of minimal} that
the weights of a minimal $G$--module $V$ lie on a single orbit of the Weyl group
and satisfy 
$\<\mu,\alpha^{\vee}\>\in\{1,0,-1\}$ for any root $\alpha$ and corresponding coroot
$\alpha^{\vee}=h_{\alpha}$.
If $\mu,\wt\mu$ are weights of $V$,
then $\wt\mu=\mu+\alpha_{1}+\cdots+\alpha_{k}$ where the $\alpha_{i}$ are possibly
repeated  roots. Since $\<\mu,\alpha_{i}\>\geq 0$ for all $i$ would imply
$\|\wt\mu\|>\|\mu\|$
there exists an $\alpha_{i}$ such that $\<\mu,\alpha_{i}^{\vee}\>=-1$ and therefore
$\|\mu+\alpha_{i}\|=\|\mu\|$. An iteration of this argument produces a permutation 
$\tau$ of $\{1,\ldots,k\}$ such that
$\|\mu+\alpha_{\tau(1)}+\cdots+\alpha_{\tau(j)}\|=\|\mu\|$ for any $j=1\ldots k$.
In representation theoretic terms, the weight spaces of a minimal representation
are one--dimensional and if $v_{\wt\mu}$, $v_{\mu}$ are eigenvectors corresponding
to the weights $\wt\mu$, $\mu$ then, up to a non--zero multiplicative constant,
$v_{\wt\mu}=e_{\alpha_{\tau(k)}}\cdots e_{\alpha_{\tau(1)}}v_{\mu}$. Indeed, by
elementary $\mathfrak{sl}_{2}(\IC)$ theory, for any $j\in\{1\ldots k\}$,
$e_{\alpha_{\tau(j)}}e_{\alpha_{\tau(j-1)}}\cdots e_{\alpha_{\tau(1)}}v_{\mu}\neq 0$
since $\<\alpha_{j}^{\vee},\mu+\alpha_{1}+\cdots+\alpha_{j-1}\>=-1$ and therefore
$v_{\wt\mu}$ and $e_{\alpha_{\tau(k)}}\cdots e_{\alpha_{\tau(1)}}v_{\mu}$ are proportional
since they lie in the same weight space.\\

If $V_{\delta},V_{\nu}$ are irreducible representations of highest weights $\delta,
\nu$ and $V_{\delta}$ is minimal, the tensor product $V_{\nu}\otimes V_{\delta}$
decomposes according to proposition \ref{ch:classification}.\ref{pr:tensor with minimal}
as
\begin{equation} \label{decompose}
V_{\nu}\otimes V_{\delta}=\bigoplus_{\wt\delta} V_{\nu+\wt\delta}
\end{equation}
where $\wt\delta$ varies among the weights of $V_{\delta}$ such that $\nu+\wt\delta$
is dominant.

\begin{lemma}\label{approx}
Let $V_{\nu}$, $V_{\delta}$ be irreducible representations with highest weights $\nu$,
$\delta$ and corresponding eigenvectors $v_{\nu}$, $v_{\delta}$. Assume $V_{\delta}$
is minimal so that $V_{\nu}\otimes V_{\delta}$ decomposes according to \eqref{decompose}.
Then, if $\wt\delta$ is a weight of $V_{\delta}$ such that $\nu+\wt\delta$ is dominant
and $v_{\wt\delta}$ is a corresponding non--zero weight vector, the orthogonal projection
of $v_{\nu}\otimes v_{\wt\delta}$ on $V_{\nu+\wt\delta}\subset V_{\nu}\otimes V_{\wt\delta}$
is non--zero.
\end{lemma}
\proof
Let $\Omega_{\nu+\wt\delta}\in V_{\nu+\wt\delta}\subset V_{\nu}\otimes V_{\wt\delta}$ be a
non--zero highest weight vector. We claim that
$(\Omega_{\nu+\wt\delta},v_{\nu}\otimes v_{\wt\delta})\neq 0$. Indeed, if the contrary
holds, we shall prove inductively that
\begin{equation}\label{vanish}
(\Omega_{\nu+\wt\delta},
 f_{\alpha_{\sigma(k)}}\cdots f_{\alpha_{\sigma(1)}}v_{\nu}\otimes
 e_{\alpha_{k}}\cdots e_{\alpha_{1}} v_{\wt\delta})=0
\end{equation}
where $\alpha_{1}\ldots\alpha_{k}$ are (possibly repeated) simple roots and $\sigma$ is any
permutation of $\{1,\ldots, k\}$. Since such vectors form a spanning set of the eigenspace
of $V_{\nu}\otimes V_{\delta}$ corresponding to the weight $\nu+\wt\delta$, it follows that
$\Omega_{\nu+\wt\delta}=0$, a contradiction. To prove the inductive claim, notice that
\begin{equation}\label{expand}
\begin{split}
 &
f_{\alpha_{\sigma(k)}}(
f_{\alpha_{\sigma(k-1)}}\cdots f_{\alpha_{\sigma(1)}}v_{\nu}\otimes
e_{\alpha_{k}}\cdots e_{\alpha_{1}}v_{\wt\delta})\\
=&
f_{\alpha_{\sigma(k)}}\cdots f_{\alpha_{\sigma(1)}}v_{\nu}\otimes
e_{\alpha_{k}}\cdots e_{\alpha_{1}} v_{\wt\delta}\\
-&
\<\alpha_{k}^{\vee},\wt\delta+(\alpha_{1}+\cdots+\alpha_{j-1})\>
f_{\alpha_{\sigma(k-1)}}\cdots f_{\alpha_{\sigma(1)}}v_{\nu}\otimes
e_{\alpha_{k}}\cdots \widehat{e_{\alpha_{j}}}\cdots e_{\alpha_{1}}v_{\wt\delta}\\
+&
f_{\alpha_{\sigma(k-1)}}\cdots f_{\alpha_{\sigma(1)}}v_{\nu}\otimes
e_{\alpha_{k}}\cdots e_{\alpha_{j}}f_{\alpha_{\sigma_{k}}}
e_{\alpha_{j-1}}\cdots e_{\alpha_{1}}v_{\wt\delta}
\end{split}
\end{equation}
where $j$ is the largest element of $\{1,\ldots,k\}$ such that $\alpha_{j}=\alpha_{\sigma(k)}$.
Now, if $f_{\alpha_{\sigma(k)}}\cdots f_{\alpha_{\sigma(1)}}v_{\nu}\otimes
e_{\alpha_{k}}\cdots e_{\alpha_{1}} v_{\wt\delta}$ is non--zero, then
$e_{\alpha_{j}}e_{\alpha_{j-1}}\cdots e_{\alpha_{1}}v_{\wt\delta}\neq 0$ and therefore by
elementary $\mathfrak{sl}_{2}(\IC)$ theory,
$f_{\alpha_{\sigma_{k}}}e_{\alpha_{j-1}}\cdots e_{\alpha_{1}}v_{\wt\delta}=0$
since $\wt\delta+(\alpha_{1}+\cdots+\alpha_{j-1})$ is a minimal weight and therefore
$\<\alpha_{j}^{\vee},\wt\delta+(\alpha_{1}+\cdots+\alpha_{j-1})\>\in\{1,0,-1\}$.
Therefore, if $f_{\alpha_{\sigma(k)}}\cdots f_{\alpha_{\sigma(1)}}v_{\nu}\otimes
e_{\alpha_{k}}\cdots e_{\alpha_{1}} v_{\wt\delta}\neq 0$, the last term on the right hand-side
of \eqref{expand} vanishes and by the inductive hypothesis and the fact that
$\Omega_{\nu+\wt\delta}$ is a highest weight vector, \eqref{vanish} holds \halmos\\

\ssubsection{Admissible intertwiners for $\Spin_{2n}$}

We begin by collecting some elementary facts about the spin representations of
$\Spin_{2n}$ which may be found in section \ref{ch:fermionic}.\ref{se:fd spinors}. The
complexified definining representation $V=\Vbox=\IC^{2n}$ has a natural bilinear,
non--degenerate symmetric form $B(\cdot,\cdot)$ which yields an isomorphism
$\so_{2n,\IC}\cong V\wedge V$ where the latter space acts on $V$ by
$u\wedge v\medspace w=uB(v,w)-v(u,w)$. If $f_{j}$, $0\neq j=-n\ldots n$ is an orthonormal
basis of $V$ satisfying $B(f_{j},f_{k})=\delta_{j+k,0}$, a basis for $\so_{2n,\IC}$
is given by the elementary matrices $F_{ij}=f_{i}\wedge f_{j}$. The Cartan subalgebra
corresponding to the block diagonal matrices in $\so_{2n,\IC}$ is then spanned by $F_{i,-i}$,
$i=1,\ldots n$. The roots of $\Spin_{2n}$ are $\theta_{k}+\theta_{l}$, $-n\leq k\neq\pm l
\leq n$ where the $\theta_{i}$, $i=1\ldots n$ are the dual basis to $F_{i,-i}$ and, by
definition $\theta_{-i}=-\theta_{i}$. The $\mathfrak{sl}_{2}(\IC)$--subalgebra
$\{e_{\alpha},f_{\alpha},h_{\alpha}\}$ of $\so_{2n,\IC}$ corresponding
to $\alpha=\theta_{k}+\theta_{l}$ is given by $\{F_{k,l},F_{k,-k}+F_{l,-l},F_{-l,-k}\}$.
The weight vectors in $V$ are exactly the $f_{j}$ since
$F_{i,-i}f_{j}=(\delta_{ij}-\delta_{i,-j})f_{i}=\theta_{j}(F_{i,-i})f_{j}$.\\

The spin representations are obtained via the representation of the Clifford algebra of $V$
generated by the $\IC$-linear symbols $\psi(u)$, $u\in V$ subject to the relations
$\psi(u)\psi(v)+\psi(v)\psi(u)=2B(u,v)$, on the exterior algebra $\Lambda V^{1,0}$ where
$V^{1,0}$ is the subspace spanned by the $f_{j}$, $j>0$. The action is given explicitly by
\begin{equation}\label{eq:action of}
\psi(u) v_{1}\wedge\cdots\wedge v_{k}=
\left\{\begin{array}{rl}
\sqrt{2}
u\wedge v_{1}\wedge\cdots\wedge v_{k}&\text{if $u\in V^{1,0}$}\\[1.2 ex]
\sqrt{2} \sum_{j} (-1)^{j+1}B(u,v_{j})
v_{1}\wedge\cdots\wedge\widehat{v_{j}}\wedge\cdots\wedge v_{k}&\text{if $u\in V^{0,1}$}
\end{array}\right.
\end{equation}

The representation is obtained by letting the Lie algebra element $u\wedge v$ act as
$\frac{1}{4}(\psi(u)\psi(v)-\psi(v)\psi(u))=\half{1}(\psi(u)\psi(v)-B(u,v))$. The even
and odd parts of the exterior algebra are clearly invariant under this action and are
irreducible. Their weights are easily read from
\begin{equation}
F_{i,-i} f_{j_{1}}\wedge\cdots\wedge f_{j_{k}}=
(\psi(f_{i})\psi(f_{-i})-\half{1})f_{j_{1}}\wedge\cdots\wedge f_{j_{k}}=
\left\{\begin{array}{rl}
 \half{1}f_{j_{1}}\wedge\cdots\wedge f_{j_{k}}&
 \text{if $i\in\{j_{1},\ldots,j_{k}\}$}\\[1.2 ex]
-\half{1}f_{j_{1}}\wedge\cdots\wedge f_{j_{k}}&
 \text{if $i\notin\{j_{1},\ldots,j_{k}\}$}
\end{array}\right.
\end{equation}
so that the vector $f_{J}=f_{j_{1}}\wedge\cdots\wedge f_{j_{k}}$ corresponds to the weight
$-\half{1}\sum\theta_{i}+\sum_{p}\theta_{j_{p}}$. The highest weight of the half of the
exterior algebra containing the top exterior power $\Lambda^{n} V^{1,0}$ has therefore
highest weight $s_{+}=\half{1}(\theta_{1}+\cdots+\theta_{n})$ and the other has highest
weight $s_{-}=\half{1}(\theta_{1}+\cdots+\theta_{n-1}-\theta_{n})$.\\

Finally, the Clifford multiplication map $V\otimes\Lambda V^{1,0}\rightarrow\Lambda V^{1,0}$
given by $v\otimes f_{J}\rightarrow\psi(v)f_{J}$ commutes with the action of $\Spin_{2n}$
and gives therefore rise to two intertwiners $V\otimes V_{s_{\pm}}\rightarrow V_{s_{\mp}}$.

\begin{lemma}\label{fd intertwiners}
Let $V_{i},V_{j},V_{k}$ be irreducible representations of $G=\Spin_{2n}$, one of which
is minimal so that $\Hom_{G}(V_{i}\otimes V_{k},V_{j})$ is at most one--dimensional. If
all are admissible at level $\ell$ and $\Hom_{G}(V_{i}\otimes V_{k},V_{j})=\IC$, there 
exist minimal $G$--modules $V_{i_{p}},V_{j_{p}},V_{k_{p}}$, $p=1\ldots\ell$ and
intertwiners $\phi_{p}\in\Hom_{G}(V_{i_{p}}\otimes V_{k_{p}},V_{j_{p}})$ such that
\begin{xalignat}{3}
V_{i}&\subset\bigotimes_{p=1}^{\ell}V_{i_{p}}&
V_{j}&\subset\bigotimes_{p=1}^{\ell}V_{j_{p}}&
V_{k}&\subset\bigotimes_{p=1}^{\ell}V_{k_{p}}
\end{xalignat}
and the corresponding restriction of $\otimes_{p}\phi_{p}$ to an intertwiner
$V_{i}\otimes V_{k}\rightarrow V_{j}$ is non--zero.
Moreover, if $V_{k}$ is minimal then one may choose $V_{i_{p}}=V_{j_{p}}$ and $V_{k_{p}}=\IC$
for $p=1\ldots\ell-1$ and $V_{k_{\ell}}=V_{k}$ so that $\otimes_{p}\phi_{p}$ is of the form
$1\otimes\cdots\otimes 1\otimes\phi_{p}$.
\end{lemma}
\proof
Up to a permutation of the modules, we may assume that $V_{k}$ is minimal and
that $\<\mu,\theta\>\leq\<\lambda,\theta\>$ where $\mu$, $\lambda$ are the
highest weights of $V_{i}$ and $V_{j}$ respectively and $\theta$ is the highest
root. If $\<\mu,\theta\>\leq\ell-1$, $V_{i}$ is contained, by lemma
\ref{ch:classification}.\ref{le:l lemma} in some tensor product
$V_{i_{1}}\otimes\cdots\otimes V_{i_{\ell-1}}$ with minimal factors. Then,
\begin{align}
V_{i}&\subset V_{i_{1}}\otimes\cdots\otimes V_{i_{\ell-1}}\otimes \IC\\[1.5 ex]
V_{k}&\subset \IC      \otimes\cdots\otimes \IC		  \otimes V_{k}\\[1.5 ex]
V_{j}&\subset V_{i_{1}}\otimes\cdots\otimes V_{i_{\ell-1}}\otimes V_{k}
\end{align}
and $\phi=1\otimes\cdots\otimes 1\otimes 1$ clearly restricts to a non--zero intertwiner.
Assume now that $\<\mu,\theta\>=\<\lambda,\theta\>=\ell\geq 2$. If $V_{k}$ is the trivial
representation then $V_{i}=V_{j}$ and, by lemma \ref{ch:classification}.\ref{le:l lemma},
$V_{i}\subset V_{i_{1}}\otimes\cdots\otimes V_{i_{\ell}}$ for some minimal $V_{i_{p}}$
and the lemma holds. We shall treat the cases when $V_{k}$ is the vector or one of the
spin representations separately.\\

$\mathbf{V_{k}=V_{s_{\pm}}}$.\\
Let $s$ and $\sbar$ be the highest weights of $V_{k}$ and of the other spin representation
so that $s=s_{\pm}=\half{1}(\theta_{1}+\cdots+\theta_{n-1}\pm\theta_{n})$ and
$\sbar=s_{\mp}$. We have $\lambda=\mu+\sigma$ where
$\sigma=\half{1}(\epsilon_{1}\theta_{1}+\cdots+\epsilon_{n}\theta_{n})$ for some
$\epsilon_{i}\in\{\pm 1\}$ is a weight of $V_{k}=V_{s}$. Up to a permutation of $V_{i}$,
$V_{j}$, we may assume that $\mu_{1}=\lambda_{1}+\half{1}$. Since
$\lambda_{1}+\lambda_{2}=\<\lambda,\theta\>=\<\mu,\theta\>=\mu_{1}+\mu_{2}$, we have
$\mu_{2}=\lambda_{2}-\half{1}$ and therefore
$\mu_{1}-\mu_{2}=\lambda_{1}-\lambda_{2}+1\geq 1$ so that $\rho=\mu-\theta_{1}$ is a
dominant weight. Let $V_{\rho}$ be the corresponding irreducible representation. Since
$\<\rho,\theta\>=\ell-1$, $V_{\rho}$ is contained in some tensor product 
$V_{i_{1}}\otimes\cdots\otimes V_{i_{\ell-1}}$ with minimal factors. Thus, by
\eqref{decompose}
\begin{align}
V_{i} \subset V_{\rho}\otimes\Vbox
     &\subset V_{i_{1}}\otimes\cdots\otimes V_{i_{\ell-1}}\otimes\Vbox\\[1.5 ex]
V_{k}&\subset\IC\otimes\cdots\otimes\IC\otimes V_{s}\\[1.5 ex]
V_{j} \subset V_{\rho}\otimes V_{\sbar}
     &\subset V_{i_{1}}\otimes\cdots\otimes V_{i_{\ell-1}}\otimes V_{\sbar}
\end{align}
since $\mu=\rho+\theta_{1}$ and $\lambda=\mu+\sigma=\rho+\sigma'$ where
$\sigma'=\theta_{1}+\sigma$ is a weight of $V_{\sbar}$.
If $\psi:\Vbox\otimes V_{s}\rightarrow V_{\sbar}$ is Clifford multiplication, we claim
that the intertwiner $\phi=1\otimes\cdots\otimes 1\otimes\psi$ has a non--zero restriction
to $V_{i}\otimes V_{k}\rightarrow V_{j}$.
To see this, notice that the highest weight vector $v_{\mu}$ in
$V_{i}\subset V_{\rho}\otimes \Vbox$
is the product $v_{\rho}\otimes v_{\theta_{1}}=v_{\rho}\otimes f_{1}$ of the corresponding
highest weight vectors.
If $v_{\sigma}\in V_{s}$ is of weight $\sigma$ so that
$v_{\sigma}=\wedge_{j:\sigma_{j}=\half{1}}f_{j}$, then
$\phi(v_{\mu}\otimes v_{\sigma})=
 v_{\rho}\otimes\psi(f_{1}\otimes v_{\sigma})=
 v_{\rho}\otimes f_{1}\wedge v_{\sigma}$. Since $f_{1}\wedge v_{\sigma}\in V_{\sbar}$
is of weight $\sigma+\theta_{1}$ and $\lambda=\rho+(\sigma+\theta_{1})$, lemma
\ref{approx} implies that $\phi(v_{\mu}\otimes v_{\sigma})$ has a non--zero projection
on $V_{j}\subset V_{\rho}\otimes V_{\sbar}$ whence the conclusion.\\

$\mathbf{V_{k}=\Vbox}$.\\
Up to a permutation of $V_{i}$ and $V_{j}$, we may assume that $\lambda$ is obtained
from $\mu$ by adding a box to the corresponding Young diagram, {\it i.e.}~ $\lambda=\mu+\theta_{j}$
where $j\geq 3$ since $\<\lambda,\theta\>=\<\mu,\theta\>$.
Let
\begin{equation}
\sigma=\half{1}(\theta_{1}+\cdots+\theta_{j-1}-\theta_{j}-\cdots-\theta_{n})
\end{equation}
so that it is a weight of $V_{s}$ with $s=s_{\pm}$ according to the parity of $n-j+1$.
As is readily verified, $\rho=\mu-\sigma$ is dominant and satisfies $\<\rho,\theta\>=\ell-1$.
Thus, if $V_{\rho}$ is the corresponding highest weight representation, then
$V_{\rho}\subset V_{i_{1}}\otimes\cdots\otimes V_{i_{\ell-1}}$ where the $V_{i_{p}}$
are minimal. It follows by \eqref{decompose} that
\begin{align}
V_{i}\subset V_{\rho} \otimes V_{s}
     &\subset V_{i_{1}}\otimes\cdots\otimes V_{i_{\ell-1}}\otimes V_{s}\\[1.5 ex]
V_{k}&\subset \IC      \otimes\cdots\otimes \IC		  \otimes \Vbox\\[1.5 ex]
V_{j}\subset V_{\rho} \otimes V_{\sbar}
     &\subset V_{i_{1}}\otimes\cdots\otimes V_{i_{\ell-1}}\otimes V_{\sbar}
\end{align}
since $\mu=\rho+\sigma$ and $\lambda=\mu+\theta_{j}=\rho+\sigma+\theta_{j}$ where
$\sigma$, $\sigma+\theta_{j}$ are weights of $V_{s}$ and $V_{\sbar}$ respectively.
Let $\psi:V_{s}\otimes\Vbox\rightarrow V_{\sbar}$ be Clifford multiplication and
$\phi=1\otimes\cdots\otimes 1\otimes\psi$.
We claim that if $v_{\mu}\in V_{i}\subset V_{\rho}\otimes V_{s}$ is the highest
weight vector then $\phi(v_{\mu}\otimes f_{j})$ is a non--zero highest weight
vector in $V_{j}\subset V_{\rho}\otimes V_{\sbar}$ so that $\phi$ restricts to
a non--zero intertwiner $V_{i}\otimes V_{k}\rightarrow V_{j}$.
To prove our claim, assume that $\phi(v_{\mu}\otimes f_{j})=0$.
Let $v_{\rho}\in V_{\rho}$ be the highest weight vector and
$v_{\sigma}=\wedge_{i=1\ldots j-1}f_{i}\in V_{s}$ a vector
of weight $\sigma$ so that $\psi(f_{-j})v_{\sigma}=0$. Then,
\begin{equation}
\begin{split}
0
&=(\phi(v_{\mu}\otimes f_{j}),v_{\rho}\otimes\psi(f_{j})v_{\sigma})\\[1.4 ex]
&=(1\otimes\psi(f_{j})v_{\mu},v_{\rho}\otimes\psi(f_{j})v_{\sigma})\\[1.4 ex]
&=(v_{\mu},v_{\rho}\otimes\psi(f_{-j})\psi(f_{j})v_{\sigma})	   \\[1.4 ex]
&=2(v_{\mu},v_{\rho}\otimes v_{\sigma})
\end{split}
\end{equation}
in contradiction with lemma \ref{approx}.
Let now $e_{\alpha}\in\so_{2n,\IC}$ be a root vector corresponding to the
simple root $\alpha$. Then,
$e_{\alpha}\phi(v_{\mu}\otimes f_{j})=\phi(v_{\mu}\otimes e_{\alpha}f_{j})$.
Since $e_{\alpha}f_{j}$ has weight $\alpha+\theta_{j}$, this does not vanish
only if $\alpha=\theta_{j-1}-\theta_{j}$ and in that case is  
proportional to $\phi(v_{\mu}\otimes f_{j-1})$. To conclude, it is therefore
sufficient to prove that $\phi(v_{\mu}\otimes f_{j-1})=0$.
To this end, notice that the weight spaces of $V_{s}$ are one--dimensional and
$v_{\mu}$ may therefore be written as
\begin{equation}
v_{\mu}=\sum_{} v_{\rho-(\sigma'-\sigma)}\otimes v_{\sigma'}
\end{equation}
where $\sigma'$ ranges over all weights of $V_{s}$ differing from $\sigma$ by
a sum of positive roots, $v_{\sigma'}=\wedge_{i:\sigma'_{i}=\half{1}}f_{i}$
and the $v_{\rho-(\sigma'-\sigma)}\in V_{\rho}$ have weight $\rho-(\sigma'-\sigma)$ .
Since such $\sigma'$ are of the form 
$\sigma'=
\half{1}(\theta_{1}+\cdots+\theta_{j-1}-\epsilon_{j}
         \theta_{j}-\cdots-\epsilon_{n}\theta_{n})$
where $\epsilon$ ranges over all even sign changes of the variables
$\theta_{j},\ldots,\theta_{n}$, we find that
$f_{j-1}\wedge v_{\sigma'}=0$ and therefore
\begin{equation}
\phi(v_{\mu}\otimes f_{j-1})=1\otimes\psi(f_{j-1})v_{\mu}=0
\end{equation}
\halmos

\ssection{Continuity of higher level primary fields}
\label{se:sobolev l}

\begin{theorem}\label{th:minimal is sobolev}
Let $\H_{i},\H_{j}$ be irreducible positive energy representations at level $\ell$ with
lowest energy subspaces $V_{i},V_{j}$ and $V_{k}$ an irreducible $\Spin_{2n}$--module
admissible at level $\ell$. If one of $V_{i},V_{j},V_{k}$ is minimal, the (projectively
unique) primary field $\psi:\hfin_{i}\otimes V_{k}[z,z^{-1}]\rightarrow \hfin_{j}$
extends to a jointly continuous operator--valued distribution
$\H_{i}^{\infty}\otimes C^{\infty}(S^{1},V_{k})\rightarrow\H_{j}^{\infty}$.
If, in addition, $V_{k}\cong\IC^{2n}$ is the
vector representation, $\psi$ extends to a bounded map
$L^{2}(S^{1},V_{k})\rightarrow\B(\H_{i},\H_{j})$.
\end{theorem}
\proof
By corollary \ref{ch:classification}.\ref{co:one dim},
$D=\dim(\Hom_{\Spin_{2n}}(V_{i}\otimes V_{k},V_{j}))\leq 1$. If $D=0$, then
$\psi=0$ and there is nothing to prove. If, on the other hand $D=1$, then by
lemma \ref{fd intertwiners}, we have
\begin{xalignat}{3}
V_{i}&\subset\bigotimes_{p=1}^{\ell} V_{i_{p}}&
V_{j}&\subset\bigotimes_{p=1}^{\ell} V_{j_{p}}&
V_{k}&\subset\bigotimes_{p=1}^{\ell} V_{k_{p}}
\end{xalignat}
where the $V_{i_{p}},V_{j_{p}},V_{k_{p}}$ are minimal and therefore admissible at
level 1. Moreover, there exist intertwiners
$\varphi_{p}\in\Hom_{\Spin_{2n}}(V_{i_{p}}\otimes V_{k_{p}},V_{j_{p}})$
such that $\otimes_{p}\varphi_{p}$ restricts to a non--zero intertwiner
$\varphi:V_{i}\otimes V_{k}\rightarrow V_{j}$.
Let $\H_{i_{p}},\H_{j_{p}}$ be the irreducible level 1 positive energy representations
with lowest energy subspaces $V_{i_{p}}$ and $V_{j_{p}}$ respectively.
The lowest energy suspaces of the level $\ell$ positive energy representations
\begin{xalignat}{2}
\H_{i_{1},\ldots,i_{\ell}}&=\bigotimes_{p=1}^{\ell}\H_{i_{p}}&
\H_{j_{1},\ldots,j_{\ell}}&=\bigotimes_{p=1}^{\ell}\H_{j_{p}}
\end{xalignat}
contain $V_{i}$ and $V_{j}$ respectively and therefore, by lemma
\ref{ch:classification}.\ref{le:generation},
$\H_{i_{1},\ldots,i_{\ell}}$ and $\H_{j_{1},\ldots,j_{\ell}}$ contain as submodules
$\H_{i}$ and $\H_{j}$ respectively. Denote by $P_{i}$ and $P_{j}$ the corresponding
orthogonal projections.
Let $\phi_{p}(z):\H_{i_{p}}\otimes V_{k_{p}}[z,z^{-1}]\rightarrow\H_{j_{p}}$ be the
level 1 primary field with initial term $\varphi_{p}$. Then
\begin{equation}
\phi(z)=P_{j}\phi_{1}(z)\otimes\cdots\otimes\phi_{\ell}(z)P_{i}
\end{equation}
is easily seen to be
a primary field of type $\H_{i}\otimes V_{k}[z,z^{-1}]\rightarrow\H_{j}$ when
$V_{k}[z,z^{-1}]$ is embedded in $V_{k_{1}}[z,z^{-1}]\otimes\cdots\otimes V_{k_{\ell}}
[z,z^{-1}]$ in the obvious way. Since the initial term of $\phi$ is $\varphi$,
$\phi\neq 0$ and therefore, up to a non--zero multiplicative constant, $\psi=\phi$.
We proceed now as in the proof of theorem \ref{th:spin sobolev}. By theorem
\ref{ch:vertex operator}.\ref{th:level 1 fields}, each $\phi_{p}$ is a product of vertex
operators with uniformly bounded modes and therefore so is $\phi$. It follows, by
proposition \ref{factorisation} that $\phi$ extends to a jointly continuous
bilinear map $C^{\infty}(S^{1},V_{k})\otimes\H_{i}^{\infty}\rightarrow\H_{j}^{\infty}$.
Finally, if $V_{k}$ is the vector representation we may, by lemma \ref{fd intertwiners},
choose $\phi(z)$ of the form
\begin{equation}
P_{j}1\otimes\cdots\otimes1\otimes\Psi(z)P_{i}
\end{equation}
where $\Psi$ is a level 1 vector primary field and is therefore continuous for
the $L^{2}$--norm by proposition \ref{ch:fermionic}.\ref{pr:L2 vector} \halmos\\

\remark Let $\phi,\phi^{*}$ be a primary field and its adjoint with corresponding vertices
$\vertex{V_{k}}{V_{j}}{V_{i}}$ and $\vertex{V_{k}^{*}}{V_{i}}{V_{j}}$. If one of $V_{i},
V_{j},V_{k}$ is minimal then, by theorem \ref{th:minimal is sobolev}, the defining identity
for $\phi^{*}$, namely
\begin{equation}
(\phi(f)\xi,\eta)=(\xi,\phi^{*}(\overline{f})\eta)
\end{equation}
where $\xi\in\hfin_{i}$, $\eta\in\hfin_{j}$ and $f\in V_{k}[z,z^{-1}]$, extends to
$\xi\in\hsmooth_{i},\eta\in\hsmooth_{j}$ and $f\in C^{\infty}(S^{1},V_{k})$. In particular,
\begin{xalignat}{3}
\phi^{*}(\overline{f})&\subseteq \phi(f)^{*}&
&\text{and}&
\phi(f)&\subseteq{\phi^{*}(\overline{f})}^{*}
\end{xalignat}
so that, for any $f\in C^{\infty}(S^{1},V_{k})$, $\phi(f)$ and $\phi^{*}(\overline{f})$ are
densely defined, closeable operators.

\ssection{Intertwining properties of smeared primary fields}
\label{se:smeared inter}

\begin{proposition}\label{pr:LG equivariance}
Let $\phi:\hfin_{i}\otimes V_{k}[z,z^{-1}]\rightarrow\hfin_{j}$ be a primary field
extending to an operator--valued distribution. Then, the following holds projectively
on $\hsmooth_{i}$, for any $\gamma\in LG\rtimes\rot$ and $f\in C^{\infty}(S^{1},V_{k})$.
\begin{equation}\label{eq:primary covariance}
\pi_{j}(\gamma)\phi(f)\pi_{i}(\gamma)^{*}=\phi(\gamma f)
\end{equation}
Moreover, if $d$ is the integrally--moded infinitesimal generator of rotations, then
\begin{equation}
e^{i\theta d}\phi(f)e^{-i\theta d}=\phi(f_{\theta})
\end{equation}
where $f_{\theta}(\phi)=f(\phi-\theta)$.
\end{proposition}
\proof
We use a standard ODE argument. Notice first that by continuity, $\phi$ intertwines the action
of $L\g$ on $\H_{i}^{\infty}$ and
$\H_{j}^{\infty}$, {\it i.e.}~
\begin{equation}
\pi_{i}(X)\phi(f)\xi-\phi(f)\pi_{j}(X)\xi=\phi(Xf)\xi
\end{equation}
for any $X\in L\g$, $f\in C^{\infty}(S^{1},V_{k})$, and $\xi\in\H_{i}^{\infty}$. Define now
$F(t)=e^{-t\pi_{j}(X)}\phi(e^{tX}f)e^{t\pi_{i}(X)}\xi$ where $X,f,\xi$ are as above. Then,
using the invariance of $\H_{i}^{\infty}$ under $LG$
(proposition \ref{ch:analytic}.\ref{invariance}) and the $C^{\infty}$--continuity of $\phi$,
\begin{equation}
\begin{split}
F(t+s)
&=
e^{-(t+s)\pi_{j}(X)}\phi(e^{tX}f+sXe^{tX}f+o(s))
\Bigl(e^{t\pi_{i}(X)}\xi+s\pi_{i}(X)e^{t\pi_{i}(X)}\xi+o(s)\Bigr)\\
&=
e^{-(t+s)\pi_{j}(X)}\Bigl(
\phi(e^{tX}f)e^{t\pi_{i}(X)}\xi+
s\phi(Xe^{tX}f)e^{t\pi_{i}(X)}\xi+s\phi(e^{tX}f)\pi_{i}(X)e^{t\pi_{i}(X)}\xi+o(s)\Bigr)\\
&=
F(t)+
s\Bigl(-\pi_{j}(X)\phi(e^{tX}f)e^{t\pi_{i}(X)}\xi+
 \phi(Xe^{tX}f)e^{t\pi_{i}(X)}\xi+s\phi(e^{tX}f)\pi_{i}(X)e^{t\pi_{i}(X)}\xi\Bigr)+o(s)\\
&=
F(t)+o(s)
\end{split}
\end{equation}
where $o(s)s^{-1}\rightarrow 0$ as $s\rightarrow 0$ in $\H^{\infty}$ and therefore in $\H$.
Thus, $F\in C^{1}(\IR,\H)$ and $\dot F\equiv 0$ so that $F\equiv F(0)=\xi$. It follows
that $\phi$ intertwines the one--parameter subgroups of $LG$ and therefore $LG$ itself.
The commutation properties with $\rot$ follow in a similar way \halmos\\


\remark\hfill
\begin{enumerate}
\item Recall that the central extensions of $L\Spin_{2n}$ corresponding to $\pi_{i},\pi_{j}$
are canonically isomorphic by proposition \ref{ch:analytic}.\ref{classify extension}
and denote either of them by $\L\Spin_{2n}$. It is not difficult to show that a more
astringent relations than \eqref{eq:primary covariance} holds, namely the (non--projective)
identity
\begin{equation}\label{eq:primary covariance 2}
\pi_{j}(\wt\gamma)\phi(f)\phi_{i}(\wt\gamma)^{*}=\phi(\wt\gamma f)
\end{equation}
for any $\wt\gamma\in\L\Spin_{2n}$.
\item Somewhat conversely, Wassermann has pointed out that it might be possible to
use the continuity of the level 1 spin primary fields and \eqref{eq:primary covariance}
to give an alternative proof of
proposition \ref{ch:analytic}.\ref{classify extension} for $G=\Spin_{2n}$.
Notice first that it is only necessary to prove this at level 1 since the level
$\ell$ representations are contained in the $\ell$--fold tensor product of the
level 1 representations. Moreover, we already know that the central extensions
corresponding to the level 1 vacuum and vector representations are isomorphic
since these arise as summands of the Neveu--Scharz Fock space
(proposition \ref{ch:fermionic}.\ref{pr:level 1 reps}). The same holds for the
level 1 spin representations since they are summands of the Ramond sector. Thus,
it is sufficient to prove the isomorphism of the central extensions corresponding
to the the vacuum representation $(\pi_{0},\H_{0})$ and one of the spin
representations $(\pi_{s},\H_{s})$ at level 1. Let
$\phi_{s}:\hfin_{0}\otimes V_{s}[z,z^{-1}]\rightarrow\hfin_{s}$ be the corresponding
primary field. By theorem \ref{th:spin sobolev} and proposition
\ref{pr:LG equivariance}, the following holds projectively for any
$f\in C^{\infty}(S^{1},V_{s})$ and $\gamma\in LG$
\begin{equation}\label{eq:correcting phase}
\pi_{s}(\gamma)\phi_{s}(f)\pi_{0}(\gamma)^{*}=\phi_{s}(\gamma f)
\end{equation}
This identity may be used to define a continuous map between the two central extensions
by associating to any lift $\wt\pi_{0}(\gamma)\in U(\H_{0})$ of $\pi_{0}(\gamma)$
the unique unitary lift $\wt\pi_{s}(\gamma)$ of $\pi_{s}(\gamma)$ such that
\eqref{eq:correcting phase} holds without any additional phase corrections. However,
the proof that it is a homomorphism requires, to the best of our knowledge, a proof
of the fact that the definition of the map is independent of the particular $f$ chosen.
\end{enumerate}


\part{Braiding of primary fields}


\newcommand {\prim}[1]{\phi_{#1}(v_{#1},z_{#1})}
\newcommand {\slc}{\mathfrak{sl}_{2}(\IC)}
\newcommand {\KZ}{Knizhnik--Zamolodchikov }
\newcommand {\derz}[1]{z_{#1}\frac{d}{dz_{#1}}}

\chapter{Braiding properties of primary fields}\label{ch:algebraic fields}

We outline the theory of the \KZ equations satisfied by the four--point
functions of primary fields and deduce from it their braiding properties.
We also compute the simplest instances of braiding. Finally, we show that,
when smeared on test functions, the primary fields satisfy the same
braiding relations provided the supports of the functions are disjoint.

\ssection{The Knizhnik--Zamolodchikov equations}
\label{se:the KZ equations}

Fix a level $\ell\in\IN$ and consider a collection of non--zero primary fields
$\phi_{n}\cdots\phi_{1}$ corresponding to the vertices
\begin{equation}
\vertex{k_{n}}{0}{i_{n}}
\vertex{k_{n-1}}{i_{n}}{i_{n-1}}
\cdots
\vertex{k_{2}}{i_{2}}{i_{1}}
\vertex{k_{1}}{i_{1}}{0}
\end{equation}
Notice that $V_{k_{1}}=V_{i_{1}}$ and $V_{k_{n}}=V_{i_{n}}^{*}$ and that, up to
a scalar the initial terms of $\phi_{n}$ and $\phi_{1}$ are the canonical
intertwiners $V_{i_{n}}\otimes V_{i_{n}}^{*}\rightarrow\IC$ and
$\IC\otimes V_{i_{1}}\rightarrow V_{i_{1}}$.
Fix $\Upsilon\in\H_{0}(0)\cong\IC$ of norm one and define the
{\it $n$--point function}
\begin{equation}
F=(\prim{n}\cdots\prim{1}\Upsilon,\Upsilon)
\end{equation}
a formal Laurent series in the variables $z_{n}\cdots z_{1}$. The $G$-equivariance
of the primary fields and the fact that $G$ fixes $\Upsilon$ imply that $F$ takes
values in $(V_{k_{n}}^{*}\otimes\cdots\otimes V_{k_{1}}^{*})^{G}$. We shall prove
that $F$ satisfies a first order partial differential equation. Attach for this
purpose an endomorphism $\Omega_{ij}$ to any pair $1\leq i,j\leq n$ by setting
\begin{equation}\label{eq:Omega operators}
\Omega_{ij}=\pi_{i}(X_{k})\pi_{j}(X^{k})
\end{equation}
where $\pi_{i}$ is the action of $\gc$ on the ith
factor of $(V_{k_{n}}^{*}\otimes\cdots\otimes V_{k_{1}}^{*})^{G}$, $\{X_{k}\},\{X^{k}\}$
are dual basis of $\gc$ with respect to the basic inner product on $\gc$ and the
summation over $k$ is implicit. Notice that $\Omega_{ij}$ is independent of the choice
of the basis $X_{k}$ and, in particular $\Omega_{ij}=\pi_{i}(X^{k})\pi_{j}(X_{k})=\Omega_{ji}$.
Moreover, $\Omega_{ii}$ acts as multiplication by $C_{k_{i}}$, the Casimir of $V_{k_{i}}$
and $V_{k_{i}}^{*}$.\\

Recall that on each $\hfin_{i_{k}}$, we have $L_{0}=d+\Delta_{i_{k}}$. Here,
as customary $d$ is the normalised generator of rotations, $L_{0}$ is given
by the Segal--Sugawara formula and $\Delta_{i_{k}}=\frac{C_{i_{k}}}{2\kappa}$
where $C_{i_{k}}$ is the Casimir of $V_{i_{k}}$ and $\kappa=\ell+\half{C_{g}}$.
Therefore, by (\ref{ch:classification}.\ref{intertwine}),
\begin{equation}\label{eq:begin}
\derz{i}\prim{i}=
[d,\prim{i}]-\Delta_{\phi_{i}}\prim{i}=[L_{0},\prim{i}]-\Delta_{k_{i}}\prim{i}
\end{equation}
Using the explicit expression of $L_{0}$,
\begin{equation}
L_{0}=
\frac{1}{\kappa}\Bigl(\half{1}X_{i}(0)X^{i}(0)+\sum_{m>0} X_{i}(-m)X^{i}(m)\Bigr)
\end{equation}
and the equivariance properties of primary fields, we find
\begin{equation}
\begin{split}
\kappa(\derz{i}+\Delta_{k_{i}})F
=&\half{1}
  (\prod_{j>i}\prim{j}[X_{k}(0)X^{k}(0),\prim{i}]\prod_{i>j}\prim{j}\Upsilon,\Upsilon)\\
+&\sum_{m>0}
  (\prod_{j>i}\prim{j}[X_{k}(-m)X^{k}(m),\prim{i}]\prod_{i>j}\prim{j}\Upsilon,\Upsilon)\\
=&\half{1}C_{k_{i}}F+\sum_{i>j} \Omega_{ij}F+
  \sum_{m>0}\biggl(
  \sum_{j>i}-\left(\frac{z_{i}}{z_{j}}\right)^{m}\Omega_{ij}F+
  \sum_{i>j} \left(\frac{z_{j}}{z_{i}}\right)^{m}\Omega_{ij}F\biggr)
\end{split}
\end{equation}
where the last equality is obtained by moving the $X_{k}(p)$, $p<0$ to the left and the
$X^{k}(q)$, $q\geq 0$ to the right and using $X(p)\Upsilon=0$ if $p\geq 0$. Thus, with
the understanding that the denominators below are to be expanded as if $|z_{n}|>\cdots>|z_{1}|$,
$F$ satisfies the {\it \KZ equations}
\begin{equation}
\frac{\partial F}{\partial z_{i}}=
\frac{1}{\kappa}
\sum_{j\neq i}\frac{\Omega_{ij}}{z_{i}-z_{j}}F
\end{equation}

Notice that $\Delta_{\phi_{1}}=0$ so that
$\prim{1}\Upsilon=\sum_{n\geq 0}z^{-n}\phi_{1}(v_{1},n)\Upsilon$
and, formally $\lim_{z_{1}\rightarrow 0}\prim{1}\Upsilon=v_{1}$. Similarly,
$\lim_{z_{n}\rightarrow\infty}(z_{n}^{\Delta_{\phi_{n}}}\prim{n}\cdots,\Upsilon)=
(\cdots,\overline{v_{n}})$ where we have used the antilinear identification
$V_{k_{n}}^{*}\rightarrow V_{k_{n}}$. Thus, the function
\begin{equation}
H(z_{2},\ldots z_{n-1})=
(\prim{n-1}\cdots\prim{2}v_{1},\overline{v_{n}})=
\lim_{\substack{z_{1}\rightarrow 0\\z_{n}\rightarrow\infty}}z_{n}^{\Delta_{\phi_{n}}}F
\end{equation}
satisfies, for $i=2\ldots n-1$ the equations
\begin{equation}
\frac{\partial H}{\partial z_{i}}=
\frac{1}{\kappa}\Bigl(\frac{\Omega_{i1}}{z_{i}}+
\sum_{j\neq 1,i,n}\frac{\Omega_{ij}}{z_{i}-z_{j}}\Bigr)H
\end{equation}

We specialise now to the case $n=4$. Set $w=z_{3}$ and $z=z_{2}$, then
$(w\partial_{w}+z\partial_{z})H=\kappa^{-1}(\Omega_{12}+\Omega_{23}+\Omega_{31})H$.
However, on $(V_{k_{1}}^{*}\otimes\cdots\otimes V_{k_{4}}^{*})^{G}$, we have
$\sum_{i=1}^{4}\pi_{i}(X)=0$ for any $X\in\g$ and therefore
\begin{equation}\label{trick}
\begin{split}
\Omega_{12}+\Omega_{23}+\Omega_{31}
&=\half{1}(\sum_{i=1}^{3}\pi_{i}(X_{k}))(\sum_{i=1}^{3}\pi_{i}(X^{k}))-
  \half{1}(\Omega_{11}+\Omega_{22}+\Omega_{33})
 =\half{1}(\Omega_{44}-\Omega_{11}-\Omega_{22}-\Omega_{33})\\
&=\kappa(\Delta_{k_{4}}-\Delta_{k_{1}}-\Delta_{k_{2}}-\Delta_{k_{3}})
 =-\kappa(\Delta_{\phi_{2}}+\Delta_{\phi_{3}})
\end{split}
\end{equation}
It follows that the {\it reduced four point function}
$\Psi(w,z)=w^{\Delta_{\phi_{3}}}z^{\Delta_{\phi_{2}}}H$ satisfies
$w\partial_{w}\Psi=-z\partial_{z}\Psi$ and therefore only contains monomials in $\zeta=zw^{-1}$.
Moreover, since $\zeta d_{\zeta}\Psi=z\partial_{z}\Psi$, $\Psi$ satisfies as a function of
one formal variable
\begin{equation} \label{kz}
\frac{d\Psi}{d\zeta}=\frac{1}{\kappa}\Bigl(
\frac{\Omega_{12}+\kappa\Delta_{\phi_{2}}}{\zeta}+\frac{\Omega_{23}}{\zeta-1}\Bigr)\Psi
\end{equation}
where $(\zeta-1)^{-1}$ is to be expanded as if $|\zeta|<1$.\\

Aside from \eqref{kz}, $\Psi$ satisfies a number of algebraic identities. To see
this, let $v_{1}\in V_{k_{1}}$ be a highest weight vector for
$\{e_{\theta},f_{\theta},h_{\theta}\}$ with $h_{\theta}v_{1}=s_{1}v_{1}$.
Since $v_{1}$ is of lowest energy,
$f_{\theta}(1)v_{1}=0$ and therefore $v_{1}$ is of highest weight for
$\{f_{\theta}(1),e_{\theta}(-1),\ell-h_{\theta}\}$ with eigenvalue
$\ell-s_{1}$. Thus, by elementary $\slc$ theory,
$e_{\theta}(-1)^{\ell-s_{1}+1}v_{1}=0$.
On the other hand, a simple induction shows that, for any $n\in\IN$
\begin{equation}
(\phi_{3}(v_{3},w)\phi_{2}(v_{2},z)e_{\theta}(-1)^{n}v_{1},v_{4})=
(-1)^{n}(w^{-1}\pi_{3}(e_{\theta})+z^{-1}\pi_{2}(e_{\theta}))^{n}
(\phi_{3}(v_{3},w)\phi_{2}(v_{2},z)v_{1},v_{4})
\end{equation}
so that the right hand--side vanishes for $n=\ell-s_{1}+1$. An additional set of identities
for $\Psi$ may be similarly derived from the equation $e_{\theta}(-1)^{\ell-s_{4}+1}v_{4}=0$
which holds whenever $v_{4}\in V_{4}$ is of highest weight for
$\{e_{\theta},f_{\theta},h_{\theta}\}$ with $h_{\theta}v_{4}=s_{4}v_{4}$. This completes the
proof of the first part of the following

\begin{proposition}[Tsuchiya-Kanie] \label{kz completeness}
Let $\phi_{3},\phi_{2}$ be primary fields with vertices $\vertex{k_{3}}{k_{4}}{j}$,
$\vertex{k_{2}}{j}{k_{1}}$ and conformal weights $\Delta_{3}$, $\Delta_{2}$ respectively.
Define the reduced four--point function
\begin{equation}
\Psi_{j}=
w^{\Delta_{3}}z^{\Delta_{2}}(\phi_{3}(v_{3},w)\phi_{2}(v_{2},z)v_{1},v_{4})=
\sum_{n\geq 0}(\psi_{3}(v_{3},n)\phi_{2}(v_{2},-n)v_{1},v_{4})\left(\frac{z}{w}\right)^{n}
\end{equation}
a formal power series in the variable $\displaystyle{\zeta=\left(\frac{z}{w}\right)}$ with
coefficents in $(V_{k_{4}}\otimes V_{k_{3}}^{*}\otimes V_{k_{2}}^{*}\otimes V_{k_{1}}^{*})^{G}$.
Then $\Psi_{j}$ converges to a holomorphic function on $|\zeta|<1$ such that 
$\Phi_{j}=\zeta^{-\Delta_{2}}\Psi_{j}$ satisfies
\begin{equation}\label{normalised kz}
\frac{d\Phi_{j}}{d\zeta}=
\frac{1}{\kappa}\Bigl(\frac{\Omega_{12}}{\zeta}+\frac{\Omega_{23}}{\zeta-1}\Bigr)\Phi_{j}
\end{equation}
and the algebraic equations
\begin{align}
(\pi_{2}(e_{\theta})+\zeta\pi_{3}(e_{\theta}))^{\ell-s_{1}+1}
\Phi_{j}(v_{1},v_{2},v_{3},v_{4})&=0 \label{boundary one}\\
\intertext{whenever $v_{1}$ is a highest weight vector for $\{e_{\theta},f_{\theta},h_{\theta}\}$
with $h_{\theta}v_{1}=s_{1}v_{1}$, and}
(\zeta\pi_{2}(f_{\theta})+\pi_{3}(f_{\theta}))^{\ell-s_{4}+1}
\Phi_{j}(v_{1},v_{2},v_{3},v_{4})&=0 \label{boundary two}
\end{align}
whenever $v_{4}$ is a highest weight vector for $\{e_{\theta},f_{\theta},h_{\theta}\}$
with $h_{\theta}v_{1}=s_{4}v_{1}$.
The set of $\Phi_{j}$ corresponding to all intermediate $V_{j}$ is a basis of the space of
solutions to the above system of equations near 0.
\end{proposition}

\remark The convergence of $\Psi_{j}$ follows at once because $\Psi_{j}$ is a formal power
series solution of \eqref{normalised kz}, which has a regular singularity at $\zeta=0$.
The completeness of the set of reduced four--point functions is stated as Proposition 4.3
of \cite{TK1}. Its proof contains a gap bridged by the treatment given in \cite{Wa2}.

\ssection{Braiding of primary fields}\label{se:alg braiding}

\begin{theorem}[Braiding]\label{th:existence of braiding}
Let $\phi_{k_{4}j}^{k_{3}},\phi_{jk_{1}}^{k_{2}}$ be primary fields corresponding to the
vertices $\vertex{k_{3}}{k_{4}}{j}$ and $\vertex{k_{2}}{j}{k_{1}}$. Then the product
$\phi_{k_{4}j}^{k_{3}}(v_{3},w)\phi_{jk_{1}}^{k_{2}}(v_{2},z)$ defines a single--valued
holomorphic function on $\{(w,z)|\frac{w}{z}\notin[0,\infty)\}$ when it is evaluated
between finite energy vectors. Moreover,
\begin{equation}\label{braiding}
\phi_{k_{4}j}^{k_{3}}(v_{3},w)\phi_{jk_{1}}^{k_{2}}(v_{2},z)=
\sum_{h}b_{h}\phi_{k_{4}h}^{k_{2}}(v_{2},z)\phi_{hk_{1}}^{k_{3}}(v_{3},w)
\end{equation}
where $b_{h}\in\IC$ and the sum on the right hand--side spans the products of all primary
fields with vertices $\vertex{k_{2}}{k_{4}}{h},\vertex{k_{3}}{h}{k_{1}}$ corresponding
to any intermediate $V_{h}$.
\end{theorem}
\proof
By the commutation properties of the primary fields with $\lpol\g$, it is sufficient to prove
\eqref{braiding} when it is evaluated between lowest energy vectors.
Let $\Psi_{j}$ be the reduced four point-function in the variable $\zeta=\frac{z}{w}$
corresponding to $\phi_{k_{4}j}^{k_{3}}\phi_{jk_{1}}^{k_{2}}$. By the previous proposition,
$\Phi_{j}=\zeta^{-\Delta_{2}}\Psi_{j}$, $\Delta_{2}=\Delta_{k_{1}}+\Delta_{k_{2}}-\Delta_{j}$
converges to a homolorphic function on $|\zeta|<1$ which may analytically be continued to a
single--valued function on $\IC\backslash[0,\infty]$ satisfying
\eqref{normalised kz}--\eqref{boundary two}.
Near $\infty$, we may rewrite \eqref{normalised kz} in terms of the local coordinate
$\eta=\zeta^{-1}$. Since
$\zeta\partial_{\zeta}=-\eta\partial_{\eta}$, $\Phi_{j}$ satisfies
\begin{equation}
\frac{d\Phi_{j}}{d\eta}=
\frac{1}{\kappa}\Bigl(
-\frac{\Omega_{12}}{\eta}+
 \frac{\Omega_{23}}{\eta(\eta-1)}\Bigr)\Phi_{j}=
 \frac{1}{\kappa}\Bigl(\frac{\Omega_{13}+\kappa(\Delta_{2}+\Delta_{3})}{\eta}
+\frac{\Omega_{23}}{\eta-1}\Bigr)\Phi_{j}
\end{equation}
where we used \eqref{trick}. On the other hand, the boundary conditions may be rewritten in
terms of $\eta$ as
\begin{align}
(\eta\pi_{2}(e_{\theta})+\pi_{3}(e_{\theta}))^{\ell-s_{1}+1}
\Phi_{j}(v_{1},v_{2},v_{3},v_{4})&=0\\
(\pi_{2}(f_{\theta})+\eta\pi_{3}(f_{\theta}))^{\ell-s_{4}+1}
\Phi_{j}(v_{1},v_{2},v_{3},v_{4})&=0
\end{align}
It follows that
$\Phi_{j}\eta^{-(\Delta_{2}+\Delta_{3})}=\Psi_{j}\zeta^{3}$
is a solution of \eqref{normalised kz}--\eqref{boundary two}
with respect $\eta$ when
$(V_{k_{4}}\otimes V_{k_{3}}^{*}\otimes V_{k_{2}}^{*}\otimes V_{k_{1}}^{*})^{G}$ and
$(V_{k_{4}}\otimes V_{k_{2}}^{*}\otimes V_{k_{3}}^{*}\otimes V_{k_{1}}^{*})^{G}$ are
identified in the obvious way. By proposition \ref{kz completeness}, it may be therefore be
written as a sum of reduced four point functions
corresponding to the vertices $\vertex{k_{2}}{k_{4}}{h}$, $\vertex{k_{3}}{h}{k_{1}}$ multiplied
by an appropriate power of $\eta$. Writing $\zeta=zw^{-1}$ and $\eta=wz^{-1}$, the powers of
$z$ and $w$ are easily seen to cancel and \eqref{braiding} is proved \halmos\\

\remark
\begin{enumerate}
\item
The braiding coefficients $b_{h}$ in \eqref{braiding} are only defined up to multiplication
by non--zero constants accounting for the multiplicative ambiguity in the choice of the
initial term of the various primary fields. A further lack of determination arises if the
$\Hom$ spaces corresponding to some of the vertices are of dimension greater or equal to 2.
However, once a definite choice of initial terms is made, the $b_{h}$ are unambiguously
defined and may be computed by expressing the analytic continuation of the reduced four--point
function corresponding to the left hand--side of \eqref{braiding} from $\zeta=0$ to
$\infty$ in terms of the reduced four point functions corresponding to the right hand--side
evaluated at $\zeta^{-1}$.
\item
The existence of braiding is perhaps best demonstrated using the \KZ equations
in their PDE form. These define a flat connection on the configuration space of
$n$ ordered points in $\IC$ for the (topologically) trivial bundle with fibre
$(V_{k_{n}}^{*}\otimes\cdots\otimes V_{k_{1}}^{*})^{G}$. The permutation of the
fields $\phi_{i}$ and $\phi_{i+1}$ is then obtained by analytic continuation
along a path connecting $(z_{n}\ldots z_{i+1},z_{i}\ldots z_{1})$ to
$(z_{n}\ldots z_{i},z_{i+1}\ldots z_{1})$.
The above one--variable approach is possible because the connection is invariant
under the diagonal action of $SL(2,\IC)$ and gives an effective method for
computing the braiding coefficients.
\end{enumerate}

\ssection{Abelian braiding}\label{se:abelian}

The simplest instance of braiding is the {\it abelian braiding} obtained from
\KZ equations with values in a one--dimensional space. We compute it below

\begin{lemma}\label{le:abelian}
Let $V_{i}$, $i=2\ldots 4$ be irreducible $G$--modules such that
$V_{4}\subset V_{2}\otimes V_{3}$ with multiplicity one and assume that the
generator of $\Hom_{G}(V_{2}\otimes V_{3},V_{4})$ is the initial term of a
primary field $\field{V_{3}}{V_{4}}{V_{2}}$. Then,
\begin{equation}\label{eq:abelian}
\field{V_{3}}{V_{4}}{V_{2}}(w)\field{V_{2}}{V_{2}}{V_{0}}(z)=
\lambda
\field{V_{2}}{V_{4}}{V_{3}}(z)\field{V_{3}}{V_{3}}{V_{0}}(w)
\end{equation}
where $\lambda\neq 0$. In fact, let the initial terms of $\field{V_{i}}{V_{i}}{V_{0}}$,
$i=2,3$ be the canonical intertwiners $\IC\otimes V_{i}\rightarrow V_{i}$
and denote those of $\field{V_{3}}{V_{4}}{V_{2}}$ and $\field{V_{2}}{V_{4}}{V_{3}}$ by
$\varphi_{3}$ and $\varphi_{2}$ respectively. Then
\begin{enumerate}
\item If $\varphi_{3},\varphi_{2}$ are normalised by setting
$\varphi_{2}=\varphi_{3}\sigma$ where
$\sigma:V_{3}\otimes V_{2}\rightarrow V_{2}\otimes V_{3}$ is permutation, then
$\lambda=e^{-i\pi\beta}$ where $\beta=\frac{1}{2\kappa}(C_{4}-C_{2}-C_{3})$ and
the $C_{i}$ are the Casimirs of the corresponding representations.
\item If $V_{3}=V_{2}$, and $\varphi_{3},\varphi_{2}$ are normalised by equating
them so that $\field{3}{4}{2}=\field{2}{4}{3}$, then
$\lambda=\varepsilon e^{-i\pi\beta}$ where $\varepsilon=\pm 1$ according to whether
$V_{4}\subset V_{2}\otimes V_{2}$ is symmetric or anti-symmetric under the action
of $\mathfrak S_{2}$.
\end{enumerate}
\end{lemma}
\proof
By assumption,
$W=(V_{4}\otimes V_{3}^{*}\otimes V_{2}^{*}\otimes\IC)^{G}
  \cong\Hom_{G}(V_{2}\otimes V_{3},V_{4})$
is one--dimensional and the \KZ equations with values in $W$ have a unique non
trivial solution $f$ given by the reduced four--point function of the left--hand
side of \eqref{eq:abelian}. The vanishing of $\lambda$ would lead to that of $f$,
a contradiction.\\
(i) $\Omega_{23}$ acts on $W$ as multiplication by $\kappa\beta$. Indeed,
\begin{equation}
\begin{split}
\Omega_{23}P
&=P\pi_{2}(X_{k})\pi_{3}(X^{k})\\
&=\half{1}P(\pi_{2}(X_{k})+\pi_{3}(X_{k}))(\pi_{2}(X^{k})+\pi_{3}(X^{k}))-
  \half{1}P(\pi_{2}(X_{k})\pi_{2}(X^{k})+\pi_{3}(X_{k})\pi_{3}(X^{k}))\\
&=\half{1}\pi_{4}(X_{k})\pi_{4}(X^{k})P-
  \half{1}P(C_{2}+C_{3})\\
&=\half{1}(C_{4}-C_{2}-C_{3})P
\end{split}
\end{equation}
Thus, the \KZ equations read $f'=\beta(\zeta-1)^{-1}f$ and, up to a constant,
$f=(\zeta-1)^{\beta}$ where, in accordance with our conventions on the definition of
braiding, the function is defined with a cut in the $\zeta-1$ plane along $[0,\infty]$
so that $f$ is continuous on $\IC\backslash[0,\infty]$.
Since the leading term of $f$ at $0$ is $\varphi_{3}$, we get
$f=e^{-\pi i\beta}(z-1)^{\beta}\varphi_{3}$.
Similarly, the leading coefficient at $\infty$ of the reduced four point-function $g$
corresponding to the right hand--side of \eqref{eq:abelian} is $\varphi_{2}$ so that
$g(\eta)=e^{-\pi i\beta}(\eta-1)^{\beta}\varphi_{2}$ where $\eta=\zeta^{-1}$.
If $\zeta\in (-\infty,0)$ so that $\arg\zeta=\pi$, then
$(\zeta-1)^{\beta}=e^{-i\pi\beta}\zeta^{\beta}(1/\zeta-1)^{\beta}$ and therefore
the analytic continuation of $f\sigma$ at $\infty$ is
$e^{-2\pi i\beta}\eta^{-\beta}(\eta-1)^{\beta}\varphi_{3}\sigma$ whence
$\lambda=e^{-2\pi i\beta}$ as claimed.\\
(ii) follows at once from (i) since $\varphi_{3}\sigma=\varepsilon\varphi_{3}$
\halmos

\ssection{Braiding of smeared primary fields}\label{se:smeared braiding}

\begin{proposition}\label{pr:smeared braiding}
Let $\phi_{kj},\phi_{ji}$, $\phi_{kh},\phi_{hi}$ be primary fields of charges $U,V,V,U$
respectively satisfying the braiding relations
\begin{equation} \label{basic braiding}
\phi_{kj}(u,w)\phi_{ji}(v,z)=
\sum_{h} b_{h}\phi_{kh}(v,z)\phi_{hi}(u,w)
\end{equation}
If all fields extend to operator--valued distributions, then for any $f\in C^{\infty}(S^{1},U)$
and $g\in C^{\infty}(S^{1},V)$ supported in $I=(0,\pi)$ and $I^{c}=(\pi,2\pi)$ respectively, the
following holds
\begin{equation} \label{smeared braiding}
\phi_{kj}(f)\phi_{ji}(g)=
\sum_{h} b_{h}\phi_{kh}(ge_{\alpha_{jh}})\phi_{hi}(fe_{-\alpha_{jh}})
\end{equation}
where $\alpha_{jh}=\Delta_{i}+\Delta_{k}-\Delta_{j}-\Delta_{h}$ and
$e_{\mu}(\theta)=e^{i\mu\theta}$ for $\theta\in(0,2\pi)$.
\end{proposition}
\proof
By continuity, we need only prove \eqref{smeared braiding} when both sides are
evaluated between finite energy vectors and by the commutation relations of the
primary fields with $\lpol\g$ it is sufficient to consider lowest energy vectors.
Moreover, we may assume that $f$ and $g$ are of the form $a\otimes u$, $b\otimes v$ where
$a,b\in C^{\infty}(S^{1},\IC)$ are supported in $I,I^{c}$ respectively and
$(u,v)\in U\times V$. Let $\xi,\eta$ be lowest energy vectors, then
\begin{equation} \label{left side}
(\phi_{kj}(f)\phi_{ji}(g)\xi,\eta)=
\sum_{n,m}     a_{n}b_{m} (\phi_{kj}(u,n)\phi_{ji}(v,m)\xi,\eta)=
\sum_{n\geq 0} a_{n}b_{-n}(\phi_{kj}(u,n)\phi_{ji}(v,-n)\xi,\eta)
\end{equation}
an absolutely convergent sum since the Fourier coefficents $a_{n},b_{m}$ of $a,b$ decrease
rapidly and, by the continuity of $\phi_{kj},\phi_{ji}$,
$\|\phi_{ji}(u,n)\xi\|\leq C(1+|n|)^{\beta_{ji}}$,
$\|\phi_{kj}(v,m)\eta\|\leq C(1+|m|)^{\beta_{kj}}$ for some
$\beta_{ji},\beta_{kj}\geq 0$.
Let
\begin{equation}
F_{j}(\zeta)=\sum_{n\geq 0}(\phi_{kj}(u,n)\phi_{ji}(v,-n)\xi,\eta)\zeta^{n}
\end{equation}
be the reduced four--point function corresponding to the product $\phi_{kj}\phi_{ji}$.
$F_{j}$ is convergent on $|\zeta|<1$ and extends to a single--valued holomorphic function
on $\IC\backslash[0,\infty)$.
The function $\check a*b(\theta)=\int_{0}^{2\pi}a(-\phi)b(\theta-\phi)\frac{d\phi}{2\pi}$
has Fourier coefficients $a_{-n}b_{n}$ and is supported away from 0 so that
$\check a*b(\zeta|\zeta|^{-1})F_{j}(\zeta)$ is smooth on $\IC\backslash\{0\}$.
Thus, \eqref{left side} is equal to
\begin{equation}\label{eq:above}
\lim_{r\nearrow 1}\int_{0}^{2\pi}\check a*b(\theta)F_{j}(re^{i\theta})\frac{d\theta}{2\pi}=
\int_{0}^{2\pi} \check a*b(\theta)F_{j}(e^{i\theta})\frac{d\theta}{2\pi}=
\lim_{r\searrow 1}\int_{0}^{2\pi}\check a*b(\theta)F_{j}(re^{i\theta})\frac{d\theta}{2\pi}
\end{equation}
Let now $G_{h}$ be the reduced four--point functions corresponding to the products
$\phi_{kh}\phi_{hi}$. From the braiding relations \eqref{basic braiding}, we have
$F_{j}(\zeta)=\sum_{h}b_{h}\zeta^{\alpha_{jh}}G_{h}(\zeta^{-1})$ where
$\alpha_{jh}=\Delta_{i}+\Delta_{k}-\Delta_{j}-\Delta_{h}$
and \eqref{eq:above} is therefore equal to
\begin{equation}
\begin{split}
 &
\sum_{h}b_{h}
\lim_{r\searrow 1}r^{\alpha_{jh}}\int_{0}^{2\pi}
\check a*b\thickspace e_{\alpha_{jh}}(\theta)
G_{h}(r^{-1}e^{-i\theta})\frac{d\theta}{2\pi}\\
=&
\sum_{h}b_{h}
\lim_{r\searrow 1}r^{\alpha_{jh}}\int_{0}^{2\pi}
e^{2\pi i\alpha_{jh}\theta}
\check b*a\thickspace e_{-\alpha_{jh}}(\theta)G_{h}(r^{-1}e^{i\theta})\frac{d\theta}{2\pi}
\end{split}
\end{equation}
Since $(\check b*a)e_{-\mu}=e^{-2\pi i\mu}\check{(be_{\mu})}*(ae_{-\mu})$, the above
yields, by the continuity of $\phi_{kh}$ and $\phi_{hi}$,
\begin{equation}
\sum_{h}b_{h}\sum_{n\geq 0}
(be_{\alpha_{jh}})_{n}(ae_{-\alpha_{jh}})_{-n}
(\phi_{kh}(v,n)\phi_{hi}(u,-n)\xi,\eta)=
\sum_{h}b_{h}(\phi_{kh}(ge_{\alpha_{jh}})\phi_{hi}(fe_{-\alpha_{jh}})\xi,\eta)
\end{equation}
as claimed \halmos




\newcommand {\skv}{S^{k}V}

\newcommand {\DF}{Dotsenko--Fateev }
\newcommand {\IP}{\mathbb P}
\newcommand{\esse}{\operatorname{s}}
\newcommand{\ci}{\operatorname{c}}
\newcommand {\w}{\omega}


\chapter{Braiding the vector primary field with its symmetric powers}
\label{ch:box/kbox braiding}

In this chapter, we compute the structure constants governing the braiding of the
$L\Spin_{2n}$--primary fields whose charges are the vector representation and one
of its symmetric (traceless) powers.
As explained in chapter \ref{ch:algebraic fields}, the braiding coefficients may
be obtained via the analytic continuation from $0$ to $\infty$ of the solutions
of the underlying \KZ equation. In section \ref{se:DF equation}, we study this
problem for a related third order Fuchsian ODE discovered by Dotsenko and Fateev.
Its solutions are expressed as generalised Euler integrals and their analytic
continuation computed by contour deformation. 
In section \ref{se:box/kbox equation}, we compute explicitly the matrices
$\Omega_{ij}$ of the \KZ equation. Finally, in section \ref{se:braiding relations}
we show that the \KZ equation reduces to the \DF equation and deduce the analytic
continuation results for the former from those of the latter.

\ssection{The Dotsenko--Fateev equation}\label{se:DF equation}

The \DF equation is the following third order Fuchsian ordinary differential
equation with regular singular points at $0,1,\infty$ \cite[eqn. (A.9)]{DF}
\begin{equation}\label{DF}
f'''+
\frac{K_{1}z+K_{2}(z-1)}{z(z-1)}f''+
\frac{L_{1}z^{2}+L_{2}(z-1)^{2}+L_{3}z(z-1)}{z^{2}(z-1)^{2}}f'+
\frac{M_{1}z+M_{2}(z-1)}{z^{2}(z-1)^{2}}f=0
\end{equation}
where
\begin{xalignat}{2}
K_{1}&=-(g+3b+3c)&
K_{2}&=-(g+3a+3c)\\
L_{1}&=(b+c)(2b+2c+g+1)&
L_{2}&=(a+c)(2a+2c+g+1)
\end{xalignat}
\begin{equation}
\begin{split}
L_{3}&=(b+c)(2a+2c+g+1)+(a+c)(2b+2c+g+1)\\
     &+(c-1)(a+b+c)+(3c+g)(a+b+c+g+1)
\end{split}
\end{equation}
\begin{xalignat}{2}
M_{1}&=-c(2b+2c+g+1)(2a+2b+2c+g+2)&
M_{2}&=-c(2a+2c+g+1)(2a+2b+2c+g+2)
\end{xalignat}
and $a,b,c,g\in\IC$ are free parameters. In \S \ref{ss:derivation of DF}, we prove
that if $a,b,c,g$ lie in a suitable range, the generalised Euler integrals
\begin{equation}\label{euler}
\int\limits_{C_{1}}\int\limits_{C_{2}}
t_{1}^{a}(t_{1}-1)^{b}(t_{1}-z)^{c}
t_{2}^{a}(t_{2}-1)^{b}(t_{2}-z)^{c}(t_{1}-t_{2})^{g}dt_{1}dt_{2}
\end{equation}
where the $C_{i}$ are contours joining one of the points $0,z,1,\infty$ to another,
yield solutions of \eqref{DF}.
Following \cite{DF}, we show in \S \ref{ss:monodromy} how different choices of the
$C_{i}$ lead to a basis of solutions diagonalising the monodromy at a given singular
point and in \S \ref{ss:continuation} compute the analytic continuation of solutions
from $0$ to $\infty$ by deforming the contours. We begin by establishing the
convergence of \eqref{euler}.

\ssubsection{Convergence of contour integrals}\label{ss:convergence}

We shall need a special case of a more general result of Selberg \cite{Sel}

\begin{proposition}[Selberg]\label{selberg}
The improper integral
\begin{equation}
J_{2}(\alpha,\beta;\gamma)=
\int_{0}^{1}\int_{0}^{1}
t_{1}^{\alpha}(1-t_{1})^{\beta}
t_{2}^{\alpha}(1-t_{2})^{\beta}|t_{1}-t_{2}|^{\gamma}dt_{1}dt_{2}
\end{equation}
is absolutely convergent if
\begin{xalignat}{3}\label{range 0}
\Re\alpha,\Re\beta,\Re\gamma &> -1&
\Re(2\alpha+\gamma) &> -2&
\Re(2\beta+\gamma)  &> -2
\end{xalignat}
and is equal to
\begin{equation}
\prod_{j=1}^{2}\frac
{\Gamma(j\half{\gamma}+1)\Gamma(\alpha+(j-1)\half{\gamma}+1)
 \Gamma(\beta+(j-1)\half{\gamma}+1)}
{\Gamma(\half{\gamma}+1)\Gamma(\alpha+\beta+j\half{\gamma}+2)}
\end{equation}
\end{proposition}

\remark If $\alpha,\beta,\gamma$ satisfy \eqref{range 0} then
$J_{2}(\alpha,\beta;\gamma)\neq 0$ since $\Gamma(\zeta)^{-1}$ vanishes iff
$\zeta\in\{0,-1,-2,\ldots\}$ and
$\Re(\half{\gamma}+1)>\half{1}$,
$\Re(\alpha+\beta+\half{\gamma}+2)=\Re(\alpha+\half{1}(2\beta+\gamma)+2)>0$,
$\Re(\alpha+\beta+\gamma+2)=\Re(\half{1}(2\alpha+\gamma)+\half{1}(2\beta+\gamma)+2)>0$.

\begin{corollary}\label{converge}
Let $a,b,c,g\in\IC$ satisfy
\begin{xalignat}{2}
\Re a,\Re b,\Re c,\Re g &> -1&
\Re (2a+g),\Re (2b+g),\Re (2c+g) &> -2
\label{range 1}\\
\Re (a+b+c+g) &< -1&
\Re (2a+2b+2c+g) &< -2
\label{range 2}
\end{xalignat}
Then, for any pair of contours $C_{i}:I\rightarrow\IP^{1}$ with interior
$C_{i}(\overset{\circ}{I})\subset\IP^{1}\setminus\{0,1,z,\infty\}$ and end--points
$C_{i}(\partial I)\in\{0,1,z,\infty\}$, the improper integral
\begin{equation}
\int\limits_{C_{1}}\int\limits_{C_{2}}
t_{1}^{a}(t_{1}-1)^{b}(t_{1}-z)^{c}
t_{2}^{a}(t_{2}-1)^{b}(t_{2}-z)^{c}(t_{1}-t_{2})^{g}dt_{1}dt_{2}
\end{equation}
is absolutely convergent, locally uniformly in $z\in\IC\setminus\{0,1\}$.
\end{corollary}
\proof This follows by a tedious case--by--case analysis. When both $C_{1}$ and
$C_{2}$ have their end--points at $0$ and $1$, the result is an immediate corollary
of proposition \ref{selberg} since $|(t_{1}-z)^{g}(t_{2}-z)^{g}|$ is bounded on
$C_{1}\times C_{2}$ locally uniformly in $z$. The other cases follow in a similar
fashion after a change of variable bringing the end--points of $C_{1}$ and $C_{2}$
onto $0,1$ \halmos

\ssubsection{A derivation of the Dotsenko--Fateev equation}\label{ss:derivation of DF}

We prove in proposition \ref{euler solves} that the integrals \eqref{euler}
are solutions of the \DF equation. The proof proceeds as follows. We consider
the multi--valued function
\begin{equation}
\Phi=t^{a}(t-1)^{b}(t-z)^{c}s^{a}(s-1)^{b}(s-z)^{c}(t-s)^{g}
\end{equation}
and associated holomorphic differential form $\w=\Phi dt\wedge ds$ on
$C_{z}=\{(t,s)\in\IC^{2}|\thinspace t,s\notin\{0,1,z\}, t\neq s\}$ and
prove that the cohomology class of $\w$ satisfies \eqref{DF}. Integration
by parts then shows that $\int_{C_{1}\times C_{2}}\w$ satisfies \eqref{DF}.\\

A word of explanation is owed to the geometrically minded reader. A more elegant
approach to the above, indeed one that yields solutions of \eqref{DF} for all
values of $a,b,c,g$, is to consider over each fibre of
$C_{z}\rightarrow C\rightarrow \IC\setminus\{0,1\}\ni z$
the flat line bundle $\L$ determined by $\Phi$.
$\w$ may then be interpreted as a section of the cohomology bundle whose fibre
at $z$ is $H^{2}(C_{z},\L)$ and correspondingly it should be integrated over
homology cycles in $H_{2}(C_{z},\L^{\vee})$ of which the $C_{1}\times C_{2}$ 
we have chosen are not elements. 
Explicit basis elements of $H_{2}(C_{z},\L^{\vee})$ are described in \cite{Ko}
but leave one with the problem of showing that these may be deformed to the more
singular contours we have been using. This is necessary for the computations of
\S \ref{ss:continuation} below critically depend upon the use of these simpler
contours. On the other hand, the deformation result requires
a substantial amount of book--keeping, see for example \cite[\S 5]{TK2} where a
somewhat simpler case is treated. We have therefore opted for a more bare--handed and
direct approach.\\

The proof of proposition \ref{euler solves} relies on the computation of the
exterior derivative $d$ of a number of differential forms of the kind $\Phi f$
where $f$ is meromorphic. These are more easily spelled out by omitting $\Phi$
and using the twisted differential $\wt d=d+d\log\Phi$ given by
\begin{equation}\label{twisted d}
\wt d=d+
(\frac{a}{t}+\frac{b}{t-1}+\frac{c}{t-z}+\frac{g}{t-s})dt+
(\frac{a}{s}+\frac{b}{s-1}+\frac{c}{s-z}+\frac{g}{s-t})ds
\end{equation}
since $d(\Phi f)=(\Phi d+\Phi d\log\Phi)f=\Phi\wt d f$. We begin with a digression.\\

The prototype of the derivation of the \DF equation is that of the hypergeometric
one satisfied by the integral
\begin{equation}
\int_{1}^{\infty} t^{a}(t-1)^{b}(t-z)^{c}dt
\end{equation}
which we briefly sketch. The basic multi-valued function is
$\Phi=t^{a}(t-1)^{b}(t-z)^{c}$ with corresponding twisted differential
$\wt d=d+(\frac{a}{t}+\frac{b}{t-1}+\frac{c}{t-z})dt$. Given our convention on omitting
$\Phi$, differentiation by $z$ is given by $\nabla=\frac{d}{dz}-\frac{c}{t-z}$ so that
the basic class $\w=\Phi dt$ and its derivatives $\w^{(k)}=\nabla^{k}\w$ are
\begin{xalignat}{3}
\w&=dt&
\w'&=-c\frac{dt}{t-z}&
\w''&=c(c-1)\frac{dt}{(t-z)^{2}}
\end{xalignat}

The hypergeometric equation may be obtained almost at one stroke by writing
\begin{equation}\label{bas}
\begin{split}
0=\wt d\Bigl(\frac{1}{t-z}\Bigr)
 &=-\frac{dt}{(t-z)^{2}}
  +(\frac{a}{t}+\frac{b}{t-1}+\frac{c}{t-z})\frac{dt}{t-z}\\
 &=c^{-1}\w''-(\frac{a}{z}+\frac{b}{z-1})c^{-1}\w'-
  \frac{a}{z}\frac{dt}{t}-\frac{b}{z-1}\frac{dt}{t-1}
\end{split}
\end{equation}

The last two terms can be eliminated from \eqref{bas} by using
\begin{equation}\label{sub 1}
0=\wt d(t-1)
 =dt+(a\frac{t-1}{t}+b+c\frac{t-1}{t-z})dt
 =(1+a+b+c)\w-(z-1)\w'-a\frac{dt}{t}
\end{equation}
and
\begin{equation}\label{sub 2}
0=\wt dt
=dt+(a+b\frac{t}{t-1}+c\frac{t}{t-z})dt
 =(1+a+b+c)\w-z\w'+b\frac{dt}{t-1}
\end{equation}

Substituting \eqref{sub 1}--\eqref{sub 2} into \eqref{bas} yields
\begin{equation}
\w''-\frac{(b+c)z+(a+c)(z-1)}{z(z-1)}\w'+\frac{c(1+a+b+c)}{z(z-1)}\w=0
\end{equation}
so that the cohomology class of $\w$ satisfies the hypergeometric equation
\begin{equation}
\w''+\frac{(1+\alpha+\beta)z-\gamma}{z(z-1)}\w'+\alpha\beta\w=0
\end{equation}
with $\alpha=-c$, $\beta=-(1+a+b+c)$ and $\gamma=-(a+c)$ \cite[\S 7.23]{In}.\\

We turn now to the differential equation satisfied by the contour integrals
\begin{equation}
\int_{C_{1}}\int_{C_{2}}t^{a}(t-1)^{b}(t-z)^{c}s^{a}(s-1)^{b}(s-z)^{c}(t-s)^{g}dtds
\end{equation}

The twisted differential is given by \eqref{twisted d} and differentiation by $z$
acts as $\nabla=\frac{d}{dz}-c(\frac{1}{t-z}+\frac{1}{s-z})$. Denoting the
anti-symmetrisation of a differential form $\phi$ with respect to the permutation
$\sigma(t,s)=(s,t)$ by $[\phi]=\phi-\sigma^{*}\phi$, the basic class and its
derivatives are
\begin{xalignat}{3}
\w   &= dt\wedge ds &
\w'  &= -c[\frac{dt\wedge ds}{t-z}] &
\w'' &= c(c-1)[\frac{dt\wedge ds}{(t-z)^{2}}]+c^{2}[\frac{dt\wedge ds}{(t-z)(s-z)}]
\end{xalignat}
and
\begin{equation}
\w''' =-c(c-1)(c-2)[\frac{dt\wedge ds}{(t-z)^{3}}]
       -3c^{2}(c-1)[\frac{dt\wedge ds}{(t-z)^{2}(s-z)}]
\end{equation}

\begin{proposition}\label{euler solves}
Let $a,b,c,g\in\IC$ lie in the range \eqref{range 1}--\eqref{range 2} so that the
integrals
\begin{equation}
\int\limits_{C_{1}}\int\limits_{C_{2}}
t^{a}(t-1)^{b}(t-z)^{c}
s^{a}(s-1)^{b}(s-z)^{c}(t-s)^{g}dtds
\end{equation}
where $C_{1}$ and $C_{2}$ are as in corollary \ref{converge}, converge and define
multi--valued holomorphic functions of $z$. Then, these satisfy the \DF equation
\eqref{DF}.
\end{proposition}
\proof The strategy is similar to that used for the hypergeometric equation above.
We start from a given cohomological identity involving $\w$ and some additional
terms and gradually eliminate these by using other cohomological identities. Write
\begin{equation}\label{basic identity}
\begin{split}
0=\wt d[\frac{ds}{(t-z)^{2}}]
&=
(c-2)[\frac{dt\wedge ds}{(t-z)^{3}}]-
g[\frac{dt\wedge ds}{(t-z)^{2}(s-z)}]+
(\frac{a}{z}+\frac{b}{z-1})[\frac{dt\wedge ds}{(t-z)^{2}}]\\
&-
(\frac{a}{z^{2}}+\frac{b}{(z-1)^{2}})[\frac{dt\wedge ds}{t-z}]+
\frac{a}{z^{2}}[\frac{dt\wedge ds}{t}]+
\frac{b}{(z-1)^{2}}[\frac{dt\wedge ds}{t-1}] \\
&=
-c^{-1}(c-1)^{-1}w'''+
c^{-1}(c-1)^{-1}(\frac{a}{z}+\frac{b}{z-1})\w''+
c^{-1}(\frac{a}{z^{2}}+\frac{b}{(z-1)^{2}})\w'\\
&-
(3c+g)[\frac{dt\wedge ds}{(t-z)^{2}(s-z)}]-
c(c-1)^{-1}(\frac{a}{z}+\frac{b}{z-1})[\frac{dt\wedge ds}{(t-z)(s-z)}]\\
&+
\frac{a}{z^{2}}[\frac{dt\wedge ds}{t}]+
\frac{b}{(z-1)^{2}}[\frac{dt\wedge ds}{t-1}]
\end{split}
\end{equation}

The last two terms may be expressed in terms of $\w$ and its derivatives by writing
\begin{equation}
\begin{split}
0=\wt d[(t-1)ds]
&=
(1+a+b+c+\half{g})[dt\wedge ds]-a[\frac{dt\wedge ds}{t}]
+c(z-1)[\frac{dt\wedge ds}{t-z}]\\
&=
2(1+a+b+c+\half{g})\w-(z-1)\w'-a[\frac{dt\wedge ds}{t}]
\end{split}
\end{equation}
and
\begin{equation}
\begin{split}
0=\wt d[tds]
&=
(1+a+b+c+\half{g})[dt\wedge ds]+b[\frac{dt\wedge ds}{t-1}]+
cz[\frac{dt\wedge ds}{t-z}]\\
&=
2(1+a+b+c+\half{g})\w-z\w'+b[\frac{dt\wedge ds}{t-1}]
\end{split}
\end{equation}

To eliminate $[(t-z)^{-1}(s-z)^{-1}dt\wedge ds]$ from \eqref{basic identity},
we use
\begin{equation}
\begin{split}
0=\wt d[\frac{ds}{t-z}]
&=
(c-1)[\frac{dt\wedge ds}{(t-z)^{2}}]-
\half{g}[\frac{dt\wedge ds}{(t-z)(s-z)}]+
(\frac{a}{z}+\frac{b}{z-1})[\frac{dt\wedge ds}{t-z}]\\
&-
\frac{a}{z}[\frac{dt\wedge ds}{t}]-
\frac{b}{z-1}[\frac{dt\wedge ds}{t-1}]\\
&=
 c^{-1}\w''
-c^{-1}(\frac{a}{z}+\frac{b}{z-1})\w'\\
&-
(c+\half{g})[\frac{dt\wedge ds}{(t-z)(s-z)}]
-\frac{a}{z}[\frac{dt\wedge ds}{t}]-
\frac{b}{z-1}[\frac{dt\wedge ds}{t-1}]
\end{split}
\end{equation}

Finally, to eliminate $[(t-z)^{-2}(s-z)^{-1}dt\wedge ds]$ from \eqref{basic identity},
write
\begin{equation}
\begin{split}
0=\wt d[\frac{ds}{(t-z)(s-z)}]
&=
(c-1)[\frac{dt\wedge ds}{(t-z)^{2}(s-z)}]+
\left(\frac{a}{z}+\frac{b}{z-1}\right)[\frac{dt\wedge ds}{(t-z)(s-z)}]\\
&-
\frac{a}{z}[\frac{dt\wedge ds}{(t-z)s}]-
\frac{b}{z-1}[\frac{dt\wedge ds}{(t-z)(s-1)}]
\end{split}
\end{equation}

This introduces the additional terms $[(t-z)^{-1}s^{-1}dt\wedge ds]$
and $[(t-z)^{-1}(s-1)^{-1}dt\wedge ds]$ which may in turn be eliminated
by using
\begin{align}
0=\wt d[\frac{ds}{s-z}]
&=
(c+\half{g})[\frac{dt\wedge ds}{(t-z)(s-z)}]+
a[\frac{dt\wedge ds}{(t-z)s}]+
b[\frac{dt\wedge ds}{(t-z)(s-1)}]\\
0=\wt d[(t-z)\frac{ds}{s-z}]
&=
(1+a+b+c+g)[\frac{dt\wedge ds}{t-z}]-
az[\frac{dt\wedge ds}{(t-z)s}]-
b(z-1)[\frac{dt\wedge ds}{(t-z)(s-1)}]
\end{align}

A straightforward check now shows that, up to coboundaries, $\w$ satisfies \eqref{DF}.
To conclude, notice that if $f$ is one of the differential forms of which we have taken
the twisted differential, then
\begin{equation}
\int\limits_{C_{1}\times C_{2}}\Phi\wt df=
\int\limits_{C_{1}\times C_{2}}d(\Phi f)=
\int\limits_{\partial(C_{1}\times C_{2})}\Phi f=0
\end{equation}
since, by \eqref{range 1}--\eqref{range 2}, $\Phi f$ vanishes on 
$\partial(C_{1}\times C_{2})=
 C_{1}(\partial I)\times C_{2}\cup C_{2}\times C_{2}(\partial I)$. Thus if
$P(\frac{d}{dz})$ is the \DF differential operator,
$\int_{C_{1}\times C_{2}}P(\frac{d}{dz})\w=0$ and it follows that
$\int_{C_{1}\times C_{2}}\w$ satisfies the \DF equation since
$\int \frac{d}{dz}\w=\frac{d}{dz}\int\w$.
This is clear if the end--points of the $C_{i}$ lie in $\{0,1,\infty\}$ and follow in
general because $\Phi$ vanishes on $\partial(C_{1}\times C_{2})$ by 
\eqref{range 1}--\eqref{range 2} \halmos

\ssubsection{Monodromy properties of contour integrals}\label{ss:monodromy}

Following \cite{DF}, we show below how suitable choices of contours in the integral
representation \eqref{euler} lead to basis of solutions of the \DF equation diagonalising
the monodromy at a given singular point $z_{0}\in\{0,1,\infty\}$. We restrict our
attention to $z_{0}=0$ and $\infty$ since these are the only cases we shall need.
For graphical simplicity, the contours are represented with $z$ lying on the negative
real axis $(-\infty,0)$.
To fix conventions, all functions of the form $w^{\lambda}$ where $w$ is a function of
$t_{1},t_{2}$ will be defined on the $w$ plane with a cut along $[-\infty,0]$ so that 
$\arg w\in(-\pi,\pi)$ and 
\begin{equation}\label{convention}
(-w)^{\lambda}=
  \begin{cases}
    w^{\lambda} e^{-i\pi\lambda}&\text{if $\Im w>0$} \\
    w^{\lambda} e^{i\pi\lambda}&\text{if $\Im w<0$}
  \end{cases}
\end{equation}
The complex powers $z^{\lambda}$ however are defined with a cut along $z\in[0,\infty]$.

\begin{proposition}\label{DF solutions}
Let the parameters $a,b,c,g$ of the \DF equation lie in the range
\eqref{range 1}--\eqref{range 2}. Then the integrals \eqref{euler} corresponding to the
double contours
\begin{xalignat}{2}
I_{0,1}		&: \I01 &
I_{\infty,1}	&: \I21 \\
\nonumber \\
I_{0,2}		&: \I02 &
I_{\infty,2}	&: \I22 \\
\nonumber \\
I_{0,3}		&: \I03 &
I_{\infty,3}	&: \I23 \\
\nonumber
\end{xalignat}
yield solutions having the following monodromy at $z=0$ and $z=\infty$ respectively.
\begin{xalignat}{2}
I_{0,1}&=z^{0}(\rho_{0,1}+zO(z))		&
I_{\infty,1}&=
\Bigl(\frac{1}{z}\Bigr)^{-2c}(\rho_{\infty,1}+\frac{1}{z}O(\frac{1}{z}))\\
I_{0,2}&=z^{1+a+c}(\rho_{0,2}+zO(z))		&
I_{\infty,2}&=
\Bigl(\frac{1}{z}\Bigr)^{-(1+a+b+2c+g)}(\rho_{\infty,2}+\frac{1}{z}O(\frac{1}{z}))\\
I_{0,3}&=z^{2(1+a+c+\half{g})}(\rho_{0,3}+zO(z))&
I_{\infty,3}&=
\Bigl(\frac{1}{z}\Bigr)^{-2(1+a+b+c+\half{g})}(\rho_{\infty,3}+\frac{1}{z}O(\frac{1}{z}))
\end{xalignat}
The leading coefficients are given by
\begin{xalignat}{2}
\rho_{0,1}	&=\half{1}(1+e^{i\pi g})J_{2}(-(2+a+b+c+g),b;g)	&
\rho_{\infty,1}	&=\half{1}(1+e^{i\pi g})J_{2}(a,b;g)		\\
\rho_{0,2}	&=B(a+1,c+1)B(-(1+a+b+c),b+1)			&
\rho_{\infty,2}	&=B(a+1,b+1)B(-(1+a+b+c),c+1)			\\
\rho_{0,3}	&=\half{1}(1+e^{i\pi g})J_{2}(a,c;g)		&
\rho_{\infty,3}	&=\half{1}(1+e^{i\pi g})J_{2}(-(2+a+b+c+g),c;g)
\end{xalignat}
where $B$ is the Euler $\beta$ function, and do not vanish if $g\neq -1$.
\end{proposition}
\proof
Consider $I_{\infty,1}$ first. The integrand in $I_{\infty,1}$ may be written,
up to a complex phase factor as $z^{2c}\psi(z)$ where
$\psi(z)=
t_{1}^{a}(t_{1}-1)^{b}(1-t_{1}z^{-1})^{c}
t_{2}^{a}(t_{2}-1)^{b}(1-t_{2}z^{-1})^{c}(t_{1}-t_{2})^{g}$
is single--valued on $|z|>1$ when $|t_{1}|,|t_{2}|<1$.
Thus, $I_{\infty,1}=\Bigl(\frac{1}{z}\Bigr)^{-2c}\Psi(\frac{1}{z})$ where $\Psi$
is holomorphic and single--valued in a neighborhood of $0$.
The leading term of $\Psi$ may be computed by deforming successively the contours
$C_{1}$, $C_{2}$ and using Selberg's formula
\begin{equation}
\begin{split}
\Psi(0)
&= \int_{C_{1}}\int_{C_{2}}
t_{1}^{a}(t_{1}-1)^{b}t_{2}^{a}(t_{2}-1)^{b}(t_{1}-t_{2})^{g}dt_{1}dt_{2}\\
&= \int_{0}^{1}\int_{C_{2}}
t_{1}^{a}(1-t_{1})^{b}t_{2}^{a}(1-t_{2})^{b}(t_{1}-t_{2})^{g}dt_{1}dt_{2}\\
&= \int_{0}^{1}\int_{C_{2}^{1}}
t_{1}^{a}(1-t_{1})^{b}t_{2}^{a}(1-t_{2})^{b}(t_{1}-t_{2})^{g}dt_{1}dt_{2}+
   \int_{0}^{1}\int_{C_{2}^{2}}
t_{1}^{a}(1-t_{1})^{b}t_{2}^{a}(1-t_{2})^{b}(t_{1}-t_{2})^{g}dt_{1}dt_{2}
\end{split}
\end{equation}
where the $C_{2}^{1}$ and $C_{2}^{2}$ lie below the real axis and join $0$ to $t_{1}$ and
$t_{1}$ to $1$ respectively. Since $\Re(t_{1}-t_{2})>0$ along $C_{2}^{1}$ we may deform it
onto the segment $(0,t_{1})$ to get a contribution
\begin{equation}
\int_{0}^{1}\int_{0}^{t_{1}}
t_{1}^{a}(1-t_{1})^{b}t_{2}^{a}(1-t_{2})^{b}|t_{1}-t_{2}|^{g}dt_{1}dt_{2}=
\half{1}\int_{0}^{1}\int_{0}^{1}
t_{1}^{a}(1-t_{1})^{b}t_{2}^{a}(1-t_{2})^{b}|t_{1}-t_{2}|^{g}dt_{1}dt_{2}
\end{equation}
On $C_{2}^{2}$, \eqref{convention} gives $(t_{1}-t_{2})^{g}=e^{i\pi g}(t_{2}-t_{1})^{g}$. 
Deforming onto the real axis we get the contribution
\begin{equation}
e^{i\pi g}\int_{0}^{1}\int_{t_{1}}^{1}
t_{1}^{a}(1-t_{1})^{b}t_{2}^{a}(1-t_{2})^{b}|t_{1}-t_{2}|^{g}dt_{1}dt_{2}=
e^{i\pi g}\half{1}\int_{0}^{1}\int_{0}^{1}
t_{1}^{a}(1-t_{1})^{b}t_{2}^{a}(1-t_{2})^{b}|t_{1}-t_{2}|^{g}dt_{1}dt_{2}
\end{equation}
And therefore, $\Psi(0)=\half{1}(1+e^{i\pi g})J_{2}(a,b;g)$ which, in view of the remark
following proposition \ref{selberg} vanishes iff $g\in 2\IZ+1$.
However, by \eqref{range 1}--\eqref{range 2}, $\Re g> -1$ and $\Re g<-1-\Re(a+b+c)<2$ and it
follows that $\Psi(0)$ is non--zero if $g\neq 1$.
$I_{\infty,3}$ is treated similarly by setting $t_{i}=z u_{i}^{-1}$, where
the $u_{i}$ lie on contours $C_{i}'$ joining $0$ to $1$. Performing the change
of variables, one finds
\begin{equation}
\begin{split}
I_{\infty,3}
=&\Bigl(\frac{1}{z}\Bigr)^{-2(1+a+b+c+\half{g})}
  \int_{C_{1}'}\int_{C_{2}'}du_{1}du_{2}\\
 &
u_{1}^{-(2+a+b+c+g)}(1-u_{1}z^{-1})^{b}(1-u_{1})^{c}
u_{2}^{-(2+a+b+c+g)}(1-u_{2}z^{-1})^{b}(1-u_{2})^{c}
(u_{1}-u_{2})^{g}
\end{split}
\end{equation}
The coefficient of the leading term at $z=\infty$ is obtained as for $I_{\infty,1}$ and is
equal to $\half{1}(1+e^{i\pi g})J_{2}(-(2+a+b+c+g),c;g)$ so that it does not vanish if
$g\neq 1$. Finally, $I_{\infty,2}$ is obtained by setting $t_{2}=zu_{2}^{-1}$ where
$u_{2}\in C_{2}'$ runs from $0$ to $1$. This yields
\begin{equation}
\begin{split}
I_{2,\infty}
=&-\Bigl(\frac{1}{z}\Bigr)^{-(1+a+b+2c+g)}\\
 &\int_{C_{1}}\int_{C_{2}'}
t_{1}^{a}(t_{1}-1)^{b}(1-t_{1}z^{-1})^{c}
u_{2}^{-(2+a+b+c+g)}(1-u_{2}z^{-1})^{b}(1-u_{2})^{c}
(1-t_{1}u_{2}z^{-1})^{g}dt_{1}du_{2}
\end{split}
\end{equation}
with leading coefficient
\begin{equation}
\int_{C_{1}}t_{1}^{a}(t_{1}-1)^{b}dt_{1}
\int_{C_{2}'}u_{2}^{-(2+a+b+c+g)}(1-u_{2})^{c}du_{2}=
B(a+1,b+1)B(-(1+a+b+c+g),c+1)\neq 0
\end{equation}
The cases $I_{0,i}$, $i=1\ldots 3$ follow similarly \halmos

\ssubsection{Connection matrices for the Dotsenko--Fateev equation}\label{ss:continuation}

We compute below the analytic continuation of the solution $I_{0,1}$ of the \DF
equation equation from $0$ to $\infty$ and re--express it in terms of the
solutions $I_{\infty,i}$ diagonalising the monodromy at $\infty$. Following
\cite{DF}, we use the technique of contour deformation.
This consists in deforming the double contour corresponding to $I_{0,1}$ 
in various ways. Combinining these with appropriate phase corrections, one obtains
a linear combination of contours corresponding to the solutions $I_{\infty,i}$.
The same method applies to all solutions $I_{0,i}$ and works equally well to compute
the analytic continuation of solutions from 0 to 1, see \cite{DF}. For later use
however, we only need the results pertaining to $I_{0,1}$.
We shall compute the analytic continuation along the negative real axis and therefore
assume that $z\in(-\infty,0)$. We adhere to the notations and conventions
of \S \ref{ss:monodromy}.

\begin{proposition}\label{continuation}
Let $I_{z_{i},j}$ be the solutions of the \DF equation given by proposition
\ref{DF solutions}. Then
\begin{equation}\label{identity}
\begin{split}
I_{0,1}
&= 
\frac{\esse(a)\esse(a+\half{g})}{\esse(a+b+\half{g})\esse(a+b+g)}I_{\infty,1}\\
&+
2e^{-\pi i(a+c+\half{g})}\ci(\half{g})
\frac{\esse(a)\esse(c)}{\esse(a+b)\esse(a+b+g)}I_{\infty,2}\\
&+
\frac{\esse(c)\esse(c+\half{g})}{\esse(a+b+\half{g})\esse(a+b+\half{g})}I_{\infty,3}
\end{split}
\end{equation}
where $\esse(\alpha)=\sin(\pi\alpha)$ and $\ci(\alpha)=\cos(\pi\alpha)$.
\end{proposition}
\proof
If we take the equator of $\IP^{1}$ as the line through $0,1,\infty$ and
$z$, then deforming the $t_{1}$ contour in the northern hemisphere or the
$t_{2}$ contour in the southern one gives

\begin{equation*}
\I01 =
\begin{cases}
 &\picturebegin
 \axis\leftuparc 0,1\leftuparc 1,1\leftuparc 2,1\rightdownarc 3,2
 \pictureend \\
 &\\
 &\picturebegin
 \axis\rightuparc 3,1\leftdownarc 0,2\leftdownarc 1,2\leftdownarc 2,2
 \pictureend
 =
 e^{\pi ig}
 \picturebegin
 \axis\leftdownarc 0,1\leftdownarc 1,1\leftdownarc 2,1\rightuparc 3,2
 \pictureend
\end{cases}
\end{equation*}
\hfill\break
where the last equality follows from the change of variables
$t_{1}\rightarrow t_{2}$, $t_{2}\rightarrow t_{1}$ (so that 
$(t_{1}-t_{2})^{g}\rightarrow(t_{2}-t_{1})^{g}=e^{\pi ig}(t_{1}-t_{2})^{g}$
by \eqref{convention} since since $\Im(t_{1}-t_{2})<0$). Now
\begin{equation*}
\begin{split}
\picturebegin
\axis\leftdownarc 0,1\leftdownarc 1,1\leftdownarc 2,1\rightuparc 3,2
\pictureend
&=
\picturebegin
\axis\leftdownarc 2,1\rightuparc 3,2
\pictureend
+
\picturebegin
\axis\leftdownarc 1,1\rightuparc 3,2
\pictureend
+
\picturebegin
\axis\leftdownarc 0,1\rightuparc 3,2
\pictureend \\
&=e^{-2\pi i(b+g)}
\picturebegin
\axis\leftuparc 2,1\rightdownarc 3,2
\pictureend
+e^{-2\pi i(a+b+g)}
\picturebegin
\axis\leftuparc 1,1\rightdownarc 3,2
\pictureend \\
&+e^{-2\pi i(a+b+c+g)}
\picturebegin
\axis\leftuparc 0,1\rightdownarc 3,2
\pictureend \\
\end{split}
\end{equation*}
and therefore
\begin{equation}
\begin{split}
\left(1-e^{2\pi i(a+b+\half{g})}\right)\I01
&=
\left(1-e^{2\pi ia}\right)
\picturebegin
\axis\leftuparc 2,1\rightdownarc 3,2
\pictureend \\
&+
\left(1-e^{-2\pi ic}\right)
\picturebegin
\axis\leftuparc 0,1\rightdownarc 3,2
\pictureend \\ \label{C11}
\end{split}
\end{equation} 
Proceeding in a similar fashion, we may deform the first summand on the right
hand--side of \eqref{C11}
\begin{equation*}
\picturebegin
\axis\leftuparc 2,1\rightdownarc 3,2
\pictureend
=\begin{cases}
&
\picturebegin
\axis\leftuparc 2,1\leftdownarc 0,2\leftdownarc 1,2\leftdownarc 2,2
\pictureend \\
& \\
e^{2\pi i(b+g)}
&\picturebegin
\axis\leftdownarc 2,1\rightuparc 3,2
\pictureend 
=
e^{2\pi i(b+g)}
\picturebegin
\axis\leftdownarc 2,1\leftuparc 0,2\leftuparc 1,2\leftuparc 2,2
\pictureend
\end{cases}
\end{equation*}
\begin{equation*}
\begin{split}
\picturebegin
\axis\leftdownarc 2,1\leftuparc 0,2\leftuparc 1,2\leftuparc 2,2
\pictureend
&=
\picturebegin
\axis\leftdownarc 2,1\leftuparc 2,2
\pictureend
+
\picturebegin
\axis\leftdownarc 2,1\leftuparc 1,2
\pictureend
+
\picturebegin
\axis\leftdownarc 2,1\leftuparc 0,2
\pictureend \\
&=
e^{-\pi ig}
\picturebegin
\axis\leftuparc 2,1\leftdownarc 2,2
\pictureend
+
e^{2\pi i(-b+b+a)}
\picturebegin
\axis\leftuparc 2,1\leftdownarc 1,2
\pictureend \\
&+
e^{2\pi i(-b+a+b+c)}
\picturebegin
\axis\leftuparc 2,1\leftdownarc 0,2
\pictureend \\
\end{split}
\end{equation*}
and get
\begin{equation}
\begin{split}
\left(1-e^{-2\pi i(a+b+g)}\right)
\picturebegin
\axis\leftuparc 2,1\leftdownarc 3,2
\pictureend
&=
\left(1-e^{-2\pi i(a+\half{g})}\right)
\I13 \\
&-
\left(1-e^{2\pi ic}\right)
\I22 \\ \label{C12}
\end{split}
\end{equation}
The second summand on the other hand yields
\begin{equation*}
\picturebegin
\axis\leftuparc 0,1\rightdownarc 3,2
\pictureend
=\begin{cases}
&
\picturebegin
\axis\leftuparc 0,1\leftdownarc 0,2\leftdownarc 1,2\leftdownarc 2,2
\pictureend \\
& \\
e^{2\pi i(a+b+c+g)}&
\picturebegin
\axis\leftdownarc 0,1\rightuparc 3,2
\pictureend 
=
e^{2\pi i(a+b+c+g)}
\picturebegin
\axis\leftdownarc 0,1\leftuparc 0,2\leftuparc 1,2\leftuparc 2,2
\pictureend
\end{cases}
\end{equation*}
\begin{equation*}
\begin{split}
\picturebegin
\axis\leftdownarc 0,1\leftuparc 0,2\leftuparc 1,2\leftuparc 2,2
\pictureend
&=
\picturebegin
\axis\leftdownarc 0,1\leftuparc 0,2
\pictureend
+
\picturebegin
\axis\leftdownarc 0,1\leftuparc 1,2
\pictureend
+
\picturebegin
\axis\leftdownarc 0,1\leftuparc 2,2
\pictureend \\
& \\
&=
e^{-\pi ig}
\picturebegin
\axis\leftuparc 0,1\leftdownarc 0,2
\pictureend
+
e^{2\pi i(-a-b-c-g+a+b)}
\picturebegin
\axis\leftuparc 0,1\leftdownarc 1,2
\pictureend \\
&+
e^{2\pi i(-a-b-c-g+b)}
\picturebegin
\axis\leftuparc 0,1\leftdownarc 2,2
\pictureend \\
\end{split}
\end{equation*}
so that
\begin{equation}\label{C13}
\begin{split}
\left(1-e^{-2\pi i(a+b)}\right)
\picturebegin
\axis\leftuparc 0,1\rightdownarc 3,2
\pictureend
&=
\left(1-e^{2\pi i(c+\half{g})}\right)
\I11 \\
&+
\left(1-e^{-2\pi ia}\right)
\picturebegin
\axis\leftuparc 0,1\leftdownarc 2,2
\pictureend \\
&=
\left(1-e^{2\pi i(c+\half{g})}\right)
\I11 \\
&-
\left(1-e^{-2\pi ia}\right)e^{2\pi i(a+c+\half{g})}
\I22 \\ 
\end{split}
\end{equation}

where the last identity follows from

\begin{equation*}
\picturebegin
\axis\leftuparc 0,1\leftdownarc 2,2
\pictureend
=
e^{\pi ig}
\picturebegin
\axis\leftdownarc 2,1\leftuparc 0,2
\pictureend
=
-e^{2\pi i(\half{g}-b+a+b+c)}
\I22
\end{equation*}
\hfill\break
Substituting \eqref{C12} and \eqref{C13} in \eqref{C11} gives finally
\begin{equation*}
\begin{split}
\left( 1-e^{2\pi i(a+b+\half{g})} \right)
\left( 1-e^{-2\pi i(a+b+g)} \right)
\left( 1-e^{-2\pi i(a+b)} \right)
I^{0}_{1} 
&= \\
\left( 1-e^{2\pi ia} \right)
\left( 1-e^{-2\pi i(a+b)} \right)
\left( 1-e^{-2\pi i(a+\half{g})} \right)
I^{\infty}_{1}
&- \\
\left(
\left( 1-e^{2\pi ia} \right)
\left( 1-e^{-2\pi i(a+b)} \right)
\left( 1-e^{2\pi ic} \right)
e^{-2\pi i(a+c+\half{g})} \right.
&+ \\
\left.
\left( 1-e^{-2\pi ic} \right)
\left( 1-e^{-2\pi i(a+b+g)} \right)
\left( 1-e^{-2\pi ia} \right) \right)
I^{\infty}_{2} 
&+ \\
\left( 1-e^{-2\pi ic} \right)
\left( 1-e^{-2\pi i(a+b+g)} \right)
\left( 1-e^{2\pi i(c+\half{g})} \right)
I^{\infty}_{3}
& 
\end{split}
\end{equation*}
which, after simplification, yields \eqref{identity} \halmos

\ssection{The Knizhnik--Zamolodchikov equations with values in
$(V_{k\theta_{1}}\otimes V_{k\theta_{1}}^{*}\otimes 
  V_{\theta_{1}} \otimes V_{\theta_{1}}^{*})^{\SO_{2n}}$}\label{se:box/kbox equation}

In the rest of this chapter, we label representations of $\Spin_{2n}$ by their highest
weight. Thus, $V_{\theta_{1}}$ and $V_{k\theta_{1}}$ denote the vector representation of
$\SO_{2n}$ and its k--fold symmetric, traceless power respectively.
In \S \ref{ss:specific omega}, we compute explicitly the residue matrices $\Omega_{ij}$
of the \KZ equations corresponding to the tensor product
$(V_{k\theta_{1}}\otimes V_{k\theta_{1}}^{*}\otimes 
  V_{\theta_{1}} \otimes V_{\theta_{1}}^{*})^{\SO_{2n}}$. The calculation depends on
some formulae of Christe and Flume \cite{CF} which are obtained in \S \ref{ss:general omega}.

\ssubsection{The $\Omega_{ij}$ matrices for
$(V_{4}\otimes V_{3}^{*}\otimes V_{2}^{*}\otimes V_{1}^{*})^{G}$}\label{ss:general omega}

Let $G$ be a compact, connected simple Lie group.
Recall from chapter \ref{ch:algebraic fields} that to any tensor product
$(V_{n}\otimes V_{n-1}^{*}\otimes\cdots\otimes V_{1}^{*})^{G}$ of irreducible
$G$--modules one associates the endomorphisms $\Omega_{ij}=\pi_{i}(X_{k})\pi_{j}(X^{k})$,
$1\leq i,j\leq n$ where $X_{k}$, $X^{k}$ are basis of $\gc$ dual with respect to the
basic inner product. We compute below some matrix entries of the $\Omega_{ij}$ when
$n=4$, in any tensor product basis corresponding to the the isomorphism
\begin{equation}\label{decomposition}
(V_{4}\otimes V_{3}^{*}\otimes V_{2}^{*}\otimes V_{1}^{*})^{G}
\cong\bigoplus_{U}\Hom_{G}(U,V_{3}^{*}\otimes V_{4})\otimes\Hom_{G}(V_{1}\otimes V_{2},U)
\end{equation}
where $U$ ranges over the irreducible summands of $V_{1}\otimes V_{2}$.
The inner product is then given by
$ (\psi_{U}\otimes\phi_{U},\psi_{W}\otimes\phi_{W})=
  \tr_{V_{1}\otimes V_{2}}(\phi_{W}^{*}\psi_{W}^{*}\psi_{U}\phi_{U})=
  \tr_{W}(\psi_{W}^{*}\psi_{U}\phi_{U}\phi_{W}^{*})$.
The following formulae are somewhat imprecisely stated in \cite{CF}

\begin{proposition}\label{cf rules}
Let $\Omega_{12},\Omega_{23}$ be the operators corresponding to the tensor
product $(V_{4}\otimes V_{3}^{*}\otimes V_{2}^{*}\otimes V_{1}^{*})^{G}$.
Then
\begin{enumerate}
\item $\Omega_{12}$ is diagonal in any basis corresponding to the decomposition
\eqref{decomposition}. Its diagonal entries are
\begin{equation}\label{conformal weight}
(\Omega_{12})_{\psi_{U}\otimes\phi_{U},\psi_{U}\otimes\phi_{U}}=
\half{1}(C_{U}-C_{1}-C_{2})=:\delta_{U}
\end{equation}
\item If $U=\IC$ is the trivial representation, then
\begin{equation}
(\Omega_{23})_{\psi_{U}\otimes\phi_{U},\psi_{U}\otimes\phi_{U}}=0
\end{equation}
\item If $\Hom_{G}(U\otimes\gc,U)=\IC$, so that $U\neq 0$ and $C_{U}>0$ then
\begin{equation}
(\Omega_{23})_{\psi_{U}\otimes\phi_{U},\psi_{U}\otimes\phi_{U}}=
-\frac{1}{4C_{U}}(C_{U}+C_{3}-C_{4})(C_{U}+C_{2}-C_{1})
\end{equation}
\item If $\Hom_{G}(U\otimes\gc,U')=0$, then
\begin{equation}
(\Omega_{23})_{\psi_{U}\otimes\phi_{U},\psi_{U'}\otimes\phi_{U'}}=
(\Omega_{23})_{\psi_{U'}\otimes\phi_{U'},\psi_{U}\otimes\phi_{U}}=
0
\end{equation}
\end{enumerate}
\end{proposition}
\proof
(i) From the $G$-equivariance of $\phi_{U}$, we get
\begin{equation}
\begin{split}
\Omega_{12}\psi_{U}\otimes\phi_{U}
&=
\psi_{U}\phi_{U}\pi_{2}(X_{k})\pi_{1}(X^{k})\\
&=
\half{1}\psi_{U}\phi_{U}(\pi_{2}(X_{k})+\pi_{1}(X_{k}))(\pi_{2}(X^{k})+\pi_{1}(X^{k}))\\
&-
\half{1}\psi_{U}\phi_{U}(\pi_{2}(X_{k})\pi_{2}(X^{k})+\pi_{1}(X_{k})\pi_{1}(X^{k}))\\
&=
\half{1}\psi_{U}\pi_{U}(X_{k})\pi_{U}(X^{k})\phi_{U}-
\half{1}(C_{1}+C_{2})\psi_{U}\phi_{U}\\
&=
\half{1}(C_{U}-C_{1}-C_{2})\psi_{U}\otimes\pi_{U}
\end{split}
\end{equation}
(ii) The map $\gc\rightarrow\IC$, $X\rightarrow\phi_{U}\pi_{2}(X)\phi_{U}^{*}$ is
$G$-invariant and therefore zero. Thus
\begin{equation}
(\Omega_{23}\psi_{U}\otimes\phi_{U},\psi_{U}\otimes\phi_{U})=
-\tr_{U}(\psi_{U}^{*}\pi_{\overline{3}}(X_{k})\psi_{U}\phi_{U}\pi_{2}(X^{k})\phi_{U}^{*})=
0
\end{equation}
(iii) The map $\gc\otimes U\rightarrow U$,
$X\otimes u\rightarrow\phi_{U}\pi_{2}(X)\phi_{U}^{*}u$ commutes with $G$ and is
therefore proportional to $u\otimes X\rightarrow\pi_{U}(X)u$. Similarly,
$\psi_{U}^{*}\pi_{\overline{3}}(X)\psi_{U}=c_{\psi}\pi_{U}(X)$ for any $X\in\g$. It
follows that
\begin{equation}
(\Omega_{23}\psi_{U}\otimes\phi_{U},\psi_{U}\otimes\phi_{U})=
-\tr_{U}(\psi_{U}^{*}\pi_{\overline{3}}(X_{k})\psi_{U}\phi_{U}\pi_{2}(X^{k})\phi_{U}^{*})=
-c_{\psi}c_{\phi}\dim(U)C_{U}
\end{equation}
To evaluate the constants $c_{\psi},c_{\phi}$, write
\begin{equation}
\pi_{U}(X^{k})\psi_{U}^{*}\pi_{\overline{3}}(X_{k})\psi_{U}=
c_{\psi}\pi_{U}(X^{k})\pi_{U}(X_{k})
\end{equation}
The right hand--side yields $c_{\psi}C_{U}$ while the left hand--side gives, as in (i)
$\half{1}(C_{4}-C_{U}-C_{3})\psi_{U}^{*}\psi_{U}$. The computation of $c_{\phi}$ is
similar. Summarising,
\begin{equation}
\begin{split}
(\Omega_{23}\psi_{U}\otimes\phi_{U},\psi_{U}\otimes\phi_{U})
&=-\frac{1}{4C_{U}}(C_{U}+C_{3}-C_{4})(C_{U}+C_{2}-C_{1})
\dim(U)\psi_{U}^{*}\psi_{U}\phi_{U}\phi_{U}^{*}\\
&=-\frac{1}{4C_{U}}(C_{U}+C_{3}-C_{4})(C_{U}+C_{2}-C_{1})
\|\psi_{U}\otimes\phi_{U}\|^{2}
\end{split}
\end{equation}
(iv) The map $\gc\otimes U\rightarrow U'$,
$X\otimes u\rightarrow\psi_{U'}^{*}\pi_{\overline{3}}(X)\psi_{U}u$ commutes with $G$
and therefore vanishes. It follows that
\begin{equation}
(\Omega_{23}\psi_{U}\otimes\phi_{U},\psi_{U'}\otimes\phi_{U'})=
 -\tr_{U'}(\psi_{U'}^{*}\pi_{\overline{3}}(X_{k})\psi_{U}
           \phi_{U}\pi_{2}(X^{k})\phi_{U'}^{*})=0
\end{equation}
The second identity follows by permuting $U$ and $U'$ and noticing that the adjoint
representation is real and therefore self-dual so that
\begin{equation}
\Hom_{G}(U'\otimes\gc,U)=\Hom_{G}(U',U\otimes\gc)=\Hom_{G}(U\otimes\gc,U')^{*}
\end{equation}
\halmos\\

\remark Clearly $\Omega_{13}$ and $\Omega_{23}$ are diagonal in any tensor product basis
of $W=(V_{4}\otimes V_{3}^{*}\otimes V_{2}^{*}\otimes V_{1}^{*})^{G}$ corresponding
respectively to the isomorphisms
\begin{align}
W&\cong
\bigoplus_{U}\Hom_{G}(U,V_{2}^{*}\otimes V_{4})\otimes\Hom_{G}(V_{1}\otimes V_{3},U)\\
\intertext{and}
W&\cong
\bigoplus_{U}\Hom_{G}(U,V_{1}^{*}\otimes V_{4})\otimes\Hom_{G}(V_{2}\otimes V_{3},U)
\end{align}
In particular, their eigenvalues may be labelled by the irreducible summands of
$V_{1}\otimes V_{3}$ and $V_{2}\otimes V_{3}$ respectively.

\ssubsection{The $\Omega_{ij}$ matrices for 
$(V_{k\theta_{1}}\otimes V_{k\theta_{1}}^{*}\otimes 
  V_{\theta_{1}} \otimes V_{\theta_{1}}^{*})^{\SO_{2n}}$}\label{ss:specific omega}

Specialising the results of \S \ref{ss:general omega} to $G=\Spin_{2n}$, we find

\begin{proposition}\label{pr:explicit}
Let $\Omega_{ij}$ be the operators corresponding to 
$W=(V_{k\theta_{1}}\otimes V_{k\theta_{1}}^{*}\otimes
    V_{\theta_{1}}\otimes V_{\theta_{1}}^{*})^{\SO_{2n}}$. Then
\begin{enumerate}
\item The eigenvalues of $\Omega_{12}$ are labelled by the summands of
\begin{equation}\label{tensor two}
V_{\theta_{1}}\otimes V_{\theta_{1}}^{*}=
V_{2\theta_{1}}\oplus V_{\theta_{1}+\theta_{2}}\oplus V_{0}
\end{equation}
and are given by
\begin{xalignat}{3}
\delta_{2\theta_{1}}&=1&
\delta_{\theta_{1}+\theta_{2}}&=-1&
\delta_{0}&=-2n+1
\end{xalignat}

\item $\Omega_{13}$ and $\Omega_{23}$ are conjugate. Their eigenvalues are
labelled by the summands of
\begin{equation}\label{tensor one}
V_{\theta_{1}}    \otimes V_{k\theta_{1}}\cong
V_{\theta_{1}}^{*}\otimes V_{k\theta_{1}}=
V_{(k+1)\theta_{1}}\oplus V_{k\theta_{1}+\theta_{2}}\oplus V_{(k-1)\theta_{1}}
\end{equation}
and are given by
\begin{xalignat}{3}
\wt\delta_{(k+1)\theta_{1}}&=k&
\wt\delta_{k\theta_{1}+\theta_{2}}&=-1&
\wt\delta_{(k-1)\theta_{1}}&=-2(n-1)-k
\end{xalignat}

\item In the basis corresponding to
\begin{equation}\label{basis}
W\cong\bigoplus_{U}
  \Hom_{\SO_{2n}}(U,V_{k\theta_{1}}\otimes V_{k\theta_{1}}^{*})\otimes
  \Hom_{\SO_{2n}}(V_{\theta_{1}}\otimes V_{\theta_{1}}^{*},U)
\end{equation}
where $U$ ranges over the summands of \eqref{tensor two}, $\Omega_{12}$, $\Omega_{23}$
and $\Omega_{13}$ are given by
\begin{xalignat}{3}\label{matrices}
\Omega_{12}&=
\left(\begin{array}{crc}
1&0&0\\0&-1&0\\0&0&-2n+1
\end{array}\right)&
\Omega_{23}&=
\left(\begin{array}{rcc}
-n   &x    &0\\
\wt x&-n+1 &y\\
0    &\wt y&0
\end{array}\right)&
\Omega_{13}&=
\left(\begin{array}{rcc}
-n    &-x    &0 \\
-\wt x&-n+1  &-y\\
0     &-\wt y&0
\end{array}\right)
\end{xalignat}
where
\begin{xalignat}{3}
x\wt x&=\frac{(n-1)}{n}(k^{2}+2k(n-1)+n(n-2))&
&\text{and}&
y\wt y&=\frac{k}{n}(k+2(n-1))
\end{xalignat}
\end{enumerate}
\end{proposition}
\proof
Let $V=V_{\theta_{1}}\cong V_{\theta_{1}}^{*}$ be the defining representation of
$\SO_{2n}$. $V$ is minimal with weights $\pm\theta_{j}$, $1\leq j\leq n$ and therefore
\eqref{tensor two} and \eqref{tensor one} follow from proposition
\ref{ch:classification}.\ref{pr:tensor with minimal}.
In particular,
\begin{equation}
W\cong\End_{\SO_{2n}}(V\otimes V_{k\theta_{1}})\cong\IC^{3}
\end{equation}
Since the Casimir $C_{\lambda}$ of a representation of highest weight $\lambda$ is
$\<\lambda,\lambda+2\rho\>$ where $2\rho=2\sum_{j}(n-j)\theta_{j}$ is the sum of the
positive roots of $\SO_{2n}$, the Casimirs of the summands of \eqref{tensor two} are
$C_{2\theta_{1}}=4n$,
$C_{\theta_{1}+\theta_{2}}=4(n-1)$ and $C_{0}=0$. It follows from proposition
\ref{cf rules} and $C_{\theta_{1}}=2n-1$ that the eigenvalues of $\Omega_{12}$
are $1,-1,-2n+1$. We claim that each has multiplicity one, {\it i.e.}~ that the spaces
$\Hom_{\SO_{2n}}(U,V_{k\theta_{1}}^{*}\otimes V_{k\theta_{1}})$ are one-dimensional for
$U=V_{2\theta_{1}}$, $U=V_{0}\cong\IC$ and $U=V_{\theta_{1}+\theta_{2}}\cong\so_{2n}$
so that $\Omega_{12}$ is given by \ref{pr:explicit}.
Since the $\Hom$ spaces are non--zero in the last two cases and $W$ is
three--dimensional, we will prove our claim by exhibiting a non--zero intertwiner
$V_{k\theta_{1}}\otimes V_{2\theta_{1}}\rightarrow V_{k\theta_{1}}$.
To this end, realise $V_{k\theta_{1}}$ as the space of traceless, symmetric $k$-tensors as
follows. Let $\skv\subset V^{\otimes k}$ be the k-fold symmetric tensor power of $V$
spanned by the tensors
\begin{equation}
[v_{1}\otimes\cdots\otimes v_{k}]=
\frac{1}{k!}\sum_{\sigma\in\mathfrak{S}_{k}}
v_{\sigma(1)}\otimes\cdots\otimes v_{\sigma(k)}
\end{equation}
Let moreover $B$ be the $\SO_{2n}$--invariant symmetric, bilinear form on $V$.
The natural trace map
\begin{equation}
\tau_{k}:V^{\otimes k}\rightarrow V^{\otimes(k-2)},
\thickspace
v_{1}\otimes\cdots\otimes v_{k}\rightarrow
\sum_{1\leq i<j\leq k}B(v_{i},v_{j})
v_{1}\otimes\cdots\otimes\wh{v_{i}}\otimes\cdots\otimes\wh{v_{j}}
     \otimes\cdots\otimes v_{k}
\end{equation}
satisfies
$\tau_{k}[v_{1}\otimes\cdots\otimes v_{k}]=[\tau_{k}(v_{1}\otimes\cdots\otimes v_{k})]$
and therefore restricts to a map $\tau_{k}:\skv\rightarrow S^{k-2}V$. The kernel of
$\left.\tau_{k}\right|_{\skv}$ is $V_{k\theta_{1}}$ and the intertwiner
$V_{k\theta_{1}}\otimes V_{2\theta_{1}}\rightarrow V_{k\theta_{1}}$ is given by
$[v_{1}\otimes\cdots\otimes v_{k}]\otimes[w_{1}\otimes w_{2}]\rightarrow
  P\tau_{k+2}[v_{1}\otimes\cdots\otimes v_{k}\otimes w_{1}\otimes w_{2}]$
where $P:\skv\rightarrow V_{k\theta_{1}}$ is the orthogonal projection.\\

Let now $i:V_{\theta_{1}}\rightarrow V_{\theta_{1}}^{*}$ be the $\SO_{2n}$ identification
and $\sigma$ the permutation on $V_{\theta_{1}}\otimes V_{\theta_{1}}$. Then,
$\Omega_{13}=F\Omega_{23}F^{-1}$ where $F=(1\otimes i)\sigma(1\otimes i^{-1})$.
Notice that $\sigma$ acts as multiplication by $1,-1,1$ respectively on the summands of
\eqref{tensor two} and therefore
\begin{equation}\label{eq:F matrix}
F=\left(\begin{array}{crc}1&0&0\\0&-1&0\\0&0&1\end{array}\right)
\end{equation}

To compute $\Omega_{23}$, notice that
the Casimirs of the summands of \eqref{tensor one} are $C_{(k+1)\theta_{1}}=k(k+2(n-1))$,
$C_{k\theta_{1}+\theta_{2}}=k^{2}+2(k+1)(n-1)-1$ and
$C_{(k-1)\theta_{1}}=(k-1)(k-1+2(n-1))$.
Since $\Omega_{23}$ acts by pre-multiplication on $V_{k\theta_{1}}\otimes V^{*}$ when
$W$ is described as
\begin{equation}
 \End_{\SO_{2n}}(V_{k\theta_{1}}\otimes V^{*})\cong
 \bigoplus_{U}\Hom_{\SO_{2n}}(U,V_{k\theta_{1}}\otimes V^{*})\otimes
              \Hom_{\SO_{2n}}(V_{k\theta_{1}}\otimes V^{*},U)
\end{equation}
it is diagonal in the basis given by the orthogonal projections onto the irreducible
summands of $V_{k\theta_{1}}\otimes V^{*}$. The corresponding eigenvalues are, by (i)
of proposition
\ref{cf rules}, $k,-1,-2(n-1)-k$. Finally, by the Brauer rules for computing tensor
products \cite[\S 24.4, Ex. 9]{Hu} and the fact that the roots of $\so_{2n}$ are
$\pm(\theta_{i}\pm\theta_{j})$, $1\leq i<j\leq n$, we find
$\Hom_{\SO_{2n}}(V_{2\theta_{1}}\otimes\so_{2n},V_{2\theta_{1}})\cong\IC$
and therefore, by (ii) and (iii) of proposition \ref{cf rules}, we find
$(\Omega_{23})_{2\theta_{1},2\theta_{1}}=-n$, $(\Omega_{23})_{0,0}=0$
and, by (iv), $(\Omega_{23})_{2\theta_{1},0}=(\Omega_{23})_{0,2\theta_{1}}=0$.
The remaining entries are simply found by requiring that
\begin{equation}
\det(\Omega_{23}-t)=(k-t)(-1-t)(-2(n-1)-k-t)
\end{equation}
The matrix giving $\Omega_{13}$ now follows from $\Omega_{13}=F\Omega_{23}F^{-1}$
and \eqref{eq:F matrix} \halmos\\

For later reference, we shall need the following immediate

\begin{corollary}\label{technical}
Let $\kappa=\ell+2(n-1)$ and $\Omega_{ij}$ be the matrices given by proposition
\ref{pr:explicit} where $1\leq k\leq \ell$. Then
\begin{enumerate}
\item The eigenvalues of $\kappa^{-1}\Omega_{12}$ do not differ by integers unless $\ell=2$.
For $\ell=2$, $\kappa^{-1}(\delta_{2\theta_{1}}-\delta_{0})=1$ and all other eigenvalues do
not differ by integers.
\item The eigenvalues of $\kappa^{-1}\Omega_{13}$ do not differ by integers unless $\ell=2k$.
For $\ell=2k$, $\kappa^{-1}(\wt\delta_{(k+1)\theta_{1}}-\wt\delta_{(k-1)\theta_{1}})=1$ and
all other eigenvalues do not differ by integers.
\end{enumerate}
\end{corollary}

\ssection{The braiding relations}\label{se:braiding relations}

\ssubsection{General strategy}\label{ss:strategy}

We consider in this section the level $\ell$ \KZ equations
\begin{equation}\label{eq:vector}
\frac{dF}{dz}=
\frac{1}{\kappa}
\Bigl(\frac{\Omega_{12}-\delta_{0}}{z}+\frac{\Omega_{23}}{z-1}\Bigr)F
\end{equation}
with values in
$(V_{k\theta_{1}}\otimes V_{k\theta_{1}}^{*}\otimes 
  V_{\theta_{1}} \otimes V_{\theta_{1}}^{*})^{\SO_{2n}}$,
where $\kappa=\ell+\half{C_{\g}}$, with $\half{C_{\g}}=2(n-1)$ the dual Coxeter number of
$\SO_{2n}$ and the matrices $\Omega_{ij}$ and $\delta_{0}$ are given by proposition
\ref{pr:explicit}.\\

We prove in \S \ref{ss:reduction} that the reduction of \eqref{eq:vector}
to a scalar, third order equation coincides with the \DF equation \eqref{DF} for a suitable
choice of the parameters $a,b,c,g$ and use this to derive the analytic continuation of
solutions of \eqref{eq:vector} from 0 to $\infty$ via proposition \ref{continuation}.
These calculations are used in \S \ref{ss:coefficients} to obtain the main result of
this chapter, namely that the structure constants governing the braiding of the primary
fields with charges $V_{\theta_{1}}$ and $V_{k\theta_{1}}$ do not vanish.
A number of technical difficulties which we describe below arise in this approach. We
circumvent them in \S \ref{ss:deformation}.\\

Unfortunately, the values of the parameters of the \DF equation \eqref{DF} for which it
coincides with the scalar reduction of \eqref{eq:vector} lie outside the range
\eqref{range 1}--\eqref{range 2} where the solutions of \eqref{DF} may be expressed by
means of generalised Euler integrals. A further difficulty arises in
identifying specific solutions of \eqref{eq:vector} with those of \eqref{DF}. This 
needs to be done with great care for we are ultimately concerned with the analytic
continuation of a given solution of \eqref{eq:vector}, namely the reduced four--point
function of the product of primary fields with charges $V_{\theta_{1}}$ and $V_{k\theta_{1}}$.
For generic values of the level $\ell\in\IN$, the identification is easily obtained by
matching the monodromic behaviour of the solutions about a given singular point. However,
when $\ell\in\{2,2k\}$, corollary \ref{technical} implies that the monodromy operator
corresponding to \eqref{eq:vector} has degenerate eigenvalues at $0$ or $\infty$. For these
values of $\ell$, the solutions of \eqref{eq:vector} are not uniquely characterised by their
monodromy about these singular points and the identification becomes ambiguous.\\

We shall circumvent both of the above difficulties in the following way. We prove
in \S \ref{ss:deformation} that four--point functions at a given level $\ell=\ell_{0}$
may be analytically continued in $\ell\in(\IC\setminus\IR)\cup D$ trough solutions of the
\KZ equations where $D$ is a disc centred at $\ell_{0}$, so that the braiding coefficients
depend holomorphically on $\ell$
\footnote{As will become apparent in \S \ref{ss:deformation} this is not true of general
solutions of the \KZ equations.}.
When $\ell$ is not real, the eigenvalues of the residue matrices
of the \KZ equations do not differ by integers and the identification of the continued
four--point functions with solutions of the \DF equation is unambiguous.
Moreover, if $\ell$ has a large negative real part, the parameters of the \DF equation
fall into the good range \eqref{range 1}--\eqref{range 2} and the analytic continuation
from $0$ to $\infty$ in the $z$ variable may therefore be computed. This yields the braiding
coefficients for complex $\ell$. The required values are obtained by evaluation at
$\ell=\ell_{0}$.

\ssubsection{Holomorphic dependence of four--point functions on the level}\label{ss:deformation}

Let $G$ be a compact, connected, simple Lie group with dual Coxeter number $\half{C_{\g}}$.

\begin{lemma}
Let $\phi_{3}$, $\phi_{2}$ be primary fields at level $\ell_{0}$ with vertices
$\vertex{V_{3}}{V_{4}}{U}$, $\vertex{V_{2}}{U}{V_{1}}$ and
$F=\sum_{n\geq 0}a_{n}z^{n}$ the formal power series expansion of the
corresponding reduced four--point function. If $\Omega_{12}(a,b)$ is the
subspace of $V=(V_{4}\otimes V_{3}^{*}\otimes V_{2}^{*}\otimes V_{1}^{*})^{G}$
corresponding to the $\Omega_{12}$-eigenvalues $\delta\in(a,b)$, then for $n\geq 1$
\begin{equation}\label{spectral}
a_{n}\in\Omega_{12}(-\infty,\delta_{U}+n(\ell_{0}+\half{C_{\g}}))
\end{equation}
\end{lemma}
\proof
By proposition \ref{cf rules}, $\Omega_{12}$ acts diagonally in the basis of
$V\cong\Hom_{G}(V_{1}\otimes V_{2},V_{3}^{*}\otimes V_{4})$ given by the elements
$\psi_{W}\otimes\phi_{W}\in
 \Hom_{G}(W,V_{3}^{*}\otimes V_{4})\otimes
 \Hom_{G}(V_{1}\otimes V_{2},W)$,
with corresponding eigenvalue $\delta_{W}=\half{1}(C_{W}-C_{1}-C_{2})$.
Since $a_{n}=\phi_{3}(\cdot,n)\phi_{2}(\cdot,-n)$, we find
\begin{equation}
(a_{n},\psi_{W}\otimes\phi_{W})=
\tr_{V_{1}\otimes V_{2}}(\phi_{W}^{*}\psi_{W}^{*}a_{n})=
\tr_{W}(\psi_{W}^{*}\phi_{3}(\cdot,n)\phi_{2}(\cdot,-n)\phi_{W}^{*})
\end{equation}
Now, $\phi_{2}(\cdot,-n)\phi_{W}^{*}$ is a $G$-intertwiner $W\rightarrow\H_{U}(n)$ and
is therefore zero unless $W\subset\H_{U}(n)$. By proposition 11.4 b) of \cite{Ka1}, the
latter condition implies that $C_{W}-C_{U}<2n(\ell_{0}+\half{C_{\g}})$ and therefore
$(a_{n},\psi_{W}\otimes\phi_{W})=0$ if
$\delta_{W}-\delta_{U}=\half{1}(C_{W}-C_{U})\geq n(\ell_{0}+\half{C_{\g}})$ \halmos\\

\remark When expressed in terms of the power series expansion $F=\sum_{n}a_{n}z^{n}$
of a reduced four point-function with intermediate representation $U$, the \KZ equation
satisfied by $F$ becomes the recurrence relation
\begin{equation}
(\Omega_{12}-\delta_{U}-n(\ell_{0}+\half{C_{\g}}))a_{n}=\Omega_{23}\sum_{m=0}^{n-1}a_{m}
\end{equation}
A priori, this does not determine $F$ completely from $a_{0}$ since
$(\Omega_{12}-\delta_{U}-n(\ell_{0}+\half{C_{\g}}))$ may fail to be invertible for some $n$.
The above corollary however shows that this is not the case since
$a_{n}\in\Ker(\Omega_{12}-(\delta_{U}+n\kappa))^{\perp}$.

\begin{proposition}\label{deform l}
Let
$\phi_{2}:\Hfin_{V_{1}}\otimes V_{2}[z,z^{-1}]\rightarrow\Hfin_{U}$,
$\phi_{3}:\Hfin_{U}\otimes V_{3}[z,z^{-1}]\rightarrow\Hfin_{V_{4}}$
be primary fields at level $\ell_{0}\in\IN$ with initial terms $\Phi_{2},\Phi_{3}$
and denote by $F_{\ell_{0}}$ the corresponding reduced four--point function with values
in $V=(V_{4}\otimes V_{3}^{*}\otimes V_{2}^{*}\otimes V_{1}^{*})^{G}$.
If the eigenvalues of $(\ell_{0}+\half{C_{\g}})^{-1}\Omega_{12}$ on $V$ do not differ
by integers larger than $1$, there exists an $\epsilon>0$ such that
\begin{enumerate}
\item For any
$\ell\in\D=
 \{\lambda|\Im\lambda\neq 0\}\cup\{\lambda|\medspace|\lambda-\ell_{0}|<\epsilon\}$
with $\ell\neq\ell_{0}$, the \KZ equation
\begin{equation} \label{deform KZ}
(\ell+\half{C_{\g}})\frac{dF_{\ell}}{dz}=
\Bigl(\frac{\Omega_{12}-\delta_{U}}{z}+\frac{\Omega_{23}}{z-1}\Bigr)F_{\ell}
\end{equation}
possesses a unique formal power series solution $F_{\ell}=\sum_{n\geq 0}a_{n}(\ell)z^{n}$
with $a_{0}(\ell)=\Phi_{3}\Phi_{2}$.
\item $F_{\ell}$ is holomorphic on $\IC\setminus[1,\infty)$ and the assignement
$(\ell,z)\rightarrow F_{\ell}(z)$ defines a holomorphic function on
$\D\times(\IC\setminus[1,\infty))$.
\end{enumerate}
\end{proposition}
\proof For any $\ell\in\IC$, a formal power series
$F_{\ell}=\sum_{n\geq 0}a_{n}(\ell)z^{n}$ satisfies \eqref{deform KZ} if, and only if
\begin{align}
 \Omega_{12}a_{0}(\ell)&=\delta_{U}a_{0}(\ell)\\
(\Omega_{12}-\delta_{U}-n(\ell+\half{C_{\g}}))a_{n}(\ell)&=
 \Omega_{23}\sum_{m=0}^{n-1}a_{m}(\ell)
 \label{deform recurrence}
\end{align}
When $\ell=\ell_{0}$ and $F_{\ell_{0}}$ is the reduced four--point function of the product
$\phi_{3}\phi_{2}$ then $a_{0}(\ell_{0})=\Phi_{3}\Phi_{2}$ and 
\eqref{spectral} yields
$a_{1}(\ell_{0})\in\Omega_{12}(-\infty,\delta_{U}+\ell_{0}+\half{C_{\g}})$. Thus, by
\eqref{deform recurrence} with $n=1$
\begin{equation}\label{spectral estimate}
\Omega_{23}\Phi_{3}\Phi_{2}\in\Omega_{12}(-\infty,\delta_{U}+\ell_{0}+\half{C_{\g}})
\end{equation}

Let now $\epsilon>0$ be such that the eigenvalues of $(\ell+\half{C_{\g}})^{-1}\Omega_{12}$
do not differ by positive integers whenever $0<|\ell-\ell_{0}|<\epsilon$ and set
$\D=(\IC\setminus\IR)\cup\{\ell|\thinspace|\ell-\ell_{0}|<\epsilon\}$. Then, for any
$\ell\in\D\setminus\{\ell_{0}\}$, the matrices $(\Omega_{12}-\delta_{U}-n(\ell+\half{C_{\g}}))$
are invertible and therefore \eqref{deform KZ} possesses a unique formal power series
solution $F_{\ell}=\sum_{n\geq 0}a_{n}(\ell)z^{n}$ with $a_{0}(\ell)=\Phi_{3}\Phi_{2}$. By
\eqref{deform recurrence}, the
$a_{n}(\ell)$ are rational functions of $\ell$ and are holomorphic on $\D\setminus\{\ell_{0}\}$.
We claim that $a_{n}(\ell)\rightarrow a_{n}(\ell_{0})$ as $\ell\rightarrow\ell_{0}$ so that
$\ell_{0}$ is a removable singularity of the $a_{n}(\ell)$.
By our assumption on the eigenvalues of $\Omega_{12}$ we need only check this for $n=1$.
In this case, the result follows from \eqref{spectral estimate} since
$a_{1}(\ell)=(\Omega_{12}-\delta_{U}-(\ell+\half{C_{\g}}))^{-1}\Omega_{23}\Phi_{3}\Phi_{2}$
converges in $\Ker(\Omega_{12}-\delta_{U}-(\ell_{0}+\half{C_{\g}}))^{\perp}$ to the unique
solution of $(\Omega_{12}-\delta_{U}-(\ell_{0}+\half{C_{\g}}))x=\Omega_{23}\Phi_{3}\Phi_{2}$
and therefore to $a_{1}(\ell_{0})$.\\

We claim now that $F_{\ell}(z)$
is holomorphic on $\D\times\{z|\thinspace|z|<1\}$. Let $\Lambda\subset\D$ be a
compact set and $M=\sup_{\ell\in\Lambda}\|a_{1}(\ell)\|+\|\Phi_{3}\Phi_{2}\|$. Then, for
$\ell\in\Lambda$ and $n\geq 2$, the norms
$n\|(\Omega_{12}-\delta_{U}-n(\ell+\half{C_{\g}}))^{-1}\Omega_{23}\|$ are uniformly bounded by
some constant $C$ and therefore, by \eqref{deform recurrence}, the quantities
$\sup_{\ell\in\Lambda}\|a_{n}(\ell)\|$, $n\geq 1$ are bounded by the solution of the recurrence 
relation
\begin{align}
\alpha_{1}&=M\\
n\alpha_{n}&=C\sum_{m=1}^{n-1}a_{m}
\end{align}
By uniqueness, these are the coefficients of the power series expansion of $M(1-z)^{-C}$
and it follows that $F_{\ell}(z)=\sum_{\n\geq 0}a_{n}(\ell)z^{n}$ is absolutely
convergent, and therefore holomorphic on $\Lambda\times\{z|\medspace|z|<1\}$.
To conclude, notice that the value of $F_{\ell}(z)$ on $|z|\geq 1$ are given by parallel
transport from those in $|z|<1$ with respect to the connection \eqref{deform KZ}. Since
the latter depends analytically on $\ell$ we deduce that $F_{\ell}(z)$ is holomorphic on
$\D\times(\IC\setminus[1,\infty))$ as claimed \halmos\\

\remark The above proposition should hold without any assumptions on the eigenvalues of
$\Omega_{12}$ and would then imply that braiding coefficients may always be computed for
generic values of the level. At present, we can only prove this for $G=\SU_{2}$.

\ssubsection{The reduction of the KZ equation to the DF equation}\label{ss:reduction}

\begin{proposition}\label{KZ to DF}
Let $\kappa\in\IC^{*}$ and consider the \KZ equation with values in
$W=(V_{k\theta_{1}}\otimes V_{k\theta_{1}}^{*}\otimes 
    V_{ \theta_{1}}\otimes V_{ \theta_{1}}^{*})^{\SO_{2n}}$ given by
\begin{equation}\label{KZ 2}
\frac{df}{dz}=
\frac{1}{\kappa}
\Bigl(\frac{\Omega_{12}-\delta_{0}}{z}+\frac{\Omega_{23}}{z-1}\Bigr)f
\end{equation}
where $\delta_{0}$ and the matrices $\Omega_{ij}$ are given by proposition
\ref{pr:explicit}. If $e_{0}$ is the eigenvector of
$\Omega_{12}$ corresponding to the eigenvalue $\delta_{0}$, the map
$f\rightarrow R(f)=(f,e_{0})(z-1)^{-k/\kappa}$ is an isomorphism of
the space of solutions of \eqref{KZ 2} onto the solutions of the \DF
equation \eqref{DF} with parameters
\begin{xalignat}{4}\label{fitting}
a&=\frac{2(n-1)+k}{\kappa}&
b&=-\frac{\kappa+1}{\kappa}&
c&=-\frac{k}{\kappa}&
g&=-2\frac{n-2}{\kappa}
\end{xalignat}
\end{proposition}
\proof We shall use the basis of $W$ where the matrices $\Omega_{ij}$ have
the form given by proposition \ref{pr:explicit} so that
$e_{0}=\begin{pmatrix}0\\0\\1\end{pmatrix}$.
Let $f$ be a solution of \eqref{KZ 2}. By setting $g=(z-1)^{-k/\kappa}f$ and
rescaling the entries of $g$, the equations may be written as
\begin{equation}\label{system}
\begin{split}
\kappa\frac{dg}{dz}
&=
\left(\begin{array}{ccc}
2n&0&0\\0&2(n-1)&0\\0&0&0
\end{array}\right)
\frac{g}{z}\\
&+
\left(\begin{array}{ccc}
-n-k  	&\frac{(n-1)}{n}(k^{2}+2k(n-1)+n(n-2))	&0\\
1	&-n-k+1					&\frac{k}{n}(k+2(n-1))\\
0	&1					&-k
\end{array}\right)
\frac{g}{z-1}
\end{split}
\end{equation}
Let $g=\begin{pmatrix}u\\ v\\ w\end{pmatrix}$. Using the last component of
\eqref{system} to express $v$ in terms of $w$ and the second to express $u$
in terms of $w$, one finds
\begin{align}
u &=
\frac{k(n-1)}{z}(2-\frac{n-k+2}{n}z)w+
\frac{\kappa(z-1)}{z}(2(n-1)+(\kappa-n+2k+1)z)w'+
\kappa^{2}(z-1)^{2}w''
\label{elim u}\\
v &= kw+\kappa(z-1)w'
\label{elim v}
\end{align}
Eliminating $u$ and $v$ from the first component of \eqref{system} yields a
third order ODE for $w=R(f)$ of the form \eqref{DF} with
\begin{xalignat}{2}
K_{1} &= \frac{2n+3(k+\kappa)-1}{\kappa}&
K_{2} &= -2\frac{2n-1}{\kappa}
\end{xalignat}
\begin{xalignat}{2}
L_{1} &= \frac{(1+k+\kappa)(\kappa+2(n+k-1))}{\kappa^{2}}&
L_{2} &= 2\frac{(n-1)(2n+\kappa)}{\kappa^{2}}
\end{xalignat}
\begin{equation}
L_{3} =-2\frac{2n^{2}+4kn+3n\kappa-2(n+k+\kappa)}{\kappa^{2}}
\end{equation}
\begin{xalignat}{2}
M_{1} &= -2\frac{k(n-1)(\kappa+2(n+k-1))}{\kappa^{3}}&
M_{2} &=  2\frac{k(n-1)(2n+\kappa)}{\kappa^{3}}
\end{xalignat}
which is indeed the Dotsenko--Fateev equation with parameters given by
\eqref{fitting}. Since $w=0$ implies $u=v=0$ by \eqref{elim u}--\eqref{elim v},
$R$ is injective and therefore surjective \halmos\\

The following elementary lemma justifies our use of complex values of $\kappa$
or equivalently of $\ell=\kappa-2(n-1)$. Indeed, for $\kappa\in\IR_{+}$,
$b=-(\kappa+1)\kappa^{-1}\in(-\infty,-1)$ and the parameters of the \DF equation
given by \eqref{fitting} lie outside the range \eqref{range 1}--\eqref{range 2} where
we have a good description of its solutions.

\begin{lemma}
The parameters $a,b,c,g$ given by \eqref{fitting} lie in the range
\begin{xalignat}{2}
\Re a,\Re b,\Re c,\Re g &> -1&
\Re (2a+g),\Re (2b+g),\Re (2c+g) &> -2
\\
\Re (a+b+c+g) &< -1&
\Re (2a+2b+2c+g) &< -2
\end{xalignat}
if, and only if
\begin{xalignat}{3}\label{kappa range}
\Re\kappa&<0&
&\text{and}&
|\kappa+\half{1}(2(n-1)+k)|^{2}&>\frac{1}{4}(2(n-1)+k)^{2}
\end{xalignat}
\end{lemma}

The next lemma describes the monodromic behaviour at $0$ or $\infty$ of the scalar
reduction $R(\cdot)$ of specific solutions of the \KZ equations \eqref{KZ 2}. When 
$\kappa$ satisfies \eqref{kappa range}, this will allow us to identify them with
the solutions of the \DF equation given in proposition \ref{DF solutions}.

\begin{lemma}\label{KZ identikit}
Fix a basis of eigenvectors of each of the matrices $\Omega_{12}$ and $\Omega_{13}$
given by proposition \ref{pr:explicit}. Denote the solutions of \eqref{KZ 2} at $0$ and
$\infty$ with initial terms equal to these eigenvectors by
$f_{2\theta_{1}},f_{\theta_{1}+\theta_{2}},f_{0}$
and
$g_{(k+1)\theta_{1}},g_{k\theta_{1}+\theta_{2}},g_{(k-1)\theta_{1}}$
respectively where the labels refer to those of the corresponding eigenvalues.
Then, up to multiplicative constants independent of $\kappa$
\begin{xalignat}{2}
&&
R(g_{(k+1)\theta_{1}})&=
\Bigl(\frac{1}{z}\Bigr)^{2k/\kappa}(1+z^{-1}O(z^{-1}))\\
R(f_{0})&=e^{-i\pi k/\kappa}(1+zO(z))&
R(g_{k\theta_{1}+\theta_{2}})&=
\Bigl(\frac{1}{z}\Bigr)^{(k-1)/\kappa}(1+z^{-1}O(z^{-1}))\\
&&
R(g_{(k-1)\theta_{1}})&=
\Bigl(\frac{1}{z}\Bigr)^{-2(n-1)/\kappa}(1+z^{-1}O(z^{-1}))
\end{xalignat}
\end{lemma}
\proof We use the basis of
$W=(V_{k\theta_{1}}\otimes V_{k\theta_{1}}^{*}\otimes 
    V_{ \theta_{1}}\otimes V_{ \theta_{1}}^{*})^{\SO_{2n}}$
given by proposition \ref{pr:explicit} so that the eigenvector $e_{0}$ of $\Omega_{12}$
corresponding to the eigenvalue $\delta_{0}$ is $\begin{pmatrix}0\\0\\1\end{pmatrix}$.
Since the Taylor series of $f_{0}$ at $0$ is $e_{0}+\sum_{n\geq 1}e_{n}z^{n}$, the
lemma follows for $R(f_{0})$.
The Taylor series of $g_{\mu}$ at $z=\infty$ is
$\Bigl(\frac{1}{z}\Bigr)^{\wt\delta_{\mu}}\sum_{n\geq 0}\xi_{n,\mu}z^{-n}$ where the
$\wt\delta_{\mu}$ and $\xi_{0,\mu}$ are the eigenvalues and corresponding eigenvectors
of the residue matrix of \eqref{KZ 2} at infinity. By (\ref{ch:algebraic fields}.\ref{trick}),
the latter is $-(\Omega_{12}-\delta_{0}+\Omega_{23})=\Omega_{13}$ and the $\wt\delta_{\mu}$
are therefore given by proposition \ref{pr:explicit}.
Rescaling the basis of $W$ so as to have
\begin{equation}
\Omega_{13}=
\left(\begin{array}{rcc}
-n   &\frac{(n-1)}{n}(k^{2}+2k(n-1)+n(n-2))	&0\\
1    &-n+1					&\frac{k}{n}(k+2(n-1))\\
0    &1						&0
\end{array}\right)
\end{equation}
the corresponding eigenvectors are
\begin{xalignat}{2}
v_{(k+1)\theta_{1}}&=
\begin{pmatrix}
k(n-1)(k+n-2)\\
-kn\\
n
\end{pmatrix}&
v_{k\theta_{1}+\theta_{2}}&=
\begin{pmatrix}
k^{2}+2k(n-1)+n(n-2)\\
-n\\
-n
\end{pmatrix}
\end{xalignat}
\begin{equation}
v_{(k-1)\theta_{1}}=
\begin{pmatrix}
(n-1)(n+k)(k+2(n-1))\\
n(k+2(n-1))\\
n
\end{pmatrix}
\end{equation}
Thus, up to a multiplicative constant
\begin{equation}
R(g_{\mu})=
(z-1)^{-k/\kappa}\Bigl(\frac{1}{z}\Bigr)^{\wt\delta_{\mu}}(1+z^{-1}O(z^{-1}))=
\Bigl(\frac{1}{z}\Bigr)^{\wt\delta_{\mu}+k/\kappa}(1+z^{-1}O(z^{-1}))
\end{equation}
as claimed \halmos

\ssubsection{Non--vanishing of the braiding coefficients}\label{ss:coefficients}

\begin{theorem}\label{th:braiding identity}
At any level $\ell_{0}\in\IN$, the following braiding identities hold
\begin{align}
\field{k\theta_{1}}{k\theta_{1}}{0}(w)
\field{ \theta_{1}}{0}{ \theta_{1}}(z)
&=
\sum_{\mu\in\{(k-1)\theta_{1},k\theta_{1}+\theta_{2},(k+1)\theta_{1}\}}
\beta_{k\theta_{1},\mu}
\field{ \theta_{1}}{k\theta_{1}}{\mu}(z)
\field{k\theta_{1}}{\mu}{ \theta_{1}}(w)
\label{braiding 1}\\
\intertext{if $k=1\ldots\ell_{0}-1$ and}
\field{\ell_{0}\theta_{1}}{\ell_{0}\theta_{1}}{0}(w)
\field{      \theta_{1}}{0}{       \theta_{1}}(z)
&=
\beta_{\ell_{0}\theta_{1},(\ell_{0}-1)\theta_{1}}
\field{    \theta_{1}}{\ell_{0}\theta_{1}}{(\ell_{0}-1)\theta_{1}}(z)
\field{\ell_{0}\theta_{1}}{(\ell_{0}-1)\theta_{1}}{    \theta_{1}}(w)
\label{braiding 2}
\end{align}
where the braiding coefficients $\beta_{\lambda,\mu}$ are non--zero.
\end{theorem}
\proof
By the tensor product rule \eqref{tensor one} and theorem
\ref{ch:algebraic fields}.\ref{th:existence of braiding}, \eqref{braiding 1} holds
for $k=1\ldots\ell_{0}-1$ for some braiding coefficients $\beta_{\lambda,\mu}$. When
$k=\ell_{0}$, the sum on the right hand--side of \eqref{braiding 1} is restricted to
$\mu=(k-1)\theta_{1}=(\ell_{0}-1)\theta_{1}$ since the irreducible representations
of $\SO_{2n}$ with highest weight given by the other values of $\mu$ are not admissible
at level $\ell_{0}$. Thus, \eqref{braiding 2} holds for some
$\beta_{\ell_{0}\theta_{1},(\ell_{0}-1)\theta_{1}}$ whose vanishing would imply that of the
four--point function of the product
$\field{k\theta_{1}}{k\theta_{1}}{0}(w)\field{ \theta_{1}}{0}{ \theta_{1}}(z)$,
a contradiction. This settles \eqref{braiding 2}.
Let now $f_{0}$ be the reduced four--point function corresponding to the left hand--side
of \eqref{braiding 1} so that $f_{0}$ satisfies \eqref{KZ 2} and
the leading term of its Taylor expansion at $0$ is the $\Omega_{12}$ eigenvector
$e_{0}$ corresponding to the eigenvalue $\delta_{0}$. Let $g_{\mu}$ be the four--point
functions corresponding to the right hand--sides of \eqref{braiding 1}
so that the initial terms of their expansion at $\infty$ are the eigenvectors of
$\Omega_{13}$.
By corollary \ref{technical}, $\Omega_{12}$ and $\Omega_{13}$ satisfy the assumptions
of proposition \ref{deform l} and we may therefore analytically continue $f_{0}$
and the $g_{\mu}$ through solutions of \eqref{KZ 2} for any $\kappa=\ell+2(n-1)$
where $\ell\in\D=(\IC\setminus\IR)\cup D$ and $D$ is a disc centred at $\ell_{0}$.
We shall abusively denote the continuation of the four--point functions by the same
symbol.\\

Take $\kappa=\ell+2(n-1)$, $\ell\in\D$ in the range \eqref{kappa range} such that
$\Im\kappa\neq 0$. Then, the characteristic exponents of the \DF equation given
by the monodromic behaviour of the solutions of proposition \ref{DF solutions}
do not differ by integers when the parameters $a,b,c,g$ are bound by \eqref{fitting}
and it follows that the solutions of the \DF equation are, up to a constant factor
uniquely determined by their monodromy about $0$ or $\infty$.
Thus, using \eqref{fitting} to compare the exponents of $R(f_{0})$ and $R(g_{\mu})$ given
by lemma \ref{KZ identikit} with those of the solutions $I_{z_{i},j}$ given by proposition
\ref{DF solutions}, we find
\begin{xalignat}{2}
&&
R(g_{(k+1)\theta_{1}})&=\rho_{\infty,1}^{-1}I_{\infty,1}\\
R(f_{0})&=e^{-i\pi k/\kappa}\rho_{0,1}^{-1}I_{0,1}&
R(g_{k\theta_{1}+\theta_{2}})&=\rho_{\infty,2}^{-1}I_{\infty,2}\\
&&
R(g_{(k-1)\theta_{1}})&=\rho_{\infty,3}^{-1}I_{\infty,3}
\end{xalignat}

and therefore, by proposition \ref{continuation}
\begin{equation}\label{almost}
\begin{split}
R(f_{0})
&=
e^{-i\pi k/\kappa}\rho_{0,1}^{-1}
\rho_{\infty,1}
\frac{\esse(a)\esse(a+\half{g})}{\esse(a+b+\half{g})\esse(a+b+g)}
R(g_{(k+1)\theta_{1}})
\\
&+
e^{-i\pi k/\kappa}\rho_{0,1}^{-1}
\rho_{\infty,2}
2e^{-\pi i(a+c+\half{g})}\ci(\half{g})
\frac{\esse(a)\esse(c)}{\esse(a+b)\esse(a+b+g)}
R(g_{k\theta_{1}+\theta_{2}})
\\
&+
e^{-i\pi k/\kappa}\rho_{0,1}^{-1}
\rho_{\infty,3}
\frac{\esse(c)\esse(c+\half{g})}{\esse(a+b+\half{g})\esse(a+b+\half{g})}
R(g_{(k-1)\theta_{1}})
\end{split}
\end{equation}
where $a,b,c,g$ are given by \eqref{fitting} and the $\rho_{z_{i},j}$ by proposition
\ref{DF solutions}.
Since $R$ is injective the same linear relation binds $f_{0}$ and the $g_{\mu}$.
Continuing back to $\kappa=\ell_{0}+2(n-1)$ therefore yields
\begin{equation}
f_{0}=
\beta_{(k+1)\theta_{1}}g_{(k+1)\theta_{1}}+
\beta_{k\theta_{1}+\theta_{2}}g_{k\theta_{1}+\theta_{2}}+
\beta_{(k-1)\theta_{1}}g_{(k-1)\theta_{1}}
\end{equation}
where the coefficients $\beta_{\mu}$ are obtained by evaluating those of \eqref{almost}
and, by direct inspection, do not vanish if $1\leq k\leq \ell_{0}-1$ \halmos \\

\remark The coefficients involved in the braiding of the primary fields with charge the
vector representation of $\SO_{2n}$ and its second exterior power $V_{\theta_{1}+\theta_{2}}$
may in principle be obtained by studying the \KZ equations with values in 
$(V_{\theta_{1}+\theta_{2}}\otimes V_{\theta_{1}+\theta_{2}}^{*}\otimes
  V_{\theta_{1}}           \otimes V_{\theta_{1}}^{*})^{\SO_{2n}}$.
A computation similar to the proof of proposition \ref{pr:explicit} shows in fact
that the corresponding operators $\Omega_{ij}$ are
\begin{xalignat}{2}\label{matrices 2}
\Omega_{12}&=
\left(\begin{array}{crc}
1&0&0\\0&-1&0\\0&0&-2n+1
\end{array}\right)&
\Omega_{23}&=
\left(\begin{array}{rcc}
-n   &x    &0\\
\wt x&-n+1 &y\\
0    &\wt y&0
\end{array}\right)
\end{xalignat}
where
\begin{xalignat}{3}
x\wt x&=\frac{(n+1)(n-2)^{2}}{n}&
&\text{and}&
y\wt y&=4\frac{n-1}{n}
\end{xalignat}
In spite of the similarity of \eqref{matrices 2} and \eqref{matrices}, the corresponding
\KZ equations, when reduced to a third order scalar equation are not of the \DF form
so that no explicit solutions are readily available.\\

\remark The computations of this chapter apply, with very little change to the group
$G=\SO_{2n+1}$ and may be used to show that the braiding relations analogous to those
of theorem \ref{th:braiding identity} only involve non--zero coefficients.\\

\remark The solutions \eqref{euler} were discovered by Dotsenko and Fateev in the
context of minimal models \cite{DF}. They were subsequently generalised
by Schechtman--Varchenko \cite{SV} and Feigin--Frenkel \cite{EdF} who used them
systematically to give solutions of all KZ equations.
It is interesting to note however that their generalised hypergeometric solutions, which
first appeared in the work of Aomoto and Gelfand \cite{Ao,GKZ}, use Euler--like contour
integrals with a number of integration variables growing linearly in $n$, even in the
specific cases considered in this chapter. This makes them intractable for computational
purposes.


\part{Fusion of positive energy representations}




\newcommand {\RRR}{\mathfrak R}
\newcommand {\Nbox}{N_{\Box}}

\newcommand {\pull}[1]{\pi^{*}_{#1}U(\H_{#1})}


\chapter{Connes fusion of positive energy representations of $L\Spin_{2n}$}
\label{ch:connes fusion}

In this final chapter, we define a tensor product operation $\fuse$ or {\it Connes
fusion} on the category $\pl$ of positive energy representations of $L\Spin_{2n}$
at a fixed level $\ell$ and study the resulting algebraic structure on $\pl$.\\

In section \ref{se:definition} we give the definition and elementary properties
of Connes fusion. In section \ref{se:level 1 fusion ring}, we use the action of
the centre $Z(\Spin_{2n})$ on $\pl$ via conjugation by discontinuous loops to
compute the fusion with representations lying on the orbit of the vacuum $\H_{0}$.
As a  simple corollary, we show that the level 1 fusion ring of $L\Spin_{2n}$ is
isomorphic to the group algebra of $Z(\Spin_{2n})$.
In section \ref{se:bounded braiding}, we derive braiding identities for bounded
intertwiners for the local loop groups $L_{I}\Spin_{2n}$ from the ones satisfied
by smeared primary fields. These are used in section \ref{se:upper bound} to give
an upper bound for the fusion with the vector representation $\Hbox$
in terms of the {\it Verlinde rules}.
The rest of the chapter is devoted to showing that this bound is actually attained.
In section \ref{se:fusion box/kbox}, we prove this for the fusion of  $\Hbox$ with
its symmetric powers by using the braiding computations
of chapter \ref{ch:box/kbox braiding}.
Section \ref{se:subring} initiates the study of the ring $\RR$ generated by the
irreducible summands of the iterated fusion products of $\Hbox$.
We prove that $\RR$ is commutative by using a {\it braiding operator} $B$ which
gives an isomorphism $\H_{1}\fuse\H_{2}\cong\H_{2}\fuse\H_{1}$ for any
$\H_{1},\H_{2}\in\RR$. 
In section \ref{se:eigenvalues of braiding} we compute the eigenvalues of the
braiding operator acting on the fusion of $\Hbox$ with itself.
Section \ref{se:DHR} describes some important algebraic properties of $\RR$ -- most
notably the existence of a {\it quantum dimension function} compatible with fusion --
which are derived from the Doplicher--Haag--Roberts theory of superselection sectors.
In section \ref{se:qdim of box}, we compute the quantum dimension of $\Hbox$
by using some arguments of Wenzl and the fusion rules obtained in section
\ref{se:fusion box/kbox}.
Section \ref{se:Verlinde} gives an alternative computation of this 
dimension under the assumption that the Verlinde rules hold. Reassuringly,
we find the value obtained in section \ref{se:qdim of box}.
Section \ref{se:main results} contains the main results of this thesis.
By using a Perron--Frobenius argument based on the coincidence of the two
computations of the quantum dimension of $\Hbox$, we show
that its fusion with the positive energy
representations of $L\Spin_{2n}$ whose lowest energy subspace is a single--valued
$\SO_{2n}$--module is given by the Verlinde rules. This implies in particular
that these representations are closed under fusion and form a commutative and
associative ring. Finally, using the action of discontinuous loops on $\pl$,
we extend the previous results to all positive energy representations of
$L\Spin_{2n}$ when the level is odd.

\ssection{Definition of Connes fusion}\label{se:definition}

This section follows \cite[\S 30]{Wa3} and rests on the results describing the von
Neumann algebras generated by local loop groups in positive energy representations
obtained in chapter \ref{ch:loc loops}.
Fix $\ell\in\IN$ and consider the set $\pl$ of positive energy representations of
$LG$ at level $\ell$. Here and in what follows, $G=\Spin_{2n}$, $n\geq 3$.\\

Let $(\pi_{0},\H_{0})\in\pl$ be the vacuum representation and denote by $\LG=\pull{0}$
and $\lig$ the corresponding central extension of $LG$ and its restriction to a local
loop group $\Lig$.
By proposition \ref{ch:classification}.\ref{classify extension}, $\LG$ is isomorphic
to $\pull{i}$ for any $(\pi_{i},\H_{i})\in\pl$ and therefore acts unitarily on $\H_{i}$.
We denote the corresponding representation by the same symbol.
The restriction of $\pi_{i}$ to $\lig$ has the simple form
$\pi_{i}(\gamma)=U\pi_{0}(\gamma)U^{*}$ where $U:\H_{0}\rightarrow\H_{i}$ is a unitary
equivalence of $\Lig$--modules provided by local equivalence. Indeed, both $\pi_{i}$
and conjugation by $U$ yield isomorphisms $\lig\cong\left.\pull{i}\right|_{\Lig}$ and
therefore differ by a character of $\Lig$ which, by lemma
\ref{ch:loc loops}.\ref{local is perfect} is necessarily trivial. It follows that the
restriction of $\pi_{i}$ to $\Lig$ extends to the canonical spatial isomorphism
$\pi_{0}(\Lig)''\cong\pi_{i}(\Lig)''$. In particular, if $x\in\Hom_{\lig}(\H_{0},\H_{i})$,
then $x a=\pi_{i}(a) x$ for any $a\in\pi_{0}(\Lig)''$.\\

Let now $\H_{i},\H_{j}\in\pl$. By locality, we may regard each as a bimodule over
the pair $(\lig,\licg)$ and form the intertwiner spaces
\begin{xalignat}{3}\label{eq:hom spaces}
\XX_{i}&=\Hom_{\licg}(\H_{0},\H_{i})&
&\text{and}&
\YY_{j}&=\Hom_{\lig} (\H_{0},\H_{j})
\end{xalignat}
where $\H_{0}$ is the vacuum representation at level $\ell$ with highest weight
vector $\Omega$. These are $\lig$ and $\licg$--modules respectively.

\begin{lemma}\label{density}
The maps $\XX_{i}\rightarrow\H_{i}$ and $\YY_{j}\rightarrow\H_{j}$ given by
$x\rightarrow x\Omega$ and $y\rightarrow y\Omega$ are embeddings with dense
image.
\end{lemma}
\proof
If $x\in\XX_{i}$ and $x\Omega=0$, then
$x\pi_{0}(\licg)''\Omega=\pi_{i}(\licg)''x\Omega=0$
and therefore, by the Reeh--Schlieder theorem, $x=0$. To prove the density of
the embedding, we may assume that $\H_{i}$ is irreducible. Pick a unitary
$u_{i}\in\Hom_{\licg}(\H_{0},\H_{i})$ using local equivalence. Then, by Haag
duality $\XX_{i}=u_{i}\pi_{0}(\licg)'=u_{i}\pi_{0}(\lig)''$ and therefore, by
the Reeh--Schlieder theorem
$\overline{\XX_{i}\Omega}=u_{i}\overline{\pi_{0}(\lig)''\Omega}=\H_{i}$. The
fact that $\YY_{j}$ embeds densely in $\H_{j}$ is proved identically \halmos\\

Consider the hermitian form $\<\cdot,\cdot\>$ on the algebraic tensor product
$\XX_{i}\otimes\YY_{j}$ defined by
\begin{equation}
\<x_{1}\otimes y_{1},x_{2}\otimes y_{2}\>=
(x_{2}^{*}x_{1}y_{2}^{*}y_{1}\Omega,\Omega)
\end{equation}
where the inner product is taken in $\H_{0}$.

\begin{lemma}
The form $\<\cdot,\cdot\>$ is positive semi--definite.
\end{lemma}
\proof 
Let $z=\sum_{p=1}^{n}x_{p}\otimes y_{p}\in\XX_{i}\otimes\YY_{j}$ so
that $\<z,z\>=\sum_{p,q}(x_{q}^{*}x_{p}y_{q}^{*}y_{p}\Omega,\Omega)$.
Let $\{e_{p}\}_{p=1}^{n}$ be an orthonormal basis of $\IC^{n}$ and set
$\Omega_{p}=\Omega\otimes e_{p}\in\H=\H_{0}\otimes\IC^{n}$. Consider
the non--negative $n\times n$ matrices $X,Y$ acting on $\H$ with entries
$x_{p}^{*}x_{q}\in\pi_{0}(\licg)'$ and $y_{q}^{*}y_{p}\in\pi_{0}(\lig)'$
respectively. By Haag duality, $X$ and $Y$ commute and therefore
$\<z,z\>=\sum_{p}(XY\Omega_{p},\Omega_{p})\geq 0$ \halmos\\

\definition
The {\it Connes fusion} $\H_{i}\fuse\H_{j}$ is the Hilbert space completion of
the quotient of $\XX_{i}\otimes\YY_{j}$ by the radical of $\<\cdot,\cdot\>$.
$\H_{i}\fuse\H_{j}$ supports commuting unitary actions of $\lig\times\licg$
given by
$(\gamma_{I},\gamma_{I^{c}})x\otimes y=
 \pi_{i}(\gamma_{I})x\otimes\pi_{j}(\gamma_{I^{c}})y$
which are easily seen to be strongly continuous.\\

\remark
The Connes fusion $\H_{i}\fuse\H_{j}$ of two irreducibles $\H_{i},\H_{j}\in\pl$
is {\it a priori} only an $(\lig,\licg)$--bimodule since the action of
$(\lig,\licg)$ need not necessarily extend to one of $\LG$, let alone
a positive energy one of $\LG\rtimes\rot$. One of our main results
(theorems \ref{th:closure of ordinary} and \ref{th:closure at odd}) states
that this action does extend, and therefore that $\H_{i}\fuse\H_{j}$ is a
positive energy representation of $LG$, when $\ell$ is odd or the lowest
energy subspaces of $\H_{i},\H_{j}$ are single--valued $\SO_{2n}$--modules.

\begin{lemma}
$\H_{0}\fuse\H_{i}\cong\H_{i}\cong\H_{i}\fuse\H_{0}$
\end{lemma}
\proof
The map $U:\XX_{0}\otimes\YY_{i}\rightarrow\H_{i}$, $x\otimes y\rightarrow yx\Omega$ is
readily seen to be an $\lig\times\licg$--equivariant isometry. By lemma \ref{density},
$\overline{U(\XX_{0}\otimes\YY_{i})}=
 \overline{\YY_{i}\XX_{0}\Omega}=
 \overline{\YY_{i}\H_{0}}=\H_{i}$
and $U$ therefore extends to a unitary $\H_{0}\fuse\H_{i}\rightarrow\H_{i}$.
Similarly, the unitary equivalence $\H_{i}\fuse\H_{0}\cong\H_{i}$ is given
by $y\otimes x\rightarrow yx\Omega$ \halmos

\begin{proposition}\label{pr:associativity}
Connes fusion is associative, that is
$\H_{i}\fuse(\H_{j}\fuse\H_{k})\cong(\H_{i}\fuse\H_{j})\fuse\H_{k}$.
\end{proposition}
\proof
We follow \cite{Lo}. Notice that for any $x_{p}\in\XX_{i}$ and $y_{q}\in\YY_{j}$,
\begin{equation}
\<x_{1}\otimes y_{1},x_{2}\otimes y_{2}\>=
(\pi_{i}(y_{2}^{*}y_{1})x_{1}\Omega,x_{2}\Omega)=
(\pi_{j}(x_{2}^{*}x_{1})y_{1}\Omega,y_{2}\Omega)
\end{equation}
where, by Haag duality $y_{2}^{*}y_{1}\in\pi_{0}(\lig)'=\pi_{0}(\licg)''$ and
$x_{2}^{*}x_{1}\in\pi_{0}(\lig)''$ so that, by local equivalence, they may be
represented on $\H_{i}$ and $\H_{j}$ respectively. It follows from lemma
\ref{density} that $\H_{i}\fuse\H_{j}$ may equivalently be defined as the
Hilbert space completion of $\H_{i}\otimes\YY_{j}$ or $\XX_{i}\otimes\H_{j}$
with respect to the forms
\begin{xalignat}{3}
(\xi_{1}\otimes y_{1},\xi_{2}\otimes y_{2})&=
(\pi_{i}(y_{2}^{*}y_{1})\xi_{1},\xi_{2})&
&\text{and}&
(x_{1}\otimes\eta_{1},x_{2}\otimes\eta_{2})&=
(\pi_{j}(x_{2}^{*}x_{1})\eta_{1},\eta_{2})
\end{xalignat}
respectively. More generally, $\H_{i}\fuse\H_{j}$ is the completion with respect
to these forms of $\K_{i}\otimes\YY_{j}$ or $\XX_{i}\otimes\K_{j}$, where the
$\K_{p}\subset\H_{p}$ are dense subspaces.
With these alternative descriptions of Connes fusion at hand, a unitary giving
the claimed isomorphism may be densely defined as the $\lig\times\licg$--equivariant
map
\begin{equation}
\XX_{i}\otimes(\H_{j}\otimes\YY_{k})\longrightarrow(\XX_{i}\otimes\H_{j})\otimes\YY_{k},
\medspace x\otimes(\eta\otimes z)\longrightarrow (x\otimes\eta)\otimes z
\end{equation}
\halmos\\

\remark
As stated, proposition \ref{pr:associativity} requires $\H_{j}\fuse\H_{k}$ and
$\H_{i}\fuse\H_{j}$ to be of positive energy and we shall use it in those cases
only. We note however that Connes fusion may be defined in the larger category
of bimodules over the III$_{1}$ factors $(\pi_{0}(\Lig)'',\pi_{0}(\Licg)'')$
where the proof of proposition \ref{pr:associativity} shows that it is an
associative operation \cite{Wa2}.

\ssection{Fusion with quasi--vacuum representations and the level 1 fusion ring}
\label{se:level 1 fusion ring}

Recall from section \ref{ch:classification}.\ref{se:disc loops} that the centre of
$G$ acts on the positive energy representations at level $\ell$ via conjugation by
discontinuous loops. If $\zeta\in\lzg$ and $(\pi,\H)\in\pl$, we denote the
corresponding conjugated representation by $(\zeta_{*}\pi,\zeta\H)$.
We call the $\zeta\H_{0}$ {\it quasi--vacuum} representations in view of the following

\begin{proposition}
Let $(\zeta_{*}\pi_{0},\zeta\H_{0})$ be a quasi--vacuum representation. Then,
Haag duality holds
\begin{equation}
\zeta_{*}\pi_{0}(\lig)'=\zeta_{*}\pi_{0}(\licg)''
\end{equation}
\end{proposition}
\proof This follows from Haag duality for the vacuum representation and the
fact that conjugation by $\zeta$ normalises the local loop groups \halmos

\begin{proposition}\label{fuse with quasi}
Let $\zeta\H_{0}$ be a quasi--vacuum representation, then
$\zeta\H_{0}\fuse\H_{i}\cong\zeta\H_{i}\cong\H_{i}\fuse\zeta\H_{0}$.
\end{proposition}
\proof
By multiplying $\zeta$ by a suitable element in $LG$, we may assume that it is
equal to 1 on $I^{c}$. Then,
$\zeta\XX_{0}=\Hom_{\licg}(\H_{0},\zeta\H_{0})=\pi_{0}(\licg)'=\pi_{0}(\lig)''$.
The map $\zeta\XX_{0}\otimes\YY_{i}\rightarrow\H_{i}$,
$x\otimes y\rightarrow yx\Omega$ is then easily seen to extend to a unitary
$\zeta\H_{0}\fuse\H_{i}\rightarrow\zeta\H_{i}$ intertwining $\lig\times\licg$.
The isomorphism $\H_{i}\fuse \zeta\H_{0}\cong\zeta\H_{i}$ follows in a similar
way by choosing $\zeta$ supported in $I^{c}$ \halmos\\

As a simple application of the above, we determine the level 1 fusion ring of
$L\Spin_{2n}$.

\begin{theorem}\label{th:level 1 ring}
The level 1 representations of $L\Spin_{2n}$ are closed under fusion. Moreover,
if $\H_{0}$ is the vacuum representation, then
\begin{equation}\label{eq:map}
z\longrightarrow z\H_{0}
\end{equation}
yields an isomorphism of the group algebra $\IC[Z(\Spin_{2n})]$ and the level
1 fusion ring of $L\Spin_{2n}$. Explicitly, if $\H_{0},\H_{v},\H_{s_{\pm}}$
are the irreducible level 1 representations, then
\begin{xalignat}{2}
\H_{v}\fuse\H_{v}&\cong\H_{0}&
\H_{v}\fuse\H_{s_{\pm}}&\cong\H_{s_{\mp}}\label{eq:1 fuse 1}\\[1.4 ex]
\H_{s_{\pm}}\fuse\H_{s_{\pm}}&\cong
	\left\{\begin{array}{ll}
	\H_{0}&\text{if $n$ even}\\
	\H_{v}&\text{if $n$ odd}
	\end{array}\right.&
\H_{s_{\pm}}\fuse\H_{s_{\mp}}&\cong
	\left\{\begin{array}{ll}
	\H_{v}&\text{if $n$ even}\\
	\H_{0}&\text{if $n$ odd}
	\end{array}\right.
\label{eq:1 fuse 2}
\end{xalignat}
\end{theorem}
\proof
The map \eqref{eq:map} is bijective since, by corollary
\ref{ch:classification}.\ref{co:level 1 Z}, $Z(\Spin_{2n})$ acts freely
and transitively on the irreducible level 1 representations. Moreover,
by proposition \ref{fuse with quasi},
\begin{equation}
z_{1}\H_{0}\fuse z_{2}\H_{0}\cong z_{1}(z_{2}\H_{0})\cong z_{1}z_{2}\H_{0}
\end{equation}
so that \eqref{eq:map} is a ring isomorphism.
Let $\H_{\lambda}$ be the irreducible, level 1 representation with highest
weight $\lambda$. As noted in \S \ref{ch:classification}.\ref{ss:level 1},
$\lambda$ is a minimal dominant weight since $\Spin_{2n}$ is simply--laced
and therefore, by theorem \ref{ch:classification}.\ref{th:Z on per} and
the explicit action of $Z(\Spin_{2n})$ on the level 1 alcove given by
proposition \ref{ch:classification}.\ref{geometric action},
$\H_{\lambda}=\exp_{T}(-2\pi i\lambda)\H_{0}$. Consequently,
$\H_{\lambda}\fuse\H_{\mu}=\H_{\nu}$ where $\nu$ is the unique minimal
dominant weight in the root lattice coset of $\lambda+\mu$. The fusion
rules \eqref{eq:1 fuse 1}--\eqref{eq:1 fuse 2} now follow from the tables
in \S \ref{ss:Z on alcove 2} of chapter \ref{ch:classification} \halmos

\ssection{Braiding properties of bounded intertwiners}\label{se:bounded braiding}

This section follows \cite{Wa2}. We derive braiding relations satisfied by {\it bounded}
intertwiners for the local loop groups $\lig$ from the corresponding ones for smeared
primary fields. These will be used in the next section to give an upper bound for the
fusion with the vector representation.\\

Using Young diagrams, denote by $\Vbox$ the vector representation of $\SO_{2n}$
and by $\Hbox$ the corresponding positive energy representation at level $\ell$.
Let $\H_{i}$ be an irreducible positive energy representation with lowest energy
subspace $V_{i}$ and consider the braiding relation
\begin{align}
\field{i}{i}{0}(w)\field{\overline{\Box}}{0}{\Box}(z)&=
\sum_{j}\lambda_{j}
\field{\overline{\Box}}{i}{j}(z)\field{i}{j}{\Box}(w)
\label{vertex braid}\\
\intertext{and the abelian braiding}
\field{i}{j}{\Box}(w)\field{\Box}{\Box}{0}(z)&=
\epsilon_{j}
\field{\Box}{j}{i}(z)\field{i}{i}{0}(w)
\label{ab braid}
\end{align}
where $j$ labels the irreducible summands $V_{j}$ of $V_{i}\otimes V_{\Box}$ which
are admissible at level $\ell$. Since $\Vbox$ is minimal, each $V_{j}$ has multiplicity
one by proposition \ref{ch:classification}.\ref{pr:tensor with minimal}.
We normalise the above as follows. The initial terms of the primary fields
$\field{i}{i}{0},\field{\Box}{\Box}{0},\field{\overline{\Box}}{0}{\Box}$ with vertices 
$\vertex{V_{i}}{V_{i}}{\IC},\vertex{\Vbox}{\Vbox}{\IC},\vertex{\Vbox^{*}}{\IC}{\Vbox}$
are the corresponding canonical $\Spin_{2n}$--intertwiners. Fix for any $j$ a generator
$\varphi_{j}\in\Hom_{\Spin_{2n}}(V_{\Box}\otimes V_{i},V_{j})\cong\IC$. This determines
initial terms for $\field{i}{j}{\Box},\field{\Box}{j}{i}$ and
$\field{\overline{\Box}}{i}{j}$ 
by taking $\varphi_{j},\varphi_{j}\sigma$ and $(\varphi_{j}\sigma)^{*}$ respectively,
where $\sigma:V_{i}\otimes V_{\Box}\rightarrow V_{\Box}\otimes V_{i}$ is permutation.
With these normalisations, $\lambda_{j}$ depends upon $\varphi_{j}$ only up to a
positive factor and, by (i) of lemma \ref{ch:algebraic fields}.\ref{le:abelian}, 
$\epsilon_{j}$ is independent of the choice of $\varphi_{j}$ and of modulus 1.\\

Let now $f\in C^{\infty}(S^{1},V_{i})$ and $g\in C^{\infty}(S^{1},V_{\Box})$
be supported in $I,I^{c}$ respectively. Since all the vertices above involve
one of the minimal representations $\IC,V_{\Box}$, all primary fields in
\eqref{vertex braid}--\eqref{ab braid} extend to operator--valued distributions
by theorem \ref{ch:sobolev fields}.\ref{th:minimal is sobolev}. We may therefore
use proposition \ref{ch:algebraic fields}.\ref{pr:smeared braiding} to smear
the braiding relations and find
\begin{align}
\field{i}{i}{0}(f)\field{\overline{\Box}}{0}{\Box}(\overline{g})&=
\sum_{j}\lambda_{j}
\field{\overline{\Box}}{i}{j}(\overline{g}e_{\alpha_{j}})
\field{i}{j}{\Box}(fe_{-\alpha_{j}})
\label{seed one}\\
\field{i}{j}{\Box}(fe_{-\alpha_{j}})\field{\Box}{\Box}{0}(g)&=
\epsilon_{j}\field{\Box}{j}{i}(g e_{-\alpha_{j}})\field{i}{i}{0}(f)
\label{seed two}
\end{align}
where $\alpha_{j}=\Delta_{i}+\Delta_{\Box}-\Delta_{j}$. With our normalisations,
$\field{\overline{\Box}}{0}{\Box}=(\field{\Box}{\Box}{0})^{*}$ and therefore
$\field{\overline{\Box}}{0}{\Box}(\overline{g})\subseteq
 (\field{\Box}{\Box}{0}(g))^{*}$ which is in fact an equality since the vector
primary field is bounded by theorem
\ref{ch:sobolev fields}.\ref{th:minimal is sobolev}. Rewriting the above more
synthetically, it follows that there exist operators
$x_{qp},y_{qp}:\hsmooth_{p}\rightarrow\hsmooth_{q}$ with the $y_{qp}$ bounded
such that
\begin{xalignat}{2}\label{seed}
x_{i0}   y_{\Box 0}^{*}	&=\sum_{j}\lambda_{j} y_{ji}^{*}x_{j\Box}&
x_{j\Box}y_{\Box 0}	&=\epsilon_{j}y_{ji}x_{i0}
\end{xalignat}
Moreover, by proposition \ref{ch:sobolev fields}.\ref{pr:LG equivariance}
$x_{qp},y_{qp}$ commute with $\licg$ and $\lig$ respectively.\\

\remark For suitable choices of $f$ and $g$, none of the above products of
operators, and {\it a fortiori} their individual factors are zero. More
precisely, if $h_{0}=\int h(\theta)\frac{d\theta}{2\pi}$ is the constant
term in the Fourier expansion of $h$, we henceforth choose $f$ and $g$ with
$f_{0},(fe_{-\alpha_{j}})_{0},
 g_{0},(ge_{-\alpha_{j}})_{0}\neq 0$ for any $j$. Then, for any
$v_{j}\in V_{j}=\H_{j}(0)$
\begin{equation} 
(y_{ji}x_{i0}\Omega,v_{j})=
(\phi_{ji}^{\Box}((ge_{-\alpha_{j}})_{0},0)\phi_{i0}^{i}(f_{0},0)\Omega,v_{j})
\end{equation}
The vanishing of the above for any $v_{j}$ implies that of the initial term
of one of the primary fields and therefore of the primary field itself, a
contradiction. Thus, $\xi_{j}=y_{ji}x_{i0}\Omega\neq 0$. The non--zeroness
of the other products follows similarly.\\

We show below that the operators in \eqref{seed} may be replaced by bounded
intertwiners without altering the braiding relations in such a way that
$x_{i0}$ and $y_{\Box0}$ become unitaries. We need a preliminary

\begin{lemma}\label{positivity}
The constants $\mu_{j}=\epsilon_{j}\lambda_{j}$ are non--negative and therefore
$\lambda_{j}\epsilon_{j}=|\lambda_{j}|$.
\end{lemma}
\proof
Consider the functional on $\lig\times\licg$ given by the four--point function
with insertions formula
\begin{equation}\label{eq:functional}
(\gamma_{I},\gamma_{I^{c}})\rightarrow
 (\pi_{i}(\gamma_{I})x_{i0}
  y_{\Box 0}^{*}\pi_{\Box}(\gamma_{I^{c}})y_{\Box 0}\Omega,x_{i0}\Omega)
\end{equation}
Since $y_{\Box 0}$ is bounded,
$y_{\Box 0}^{*}\pi_{\Box}(\gamma_{I^{c}})y_{\Box 0}\in\pi_{0}(\lig)'=\pi_{0}(\licg)''$
and we may rewrite \eqref{eq:functional} as
\begin{equation}
(\pi_{i}(\gamma_{I})\pi_{i}(y_{\Box 0}^{*}\pi_{\Box}(\gamma_{I^{c}})y_{\Box 0})
x_{i0}\Omega,x_{i0}\Omega)=
(\pi_{i}(y_{\Box 0}^{*}\pi_{\Box}(\gamma_{I})\pi_{\Box}(\gamma_{I^{c}})y_{\Box 0})
x_{i0}\Omega,x_{i0}\Omega)
\end{equation}
so that it extends to a positive linear functional on the algebraic tensor product
$\A=\pi_{0}(\lig)''\otimes\pi_{0}(\licg)''$. On the other hand, using the braiding
relations \eqref{seed} we may write \eqref{eq:functional} as
\begin{equation}
\sum_{j}\lambda_{j}
(\pi_{j}(\gamma_{I})
 \pi_{j}(\gamma_{I^{c}})x_{j\Box}y_{\Box 0}\Omega,y_{ji}x_{i0}\Omega)=
\sum_{j}\lambda_{j}\epsilon_{j}
(\pi_{j}(\gamma_{I})\pi_{j}(\gamma_{I^{c}})y_{ji}x_{i0}\Omega,y_{ji}x_{i0}\Omega)
\end{equation}
so that the functional
$\phi:a_{I}\otimes b_{I^{c}}
\rightarrow\sum_{j}\epsilon_{j}\lambda_{j}
(\pi_{j}(a_{I})\pi_{j}(b_{I^{c}})\xi_{j},\xi_{j})$, where $\xi_{j}=y_{ji}x_{i0}\Omega$,
is positive on $\A$.
By von Neumann density, the $\H_{j}$ involved are irreducible and inequivalent
$\A$-modules, and therefore $(\oplus_{j}\pi_{j}(\A))''=\bigoplus_{j}\B(\H_{j})$.
In particular, any $T\in\B(\H_{j})$ is a strong limit of elements of the form
$\oplus_{j}\pi_{j}(a)$ so that $T^{*}T$ is a weak limit of elements
$\oplus_{j}\pi_{j}(a^{*}a)$ and therefore $\phi(T^{*}T)\geq 0$.
Choosing $T=1_{j}$, we find $\phi(T^{*}T)=\lambda_{j}\epsilon_{j}\|\xi_{j}\|^{2}\geq 0$
which, in view of the remark following \eqref{seed}, implies
$\lambda_{j}\epsilon_{j}\geq 0$ \halmos

\begin{lemma}\label{phase}
Let $\hsmooth_{p}$ be dense subspaces of the Hilbert spaces $\H_{p}$, $p=1\ldots 4$
and consider operators
\begin{xalignat}{2}
&
\begin{diagram}[height=2em,width=2.2em]
\hsmooth_{1}&\rTo^{T}&\hsmooth_{2}\\
\dTo^{u}    &        &\dTo_{v}     \\
\hsmooth_{3}&\rTo_{S}&\hsmooth_{4}
\end{diagram}
&&
\begin{diagram}[height=2em,width=2.2em]
\hsmooth_{1}&\rTo^{T}&\hsmooth_{2}\\
\uTo^{u^{*}}&        &\uTo_{v^{*}}\\
\hsmooth_{3}&\rTo_{S}&\hsmooth_{4}
\end{diagram}
\end{xalignat}
with $T,S$ closeable and $u,v$ bounded. If the above diagrams are commutative, they
remain so when $T$ and $S$ are replaced by their phases.
\end{lemma}
\proof
By the boundedness of $u,v$, we get
$v\overline{T}		\subseteq \overline{S}u$ and, in particular
$u\D(\overline{T})	\subseteq \D(\overline{S})$. Similarly,
$v^{*}\overline{S}	\subseteq \overline{T}u^{*}$ whence, taking adjoints,
$uT^{*}			\subseteq S^{*}v$ and
$v\D(T^{*})		\subseteq \D(S^{*})$.
We claim that $uT^{*}\overline{T}\subseteq S^{*}\overline{S}u$. To see this, let
$\xi\in\D(T^{*}\overline{T})$. Then $\xi\in\D(\overline{T})$ so that
$u\xi\in\D(\overline{S})$ and $\overline{S}u\xi=v\overline{T}\xi\in\D(S^{*})$ since
$\overline{T}\xi\in\D(T^{*})$. It follows that
$S^{*}\overline{S}u\xi=S^{*}v\overline{T}\xi=uT^{*}\overline{T}\xi$ as claimed. By
functional calculus, $uf(T^{*}\overline{T})=f(S^{*}\overline{S})u$ for any bounded
measurable function $f$. In particular,
\begin{equation}
v\overline{T}(T^{*}\overline{T}+\epsilon)^{-\half{1}}=
\overline{S}u(T^{*}\overline{T}+\epsilon)^{-\half{1}}=
\overline{S}(S^{*}\overline{S}+\epsilon)^{-\half{1}}u
\end{equation}
which, in the limit $\epsilon\rightarrow 0$ yields the commutativity of the first
diagram when $T$ and $S$ are replaced by their phases. That of the second follows
by the above argument by permuting $(T,S)$, $(u,u^{*})$, $(v,v^{*})$,
$(\H_{1},\H_{3})$ and $(\H_{2},\H_{4})$ \halmos

\begin{proposition}\label{replace}
There exist non--zero bounded operators $x_{qp},y_{qp}:\H_{p}\rightarrow\H_{q}$
commuting with $\licg$ and $\lig$ respectively such that
\begin{xalignat}{2}\label{bounded braiding}
x_{i0}   y_{\Box 0}^{*}	&=\sum_{j}\lambda_{j} y_{ji}^{*}x_{j\Box}&
x_{j\Box}y_{\Box 0}	&=\epsilon_{j}y_{ji}x_{i0}
\end{xalignat}
where $x_{i0},y_{\Box0}$ are unitary and $\lambda_{j},\epsilon_{j}$ are the braiding
coefficients \eqref{vertex braid}--\eqref{ab braid}.
\end{proposition}
\proof We start from the relations \eqref{seed} satisfied by the non--zero operators
$x_{qp},y_{qp}$ and modify the $x_{qp}$ in various steps without altering the $y_{qp}$
and the braiding relations \eqref{seed}.\\

{\it 1st step.} Applying lemma \ref{phase} to
\begin{equation}\label{dgr}
\begin{diagram}[height=2.5em,width=2.8em]
\hsmooth_{0}	 &\rTo^{x_{i0}}	&\hsmooth_{i}\\
\dTo^{y_{\Box 0}}& &\dTo_{\oplus_{j}|\lambda_{j}|^{\half{1}}\epsilon_{j}y_{ji}}\\
\hsmooth_{\Box}	 &\rTo_{\oplus_{j}|\lambda_{j}|^{\half{1}}x_{j\Box}}
				&\bigoplus_{j}\hsmooth_{j}
\end{diagram}
\end{equation}
we may replace the $x_{qp}$ in \eqref{seed} by bounded operators. This yields
a non--zero $x_{i0}$, since it is simply the phase of the former $x_{i0}$.\\

{\it 2nd step.}
We wish to modify $x_{i0}$ to make it injective with dense image. Another
application of lemma \ref{phase} will then allow its replacement by a unitary phase.
Although we shall modify the range $R(x_{i0})$ first, notice that $x_{i0}$ may
always be assumed to be injective. Indeed, the projection onto $\Ker(x_{i0})^{\perp}$
lies in $\pi_{0}(\licg)'=\pi_{0}(\lig)''$. Since the latter is a type III factor,
there exists a partial isometry $u\in\pi_{0}(\lig)''$ with initial space $\H_{0}$ and
final space $\Ker(x_{i0})^{\perp}$. Replacing each $x_{pq}$ by $x_{pq}\pi_{q}(u)$
yields an injective $x_{i0}$ without altering the braiding relations.\\

{\it 3rd step.}
Making $x_{i0}$ surjective is a trifle more involved since in general
$\pi_{i}(\licg)'\subsetneq\pi_{i}(\lig)''$ and the previous
device does not apply. We resort to an averaging procedure relying on von
Neumann density or equivalently the fact that $\pi_{i}(\licg)'\cap\pi_{i}(\lig)'=\IC$.
Let $\{g_{n}\}$ be a dense, countable subgroup of $U(\pi_{i}(\lig)'')$ and
$u_{m}\in\pi_{0}(\lig)''$ partial isometries satisfying $u_{m}^{*}u_{n}=\delta_{mn}$.
We may replace each $x_{qp}$ by the norm convergent sum
$\sum_{n}2^{-n}\pi_{q}(g_{n})x_{qp}\pi_{p}(u_{n})$ without altering the braiding
relations. If $p\in\pi_{i}(\licg)'$ is the projection onto
$R(x_{i0})=\Ker(x_{i0}x_{i0}^{*})^{\perp}$, the corresponding projection for
$\wt x_{i0}=\sum_{n}2^{-n}\pi_{i}(g_{n})x_{i0}u_{n}$ is $\wt p=\bigvee_{n}g_{n}pg_{n}^{*}$
since $\wt x_{i0}\wt x_{i0}^{*}\xi=0$ iff $g_{n}x_{i0}x_{i0}^{*}g_{n}\xi=0$ for all
$n$. $\wt p$ commutes with the $g_{n}$ and therefore with $\pi_{i}(\lig)'$.
Thus, $\wt p\in\pi_{i}(\lig)'\cap\pi_{i}(\lig)'=\IC$ whence $\wt p=1$ since
$\wt p\neq 0$.\\

{\it 4th step.} It follows that $x_{i0}$ may be chosen to have dense range and
by our previous argument, to be injective. A further application of lemma \ref{phase} then
yields a unitary $x_{i0}$. Notice that the modified $x_{j\Box}$ are non--zero. Indeed,
by the second braiding relation \eqref{bounded braiding} and the unitarity of $x_{i0}$,
the vanishing of $x_{j\Box}$ implies that of the original $y_{ji}$, a contradiction.\\

Steps 3 and 2 may now be applied to $y_{\Box0}$ yielding an injective operator with
dense range. A final application of lemma \ref{phase} to the diagram \eqref{dgr} reflected
across the NW--SE diagonal then allows $y_{\Box0}$ to be replaced by its unitary phase.
The modified $y_{ji}$ are non--zero since
$y_{ji}=\epsilon_{j}^{*}x_{j\Box}y_{\Box0}x_{i0}^{*}$ \halmos

\ssection{Upper bounds for fusion with the vector representation}
\label{se:upper bound}

\begin{proposition}\label{pr:upper bound}
Let $\H_{i},\H_{\Box}$ be the positive energy representations at level $\ell$
with lowest energy subspaces $V_{i},V_{\Box}$. Then,
\begin{equation}
\H_{i}\fuse\Hbox=
\bigoplus_{j} N_{i\thinspace \Box}^{j}\H_{j}
\end{equation}
where $\H_{j}(0)=V_{j}$ ranges over the summands of $U\otimes V_{\Box}$ which are
admissible at level $\ell$ and $1\geq N_{i\thinspace \Box}^{j}\geq 0$. Moreover,
$N_{i\thinspace \Box}^{j}$ vanishes if, and only if the corresponding braiding
coefficient $\lambda_{j}$ in \eqref{vertex braid} does.
\end{proposition}
\proof
Let $x_{qp},y_{qp}:\H_{p}\rightarrow\H_{q}$ be the intertwiners given by proposition
\ref{replace}. Since $x_{i0}$ and $y_{\Box0}$ are unitaries, we have
$\XX_{i}=\Hom_{\licg}(\H_{0},\H_{i})=x_{i0}\pi_{0}(\licg)'=x_{i0}\pi_{0}(\lig)''$
and similarly $\YY_{\Box}=y_{\Box0}\pi_{0}(\lig)''$.
Moreover, if $a_{1},a_{2}\in\pi_{0}(\lig)''$ and $b_{1},b_{2}\in\pi_{0}(\licg)''$,
we find by \eqref{bounded braiding}
\begin{equation}
\<x_{i0}a_{1}\otimes y_{\Box0}b_{1},x_{i0}a_{2}\otimes y_{\Box0}b_{2}\>=
\sum_{j}\lambda_{j}\epsilon_{j}
(y_{ji}\pi_{i}(b_{1})x_{i0}a_{1}\Omega,y_{ji}\pi_{i}(b_{2})x_{i0}a_{2}\Omega)
\end{equation}
It follows from lemma \ref{positivity} that the map
$U:\XX_{i}\otimes\YY_{\Box}\rightarrow\bigoplus_{j}\H_{j}$,
$x\otimes y\rightarrow
 \bigoplus_{j}|\lambda_{j}|^{\half{1}}y_{ji}\pi_{i}(y_{\Box0}^{*}y)x\Omega$
extends to an isometry $\H_{i}\fuse\H_{\Box}\rightarrow\bigoplus_{j}\H_{j}$
which is easily seen to commute with the action of $\lig\times\licg$.
Notice that the range of $U$ intersects non--trivially all $\H_{j}$ such that
$\lambda_{j}\neq 0$. Indeed, for any such $\lambda_{j}$,
$U(\XX_{i}\otimes y_{\Box0})\cap\H_{j}=y_{ji}\XX_{i}\Omega$ which is a non--zero
subspace of $\H_{j}$ since $y_{ji}\neq 0$ and $\overline{\XX_{i}\Omega}=\H_{i}$.
Thus, by the irreducibility of the $\H_{j}$ under $\lig\times\licg$, the image
of $U$ is precisely $\bigoplus_{j:\lambda_{j}\neq 0}\H_{j}$ \halmos

\ssection{Fusion of the vector representation with its symmetric powers}
\label{se:fusion box/kbox}

An important corollary of proposition \ref{pr:upper bound} and the 
computations of chapter \ref{ch:box/kbox braiding} is the following

\begin{proposition}\label{pr:box/kbox fuse}
Let $\H_{\theta_{1}}=\H_{\Box}$ and $\H_{k\theta_{1}}$ be the positive energy
representations at level $\ell$ corresponding to the defining representation
$V_{\theta_{1}}$ of $\SO_{2n}$ and its $k$--fold symmetric, traceless power
$V_{k\theta_{1}}$. Then,
\begin{align}
\H_{k\theta_{1}}\fuse\H_{\theta_{1}}&=
\H_{(k-1)\theta_{1}}\oplus
\H_{k\theta_{1}+\theta_{2}}\oplus
\H_{(k+1)\theta_{1}}
\label{eq:sym 1}\\
\intertext{if $1\leq k\leq\ell-1$ and}
\H_{\ell\theta_{1}}\fuse\H_{\theta_{1}}&=
\H_{(\ell-1)\theta_{1}}
\label{eq:sym 2}
\end{align}
\end{proposition}
\proof
By theorem \ref{ch:box/kbox braiding}.\ref{th:braiding identity}, the constants
governing the braiding relations
\begin{align}
\field{k\theta_{1}}{k\theta_{1}}{0}(w)
\field{ \theta_{1}}{0}{ \theta_{1}}(z)
&=
\sum_{\mu\in\{(k-1)\theta_{1},k\theta_{1}+\theta_{2},(k+1)\theta_{1}\}}
\beta_{k\theta_{1},\mu}
\field{ \theta_{1}}{k\theta_{1}}{\mu}(z)
\field{k\theta_{1}}{\mu}{ \theta_{1}}(w)\\
\intertext{if $k=1\ldots\ell-1$ and}
\field{\ell\theta_{1}}{\ell\theta_{1}}{0}(w)
\field{      \theta_{1}}{0}{       \theta_{1}}(z)
&=
\beta_{\ell\theta_{1},(\ell-1)\theta_{1}}
\field{    \theta_{1}}{\ell\theta_{1}}{(\ell-1)\theta_{1}}(z)
\field{\ell\theta_{1}}{(\ell-1)\theta_{1}}{    \theta_{1}}(w)
\end{align}
are non--zero. The fusion rules \eqref{eq:sym 1}--\eqref{eq:sym 2} now follow
from proposition \ref{pr:upper bound} \halmos

\ssection{The fusion ring $\mathcal R_{0}$}\label{se:subring}

Let $\H_{\theta_{1}}=\Hbox$ be the vector representation of $LG$ at a fixed
level $\ell$. We consider below the vector space $\RR$ generated by the
irreducible summands of the iterated fusion products $\Hbox^{\fuse k}$ and
show that it is closed under fusion and forms an associative and commutative
ring. The latter property is a direct consequence of the existence of a
{\it braiding} operator $B$ giving a unitary map
$\H_{i}\fuse\H_{j}\rightarrow\H_{j}\fuse\H_{i}$ intertwining $LG$. A number
of other important properties of $\RR$, most notably the existence of a
positive character or quantum dimension function, will be established in
section \ref{se:DHR}.\\

Let $L_{0}$ be the infinitesimal generator of rotations given by the Segal--Sugawara
formula (\ref{ch:analytic}.\ref{sugawara})

\begin{lemma}[\cite{Wa3}]\label{le:braiding operator}
Let $\H_{i},\H_{j}\in\pl$ be such that $\H_{j}\fuse\H_{i}$ is of positive energy.
Then, if $\XX_{i},\YY_{j}$ are the spaces defined by \eqref{eq:hom spaces}, the
operator $B:\XX_{i}\otimes\YY_{j}\longrightarrow\H_{j}\fuse\H_{i}$ given by
\begin{equation}\label{eq:definition of B}
B(x\otimes y)= e^{-i\pi L_{0}}\medspace y^{\pi}\otimes x^{\pi}
\end{equation}
where $z^{\pi}=e^{i\pi L_{0}}z e^{-i\pi L_{0}}$ and the $e^{-i\pi L_{0}}$
on the right hand--side of \eqref{eq:definition of B} refers to the positive
energy action of $\rot$ on $\H_{j}\fuse\H_{i}$, extends to a unitary
$\H_{i}\fuse\H_{j}\rightarrow\H_{j}\fuse\H_{i}$ intertwining $\LG$.
\end{lemma}
\proof
Notice that $B$ is well--defined since $y^{\pi}\in\XX_{j}$ and $x^{\pi}\in\YY_{i}$.
Moreover, $L_{0}\Omega=0$ and therefore
\begin{equation}
\|y^{\pi}\otimes x^{\pi}\|^{2}=
((y^{\pi})^{*}y^{\pi}(x^{\pi})^{*}x^{\pi}\Omega,\Omega)=
(y^{*}yx^{*}x e^{-i\pi L_{0}}\Omega,e^{-i\pi L_{0}}\Omega)=
\|x\otimes y\|^{2}
\end{equation}
so that $B$ is norm--preserving and extends to a unitary since it has dense range. As
is readily verified, $B$ intertwines the actions of $\lig\times\licg$. It follows that
$\H_{i}\fuse\H_{j}$ is of positive energy and, by von Neumann density that $B$
intertwines $\LG$ \halmos\\

\definition By proposition \ref{pr:upper bound}, the iterated fusion products of $\Hbox$
are finitely reducible positive energy representations. Their irreducible summands
therefore generate a vector space $\RR$ which is closed under fusion and therefore
is, by proposition \ref{pr:associativity} and lemma \ref{le:braiding operator} an
associative and commutative ring.

%

\ssection{Eigenvalues of the braiding operator}
\label{se:eigenvalues of braiding}

We compute below the eigenvalues of the braiding operator $B$ on $\Hbox\fuse\Hbox$.
For convenience, we adopt a graphical notation and label the irreducible
$\SO_{2n}$--modules and corresponding positive energy representations by
the Young diagram describing their highest weight. Thus, the second and
third symmetric, traceless powers of the vector representation $\Vbox$
of $\SO_{2n}$ will be denoted by $\Vsym,\Vsymm$ and its second and third
exterior powers by $\Valt,\Valtt$. The corresponding positive energy
representations are $\Hsym,\Hsymm,\Halt$ and $\Haltt$ respectively.
In particular, if the level $\ell$ is greater or equal to 2, proposition
\ref{pr:box/kbox fuse} yields
\begin{equation}\label{eq:summands}
\Hbox\fuse\Hbox=
\Hsym\oplus\Halt\oplus\H_{0}
\end{equation}
The calculation of the eigenvalues of $B$ on \eqref{eq:summands} relies on a version
of the braiding relations of proposition \ref{replace} which is equivariant with
respect to rotation by $\pi$ in the sense that
$y_{qp}=x_{qp}^{\pi}=e^{i\pi L_{0}}x_{qp}e^{-i\pi L_{0}}$.

\begin{lemma}\label{le:equi braiding}
Let $\ell\geq 2$. Then, there exist non--zero bounded operators
$x_{qp}\in\Hom_{\licg}(\H_{p},\H_{q})$ with $x_{\Box 0}$ unitary satisfying
\begin{xalignat}{2}\label{equivariant braiding}
x_{\square 0}(x_{\square 0}^{\pi})^{*}&=
\sum_{j}\lambda_{j}(x_{j \square}^{\pi})^{*}x_{j \square}&
x_{j \square}x_{\square 0}^{\pi}&=
\epsilon_{j}x_{j \square}^{\pi}x_{\square 0}
\end{xalignat}
where $j$ labels the summands of \eqref{eq:summands} and $\lambda_{j}\epsilon_{j}>0$.
Moreover,
\begin{equation}\label{eq:epsilon j}
\epsilon_{j}=\sigma_{j}e^{i\pi(2\Delta_{\square}-\Delta_{j})}
\end{equation}
where $\sigma_{j}=\pm 1$ according to whether $\H_{j}(0)\subset\Vbox\otimes\Vbox$ is
symmetric or anti--symmetric under $\mathfrak S_{2}$.
\end{lemma}
\proof
The proof is almost identical to the discussion of section \ref{se:bounded braiding}
and the proof of proposition \ref{replace} and we simply indicate the points where
they differ.
We start from the smeared braiding relations \eqref{seed one}--\eqref{seed two} where
now $i=\square$ labels the vector representation $\Vbox$ and $j$ the summands of
\begin{equation}
\Vbox\otimes\Vbox=\Vsym\oplus\Valt\oplus\IC
\end{equation}
Let $f\in C^{\infty}(S^{1},\Vbox)$ be supported in $I=(0,\pi)$ and set $g=f^{\pi}$
where $f^{\pi}(\theta)=f(\theta-\pi)$. Then, \eqref{seed one}--\eqref{seed two} read
\begin{align}
\field{\square}{\square}{0}(f)\field{\overline{\square}}{0}{\square}(\overline{f^{\pi}})&=
\sum_{j}\lambda_{j}
\field{\overline{\square}}{\square}{j}(\overline{f^{\pi}e_{-\alpha_{j}}})
\field{\square}{j}{\square}(fe_{-\alpha_{j}})\\
\field{\square}{j}{\square}(fe_{-\alpha_{j}})\field{\square}{\square}{0}(f^{\pi})&=
\epsilon_{j}
\field{\square}{j}{\square}(f^{\pi}e_{-\alpha_{j}})\field{\square}{\square}{0}(f)
\label{eq:new norm}
\end{align}
where $\alpha_{j}=2\Delta_{\square}-\Delta_{j}$. Here, we normalise \eqref{eq:new norm}
by equating the initial terms of the left--most primary fields so that, by (ii) of lemma
\ref{ch:algebraic fields}.\ref{le:abelian}, $\epsilon_{j}$ is given by \eqref{eq:epsilon j}.
By theorem \ref{ch:box/kbox braiding}.\ref{th:braiding identity}, $\lambda_{j}\neq 0$ and
therefore, by lemma \ref{positivity}, $\lambda_{j}\epsilon_{j}>0$.
Notice that $f^{\pi}e_{-\alpha_{j}}=e^{-i\pi\alpha_{j}}(fe_{-\alpha_{j}})^{\pi}$ since
$f$ is supported in $I$. Moreover, by the equivariance of primary fields under the
integrally--moded action $R_{\theta}=e^{i\theta d}$ of $\rot$
(proposition \ref{ch:sobolev fields}.\ref{pr:LG equivariance}),
\begin{equation}
 \phi_{qp}(g^{\pi})=
 R_{\pi}\phi_{qp}(g)R_{-\pi}=
 e^{i\pi(\Delta_{p}-\Delta_{q})}e^{i\pi L_{0}}\phi_{qp}(g)e^{-i\pi L_{0}}=
 e^{i\pi(\Delta_{p}-\Delta_{q})}\phi_{qp}(g)^{\pi}
\end{equation}
and therefore
\begin{align}
\field{\square}{\square}{0}(f)
\field{\overline{\square}}{0}{\square}(\overline{f})^{\pi}&=
\sum_{j}\lambda_{j}
\field{\overline{\square}}{\square}{j}(\overline{fe_{-\alpha_{j}}})^{\pi}
\field{\square}{j}{\square}(fe_{-\alpha_{j}})\\
\field{\square}{j}{\square}(fe_{-\alpha_{j}})\field{\square}{\square}{0}(f)^{\pi}&=
\epsilon_{j}
\field{\square}{j}{\square}(fe_{-\alpha_{j}})^{\pi}\field{\square}{\square}{0}(f)
\end{align}
Using the adjunction property
$\field{\overline{\square}}{q}{p}(\overline{g})=\field{\square}{p}{q}(g)^{*}$,
we find operators $x_{qp}$ satisfying \eqref{equivariant braiding}.
The $x_{qp}$ may now be modified in various stages as in the proof of proposition
\ref{replace} with the additional requirement that the $x_{qp}^{\pi}$ be altered
accordingly. This yields a new set of operators satisfying the same braiding
relations but with $x_{\square 0}$ unitary \halmos

\begin{proposition}\label{eigenvalues}
Let $\Hbox$ be the vector representation at level $\ell\geq 2$ and
$B\in\End_{\LG}(\Hbox\fuse\Hbox)$ the braiding operator defined by lemma
\ref{le:braiding operator}. Then, the eigenvalues of $B$ corresponding
to the summands of \eqref{eq:summands} are distinct and given by
\begin{xalignat}{3}\label{evals}
\beta_{\sym}	&= q	 &
\beta_{\alt}	&=-q^{-1}&
\beta_{0}	&= r^{-1}
\end{xalignat}
where $q=e^{-\frac{i\pi}{\kappa}}$ with $\kappa=\ell+2(n-1)$ and $r=q^{2n-1}$.
\end{proposition}
\proof
Let $j$ label the summands of \eqref{eq:summands}. By lemma \ref{le:equi braiding}
and the proof of proposition \ref{pr:upper bound}, the isomorphism
$\Hbox\fuse\Hbox\rightarrow\bigoplus_{j}\H_{j}$ is given by
\begin{equation}\label{eq:isom}
a\otimes b
\longrightarrow\bigoplus_{j}|\lambda_{j}\epsilon_{j}|^{\half{1}}
x_{j\Box}^{\pi}\pi_{\Box}((x_{\Box 0}^{\pi})^{*}b)a\Omega
\end{equation}
for any
$a\in\XX_{\Box}=\Hom_{\licg}(\H_{0},\Hbox)$ and
$b\in\YY_{\Box}=\Hom_{\lig} (\H_{0},\Hbox)$. Thus, by \eqref{eq:definition of B},
$B(a\otimes b)$ is mapped to
\begin{equation}\label{eq:first attempt}
\begin{split}
  \bigoplus_{j}|\lambda_{j}\epsilon_{j}|^{\half{1}}
  e^{-i\pi L_{0}}x_{j\Box}^{\pi}\pi_{\Box}((x_{\Box 0}^{\pi})^{*}a^{\pi})b^{\pi}\Omega
&=\bigoplus_{j}|\lambda_{j}\epsilon_{j}|^{\half{1}}
  e^{-i\pi L_{0}}x_{j\Box}^{\pi}b^{\pi}\pi_{0}((x_{\Box 0}^{\pi})^{*}a^{\pi})\Omega\\
&=\bigoplus_{j}|\lambda_{j}\epsilon_{j}|^{\half{1}}
  x_{j\Box}bx_{\Box 0}^{*}a\Omega\\
&=\bigoplus_{j}|\lambda_{j}\epsilon_{j}|^{\half{1}}
  x_{j\Box}x_{\Box 0}^{\pi}(x_{\Box 0}^{\pi})^{*}bx_{\Box 0}^{*}a\Omega\\
&=\bigoplus_{j}|\lambda_{j}\epsilon_{j}|^{\half{1}}
  x_{j\Box}x_{\Box 0}^{\pi}x_{\Box 0}^{*}\pi_{\Box}((x_{\Box 0}^{\pi})^{*}b)a\Omega
\end{split}
\end{equation}
which, using \eqref{equivariant braiding} yields
\begin{equation}\label{eq:braided}
\bigoplus_{j}\epsilon_{j}|\lambda_{j}\epsilon_{j}|^{\half{1}}
  x_{j\Box}^{\pi}\pi_{\Box}((x_{\Box 0}^{\pi})^{*}b)a\Omega
\end{equation}
Comparing \eqref{eq:braided} with \eqref{eq:isom} we find that the eigenvalues of $B$
are the $\epsilon_{j}$ given by \eqref{eq:epsilon j} where, as customary
$\Delta_{j}=\frac{C_{j}}{2\kappa}$. Since the Casimir of a representation of highest
weight $\lambda$ is $\<\lambda,\lambda+2\rho\>$ where
$2\rho=2\sum_{j=1}^{n}(n-j)\theta_{j}$ is the sum of the positive roots of $\Spin_{2n}$,
we have $C_{0}=0$, $C_{\Box}=2n-1$, $C_{\sym}=4n$ and $C_{\alt}=4(n-1)$ and therefore
\eqref{evals} holds \halmos\\

\remark The eigenvalues of the braiding operator agree with those of the generators of the
braid group $B_{m}$ on $m$ strings when the latter is mapped into the algebra $C_{m}(q,r)$
defined by Wenzl \cite{We2}. Indeed, the latter is the
quotient of the group algebra $\IC B_{m}$ with generators $g_{i}$, $i=1\ldots m-1$
corresponding to the permutation of two succesive strings by a set of relations comprising
$(g_{i}-r^{-1})(g_{i}+q^{-1})(g_{i}-q)=0$. The underlying reason for this is that, for any
$m$, $\End_{\LG}(\H_{\square}^{\boxtimes m})$ is isomorphic to $C_{m}(q,q^{2n-1})$,
though the proof of this fact requires the knoweldge of the fusion rules with $\H_{\Box}$.

\ssection{The Doplicher--Haag--Roberts theory of superselection sectors}\label{se:DHR}

We require some consequences of the theory of superselection sectors, most importantly
the existence of a positive homomorphism or {\it quantum dimension} function on the ring
$\RR$ of section \ref{se:subring}, which were originally obtained by
Doplicher--Haag--Roberts \cite{DHR1,DHR2} and Fredenhagen--Rehren--Schroer \cite{FRS}.
The precise relation between the predictions of these theories and the bimodule
framework for loop groups was first explained in \cite{Wa4}. A brief account follows.\\

Let $\H\in\pl$ be an irreducible positive energy representation of $LG$ and
assume that $\H$ has a conjugate, {\it i.e.}~ an irreducible $\overline{\H}\in\pl$
such that $\H\fuse\overline{\H}\supseteq\H_{0}$ where $\H_{0}$ is the vacuum
representation. Then, $\overline{\H}$ is unique and $\H_{0}$ is contained in
$\H\fuse\overline{\H}$ with multiplicity one. Assume further that all iterated
fusion products $\H^{\fuse m}\fuse\overline{\H}^{\fuse n}$ are of positive energy
and let $\RRR$ be the ring additively generated by their irreducible summands
$\H_{i}$. Then,

\begin{enumerate}
\item[$\bullet$] There exists a unique faithful trace on each
$\End_{\LG}(\H_{i_{1}}\fuse\cdots\fuse\H_{i_{m}})$ normalised by
$\tr(1)=1$ and consistent with the inclusions
\begin{xalignat}{2}
\End_{\LG}(\H_{i_{1}}\fuse\cdots\fuse\H_{i_{m}})&\hookrightarrow
\End_{\LG}(\H_{i_{1}}\fuse\cdots\fuse\H_{i_{m}}\fuse H_{j})&
x&\rightarrow x\fuse 1\\
\End_{\LG}(\H_{i_{1}}\fuse\cdots\fuse\H_{i_{m}})&\hookrightarrow
\End_{\LG}(H_{j}\fuse\H_{i_{1}}\fuse\cdots\fuse\H_{i_{m}})&
x&\rightarrow 1\fuse x
\end{xalignat}
\item[$\bullet$] (Jones' relations)
Each $\H_{i}$ possesses a (necessarily unique) conjugate. If
$e_{k}\in\End_{\LG}(\H_{i}\fuse\overline{\H_{i}}\fuse\H_{i})$, $k=1,2$
are the (Jones) projections onto
$\H_{0}\fuse\H_{i}\subseteq(\H_{i}\fuse\overline{\H_{i}})\fuse\H_{i}$ and
$\H_{i}\fuse\H_{0}\subseteq\H_{i}\fuse(\overline{\H_{i}}\fuse\H_{i})$ respectively,
then
\begin{xalignat}{3}\label{eq:Jones}
e_{1}e_{2}e_{1}&=\tau e_{1}&&\text{where}&\tau=\tr(e_{1})
\end{xalignat}
\item[$\bullet$](Markov property I)
Let $e\in\End_{\LG}(\H_{i}\fuse\overline{\H_{i}})$ be the Jones projection onto $\H_{0}$
and $a\in\End_{\LG}(\H_{j}\fuse\H_{i})$. Then $a_{1}=a\fuse 1$ and
$e_{2}=1\fuse e\in\End_{\LG}(\H_{j}\fuse\H_{i}\fuse\overline{\H_{i}})$ satisfy
\begin{equation}\label{eq:markov 1}
\tr(a_{1}e_{2})=\tr(a_{1})\tr(e_{2})
\end{equation}
\end{enumerate}

The above properties are in fact true for any system of bimodules over a type
III factor. In our setting however, the following additional properties hold

\begin{enumerate}
\item[$\bullet$](Statistics or braiding operators) There exists a canonical
 unitary $g\in\End_{\LG}(\H_{i}\fuse\H_{i})$ and the operators
$g_{1}=g\fuse 1,g_{2}=1\fuse g\in\End_{\LG}(\H_{i}\fuse\H_{i}\fuse\H_{i})$
satisfy
\begin{equation}\label{eq:braid relations}
g_{1}g_{2}g_{1}=g_{2}g_{1}g_{2}
\end{equation}
Moreover, $g$ coincides with the brading operator $B$ defined by lemma
\ref{le:braiding operator}.
\item[$\bullet$](Markov property II)
Let $g_{2}=1\fuse g\in\End_{\LG}(\H_{j}\fuse\H_{i}\fuse\H_{i})$ be the braiding map
corresponding to the last two factors. Then, for any $a\in\End_{\LG}(\H_{j}\fuse\H_{i})$
and $a_{1}=a\fuse 1\in\End_{\LG}(\H_{j}\fuse\H_{i}\fuse\H_{i})$,
\begin{equation}
\tr(a_{1}g^{\pm 1}_{2})=\tr(a_{1})\tr(g^{\pm 1}_{2})
\end{equation}
\item[$\bullet$](Compatibility)
If $g\in\End_{\LG}(\H_{i}\fuse\H_{i})$ and $e\in\End_{\LG}(\H_{i}\fuse\overline{\H_{i}})$
are the braiding operator and Jones projection respectively, and
$g_{1}=g\fuse 1,e_{2}=1\fuse e\in\End_{\LG}(\H_{i}\fuse\H_{i}\fuse\overline{\H_{i}})$,
then
\begin{xalignat}{3}
e_{2}g_{1}e_{2}&=\lambda e_{2}&&\text{where}&|\lambda|&=|\tr(g_{2})|
\end{xalignat}
\item[$\bullet$](Statistical or quantum dimension)
The map
\begin{equation}\label{eq:defn of dim}
d(\H_{i})=1/\sqrt{\tr(e)}
\end{equation}
where $e\in\End_{\LG}(\H_{i}\fuse\overline{\H_{i}})$ is the Jones projection, extends
to a positive homomorphism of the ring $\RRR$ into $\IR$.
\end{enumerate}

The foregoing discussion and the fusion rules of proposition \ref{pr:box/kbox fuse} imply

\begin{corollary}\label{co:observation}
Let $\RR$ be the ring generated by the irreducible summands of the iterated fusion
products of $\Hbox$ defined in section \ref{se:subring}. Then, there exists a
positive homomorphism or quantum dimension function $d:\RR\rightarrow\IR$.
\end{corollary}
\proof
By proposition \ref{pr:box/kbox fuse}, we have $\Hbox\fuse\Hbox\supseteq\H_{0}$ so
that $\Hbox$ is self--conjugate and $\RR$ possesses a positive character \halmos

\ssection{The quantum dimension of $\H_{\square}$ by Wenzl's lemmas}\label{se:qdim of box}

Assuming $\ell\geq 2$, we compute in this section the quantum dimension $d(\Hbox)$
of the vector representation by adapting a calculation of Wenzl \cite{We3} to our
framework. We proceed as follows.
By \eqref{eq:defn of dim} and \eqref{eq:Jones}, $d(\Hbox)$ may be derived from the
Jones relations $e_{1}e_{2}e_{1}=\tau e_{1}$ satisfied by the projections $e_{i}$.
To this end, we use a representation of
$\A_{3}=\End_{\LG}(\Hbox\fuse\Hbox\fuse\Hbox)\ni e_{1},e_{2}$ where the braiding
operators $g_{1},g_{2}$ of \eqref{eq:braid relations} have a simple form. This
gives matrix representatives for the $e_{i}$ since, by virtue of proposition
\ref{eigenvalues}, the $g_{i}$ have distinct eigenvalues and the $e_{i}$ therefore
are spectral projections of the $g_{i}$. $d(\Hbox)$ is then obtained by comparing
$e_{1}e_{2}e_{1}$ and $e_{1}$.\\

We begin by using our partial knowledge of the fusion rules to describe the
structure of $\A_{3}$. By proposition \ref{pr:box/kbox fuse},
\begin{equation}
\Hbox\fuse\Hbox = \Hsym\oplus\Halt\oplus\H_{0}\label{box}
\end{equation}
and
\begin{xalignat}{4}
\Hbox\fuse\Hsym &= \Hsymm\oplus\Htwone\oplus\Hbox&
&\text{if $\ell\geq 3$}&
\Hbox\fuse\Hsym &= \Hbox&
&\text{if $\ell=2$}
\label{sym}\\
\intertext{Moreover, by proposition \ref{pr:upper bound},}
\Hbox\fuse\Halt &\subseteq \Htwone\oplus\Haltt\oplus\Hbox&
&\text{if $\ell\geq 3$}&
\Hbox\fuse\Halt &\subseteq \Haltt\oplus\Hbox&
&\text{if $\ell=2$}
\label{alt}
\end{xalignat}
since $\Vtwone$ is not admissible at level 2. In fact,

\begin{lemma}
If $\ell\geq 2$, then $\displaystyle{\Hbox\fuse\Halt\supseteq\Hbox}$
\end{lemma}
\proof
From \eqref{box} we see that $\Hbox$ is self--conjugate and therefore that all summands
of \eqref{box} possess conjugates. Conjugating, we find
$\overline{\Halt}\subseteq\Hbox^{\fuse 2}$ and therefore
\begin{equation}
\Hbox\fuse(\Hbox\fuse\Halt)=
(\Hbox\fuse\Hbox)\fuse\Halt
\supseteq\H_{0}
\end{equation}
so that, by uniqueness of conjugates, $\Hbox=\overline{\Hbox}\subseteq\Hbox\fuse\Halt$
\halmos\\

The above fusion rules may be used to describe the inclusion of finite--dimensional
algebras $\A_{1}\subseteq\A_{2}\subseteq\A_{3}$ where $\A_{k}=\End_{\LG}(\Hbox^{\fuse k})$.
This is best encoded in a {\it Bratteli diagram} as follows. Each $\A_{k}$ is a direct sum
of matrix algebras labelled by the inequivalent irreducible summands of $\Hbox^{\fuse k}$.
The corresponding labels are drawn in a row. Lines are drawn from the row corresponding
to an algebra to the next, connecting the matrix block of the smaller algebra to each
block of the larger which contains it. We therefore get

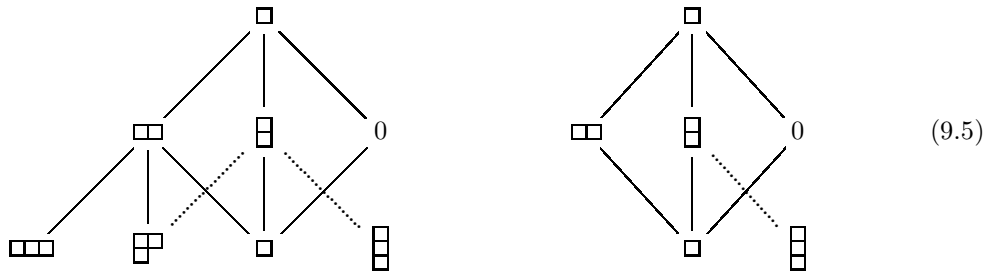
\newarrow{Ddot}.....
\begin{xalignat}{2}\label{bratteli}
\begin{diagram}[height=2.2em,width=2.2em]
	&    	& 	&     		&\Box	&     		&	\\
	&    	& 	&\ldLine	&\vLine	&\rdLine	&	\\
	&	&\sym	&     		&\alt	&     		&0	\\
	&\ldLine&\vLine	&\rdLine\ldDdot	&\vLine	&\rdDdot\ldLine	&	\\
\symm	&    	&\twone	&		&\Box	&     		&\altt	\\
\end{diagram}&&
\begin{diagram}[height=2.2em,width=2.0em]
 	&     		&\Box	&     		&	\\
 	&\ldLine	&\vLine	&\rdLine	&	\\
\sym	&     		&\alt	&     		&0	\\
	&\rdLine	&\vLine	&\rdDdot\ldLine	&	\\
	&		&\Box	&     		&\altt	\\
\end{diagram}&
\end{xalignat}

according to whether $\ell\geq 3$ or $\ell=2$. The dotted lines indicate that $\Htwone$
or $\Haltt$ may fail to appear in $\Halt\fuse\Hbox$, in accordance with \eqref{alt}.
Let now $g$ be the braiding operator on $\Hbox\fuse\Hbox$ and $e$ the Jones projection
onto $\H_{0}\subset\Hbox\fuse\Hbox$. By proposition \ref{eigenvalues}, $g=B$ has three
distinct eigenvalues $q,-q^{-1},r^{-1}$ so that $e$ is the spectral projection of $g$
corresponding to $r^{-1}$ and therefore
\begin{equation}\label{eq:spectral}
e=\frac{(g-q)(g+q^{-1})}{(r^{-1}-q)(r^{-1}+q^{-1})}
\end{equation}
The same relation binds $g_{i},e_{i}\in\A_{3}$ where $e_{1}=e\fuse 1$ and $e_{2}=1\fuse e$
are the Jones projections onto
$\H_{0}\fuse\Hbox\subset(\Hbox\fuse\Hbox)\fuse\Hbox$ and
$\Hbox\fuse\H_{0}\subset\Hbox\fuse(\Hbox\fuse\Hbox)$ respectively and
$g_{1}=g\fuse 1,g_{2}=1\fuse g$ are the braiding operators on $\Hbox\fuse\Hbox\fuse\Hbox$.
Since $\H_{0}\fuse\Hbox\cong\Hbox\cong\Hbox\fuse\H_{0}$ is irreducible, the $e_{i}$ are
minimal projections in $\A_{3}$ and $e_{i}\leq p_{\Box}$ where the latter is the
minimal central projection of $\A_{3}$ corresponding to the matrix block labelled by the
vector representation.
By minimality of $e_{2}$, the left action of $\A_{3}$ on $\A_{3}e_{2}$ is an irreducible
representation of $\A_{3}$ which is faithful on $p_{\Box}\A_{3}p_{\Box}$ and is 
therefore three--dimensional. We shall use it to find explicit matrix representatives of the
$e_{i}$. The following gives a convenient basis of $\A_{3}e_{2}$

\begin{lemma}[Wenzl,\cite{We3}]\label{le:basis of Ae}
If $\ell\geq 2$ then $e_{2},e_{1}e_{2},g_{1}e_{2}\in\A_{3}=\End_{\LG}(\Hbox^{\fuse 3})$
are linearly independent and therefore form a basis of the left $\A_{3}$--module
$\A_{3}e_{2}$.
\end{lemma}
\proof
Notice first that $1,e_{1},g_{1}$ are linearly independent since $g_{1}$ has three
dictinct eigenvalues $q,-q^{-1},r^{-1}$ and $e_{1}$ is the spectral projection
corresponding to $r^{-1}$. Let now $\alpha,\beta,\gamma\in\IC$ be such that
$\alpha e_{2}+\beta e_{1}e_{2}+\gamma g_{1}e_{2}=0$. Multipliplying by
$e_{2}(\alpha+\beta e_{1}+\gamma g_{1})^{*}$, we find
\begin{equation}
e_{2}(\alpha+\beta e_{1}+\gamma g_{1})^{*}
     (\alpha+\beta e_{1}+\gamma g_{1})e_{2}=0
\end{equation}
taking traces and using the Markov property \eqref{eq:markov 1} yields
\begin{equation}
\begin{split}
\tr(e_{2}(\alpha+\beta e_{1}+\gamma g_{1})^{*}
         (\alpha+\beta e_{1}+\gamma g_{1})e_{2})
&=
\tr((\alpha+\beta e_{1}+\gamma g_{1})^{*}
    (\alpha+\beta e_{1}+\gamma g_{1})e_{2})\\
&=
\tr((\alpha+\beta e_{1}+\gamma g_{1})^{*}
    (\alpha+\beta e_{1}+\gamma g_{1}))\tr(e_{2})\\
&=0
\end{split}
\end{equation}
and therefore, by faithfulness of the trace $\alpha+\beta e_{1}+\gamma g_{1}=0$
whence $\alpha=\beta=\gamma=0$ \halmos

\begin{proposition}[Wenzl, \cite{We3}]
If $\ell\geq 2$, then $e_{1}e_{2}e_{1}=\tau e_{1}$ where
$\tau=\displaystyle{\bigl(1+\frac{r-r^{-1}}{q-q^{-1}}\bigr)^{-2}}$.
\end{proposition}
\proof
We begin by computing the matrices representing the braiding operators $g_{1},g_{2}$ on
$\A_{3}e_{2}$ in the basis $e_{2},g_{1}e_{2},e_{1}e_{2}$ of lemma \ref{le:basis of Ae}.
The first columns of $g_{1}$ and $g_{2}$ and the last column of $g_{1}$ are obvious since
$g_{i}e_{i}=r^{-1}e_{i}$. The second column of $g_{1}$ may be computed by expressing
$g_{1}^{2}$ in terms of $e_{1},g_{1}$ via \eqref{eq:spectral}. We find
\begin{xalignat}{3}\label{eq:spectral II}
g_{i}^{2}&=1+g_{i}(q-q^{-1})-zr^{-1}(q-q^{-1})e_{i}&
&\text{where}&
z=&\Bigl(1+\frac{r-r^{-1}}{q-q^{-1}}\Bigr)
\end{xalignat}
It will be more convenient to work with $c_{i}=z e_{i}$. It follows from the above that in
the basis $c_{2},g_{1}c_{2},c_{1}c_{2}$ of $\A_{3}e_{2}$, $g_{1}$ and $g_{2}$ are given by
\begin{xalignat}{3}
g_{1}&=
\left(\begin{array}{ccc}
0&1			&0\\
1&q-q^{-1}		&0\\
0&-r^{-1}(q-q^{-1})	&r^{-1}
\end{array}\right)&
&\text{and}&
g_{2}&=
\left(\begin{array}{ccc}
r^{-1}	&a_{12}&a_{13}\\
0	&a_{22}&a_{23}\\
0	&a_{32}&a_{33}
\end{array}\right)
\end{xalignat}
The $a_{ij}$ may be determined from the braid relation $g_{1}g_{2}g_{1}=g_{2}g_{1}g_{2}$.
Indeed, comparing the first columns of the two products we find $a_{12}=a_{22}=0$. Since
the braid relations imply that $g_{1}$ and $g_{2}$ are conjugate, we find by equating their
traces and determinants that $a_{33}=q-q^{-1}$ and $a_{23}a_{32}=1$. The other entries are
now easily found by imposing the braid relations and yield
\begin{equation}
g_{2}=
\left(\begin{array}{ccc}
r^{-1}	&0&-(q-q^{-1})\\
0	&0&1\\
0	&1&q-q^{-1}
\end{array}\right)
\end{equation}
From \eqref{eq:spectral II} and $c_{i}=ze_{i}$ we therefore find
\begin{xalignat}{3}
c_{1}&=
\left(\begin{array}{ccc}
0&0	&0\\
0&0	&0\\
1&r^{-1}&z
\end{array}\right)&
&\text{and}&
c_{2}&=
\left(\begin{array}{ccc}
z&r&1\\
0&0&0\\
0&0&0
\end{array}\right)
\end{xalignat}
so that $c_{1}c_{2}c_{1}=c_{1}$ and therefore $e_{1}e_{2}e_{1}=z^{-2}e_{1}$ as claimed
\halmos

\begin{corollary}\label{dimension}
If $\ell\geq 2$, the quantum dimension of $\Hbox$ is equal to
$\displaystyle{1+\frac{\sin(\frac{(2n-1)\pi}{\kappa})}{\sin(\frac{\pi}{\kappa})}}$
where $\kappa=\ell+2(n-1)$.
\end{corollary}
\proof
From \eqref{eq:defn of dim}, \eqref{eq:Jones} and the previous proposition, we find
$d(\Hbox)=\displaystyle{\bigl(1+\frac{r-r^{-1}}{q-q^{-1}}\bigr)}$ where $q,r$ are
given by proposition \ref{eigenvalues} \halmos\\

\remark Observe that as $\kappa\rightarrow\infty$, $d(\Hbox)\rightarrow 2n$ which is
the dimension of the vector representation of $\SO_{2n}$.

\ssection{The Verlinde algebra}\label{se:Verlinde}

We give below an algebraic characterisation of the {\it Verlinde rules} which
conjecturally describe the fusion of a positive energy representation of
$L\Spin_{2n}$ with the vector representation $\Hbox$. This will be used in the
next section to prove that these rules actually hold.\\

We begin with a heuristic digression.
Let $\pl$ be the set of positive energy representations of $LG$ at level $\ell$
and $\al$ the alcove parametrising the irreducibles in $\pl$. We
define below a set of functions $\phi_{\mu}$ on $\pl$ indexed by $\mu\in\al$
which are to be thought of as the characters of the putative $*$-algebra $\IC[\pl]$
whose multiplication and involution are defined by fusion and conjugation. If
$\lambda\in\al$ and $N_{\lambda}$ is the matrix giving the fusion with
$\H_{\lambda}$, {\it i.e.}~
\begin{equation}\label{putative}
\H_{\lambda}\fuse\H_{\mu}=\bigoplus_{\nu} N_{\lambda\mu}^{\nu}\H_{\nu}
\end{equation}
then applying a character $\phi$ of $\IC[\pl]$ to both sides of \eqref{putative},
we find
$\sum_{\nu}N_{\lambda \mu}^{\nu}\phi(\H_{\nu})=\phi(\H_{\lambda})\phi(\H_{\mu})$
so that the vector $\phi$ with entries $\phi(\H_{\cdot})$ is an eigenvector of
$N_{\lambda}$ with eigenvalue $\phi(\H_{\lambda})$. We now prove this assertion
when $\lambda$ is a minimal weight and $N_{\lambda}$ is given by the
{\it Verlinde rules}

\begin{proposition}\label{vector}
For any integral weight $\mu$ in the level $\ell$ alcove, define a function
$\phi_{\mu}:\pl\rightarrow\IC$ by $\phi_{\mu}(\H_{\nu})=\cchi_{\nu}(S_{\mu})$
where $\cchi_{\nu}$ is the character of the representation $V_{\nu}$ and
$S_{\mu}=\exp_{T}(2\pi i\frac{\mu+\rho}{\kappa})$ where $\kappa=\ell+\half{C_{\g}}$.
If $\lambda$ is a dominant minimal weight and $N_{\lambda}$ is the matrix defined by
\begin{equation}
N_{\lambda\medspace\delta}^{\nu}=
\left\{\begin{array}{cl}
1&\text{if $V_{\nu}\subseteq V_{\lambda}\otimes V_{\delta}$ and $\<\nu,\theta\>\leq\ell$}\\
0&\text{otherwise}
\end{array}\right.
\end{equation}
then, for any $\delta\in\al$,
$\sum_{\nu}N_{\lambda\delta}^{\nu}\phi_{\mu}(\H_{\nu})=
 \phi_{\mu}(\H_{\lambda})\phi_{\mu}(\H_{\delta})$
so that the vector with entries $\phi(\H_{\cdot})$ is an eigenvector of $N_{\lambda}$
with eigenvalue $\phi_{\mu}(\H_{\lambda})$.
\end{proposition}
\proof
Let $\delta\in\al$, then $\cchi_{\lambda}\cchi_{\delta}=\sum_{\nu}\cchi_{\nu}$ where $\nu$
spans the highest weights of the irreducible summands
$V_{\nu}\subset V_{\lambda}\otimes V_{\delta}$. By the tensor product rules with minimal
representations (proposition \ref{ch:classification}.\ref{pr:tensor with minimal}), these
are of the form $\lambda'+\delta$ where $\lambda'\in W\lambda$ is a weight of $V_{\lambda}$
and therefore satisfies
$|\<\lambda',\theta\>|\leq 1$ so that $\<\nu,\theta\>\leq\ell+1$. We need to prove that
$\cchi_{\nu}(S_{\mu})=0$ whenever $\<\nu,\theta\>=\ell+1$. By the Weyl character formula,
$\cchi_{\nu}=\frac{A_{\nu+\rho}}{A_{\rho}}$ where $2\rho$ is the sum of the positive roots
of $\gc$ and $A_{\alpha}(\exp(\beta))=\sum_{w}(-1)^{w}e^{\<w\alpha,\beta\>}$. If
$\<\nu,\theta\>=\ell+1$, then $\nu+\rho$ lies on the affine hyperplane
$H_{\kappa}=\{h|\<h,\theta\>=\kappa\}$ since 
$\kappa=\ell+\half{C_{\g}}=\ell+1+\<\rho,\theta\>$. The affine reflection corresponding to
$H_{\kappa}$ is $t\rightarrow\sigma_{\theta}(t)+\kappa\theta$ where $\sigma_{\theta}$ is the
orthogonal reflection determined by $\theta$, and therefore, by the integrality of $\mu,\rho$
\begin{equation}
\begin{split}
A_{\nu+\rho}(\exp(2\pi i\frac{\mu+\rho}{\kappa}))
&=
A_{\sigma_{\theta}(\nu+\rho)+\kappa\theta}(\exp(2\pi i\frac{\mu+\rho}{\kappa}))\\
&=
\sum_{w\in W}(-1)^{w}
e^{2\pi i\frac{\<w\sigma_{\theta}(\nu+\rho),\mu+\rho\>}{\kappa}}
e^{2\pi i\<w\theta,\mu+\rho\>}\\
&=
(-1)^{\sigma_{\theta}}A_{\nu+\rho}(\mu+\rho)
\end{split}
\end{equation}
so that $\chi_{\nu}(S_{\mu})=0$ as claimed \halmos\\

When $G=\Spin_{2n}$ and $\lambda$ is the highest weight of the vector representation
of $\SO_{2n}$ with corresponding Young diagram $\Box$, we find

\begin{corollary}\label{other dimension}
Let $V_{\Box}\cong\IC^{2n}$ be the vector representation of $\SO_{2n}$ and $\Nbox$
the matrix with entries labelled by the irreducible $\Spin_{2n}$--modules which are
admissible at level $\ell$ with $N_{U\Box}^{W}=1$ if $W\subset U\otimes\Vbox$
and $0$ otherwise. Then, $\phi_{0}$ is an eigenvector of $\Nbox$ with eigenvalue
$\displaystyle{1+\frac{\sin{\frac{(2n-1)\pi}{\kappa}}}{\sin{\frac{\pi}{\kappa}}}}$
where $\kappa=\ell+2(n-1)$.
\end{corollary}
\proof
The dual Coxeter number $\half{C_{\g}}$ of $\SO_{2n}$ is $2(n-1)$ and we simply
need to compute $\cchi_{\Box}(\exp(2\pi i\frac{\rho}{\kappa}))$. Since the weights
of $V_{\Box}$ are $\pm\theta_{j}$, $j=1\ldots n$ and $\rho=\sum_{j=1}^{n}(n-j)\theta_{j}$,
we get
\begin{equation}\label{eq:vector dimension}
\begin{split}
\cchi_{\Box}(\exp(2\pi i\frac{\rho}{\kappa}))
&= \sum_{j=1}^{n}\Bigl(e^{2\pi i\frac{j-n}{\kappa}}+e^{-2\pi i\frac{j-n}{\kappa}}\Bigr)
 = 2\sum_{j=0}^{n-1}\Re e^{2\pi i\frac{j}{\kappa}}\\
&= 2\cos(\frac{(n-1)\pi}{\kappa})
    \frac{\sin(\frac{n\pi}{\kappa})}{\sin(\frac{\pi}{\kappa})}
 = 1+\frac{\sin(\frac{(2n-1)\pi}{\kappa})}{\sin(\frac{\pi}{\kappa})}
\end{split}
\end{equation}
\halmos\\

We shall need the following well--known formula, see \cite{Ka2} for instance
\begin{lemma}\label{positive dim}
If $\phi_{0}$ is defined as in proposition \ref{vector} by
$\phi_{0}(\H_{\lambda})=\cchi_{\lambda}(\exp(2\pi i\frac{\rho}{\kappa}))$, then
\begin{equation}\label{eq:qdim of reps}
\phi_{0}(\H_{\lambda})=
\prod_{\alpha>0}
\frac{\sin(\frac{\pi\<\lambda+\rho,\alpha\>}{\kappa})}
     {\sin(\frac{\pi\<\rho,\alpha\>}{\kappa})}>0
\end{equation}
\end{lemma}
\proof
By the Weyl character formula,
$\cchi_{\mu}=
 \displaystyle{
 \frac{A_{\mu+\rho}}{\prod_{\alpha>0}(e_{\half{1}\alpha}-e_{-\half{1}\alpha})}}$
where
\begin{equation}
A_{\nu}(\exp(\delta))=\sum_{w\in W}e^{\<w\nu,\delta\>}=A_{\delta}(\exp(\nu))
\end{equation}
and $e_{\beta}(\exp(\delta))=e^{\<\beta,\delta\>}$. Setting
$\mu=0$, we find Weyl's denominator formula
$A_{\rho}(\exp(2\pi it))=\prod_{\alpha>0}2i\sin(\pi\<t,\alpha\>)$. Thus,
\begin{equation}
\phi_{0}(\H_{\lambda})=
\cchi_{\lambda}(\exp(2\pi i\frac{\rho}{\kappa}))=
\frac{A_{\lambda+\rho}}
     {A_{\rho}}(\exp(2\pi i\frac{\rho}{\kappa}))=
\frac{A_{\rho}(\exp(2\pi i\frac{\lambda+\rho}{\kappa}))}
     {A_{\rho}(\exp(2\pi i\frac{\rho}{\kappa}))}=
\prod_{\alpha>0}
\frac{\sin(\frac{\pi\<\lambda+\rho,\alpha\>}{\kappa})}
     {\sin(\frac{\pi\<\rho,\alpha\>}{\kappa})}
\end{equation}
to conclude, notice simply that
\begin{equation}
\<\rho,\alpha\>\leq
\<\rho+\lambda,\alpha\>\leq
\<\rho+\lambda,\theta\><
1+\<\rho,\theta\>+\ell=\kappa
\end{equation}
and therefore \eqref{eq:qdim of reps} holds \halmos

\ssection{Main results}
\label{se:main results}

Recall the following well--known result \cite[chap. 13]{GM} whose use is typical in
subfactor theory,

\begin{theorem}[Perron--Frobenius]\label{th:PF}
Let $N$ be an irreducible matrix with non--negative entries. Then
\begin{enumerate}
\item $N$ has an eigenvector with positive entries. Up to multiplication by a positive
constant, the latter is the unique eigenvector of $N$ with non--negative entries. The
corresponding (positive) eigenvalue is called the Perron--Frobenius eigenvalue of $N$.
\item Let $\wt N$ be a matrix with non--negative entries bounded above by those
of $N$. Then, any positive eigenvalue of $\wt N$ is bounded above by the
Perron--Frobenius eigenvalue of $N$. Equality holds only if $\wt N=N$.
\end{enumerate}
\end{theorem}

\begin{theorem}\label{th:fusion rule with box}
Let $\Hbox$ be the vector representation of $L\Spin_{2n}$ at level $\ell$. Then,
for any irreducible, positive energy representation $\H_{U}$ of $L\Spin_{2n}$
whose lowest energy subspace $U$ is a single--valued representation of $\SO_{2n}$,
the following holds
\begin{equation}\label{eq:fuse with box}
\H_{U}\fuse\H_{\Box}=\bigoplus_{W\subset U\otimes V_{\Box}}N_{U\Box}^{W}\H_{W}
\end{equation}
where $N_{U\Box}^{W}$ is 1 if $W$ is admissible at level $\ell$ and 0 otherwise.
Any such $\H_{U}$ appears as a summand in one of the iterated fusion powers
$\Hbox^{\fuse k}$. In particular, it has a conjugate representation and its
quantum dimension is given by
\begin{equation}\label{eq:explicit dimension}
d(\H_{U})=
\cchi_{U}(\exp(\frac{2\pi i\rho}{\kappa}))=
\prod_{\alpha>0}
\frac{\sin(\pi\frac{\<\lambda+\rho,\alpha\>}{\kappa})}
     {\sin(\pi\frac{\<\rho,\alpha\>}{\kappa})}
\end{equation}
where $\cchi_{U}$ and $\lambda$ are the character of $U$ and its highest weight,
the $\alpha$ are the positive roots of $\Spin_{2n}$, $2\rho=\sum_{\alpha>0}\alpha$
and $\kappa=\ell+2(n-1)$.
\end{theorem}
\proof
If $\ell=1$, $\H_{U}$ is either $\H_{0}$ or $\Hbox$ and \eqref{eq:fuse with box}
follows from theorem \ref{th:level 1 ring}. The existence of a conjugate is
obvious from the fusion rules \eqref{eq:1 fuse 1} which show that the
corresponding Jones projection is 1 so that $d(\H_{U})=1$. This clearly agrees
 with \eqref{eq:explicit dimension} when $U$ is the trivial representation. When
$U=\Vbox$, \eqref{eq:vector dimension} shows that the right--hand side of
\eqref{eq:explicit dimension} is 1 since $\kappa=2n-1$ and the theorem is
proved.\\

Let now $\ell\geq 2$. By proposition \ref{pr:upper bound}, \eqref{eq:fuse with box} holds
when the $N_{U\Box}^{W}$ are replaced by some $0\leq\wt N_{U\Box}^{W}\leq N_{U\Box}^{W}$,
since part of the braiding coefficients involved in the computation of fusion
might vanish.
Let $\Nbox$ be the matrix whose entries are the $N_{U\Box}^{W}$, where $U$
and $W$ are single--valued representations of $\SO_{2n}$ which are admissible
at level $\ell$, or equivalently the diagonal block of the matrix defined in
corollary \ref{other dimension} labelled by these representations.
$\Nbox$ is non--negative and irreducible. Indeed, if $U,W$ are
two such representations with highest weights $\lambda,\mu$ so that
$\<\lambda,\theta\>,\<\mu,\theta\>\leq\ell$, there exists a sequence of
dominant integral weights of $\SO_{2n}$, $\lambda=\nu_{1}\cdots\nu_{k}=\mu$
such that $\<\nu_{j},\theta\>\leq\ell$ and each $\nu_{j+1}$ is obtained by
either adding or removing a box from the Young diagram of $\mu_{j}$. By the
tensor product rules with the vector representation $\Vbox$
(proposition \ref{ch:classification}.\ref{pr:tensor with minimal}), this
implies that $W\subseteq U\otimes\Vbox^{\otimes (k-1)}$ and therefore that
$\Nbox$ is irreducible.
By theorem \ref{th:PF}, $\Nbox$ has, up to a multiplicative constant, a
unique eigenvector with positive entries. On the other hand, by proposition
\ref{vector} and lemma \ref{positive dim}, the numbers
$\phi_{0}(\H_{U}):=\phi_{0}(\H_{\lambda})$ where $\lambda$ is the highest weight
of $U$, are positive and obey
\begin{equation}
\sum_{W}N_{\Box\medspace U}^{W}\phi_{0}(\H_{W})=
\phi_{0}(\H_{\Box})\phi_{0}(\H_{U})
\end{equation}
It follows that the Perron--Frobenius eigenvector of $\Nbox$ has entries
$\phi_{0}(\H_{\cdot})$ and corresponding eigenvalue $\phi_{0}(\Hbox)$.\\

A similar statement about the matrix $\wt N_{\Box}$ may be obtained by using the
quantum dimension $d$ of corollary \ref{co:observation} on the fusion ring $\RR$
generated by the irreducible summands of the fusion powers $\Hbox^{k}$. Applying
$d$ to \eqref{eq:fuse with box}, with $N_{U\Box}^{W}$ replaced by $\wt N_{U\Box}^{W}$
yields        
\begin{equation} \label{eq:almost pf}
\sum_{W\subset V_{\Box}\otimes U}
\wt N_{\Box\medspace U}^{W}d(\H_{W})=
d(\Hbox)d(\H_{U})
\end{equation}
whenever $\H_{U}\in\RR$. We do not know $\wt N_{\Box}$ to be irreducible, nor is
\eqref{eq:almost pf} a statement about $\wt N_{\Box}$ since some $\H_{U}$ may not lie
in $\RR$ but if $M_{\Box}$ is the matrix whose entries are $\wt N_{U\Box}^{W}$ if
$\H_{U},\H_{W}\in\RR$ and zero otherwise, we clearly have
$M_{\Box}\leq\wt N_{\Box}\leq\Nbox$ entry--wise and therefore, by theorem \ref{th:PF}
$d(\Hbox)\leq\phi_{0}(\H_{\Box})$ with equality only if $M_{\Box}=\wt N_{\Box}=\Nbox$.
However, by corollaries \ref{dimension} and \ref{other dimension}
\begin{equation}\label{eq:equality}
d(\Hbox)=
1+\frac{\sin{\frac{(2n-1)\pi}{\kappa}}}{\sin{\frac{\pi}{\kappa}}}=
\phi_{0}(\Hbox)
\end{equation}
whence $\wt N_{\Box}=\Nbox$ and \eqref{eq:fuse with box} holds. The fact that any
$\H_{U}$ appears in some $\Hbox^{\fuse k}$ now follows from the irreducibility
of $\wt N_{\Box}=\Nbox$.
Finally, the claim on the quantum dimension of the $\H_{U}$ is a direct
consequence of the uniqueness of the Perron--Frobenius eigenvector of $\Nbox$,
\eqref{eq:equality} and \eqref{eq:qdim of reps} \halmos

\begin{theorem}\label{th:closure of ordinary}
The positive energy representations of $L\Spin_{2n}$ at level $\ell$ whose lowest
energy subspace is a single--valued representation of $\SO_{2n}$ are closed under
fusion. They form a commutative and associative ring.
\end{theorem}
\proof
By theorem \ref{th:fusion rule with box}, the positive energy representations of
$L\Spin_{2n}$ corresponding to single--valued representations of $\SO_{2n}$ are
exactly the irreducible summands of the iterated fusion powers $\Hbox^{\fuse k}$
of the vector representation and therefore the generators of the commutative and
associative ring $\RR$ of section \ref{se:subring} \halmos\\

The restriction to single--valued representations of $\SO_{2n}$ in theorems
\ref{th:fusion rule with box} and \ref{th:closure of ordinary} is technical
rather than conceptual and these results conjecturally hold for all positive
energy representations of $L\Spin_{2n}$. At {\it odd} level, this can be
shown by using the discontinuous loops of section
\ref{ch:classification}.\ref{se:disc loops}.

\begin{theorem}\label{th:main 3}
Let $\Hbox$ be the vector representation of $L\Spin_{2n}$ at odd level $\ell$.
Then, for any irreducible positive energy representation $\H_{U}$ with lowest
energy subspace $U$
\begin{equation}\label{eq:extended rules}
\H_{U}\fuse\H_{\Box}=\bigoplus_{W\subset U\otimes\Vbox}N_{U\medspace\Box}^{W}\H_{W}
\end{equation}
where $N_{U\medspace\Box}^{W}$ is 1 if $W$ is admissible at level $\ell$ and
zero otherwise.
\end{theorem}
\proof
Notice that the fusion rules \eqref{eq:extended rules} for $\ell=1$ were
established in theorem \ref{th:level 1 ring}. Let now $U$ be a
single--valued $\SO_{2n}$--module so that \eqref{eq:extended rules} holds
by theorem \ref{th:fusion rule with box}. If $\zeta$ is a discontinuous
loop, proposition \ref{fuse with quasi} yields
$\zeta(\H_{U}\fuse\Hbox)\cong
(\zeta\H_{0}\fuse\H_{U})\fuse\Hbox\cong
\H_{\zeta(U)}\fuse\Hbox$
where the notation refers to the action of $Z(\Spin_{2n})$ on the irreducible
$\Spin_{2n}$--modules admissible at level $\ell$ given by proposition
\ref{ch:classification}.\ref{geometric action}. Conjugating both sides
of \eqref{eq:extended rules}, we find
\begin{equation}\label{eq:conjugated extended}
\H_{\zeta(U)}\fuse\H_{\Box}=
\bigoplus_{\substack{W\subset U\otimes\Vbox\\ W\medspace\text{$\ell$-admissible}}}
\H_{\zeta(W)}
\end{equation}
We claim that $W\subseteq U\otimes\Vbox$ if, and only if
$\zeta(W)\subseteq\zeta(U)\otimes\Vbox$. Assuming this for a moment,
\eqref{eq:conjugated extended} may be rewritten as
\begin{equation}
\H_{\zeta(U)}\fuse\H_{\Box}=
\bigoplus_{\substack{X\subset\zeta(U)\otimes\Vbox\\ X\medspace\text{$\ell$-admissible}}}
\H_{X}
\end{equation}
so that \eqref{eq:extended rules} holds for $\zeta\H_{U}$. Since $\ell$ is odd,
corollary \ref{ch:classification}.\ref{co:exchange} implies that any irreducible
representation of $L\Spin_{2n}$ at level $\ell$ is of the form $\zeta\H_{U}$
for a suitable $\zeta$ and a single--valued $\SO_{2n}$--module $U$ and theorem
\ref{th:main 3} is proved.
Returning to our claim, let $\lambda_{W},\lambda_{\Box},\lambda_{U}$ be the
highest weights of $W,\Vbox$ and $U$ respectively. By the tensor product rules
with minimal representations
(proposition \ref{ch:classification}.\ref{pr:tensor with minimal}),
$W\subseteq U\otimes\Vbox$ iff $\lambda_{W}=\lambda_{U}+w\lambda_{\Box}$ for
some $w\in W$. By proposition \ref{ch:classification}.\ref{geometric action},
$\zeta(\lambda)=w_{i}\lambda+\ell\cow{i}$, and we find that this is
the case iff $\zeta(\lambda_{W})=\zeta(\lambda_{U})+w'\lambda_{\Box}$ for
some $w'\in W$ and therefore iff $\zeta(W)\subseteq \zeta(U)\otimes\Vbox$
\halmos

\begin{theorem}\label{th:closure at odd}
At odd level, the positive energy representations of $L\Spin_{2n}$ are closed under fusion
and have conjugates. They form a commutative and associative ring.
\end{theorem}
\proof
Let $U_{i}$, $i=1,2$ be single--valued representations of $\SO_{2n}$. By theorem
\ref{th:closure of ordinary}, $\H_{U_{1}}\fuse\H_{U_{2}}$ is of positive energy.
Applying discontinuous loops $\zeta_{1},\zeta_{2}$, we find that
\begin{equation}
\zeta_{1}\H_{U_{1}}\fuse\zeta_{2}\H_{U_{2}}\cong
\zeta_{1}\H_{0}\fuse(\H_{U_{1}}\fuse\zeta_{2}\H_{U_{2}})\fuse\zeta_{2}\H_{0}\cong
\zeta_{1}\zeta_{2}(\H_{U_{1}}\fuse\zeta_{2}\H_{U_{2}})
\end{equation}
is of positive energy. By corollary \ref{ch:classification}.\ref{co:exchange},
any irreducible representation at odd level is of the form $\zeta\H_{U}$ for
some $\zeta$ and $U$ single--valued and our first claim is
proved. The existence of conjugates follows similarly for
$\H_{U}\fuse\overline{\H_{U}}\supseteq\H_{0}$ for $U$ single--valued implies that
$\zeta\H_{U}\fuse\zeta^{-1}\overline{\H_{U}}\supset\H_{0}$ and therefore
that $\zeta\H_{U}$ has a conjugate. The associativity and commutativity are
immediate consequences of the associativity of fusion and the existence of braiding
\halmos\\

\remark
It is readily verified that both sides of \eqref{eq:explicit dimension} are
invariant under the action of $Z(\Spin_{2n})$ and therefore that the quantum
dimension of all irreducible positive representations of $L\Spin_{2n}$ at odd
level is given by the formula \eqref{eq:explicit dimension}.




\newcommand {\Vspin}{V_{s_{\pm}}}
\newcommand {\Hspin}{\H_{s_{\pm}}}
\newcommand {\Hext}{\H_{\Lambda^{k}V}}
 

\chapter*{Concluding remarks}

\setcounter{section}{0}

\ssection{Computation of the fusion ring of $L\Spin_{2n}$}

The structure of the level $\ell$ fusion ring $R_{\ell}(LG)$ has been known
conjecturally for some time. It should be isomorphic to the quotient of the
representation ring $R(G)$ by a certain 'holomorphic induction' map $i_{\ell}$
taking $V$ to $\H_{V}$ whenever $V$ is admissible at level $\ell$.
One way to prove this is to show that
\begin{equation}\label{eq:to show}
i_{\ell}(V)\fuse i_{\ell}(U)=i_{\ell}(V\otimes U)
\end{equation}
for all fundamental representations $V$ since these generate $R(G)$ as a ring.
For $G=\Spin_{2n}$, they are the exterior powers $\Lambda^{k}\Vbox$ of the
vector representation $\Vbox$ with $k=1\ldots n-2$ and the spin modules
$V_{s_{\pm}}$. The fusion rule \eqref{eq:to show} for $V=\Vbox$ was proved in
chapter \ref{ch:connes fusion}. In fact, for {\it odd} $\ell$
\begin{equation}\label{fusion with box}
\Hbox\boxtimes\H_{U}=\bigoplus_{W}N_{\Box U}^{W}\H_{W}
\end{equation}
where $W$ spans the summands of $\Vbox\otimes U$ and $N_{\Box U}^{W}$ is 1 or
0 according to whether $W$ is admissible at level $\ell$ or not. For even
$\ell$, \eqref{fusion with box} was shown only for single--valued
$\SO_{2n}$--modules $U$. This restriction will be lifted in \S \ref{ss:with spin}.\\

A.~Wassermann has outlined a method for computing the remaining fusion rules,
thereby completing the calculation of $R_{\ell}(L\Spin_{2n})$. We reproduce
it below.

\ssubsection{Fusing with the spin representations $\Hspin$.}\label{ss:with spin}

Since the spin modules $\Vspin$ are minimal representations of $\Spin_{2n}$,
the corresponding Verlinde rules \eqref{eq:to show} are of the same form as
\eqref{fusion with box}, namely
\begin{equation}\label{fuse with spin}
\Hspin\boxtimes\H_{U}=
\bigoplus_{W\subset\Vspin\otimes U} N_{s_{\pm} U}^{W}\H_{W}
\end{equation}
where $N_{s_{\pm} U}^{W}$ is 1 if $W$ admissible at level $\ell$ and 0 otherwise.\\

The proof of \eqref{fuse with spin} parallels that of \eqref{fusion with box} given
in chapter \ref{ch:connes fusion}. Briefly, the braiding relations of the spin
primary fields with a general primary field are converted, by smearing on test
functions with disjoint supports and taking phases, into commutation relations
of bounded intertwiners for the local loop groups. Most of the required analysis 
has been carried out in chapter \ref{ch:sobolev fields} where all the primary
fields needed for fusing a general representation with $\Hspin$ were shown to
extend to operator--valued distributions. A technical point needs to be addressed
however. At present, the procedure of taking phases to produce bounded intertwiners
applies only if one of the two primary fields extends to a bounded
operator--valued distribution. This is the case for the vector primary field,
since it is essentially a Fermi field but does not hold for the spin fields. Once
this problem is bypassed, the bounded intertwiners give an $L\Spin_{2n}$--equivariant
isometry of the left hand--side of \eqref{fuse with spin} into the right hand--side
and therefore the upper bound
\begin{equation}\label{bound}
0\leq\wt N_{s_{\pm}U}^{W}\leq N_{s_{\pm}U}^{W}
\end{equation}
where the $\wt N_{s_{\pm} U}^{W}$ are the structure constants of fusion. This shows
that fusion with the spin representations is always of positive energy and that these
representations have (necessarily unique) conjugates.
Indeed, by \eqref{bound}, $\Hspin\boxtimes\Hspin$ contains only representations whose
lowest energy subspaces are single--valued $\SO_{2n}$--modules. As shown in chapter
\ref{ch:connes fusion}, these representations have conjugate bimodules so that for
some $\H_{W}$, 
\begin{equation}
\Hspin\boxtimes(\Hspin\boxtimes\H_{W})=
(\Hspin\boxtimes\Hspin)\boxtimes\H_{W}
\supset\H_{0}
\end{equation}
whence $\H_{s_{\pm}}$ have conjugates.\\

We now use the method of chapter \ref{ch:connes fusion} to prove that
\eqref{fusion with box} holds for any representation $U$ of $\Spin_{2n}$ admissible
at level $\ell$. Using the Doplicher--Haag--Robert theory, the quantum dimension $d$
may be extended from the ring $\RR$ generated by the iterated fusion powers of $\Hbox$
to the ring $\wt\RR$ generated by the fusion powers of $\Hbox$ and $\Hspin$. Consider
the matrix $N_{\Box}$, whose entries $N_{\Box U}^{W}$ are now indexed by two--valued
representations of $\SO_{2n}$. $N_{\Box}$ is non--negative, irreducible and, by simple
inspection, has a Perron--Frobenius eigenvector whose entries are the 'dimensions'
\begin{equation}
\phi_{0}(\H_{U})=\cchi_{U}(\exp(\frac{2\pi i\rho}{\kappa}))
\end{equation}
given in \cite[\S 13.8]{Ka1} and studied in chapter \ref{ch:connes fusion}. The
corresponding eigenvalue is $\phi_{0}(\Hbox)$. The fusion rules with $\Hbox$ read
\begin{equation}
\Hbox\fuse\H_{U}=\bigoplus_{W\subset\Vbox\otimes U}\wt N_{\Box U}^{W}\H_{W}
\end{equation}
where, as shown in chapter \ref{ch:connes fusion}, $\wt N_{\Box U}^{W}\leq N_{\Box U}^{W}$
with equality if $U$ is a single--valued $\SO_{2n}$--module.
Let the matrix $M_{\Box}$ have entries $M_{\Box U}^{W}=\wt N_{\Box U}^{W}$ if $U,W\in\wt\RR$
and $0$ otherwise so that $M_{\Box}\underline{d}=d(\H_{0})\underline{d}$ where $d$ is
the (non--zero) vector with components $d(\H_{U})$. Since
$M_{\Box}\leq\wt N_{\Box}\leq N_{\Box}$ and
$d(\H_{0})=\phi_{0}(\Hbox)$, we have $M_{\Box}=N_{\Box}$. Therefore \eqref{fusion with box}
holds and all positive energy representations of level $\ell$ lie in $\wt\RR$. Moreover,
by uniqueness of Perron--Frobenius eigenvectors, the quantum dimension of those
corresponding to two--valued representations of $\SO_{2n}$ is given by
$d(\H_{U})=\gamma \phi_{0}(\H_{U})$ for some $\gamma>0$ while, by the results of chapter
\ref{ch:connes fusion}, $d(\H_{U})=\phi_{0}(\H_{U})$ whenever $U$ is a single--valued
$\SO_{2n}$--module.\\

Let now $U$ be a single--valued representation of $\SO_{2n}$, then
\begin{equation} \label{fufuse with spin}
\Hspin\fuse\H_{U}=\bigoplus_{W\subset\Vspin\otimes U}\wt N_{s_{\pm}U}^{W}\H_{W}
\end{equation}
where all $W$ are two--valued on $\SO_{2n}$. Since $d(\Hspin)=\gamma \phi_{0}(\Hspin)$,
$d(\H_{U})=\phi_{0}(\H_{U})$ and $d(\H_{W})=\gamma \phi_{0}(\H_{W})$, applying $d$ to
\eqref{fufuse with spin} we get
\begin{equation} \label{ecco}
\sum_{W\subset\Vspin\otimes U}\wt N_{s_{\pm}U}^{W}\phi_{0}(\H_{W})=
\phi_{0}(\Hspin)\phi_{0}(\H_{U})
\end{equation}
On the other hand, the dimensions $\phi_{0}(\H_{\cdot})$ satisfy \eqref{ecco} when
the $\wt N_{s_{\pm}U}^{W}$ are replaced by $N_{s_{\pm}U}^{W}$ and therefore the
two numbers must be equal by \eqref{bound}. It follows that \eqref{fuse with spin}
holds whenever $U$ is a single--valued representation of $\SO_{2n}$.
To conclude, assume that the conjugate bimodules of $\H_{s}$ is
$\overline{\H_{s}}=\H_{\overline{s}}$ where $s$ is $s_{\pm}$ and $\overline{s}$ is the
corresponding dual highest weight, so that $\overline{s}$ is $s_{\pm}$ if $n$ is even
and $s_{\mp}$ if $n$ is odd. Then, if $U$ is a single--valued $\SO_{2n}$--module, we
get by Frobenius reciprocity
\begin{equation}
N_{s U}^{W}=
\wt N_{s U}^{W}=
\dim\Hom_{L\Spin_{2n}}(\H_{s}\fuse\H_{U},\H_{W})=
\dim\Hom_{L\Spin_{2n}}(\H_{U},\H_{\overline{s}}\fuse\H_{W})=
\wt N_{\overline{s} W}^{U}
\end{equation}
and therefore, using the symmetry of the Verlinde numbers,
$\wt N_{\overline{s} W}^{U}=N_{s U}^{W}=N_{\overline{s} W}^{U}$ which concludes the
proof of \eqref{fuse with spin}. To summarise, the above method relies on
\begin{enumerate}
\item Extending the procedure of taking phases to produce the braiding of
bounded intertwiners from that of primary fields to the case where both
primary fields are unbounded operator--valued distributions.
\item Proving that $\overline{\H_{s}}=\H_{\overline{s}}$.
\end{enumerate}

\ssubsection{Fusing with the exterior powers $\Hext$.}

The fusion rules with the exterior powers $\Hext$, $k=2\ldots n-2$ are more delicate
to handle since they involve multiplicities other than 0 or 1. This difficulty can
be circumvented by realising the multiplicity spaces as representations of the
Birman--Wenzl algebra $BW_{k}$ on $k$ strings.
The latter is relevant because a generalisation of Brauer--Weyl duality \cite{Br} holds,
namely the commutant of $L\Spin_{2n}$ on the $k$--fold fusion $\Hbox^{\boxtimes k}$ is
generated by the action of the braid group $B_{k}$ and this action factors through
$BW_{k}$. Moreover, $\Hext$ is the isotypical summand of $\Hbox^{\boxtimes k}$
corresponding to the quantum sign representation $\varepsilon_{k}$ of $B_{k}$.
Thus, $\Hext\fuse\H_{U}$ may be computed as the subspace of $\Hbox^{\fuse k}\fuse\H_{U}$
transforming like $\varepsilon_{k}$ under the natural action of $B_{k}$. This may be achieved
by realising the action more explicitly, by analogy with the $L\SU_{n}$ case treated in
\cite{Wa3} on chains of $k$ vector primary fields of the form
\begin{equation}
\phi_{V_{k}V_{k-1}}(f_{k})\cdots\phi_{V_{1}U}(f_{1})
\end{equation}
smeared on functions supported in $k$ consecutive subintervals of $I=(0,2\pi)$,
where the action of $B_{k}$ is given by braiding. As for the previously
computed fusion rules, this analysis gives an upper bound for $\Hext\fuse\H_{U}$
in terms of the Verlinde rules. A Perron--Frobenius argument based on the computation
of the quantum dimension of $\Hext$ given in chapter \ref{ch:connes fusion} then
shows that the bound is attained.

\ssection{Loop groups of other Lie groups}

The loop groups of Lie groups other than $\SU_{n}$ and $\Spin_{2n}$ need to be
investigated. Our results extend almost {\it verbatim} to $\Spin_{2n+1}$. Moreover,
the analysis of the continuity properties of primary fields done in chapters
\ref{ch:vertex operator} and \ref{ch:sobolev fields} covers the level 1 and a number
of higher level fields for all simply--laced groups. The primary fields for the
remaining $\Sp_{n},F_{4},G_{2}$ may probably be studied via the generalised vertex
operator construction \cite{GNOS}.

\ssection{Loop groups, quantum invariant theory and subfactors}

Any positive energy representation $(\pi,\H)$ of $LG$ gives rise to two subfactors.
The first is obtained via the quantum invariant theory inclusion
\begin{equation}\label{II1}
(\bigcup_{m} \IC\otimes\End_{LG}(\H^{\fuse m}))''\subset
(\bigcup_{m}           \End_{LG}(\H^{\fuse m}))''
\end{equation}
and is of type II$_{1}$ with the trace given by the Doplicher--Haag--Roberts theory.
The second is the Jones--Wassermann inclusion
\begin{equation}\label{III1}
\pi(L_{I}G)''\subset\pi(L_{I^{c}}G)'
\end{equation}
where $L_{I}G$ is the group of loops supported in the proper interval $I\subset S^{1}$.
As shown by Jones and Wassermann \cite{Wa1} and explained in chapter \ref{ch:loc loops},
\eqref{III1} is an irreducible inclusion of hyperfinite factors of type III$_{1}$ whenever
$\pi$ is irreducible. Using results of Popa \cite{Po1,Po2},
Wassermann proved that for $G=\SU_{n}$, \eqref{III1} is isomorphic to the tensor product
of the hyperfinite factor of type III$_{1}$ with a Jones-Wenzl subfactor \cite{Wa1}. The
latter is defined using the Hecke algebras $H_{m}(q)$ with $q=e^{\pi i/(N+\ell)}$ and the
representation of $\SU_{n}$ corresponding to the lowest energy subspace of $\H$ \cite{We1}.
Interestingly, it may also be defined as the quantum invariant theory inclusion 
corresponding to $\H$. Since the positive energy representations of $LG$ and the subfactors
defined in \cite{We2} are characterised by the same data for $G$ of type $B$, $C$ or $D$,
a natural conjecture is that a similar result holds in these cases.
For $G=\Spin_{2n}$, this should follow from a computation of the higher relative
commutants by showing that $\End_{L\Spin_{2n}}(\Hbox^{\fuse n})$ is isomorphic to the
Birman-Wenzl algebra on $n$ strings.
In fact, only the II$_{1}$ subfactors corresponding to the single--valued representations
of $SO_{2n}$ are constructed in \cite{We2}. The quantum invariant theory
inclusion on the other hand gives subfactors for {\it any} representation of $\Spin_{2n}$.
We intend to study the corresponding spin subfactors, the structure of which should be
linked to the delicate combinatorics of the tensor product rules for the two--valued
representations of $\Spin_{2n}$.

\ssection{Modular categories and invariants of $3$-manifolds}

An important class of invariants of $3$--manifolds is obtained via surgery on
framed links in the $3$-sphere. These use {\it modular tensor categories} which
are traditionally produced as the semi--simple quotients of the representation theory
of quantum groups at roots of unity \cite{Tu}. The analysis involved relies on the
combinatorics of the tensor products of the representations of the underlying Lie
algebras \cite{TuWe}.
So far, two main difficulties have prevented the construction of the invariants
corresponding to $\Spin_{2n}$. On the one hand, the intricacies of the combinatorics
of two--valued representations of $\SO_{2n}$ have led researchers to consider the
(truncated) representation ring of $\SO_{2n}$ only. On the other, the latter gives
rise to a braided tensor category which is not modular, in that the corresponding
$S$ matrix is not invertible. As a result, only the invariants corresponding to
$\PSO_{2n}=\SO_{2n}/\IZ_{2}$ have been defined.
The category $\pl$ of positive energy representations of $L\Spin_{2n}$ at level
$\ell$, endowed with Connes fusion should provide a natural solution to this problem.
Indeed, the Doplicher--Haag--Roberts theory implies that $\pl$ is a braided tensor
category and there remains to verify that $\pl$ is modular.

\ssection{Knizhnik--Zamolodchikov equations}

In chapter \ref{ch:box/kbox braiding}, we solved the KZ differential equations
corresponding to a number of important fusion rules for $L\Spin_{2n}$ by using
the contour integrals of Dotsenko and Fateev \cite{DF}. These have been greatly
generalised by Schechtman--Varchenko \cite{SV} and Feigin--Frenkel \cite{EdF}
who used them systematically to give all solutions of the KZ equations. It is
interesting to notice however that their generalised hypergeometric solutions,
which first appeared in the work of Aomoto and Gelfand \cite{Ao,GKZ}, use
Euler--like contour integrals with a large number of integration variables 
(typically growing linearly in $n$) which make them intractable for computational
purposes. Our solutions on the other hand use 2 integration variables for any $n$.
There might therefore exist a universal cohomological simplification mechanism
applicable to the contour integrals of Schechtman--Varchenko and Feigin--Frenkel
which could well lead to a proof of the symmetry properties of the braiding
coefficents conjectured by Witten in the context of restricted solid--on--solid
statistical mechanical models.


\end{document}